\documentclass[11pt]{amsart}
\usepackage[alphabetic]{amsrefs} 
\usepackage[letterpaper,margin=1in]{geometry}

\usepackage{amsmath}
\usepackage{amssymb,mathrsfs,url}
\usepackage{amsfonts}
\usepackage{amsthm}
\usepackage[utf8]{inputenc}
\usepackage[pdfencoding=auto, psdextra,colorlinks=true,linkcolor=blue, citecolor=teal!50!green,urlcolor=cyan]{hyperref}
\usepackage{verbatim}

\usepackage{tikz-cd}
\usepackage{array,booktabs}
\usepackage{xspace}
\usepackage[all,cmtip]{xy}
\usepackage{dynkin-diagrams}
\usepackage{enumerate}
\usepackage{bm}
\usepackage{enumitem}
\usepackage{xcolor}
\usepackage{xr}
\usepackage{multicol}

\usepackage{relsize}
\usepackage[bbgreekl]{mathbbol}
\usepackage{amsfonts}
\DeclareSymbolFontAlphabet{\mathbb}{AMSb} 
\DeclareSymbolFontAlphabet{\mathbbl}{bbold}
\newcommand{\Prism}{{\mathlarger{\mathbbl{\Delta}}}}
\newcommand{\N}{\ensuremath{\mathbb{N}}} 
\newcommand{\Z}{\ensuremath{\mathbb{Z}}} 
\newcommand{\bZ}{\ensuremath{\breve{\Z}}} 
\newcommand{\bE}{\ensuremath{\breve{E}}} 
\newcommand{\Qp}{\ensuremath{\mathbb{Q}_p}} 
\newcommand{\rQ}{\ensuremath{\mathbb{Q}}} 
\newcommand{\bQ}{\ensuremath{\breve{\rQ}_p}} 

\newcommand{\R}{\ensuremath{\mathbb{R}}} 

\newcommand{\CC}{\ensuremath{\mathbb{C}}} 
\newcommand{\F}{\ensuremath{\mathbb{F}}} 

\renewcommand{\emptyset}{\ensuremath{\varnothing}}	

\newcommand{\pre}{\mathrm{pre}}
\newcommand{\Res}{\operatorname{Res}}

\newcommand{\Rep}{\operatorname{Rep}} 
\newcommand{\perf}{\ensuremath{\textrm{perf}}} 
\newcommand{\Isoc}{\operatorname{Isoc}}
\newcommand{\Stab}{\ensuremath{\mathrm{Stab}}} 

\newcommand{\Lie}{\operatorname{Lie}}
\newcommand{\twolim}{\operatornamewithlimits{2-lim}} 
\newcommand{\Spec}{\operatorname{Spec}}

\newcommand{\Gal}{\operatorname{Gal}}  
\newcommand{\Frob}{\operatorname{Frob}} 
\newcommand{\OO}{\ensuremath{\mathcal{O}}} 

\newcommand{\lrangle}[1]{\ensuremath{\langle #1 \rangle}} 
\newcommand{\lrbracket}[1]{\ensuremath{\{ #1 \}}} 
\newcommand{\wdt}[1]{\ensuremath{\widetilde{#1}}} 
\newcommand{\wdh}[1]{\ensuremath{\widehat{#1}}} 
\newcommand{\ovl}[1]{\ensuremath{\overline{#1}}} 
\newcommand{\dottimes}{\operatorname{\dot{\times}}} 

\newcommand{\identity}{\ensuremath{\mathrm{id}}}

\newcommand{\rightiso}{\ensuremath{\stackrel{\sim}{\rightarrow}}}

\newcommand{\Isom}{\operatorname{\underline{\mathrm{Isom}}}}

\newcommand{\intg}{\ensuremath{\mathrm{int}}}

\newcommand{\Int}{\operatorname{\mathrm{Int}}}

\newcommand{\EE}{\ensuremath{\mathcal{E}}}
\newcommand{\Ker}{\operatorname{ker}}

\newcommand{\Image}{\operatorname{im}}

\newcommand{\mT}{\ensuremath{T}} 
\newcommand{\mG}{\ensuremath{G}} 
 
\newcommand{\mQ}{\ensuremath{\mathcal{Q}}}

\newcommand{\mult}{\ensuremath{\mathrm{mult}}} 

\newcommand{\dR}{\ensuremath{\mathrm{dR}}} 
\newcommand{\st}{\ensuremath{\mathrm{st}}} 
\newcommand{\proet}{\text{pro{\'e}t}} 
\newcommand{\fket}{\text{fk{\'e}t}} 
\newcommand{\ket}{\text{k{\'e}t}} 
\newcommand{\proket}{\text{prok{\'e}t}} 

\newcommand{\Spf}{\operatorname{Spf}}
\newcommand{\Spa}{\operatorname{Spa}}
\newcommand{\Spd}{\operatorname{Spd}}

\newcommand{\spe}{\operatorname{sp}} 
\newcommand{\Dia}{\ensuremath{\Diamond}} 
\newcommand{\dia}{\ensuremath{\diamond}} 
\newcommand{\Perf}{\operatorname{Perf}} 
\newcommand{\Perfd}{\operatorname{Perfd}}  
\newcommand{\PCAlg}{\operatorname{PCAlg^{op}}} 
\newcommand{\SchPerf}{\operatorname{SchPerf}}  

\newcommand{\gp}{\ensuremath{\mathrm{gp}}} 
\newcommand{\sat}{\ensuremath{\mathrm{sat}}} 
\newcommand{\md}[1]{\mathsf{#1}} 

\newcommand{\Bun}{\operatorname{Bun}} 
\newcommand{\BL}{\operatorname{BL}} 
\newcommand{\Adm}{\operatorname{Adm}} 

\newcommand{\Gra}[1]{\operatorname{Gr}_{#1}} 
\newcommand{\lcM}{\ensuremath{M^{\loc}}} 
\newcommand{\vM}{\ensuremath{\mathbb{M}}} 
\newcommand{\fl}{\ensuremath{\mathcal{F}l}} 
\newcommand{\Fil}{\ensuremath{\mathrm{Fil}}} 
\newcommand{\loc}{\ensuremath{\mathrm{loc}}} 
\newcommand{\Sht}{\ensuremath{\mathrm{Sht}}} 
\newcommand{\HT}{\ensuremath{\mathrm{HT}}} 
\newcommand{\ls}{\ensuremath{\mathbb{L}}} 
\newcommand{\lp}{\ensuremath{\mathbb{P}}} 
\newcommand{\VVs}{\ensuremath{\mathscr{V}}} 
\newcommand{\PPs}{\ensuremath{\mathscr{P}}} 
\newcommand{\PPp}{\ensuremath{\mathbb{P}}} 
\newcommand{\DRT}{\operatorname{DRT}} 
\newcommand{\Loc}{\ensuremath{\mathrm{Loc}}} 

\newcommand{\Gm}{\ensuremath{\mathbb{G}_\mathrm{m}}}
\newcommand{\Ga}{\ensuremath{\mathbb{G}_\mathrm{a}}}
\newcommand{\GL}{\operatorname{GL}}

\newcommand{\GSp}{\operatorname{GSp}}
\newcommand{\GSP}{\operatorname{\mathcal{GSP}}}

\newcommand{\cl}{\ensuremath{\mathrm{cl}}} 
\newcommand{\diag}{\operatorname{diag}} 

\newcommand{\DS}{\ensuremath{\mathbb{S}}} 
\newcommand{\Shum}[1]{\ensuremath{\mathscr{S}_{#1}}} 
\newcommand{\shu}[1]{\ensuremath{\mathrm{Sh}_{#1}}} 
\newcommand{\shuc}[2]{\ensuremath{\mathrm{Sh}_{#1}^{#2}}}

\newcommand{\Shumc}[2]{\ensuremath{\mathscr{S}_{#1}^{#2}}}
\newcommand{\Shumm}[1]{\ensuremath{\mathscr{S}_{#1}^{\min}}}
\newcommand{\Cusp}{\operatorname{Cusp}} 

\newcommand{\Zb}{\ensuremath{\mathrm{Z}}} 

\newcommand{\NE}{\ensuremath{\mathcal{N}}} 
\newcommand{\CE}{\ensuremath{\mathcal{C}}} 
\newcommand{\KR}{\ensuremath{\mathrm{KR}}} 
\newcommand{\EKOR}{\ensuremath{\mathrm{EKOR}}} 
\newcommand{\ZIP}{\ensuremath{\mathrm{Zip}}} 

\newcommand{\MM}{\ensuremath{\mathcal{M}}}

\newcommand{\FFC}{\ensuremath{\mathrm{FF}}}

\newcommand{\Bui}{\operatorname{\mathcal{B}}} 
\newcommand{\FF}{\ensuremath{\mathcal{F}}} 

\newcommand{\QQc}{\ensuremath{\mathcal{Q}^{\circ}}}

\newcommand{\PPc}{\ensuremath{\mathcal{P}^{\circ}}}
\newcommand{\GGc}{\ensuremath{\mathcal{G}^{\circ}}} 

\newcommand{\KK}{\ensuremath{\Breve{K}}} 

\newcommand{\ab}{\ensuremath{\mathcal{A}}} 
\newcommand{\UU}{\ensuremath{\mathcal{U}}} 
\newcommand{\dd}{\ensuremath{\mathrm{\ddagger}}} 

\newcommand{\GG}{\ensuremath{\mathcal{G}}} 
\newcommand{\HH}{\ensuremath{\mathcal{H}}} 
\newcommand{\LL}{\ensuremath{\mathcal{L}}}
\newcommand{\PP}{\ensuremath{\mathcal{P}}} 
\newcommand{\QQ}{\ensuremath{\mathcal{Q}}} 

\newcommand{\bigsur}{\ensuremath{\star}} 

\newcommand{\bb}[1]{\mathbb{#1}}

\newcommand{\ca}[1]{\mathcal{#1}}

\newcommand{\und}[1]{\und{#1}}
\newcommand{\sbst}{\subset}

\newcommand{\lra}{\longrightarrow}

\newcommand{\mbf}[1]{\mathbf{#1}}

\newcommand{\iso}{\cong}
\newcommand{\mrm}[1]{\mathrm{#1}}
\newcommand{\wat}{\widehat}

\newcommand{\et}{\text{{\'e}t}}

\newcommand{\dr}{\mathrm{dR}}

\newcommand{\der}{\mrm{der}}
\newcommand{\ad}{\mrm{ad}}

\newcommand{\sh}{\mrm{Sh}}
\newcommand{\wdtd}{\widetilde}

\newcommand{\A}{\bb{A}_f}
\newcommand{\bss}{\backslash}
\newcommand{\disju}{\coprod}
\newcommand{\stb}{\mrm{Stab}}

\DeclareMathOperator{\lie}{\mrm{Lie}}

\newcommand{\mmin}{\mrm{min}}
\newcommand{\cusp}{\mrm{Cusp}}

\newcommand{\Ap}{\bb{A}_f^p}

\DeclareMathOperator{\spec}{\mrm{Spec}}
\newcommand{\gal}{\mrm{Gal}}

\newcommand{\Aut}{\mrm{Aut}}

\newcommand{\Hom}{\mrm{Hom}}

\newcommand{\G}{\ca{G}}

\DeclareMathOperator{\gr}{\mrm{Gr}}

\DeclareMathOperator{\im}{\mrm{Im}}

\newcommand{\ul}[1]{\underline{#1}}

\usepackage[OT2,T1]{fontenc}
\DeclareSymbolFont{cyrletters}{OT2}{wncyr}{m}{n}
\DeclareMathSymbol{\Sha}{\mathalpha}{cyrletters}{"58}

\newcommand{\Ad}{\mrm{Ad}}

\newcommand{\can}{\mrm{can}}

\DeclareMathOperator{\ob}{\mrm{Ob}}

\newcommand{\red}{\mrm{red}}

\newcommand{\cpl}[2]{{(#1)}{}^{\wedge{\ }}_{#2}}

\newcommand{\K}{\wdtd{K}}

\newtheorem{lem}{Lemma}[section]
\newtheorem{definition}[lem]{Definition}
\newtheorem{thm}[lem]{Theorem}
\newtheorem{rk}[lem]{Remark}
\newtheorem{prop}[lem]{Proposition}

\newtheorem{cor}[lem]{Corollary}
\newtheorem{convention}[lem]{Convention}
\newtheorem{construction}[lem]{Construction}
\newtheorem{assumption}[lem]{Assumption}
\newtheorem{axiom}[lem]{Axiom}
\newtheorem{Definition and Lemma}[lem]{Definition/Lemma}
\newtheorem{Definition and Proposition}[lem]{Definition/Proposition}
\numberwithin{equation}{section}

\newtheorem*{rkS}{Remark}
\newtheorem*{propS}{Proposition}

\newtheorem*{thmA}{Theorem A}
\newtheorem*{thmB}{Theorem B}
\newtheorem*{thmC}{Theorem C}
\newtheorem*{thmD}{Theorem D}
\newtheorem*{thmE}{Theorem E}
\newtheorem*{thmF}{Theorem F}
\newtheorem*{corA}{Corollary A}
\newtheorem*{corB}{Corollary B}

\newtheorem*{defA}{Definition A}
\newtheorem*{defB}{Definition B}
\newenvironment{proofof}[1][]{{\noindent\it Proof of  {#1}.}\hspace{.5em}}{\hfill $\square$\par}

\title[Canonical extension of $p$-adic shtuka]{Canonical extensions of $p$-adic shtukas on toroidal compactifications of Shimura varieties}
\author{Shengkai Mao}
\address{Morningside Center of Mathematics, Beijing 100190, China}
\email{maoshengkai@amss.ac.cn}
\author{Peihang Wu}
\address{BICMR, Peking University, Beijing 100871, China}
\email{wuph@pku.edu.cn} 
\date{\today}

\begin{document}
\begin{abstract}
We construct canonical extensions of $p$-adic shtukas on integral models of toroidal compactifications of abelian-type Shimura varieties with quasi-parahoric levels at any prime number $p$. 
More precisely, we define the notion of a log diamond as a $v$-sheaf associated with a log scheme over $\bb{Z}_p$ and construct a $p$-adic log shtuka over the log diamond of an integral toroidal compactification of an abelian-type Shimura variety by studying the ``degeneration'' of the shtuka at the boundary. 
Moreover, we provide a definition of canonical integral models of toroidal and minimal compactifications in the sense of Pappas and Rapoport, and verify it in the same generality as above.\par
Applications include the canonicity and functoriality of integral toroidal compactifications, as well as an axiomatic proof of the well‑positionedness of all well‑known stratifications on the special fiber.
\end{abstract}
\maketitle
\tableofcontents

\section*{Introduction}
Let $(G,X)$ be a Shimura datum with reflex field $\bb{E}(G,X)$. Denote by $\{\shu K(G,X)\}_{K\subset G(\A),\text{neat}}$ the associated Shimura varieties. Fix a prime number $p$ and a place $v|p$ of $\bb{E}$.\par
In \cite{PR24}, building on Scholze’s idea and theory of $p$-adic shtukas as substitutes for ``motives'' over $p$-adic fields, Pappas and Rapoport developed a general theory of ($p$-adic) shtukas and applied it to local and global Shimura varieties. They constructed canonical integral models over $\Spec\OO_{\bb{E}_v}$, proving their uniqueness and functoriality, and treated Hodge-type Shimura varieties with stabilizer parahoric level. This was extended to general parahoric levels in \cite{daniels2024conjecture}, and to abelian-type Shimura varieties in \cite{DY25} using the Kisin–Pappas–Zhou models \cite{KPZ24}.\par
Toroidal compactifications of Shimura varieties over the complex numbers were constructed in \cite{AMRT10}, and for mixed Shimura varieties over reflex fields in \cite{Pin89}. For a cone decomposition $\Sigma$, the compactification $\shuc K\Sigma(G,X)$ has boundary charts described by certain toric embeddings of mixed Shimura varieties indexed by cusp labels $[(\Phi,\sigma)]\in\Cusp_K(G,X,\Sigma)$.\par
Integral versions of this picture were developed in \cite{FC90}, \cite{Lan13}, \cite{Mad19}, and \cite{Wu25}. In the PEL-type setting, boundary strata and their gluings are described in terms of degenerations of abelian schemes. For Shimura varieties of Hodge type and of abelian type, the strategy is different: one first obtains candidate integral models from integral toroidal compactifications of Siegel type via normalizations and passing to quotients, and the main technical difficulties arise in verifying that the boundary components of these models satisfy the required structural properties.\par
A natural question thus arises:\par 
\textbf{Q}: If one interprets the integral models constructed in the framework of \cite{PR24} as moduli spaces of shtukas, then for suitably defined integral models of their toroidal compactifications, one should, at least in principle, be able to provide a precise description of their degeneration behavior along the boundary.\par 
In this paper, we discuss our understanding of this question. The goals of this paper are threefold.
\begin{enumerate}
\item We develop a theory of log shtukas over log diamonds over $\Spd\Z_p$, extending the theory of shtukas over diamonds in \cite{PR24}. On the generic fiber, log shtukas correspond to certain pro-Kummer étale local systems with Hodge–Tate period maps, while on the special fiber, they recover the usual shtukas. We also prove the rigidity of extensions from the generic fiber to the integral base.

\item We construct log shtukas on toroidal compactifications of integral models of abelian-type Shimura varieties with quasi-parahoric level (which slightly generalizes Kisin--Pappas--Zhou models), and provide an axiomatic framework that is applicable to general Shimura varieties, in which  
the extension of shtukas to toroidal compactifications does not rely on the theory of abelian schemes or $p$-divisible groups. Hence, we expect that this picture should also hold in general, although canonical integral models have not yet been constructed for general parahoric groups at any $p$. In particular, this allows
us to prove the uniqueness and functoriality of toroidal compactifications of integral models, including compatibility with changes of parahoric levels, generalizing results in \cite{Wu25} and
\cite{Mao25b}.

\item We provide a description of how the shtukas supported on integral models degenerate along the boundary. Along the way, we show that various stratifications on the special fibers, including central leaves, Newton strata, Kottwitz–Rapoport (KR)
strata, and Ekedahl–Kottwitz–Oort–Rapoport (EKOR) strata defined using shtukas have good properties along the boundary. This reproves and generalizes some results in \cite{boxer2015torsion}, \cite{lan2018compactifications} and \cite{Mao25}.
\end{enumerate}

\subsubsection*{Log shtukas on log diamonds}
Let $E$ be a local field of characteristic $0$, let $G$ be a reductive group over $\rQ_p$, let $\mu$ be a minuscule cocharacter, and let $\GG$ be a parahoric group scheme of $G$ over $\Z_p$.\par
Recall that, in \cite{PR24}, a theory of shtukas over diamonds was developed. 
Given a scheme $X$ that is separated, of finite type over $\Z_p$, there are different ways to attach certain diamonds to $X$, e.g., $X^{\Dia}$, $X^{\dia}$, and $X^{\Dia/}$ (see \cite[\S 2.2]{anschutz2022p} and \cite[\S 2.1]{PR24}).
One can define a family of shtukas $(\PPs, \phi_{\PPs})$ on $X^{\Dia}$ (with one leg bounded by $\mu$) as a $1$-morphism between $v$-stacks $X^{\Dia} \to \Sht_{\GG, \mu}$. 
For a locally noetherian adic space $X$ over $\Spa(E,\OO_E)$, the authors proved that there is an equivalence between $\GG$-shtukas over $X^\Diamond$ and pairs $(\bb{P},H)$ consisting of pro-{\'e}tale $\ul{\ca{G}(\bb{Z}_p)}$-torsors $\bb{P}$ over $X^\Diamond$ and Hodge-Tate maps $H$ (see \cite[Prop. 2.5.3]{PR24}). 
An important feature is that, given $\mathscr{X}$ a separated normal scheme of finite type over $\OO_E$, a morphism of $\GG$-shtukas on the generic fiber $(\mathscr{X}_{\eta})^{\Dia}$ extends uniquely over the integral base $\mathscr{X}^{\Dia}$ (see \cite[Thm. 2.7.7 and Cor. 2.7.9]{PR24}). In characteristic $p$, shtukas are crystalline in nature (see \cite[Thm. 2.3.8 and Ex. 2.4.9]{PR24}).

To define the extension of shtukas on toroidal compactifications, we need to generalize the theory of diamonds and shtukas to the logarithmic setting. 
Let $(X,\ca{M}_X)$ be a locally Noetherian fs log adic space over $\Spa(\bb{Z}_p)$. Recall that in \cite{DLLZ23}, a log structure $(\ca{M}_X, \alpha_X)$ on $X$ is a pair consisting of $\ca{M}_X$, a sheaf of monoids over $X_{\et}$, together with a morphism of sheaves of monoids $\alpha_X: \ca{M}_X \to \OO_{X_{\et}}$ such that $\alpha_X^{-1}(\OO_{X_{\et}}^{\times}) \to \OO_{X_{\et}}^{\times}$ is an isomorphism.

Recall that the diamond $X^\Diamond$ associated with $X$ is a $v$-sheaf, and we should associate a $v$-sheafy diamond with the log adic space $(X, \ca{M}_X)$ while preserving the log-structure information coming from the {\'e}tale topology.

\begin{defA}[{Definition \ref{def-log-diamonds-adic-space}}]
Let $S \in \Perf$. Then $(X,\ca{M}_X)^{\Dia}(S)$ consists of isomorphism classes of pairs $((S^{\sharp}, \ca{M}_{S^{\sharp}}), f)$, where $S^{\sharp}$ is an untilt of $S$, $\ca{M}_{S^{\sharp}}$ is a \emph{saturated and fine perfectoid} log structure on $S^{\sharp}$, and $f: (S^{\sharp}, \ca{M}_{S^{\sharp}}) \to (X, \ca{M}_X)$ is a morphism between log adic spaces.
\end{defA}

Here the adjective ``\emph{fine perfectoid}'' will be introduced in Definition \ref{def-perfectoid-log-structure}. Roughly speaking, we require that $\ca{M}_{S^{\sharp}}$ has local charts that are uniquely $p$-divisible and are ``perfections'' of finitely generated monoids in an appropriate sense. We show that
\begin{thmA}[{Theorem \ref{thm-log-diamond-v-sheaves}}]
$(X,\ca{M}_X)^{\Dia}$ is a $v$-sheaf on $\Perf$.
\end{thmA}
For separated schemes of finite type, one can define log diamonds $(X, \ca{M}_X)^{\Dia}$, $(X, \ca{M}_X)^{\dia}$, and $(X, \ca{M}_X)^{\Dia/}$ as $v$-sheaves in a similar way. 
\begin{defB}[{Definition \ref{def-log-shtuka}}]
The groupoid of log $\GG$-shtukas on $(X,\ca{M}_X)^\Dia$ (with one leg bounded by $\mu$) is defined as (see Definition \ref{def-2-lim})
$$\Sht_{\GG, \mu}(X,\ca{M}_X):=\twolim\limits_{(S^{\sharp},\ca{M}_{S^\sharp},f_{S^\sharp})\in ((X,\ca{M}_X)^{\Dia})^{\mrm{op}}}\Sht_{\GG, \mu}(S^\sharp).$$
  \end{defB}
In other words, giving a log shtuka amounts to giving a $1$-morphism between $v$-stacks $$(X, \ca{M}_X)^{\Dia} \to \Sht_{\GG, \mu}.$$

Let us now focus on the generic fiber. In this situation, $X$ is a locally Noetherian fs log adic space over $\Spa \rQ_p$.
We consider the category of pro-Kummer-{\'e}tale $\Z_p$-local systems on $X$, which is the generalization of the category of pro-{\'e}tale $\Z_p$-local systems to the log case. To formulate an equivalence of categories, we need to restrict the notion to a so-called \emph{pro-$p$-Kummer-{\'e}tale} $\Z_p$-local system, which requires that the associated representation from the Kummer \'etale fundamental group $\pi_1^{\ket}(X, \zeta)$ of $(X, \ca{M}_X)$ at each log geometric point $\zeta$ has no contribution from the prime-to-$p$ part; see Definition \ref{def-p-neat-fin-loc-sys}. In particular, when the local system has unipotent monodromy along the boundaries, it is pro-$p$-Kummer.  
\begin{thmB}[{Theorem \ref{thm-equi-cat-gen}}]
Let $(X,\ca{M}_X)$ be a locally Noetherian fs log adic space over $\Spa(\rQ_p, \Z_p)$, and let $(\GG, \mu)$ be as above. Then there exists an equivalence of categories between:
\begin{enumerate}
\item log $\GG$-shtukas over $(X, \ca{M}_X)^{\Dia} \to \Spd E$ with one leg bounded by $\mu$;
\item pairs $(\PPp, H)$ where $\PPp$ is a pro-$p$-Kummer-\'etale $\ul{\ca{G}(\bb{Z}_p)}$-torsor defined over $X^{\log\Dia}$ and $H: \PPp \to \ca{F}_{G,\mu^{-1}}^\Diamond$ is a $\underline{\GG(\Z_p)}$-equivariant map of $v$-sheaves over $\Spd E$.
\end{enumerate}
\end{thmB}
For a related theorem in the log prismatic theory, see \cite[Thm. 7.36]{KY25} and an update \cite{IKY26} by Inoue-Koshikawa-Yao. In fact, our proof is similar to the one presented there. This theorem also tells us that $(\sh_K^\Sigma)^{\log \Diamond}\to \Sht_{\G^c,\mu^c}$ should not factor through $\sh_K^{\mmin,\Diamond}$, even though the Hodge-Tate period map on $\sh_{K^p}^{\Sigma, \Diamond}$ should factor through $\sh_{K^p}^{\mmin, \Diamond}$.\par

Now, let $(X, \ca{M}_X)$ be an fs log scheme where $X$ is a separated, normal, and flat scheme of finite type over $\Spec \Z_p$. There are different ways to associate a log adic space to $(X, \ca{M}_X)$ and produce diamonds $(X, \ca{M}_X)^{\Dia}$, $(X, \ca{M}_X)^{\dia}$, $(X, \ca{M}_X)^{\Dia/}$. The log generalization of \cite[Thm. 2.7.7]{PR24} is as follows:
\begin{thmC}[{Theorem \ref{thm-ext-shu-gen}}]
The restriction functor
$$\Res^{\mathscr{X}}_{\mathscr{X}_{\eta}}: \Sht_{\ca{G},\mu}((\mathscr{X}, \ca{M}_{\mathscr{X}})^{\Diamond/})\lra\Sht_{\ca{G},\mu}((\mathscr{X}_{\eta}, \ca{M}_{\mathscr{X}_{\eta}})^{\Dia})$$
is fully faithful.
\end{thmC}
In addition, we show that log shtukas on special fibers are non-log shtukas (see Corollary \ref{cor-special-fiber-nonlog}). That is, in characteristic $p$, a log $\GG$-shtuka (with one leg bounded by $\mu$) on $(\mathscr{X}, \ca{M}_{\mathscr{X}})^{\dia}$ is equivalent to an actual $\GG$-shtuka (with one leg bounded by $\mu$) on $\mathscr{X}^{\dia}$.

\subsubsection*{Canonical integral models}

Let $(G, X)$ be a Shimura datum. Fix an open compact subgroup $K\sbst G(\A)$. Choose an admissible (rational polyhedral), smooth, projective cone decomposition $\Sigma$ (without self-intersections).
Denote the toroidal compactification by $\sh_K^{\Sigma}:=\sh_K^{\Sigma}(G,X)$ and the minimal compactification by $\sh_K^\mmin:=\sh_K^\mmin(G,X)$; they are defined over the reflex field $\mathbb{E}$. There is a proper morphism $\oint_{K}^{\Sigma}: \sh_K^\Sigma \to \sh_K^\mmin$ that is compatible with the stratifications on the source and the target (see \cite{AMRT10} and \cite{Pin89}).

Fix a prime number $p$ and a place $v$ of $\mathbb{E}$ over $p$. Let $E = \mathbb{E}_v$ be the completion. Following \cite[Prop. 2.1.2]{lan2018compactifications} and \cite[Thm. 4.1.5]{Mad19}, one expects that there exist normal, proper, flat models $\mathscr{S}_K^\Sigma$ and $\mathscr{S}_K^\mmin$ for $\sh^\Sigma_{K,E}$ and $\sh_{K,E}^\mmin$ over $\ca{O}_{E}$, respectively, such that a list of qualitative properties stated in Axiom \ref{axiom-good-compactification} is satisfied.

Let us explain some of the axioms. Fix a cusp label $[\Phi] = [(Q_{\Phi}, X_{\Phi}^+, g_{\Phi})] \in \Cusp_K(G, X)$. We have a tower of integral models of mixed Shimura varieties:
$$\mathscr{S}_{K_\Phi} \to \mathscr{S}_{\ovl{K}_\Phi} \to \mathscr{S}_{K_{\Phi,h}}\quad (:= \mathscr{S}_{K_\Phi}(P_\Phi,D_\Phi)\to \mathscr{S}_{\ovl{K}_\Phi}(\overline{P}_\Phi,\overline{D}_\Phi)\to \mathscr{S}_{K_{\Phi,h}}(G_{\Phi,h},D_{\Phi,h}))$$
which are normal flat schemes of finite type over $\ca{O}_{E}$ that extend those on the generic fiber. The first map is also a torsor under the same split torus $\mbf{E}_{K_\Phi}$. The second map is proper and surjective. There is a $\Delta_{\Phi, K}$-action on the whole tower that factors through a finite quotient on $\mathscr{S}_{K_{\Phi,h}}$.

Pick $\sigma \in \Sigma^+(\Phi) \subset X_*(\mbf{E}_{K_\Phi}) \otimes \R$. There is a normal subgroup $\Delta_{\Phi, K}^{\circ}\triangleleft \Delta_{\Phi, K}$ that stabilizes $\sigma$. The boundary of $\mathscr{S}_K^\Sigma$ (resp. $\mathscr{S}_K^\mmin$) is stratified by $\ca{Z}_{[(\Phi,\sigma)],K}\iso \Delta^\circ_{\Phi,K}\bss \mathscr{S}_{K_\Phi,\sigma}$ (resp. $\ca{Z}_{[\Phi],K}\iso \Delta_{\Phi,K}\bss \mathscr{S}_{K_{\Phi,h}}$). \'Etale locally on $\mathscr{S}_K^\Sigma(G, X)$ at the boundary stratum $\ca{Z}_{[(\Phi,\sigma)],K}$, the morphism $\mathscr{S}_K \hookrightarrow \mathscr{S}_K^\Sigma$ is a toroidal embedding $\Delta^\circ_{\Phi,K}\bss\Shum{K_{\Phi}} \hookrightarrow \Delta^\circ_{\Phi,K}\bss\Shum{K_{\Phi}}(\sigma)$.

Moreover, there is a stronger strata-preserving isomorphism
$$\mathfrak{X}_{\sigma,K}^\circ\iso \Delta^\circ_{\Phi,K}\bss \cpl{\mathscr{S}_{K_\Phi}(\sigma)}{\mathscr{S}_{K_\Phi,\sigma}^+}, \quad \mathscr{S}^+_{K_\Phi,\sigma}:= \cup_{\tau\subset \sigma,\tau\in \Sigma^+(\Phi)}\mathscr{S}_{K_\Phi,\tau},$$
where
$$\mathfrak{X}^\circ_{\sigma,K}:=\cpl{\mathscr{S}_K^\Sigma}{\ca{Z}^+_{\sigma,K}},\quad \ca{Z}^+_{\sigma,K}:=\cup_{\tau\subset \sigma,\tau\in\Sigma^+(\Phi)}\ca{Z}_{[(\Phi,\tau)],K}.$$
Let $(G, X)$ be an abelian-type Shimura datum, let $p$ be any prime, and let $K=K_pK^p$ be any neat level. By a series of works \cite{FC90}, \cite{Lan16b}, \cite{Mad19}, \cite{Wu25}, and others, there exist good toroidal and minimal compactifications of integral models satisfying Axiom \ref{axiom-good-compactification} (see Theorem \ref{thm-abelian-type-axiom}).

On the other hand, in \cite{PR24}, \cite{daniels2024conjecture}, and \cite{DY25}, when $K_p$ is parahoric, a conjectural framework of \emph{canonical integral models} of the pro-system $\lrbracket{\shu{K_pK^p}(G, X)}_{K^p \subset G(\A^p)}$ was also introduced and constructed essentially in the abelian-type case. 

In order to define and construct \emph{canonical integral models} of toroidal compactifications, we combine the two above-mentioned aspects and study the canonical integral models of \emph{boundary mixed} Shimura varieties. For a mixed Shimura datum $(P, \mathcal{X})$, $P$ is a non-reductive group. We generalize the notions of local models and shtukas with one leg bounded by $\mu$ to non-reductive groups, and prove many functoriality results that are used in later sections. The definitions and propositions are parallel to those defined using reductive groups.

\begin{itemize}

\item We formulate a conjectural framework \ref{def: canonical model for mixed Shimura data} for \emph{canonical integral models} of the pro-system of mixed Shimura varieties arising from the boundary $\lrbracket{\shu{K_{\Phi, p}K_{\Phi}^p}(P_{\Phi}, D_{\Phi})}_{K_{\Phi}^p \subset P_{\Phi}(\A^p)}$, following \cite[Conj. 4.2.2]{PR24}, \cite[Def. 4.1.2]{daniels2024conjecture}, and \cite[Def. 4.3]{DY25}. In particular, we require an extension of the $\PP_{\Phi}^*$-shtuka with one leg bounded by $\mu_{\Phi}^*$ from $\shu{K_{\Phi}}(P_{\Phi}, D_{\Phi})$ to the integral model $\Shum{K_{\Phi}}(P_{\Phi}, D_{\Phi})$. Here we use a slightly different group $P_{\Phi}^*$ than $P_{\Phi}^c$; see Remark \ref{rk: why use star not c}. We show the uniqueness and functoriality of such canonical integral models; see Proposition \ref{prop: morphisms extend to integral models} and Corollary \ref{cor: uniqueness of canonical mixed}.

\item We formulate a conjectural framework for \emph{canonical integral models} of the pro-system of toroidal compactifications $\lrbracket{\Shum{K}^{\Sigma}}_{K^p}$ and minimal compactifications $\lrbracket{\Shum{K}^{\min}}_{K^p}$ (see Definition \ref{def-PR-int-mod}). We show the uniqueness and functoriality of such canonical integral models (see Proposition \ref{prop-full-functoriality}). We also obtain analogous statements for integral models of minimal compactifications.
\end{itemize}

\begin{thmD}[{Theorem \ref{thm-can-int-mod-cpt-summary} and Theorem \ref{thm-ext-cim-ab}}]
Let $(G, X)$ be an abelian-type Shimura datum, $p$ be any prime, and $K_p$ be any quasi-parahoric level.
\begin{enumerate}
    \item There exist canonical integral models $\lrbracket{\Shum{K}(G, X)}_{K^p}$.
    \item There exist canonical integral models $\lrbracket{\Shum{K_{\Phi, p}K_{\Phi}^p}(P_{\Phi}, D_{\Phi})}_{K_{\Phi}^p \subset P_{\Phi}(\A^p)}$ for each $[\Phi]$.
    \item There exist canonical integral models $\lrbracket{\Shum{K}^{\Sigma}}_{K^p}$ for a final collection of smooth projective cone decompositions $\Sigma$, and canonical integral models $\lrbracket{\Shum{K}^{\min}}_{K^p}$.
\end{enumerate}

Such canonical integral models are unique and functorial.
\end{thmD}

\begin{rkS}
We remark that, even for compact Shimura varieties, the first item already extends the known cases in the literature (see \cite{PR24}, \cite{daniels2024conjecture} and \cite{DY25}) to cases allowing all primes and all quasi-parahoric levels for all abelian-type Shimura data. 
In \cite{daniels2024conjecture}, the authors showed the first item of Theorem D in the Hodge-type case for any $p$ and any quasi-parahoric $K_p$; in \cite{DY25}, the authors proved the result in the abelian-type case when $p>2$ and $K_p$ is parahoric. \par
In fact, even if one just wants to generalize these results to compactifications, it is still much more convenient to work with Bruhat-Tits stabilizer levels than parahoric levels. This fact forced us to prove the theorem in the above generality.
\end{rkS}
\begin{rkS}[{See also Remark \ref{rk-local-model-diagram}}]
To show the theorem above, it is crucial for us to handle non-$R$-smooth (cf. \cite[2.1.4]{KPZ24} for the definition of $R$-smoothness) Shimura data $(G,X)$ with parahoric levels $K_p^\circ$ whose intersections with $G^\der(\bb{Q}_p)$ are \emph{not} parahoric. For this, we need to construct Hodge-type liftings that are \textbf{accessible} (see Definition \ref{def-compatible-in-derived-part}) to the abelian-type ones in the quasisplit nonsplit $D^{\bb{H}}$-type case at all primes. This calls for a certain group-theoretic refinement in the construction.\par
However, our construction itself cannot deduce new cases of local model diagrams. When $p>2$, this has been established in \cite{KPZ24}; when $p=2$, for the progress in this direction, cf. Jie Yang's work \cite{Yan25}. One can show that, with inputs in \cite{daniels2024conjecture}, the $\G^{\ad,\circ}$-local model diagrams in \cite{KPZ24} are schematic local model diagrams in the sense of \cite[Def. 4.9.1]{PR24}; for this, see Proposition \ref{prop: schematic local model diagram, 2}.
\end{rkS}
An essential feature of canonical integral models of compactifications is the existence of a log $\G^c$-shtuka $\mathscr{P}^\can$ with one log bounded by $\mu$ on $\mathscr{S}_K^\Sigma$ that extends the given shtuka on $\mathscr{S}_K$; that is, for a canonical integral model $\mathscr{S}_K^\Sigma$, there is a morphism
$$(\mathscr{S}_K^\Sigma)^{\log \Diamond}\to \Sht_{\G^c,\mu^c}.$$
Our method here provides a description of these extensions at boundaries, which will be given below in more detail; and an important advantage is that the extension step only uses (log) $\G^c$-shtukas, which makes it possible to be applied to cases beyond abelian-type Shimura varieties (see Theorem \ref{thm-ext-toric-equiv} and Lemma \ref{lem-gluing-shtukas}).
\subsubsection*{Shtukas at the boundary}

From this subsection onward, fix any Shimura datum $(G, X)$, and let $K_p$ be a quasi-parahoric level subgroup. We work with canonical integral models $\lrbracket{\Shum{K}^{\Sigma}}_{K^p}$ and $\lrbracket{\Shum{K}^{\min}}_{K^p}$, and we prove many results without the use of abelian schemes and $p$-divisible groups\footnote{In fact, our method works in a more general context; see the discussion in \S\ref{subsec-can-mod-II}; indeed, we only need a weaker assumption \ref{ass-well-position}.}. This framework allows us to work with more general types of Shimura varieties once certain axioms can be verified. 

 Let $\mathfrak{W} = \Spf( R,I) \subset \mathfrak{X}_{\sigma}^{\circ}$ be an affine open formal subscheme. We can consider the flat morphisms 
 $W = \Spec R \to \Shumc{K}{\Sigma}$ and $W \to \Delta_{\Phi,K}^\circ\backslash\Shum{K_\Phi}(\sigma)$.
Let $W^0 \subset W$ be the common open subscheme associated with $\Shum{K} \subset \Shumc{K}{\Sigma}$ and $\Delta_{\Phi,K}^\circ\bss\Shum{K_\Phi} \subset \Delta_{\Phi,K}^\circ\bss\Shum{K_\Phi}(\sigma)$.

\begin{propS}[{Proposition \ref{prop: shtuka comparison, over W, integral base}}]
We have the following important diagram:
\begin{equation}\label{intro-prop-shtuka comparison}
\begin{tikzcd}
	{\Shum{K}(G, X)^{\dia}} & {W^{0, \dia}} & {\Delta_{\Phi, K}^{\circ}\backslash\Shum{K_{\Phi}}(P_{\Phi}, D_{\Phi})^{\dia}} & {\Delta_{\Phi, K}^{\circ}\backslash\Shum{K_{\Phi, h}}(G_{\Phi, h}, D_{\Phi, h})^{\dia}} \\
	{\Sht_{\GG^c, \mu^c, \delta = 1}} & {} & {\Sht_{\PP_{\Phi}^*, \mu_{\Phi}^*, \delta = 1}} & {\Sht_{\GG_{\Phi, h}^*, \mu_{\Phi, h}^*, \delta = 1}.}
	\arrow[from=1-1, to=2-1]
	\arrow[from=1-2, to=1-1]
	\arrow[from=1-2, to=1-3]
	\arrow[from=1-3, to=1-4]
	\arrow[from=1-3, to=2-3]
	\arrow[from=1-4, to=2-4]
	\arrow["{\Int(g_{\Phi}^{-1})}"', from=2-3, to=2-1]
	\arrow[from=2-3, to=2-4]
\end{tikzcd}
\end{equation}
\end{propS}

This tells us how shtukas degenerate along the boundary. One should compare this result with the degeneration of abelian schemes with extra structures introduced in \cite{FC90} and \cite{Lan13}.

On the other hand, we can extend the shtukas over the boundary. Constructing log shtukas on $(\Shum{K_{\Phi}}(\sigma), \ca{M}_{\Shum{K_{\Phi}}(\sigma)})$ is more approachable than constructing log shtukas on $(\mathscr{S}^\Sigma_{K}, \ca{M}_{\mathscr{S}^\Sigma_{K}})$, since we can work with an actual torsor under a split torus. By gluing log shtukas on $(\Shum{K_{\Phi}}(\sigma), \ca{M}_{\Shum{K_{\Phi}}(\sigma)})$ for $[(\Phi, \sigma)] \in \Cusp_K(G, X, \Sigma)$, we get a log shtuka on $(\mathscr{S}^\Sigma_{K}, \ca{M}_{\mathscr{S}^\Sigma_{K}})$:

\begin{thmE}[{Corollary \ref{cor-ext-toric-emb}, Theorem \ref{prop-degeneration-int}}]\leavevmode
  \begin{enumerate}
      \item For each $[(\Phi, \sigma)] \in \Cusp_K(G, X, \Sigma)$, there is a unique log $\PP_{\Phi}^*$-shtuka on $\Delta_{\Phi, K}^{\circ}\backslash \Shum{K_{\Phi}}(\sigma)$ with one leg bounded by $\mu_{\Phi}^*$ extending the one on the generic fiber, i.e., there exists a morphism $\Shum{K_{\Phi}}(\sigma)^{\log\Dia/} \to \Sht_{\PP_{\Phi}^*, \mu_{\Phi}^*, \delta = 1}$.
      \item There exists a unique log $\GG^c$-shtuka on $\mathscr{S}^\Sigma_{K}$ with one leg bounded by $\mu^c$ extending the one on the generic fiber, i.e., there exists a morphism $(\mathscr{S}^\Sigma_{K})^{\log\Dia/} \to \Sht_{\GG^c, \mu^c, \delta = 1}$.
  \end{enumerate}
\end{thmE}
Since $\mathscr{S}^\Sigma_{K}(G, X)$ is proper, we have $(\mathscr{S}^\Sigma_{K}, \ca{M}_{\mathscr{S}^\Sigma_{K}})^{\Dia/} = (\mathscr{S}^\Sigma_{K}, \ca{M}_{\mathscr{S}^\Sigma_{K}})^{\Dia} = (\mathscr{S}^\Sigma_{K}, \ca{M}_{\mathscr{S}^\Sigma_{K}})^{\dia}$. We restrict it to the special fiber and obtain a morphism
$(\mathscr{S}^\Sigma_{K, \Bar{s}}, \ca{M}_{\mathscr{S}^\Sigma_{K, \Bar{s}}})^{\dia} \to \Sht_{\GG^c, \mu^c}$, where $\Bar{s} = \Spec \ovl{\F}_p$. Let $\Sht^W_{\GG^c, \mu^c}$ be the Witt vector $\GG^c$-shtuka with one leg bounded by $\mu^c$ introduced in \cite{xiao2017cycles} (cf. \cite{shen2021ekor}; see \cite[Rmk. 3.1.8]{DvHKZ24ig} for the sign convention). Combining with Corollary \ref{cor-special-fiber-nonlog}, we obtain:
\begin{corA}[{Corollary \ref{cor-degeneration-sp}}]\label{intro-thm-toroidal}
We have a $\GG^c$-shtuka with one leg bounded by $\mu^c$ on the special fiber of $\mathscr{S}^\Sigma_{K}(G, X)$, i.e., a morphism $\mathscr{S}^\Sigma_{K}(G, X)_{\Bar{s}}^{\perf} \to \Sht^W_{\GG^c, \mu^{c}, \delta = 1}$.
\end{corA}
\begin{corB}[{Corollary \ref{cor-big-diamond}}]
We have a $\GG^c$-shtuka with one leg bounded by $\mu^c$ on the big diamond $\mathscr{S}_{K}(G, X)^{\Dia}$, i.e., a morphism $\mathscr{S}_{K}(G, X)^{\Dia} \to \Sht_{\GG^c, \mu^{c}, \delta = 1}$.
\end{corB}
Neither Corollary A nor Corollary B appears previously in the literature, even in cases such as modular curves. For cases with good moduli interpretations, we expect that these results could be proved using log $p$-divisible groups.
\subsubsection*{Well-positioned subschemes}
Using the Witt vector shtuka $\Shum{K}(G, X)_{s}^{\perf} \to \Sht^W_{\GG^c, \mu^c, \delta = 1}$, and following \cite{PR24}, we can construct Newton strata, central leaves, KR strata, and EKOR strata on $\Shum{K}(G, X)_{\Bar{s}}$.

In \cite{boxer2015torsion} and \cite{lan2018compactifications}, the authors formulate the notion of \emph{well-positioned subschemes}. Well-positioned strata have good (partial) toroidal and minimal compactifications along the boundary, with many good qualitative properties described in Axiom \ref{axiom-good-compactification}, as if they were Shimura varieties in characteristic $0$. It was proved in \cite{lan2018compactifications} that, for PEL-type Shimura varieties, Newton strata, central leaves, KR strata, and EKOR strata are well positioned. These results were generalized to the Hodge-type case in \cite{Mao25}, where it was shown that their partial minimal compactifications are again of the same type of strata (e.g., the boundary of the partial minimal compactification of a Newton stratum is stratified by Newton strata on the boundary). The proofs in \cite{lan2018compactifications} essentially used the degeneration of $p$-divisible groups along the boundary. Under the degeneration diagram (\ref{intro-prop-shtuka comparison}) (cf. Proposition \ref{prop: shtuka comparison, over W, integral base}), and since the shtukas in characteristic $p$ are crystalline in nature, we find that the well-positioned property of all these strata can be naturally explained using the degeneration of shtukas.

Assume there exist canonical integral models $\lrbracket{\Shum{K}^{\Sigma}}_{K^p}$ and $\lrbracket{\Shum{K}^{\min}}_{K^p}$ in the sense of Definition \ref{def-PR-int-mod} (for example, in the setting of Theorem D).

\begin{thmF}\leavevmode
\begin{enumerate}
\item Newton strata (resp. central leaves, KR strata, and EKOR strata), their connected components, and their closures are well positioned locally closed subschemes in the sense of Definition \ref{def: well-positioned}.
\item The boundary of the partial minimal compactifications of Newton strata (resp. central leaves, KR strata) is stratified by Newton strata (resp. central leaves, KR strata) on the boundary; see Proposition \ref{prop: Newton strata are well-positioned}, \ref{prop: central leaves are well-positioned}, and \ref{prop: KR strata are well-positioned} for details.
\end{enumerate}
\end{thmF}

Besides generalizing the well-positioned property and the boundary descriptions of partial minimal compactifications from the Hodge-type case to a more general framework (at least to the abelian-type case), we also provide a new description of partial toroidal compactifications of these strata.

\begin{propS}[{Proposition \ref{prop: toroidal of Newton strata}, \ref{prop: toroidal of central leaves}, \ref{prop: toroidal of KR strata}, \ref{prop: toroidal of EKOR strata}}]
The partial toroidal compactifications of Newton strata (resp. central leaves, KR strata, and EKOR strata) coincide with the Newton strata (resp. central leaves, KR strata, and EKOR strata) defined using $\Shum{K}^{\Sigma}(G, X)_{\Bar{s}}^{\perf} \to \Sht^W_{\GG^c, \mu^c, \delta = 1}$.
\end{propS}

\begin{rkS}
In a recent paper \cite{Inoue2025LogPD}, Kentaro Inoue constructed a log prismatic realization in the Hodge-type case for smooth integral models of toroidal compactifications by developing a log prismatic Dieudonn{\'e} theory using log $p$-divisible groups. 
Combining with it, the method provided here can be extended to the log prismatic context for more general $(G,X)$, and may also give the log prismatic realizations finer boundary descriptions. In fact, we can construct a realization functor from log prismatic $F$-crystals to log shtukas that is compatible with the \'etale realization functor. Together with Inoue, we plan to explore this direction.
\end{rkS}
\subsubsection*{Structure of the article}
Section \ref{sec-shtukas-nonreductive} can be viewed as a preliminary section of this paper, where we start by generalizing the definition of some objects in \cite{SW20} to non-reductive groups and proving basic properties. 
In Section \ref{sec-log-diamond-shtuka}, we define and study log diamonds and shtukas, generalizing and modifying the previous work \cite{PR24} and \cite{KY25}. Some results analogous to \cite[Sec. 2]{PR24} will be shown in the log setting.
In Section \ref{sec-shtukas-mixed-sh}, we study the $p$-adic local systems on mixed Shimura varieties.
In Section \ref{sec-canonical}, we define canonical integral models for boundary mixed Shimura varieties. Some basic properties of these canonical integral models will be proved.
In Section \ref{sec-can-ext}, we show the existence of canonical extensions of shtukas (or rather, log shtukas) in an axiomatic setup. We also define canonical integral models for toroidal/minimal compactifications.
In Section \ref{sec-can-ext-ab}, we demonstrate that canonical integral models for compactifications exist for all abelian-type Shimura data, all primes, and all quasi-parahoric levels. In particular, the canonical integral models for (boundary mixed) Shimura varieties also exist in the same generality.
In Section \ref{sec-well-position}, we prove the well-positionedness of all well-known stratifications of special fibers of Shimura varieties under the axioms in Section \ref{sec-can-ext}.\par
Here is a workflow diagram:\par
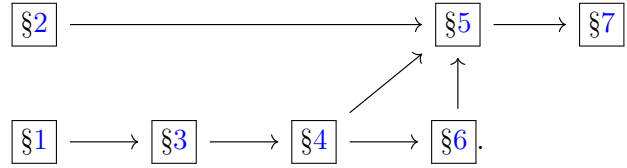
\begin{equation*}
    \begin{tikzcd}
    \fbox{\S \ref{sec-log-diamond-shtuka}}\arrow[rrr]&&&\fbox{\S \ref{sec-can-ext}}\arrow[r]&\fbox{\S \ref{sec-well-position}}\\
    \fbox{\S \ref{sec-shtukas-nonreductive}}\arrow[r]&\fbox{\S \ref{sec-shtukas-mixed-sh}}\arrow[r]&\fbox{\S \ref{sec-canonical}}\arrow[ur]\arrow[r]&\fbox{\S \ref{sec-can-ext-ab}}.\arrow[u]&
    \end{tikzcd}
\end{equation*}
\subsubsection*{Acknowledgments}
The authors thank Heng Du, Kentaro Inoue, and Yupeng Wang for helpful discussions on log geometry and $p$-adic Hodge theory during the preparation of this work. The authors thank Jingren Chi, Pol van Hoften, Dongryul Kim, Wansu Kim, Kai-Wen Lan, Sian Nie, Shen Xu, Jie Yang and Alexander Youcis for the inspiring discussions.\par
We would like to thank MCM and BICMR for a great academic environment.
\subsubsection*{Notation and conventions}
Fix a prime number $p>0$. All monoids are commutative. All rings have identities. For the conventions in $p$-adic geometry (and perfectoid geometry), we follow \cite{SW20} and \cite{scholze2017etale}. For the definition of an adic space being \emph{locally Noetherian}, we follow \cite[p.4]{DLLZ23}.  For conventions on (mixed) Shimura varieties and compactifications, we mainly follow \cite{Wu25} (we refer readers to the latest version \href{https://peihang-wu.github.io/files/cpt-ab.pdf}{here}, which will be updated on arXiv in the future). For the conventions in log geometry, we mainly follow \cite{Kat89}, \cite{Ogu18}, and \cite{DLLZ23}. A constant sheaf of monoids on $X$ with values in $\md{P}$ is denoted by $\md{P}_X$ instead of $\ul{\md{P}}$ to avoid awkward conventions such as ``$\overline{\ul{\md{P}}^a}$''. Log structures are written multiplicatively, but charts are written additively. 
For a Shimura datum $(G, X)$, $G$ is a reductive group over $\rQ$. We also write $G \otimes_{\rQ} \rQ_p$ as $G$ when it is clear in the context.
There are several ``$E$'' in different fonts: In general, ``$\bb{E}$'' denotes global fields; ``$E$'' denotes local fields; ``$\mbf{E}$'' denotes tori.
For a smooth affine group scheme $\GG$ over $\Z_p$, we usually denote $K_p = \GG(\Z_p)$ and $\KK_p = \GG(\bZ_p)$. Here $\breve{(\ast)}$ means the set of $\bZ_p$-points.
Let $\Perf$ denote the category of perfectoid spaces of characteristic $p$, and $\Perfd$ denote the category of all perfectoid spaces. Let $\PCAlg$ denote the category of affine perfect schemes. 
Some Galois groups: $\Gamma = \Gal(\ovl{\rQ}_p|\rQ_p)$, $I = \Gal(\ovl{\rQ}_p|\bQ)$, $\Sigma_0 = \lrangle{\sigma} = \Gal(\bQ|\rQ_p)$.\par
The group-theoretic conventions in \S\ref{sec-shtukas-nonreductive} and a few lemmas and definitions in \S\ref{sec-shtukas-mixed-sh} and \S\ref{sec-canonical} are slightly different from other places in the article. In \S\ref{sec-shtukas-nonreductive}, $G$ usually denotes the Levi quotient of a linear algebraic group $P$, while an embedding of $P$ into a reductive group is denoted by $P\hookrightarrow G'$. In other places, $G$ usually denotes the group that $P$ (or rather, $P_\Phi$) maps into.
\section{Shtukas for non-reductive groups}\label{sec-shtukas-nonreductive}
To understand the geometry of toroidal compactifications of Shimura varieties, we study mixed Shimura varieties. As a first step, we discuss local models and shtukas for non-reductive groups in detail. Many properties are analogous to those for reductive groups; consequently, we present only the properties needed in later sections.

\subsection{The $B_{\dR}^+$-affine Grassmannian}

\subsubsection{Setup}

Let $\PP$ be a smooth affine group scheme over $\Z_p$.

\begin{Definition and Proposition}[{\cite[\S 19, 20]{SW20}}]
     Recall that the $B_{\dR}^+$-affine Grassmannian $\Gra{\PP} \to \Spd \Z_p$ is a $v$-sheaf on $\Perf$ admitting the following equivalent descriptions.
     \begin{enumerate}
         \item It is the \'etale sheafification of the presheaf quotient $L\PP/L^+\PP$. Recall that the loop group functor $L\PP$ and the positive loop group subfunctor $L^+\PP$ are defined on $\Spd \Z_p$: Given $S = \Spa(R, R^+) \in \Perf$ equipped with an untilt $S^{\sharp} = \Spa(R^{\sharp}, R^{\sharp+})$ over $\Spa(\Z_p)$, we have
         \[ L\PP: (R^{\sharp}, R^{\sharp+}) \mapsto \PP(B_{\dR}(R^{\sharp})),\quad L^+\PP: (R^{\sharp}, R^{\sharp+}) \mapsto \PP(B^+_{\dR}(R^{\sharp})). \]
         \item It is a functor whose $S = \Spa(R, R^+)$-points parametrize untilts $S^{\sharp} = \Spa(R^{\sharp}, R^{\sharp+})$ over $\Spa(\Z_p)$ together with a $\PP$-torsor $\EE$ on $\Spec B_{\dR}^+(R^{\sharp})$ and a trivialization of it over $\Spec B_{\dR}(R^{\sharp})$.
         \item It is a functor whose $S = \Spa(R, R^+)$-points parametrize untilts $S^{\sharp} = \Spa(R^{\sharp}, R^{\sharp+})$ over $\Spa(\Z_p)$ together with a $\PP$-torsor $\EE$ on $S \dot{\times} \Z_p$ and a trivialization of $\EE|_{S \dot{\times} \Z_p\setminus S^{\sharp}}$ that is meromorphic along $S^{\sharp}$.
     \end{enumerate}
\end{Definition and Proposition}
\begin{proof}
    The equivalence of the second and third descriptions is proved using Beauville--Laszlo gluing, which is compatible with the Tannakian formalism; see, for example, \cite[Lem. 7.2]{anschutz2022extending}. The equivalence of the first and second descriptions is given by \cite[Prop. 19.1.2, 19.5.3]{SW20}. The proofs do not need the group to be reductive.
\end{proof}

We recall some simple facts.
\begin{lem}[{\cite[Thm. 1.4]{pappas2008twisted}}]\label{lem: quotient affine, closed embedding}
Let $\PP_1 \to \PP_2$ be a morphism of smooth affine group schemes. Assume that $\PP_1$ is the smoothing of the closure of $P_1$ in $\PP_2$. If $\PP_2/\PP_1$ is affine, then $\Gra{\PP_1} \to \Gra{\PP_2}$ is a closed embedding. 
\end{lem}
\begin{proof}
    This is essentially \cite[Thm. 1.4]{pappas2008twisted} (cf.\ \cite[Lem. 19.1.5]{SW20}, where we only need the input that $\PP_2/\PP_1$ is affine). Note that the proof in \cite[Lem. 19.1.5]{SW20} does not require that $\PP_1 \to \PP_2$ is a closed embedding; it only requires that, for any $\Spa(R, R^+) \in \Perf$ with an untilt $\Spa(R^{\sharp}, R^{\sharp+})$, we have $\PP_1(B_{\dR}^+(R^{\sharp})) = P_1(B_{\dR}(R^{\sharp})) \cap \PP_2(B_{\dR}^+(R^{\sharp}))$. This is always true under our assumption; see \cite[Lem. 4.6.1]{PR24}.
\end{proof}

\begin{lem}\label{lem: closed embedding of loop groups}
Let $\mathcal{X}$ be an affine scheme of finite type over $\Z_p$, and let $\mathcal{X} \hookrightarrow \mathbb{A}^n$ be a closed embedding. Then the natural inclusion functor $L^+\mathcal{X} \to L^+\mathbb{A}^n$, as well as $L\mathcal{X} \to L\mathbb{A}^n$, is a closed embedding. Moreover, $L^+\mathcal{X} \to L\mathcal{X}$ is a closed embedding.
\end{lem}
\begin{proof}
We show that $L\mathcal{X} \to L\mathbb{A}^n$ is a closed embedding. In the equal-characteristic case over geometric points, see the descriptions in \cite[\S 1]{pappas2008twisted}. In the mixed-characteristic case, see \cite[Prop. 1.1]{zhu2017affine}. In general, given an $S = \Spa(R, R^+)$-point $x \in L\mathbb{A}^n(R, R^+) = B_{\dR}(R^{\sharp+})^n$, we need to prove that the locus where this point lies in the image of $L\mathcal{X}$ is closed. The result follows from the proof of \cite[Lem. 19.1.4]{SW20}: by induction, assume $x \in B_{\dR}^+(R^{\sharp+})^n/\xi B_{\dR}^+(R^{\sharp+})^n$; then $x$ lies in the image of $L\mathcal{X}$ if and only if its reduction modulo $\xi$ lies in the image of $\mathcal{X}$. Similarly, $L^+\mathcal{X} \to L^+\mathbb{A}^n$ is a closed embedding, and therefore $L^+\mathcal{X} \to L\mathcal{X}$ is a closed embedding.
\end{proof}

Also, let us introduce a useful lemma. Given a positive integer $N$, recall that we have a notion of congruence quotient group $L^{\leq N}\PP = L^+\PP/L^+_N\PP$, where $L_N^+\PP := \ker(L^+\PP \to \PP(B_{\dR}^+/\xi^N))$.
\begin{lem}\label{lem: congruence subgroup}
    Let $X \subset \Gra{\PP}$ be a quasi-compact and quasi-separated sub-sheaf, then the natural action $a_X: L^+\PP \times X \to \Gra{\PP}$ factors through $a_{X, N}: L^{\leq N}\PP \times X \to \Gra{\PP}$ for some large $N$. Moreover, $a_{X, N}$ is proper.
\end{lem}
\begin{proof}
    Let us show the first statement. Fix a faithful representation $\rho: \PP \to \GL_n$ such that $\GL_n/\PP$ is quasi-affine (such an embedding exists; see \cite[Prop. 1.3]{pappas2008twisted}; see also \cite[Cor. 3.8]{richarz2020basics} for a more general statement), and then apply the arguments in \cite[Lem. 19.1.5]{SW20} (cf. \cite[Prop. 1.20]{zhu2017affine}) to show that $\Gra{\PP} \to \Gra{\GL_n}$ is a locally closed embedding. The image of $X$ in $\Gra{\GL_n}$ is quasi-compact and quasi-separated, and $L_N^+\PP = L^+\PP \cap L_N^+\GL_n$. It suffices to work with $\PP = \GL_n$, where the result is well-known. Now, $a_{X, N}$ is quasi-compact since both $L^{\leq N}\PP$ and $X$ are quasi-compact, is partially proper since $L^{\leq N}\PP$, $X$ and $\Gra{\PP}$ are partially proper by construction. In particular, $a_{X, N}$ is proper by \cite[Cor. 17.4.8]{SW20}.
\end{proof}

\subsubsection{Diamond functors}
For a scheme $X$ over $\bb{Z}_p$, there are three different notions of diamonds associated with $X$. See \cite[\S2.2]{anschutz2022p} and \cite[\S2.1]{PR24}.
\begin{itemize}
    \item ``big diamond'' $X^\Diamond$: A $v$-sheaf $X^\Diamond$ sending any affinoid perfectoid $S\in \ob\Perf$ to the isomorphism classes of tuples $(S^\sharp,\iota,f)$, where $(S^\sharp,\iota)$ is an untilt and $(f: S^\sharp\to X)\in X(S^\sharp)$. 
    \item ``small diamond'' $X^\diamond$: A $v$-sheaf $X^\diamond$ sending any affinoid perfectoid $S=\Spa(A,A^+)\in \ob\Perf$ to the isomorphism classes of tuples $(S^\sharp=\Spa(A^\sharp,A^{\sharp,+}),\iota,f)$, where $(S^\sharp,\iota)$ is an untilt and $f\in X(A^{\sharp,+})$.
    \item ``slashed diamond'' $X^{\Diamond/}:= X^{\diamond}\disju_{X^{\diamond}\times_{\Spd \bb{Z}_p}\Spd \bb{Q}_p}(X_{\bb{Q}_p})^\Diamond$.
\end{itemize}

If $X$ is proper over $\bb{Z}_p$, the three notions above are the same by the valuative criterion of properness.\par
For a pre-adic space $X$ over $\bb{Z}_p$, the similarly-defined functor $X^\Diamond$ sending $S\in\ob \Perf$ to $(S^\sharp, \iota,f:S^\sharp\to X)$ is also a $v$-sheaf over $\Perf$ by \cite[Lem. 18.1.1]{SW20}.

\subsubsection{Reduction}

Recall the reduction functor introduced in \cite{gleason2025specialization}. Let $\wdt{\SchPerf}$ denote the topos on $\PCAlg$, endowed with the $v$-topology. Let $\wdt{\Perf}$ denote the topos on $\Perf$ equipped with the $v$-topology. Then $(\dia, \red)$ form a pair of adjoint functors between the topoi $(\wdt{\SchPerf}, \wdt{\Perf})$. Given $\Spec A \in \SchPerf$ and a $v$-sheaf $\FF$, we denote $\FF^{\red}(A) = \Hom(\Spd A, \FF)$. 
The reduction functor $\red$ commutes with finite limits.

Let $\FF$ be a $v$-sheaf over $\Spd \Z_p$, let $\FF_{\F_p}:= \FF_{\Spd \F_p} \to \Spd \F_p$ be the base change of $\FF \to \Spd \Z_p$. We say that $\FF_{\F_p}$ is the special fiber of $\FF$. Note that $(\Spd \Z_p)^{\red}$ is represented by $\Spec \F_p$. We say that $\FF$ is formally $p$-adic (in the sense of \cite[Def 3.20]{gleason2025specialization}), if $\FF_{\Spd \F_p} = (\FF^{\red})^{\dia}$.

Recall the Witt affine Grassmannian $\Gra{\PP}^W$ defined in \cite{zhu2017affine}.
\begin{definition}
The Witt affine Grassmannian $\Gra{\PP}^W$ is the functor defined over $\PCAlg$, sending $\Spec R \in \PCAlg$ to the set of $\PP$-torsors on $\Spec W(R)$ together with a trivialization over $\Spec W(R)[1/p]$.
\end{definition}

It follows from definition that $\Gra{\PP, \F_p} = (\Gra{\PP}^W)^{\dia}$, $\Gra{\PP}^W$ is the reduction of $\Gra{\PP}$. $\Gra{\PP}^W$ is the \'etale sheafification of $L_{\F_p}\PP/L_{\F_p}^+\PP$, where $L_{\F_p}\PP$ and $L_{\F_p}^+\PP$ are the Witt loop group functor and the Witt positive loop group functor respectively. By \cite[Cor. 9.6]{bhatt2017projectivity}, $\Gra{\PP}^W$ is represented by an increasing union of perfections of quasi-projective schemes over $\F_p$ along closed immersions. Here we only need $\GL_n/\PP$ to be quasi-affine for some faithful embedding $\PP \to \GL_n$ over $\Z_p$.

\subsubsection{Over $C$}

\emph{From now on}, let $P$ be a linear algebraic group over $\rQ_p$ with a Levi decomposition $P = U \rtimes G$, where $U = R_u P$ and $G$ is the reductive quotient. Let $\GG$ be a smooth affine model of $G$ over $\Z_p$, and let $\UU$ be a smooth affine model of $U$ over $\Z_p$ such that the conjugation action of $\GG(\bZ_p)$ stabilizes $\UU(\bZ_p)$. By \cite[1.7.3 and Prop. 1.7.6]{BT84}, we can associate a smooth affine model $\PP := \UU \rtimes \GG$ to $P$ over $\Z_p$. Throughout this article, we work with such $\PP$.

Let $C$ be a perfectoid field of characteristic $0$. We write the base changes of $\Gra{\PP}$ and $\Gra{\GG}$ as $\Gra{P, C}$ and $\Gra{G, C}$, respectively. Let $\mu: \Gm \to P$ be a cocharacter over $\ovl{\rQ}_p$. We define $\Gra{P, \mu} \subset \Gra{P, C}$ as the $L^+P$-orbit of $\mu: L\Gm \to LP$, i.e., the $v$-sheaf-theoretic image of $L^+P \times \mu \to \Gra{P, C}$, and define $\Gra{P, \leq \mu}$ as the $v$-sheaf-theoretic closure of $\Gra{P, \mu}$ in $\Gra{P, C}$.

\begin{rk}\label{rk: orbital definition and usual definition}
    Recall that $\Gra{G, \mu}$ (resp. $\Gra{G, \leq \mu}$) has two definitions: one uses such orbital image (resp. and its closure) as in \cite[\S 4]{anschutz2022p}, and the other uses the pointwise boundedness condition as in \cite[\S 19]{SW20}. These two definitions coincide; see \cite[Cor. 4.6]{anschutz2022p}.
\end{rk}

\begin{lem}\label{lem: P and G associated, Gr}
    Fix a section $G \to P$ such that $\mu$ lifts to $G$. We then have a closed embedding $\Gra{G, C} \to \Gra{P, C}$, and $\Gra{P, (\leq)\mu} \subset \Gra{P, C}$ is the $L^+P$-orbit of $\Gra{G, (\leq)\mu} \subset \Gra{G, C}$, i.e., the $v$-sheaf-theoretic image of $L^+P \times \Gra{G, (\leq)\mu} \to \Gra{P, C}$. As a result, $\Gra{P, \mu}$ (resp. $\Gra{P, \leq \mu}$) agrees with the subfunctor of all maps $S \to \Gra{P, C}$ such that, for every geometric point $x: \Spa(D, D^+) \to S$, we have $x \in P(B_{\dR}^+(D^{\sharp}))\xi^{\mu}P(B_{\dR}^+(D^{\sharp}))/P(B_{\dR}^+(D^{\sharp}))$ (resp. $x \in P(B_{\dR}^+(D^{\sharp}))\xi^{\mu'}P(B_{\dR}^+(D^{\sharp}))/P(B_{\dR}^+(D^{\sharp}))$ for some $\mu' \leq \mu$). Here we identify $B_{\dR}(D^{\sharp}) \cong D^{\sharp}((\xi))$.
\end{lem}
\begin{proof}
 This follows mainly from Remark \ref{rk: orbital definition and usual definition}. Only need to show that the image of $L^+P \times \ovl{\Gra{G, \mu}}$ in $\Gra{P, C}$ is $\ovl{\Gra{P, \mu}}$, here $\ovl{(\ast)}$ means the closure under the $v$-topology. Since $\Gra{G, \leq \mu}$ is proper, thus quasi-compact and quasi-separated, by \cite[Prop. 20.2.3]{SW20}, the action $a: L^+P \times \Gra{G, \leq \mu} \to \Gra{P, C}$ factors through $a_N: L^{\leq N}P \times \Gra{G, \leq \mu} \to \Gra{P, C}$ for some large $N$ and $a_N$ is proper (see Lemma \ref{lem: congruence subgroup}). We then apply Lemma \ref{lem: cl is surj}.
\end{proof}
\begin{lem}\label{lem: cl is surj}
    Let $f: X_1 \to X_2$ be a universally closed morphism of small $v$-stacks, and let $Y_1$ and $Y_2$ be sub-$v$-stacks of $X_1$ and $X_2$, respectively. Assume that $f|_{Y_1}$ factors through $Y_2$ and $f: Y_1 \to Y_2$ is surjective. Then $f: Y_1^{\cl} \to Y_2^{\cl}$ is surjective, where $Y^{\cl}_i$ is the v-closure of $Y_i$ in $X_i$. 
\end{lem}
\begin{proof}
    By \cite[Prop. 2.8]{anschutz2022p}, $Y_i^{\cl} = |Y_i|^{\textit{wgc}} \times_{\underline{|Y_i|}} X_i$, where $|Y_i|^{\textit{wgc}}$ is the weakly generalizing closure of $|Y_i| \subset |X_i|$. Since $Y_1 \to Y_2$ is surjective, $|Y_1| \to |Y_2|$ is surjective (see \cite[Prop. 12.9]{scholze2017etale}). By \cite[Lem. 2.4]{anschutz2022p}, since $f$ is universally closed, $|Y_2|^{\textit{wgc}} \subset |f|(|Y_1|^{\textit{wgc}})$, and $f: Y_1^{\cl} \to Y_2^{\cl}$ is surjective.
\end{proof}

\begin{prop}[{\cite[Prop. 19.2.3]{SW20}}]\label{prop: leq mu is a closed embedding}
Fix a maximal torus and a Borel $T \subset B \subset G$. Fix a dominant cocharacter $\mu \in X_*(T)^+$. Then $\Gra{P, \leq \mu} \subset \Gra{P, C}$ is a closed subfunctor, and $\Gra{P, \mu} \subset \Gra{P, \leq \mu}$ is an open subfunctor. Moreover, $\Gra{P, \leq \mu}$ is a locally spatial diamond.
\end{prop}
\begin{proof}
  Let $\Gra{P, \leq \mu}' \subset \Gra{P, C}$ be the pullback of $\Gra{G, \leq \mu} \subset \Gra{G, C}$ under the projection $\Gra{P, C} \to \Gra{G, C}$. Then we have $f: \Gra{P, \leq \mu} \xrightarrow{\subset} \Gra{P, \leq \mu}'$. Since $\Gra{G, \leq \mu} \subset \Gra{G, C}$ is closed by \cite[Prop. 19.2.3]{SW20}, it suffices to show that $f$ is closed.

  We take a strictly totally disconnected perfectoid cover $X$ of the small $v$-sheaf $\Gra{G, \leq \mu}$ and claim that $f_X$ is closed. Note that closedness of a morphism descends along a $v$-cover (\cite[Prop. 10.11]{scholze2017etale}). Take a $v$-cover $\bigsqcup_{\mu' \leq \mu} X_{\mu'} \to X$, where $X_{\mu'}$ is a strictly totally disconnected perfectoid space that is a $v$-cover of $X \times_{\Gra{G, \leq \mu}} \Gra{G, \mu'}$. Fix a point $x\in X_{\mu'}$, and let $U = \Spd(R, R^+) \to X_{\mu'}$ be an \'etale neighborhood of $x$ such that $\Gra{G, C}(U) = G(B_{\dR}(R^{\sharp}))/G(B_{\dR}^+(R^{\sharp}))$. Since every \'etale covering of $X_{\mu'}$ splits, $U$ is again a strictly totally disconnected space. 
  Since the $v$-sheaves involved are partially proper, for simplicity we work with rank-$1$ points and assume $R^+ = R^{\circ}$. Set $\wdt{\OO} = B_{\dR}^+(R^{\sharp})$ and $\wdt{F} = B_{\dR}(R^{\sharp})$. By Lemma \ref{lem: fiberness result} below, the fiber $f_U$ is isomorphic to
  \[
      U(\wdt{\OO})/ U(\wdt{\OO}) \cap \xi^{\mu'}U(\wdt{\OO})\xi^{-\mu'} \hookrightarrow U(\wdt{F})/\xi^{\mu'}U(\wdt{\OO})\xi^{-\mu'},
  \]
  It is a closed embedding since $U(\wdt{\OO}) \hookrightarrow U(\wdt{F})$ is a closed embedding by Lemma \ref{lem: closed embedding of loop groups}.
  
  To prove the last statement, we take a faithful representation $\rho: P \to \GL_n$ such that $\GL_n/P$ is quasi-affine. Then $\Gra{P} \to \Gra{\GL_n}$ is a locally closed embedding. Since $\Gra{P, \leq \mu} \subset \Gra{P}$ and $\Gra{\GL_n, \leq\rho(\mu)} \subset \Gra{\GL_n}$ are closed, it follows that $\Gra{P, \leq \mu} \to \Gra{\GL_n, \leq\rho(\mu)}$ is locally closed. Since $\Gra{\GL_n, \leq\rho(\mu)}$ is a spatial diamond by \cite[Thm. 19.2.4]{SW20}, we conclude that $\Gra{P, \leq \mu}$ is a locally spatial diamond (see \cite[Cor. 11.28]{scholze2017etale}).
\end{proof}

\begin{lem}\label{lem: fiberness result}
    Let $S = \Spa(R, R^{\circ})$ be a strictly totally disconnected perfectoid space. Let $x$ be an $S$-point of $\Gra{G, \mu} \to \Spd(C, \OO_C)$, associated with the untilt $S^{\sharp} = \Spa(R^{\sharp}, R^{\sharp\circ})$. Up to replacing $S$ with an open subspace, we have
    \[
        x = g\xi^{\mu}G(B_{\dR}^+(R^{\sharp}))/G(B_{\dR}^+(R^{\sharp}))
    \]
    for some $g \in G(B_{\dR}^+(R^{\sharp}))$.
\end{lem}
\begin{proof}
    Let $L^+G(\mu) := \ker(L^+G \to \Gra{G, \mu})$ be the stabilizer subgroup. Then $L^+G \to \Gra{G, \mu}$ is an $L^+G(\mu)$-torsor. It suffices to show that the orbit map $L^+G \to \Gra{G, \mu}$ has a section over $S$ after replacing $S$ with an open subspace. Consider the reduction map $\red: L^+G \to G$, $g \mapsto \bar{g}$, modulo $\xi$. Let $L^1G(\mu)$ be the kernel of the restriction of $\red$ to $L^+G(\mu)$, and let $G_{\mu} \subset G$ be the image of $L^+G(\mu)$. We have an exact sequence $1 \to L^1G(\mu) \to L^+G(\mu) \to G_{\mu} \to 1$. To show that $H^1_v(S, L^+G(\mu)) = 0$ after replacing $S$ with an open subspace, it suffices to show that $H^1_v(S, L^1G(\mu)) = 0$ and $H^1_v(S, G_{\mu}) = 0$ after replacing $S$ with an open subspace.

    (1) $H^1_v(S, G_{\mu}) = 0$: Note that $G_{\mu}$ is a parabolic subgroup of $G$ determined by $\mu$; in particular, it is smooth. A $G_{\mu}$-torsor on a perfectoid space $S$ for the $v$-topology is a torsor for the \'etale topology. A general result for rigid groups can be found in \cite{heuer2022g}. Since every \'etale covering of $S$ splits, after replacing $S$ with an open subspace, we have $H^1_v(S, G_{\mu}) = 0$.

    (2) $H^1_v(S, L^1G(\mu)) = 0$: Recall that $G$ is split over $C$, and $L^1G(\mu)$ is a subgroup of the pro-unipotent kernel $\ker(L^+G \to G)$. More precisely, $L^1G(\mu)$ admits a descending filtration $L^{\bullet}G(\mu)$ such that $L^nG(\mu)/L^{n+1}G(\mu)$ is isomorphic to some vector group $\mathfrak{g}^{(n)} \otimes \xi^n/\xi^{n+1}$, where $\mathfrak{g}^{(1)} \subset \mathfrak{g}^{(2)} \subset \cdots \mathfrak{g} = \Lie G$ is an increasing filtration. Since $H^1_v(S, B_{\dR}^+/(\xi)) = H^1_v(S, \OO_S) = 0$, and the transition morphisms in the inverse system of vector groups $\lrbracket{L^1G(\mu)/L^{n}G(\mu)}_n$ are surjective, we conclude that $H^1_v(S, L^1G(\mu)) = 0$.\end{proof}
\subsubsection{Over $\Qp$}
Let $\lrbracket{\mu}$ be the conjugacy class of $\mu$ in $P$, and $E$ be the defining field of $\lrbracket{\mu}$. Note that the defining field of $\lrbracket{\mu}$ in $P$ is the same as the one of $\lrbracket{\mu}$ in $G$, since we can fix a section $G \to P$ such that $\mu$ lifts to $G$. Since $\Gra{G, \mu}$ and $\Gra{G, \leq \mu}$ are descended to $\Spd E$ (which we again denote by $\Gra{G, \mu}$ and $\Gra{G, \leq \mu}$ respectively), $\Gra{P, \mu}$ and $\Gra{P, \leq \mu}$ are descended to $\Spd E$  (which we again denote by $\Gra{P, \mu}$ and $\Gra{P, \leq \mu}$ respectively), by the descriptions in Lemma \ref{lem: P and G associated, Gr}. Moreover, $\Gra{P, \mu} \hookrightarrow \Gra{P, \leq \mu}$ is open and $\Gra{P, \leq \mu} \subset \Gra{P, E}$ is closed. There is a section $\Gra{G, \leq \mu} \to \Gra{P, \leq \mu}$ to the projection $\Gra{P, \leq \mu} \to \Gra{G, \leq \mu}$.

\subsubsection{Represented by rigid analytic spaces}
Let
\[
    P_{\mu} := \lrbracket{g \in P \mid \lim_{t \to \infty} \mu(t)g\mu(t)^{-1}\ \textit{exists}},
    \quad
    G_{\mu} := \lrbracket{g \in G \mid \lim_{t \to \infty} \mu(t)g\mu(t)^{-1}\ \textit{exists}}.
\]
Let $\fl_{P, \mu}:= P/P_{\mu}$ and $\fl_{G, \mu} = G/G_{\mu}$. Note that $\fl_{G, \mu}$ is the usual flag variety. The space $\fl_{P, \mu}$ can be also viewed as a smooth quasi-projective scheme over $E$.

\begin{definition}\leavevmode\label{def-minuscule-parabolic}
    \begin{enumerate}
        \item We say $\mu$ is \textbf{quotient-minuscule} in $P$ if $\mu$ is minuscule in the reductive quotient $G$.
        \item We say $\mu$ is \textbf{minuscule} in $P$ if the weights of $\mu$ on $\Lie P$ are $\lrbracket{-1, 0, 1}$.
    \end{enumerate}
\end{definition}
It follows from the definition that $\mu$ is minuscule implies that $\mu$ is quotient-minuscule.
\begin{lem}[{\cite[Prop. 19.4.2]{SW20}}]\label{lem: BB map}
    There is a Bialynicki-Birula map $\Gra{P, \mu} \to \fl_{P, \mu}^{\Dia}$, which is an isomorphism when $\mu$ is minuscule in $P$.
\end{lem}
\begin{proof}
    Consider the commutative diagram
    \[
\begin{tikzcd}
	{\Gra{P, \mu}} & {\fl_{P, \mu}^{\Dia}} \\
	{\Gra{G, \mu}} & {\fl_{G, \mu}^{\Dia},}
	\arrow[from=1-1, to=1-2]
	\arrow[from=1-1, to=2-1]
	\arrow[from=1-2, to=2-2]
	\arrow[from=2-1, to=2-2]
\end{tikzcd}
    \]
    where the horizontal maps are defined using the Tannakian formalism and \cite[Prop. 19.4.2]{SW20}, and the bottom map is an isomorphism. Since $\fl_{P, \mu}^{\Dia}$ and $\Gra{P, \mu}$ are diamonds, to check the qcqs map $\Gra{P, \mu} \to \fl_{P, \mu}^{\Dia}$ is an isomorphism, it suffices to check the bijectivity on geometric points by \cite[Lem. 11.11, Prop. 11.15]{scholze2017etale}. 
    In particular, we check the bijectivity of the fibers at every geometric point of $\Gra{G, \mu} \rightiso \fl_{G, \mu}^{\Dia}$. Let $\Spa(C, C^+)$ be a geometric point of characteristic $0$, denote $\OO := B_{\dR}^+(C) \cong C[[\xi]]$ and $F = B_{\dR}(C) \cong C((\xi))$. We have identifications:
    \[  \Gra{P, \mu}(C, C^+) = P(\OO)\xi^{\mu}P(\OO)/P(\OO) = P(\OO)/P(\OO) \cap \xi^{\mu}P(\OO)\xi^{-\mu}. \]
    The natural morphism $\Gra{P, \mu}(C, C^+) \to \fl_{P, \mu}^{\Dia}(C, C^+)$ is induced by the reduction map
    \[ P(\OO) \to P(C),\quad q\mapsto \bar{q} \mod \xi. \]
    Let $g \in G(C)$, the fiber of $\Gra{P, \mu} \to \Gra{G, \mu}$ at  $g$ (more explicity, over $g\xi^{\mu}G(\OO)$) is
\begin{equation}\label{eq: fibers at grassmannian}
  U(\OO)/(U(\OO)\cap g\xi^{\mu}U(\OO) \xi^{-\mu}g^{-1}) \cong U(\OO)/(U(\OO)\cap \xi^{\mu}U(\OO)\xi^{-\mu}).
\end{equation}
    and the fiber of $\fl_{P, \mu}^{\Dia} \to \fl_{G, \mu}^{\Dia}$ at $g$ is
    $$ U(C)/gU_{\mu}(C)g^{-1} \cong U(C)/U_{\mu}(C), $$
    where $\Lie U_{\mu}$ consists of $\mu$-weight spaces of $\Lie U$ of weights $\leq 0$. By the minuscule condition on $\mu$, the reduction map gives the bijection
    $$ U(\OO)/(U(\OO)\cap \xi^{\mu}U(\OO)\xi^{-\mu}) \rightiso U(C)/U_{\mu}(C). $$
\end{proof}

\begin{definition}\label{def: P comes from boundary}
    We say that $(P, \mu)$ \textbf{comes from the boundary} if there exists a reductive $\rQ_p$-group $G'$ with a parabolic $\rQ_p$-subgroup $Q \subset G'$ such that 
    \begin{enumerate}
        \item $\mu: \Gm \stackrel{\mu}{\to} P \to G'$ is a minuscule cocharacter in $G'$ defined over a field extension of $E$,
        \item $\lrbracket{ U_{\alpha}|\ \lrangle{\mu, \alpha} > 0, \alpha \in \Phi'} \subset Q$, where $\Phi'$ is the absolute root datum of $G'$.
        \item $P$ is a normal subgroup of $Q$ with $R_uP = R_uQ$.
    \end{enumerate}
    In this case, we denote the Levi quotient of $Q$ by $L$.
\end{definition}
\begin{lem}\label{lem: boundary Shimura datum are from boundary}
    Let $(G', X')$ be a pure Shimura datum and let $(P, D)$ be a rational boundary component. Let $x' \in D$, and let $\mu_{x'}$ be the Hodge cocharacter associated with $x'$. Then $(P_{\rQ_p}, \mu_{x', \ovl{\rQ}_p})$ comes from the boundary in the sense of Definition~\ref{def: P comes from boundary}.
\end{lem}
\begin{proof}
    The groups involved are defined over $\rQ$, we work over the base field $\CC$ in place of $\ovl{\rQ}_p$. Let $H_0$ be the reference group defined in \cite[\S 4.3]{Pin89}. It is defined over $\R$ and equipped with two homomorphisms $h_0: \DS \to H_0$ and $h_{\infty}: \DS_{\CC} \to \DS_{\CC} \times H_{0, \CC}$. Let $Q$ be an admissible parabolic subgroup of $G'$. For any $x = h_x \in X'$ that determines a point $x' = h_{x'} \in D$ in the sense of \cite[\S 4.11]{Pin89}, one can find a morphism from the reference group, $\omega_{x}: H_{0} \to G'_{\R}$, such that $h_{x} = \omega_{x}\circ h_0: \DS \to G_{\R}'$ and $h_{x'} = \omega_{x}\circ h_{\infty}: \DS_{\CC} \to P_{\CC}$; see \cite[Prop. 4.6]{Pin89}.
    
    Let $w: \Gm \to \DS$ and $\mu: (\Gm)_{\CC} \to \DS_{\CC}$ be the weight and Hodge homomorphisms, respectively. Let $\mu_{x} = h_{x}\circ \mu$ and $\mu_{x'} = h_{x'} \circ \mu$ be the Hodge cocharacters associated with $x \in X'$ and $x' \in D$, respectively. Then $\mu_{x}$ and $\mu_{x'}$ are conjugate under $G'$ (indeed, they are already conjugate in $H_0$); see \cite[Prop. 12.1, Lem. 12.2]{Pin89}. In particular, $\mu_{x'}$ is minuscule.
    
    Since the Hodge structure on $\Lie G'$ induced by $h_{x}$ is pure of type $\lrbracket{(-1, 1), (0, 0), (1, -1)}$, \cite[Lem. 4.4]{Pin89} implies that the Hodge structure induced by $h_{x'}$ is mixed, with types
    \begin{equation}\label{eq: type of weights on Lie G'}
        \lrbracket{(1, 1), (1, 0), (0, 1), (1, -1), (0, 0), (-1, 1), (0, -1), (-1, 0), (-1, -1)}.
    \end{equation}
    The adjoint action of $\mu_{x'}: \Gm \to P_{\CC} \to G'_{\CC}$ on $\Lie U_{\alpha}$ is given by multiplication by $t^{\lrangle{\mu_{x'}, \alpha}}$. Thus $\Lie U_{\alpha} \subset F^0\Lie G'$ ($F^0$ is the Hodge filtration induced by $h_{x'}$) if and only if $\lrangle{\mu_{x'}, \alpha} \leq 0$ (see conventions in \cite[\S 1.3]{Pin89}). In particular, if $\lrangle{\mu_{x'}, \alpha}>0$, then the weight of $\Lie U_{\alpha}$ under $h_{x'}\circ w$ has to be at least $0$ by \eqref{eq: type of weights on Lie G'}. Hence $\Lie U_{\alpha} \subset (\Lie Q)_{\CC}$, since $(\Lie Q)_{\CC}$ is the direct sum of all nonnegative weight spaces in $(\Lie G')_{\CC}$ under the adjoint action of $h_{x'}\circ w$; see \cite[4.1 and Prop. 4.6]{Pin89}.
\end{proof}

\begin{lem}\label{lem: flag varieties have same dimensions}
    If $(P, \mu)$ comes from the boundary, then the immersions $\fl_{P, \mu} \to \fl_{Q, \mu} \to \fl_{G', \mu}$ are open immersions; in particular, they have the same dimension.
\end{lem}
\begin{proof}
    This follows from a calculation of the tangent spaces at the identity. Since the construction of $\fl$ is compatible with base change, we base change everything to an algebraically closed field. We take a maximal torus $T$ inside $Q$,
    \[ \Lie G/\Lie G_{\mu} = \oplus_{\alpha \in \Phi, \lrangle{\mu, \alpha} > 0} \Lie U_{\alpha} = \Lie Q/\Lie Q_{\mu}. \]
    Let us denote $W:=R_uP = R_uQ$, $G_h:=P/W$, $L:=Q/W$, and $G_l:= Q/P$. We have an exact sequence of reductive groups $1 \to G_h \to L \to G_l \to 1$, and the root data of $G_h$ and $G_l$ are orthogonal. In particular, $\Lie L/\Lie L_{\mu} = \Lie G_h/\Lie G_{h, \mu}$; thus $\Lie Q/\Lie Q_{\mu} = \Lie P/\Lie P_{\mu}$.
\end{proof}
\begin{rk}\label{rk: same dimension for mixed Shimura varieties}
    Under the setting of Lemma~\ref{lem: boundary Shimura datum are from boundary}, let us give a more geometric explanation of Lemma~\ref{lem: flag varieties have same dimensions}. Consider the map $\tau: X \to \pi_0(X) \times \Hom(\DS_{\CC}, P_{\CC})$ that sends $x$ to $([x], \omega_x\circ h_{\infty})$. Let $(P, D)$ be a rational boundary component. Recall that $D \subset \pi_0(X) \times \Hom(\DS_{\CC}, P_{\CC})$ is a $P(\R)U(\CC)$-orbit containing $\tau(x)$ for some $x \in X$. Let $\phi: D \to \phi(D) \subset \Hom(\DS_{\CC}, P_{\CC})$ be the natural projection with finite fibers. Fix $h \in \phi(D)$, let $C(h)$ be the centralizer of $h$ in $P(\R)W(\CC)$. Let $M$ be a faithful representation of $P$. By the proof of \cite[Prop. 1.7]{Pin89}, we have
    \[
        \phi(D) = P(\R)W(\CC)/C(h) \hookrightarrow P(\CC)/\mathrm{exp}(F^0(\Lie P)_{\CC}) = P(\CC) / P_{\mu}(\CC) \hookrightarrow \mathrm{Grass}(M)(\CC),
    \]
    where the first morphism is an open embedding, and the last morphism is the closed embedding given by the Hodge filtration associated with the Hodge cocharacter $\mu$ determined by $h$. In particular, the dimension of $D$ agrees with that of $\fl_{P, \mu}$ as complex manifolds. This is analogous to the Borel embedding $X' \hookrightarrow \Tilde{X}' = G'/G'_{\mu'} = \fl_{G', \mu'}$, which also shows that the dimension of $X'$ equals that of $\fl_{G', \mu'}$. Let $X' \subset X$ be the preimage of $D$ under $\tau$, which is a union of connected components of $X$. Then $X' \to D$ is an open embedding; see \cite[Prop. 4.15(a)]{Pin89}. Thus $X'$ and $D$ have the same dimension. Finally, since $\mu$ and $\mu'$ are conjugate in $G'$, we conclude that $\fl_{P, \mu}$ and $\fl_{G', \mu}$ have the same dimension.
\end{rk}

\begin{lem}\label{lem: BB map for general group}
    If $(P, \mu)$ comes from the boundary, then $\Gra{P, \mu} \to \fl_{P, \mu}^{\Dia}$ and $\Gra{Q, \mu} \to \fl_{Q, \mu}^{\Dia}$ are isomorphisms. In particular, $\Gra{P, \mu}$ and $\Gra{Q, \mu}$ are represented by rigid analytic spaces over $\Spd E$. Moreover, $\Gra{P, \mu} \rightiso \Gra{Q, \mu}$ and $\fl_{P, \mu} \rightiso \fl_{Q, \mu}$.
\end{lem}
\begin{proof}

    The first two statements follow from Lemma \ref{lem: BB map}. Let us show the last statement. Given any $q \in Q(\OO)$, we can factor $q = p\prod_{\alpha}u_{\alpha}$ with $p \in P(\OO)$ and $u_{\alpha} \in U_{\alpha}(\OO)$, where $\alpha$ runs over roots in $Q\setminus P$. Since $\mu$ factors through $P$, we have $\lrangle{\alpha, \mu} = 0$ for all such $\alpha$. In particular, $U_{\alpha}(\OO) \subset Q(\OO) \cap \xi^{\mu}Q(\OO)\xi^{-\mu}$, and hence $\Gra{P, \mu}(C, C^+) = \Gra{Q, \mu}(C, C^+)$, $\Gra{P, \mu} \rightiso \Gra{Q, \mu}$. Since $\fl_{P, \mu}$ and $\fl_{Q, \mu}$ are smooth, and the $\Dia$ functor is fully faithful on seminormal rigid analytic spaces (\cite[Prop. 10.2.3]{SW20}), $\fl_{P, \mu} \rightiso \fl_{Q, \mu}$.
\end{proof}

\subsubsection{Over $\Spd \Z_p$}

We move on to the integral base and recall the definition of $v$-sheaf theoretical local models.

\begin{definition}\label{def: quasi-parahoric for non-reductive group}
    We say that $\PP$ is a \textbf{stabilizer quasi-parahoric group scheme} (resp. \textbf{quasi-parahoric group scheme}, \textbf{parahoric group scheme}) of $P$ if $\PP = \UU \rtimes \GG$, $\GG$ is a stabilizer quasi-parahoric group scheme (resp. quasi-parahoric group scheme, parahoric group scheme) of $G$ and $\UU$ is a smooth affine group scheme over $\Z_p$ with connected fibers. In these cases, denote $\PP^{\circ} = \UU \rtimes \GG^{\circ}$.
\end{definition}

Let $Q \subset G'$ be a parabolic subgroup, and let $T \subset Q$ be a maximal $\bQ$-split torus of $G'$ defined over $\rQ_p$. Let $U = R_uQ$, and let $L$ be the standard Levi subgroup of $Q$ containing $T$. 

Let $\FF \subset A(G', T)$ be a facet in the reduced apartment. Under the natural projection $A(G', T) \to A(L, T)$, there is a unique facet $\FF_L \subset A(L, T)$ containing $\FF$. Let $\LL_{\FF_L}$ (resp. $\GG'_{\FF}$) be the parahoric group scheme of $L$ (resp. $G'$) with respect to $\FF_L$ (resp. $\FF$). Then $\LL_{\FF_L}(\bZ_p) = L(\bQ) \cap \GG'_{\FF}(\bZ_p)$; see \cite[Lem. 4.1.1]{haines2010satake}. Moreover, $\LL_{\FF_L} \to \GG'_{\FF}$ is a closed embedding. To see this, let $\LL_{\FF}$ be the closure of $L$ in $\GG'_{\FF}$. Then $\LL_{\FF_L} \to \LL_{\FF}$ is the smoothing; i.e., $\LL_{\FF_L}$ is the unique smooth affine group scheme of $L$ such that $\LL_{\FF_L}(\bZ_p) = \LL_{\FF}(\bZ_p)$. Since $\LL_{\FF}$ contains the open big cell of $\LL_{\FF_L}$, we have $\LL_{\FF_L} \rightiso \LL_{\FF}$. Let $\UU_{\FF}$ be the closure of $U$ in $\GG'_{\FF}$; it is a smooth affine group scheme with connected fibers, since its special fiber is a product of root groups $U_{\alpha}$ inside the linear algebraic group $\GG'_{\FF, \F_p}$. The closure of $Q$ in $\GG'_{\FF}$ is $\QQ_{\FF} := \UU_{\FF} \rtimes \LL_{\FF}$, and $\QQ_{\FF}$ is parahoric.

\begin{definition}\label{def: parabolic embedding}
    We say that an embedding of parahoric group schemes $\QQ \subset \GG'$ is a \emph{parabolic embedding} if there exists a facet $\FF \subset A(G', T)$ for some maximal $\bQ$-split torus $T$ of $G'$ contained in $Q$ such that $\QQ = \QQ_{\FF}$ and $\GG' = \GG'_{\FF}$.
\end{definition}

Let $\Gra{\GG', \leq \mu}$ be the closure of $\Gra{G', \leq \mu}$ in $\Gra{\GG'}$; it is proper and is represented by a spatial diamond (see \cite[Prop. 20.3.6]{SW20}).

\begin{lem}\label{lem: v-sheaf local model, Q to G, locally closed}
     Let $\QQ \subset \GG'$ be a parabolic embedding of parahoric group schemes, and let $\mu$ be a cocharacter of $G'$ factoring through $Q$. Then $\Gra{\QQ, \leq\mu} \to \Gra{\GG', \leq\mu}$ is a locally closed embedding, where $\Gra{\QQ, \leq\mu}$ is the closure of the closed subfunctor $\Gra{Q, \leq\mu} \subset \Gra{Q}$ (see Proposition~\ref{prop: leq mu is a closed embedding}) in $\Gra{\QQ, \OO_E}$, and $E$ is the defining field of $\lrbracket{\mu}$.
\end{lem}
\begin{proof}
    Denote $G'$ by $G$. We apply the arguments in \cite[\S 5.1]{anschutz2022p}. Let $T \subset Q \subset G$ be a maximal $\bQ$-split torus defined over $\rQ_p$ such that $\FF \subset A(G, T)$. Let $\mathcal{T}^{\circ}$ be the connected N\'eron model of $T$, and let $\lambda: \Gm \to \mathcal{T}^{\circ}$ be a cocharacter defined over $\Z_p$ such that $Q = G_{\lambda}$ and $L = G_{\lambda} \cap G_{\lambda^{-1}}$. As in the proof of \cite[Lem. 4.5]{HainesRicharz2021Test}, $\lambda$ defines $\QQ = \GG_{\lambda}$ and $\mathcal{L} = \GG_{\lambda} \cap \GG_{\lambda^{-1}}$. Choose an embedding $\GG \hookrightarrow \GL_{n, \Z_p}$ such that $\GL_{n, \Z_p}/\GG$ is quasi-affine. Let $T' \subset G':= \GL_n$ be a maximal $\bQ$-split torus defined over $\rQ_p$ that contains $T$, and let $Q':= G'_{\lambda'}$, where $\lambda'$ is the composition of $\Gm \stackrel{\lambda}{\to} \mathcal{T}^{\circ} \to \mathcal{T}^{\prime\circ}$. Then $Q' = U' \rtimes L'$ extends to $\QQ' = \UU' \rtimes \LL'$ where $\QQ'$ is again a parahoric group scheme. As in the arguments in the proof of \cite[Thm. 5.2]{anschutz2022p}, $\QQ'/\QQ$ is quasi-affine; thus $\Gra{\QQ} \to \Gra{\QQ'}$ is a locally closed immersion. By the discussion above Lemma~\ref{lem: v-sheaf local model, Q to G, locally closed}, it suffices to work with $\QQ' \subset \GG'$. To save notation, we denote by $(\QQ', \GG', \mu') = (\QQ, \GG, \mu)$. Now $\QQ = \UU \rtimes \LL$, and $\LL$ and $\GG$ are reductive schemes over $\Z_p$ since $T$ is a split maximal torus. Let $T_L$ be the cocenter of $L$, and let $\nu$ be the image of $\mu$ under $L \to T_L$. We use results from \cite[Prop. III.3.6]{fargues2021geometrization} (the proof there works over the integral base, since \cite[Lem. VI.3.2, VI.3.3]{fargues2021geometrization} do.) Recall that we have a locally constant Kottwitz map $\Gra{\LL} \to \pi_1(L)$. Let $\Gra{\LL}^{\nu} \subset \Gra{\LL}$ be the corresponding open and closed sub-sheaf, and let $\Gra{\QQ}^{\nu} \subset \Gra{\QQ}$ be the preimage of $\Gra{\LL}^{\nu} \subset \Gra{\LL}$. By \cite[Prop. VI.3.1]{fargues2021geometrization}, $\Gra{\QQ}^{\nu} \subset \Gra{\GG}$ is locally closed. Since $\Gra{\LL, \leq \mu} \subset \Gra{\LL}^{\nu}$, $\Gra{\QQ, \leq \mu} \subset \Gra{\QQ}$ and $\Gra{\GG, \leq \mu} \subset \Gra{\GG}$ are closed embeddings, $\Gra{\QQ, \leq\mu} \to \Gra{\GG, \leq\mu}$ is locally closed.
     \end{proof}
     
Let $G_h \subset L$ be a normal subgroup with quotient $G_l$. Then $G_h$ and $G_l$ are reductive. Moreover, $L$ is an almost product of $G_h$ and $G_l'$, where $G_l' \to G_l$ is an isogeny. Recall that we have a $G_h(\bQ) \times G_l'(\bQ)$-$\sigma$-equivariant bijection between reduced Bruhat--Tits buildings that preserves the polysimplicial structure:
\begin{equation}\label{eq: bijection on buildings}
   \Bui_{\red}(G_h, \bQ) \times \Bui_{\red}(G_l, \bQ) \to \Bui_{\red}(L, \bQ),\quad (\FF_h, \FF_l) \leftrightarrow \FF_L, \quad   (x_h, x_l) \leftrightarrow x_L.
\end{equation}

In particular, $\GG_{h, x_h}(\bZ_p) = G_h(\bQ) \cap \LL_{x_L}(\bZ_p)$. However, the parahoric group schemes $\GG_{h, \FF_h}$ and $\LL_{\FF_L}$ might not satisfy $\GG_{h, \FF_h}(\bZ_p) = G_h(\bQ) \cap \LL_{\FF_L}(\bZ_p)$; see \cite[\S 2.4.3]{Mao25b} (especially \cite[Example 2.61]{Mao25b}) for detailed discussions. Moreover, the closure of $G_h$ in $\LL_{x_L}$ might not be smooth. It is smooth if there exists a maximal $\bQ$-split torus of $G_h$ that is $R$-smooth; see \cite[Lem. 2.77]{Mao25b}.

We make the following definition, which will be used when considering mixed Shimura varieties at the boundary.
\begin{definition}\label{def: PP, mu comes from boundary}
    Assume that $(P, \mu)$ comes from the boundary, with $P \to Q \to G'$ as usual. Let $\PP \to \QQ \to \GG'$ be quasi-parahoric group schemes of $P \to Q\to G'$. We say that $(\PP, \mu)$ comes from the boundary associated with $\PP \to \QQ \to \GG'$ if
    \begin{enumerate}
        \item $\QQ$ (resp. $\PP$) is the smoothing of the closure of $Q$ (resp. $P$) in $\GG'$.
        \item $\PPc \to \QQc \to \GG^{\prime\circ}$ is of the form $\PP_{\FF_h} \to \QQ_{\FF_L} \to \GG^\prime_{\FF}$, with facets $\FF_h$, $\FF_L$, $\FF$ related in the way discussed above.
    \end{enumerate}
\end{definition}

\subsubsection{Pass to parahoric levels}\label{subsubsec-pass-to-parahoric}

In the remainder of this section, we always assume that $\PP = \UU \rtimes \GG$ is a quasi-parahoric group scheme.

\begin{lem}\label{lem: Grassmannian, cartesian diagram}
\[
\begin{tikzcd}
	{\Gra{\PP^{\circ}, \Spd \Z_p}} & {\Gra{\PP, \Spd \Z_p}} & {\Gra{\PP^{\circ}}^W} & {\Gra{\PP}^W} \\
	{\Gra{\GG^{\circ}, \Spd \Z_p}} & {\Gra{\GG, \Spd \Z_p},} & {\Gra{\GG^{\circ}}^W} & {\Gra{\GG}^W.}
	\arrow[from=1-1, to=1-2]
	\arrow[from=1-1, to=2-1]
	\arrow[from=1-2, to=2-2]
	\arrow[from=1-3, to=1-4]
	\arrow[from=1-3, to=2-3]
	\arrow[from=1-4, to=2-4]
	\arrow[from=2-1, to=2-2]
	\arrow[from=2-3, to=2-4]
\end{tikzcd}
\]
    are Cartesian diagrams.
\end{lem}
\begin{proof}
The commutativity of both diagrams follows from the fact that the sheafification functor commutes with fiber products, and the diagram of presheaves is clearly Cartesian (since $\PP \to \GG$ is surjective).
\end{proof}

\begin{lem}\label{lem: fibers of Gra P to G are connected}
    The fibers of $\Gra{\PP}^W \to \Gra{\GG}^W$ are connected. Moreover, $\Gra{\PP}^W \to \Gra{\GG}^W$ induces a bijection $\pi_0(\Gra{\PP}^W) \to \pi_0(\Gra{\GG}^W)$.
\end{lem}
\begin{proof}
    Recall that $L\UU = LU$ is connected as an ind-scheme by fixing some $U \cong \mathbb{A}^n$, and $L^+\UU$ is also connected when $\UU$ is connected. Hence, $\Gra{\UU}^W$ is connected as an ind-scheme. Let $x = \Spec l \to \Gra{\GG}^W$ be a geometric point, and write $x = g\GG(W(l))/\GG(W(l))$. Then the fiber of $\Gra{\PP}^W \to \Gra{\GG}^W$ at $x$ is $U(W(l)[1/p])/g\UU(W(l)[1/p])g^{-1}$, which is conjugate to $\Gra{\UU}^W$. Moreover, fixing a section $\GG \to \PP$, since the action $LU \times \Gra{\GG}^W \to \Gra{\PP}^W$ is surjective, $\Gra{\PP}^W \to \Gra{\GG}^W$ induces a bijection $\pi_0(\Gra{\PP}^W) \to \pi_0(\Gra{\GG}^W)$.
\end{proof}

\begin{cor}[{\cite[Prop. 1.21]{zhu2017affine}}]\label{cor: structure of Gra PPcirc to Gra PP}
    $\Gra{\PP^{\circ}}^W \to \Gra{\PP}^W$ is a $\PP/\PP^{\circ}$-torsor that maps connected components isomorphically onto connected components.
\end{cor}
\begin{proof}
    This statement is true when $\PP = \GG$; see \cite[Prop. 1.21]{zhu2017affine} (cf. \cite[Prop. 21.1.4]{SW20}). In general, $\Gra{\PP^{\circ}}^W \to \Gra{\PP}^W$ is a $\PP/\PP^{\circ}$-torsor since $\Gra{\GG^{\circ}}^W \to \Gra{\GG}^W$ is a $\GG/\GG^{\circ}$-torsor and $\PP/\PP^{\circ} = \GG/\GG^{\circ}$ by Lemma \ref{lem: Grassmannian, cartesian diagram}. The last statement then follows from Lemma \ref{lem: fibers of Gra P to G are connected}.
\end{proof}

\subsection{Local models}

In the remainder of this section, we always assume that $\mu$ is \textbf{quotient-minuscule} in the sense of Definition \ref{def-minuscule-parabolic}. We define 
\begin{equation}\label{def: local model}
    \vM_{\PP, \mu}^v \subset \Gra{\PP, \Spd \OO_E} = \Gra{\PP, \Spd \Z_p} \times_{\Spd \Z_p} \Spd \OO_E
\end{equation}
as the closure of $\Gra{P, \mu}$. The closed embedding $\Gra{\GG} \to \Gra{\PP}$ induces $\vM_{\GG, \mu}^v \to \vM_{\PP, \mu}^v$, which gives a section of the projection $\vM_{\PP, \mu}^v \to \vM_{\GG, \mu}^v$, thus $\vM_{\PP, \mu}^v$ projects onto $\vM_{\GG, \mu}^v$.

\begin{lem}[{\cite[Prop. 21.4.3]{SW20}}]\label{lem: isomorphism of local model, quasi-parahoric}
    The natural morphism $\vM_{\PP^{\circ}, \mu}^v \to \vM_{\PP, \mu}^v$ is an isomorphism.
\end{lem}
\begin{proof}
    This statement is true when $\PP = \GG$; see \cite[Prop. 21.4.3]{SW20}. In general, we use Lemma \ref{lem: Grassmannian, cartesian diagram}. Since $\Gra{\GG^{\circ}, \Spd \Z_p} \to \Gra{\GG, \Spd \Z_p}$ is proper, $\Gra{\PP^{\circ}, \Spd \Z_p} \to \Gra{\PP, \Spd \Z_p}$ is also proper. In particular, $\vM_{\PP^{\circ}, \mu}^v \to \vM_{\PP, \mu}^v$ is proper. Applying the arguments in the proof of \cite[Prop. 21.4.3]{SW20}, and using Corollary \ref{cor: structure of Gra PPcirc to Gra PP}, it suffices to show that the reduction of $\vM_{\PP^{\circ}, \mu}^v$ (see Lemma \ref{lem: reduction of local model}) is completely contained in one connected component of $\Gra{\PP^{\circ}}^W$. Again, this follows from Corollary \ref{cor: structure of Gra PPcirc to Gra PP} and the fact that $\vM_{\PP^{\circ}, \mu}^v$ projects onto $\vM_{\GG^{\circ}, \mu}^v$.
\end{proof}

\begin{prop}
    Assume $(\PP, \mu)$ comes from the boundary with $\PP \to \QQ \to \GG'$ (see Definition \ref{def: PP, mu comes from boundary}); then the natural embedding $\vM_{\PP, \mu}^v \rightiso \vM_{\QQ, \mu}^v \to \vM_{\GG', \mu}^v$ is a locally closed embedding.
\end{prop}
\begin{proof}
    Lemma \ref{lem: BB map for general group} says that $\Gra{P, \mu} \to \Gra{Q, \mu}$ is an isomorphism. Since $\QQ^{\circ}/\PP^{\circ}$ is affine, $\Gra{\PP^{\circ}} \to \Gra{\QQ^{\circ}}$ is a closed embedding, see Lemma \ref{lem: quotient affine, closed embedding}. In particular, $\vM_{\PP^{\circ}, \mu}^v \rightiso \vM_{\QQ^{\circ}, \mu}^v$. By Lemma \ref{lem: isomorphism of local model, quasi-parahoric}, we have $\vM_{\PP, \mu}^v \rightiso \vM_{\QQ, \mu}^v$. By Lemma \ref{lem: v-sheaf local model, Q to G, locally closed} and Lemma \ref{lem: isomorphism of local model, quasi-parahoric}, $\vM_{\QQ, \mu}^v \to \vM_{\GG', \mu}^v$ is a locally closed embedding.
\end{proof}

\subsubsection{Orbital descriptions}
\begin{lem}[{\cite[Prop. 4.13]{anschutz2022p}}]\label{lem: orbit action closed}
    The local model $\vM_{\PP, \mu}^v$ is stable under $L^+\PP$-action. Moreover, $\vM_{\PP, \mu}^v$ is the $L^+\PP$-orbit of $\vM_{\GG, \mu}^v$ under the fixed closed embedding $\Gra{\GG} \subset \Gra{\PP}$. 
\end{lem}
\begin{proof}
Choosing a faithful embedding $\rho: \PP \to \GL_n$ such that $\GL_n/\PP$ is quasi-affine, we have a locally closed embedding $\Gra{\PP} \to \Gra{\GL_n}$. 
Let $X \subset \Gra{\GL_n, \OO_E}$ be a quasi-compact and quasi-separated closed subsheaf that contains $\vM_{\PP, \mu}^v$ (for example, take $X = \Gra{\GL_n, \leq \rho(\mu), \OO_E}$). Since $\vM_{\PP, \mu}^v \to X$ is a monomorphism, as in Lemma \ref{lem: congruence subgroup}, $a: L_{\OO_E}^+\PP \times \vM_{\PP, \mu}^v \to \Gra{\PP, \OO_E}$ factors through $a_N: L_{\OO_E}^{\leq N}\PP \times \vM_{\PP, \mu}^v \to \Gra{\PP, \OO_E}$ for some large $N$, where $N$ is determined by $X$. Let $X_1 = X \cap \Gra{\PP, \OO_E} \hookrightarrow \Gra{\PP, \OO_E}$, it is a closed sub-sheaf, $L_{\OO_E}^+\PP$ acts on $X_1$ and factors through $L_{\OO_E}^{\leq N}\PP$. Since $L_{\OO_E}^{\leq N}\PP \to \Spd \OO_E$ is partially proper and open, by \cite[Cor. 2.9]{anschutz2022p} (which says that taking closure commutes with base change along partially proper and open morphisms), base change along the projection $L_{\OO_E}^{\leq N}\PP \times_{\Spd \OO_E} X_1 \to X_1$, we have
    \begin{equation*}
       L_{\OO_E}^{\leq N}\PP \times_{\Spd \OO_E} \vM_{\PP, \mu}^v = (L_E^{\leq N}\PP \times_{\Spd E} \Gra{P, \mu})^{\cl},
    \end{equation*}   
    as a closed subsheaf of $L_{\OO_E}^{\leq N}\PP \times_{\Spd \OO_E} X_1$. Since $L_E^{\leq N}\PP \times_{\Spd E} \Gra{P, \mu} \to X_{1, E}$ factors through $\Gra{P, \mu}$ by construction, combining with the equality displayed, we see that $L_{\OO_E}^{\leq N}\PP \times_{\Spd \OO_E} \vM_{\PP, \mu}^v \to X_1$ factors through $\vM_{\PP, \mu}^v$.

    Let us show the second statement. Since $\vM_{\GG, \mu}^v \subset \vM_{\PP, \mu}^v$, $\vM_{\PP, \mu}^v$ is stable under $L^+\PP$-action implies that the morphism $a': L_{\OO_E}^+\PP \times_{\Spd \OO_E} \vM_{\GG, \mu}^v \to \Gra{\PP, \Spd \OO_E}$ factors through $\vM_{\PP, \mu}^v$. We need to show it is surjective. Since $\vM_{\GG, \mu}^v$ is quasi-compact and quasi-separated, apply Lemma \ref{lem: congruence subgroup}, $a'$ factors through some $a'_N$ for large $N$. Apply \cite[Cor. 2.9]{anschutz2022p} again. Note that $\Gra{G, \mu} \subset \Gra{P, \mu}$, as well as $\vM_{\GG, \mu}^v \subset \vM_{\PP, \mu}^v$, is a closed embedding. Base-changing along the projection $L^{\leq N}_{\OO_E}\PP \times_{\Spd \OO_E} \vM_{\PP, \mu}^v \to \vM_{\PP, \mu}^v$, we have
    \begin{equation*}
       L^{\leq N}_{\OO_E}\PP \times_{\Spd \OO_E} \vM_{\GG, \mu}^v = (L^{\leq N}_E\PP \times_{\Spd E} \Gra{G, \mu})^{\cl},
    \end{equation*}
    as a closed subsheaf of $L_{\OO_E}^{\leq N}\PP \times_{\Spd \OO_E} \vM_{\PP, \mu}^v$. Since $L_{\OO_E}^{\leq N}\PP \times_{\Spd E} \Gra{G, \mu}$ surjects onto $\Gra{P, \mu}$ by construction, we apply Lemma \ref{lem: cl is surj} to the action $a_N'$.
\end{proof}

\subsubsection{Over $\Spd \F_p$}

Let $\GG$ be parahoric. Let $\lcM_{\GG, \mu} \subset \Gra{\GG, k_E}^W$ be the reduction of $\vM_{\GG, \mu}^v \subset \Gra{\GG, \OO_E}$. We have $\lcM_{\GG, \mu} = \ab_{\GG, \mu} \subset \Gra{\GG, k_E}^W$ (see \cite[Thm. 1.5]{anschutz2022p}). Here $\ab_{\GG, \mu}$ is the $\mu$-admissible locus in the sense of Kottwitz–Rapoport, which is the $k_E$ descent of
\[ \ab_{\GG, \mu, k} = \bigcup_{w \in \Adm_{G}(\lrbracket{\mu^{-1}})_{\GG(\bZ_p)}} \Gra{\GG, w, k}^W, \]
where $k = \bar{k}_E$, $\Gra{\GG, w, k}^W$ is the usual Schubert cell, it is the perfection of a quasi-projective scheme. 

Let $\PP$ be parahoric and $\wdt{W}$ be the Iwahori-Weyl group of $G$. For any $w \in \wdt{W}$, let $\Gra{\PP, w}^W \subset \Gra{\PP, k}^W$ be the orbit of $L^+_{k}\PP$ of $\Gra{\GG, w}^W \to \Gra{\GG, k}^W \to \Gra{\PP, k}^W$. By calculating the stabilizer of the orbit of $L^+_{k}\PP$ at $w$, we see that $\Gra{\PP, w}^W$ is locally perfectly of finite type. Let us define the $\mu$-admissible locus $\ab_{\PP, \mu} \subset \Gra{\PP, k_E}^W$ as the $k_E$-descent of
\begin{equation}\label{eq: P-grassmanian factorization}
    \ab_{\PP, \mu, k} = \bigcup_{w \in \Adm_P(\lrbracket{\mu^{-1}})_{\PP(\bZ_p)}} \Gra{\PP, w}^W, \quad \Adm_P(\lrbracket{\mu^{-1}})_{\PP(\bZ_p)}:= \Adm_G(\lrbracket{\mu^{-1}})_{\GG(\bZ_p)}.
\end{equation}

Let $\KK = \GG(\Z_p)$, and let $\prescript{}{\KK}{\wdt{W}}^{\KK} \subset \wdt{W}$ be the index set introduced in \cite{Richarz2013} (cf. \cite[\S 1.2.6]{shen2021ekor}).
\begin{lem}\label{lem: closure on witt grassmannian}
    Given $w \in \prescript{}{\KK}{\wdt{W}}^{\KK}$, let $\ovl{\Gra{\PP, w}^W}$ be the closure of $\Gra{\PP, w}^W \subset \Gra{\PP, k}^W$, then 
    \[ \ovl{\Gra{\PP, w}^W} = \Gra{\PP, \leq w}^W := \bigcup_{w' \leq w, w' \in \prescript{}{\KK}{\wdt{W}}^{\KK}} \Gra{\PP, w'}^W. \]
    In particular, $\ab_{\PP, \mu} \subset \Gra{\PP, k_E}^W$ is a closed subfunctor.
\end{lem}
\begin{proof}
    Since the closure of $\Gra{\GG, w}^W$ in $\Gra{\GG}^W$ is $\Gra{\GG, \leq w}^W$, the usual affine Schubert variety. Use the orbit map $L^+_{k}\PP \times \Gra{\GG, k}^W \to \Gra{\PP, k}^W$, we see that $\Gra{\PP, \leq w}^W\subset\ovl{\Gra{\PP, w}^W}$. On the other hand, by Lemma \ref{lem: congruence subgroup}, the action $a: L^+_k\PP \times \Gra{\GG, \leq w}^W \to \Gra{\PP, k}^W$ factors through $a_N$ for some large $N$, and $a_N$ is proper. In particular, $L^+_k\PP \times \Gra{\GG, \leq w}^W$ has closed image in $\Gra{\PP}$ that contains $\Gra{\PP, w}$, thus $\ovl{\Gra{\PP, w}^W} \subset \Gra{\PP, \leq w}^W$.
\end{proof}

\begin{lem}\label{lem: reduction of local model}
    Let $\lcM_{\PP, \mu} \subset \Gra{\PP, k_E}^W$ be the reduction of $\vM_{\PP, \mu}^v \subset \Gra{\PP, \OO_E}$. Then $\lcM_{\PP, \mu} \subset \Gra{\PP, k_E}^W$ is representable and closed.
\end{lem}
\begin{proof}
    We need to show that $\vM_{\PP, \mu}^v$ is $\varpi$-adic formal (i.e. the structure morphism to $\Spd \OO_E$ is formally adic in the sense of \cite[Def. 3.20]{gleason2025specialization}); then $\lcM_{\PP, \mu} \subset \Gra{\PP, k_E}^W$ represents the base change of $\vM_{\PP, \mu}^v \subset \Gra{\PP, \OO_E}$ along $\Spd \F_q \to \Spd \OO_E$. To show the claim, we follow the proof in \cite[Prop. 4.14]{anschutz2022p}: We need to show that the special fiber of $\vM_{\PP, \mu}^v$ is represented by a perfect scheme. This follows from Lemma \ref{lem: orbit action closed} and in particular (\ref{eq: P-grassmanian factorization}).
\end{proof}

By checking the set of $k$-points of these two closed functors, we have
\begin{cor}\label{cor: special fiber of local model}
    Let $\PP$ be parahoric. Then $\lcM_{\PP, \mu} = \ab_{\PP, \mu}$.
\end{cor}

\begin{rk}\label{rk: v-formalizing possibly fails}
    Recall that, $\vM_{\GG, \mu}^v$ is a proper topologically normal rich kimberlite, due to results in \cite{gleason2025specialization}, \cite{anschutz2022p}, \cite{gleason2024tubular}. One might expect that $\vM_{\PP, \mu}^v$ is also a topologically normal rich kimberlite, though not necessarily proper. Most arguments in the cited papers work here, except for one crucial ingredient, the purity of $\GG$-torsor proved in \cite{anschutz2022extending}. This purity result simply fails for non-reductive group, in particular, $\Gra{\PP}$ might not be $v$-formalizing, thus might not be pre-kimberlite.
\end{rk}

Let $\GG_0$ be the special fiber of $\GG$. The $L^+_{k_E}\GG$-action on $\lcM_{\GG, \mu}$ factors through $\GG_0$ (see \cite[\S 4.9.1]{PR24}). In particular, we have $|\GG_0\backslash \lcM_{\GG, \mu}| = \Adm_G(\lrbracket{\mu^{-1}})_{\GG(\bZ_p)}$. Similarly, let $\PP_0$ be the special fiber of $\PP$. Assume $(\PP, \mu)$ comes from the boundary, associated with $\PP \to \QQ \to \GG'$. Use the monomorphisms $\PP \to \GG'$ and $\vM_{\PP, \mu}^v \to \vM_{\GG', \mu}^v$, we have
\begin{cor}\label{cor: comes from boundary, action factors through G_0}
     Assume $(\PP, \mu)$ comes from the boundary. Then the $L^+_{k_E}\PP$-action on $\lcM_{\PP, \mu}$ factors through $\PP_0$. In particular, $|\PP_0\backslash \lcM_{\PP, \mu}| = \Adm_P(\lrbracket{\mu^{-1}})_{\PP(\bZ_p)}$.
\end{cor}

\subsubsection{Functoriality}

\begin{definition}\label{def: compatible, P and G}
    Let $P_1 = U_1 \rtimes G_1$, $P_2 = U_2 \rtimes G_2$. We say $P_1 \to P_2$ is compatible with $U_1 \rtimes G_1 \to U_2 \rtimes G_2$ if $P_1 \to P_2$ induces $U_1 \to U_2$, $G_1 \to G_2$, and there exist sections $G_1 \to P_1$ and $G_2 \to P_2$ such that $G_1 \to P_1 \to P_2$ factors through $G_2$.
\end{definition}
\begin{rk}\leavevmode
    \begin{itemize}
        \item In this case, $P_1 \to P_2$ is uniquely determined by the maps $U_1 \to U_2$ and $G_1 \to G_2$.
        \item When $P_1 \to P_2$ induces $U_1 \to U_2$, $G_1 \to G_2$, then one can always find such a pair of sections $G_1 \to P_1$ and $G_2 \to P_2$ in characteristic $0$. 
    \end{itemize}
\end{rk}

Let $P_1 \to P_2$ be compatible with $U_1 \rtimes G_1 \to U_2 \rtimes G_2$, and let $\PP_1 \to \PP_2$ be a morphism of quasi-parahoric group schemes that extends $P_1 \to P_2$; then $\PP_1 \to \PP_2$ induces $\UU_1 \to \UU_2$ and $\GG_1 \to \GG_2$. 
\begin{definition}\label{def: compatible, PP and GG}
    We say $\PP_1 \to \PP_2$ is compatible with $\GG_1 \to \GG_2$ if there exist sections $\GG_1 \to \PP_1$, $\GG_2 \to \PP_2$ such that $\GG_1 \to \PP_1 \to \PP_2$ factors through $\GG_2$.
\end{definition}

In this subsection, we consider the following case. Let $P_1 = U \rtimes G_1 \to P_2 = U \rtimes G_2$ be compatible with $U \rtimes G_1 \stackrel{(\identity, f)}{\to} U \rtimes G_2$. Let $\PP_1 = \UU \rtimes \GG_1$ and $\PP_2 = \UU \rtimes \GG_2$ be parahoric group schemes of $P_1$ and $P_2$, respectively, such that $P_1 \to P_2$ induces $\PP_1 \to \PP_2$ compatible with $\GG_1 \to \GG_2$.

Recall that the sheafification functor commutes with fiber products. Since $\PP_2 \to \GG_2$ is surjective with a section $\GG_2 \to \PP_2$, we have:
\begin{lem}
    The following diagrams are Cartesian:
\[
\begin{tikzcd}
	{\Gra{P_1}} & {\Gra{P_2}} & {\Gra{\PP_1}} & {\Gra{\PP_2}} & {\Gra{\PP_1}^W} & {\Gra{\PP_2}^W} \\
	{\Gra{G_1}} & {\Gra{G_2},} & {\Gra{\GG_1}} & {\Gra{\GG_2},} & {\Gra{\GG_1}^W} & {\Gra{\GG_2}^W.}
	\arrow[from=1-1, to=1-2]
	\arrow[from=1-1, to=2-1]
	\arrow[from=1-2, to=2-2]
	\arrow[from=1-3, to=1-4]
	\arrow[from=1-3, to=2-3]
	\arrow[from=1-4, to=2-4]
	\arrow[from=1-5, to=1-6]
	\arrow[from=1-5, to=2-5]
	\arrow[from=1-6, to=2-6]
	\arrow[from=2-1, to=2-2]
	\arrow[from=2-3, to=2-4]
	\arrow[from=2-5, to=2-6]
\end{tikzcd}
\]
\end{lem}

Assume that the cocharacter $\mu$ of $P_1$ is quotient-minuscule in both $P_1$ and $P_2$.

\begin{lem}\label{lem: generic fiber, same U, cartesian}
    We have a Cartesian diagram:
\[
\begin{tikzcd}
	{\Gra{P_1, \mu}} & {\Gra{P_2, \Spd E_1, \mu}} \\
	{\Gra{G_1, \mu}} & {\Gra{G_2, \Spd E_1, \mu}.}
	\arrow[from=1-1, to=1-2]
	\arrow[from=1-1, to=2-1]
	\arrow[from=1-2, to=2-2]
	\arrow[from=2-1, to=2-2]
\end{tikzcd}
\]
    In particular, $\Gra{P_1, \mu} \to \Gra{P_2, \Spd E_1, \mu}$ is proper.
\end{lem}
\begin{proof}

For $i = 1, 2$, let $\Gra{P_i, \mu, \Spd E_1}' \subset \Gra{P_i, \Spd E_1}$ be the base change of $\Gra{G_i, \Spd E_1, \mu} \subset \Gra{G_i, \Spd E_1}$. Consider the commutative diagram with vertical diagrams being Cartesian:
\[
\begin{tikzcd}[sep=tiny]
	{\Gra{P_1, \mu}'} &&& {\Gra{P_2, \Spd E_1, \mu}'} \\
	& {\Gra{P_1, \Spd E_1}} & {\Gra{P_2, \Spd E_1}} \\
	{\Gra{G_1, \mu}} &&& {\Gra{G_2, \Spd E_1, \mu}} \\
	& {\Gra{G_1, \Spd E_1}} & {\Gra{G_2, \Spd E_1}.}
	\arrow[from=1-1, to=1-4]
	\arrow[from=1-1, to=2-2]
	\arrow[from=1-1, to=3-1]
	\arrow[from=1-4, to=2-3]
	\arrow[from=1-4, to=3-4]
	\arrow[from=2-2, to=2-3]
	\arrow[from=2-2, to=4-2]
	\arrow[from=2-3, to=4-3]
	\arrow[from=3-1, to=3-4]
	\arrow[from=3-1, to=4-2]
	\arrow[from=3-4, to=4-3]
	\arrow[from=4-2, to=4-3]
\end{tikzcd}
\]
Both $\Gra{G_1, \mu}$ and $\Gra{G_2, \Spd E_1, \mu}$ are proper over $\Spd E_1$, thus $\Gra{G_1, \mu} \to \Gra{G_2, \Spd E_1, \mu}$ is proper. By Proposition \ref{prop: leq mu is a closed embedding}, $\Gra{P_i, \Spd E_1, \mu} \to \Gra{P_i, \Spd E_1, \mu}'$ are closed embeddings for $i = 1, 2$. In particular, $\Gra{P_1, \mu} \to \Gra{P_2, \Spd E_1, \mu}$ is proper. In order to show the proper morphism between diamonds
\begin{equation}\label{eq: generic fiber, same U, cartesian}
    \Gra{P_1, \mu} \to \Gra{P_2, \Spd E_1, \mu} \times_{\Gra{G_2, \Spd E_1, \mu}} \Gra{G_1, \mu}
\end{equation}
is an isomorphism, it suffices to check it is an isomorphism over any $\Spa(C, C^+)$-point in $\Gra{G_1, \mu}$. Let $\OO = B_{\dR}^+(C^{\sharp})$, $F = B_{\dR}(C^{\sharp})$, $g\xi^{\mu}G(\OO) \in \Gra{G_1, \mu}(C, C^+)$, with $g \in G_1(\OO)$. Since $g$ normalizes $U(\OO) \rtimes G_i(\OO)$, then the fibers of both sides of (\ref{eq: generic fiber, same U, cartesian}) over $g\xi^{\mu}G(\OO)$ are given by the same set (\ref{eq: fibers at grassmannian}).
\end{proof}

\begin{lem}\label{lem: integral, same U, cartesian}
    We have cartesian diagrams:
    \begin{equation}\label{eq: integral, same U, cartesian}
\begin{tikzcd}
	{\lcM_{\PP_1, \mu}} & {\lcM_{\PP_2, \mu, k_{E_1}}} & {\mathbb{M}_{\PP_1, \mu}^v} & {\mathbb{M}_{\PP_2, \mu}^v \times_{\Spd \OO_{E_2}} \Spd \OO_{E_1}} \\
	{\lcM_{\GG_1, \mu}} & {\lcM_{\GG_2, \mu, k_{E_1}},} & {\mathbb{M}_{\GG_1, \mu}^v} & {\mathbb{M}_{\GG_2, \mu}^v \times_{\Spd \OO_{E_2}} \Spd \OO_{E_1}.}
	\arrow[from=1-1, to=1-2]
	\arrow[from=1-1, to=2-1]
	\arrow[from=1-2, to=2-2]
	\arrow[from=1-3, to=1-4]
	\arrow[from=1-3, to=2-3]
	\arrow[from=1-4, to=2-4]
	\arrow[from=2-1, to=2-2]
	\arrow[from=2-3, to=2-4]
\end{tikzcd}
\end{equation}
\end{lem}
\begin{proof}
    Apply similar arguments appeared in the proof of Lemma \ref{lem: generic fiber, same U, cartesian}, $\mathbb{M}_{\PP_1, \mu}^v \to \mathbb{M}_{\PP_2, \mu}^v \times_{\Spd \OO_{E_2}} \Spd \OO_{E_1}$ is proper. In particular, to show the proper morphism 
    \begin{equation*}
    \mathbb{M}_{\PP_1, \mu}^v \to \mathbb{M}_{\PP_2, \mu, \Spd \OO_{E_1}}^v \times_{\mathbb{M}_{\GG_2, \mu, \Spd \OO_{E_1}}^v}\mathbb{M}_{\GG_1, \mu}^v
\end{equation*}
    is an isomorphism, it suffices to check over geometric points (see \cite[Cor. 17.4.8]{SW20}). Over generic fiber, this is Lemma \ref{lem: generic fiber, same U, cartesian}. Over special fiber, this is the proper morphism
    \begin{equation*}
    \lcM_{\PP_1, \mu} \to \lcM_{\PP_2, \mu, k_{E_1}} \times_{\lcM_{\GG_2, \mu, k_{E_1}}}\lcM_{\GG_1, \mu}
\end{equation*}
    induced by the left diagram in (\ref{eq: integral, same U, cartesian}). We apply the same arguments as in the proof of Lemma \ref{lem: generic fiber, same U, cartesian} and compute the fibers of both sides over $\lcM_{\GG_1, \mu}$, with the help of Corollary \ref{cor: special fiber of local model} (cf. Lemma \ref{lem: orbit action closed}) and Lemma \ref{lem: isomorphism of local model, quasi-parahoric}.
\end{proof}

\subsubsection{Devissage}

 Let $\wdt{P} \to P$ be a surjective morphism with kernel a central multiplicative group $Z \subset \wdt{P}$. Let $\wdt{G}$ be the Levi quotient of $\wdt{P}$. Note that $Z$ maps isomorphically into a central multiplicative group in $\wdt{G}$. Let $\wdt{\PP} = \wdt{\UU} \rtimes \wdt{\GG}$ (resp. $\PP = \UU \rtimes \GG$) be a quasi-parahoric group scheme of $\wdt{P}$ (resp. $P$). Assume $\wdt{P} \to P$ induces $\wdt{\PP} \to \PP$ and $\wdt{\UU} = \UU$.
\begin{lem}[{\cite[Prop. 21.5.2]{SW20}}]\label{lem: Devissage}
    Let $\wdt{\mu}$ be a quotient-minuscule cocharacter of $\wdt{P}$ and $\mu$ be its projection in $P$ defined over $E \subset \wdt{E}$, then the natural morphisms are isomorphisms:
    \[ \mathbb{M}_{\wdt{\GG}, \wdt{\mu}}^v \rightiso \mathbb{M}_{\GG, \mu}^v \otimes_{\Spd \OO_E} \OO_{\wdt{E}},\quad \mathbb{M}_{\wdt{\PP}, \wdt{\mu}}^v \rightiso \mathbb{M}_{\PP, \mu}^v \otimes_{\Spd \OO_E} \OO_{\wdt{E}}. \]
\end{lem}
\begin{proof}
    See \cite[Prop. 21.5.2]{SW20}, Lemma \ref{lem: integral, same U, cartesian} and \ref{lem: isomorphism of local model, quasi-parahoric}.
\end{proof}

\subsection{Shtukas with one leg bounded by $\mu$}\label{subsec-shtuka-nonred}

Let $\PP$ be a quasi-parahoric group scheme, and $\mu:\bb{G}_{m,\overline{\bb{Q}}_p}\to P_{\overline{\bb{Q}}_p}$ be a quotient-minuscule cocharacter.
\begin{definition}[{\cite[Def. 2.4.3]{PR24}}]
    Let $S$ be a perfectoid space over $k = \ovl{\mathbb{F}}_p$, and let $S^{\sharp}$ be an untilt of $S$ over $\OO_{\bE}$. A $\PP$-shtuka over $S$ with one leg at $S^{\sharp}$ is a pair $(\PPs, \phi_{\PPs})$, where
    \begin{enumerate}
        \item $\PPs$ is a $\PP$-torsor over the adic space $S\dottimes \Z_p$,
        \item $\phi_{\PPs}$ is a $\PP$-torsor isomorphism
        \begin{equation*}
            \phi_{\PPs}: \Frob_{S}^*(\PPs)|_{S\dottimes \Z_p \setminus S^{\sharp}} \rightiso \PPs|_{S\dottimes \Z_p \setminus S^{\sharp}},
        \end{equation*}
        which is meromorphic along the closed Cartier divisor $S^{\sharp} \subset S\dottimes \Z_p$.
    \end{enumerate}
    We say $(\PPs, \phi_{\PPs})$ is bounded by $\mu$ if the relative position $\phi_{\PPs}(\Frob_{S}^*(\PPs))$ and $\PPs$ at $S^{\sharp}$ is bounded by $\vM_{\PP, \mu}^v$ (\ref{def: local model}).
\end{definition}

\begin{definition}[{\cite[Def. 2.4.8]{PR24}}]
    Given a $v$-sheaf $\FF$, a $\PP$-shtuka over $\FF$ (resp. with one leg bounded by $\mu$) is a section of the $v$-stack (\cite[Prop. 19.5.3]{SW20}) given by the groupoid of $\PP$-shtukas with one leg (resp. with one leg bounded by $\mu$) over $\FF$. 
\end{definition}

\subsubsection{Functoriality}

Let $P_1 \to P_2$ be a morphism, $\PP_1$, $\PP_2$ be quasi-parahoric models of $P_1$, $P_2$ respectively, such that $P_1 \to P_2$ induces a morphism $\PP_1 \to \PP_2$. Let $\mu_1$ be a cocharacter of $P_1$ with reflex field $E_1$, and let $\mu_2$ be the composition of $\mu_1$ with $P_1 \to P_2$, with reflex field $E_2$. Assume $\mu_1$ and $\mu_2$ are quotient-minuscule in $P_1$ and $P_2$ respectively. Then the natural functor $\Gra{\PP_1} \to \Gra{\PP_2}$ induces
$$  \Gra{P_1, \mu_1} \to \Gra{P_2, \Spd E_1, \mu_2}, \quad \vM_{\PP_1,  \mu_1}^v \to \vM_{\PP_2,  \mu_2}^v \times_{\Spd \OO_{E_2}} \Spd \OO_{E_1}.$$
We also have
$$\Sht_{\PP_1} \to \Sht_{\PP_2},\quad \Sht_{\PP_1,  \mu_1} \to \Sht_{\PP_2,  \mu_2},$$
which come from pushing out the $\PP_1$-torsor to a $\PP_2$-torsor and naturally extending the isomorphism of $\PP_1$-torsors $\phi_{\PP_1}$ to an isomorphism of $\PP_2$-torsors $\phi_{\PP_2}$. 

The projection $\Sht_{\PP, \mu} \to \Sht_{\GG, \mu}$ has a section $\Sht_{\GG, \mu} \to \Sht_{\PP, \mu}$. The substack $\Sht_{\PP, \mu} \subset \Sht_{\PP}$ is contained in the preimage of $\Sht_{\GG, \mu} \subset \Sht_{\GG}$, but they are not equal.

\subsubsection{Trivial cocharacter}

\begin{definition}\label{def-tri-shtuka}
    We say $(\PPs, \phi_{\PPs})$ is a trivial $\PP$-shtuka over $\FF$ if for each perfectoid space $S \to \FF$, $\PPs_S$ is a trivial $\PP$-torsor over $S \dottimes \Z_p$, and $\phi_{\PPs_S}$ is the isomorphism induced by the Frobenius on the base $S$.
\end{definition}

In the following, we assume that $X$ is a normal scheme, flat and of finite type over $\OO_E$. By \cite[Cor. 2.7.10]{PR24}, a $\PP$-shtuka over $X^{\Dia/}$ is uniquely determined by its generic fiber over $X_{\eta}^{\Dia}$. In particular, a $\PP$-shtuka over $X^{\Dia/}$ is trivial if it is trivial over $X_{\eta}^{\Dia}$.

Let $1 \to P_1 \to P_2 \to P_3 \to 1$ be an exact sequence of linear algebraic groups which induces exact sequences
\[ 1 \to U_1 \to U_2 \to U_3 \to 1,\quad 1 \to G_1 \to G_2 \to G_3 \to 1. \]
Let $1 \to \PP_1 \to \PP_2 \to \PP_3 \to 1$ be an exact sequence of quasi-parahoric group schemes extending the generic fiber which induces
\[ 1 \to \UU_1 \to \UU_2 \to \UU_3 \to 1,\quad 1 \to \GG_1 \to \GG_2 \to \GG_3 \to 1. \]
In other words, in the decomposition $\PP_i = \UU_i \rtimes \GG_i$, we assume there exist compatible sections $\GG_i \to \PP_i$. Let $\mu$ be a quotient-minuscule cocharacter of $P_1$, then it is automatically quotient-minuscule in $P_2$.
\begin{prop}\label{prop: generic fiber reduce, then integral reduce}
   In the above setting, let $(\PPs_2, \phi_{\PPs_2})$ be a $\PP_2$-shtuka over $X^{\Dia/}$ (resp. with one leg bounded by $\mu$). Assume that its base change over $X_{\eta}^{\Dia}$ reduces to a $\PP_1$-shtuka $(\PPs_{1, \eta}, \phi_{\PPs_{1, \eta}})$ over $X_{\eta}^{\Dia}$ (resp. with one leg bounded by $\mu$). Then $(\PPs_2, \phi_{\PPs_2})$ reduces uniquely to a $\PP_1$-shtuka $(\PPs_1, \phi_{\PPs_1})$ (resp. with one leg bounded by $\mu$) over $X^{\Dia/}$ whose base change over $X_{\eta}^{\Dia}$ is $(\PPs_{1, \eta}, \phi_{\PPs_{1, \eta}})$.
\end{prop}
\begin{proof}
Consider the pushout $X^{\Dia/} \to \Sht_{\PP_2} \to \Sht_{\PP_3}$. The pushout $\PP_3$-shtuka $(\PPs_3, \phi_{\PPs_3})$ is trivial over $X_{\eta}^{\Dia}$, and thus is itself trivial over $X^{\Dia/}$. In other words, given any perfectoid space $S \to X^{\Dia/}$, the $\PP_2$-torsor $\PPs_2$ over $S\dottimes \Z_p$ reduces to a $\PP_1$-torsor $\PPs_1$, and the isomorphism
    \begin{equation*}
        \Phi_{\PPs_{2, S}}: \Frob_S^* \PPs_{2, S}|_{S\dottimes \Z_p \backslash S^{\sharp}} \cong  \PPs_{2, S}|_{S\dottimes \Z_p \backslash S^{\sharp}}
    \end{equation*}
    reduces to an isomorphism of $\PP_1$-torsor:
    \begin{equation*}
        \Phi_{\PPs_{1, S}}: \Frob_S^* \PPs_{1, S}|_{S\dottimes \Z_p \backslash S^{\sharp}} \cong  \PPs_{1, S}|_{S\dottimes \Z_p \backslash S^{\sharp}}.
    \end{equation*}
    Let us check the boundedness condition.
    \begin{enumerate}
        \item When $P_3$ is a reductive group, equivalently, we have that $U_1 = U_2$. By \cite[Prop. 21.5.1]{SW20}, $\mathbb{M}_{\GG_1, \mu}^v \rightiso \mathbb{M}_{\GG_2, \mu}^v \otimes \Spd \OO_{E_1}$ is a canonical isomorphism. Lemma \ref{lem: integral, same U, cartesian} gives $\mathbb{M}_{\PP_1, \mu}^v \rightiso \mathbb{M}_{\PP_2, \mu}^v \otimes \Spd \OO_{E_1}$. Hence, $(\PPs_1, \phi_{\PPs_1})$ is bounded by $\mu$.
        \item When $P_3$ is a unipotent group, equivalently, we have that $G_1 = G_2$, we claim that the following commutative diagram of $v$-sheaves is Cartesian. In particular, $(\PPs_1, \phi_{\PPs_1})$ is bounded by $\mu$.
\[
\begin{tikzcd}
	{\mathbb{M}_{\PP_1, \mu}^v} & {\mathbb{M}_{\PP_2, \mu}^v \times_{\Spd \OO_{E_2}} \Spd \OO_{E_1}} \\
	{\Gra{\PP_1, \Spd \OO_{E_1}}} & {\Gra{\PP_2, \Spd \OO_{E_1}}.}
	\arrow[from=1-1, to=1-2]
	\arrow[from=1-1, to=2-1]
	\arrow[from=1-2, to=2-2]
	\arrow[from=2-1, to=2-2]
\end{tikzcd}
\]
    The bottom morphism is a closed embedding since $\PP_2/\PP_1 = \UU$ is affine; then the top morphism is a closed embedding by construction. Then the induced morphism 
    \begin{equation}\label{eq: generic fiber reduce, then integral reduced, 1}
        \mathbb{M}_{\PP_1, \mu}^v \to \mathbb{M}_{\PP_2, \mu, \Spd \OO_{E_1}}^v \times_{\Gra{\PP_2, \Spd \OO_{E_1}}} \Gra{\PP_1, \Spd \OO_{E_1}}
    \end{equation}
    is also closed. To show that it is an isomorphism, it suffices to check the surjectivity on fibers over the $\Spa(C, C^+)$-points in $\Gra{\GG, \Spd \OO_{E_1}}$.
    
    Over the generic fiber, (\ref{eq: generic fiber reduce, then integral reduced, 1}) becomes
    \begin{equation}\label{eq: generic fiber reduce, then integral reduced, 2}
        \Gra{P_1, \mu} \to \Gra{P_2, \Spd E_1, \mu} \times_{\Gra{P_2, \Spd E_1}} \Gra{P_1, \Spd E_1}. 
    \end{equation}
    Let $\OO = B_{\dR}^+(C^{\sharp})$, $F = B_{\dR}(C^{\sharp})$, fix a point $g\xi^{\mu}G(\OO) \in \Gra{G, \mu}(C, C^+)$, where $g \in G(\OO)$, let $h = g\xi^{\mu}$, the fiber of the right hand side of (\ref{eq: generic fiber reduce, then integral reduced, 2}) over $h$ is
   \[U_2(\OO)/(U_2(\OO) \cap \xi^{\mu}U_2(\OO)\xi^{-\mu}) \times_{U_2(F)/\xi^{\mu}U_2(\OO)\xi^{-\mu}} U_1(F)/\xi^{\mu}U_1(\OO)\xi^{-\mu},\]
   as subset of $U_1(F)/\xi^{\mu}U_1(\OO)\xi^{-\mu}$, equals the fiber of the left hand side of (\ref{eq: generic fiber reduce, then integral reduced, 2}) over $h$:
   \[ U_1(\OO)/(U_1(\OO) \cap \xi^{\mu}U_1(\OO)\xi^{-\mu}), \]
   since $U_1(F) \subset U_2(F)$ (resp. $U_1(\OO_F) \subset U_2(\OO_F)$) is stabilized by the conjugation of $P_1(F)$ (resp. $P_1(\OO_F)$). Over the special fiber, the isomorphism can be checked similarly, with the help of Corollary \ref{cor: special fiber of local model} (cf. Lemma \ref{lem: orbit action closed}) and Lemma \ref{lem: isomorphism of local model, quasi-parahoric}.
    \item In general, we can insert an auxiliary $\PP'$ into
    \[ \PP_1 = \UU_1 \rtimes \GG_1 \to \PP' = \UU_2 \rtimes \GG_1 \to \PP_2 = \UU_2 \rtimes \GG_2, \]
    where $\GG_1$ acts on $\UU_2$ through $\GG_1 \to \GG_2$, and apply the above two steps successively.
    \end{enumerate} 
\end{proof}

\subsubsection{$\Bun_P$ and $P\textit{-}\Isoc$}\label{subsec: isocrystals}

Recall that $P\textit{-}\Isoc$ is a stack on $\PCAlg$ that maps $S = \Spec R$ to the groupoid of pairings $(\EE, \beta)$, where $\EE$ is a $P$-torsor over $\Spec W(R)[1/p]$, and $\beta: \sigma^*\EE \to \EE$ is an isomorphism of $P$-torsors. By \cite[Thm 7.14, Prop. 4.7, Prop. 6.3]{gleason2023meromorphic}, $P\textit{-}\Isoc$ is a $v$-stack, and $P\textit{-}\Isoc \cong \Bun_P^{\red}$, the proof of this part does not need $P$ to be reductive. By definition there is a bijection of the underlying topological spaces $|P\textit{-}\Isoc| \cong B(P)$.

Let $\mu$ be a quotient-minuscule cocharacter of $P$, also denote by $\mu$ its projection on $G$. Recall that one defines $G\textit{-}\Isoc_{\mu} \subset G\textit{-}\Isoc$ as the closed substack corresponding to $B(G, \mu) \subset B(G)$, following the main result in \cite{rapoport1996classification}. Since the natural projection $B(P) \to B(G)$ is a bijection (\cite[\S 3.6]{kottwitz1997isocrystals}), we let $P\textit{-}\Isoc_{\mu} \subset P\textit{-}\Isoc$ be the pullback of $G\textit{-}\Isoc_{\mu} \subset G\textit{-}\Isoc$, which is again a closed substack.

Let $\Sht_{\ca{P}}^W$ be the functor sending any perfect algebra $R$ to the groupoid of $(\ca{E},\beta)$, where $\ca{E}$ is a $\ca{P}$-torsor on $\Spec W(R)$ and $\beta:\sigma^*\ca{E}\dashrightarrow\ca{E}$ a modification. By \cite[Thm. 2.3.8]{PR24} and the Tannakian formalism, there is a natural isomorphism $\Sht^{\red}_{\PP} \cong \Sht^W_{\PP}$; cf. \cite[Lem. 3.1.5]{DvHKZ24ig}. When $\PP$ is parahoric, one can define $\Sht^W_{\PP, \mu}$ by requiring the relative position to be controlled by $\Gra{\PP, \mu}^W$. Then, under $\Sht^{\red}_{\PP} \cong \Sht^W_{\PP}$, we have $\Sht^{\red}_{\PP, \mu} \cong \Sht^W_{\PP, \mu}$ by \cite[Lem. 3.1.7]{DvHKZ24ig} (with the help of Corollary \ref{cor: special fiber of local model}). When $\PP$ is quasi-parahoric, we let $\Sht^W_{\PP, \mu} \subset \Sht^W_{\PP}$ be the reduction of $\Sht_{\PP, \mu} \subset \Sht_{\PP}$; this is compatible with the parahoric case.

Recall that the natural projection $\Sht^W_{\PP} \to P\textit{-}\Isoc$ defined by $(\EE, \beta) \mapsto (\EE|_{W(\ast)[1/p]}, \beta)$ is a $v$-cover. By construction, $\Sht^W_{\PP} \to P\textit{-}\Isoc \to G\textit{-}\Isoc$ factors through $\Sht^W_{\GG}$. By \cite[\S 3.2.2]{DvHKZ24ig}, the projection $\Sht^W_{\GG, \mu} \to G\textit{-}\Isoc$ factors through $G\textit{-}\Isoc_{\mu^{-1}}$ when $\GG$ is parahoric, then $\Sht^W_{\PP, \mu} \to P\textit{-}\Isoc$ factors through $P\textit{-}\Isoc_{\mu^{-1}}$ when $\PP$ is parahoric.

Recall that one has the Beauville--Laszlo morphism $\BL: \Sht_{\PP} \to \Bun_P$ as follows: let $(\PPs, \phi_{\PPs})$ be a $\PP$-shtuka on $S$ with one leg at $S^{\sharp}$. Restrict $(\PPs, \phi_{\PPs})$ to $\mathcal{Y}_{[r, \infty)}(S)$ for sufficiently large $r$ so that it excludes the leg at $S^{\sharp}$, and descend it to $X_{\FFC, S}$; we obtain a $\PP$-torsor $\EE(\PPs, \phi_{\PPs})$ on $X_{\FFC, S}$. By taking the reduction, we have a morphism $\BL^{\red}: \Sht_{\PP}^{\red} \to \Bun_P^{\red}$. This is compatible with the projection $\Sht^W_{\PP} \to P\textit{-}\Isoc$ together with the identification $P\textit{-}\Isoc \cong \Bun_P^{\red}$.

On the other hand, we have a functor $P\textit{-}\Isoc \to \Bun_P$, defined by composing the exact $\otimes$-functor $\Rep_{\Qp} P \to \Isoc$ with the exact $\otimes$-functor $\Isoc \to \Bun$ (see \cite[\S III.2]{fargues2021geometrization}, \cite{anschutz2019reductive}). By results in \cite{fargues2020g} and \cite{anschutz2019reductive}, when $P = G$ is a reductive group, such a functor induces a bijection on the underlying topological spaces (and moreover a homeomorphism when $B(G)$ is endowed with the topology induced by the Bruhat order reversed from \cite{rapoport1996classification}; see \cite{viehmann2024newton}). In general, we also have
\begin{lem}\label{lem: bijection B(P) to P-Isoc}
    The functor $P\textit{-}\Isoc \to \Bun_P$ induces a bijection $B(P) \to \Bun_P(C)/\cong$.
\end{lem}
\begin{proof}
    Since $B(P) \to B(G)$ and $B(G) \to \Bun_G(C)/\cong$ are bijections, it suffices to show that $H^1(X_{\FFC, C}, U) = 0$. In characteristic $0$, the unipotent group $U$ is split and is a successive extension of $\Ga$. Since $H^1(X_{\FFC, C}, \Ga) = H^1(X_{\FFC, C}, \OO_{X_{\FFC, C}}) = 0$ (see \cite[Prop. II.2.5(ii)]{fargues2021geometrization}), we have $H^1(X_{\FFC, C}, U) = 0$.
\end{proof}

\subsubsection{Shtukas and local systems}\label{subsec: shtukas and local systems}

Using the Tannakian formalism and Beauville--Laszlo gluing, and following the proofs of \cite[Prop. 12.4.6]{SW20} and \cite[Prop. 2.5.1]{PR24} (\cite[Thm. 22.5.2, Prop. 22.6.1]{SW20} works for any smooth affine model $\PP$ with connected fibers), we have the following:
\begin{prop}[{\cite[Prop. 2.5.1]{PR24}}]\label{prop: shtuka and loc system, general}
    Let $\PP$ be parahoric. The construction
\begin{equation}\label{eq: shtuka and local system}
    (\PPs, \phi_{\PPs}) \mapsto (\PPp, \DRT(\PPs))
\end{equation}
gives an equivalence of categories between
\begin{enumerate}
    \item $\PP$-shtukas over $S/\Spd E$ (resp.  with one leg bounded by $\mu$),
    \item pairs $(\PPp, D)$, where $\PPp$ is a pro-\'etale $\underline{\PP(\Z_p)}$-torsor over $S$ and $D: \PPp \to \Gra{P, \Spd E}$ (resp. $D: \PPp \to \Gra{P, \mu^{-1}}$) is a $\underline{\PP(\Z_p)}$-equivariant morphism over $\Spd E$.
\end{enumerate}
   This result can be generalized to the case where the base $S/\Spd E$ is replaced by a $v$-sheaf $\FF$.
\end{prop}

Finally, as in \cite[\S 2.6]{PR24}, using the Tannakian formalism, given a de-Rham pro-\'etale torsor $\PPp$ under $\underline{\PP(\Z_p)}$ over $X^{\Dia}/\Spd E$ (resp. we moreover assume that for all classical points $x\in X$, the filtration $\Fil^{\bullet}D_{\dR}(\ls_{\rho, \bar{x}})$ has constant conjugacy class), then $\PPp$ is equipped with a canonical Hodge--Tate period map $D: \PPp \to \Gra{P, \Spd E}$ (resp. $D: \PPp \to \Gra{P, \mu^{-1}}$), and one can attach to it a $\PP$-shtuka $(\PPs, \phi_{\PPs})$ over $X^{\Dia}/\Spd E$ (resp. with one leg bounded by $\mu$).

\subsubsection{Quasi-parahoric group schemes}

\begin{lem}\label{lem: fiber product in shtuka}
    Let $\PP = \UU \rtimes \GG$ be a quasi-parahoric group scheme. Then the following $2$-commutative diagrams are $2$-Cartesian:
\[
\begin{tikzcd}
	{\Sht_{\PP^{\circ}}} & {\Sht_{\PP}} & {\Sht_{\PP^{\circ}, \mu}} & {\Sht_{\PP, \mu}} \\
	{\Sht_{\GG^{\circ}}} & {\Sht_{\GG},} & {\Sht_{\GG^{\circ}, \mu}} & {\Sht_{\GG, \mu}.}
	\arrow[from=1-1, to=1-2]
	\arrow[from=1-1, to=2-1]
	\arrow[from=1-2, to=2-2]
	\arrow[from=1-3, to=1-4]
	\arrow[from=1-3, to=2-3]
	\arrow[from=1-4, to=2-4]
	\arrow[from=2-1, to=2-2]
	\arrow[from=2-3, to=2-4]
\end{tikzcd}
\]
\end{lem}
\begin{proof}
    Over $S \dot{\times}\Z_p$, it is easy to show that the groupoid of $\PP^{\circ}$-torsors is isomorphic to the groupoid of $\PP$-torsors whose reduction to $\GG$-torsors comes from $\GG^{\circ}$-torsors. To check the boundedness condition, we use Lemma \ref{lem: isomorphism of local model, quasi-parahoric}.
\end{proof}

Consider the Kottwitz map $\kappa_G: |\Bun_G| \to \pi_1(G)_{\Gamma}$, where $\Gamma = \Gal(\ovl{\rQ}_p|\Qp)$ and $\pi_1(G)$ is the algebraic fundamental group of $G$. Then $\kappa_G$ is locally constant and maps $\Bun_{G, \mu^{-1}}$ to $-\mu^{\natural}$. Let $\GG$ be a quasi-parahoric group scheme of $G$, and let $\Sht_{\GG, \mu}^{\kappa = -\mu^{\natural}} \subset \Sht_{\GG, \mu}$ be the open and closed substack that maps to $-\mu^{\natural}$ under the composition of the Beauville--Laszlo map $\BL^{\circ}: \Sht_{\GG, \mu} \to \Bun_G$ with $\kappa_G$. The map $\BL^{\circ}$ does not factor through $\Bun_{G, \mu^{-1}}$; nevertheless, the restriction $\BL^{\circ}: \Sht_{\GG, \mu}^{\kappa = -\mu^{\natural}} \to \Bun_{G}$ does factor through $\Bun_{G, \mu^{-1}}$; see \cite[Prop. 3.1.10]{daniels2024conjecture}. 

Recall that we have a short exact sequence
\[ 1 \to \GG^{\circ}(\Z_p) \to \GG(\Z_p) \to \pi_0(\GG)^{\phi} \to 1. \]
In \cite[\S 3.2]{daniels2024conjecture}, there is an action of $\pi_0(\GG)^{\phi}$ on $\Sht_{\GG^{\circ}, \mu}$, and the pushout $\EE \mapsto \GG \times^{\GG^{\circ}} \EE$ naturally induces an open and closed substack $[\Sht_{\GG^{\circ}, \mu}/\pi_0(\GG)^{\phi}] \to \Sht_{\GG, \mu}$; see \cite[Thm. 3.3.5]{daniels2024conjecture}. Let $\Sht_{\GG, \mu, \delta=1}$ be the image of $\Sht_{\GG^{\circ}, \mu}$. Then $\Sht_{\GG^{\circ}, \mu} \to \Sht_{\GG, \mu, \delta=1}$ is a torsor under the finite abelian group $\pi_0(\GG)^{\phi}$.

Let $\PP = \UU \rtimes \GG$ be a quasi-parahoric group scheme. We define 
\[ \Sht_{\PP, \mu, \delta = 1} \subset \Sht_{\PP, \mu}^{\kappa = -\mu^{\natural}} \subset \Sht_{\PP, \mu}\quad (\textit{resp}.\ \Sht_{\PP, \mu, \delta = 1}^W \subset \Sht_{\PP, \mu}^{W, \kappa = -\mu^{\natural}} \subset \Sht_{\PP, \mu}^W)\] 
as the open and closed substacks that are the pullbacks of
\[ \Sht_{\GG, \mu, \delta = 1} \subset \Sht_{\GG, \mu}^{\kappa = -\mu^{\natural}} \subset \Sht_{\GG, \mu}\quad (\textit{resp}.\ \Sht_{\GG, \mu, \delta = 1}^W \subset \Sht_{\GG, \mu}^{W, \kappa = -\mu^{\natural}} \subset \Sht_{\GG, \mu}^W)\]
under the natural projection $\Sht_{\PP, \mu} \to \Sht_{\GG, \mu}$ (resp. $\Sht_{\PP, \mu}^W \to \Sht_{\GG, \mu}^W$). Then $\Sht_{\PP, \mu}^{W, \kappa = -\mu^{\natural}} \to P\textit{-}\Isoc$ factors through $P\textit{-}\Isoc_{\mu^{-1}}$.

Also, we have a short exact sequence
\begin{equation}\label{eq: P and P circ exact sequence}
    1 \to \PP^{\circ}(\Z_p) \to \PP(\Z_p) \to \pi_0(\PP)^{\phi} \to 1.
\end{equation}
We define an action of $\pi_0(\PP)^{\phi}$ on $\Sht_{\PP^{\circ}, \mu}$ as in \cite[\S 3.2]{daniels2024conjecture}, and the pushout functor naturally induces $[\Sht_{\PP^{\circ}, \mu}/\pi_0(\PP)^{\phi}] \to \Sht_{\PP, \mu}$.

\begin{cor}[{cf. \cite[Cor. 3.3.8, 3.3.9, Rmk. 4.1.5]{daniels2024conjecture}}]\label{cor: quasi-parahoric}
    \begin{enumerate}
        \item $[\Sht_{\PP^{\circ}, \mu}/\pi_0(\PP)^{\phi}]$ is open and closed in $\Sht_{\PP, \mu}$ with image $\Sht_{\PP, \mu, \delta = 1}$. In particular, $\Sht_{\PP^{\circ}, \mu} \to \Sht_{\PP, \mu, \delta = 1}$ is a torsor under the finite abelian group $\pi_0(\PP)^{\phi}$.
        \item There is a natural isomorphism
    \[ \Sht_{\PP, \mu, \delta = 1} \times_{\Spd \OO_E} \Spd E \cong [\Gra{P, \mu^{-1}}/\underline{\PP(\Z_p)}]. \]
        \item Let $X$ be a normal scheme of finite type and flat over $\Z_p$. Assume that we have $X^{\Dia/} \to \Sht_{\PP, \mu}$, and $X_{\eta}^{\Dia} \to \Sht_{\PP, \mu, \Spd E}$ factors through $\Sht_{\PP, \mu, \delta = 1, \Spd E}$. Then $X^{\Dia/} \to \Sht_{\PP, \mu}$ factors through $\Sht_{\PP, \mu, \delta = 1}$.
    \end{enumerate}
\end{cor}
\begin{proof}
    \begin{enumerate}
        \item By definition, $\pi_0(\GG) = \pi_0(\PP)$. This corollary directly follows from Lemma \ref{lem: fiber product in shtuka}.
        \item For parahoric subgroups, this follows from \cite[Prop. 11.16]{zhang2023pel}; note that \cite[Thm. 22.5.2, Prop. 22.6.1]{SW20} works for any smooth affine model $\PP$ with connected fibers. For quasi-parahoric subgroups, apply the results in Part $(1)$.
        \item See \cite[Rmk. 4.1.5]{daniels2024conjecture}.
    \end{enumerate}
\end{proof}

\begin{lem}\label{lem: fiber product of shtukas}
    Let $P_1 = U \rtimes G_1$, $P_2 = U \rtimes G_2$, $P_1 \to P_2$ be a morphism compatible with $U \rtimes G_1 \stackrel{(\identity, f)}{\to} U \rtimes G_2$ in the sense of \ref{def: compatible, P and G}. Let $\PP_i = \UU \rtimes \GG_i$ be parahoric group schemes of $P_i$. Assume $P_1 \to P_2$ induces $\PP_1 \to \PP_2$ that is compatible with $\GG_1 \to \GG_2$ in the sense of \ref{def: compatible, PP and GG}. Assume the cocharacter $\mu$ of $P_1$ is quotient-minuscule in both $P_1$ and $P_2$. Then the following $2$-commutative diagrams are $2$-Cartesian:
\[
\begin{tikzcd}
	{\Sht_{\PP_1}} & {\Sht_{\PP_2}} & {\Sht_{\PP_1, \mu}} & {\Sht_{\PP_2, \mu}} \\
	{\Sht_{\GG_1}} & {\Sht_{\GG_2},} & {\Sht_{\GG_1, \mu}} & {\Sht_{\GG_2, \mu}.}
	\arrow[from=1-1, to=1-2]
	\arrow[from=1-1, to=2-1]
	\arrow[from=1-2, to=2-2]
	\arrow[from=1-3, to=1-4]
	\arrow[from=1-3, to=2-3]
	\arrow[from=1-4, to=2-4]
	\arrow[from=2-1, to=2-2]
	\arrow[from=2-3, to=2-4]
\end{tikzcd}
\]
\end{lem}
\begin{proof}
    The left diagram is Cartesian; see \cite[Cor. 2.17]{DY25}. Note that $\PP_1 = \GG_1 \times_{\GG_2} \PP_2$ by assumption. The right diagram is Cartesian by Lemma \ref{lem: integral, same U, cartesian}.    
\end{proof}

\begin{cor}\label{cor: quasi-parahoric, fiber product}
     Keep the notation and assumptions from Lemma \ref{lem: fiber product of shtukas}, with the word \textbf{parahoric} replaced by \textbf{quasi-parahoric}. Then the following $2$-commutative diagram is $2$-Cartesian:
    \begin{equation}\label{eq: quasi-parahoric, fiber product}
\begin{tikzcd}
	{\Sht_{\PP_1, \mu, \delta=1}} & {\Sht_{\PP_2, \mu, \delta=1}} \\
	{\Sht_{\GG_1, \mu, \delta=1}} & {\Sht_{\GG_2, \mu, \delta=1}.}
	\arrow[from=1-1, to=1-2]
	\arrow[from=1-1, to=2-1]
	\arrow[from=1-2, to=2-2]
	\arrow[from=2-1, to=2-2]
\end{tikzcd}
    \end{equation}
\end{cor}
\begin{proof}
    First, assume $G_1 = G_2$ and $\GG_1 \to \GG_2$ induces $\GG_1^{\circ} = \GG_2^{\circ}$. By \cite[Cor. 3.3.11]{daniels2024conjecture} and the first part of Corollary \ref{cor: quasi-parahoric}, the diagram (\ref{eq: quasi-parahoric, fiber product}) is $2$-Cartesian, where both horizontal morphisms are torsors under finite abelian groups
    \[\pi_0(\PP_2)^{\phi}/\pi_0(\PP_1)^{\phi} = \pi_0(\GG_2)^{\phi}/\pi_0(\GG_1)^{\phi}. \]
    In general, given a morphism of quasi-parahoric group schemes $\GG_1 \to \GG_2$, then $\GG_1^{\circ} \to \GG_1 \to \GG_2$ factors through $\GG_2^{\circ}$, this can be seen by the equation $\GG_i^{\circ}(\bZ_p) = \GG_i(\bZ_p) \cap \Ker \Tilde{\kappa}_{G_i}$, where $\Tilde{\kappa}_{G_i}: G_i(\bQ) \to \pi_1(G_i)_I$ is the functorial Kottwitz map. Consider the diagram
    \[
\begin{tikzcd}[sep=small]
	{\Sht_{\PP_1^{\circ}, \mu}} && {\Sht_{\PP_2^{\circ}, \mu}} & \\
	& {\Sht_{\PP_1, \mu, \delta=1}} && {\Sht_{\PP_2, \mu, \delta=1}} \\
	{\Sht_{\GG_1^{\circ}, \mu}} && {\Sht_{\GG_2^{\circ}, \mu}} \\
	& {\Sht_{\GG_1, \mu, \delta=1}} && {\Sht_{\GG_2, \mu, \delta=1}.}
	\arrow[from=1-1, to=1-3]
	\arrow[from=1-1, to=2-2]
	\arrow[from=1-1, to=3-1]
	\arrow[from=1-3, to=2-4]
	\arrow[from=1-3, to=3-3]
	\arrow[from=2-2, to=2-4]
	\arrow[from=2-2, to=4-2]
	\arrow[from=2-4, to=4-4]
	\arrow[from=3-1, to=3-3]
	\arrow[from=3-1, to=4-2]
	\arrow[from=3-3, to=4-4]
	\arrow[from=4-2, to=4-4]
\end{tikzcd}
    \]
    Since the left, right and the back diagrams are 2-Cartesian by above discussions and Lemma \ref{lem: fiber product of shtukas}, and since $\Sht_{\GG_1^{\circ}, \mu} \to \Sht_{\GG_1, \mu, \delta=1}$ is a surjection of $v$-stacks, the front diagram (\ref{eq: quasi-parahoric, fiber product}) is 2-Cartesian.
\end{proof}

\section{Log diamonds and log shtukas}\label{sec-log-diamond-shtuka}

We discuss some theorems in \cite[Sec. 2]{PR24} in the log setting. 
To start with, we also introduce some basic results about log diamonds defined in the $v$-topology for log adic spaces, schemes and formal schemes. Since we will deal with many slightly different cases simultaneously, the readers will be frequently reminded about the cases under consideration before seeing the results. 
\subsection{Log diamonds associated with log schemes}\label{subsec-log-diamond}
Our goal is to develop results on log diamonds associated with fs log schemes over $\bb{Z}_p$. This may be viewed as a common generalization of some definitions and propositions in \cite[\S18]{SW20} and \cite[\S7]{KY25}.\par
\subsubsection{}\label{subsubsec-diamond}
Let $X$ be a scheme locally of finite type over $\bb{Z}_p$. We write $X^\ad$ as the adic space that represents the fiber product
$$X\times_{\Spec\bb{Z}_p}\Spa \bb{Z}_p.$$
Let $\wat{X}$ be the formal scheme defined by the $p$-adic completion of $X$ and $\wat{X}^\ad$ the adic space associated with it. We can also define $X^\ad$ similarly when $X$ is of finite type over $\bb{Q}_p$. \par 
When $X$ is separated and of finite type, there is an open embedding $\wat{X}^\ad\hookrightarrow X^\ad.$
Now let $X$ be a separated scheme of finite type over $\ca{O}_E$, where $E$ is a finite field extension of $\bb{Q}_p$. We list some easy facts that will be used later:
\begin{enumerate}
    \item Let $X_{\eta}$ be the generic fiber of $X$, we can also define $(X_{\eta})^{\ad} = X_{\eta} \times_{\Spec E} \Spa(E, \OO_E)$ in the sense of \cite[Proposition 3.8]{Hub94}, then
\begin{equation}\label{eq-dia-asso-with-scheme-and-adic-space}
    (X_{\eta})^{\ad} = (X \times_{\Spec \OO_E} \Spa \OO_E) \times_{\Spa \OO_E} \Spa(E, \OO_E) = (X^{\ad})_{\eta}.
\end{equation}
    \item Let $X \to Y \leftarrow Z$ be schemes separated locally of finite type over $E$, then $X^{\ad}$, $Y^{\ad}$, $Z^{\ad}$ are Tate. In particular, $X^{\ad} \to Y^{\ad} \leftarrow Z^{\ad}$ are adic, the fiber product exists. Checking the functoriality in \cite[Proposition 3.8]{Hub94}, we have
    \begin{equation}\label{eq-anal-comm-with-fiber-product}
        (X \times_Y Z)^{\ad} = X^{\ad} \times_{Y^{\ad}} Z^{\ad}
    \end{equation}
    \item Let $X$ be a smooth variety over $E$, then $X^{\ad}$ is a sousperfectoid analytic adic space.
    \item Let $\wdh{X} \to \wdh{Y} \leftarrow \wdh{Z}$ be formal schemes separated locally of finite type over $\OO_E$, then (\ref{eq-anal-comm-with-fiber-product}) also holds:
    \begin{equation}\label{eq-anal-of-formal-scheme-commutes-with-fiber-product}
        (\wdh{X} \times_{\wdh{Y}} \wdh{Z})^{\ad} = \wdh{X}^{\ad} \times_{\wdh{Y}^{\ad}} \wdh{Z}^{\ad}.
    \end{equation} 
\item
The functor $X^{\Dia}$ (resp. $X^{\dia}$) can be constructed via taking the diamond functor for adic spaces for $X^\ad$ (resp. $\wat{X}^\ad$); both $(-)^\Dia$ and $(-)^\dia$ commute with fiber products in the category of separated schemes of finite type over $\OO_E$, and are compatible with (\ref{eq-anal-comm-with-fiber-product}) and (\ref{eq-anal-of-formal-scheme-commutes-with-fiber-product}) under (\ref{eq-dia-asso-with-scheme-and-adic-space}).
\end{enumerate}
\subsubsection{}Let us recall some terminology in log geometry. For details, we refer the readers to \cite{Kat89}, \cite{Ogu18}, \cite{DLLZ23} and \cite{KY25}. \par
Let $\md{P}$ be a monoid. One can associate a group $\md{P}^{\gp}$ with a natural homomorphism between monoids $\gp_{\md{P}}:\md{P}\to \md{P}^{\gp}$; in fact, $\md{P}\mapsto \md{P}^{\gp}$ is the left adjoint functor of the natural inclusion of the category of groups into the category of monoids.\par
We say that $\md{P}$ is \emph{integral} if $\gp_{\md{P}}$ is injective, say that $\md{P}$ is \emph{saturated}, if $\md{P}$ is integral and any $x\in \md{P}^{\gp}$ such that $nx\in \md{P}$ for some $n\in \bb{Z}_{>0}$ is in $\md{P}$, and say that $\md{P}$ is \emph{fine} if $\md{P}$ is integral and finitely generated. We say that $\md{P}$ is \emph{fs} if $\md{P}$ is both fine and saturated. Denote $\overline{\md{P}}:=\md{P}/\md{P}^\times$.\par
We say an adic space $X$ is \emph{{\'e}tale sheafy} if $X_{\et}$ is a site and if $\ca{O}_{X_\et}$ is a sheaf. This is satisfied, for example, when $X$ is analytic stably sheafy or is Noetherian affinoid (cf. \cite[Cor. A.11]{DLLZ23}).
\begin{definition}[{\cite[Def. 2.2.2]{DLLZ23}}]\label{def-log-adic-space}
Let $X$ be an {\'e}tale sheafy adic space. A \emph{prelog structure} on X is a pair $(\ca{M}_X,\alpha)$, where $\ca{M}_X$ is a sheaf of monoids over $X_{\et}$ and $\alpha:\ca{M}_X\to\ca{O}_{X_{\et}}$ is a morphism of sheaves of monoids, and such a pair is a \emph{log structure} if the induced morphism $\alpha^{-1}(\ca{O}^\times_{X_\et})\to\ca{O}^\times_{X_\et}$ is an isomorphism. \par
A \emph{log adic space} is a triple $(X,\ca{M}_X,\alpha)$ consisting of an {\'e}tale sheafy adic space $X$ and a log structure $(\ca{M}_X,\alpha)$ as above.\par
A \emph{morphism} $f:(X,\ca{M}_X,\alpha_X)\to (Y,\ca{M}_Y,\alpha_Y)$ of log adic spaces is a morphism $f: X\to Y$ of underlying adic spaces $X$ and $Y$ with a morphism of sheaves of monoids $f^\sharp: f^{-1}\ca{M}_Y\to \ca{M}_X$ such that the diagram
\begin{equation}
    \begin{tikzcd}
    f^{-1}\ca{M}_Y\arrow[rr,"f^\sharp"]\arrow[d]&& \ca{M}_X\arrow[d]\\
    f^{-1}\ca{O}_{Y_\et}\arrow[rr]&&\ca{O}_{X_\et}
    \end{tikzcd}
\end{equation}
commutes.
\end{definition}
We call such a log structure $(\ca{M}_X,\alpha)$ a log structure on $X_{\et}$ (or simply on $X$). We omit $\alpha$ if it is not important or is clear in the context. \par
\begin{definition}\label{def-perfectd-monoid}
\begin{enumerate}
\item Let $\ca{M}$ be a log structure on $X_{\et}$.
We say $\ca{M}_X$ is \emph{integral} (resp. \emph{saturated}) if it is a sheaf of integral (resp. saturated) monoids. Let $f:X\to Y$ be a morphism of adic spaces. The \emph{inverse image} $f^*\ca{M}_Y$ of a log structure $\ca{M}_Y$ of $Y$ is defined to be the log structure associated with the prelog structure on $X$ given by $f^{-1}\ca{M}_Y\to f^{-1}\ca{O}_{Y_{\et}}\to \ca{O}_{X_{\et}}$ (cf. \cite[(1.4)]{Kat89}).
\item For a sheaf of monoids $\ca{M}$ on $X_{\et}$, denote by $\ca{M}^\times$ its sheaf of invertible elements. Denote by $\ca{M}^\flat$ its tilt $\varprojlim_{x\mapsto x^p}\ca{M}$. 
Following \cite[Def. 2.19]{KY25}, we say a sheaf of monoids $\ca{M}$ on $X_{\et}$ is \emph{perfectoid} if $\overline{\ca{M}^\flat}:=\ca{M}^\flat/\ca{M}^{\flat,\times}\to \ca{M}/\ca{M}^\times=:\overline{\ca{M}}$ is an isomorphism. In particular, $\overline{\ca{M}}$ is \emph{uniquely $p$-divisible} if $\ca{M}$ is perfectoid (cf. \cite[Rmk. 2.23]{KY25}).
\end{enumerate}
\end{definition}
\begin{definition}[{\cite[Def. 2.3.1]{DLLZ23}}]\label{def-charts}
Let $\md{P}$ be a monoid. A \emph{chart of $(X,\ca{M}_X)$ modeled on $\md{P}$} is a morphism $\theta: \md{P}_X\to \ca{M}_X$ of {\'e}tale sheaves of monoids such that $\alpha\circ\theta$ factors through $\ca{O}^+_{X_\et}$ and the associated log structure of $\alpha\circ\theta$ is canonically isomorphic to that of $\alpha$. We say the chart is fs if $\md{P}$ is fs.
\end{definition}
A log adic space $X$ is \emph{fs} (resp. \emph{fine}) if $X$ {\'e}tale locally admits fs (resp. fine) charts. Following the conventions in \cite{Kat89}, for a prelog ring/space $(A,\md{P})$, one can associate a log ring/space denoted by $(A,\md{P})^a$ (see also \cite[Def. 2.2.2(6)]{DLLZ23}); the functor $(-)^a$ is the left adjoint of the natural inclusion of the category of log adic spaces into the category of prelog adic spaces.
\subsubsection{The cases we consider}\label{subsubsec-log-diamonds-setup}
From now on, the pair $(X,\ca{M}_X)$ can be the following cases:
\begin{enumerate}
\item \label{case-gen-log-adic-space} an fs log adic space over $\Spa \bb{Z}_p$;
    \item \label{case-gen-log-sch} an fs log scheme over $\Spec \bb{Z}_p$;
\item \label{case-gen-for-sch} an fs log formal scheme over $\Spf \bb{Z}_p$.    
\end{enumerate}
As a (pre-)adic space might not be {\'e}tale sheafy in general, sometimes it is worth separating the last two cases.\par
Let $X$ be a scheme over $\bb{Z}_p$ or $\bb{Q}_p$. 
In some situations, we assume that
\begin{description}
    \item[(SF)]\label{ass-gen} Let $(X,\ca{M}_X)$ be an fs log scheme over $\bb{Z}_p$ or $\bb{Q}_p$ with $X$ a separated and of finite type over $\bb{Z}_p$ or $\bb{Q}_p$.
\end{description}
We do not assume the properness of $X\to \Spec \bb{Z}_p$ here. \par
By \cite[Cor. A.11]{DLLZ23}, for log schemes satisfying (SF),  $X^\ad$ and $\wat{X}^\ad$ are both {\'e}tale sheafy.\par
There is a natural morphism of {\'e}tale sites 
$\nu^\ad:=\nu^\ad_X: X^\ad_{\et}\to X_{\et}$ (see \cite[3.2.8, p.180]{Hub96}). Also, there is a natural morphism of {\'e}tale sites $\wat{\nu}:=\wat{\nu}_X:\wat{X}^\ad_\et\to X_\et$ defined by composing with $\wat{X}_{\et}^\ad\to X^\ad_{\et}$.
\subsubsection{}\label{subsubsec-fine-perfd-log-structure}To define log diamonds, it is crucial to study a certain class of log structures on perfectoid spaces.\par
Let $Y$ be a perfectoid space. Denote by $Y_{\et}$ (resp. $Y_v$) the {\'e}tale site (resp. $v$-site) of $Y$. There is a natural projection of sites 
$$\nu:Y_{v}\to Y_{\et}.$$
By \cite[Thm. 8.7]{scholze2017etale}, $\ca{O}_{Y_v}$ is a $v$-sheaf. So it makes sense to define log structures on $Y_v$.
\begin{definition}\label{def-log-str-v}
A prelog structure on $Y_v$ is a pair $(\ca{M}_Y,\alpha)$ consisting of a $v$-sheaf of monoids $\ca{M}_Y$ in $Y_v^\sim$ and a morphism between $v$-sheaves $\alpha: \ca{M}_Y\to \ca{O}_{Y_v}$. It is a log structure on $Y_v$ if $\alpha^{-1}(\ca{O}^\times_{Y_v})\to \ca{O}^\times_{Y_v}$ is an isomorphism.\par
The definitions of morphisms and charts are also similar to those defined for log structures on $Y_{\et}$ (cf. Definition \ref{def-log-adic-space} and \ref{def-charts}).
\end{definition}
From now on, we let $Y_\tau$ stand for $Y_{\et}$ (resp. $Y_v$) if $\tau={\et}$ (resp. $\tau=v$).
We say a sheaf of monoids $\ca{M}$ on $Y_\tau$ is uniquely $p$-divisible (resp. perfectoid) if 
$\ca{M}\xrightarrow{x\mapsto x^p} \ca{M}$ (resp. 
$\overline{\ca{M}^\flat}\to\overline{\ca{M}}$) is an isomorphism. That is, for any $U\in Y_{\tau}$, $\ca{M}(U)$ is uniquely $p$-divisible (resp. $\overline{\ca{M}^\flat}(U)\to \overline{\ca{M}}(U)$ is an isomorphism).\par 
The following property is weaker than being uniquely $n$-divisible but stronger than being $n$-torsion-free.
\begin{definition}\label{def-p-unique}
Let $n$ be a positive integer and $\md{P}$ a monoid. We say that $\md{P}$ is \textbf{$n$-unique} if $\md{P}\xrightarrow{n\cdot} \md{P}$ is injective. Similarly, a sheaf of monoids $\ca{M}$ on a site $\ca{C}$ is $n$-unique if $\ca{M}(U)$ is $n$-unique for any $U\in \ob \ca{C}$. 
\end{definition}
\begin{definition}\label{def-perfectoid-log-structure}
Let $Y$ be a perfectoid space. Set $\tau=\et$ or $v$. Let $\bb{N}[\frac{1}{p}]\sbst \bb{Q}$ be the monoid consisting of elements of the form $\frac{a}{b}$ where $a\in \bb{N}$ and $b\in \{p^n\}$ for $n\geq 0$.
\begin{enumerate}
\item A uniquely $p$-divisible monoid $\md{P}$ is called \textbf{$p$-finitely generated} if $\md{P}$ is $n$-unique for all positive integers $n$, and there is a finite set of elements $S\sbst \md{P}$ such that the set of all $p^i$-th roots of elements of $S$ generates $\md{P}$ for integers $i>0$. In other words, there is a surjection of monoids $\bb{N}[\frac{1}{p}]^{\oplus n}\twoheadrightarrow \md{P}$.
\item A uniquely $p$-divisible monoid $\md{P}$ is called \textbf{$p$-weakly finitely generated} if $\md{P}$ is $n$-unique for all positive integers $n$, and there is a uniquely $p$-divisible, $p$-finitely generated submonoid $\md{P}'\sbst \md{P}$ such that $\md{P}'\sbst \md{P}\sbst \bb{Q}_{\geq 0}\md{P'}$. Note that the notation $\bb{Q}_{\geq 0}\md{P}':=\varinjlim_{a\mapsto na, a\in\md{P},n\geq 1}\md{P}$ and the inclusions make sense because of the $n$-uniqueness assumption.
\item A perfectoid log structure $(\ca{M}_Y,\alpha)$ on $Y_\tau$ is \textbf{$p$-coherent} (resp. \textbf{$p$-weakly coherent}) if, for any geometric point $\overline{x}$ of $Y_{\et}$, there is an {\'e}tale neighborhood $U_{\overline{x}}$ such that $\ca{M}_{Y}|_{U_{\overline{x}}}$ admits a $p$-finitely generated (resp. $p$-weakly finitely generated) chart $\md{P}_{U_{\overline{x}}}$ (of log structures on $(U_{\overline{x}})_\tau$).
\item A log structure $(\ca{M}_Y,\alpha)$ on $Y_{\tau}$ is called \textbf{fine perfectoid} (resp. strongly fine perfectoid) if $\ca{M}_Y$ is integral, perfectoid, and $p$-weakly coherent (resp. $p$-coherent).
\end{enumerate}
\end{definition}
\begin{rk}\label{rk-fine-perfectoid-pushforward}
A log structure $\ca{M}$ on $Y_v$ is fine perfectoid implies that $\nu_*\ca{M}$ on $Y_{\et}$ is fine perfectoid. In fact, to see that $\nu_*\ca{M}$ is perfectoid, note that $\overline{\nu_*\ca{M}}\iso \nu_*\overline{\ca{M}}$ and $\overline{\nu_*\ca{M}^\flat}\iso \nu_*\overline{\ca{M}^\flat}$; this follows from Lemma \ref{lem-nubarM-iso}. Then the isomorphism of sheaves holds after pushforward.
\end{rk}
\begin{rk}\label{rk-condition-log-necessary}
We remark that the integrality is crucial in the proof of Theorem \ref{thm-log-diamond-v-sheaves}, while perfectoidness and (weak) $p$-finite generation will be used in the proof of Theorem \ref{thm-equi-cat-gen}. 
\end{rk}
We discuss some basic properties of the log structures defined above.
\begin{lem}[{cf. \cite[Lem. 2.1.10]{DLLZ23}}]\label{lem-unique-p-divisible-section}
Let $\lambda:\md{M}\to\md{P}'$ be a surjective morphism between integral monoids such that $\ker \lambda^{\gp}\sbst \md{M}$. Suppose that $\md{P}'$ is additionally saturated, sharp, uniquely $p$-divisible, and $p$-finitely generated monoid. Then $\lambda$ admits a section. 
\end{lem}
\begin{proof}
By integrality, there are injections $\md{M}\hookrightarrow\md{M}^\gp$ and $\md{P}'\hookrightarrow\md{P}^{\prime,\gp}$. Since $\md{P}'$ is saturated and sharp, $\md{P}^{\prime,\gp}$ is torsion-free. As there is a surjection $\bb{N}[\frac{1}{p}]^N\to \md{P}'$ for some positive integer $N$, we know that there is a surjection $\bb{Z}[\frac{1}{p}]^N\to \md{P}^{\prime,\gp}$. So $\md{P}^{\prime,\gp}\iso \bb{Z}[\frac{1}{p}]^n$ for some $n\leq N$ as $\bb{Z}{[1/p]}$ is a PID. So $\lambda^\gp:\md{M}^\gp\to \md{P}^{\prime,\gp}$ admits a section. Since $\ker\lambda^{\gp}\sbst \md{M}$, we see this section restricted to $\md{P}'$ must be contained in $\md{M}$. In fact, for any $s\in \md{M}$, $\lambda^{\gp,-1}(\lambda(s))=\lambda^{-1}(\lambda(s))\sbst \md{M}$.
\end{proof}
\begin{lem}\label{lem-satfineperf-satfineperfbar}
Let $(\ca{M}_Y,\alpha)$ be a saturated and fine perfectoid log structure on $Y_{\tau}$. Then $\overline{\ca{M}}_Y$ is saturated, uniquely $p$-divisible, and $p$-weakly finitely generated.  
\end{lem}
\begin{proof}
We only show it when $\tau=\et$. The case when $\tau=v$ is similar. Saturatedness and unique $p$-divisibility follow from the saturatedness and the fine-perfectoidness of $\ca{M}_Y$, respectively. The question is {\'e}tale local, so we assume that $\ca{M}_Y$ admits a global (uniquely $p$-divisible and) $p$-weakly finitely generated chart $\md{P}$. From the definition above, there is a $p$-finitely generated submonoid $\md{P}'$ and inclusions $\md{P}'\sbst \md{P}\sbst \bb{Q}_{\geq 0}\md{P}'$. Fix a surjection $\bb{N}[\frac{1}{p}]^n\twoheadrightarrow\md{P}'$ with a commutative diagram
\begin{equation*}
    \begin{tikzcd}
        \bb{N}[\frac{1}{p}]^n\arrow[r,hook]\arrow[d,two heads]&\md{N}\arrow[r,hook]\arrow[d,two heads]& \bb{Q}_{\geq 0}^n\arrow[d,two heads]\\
        \md{P}'\arrow[r,hook]&\md{P}\arrow[r,hook]&\bb{Q}_{\geq 0}\md{P}',
    \end{tikzcd}
\end{equation*}
where $\md{N}$ is the pullback of $\md{P}$ via $\bb{Q}^n_{\geq 0}\to \bb{Q}_{\geq 0}\md{P}'$.\par 
Let $\overline{\ca{M}}_Y'$ be the image of $\md{P}'_Y$ via $\md{P}_Y\to \overline{\ca{M}}_Y$. The diagram above extends to a commutative diagram
\begin{equation}\label{diag-p-w-fg-quotient}
    \begin{tikzcd}
        \bb{N}[\frac{1}{p}]^n_Y\arrow[r,hook]\arrow[d,two heads]&\md{N}_Y\arrow[r,hook]\arrow[d,two heads]& \bb{Q}_{\geq 0,Y}^n\arrow[d,two heads]\\
        \md{P}'_Y\arrow[r,hook]\arrow[d,two heads]&\md{P}_Y\arrow[r,hook]\arrow[d,two heads]&(\bb{Q}_{\geq 0}\md{P}')_Y\\
        \overline{\ca{M}}_Y'\arrow[r,hook]&\overline{\ca{M}}_Y. &
    \end{tikzcd}
\end{equation}
Fix $U\in Y_{\et}$ and $x_1,x_2\in \overline{\ca{M}}_Y(U)$ such that $nx_1=nx_2$ for some positive integer $n$. Replacing $U$ with a cover, assume that $x_1$ (resp. $x_2$) lifts to $y_1$ (resp. $y_2$) in $\md{N}$. Taking $(*)^\gp$ to the left diagram, we have
\begin{equation*}
\begin{tikzcd}
    \bb{Z}[\frac{1}{p}]^n\arrow[r,hook]\arrow[d,two heads,"\beta_p"]&\md{N}^{\gp}\arrow[d,two heads,"\beta"]\\
    (\overline{\ca{M}}'_Y)^{\gp}\arrow[r,hook]&\overline{\ca{M}}_Y^\gp.
    \end{tikzcd}
\end{equation*}
Since $\overline{\ca{M}}_Y$ is saturated and sharp, both $\overline{\ca{M}}_Y$ and $\overline{\ca{M}}_Y^{\gp}$ are torsion-free. Hence, if $y\in \md{N}^{\gp}$ such that $ny\in \ker \beta_p$, then $y\in \ker \beta$. So $y_1-y_2\in \ker \beta$. This implies that $x_1=x_2$ in $\overline{\ca{M}}_Y^{\gp}(U)$ and, \emph{a priori}, in $\overline{\ca{M}}_Y(U)$ by integrality. Now we have proved the $n$-uniqueness of $\overline{\ca{M}}_Y$. Other properties follow from the diagram (\ref{diag-p-w-fg-quotient}).
\end{proof}
\begin{lem}\label{lem-fs-submonoid}
Let $\md{P}'$ be a saturated, sharp, uniquely $p$-divisible, and $p$-finitely generated monoid. Then there is an fs sharp submonoid $\md{P}_0\sbst\md{P}'$ such that $\md{P}'=\md{P}_0[\frac{1}{p}]$, where $\md{P}_0[\frac{1}{p}]:=\varinjlim_n\frac{\md{P}_0}{p^n}\sbst \md{P}^{\prime,\gp}$.
\end{lem}
\begin{proof}
Since $\md{P}'$ is sharp and saturated, we embed $\md{P}'$ into $\md{P}^{\prime,\gp}\iso \bb{Z}[\frac{1}{p}]^n$ (see the proof of Lemma \ref{lem-unique-p-divisible-section}). Then, by the definition of being $p$-finitely generated, there is a set $S_0$ of $\bb{N}[\frac{1}{p}]$-generators of $\md{P}'$. Let $\md{P}_0^{\gp}$ be the $\bb{Z}$-module generated by $S_0$ in $\md{P}^{\prime,\gp}$, and $\md{M}_0$ be the ($\bb{N}$-)monoid generated by $S_0$ in $\md{P}_0^{\gp}\cap\md{P}'$. 
Let $\md{P}_0$ be the saturation of $\md{M}_0$ (in $\md{P}_0^{\gp}$); this monoid is saturated, finitely generated by \cite[I. Cor. 2.1.11]{Ogu18}, and integral by construction. \par
Note that the saturation $\md{P}_0$ of $\md{M}_0\sbst \md{P}_0^{\gp}\cap\md{P}'$ is in $\md{P}'$, as $\md{P}'$ is saturated. 
Hence, $\md{P}_0$ is fine by the above paragraph, and is sharp since it is included in a sharp monoid $\md{P}'$.
\end{proof}
\begin{lem}\label{lem-satfineperf-implies-chart}
Let $(\ca{M}_Y,\alpha)$ be a saturated and strongly fine perfectoid log structure on $Y_{\et}$.\par 
Then, for any geometric point $\overline{x}$ of $Y_{\et}$, there is an {\'e}tale neighborhood $U_{\overline{x}}$ such that $\ca{M}_Y$ admits a chart modeled on $\overline{\ca{M}}_{Y,\overline{x}}$.
\end{lem}
\begin{proof}
The proof is similar to that of \cite[Prop. 2.3.12-2.3.13]{DLLZ23}, but it is much more involved here, as we do not assume that the log structure is fs. Since this is a local question, we assume that there is an $p$-finitely generated $\md{P}$ together with $\theta:\md{P}\to\ca{M}_{Y}$ factoring through $\ca{O}^+_{Y_{\et}}$, such that $\theta^a:(\md{P}_X)^a\to\ca{M}_Y$ is an isomorphism. Note that this implies that $\overline{\theta}_{\overline{x}}:\md{P}\to(\md{P}_X)^a/(\alpha\circ\theta)^{-1}(\ca{O}^\times_{Y_{\et},\overline{x}})\to \overline{\ca{M}}_{Y,\overline{x}}$ is surjective (see \cite[Rmk. 2.3.4]{DLLZ23}). 
Our goal is reduced to finding a $p$-finitely generated chart $\md{P}':=\overline{\ca{M}}_{Y,\overline{x}}$ in an {\'e}tale neighborhood of $\overline{x}$.\par 
Firstly, we show that there is a section at the geometric point $\overline{x}$.
There is a natural surjective map between integral monoids 
$$\ca{M}^\flat_{Y,\overline{x}}\to\overline{\ca{M}^\flat_{Y,\overline{x}}}=\overline{\ca{M}}_{Y,\overline{x}}.$$
The last equality follows from the definition of perfectoid log structures. 
By Lemma \ref{lem-unique-p-divisible-section} above, there is a section from $\overline{\ca{M}}_{Y,\overline{x}}$ to $\ca{M}^\flat_{Y,\overline{x}}$. 
Let $\md{P}':=\overline{\ca{M}}_{Y,\overline{x}}$ and fix a section $\md{s}^\flat:\md{P}'\to\ca{M}^\flat_{Y,\overline{x}}$. Composing $\md{s}^\flat$ with the natural projection $p_0: \ca{M}^\flat_{Y,\overline{x}}\to \ca{M}_{Y,\overline{x}}$, we obtain a section $$\md{s}:\md{P}'\iso(\overline{\ca{M}^\flat_Y})_{\overline{x}}\xrightarrow{\md{s}^\flat}\ca{M}^\flat_{Y,\overline{x}}\xrightarrow{p_0}\ca{M}_{Y,\overline{x}}.$$

Denote by $\md{P}_0$ the fs sharp submonoid in $\md{P}'$ as constructed in Lemma \ref{lem-fs-submonoid}, and denote by $S_0'$ the finite set of generators of $\md{P}_0$. By the proof of both Lemma \ref{lem-unique-p-divisible-section} and Lemma \ref{lem-fs-submonoid}, there is a commutative diagram
\begin{equation*}
    \begin{tikzcd}
    \md{P}_0\arrow[rr,hook]\arrow[d,hook]&&\md{P}'\arrow[d,hook]\\
    \md{P}_0^\gp\arrow[rr,hook]&&\md{P}^{\prime,\gp},
    \end{tikzcd}
\end{equation*}
where $\md{P}_0^\gp$ (resp. $\md{P}^{\prime,\gp}$) is a finite free $\bb{Z}$- (resp. $\bb{Z}[\frac{1}{p}]$-) module. 
Since $\ca{M}_Y$ is perfectoid, $\ca{M}_{Y,\overline{x}}$ is also perfectoid. Since $\overline{\ca{M}^\flat_Y}_{\overline{x}}\to\overline{\ca{M}}_{Y,\overline{x}}$ is an isomorphism, we can lift any $s\in\overline{\ca{M}}_{Y,\overline{x}}$ to $\wdtd{s}=(s_0,s_1,\dots)\in \ca{M}^\flat_{Y,\overline{x}}$ via $\md{s}$. The $p^n$-th root of $\wdtd{s}$ taken in $\ca{M}^\flat_{Y,\overline{x}}$ is given by shifting the limit representing $\wdtd{s}$ to the left by $n$.\par
Consider the map $\alpha\circ\md{s}: \md{P}'\to\ca{M}^\flat_{Y,\overline{x}}\to \ca{M}_{Y,\overline{x}}\xrightarrow{\alpha}\ca{O}_{Y,\overline{x}}$. We can make this composition factor through $\ca{O}^+_{Y,\overline{x}}$ by adjusting $\md{s}$, which will be explained now.
Since $Y$ is perfectoid (and Tate), we have that $\{f\in\ca{O}_{Y,\overline{x}}||f(\overline{x})|>1\}\sbst \ca{O}^\times_{Y,\overline{x}}$ and $\{f\in\ca{O}_{Y,\overline{x}}||f(\overline{x})|\leq 1\}=\ca{O}^+_{Y,\overline{x}}$. If there is $c\in S_0'$ such that $|\alpha\circ\md{s}(c)(\overline{x})|>1$, then $f_c:=\alpha\circ\md{s}(c)\in \ca{O}^\times_{Y,\overline{x}}$. 
Assume that $|f_c(\overline{x})|$ is maximal among all $c_i\in S_0'$ (noting that $|S_0'|$ is finite). Write $S_0'=\{c_i\}_{i\in I}$. 
By \cite[I. Cor. 2.2.7]{Ogu18} (see also \cite[Lem. 2.3.12]{DLLZ23}), $\md{P}_0$ is embedded in some $\bb{N}^m$, where $c_i$ can be uniquely represented by a tuple $\ul{c}_i=(n^i_1,\cdots,n^i_m)\neq (0,0,\cdots,0)$ of natural numbers. Denote by $||\ul{c}_i||$ the sum of the coefficients of $\ul{c}_i$ representing $c_i$; it is not less than $1$. 
Define a map $$\mathsf{s}^\flat_{/f_c}: \sum_{i\in I} a_i c_i\mapsto \prod_{i\in I}(\md{s}^\flat(c_i)\cdot p_0^{-1}[(\md{s}(c))^{-||\ul{c}_i||}])^{a_i},$$
where $a_i\in \bb{N}[\frac{1}{p}]$. The function $(-)^{a_i}$ makes sense, as explained in the last paragraph. Note that $\mathsf{s}^\flat_{/f_c}$ is well-defined: On $\md{P}_0$, it is well-defined because $\ul{c}_i$ is uniquely defined; on $\md{P}'$, note that any element lies in $\md{P}_0$ by multiplying some $p^k$ and that $\ca{M}_{Y,\overline{x}}^\flat$ is uniquely $p$-divisible by the explanation in the last paragraph.\par 
Then $\mathsf{s}^\flat_{/f_c}:\md{P}'\to \ca{M}^\flat_{Y,\overline{x}}$ is a $\bb{N}[\frac{1}{p}]$-equivariant morphism. Set $\md{s}_{/f_c}:=p_0\circ\md{s}^\flat_{/f_c}$. Then $\alpha\circ\md{s}_{/f_c}$ factors through $\ca{O}^+_{Y,\overline{x}}$, as desired. \par
Set $\theta'_{\overline{x}}:=\md{s}_{/f_c}$ and $\theta^{\prime,a}_{\overline{x}}:(\md{P}')^a\to \ca{M}_{Y,\overline{x}}$. 
There is an {\'e}tale neighborhood $U$ of $\overline{x}$ such that $\theta'_{\overline{x}}$ extends to $\theta':\md{P}'_U\to \ca{M}_Y|_U$, since $\md{s}^\flat:\md{P}'\to \ca{M}^\flat_{Y,\overline{x}}$ extends to an {\'e}tale neighborhood. Consequently, there is an extension $\theta^{\prime,a}:=(\theta')^a$ for $\theta^{\prime,a}_{\overline{x}}$ on $U$. We only need to show that $\theta^{\prime,a}$ is an isomorphism up to shrinking the {\'e}tale neighborhood $U$.\par
Indeed, $\theta^{\prime,a}$ is surjective over some $U_1\in Y_{\et}\to U$, since $\ca{M}_{Y}|_U/\theta^{\prime,a}(\md{P}^{\prime,a}_U)=\theta^a(\md{P}^a_U)/\theta^{\prime,a}(\md{P}^{\prime,a}_U)$ is trivial at $\overline{x}$, and both $\md{P}$ and $\md{P}'$ are uniquely $p$-divisible and $p$-finitely generated. \par
Moreover, this map is injective. The congruence relation $\md{R}_0$ in $\md{P}_0\times \md{P}_0$ induced by $\md{P}_0\to \overline{\ca{M}}_{Y}$ is finitely generated by \cite[I. Lem. 2.1.9]{Ogu18} and is trivial at $\overline{x}$. Hence, it is trivial over some {\'e}tale neighborhood $U_2$ of $\overline{x}$, which we assume there is an {\'e}tale morphism $U_2\to U_1$. As both $\md{P}'$ and $\overline{\ca{M}}_Y$ are integral and uniquely $p$-divisible, the congruence relation $\md{R}\sbst \md{P}'\times\md{P}'$ induced by $\md{P}'\to \overline{\ca{M}}_Y$ is trivial over $U_{2,\et}$. Indeed, for $a,b\in \md{P}'$, if $a\sim b$ in $\overline{\ca{M}}_Y$ over $U_{2,\et}$, so is $p^i a$ and $p^i b$. Conversely, we can divide relations in $\md{R}_0$ by $p$: if $a\sim b$, let $a'=\frac{a}{p},b'=\frac{b}{p}$ in $\md{P}'$, then the images $c_1,c_2$ of $a'$ and $b'$ in $\overline{\ca{M}}_{Y}$ are equal since $pc_1=pc_2$ by the unique $p$-divisibility of $\overline{\ca{M}}_Y$. So $\md{R}$ is trivial if and only if $\md{R}_0$ is trivial.
This completes the proof.
\end{proof}
\subsubsection{}\label{subsubsec-log-diamonds}
We define log diamonds for the cases in \S\ref{subsubsec-log-diamonds-setup}.\par
In Case \ref{case-gen-log-adic-space},  
\begin{definition}\label{def-log-diamonds-adic-space}
Let $(X,\ca{M}_X,\alpha)$ be an ({\'e}tale sheafy) fs log adic space over $\Spa\bb{Z}_p$. Then the log diamond $(X,\ca{M}_X)^\Diamond$ (or $X^{\log \Diamond}$) is a functor sending $S=\Spa(A,A^+)\in \ob\Perf$ to the isomorphism classes $\{(S^\sharp,\iota,\ca{M}_{S^\sharp},f)\}/\simeq$, where $(S^\sharp,\iota)$ denotes an untilt as in \S\ref{subsubsec-diamond}, $\ca{M}_{S^\sharp}$ is a saturated and fine perfectoid log structure, and $f:S^\sharp\to X$ is a morphism between adic spaces that induces a morphism between log structures $f^*\ca{M}_X\to \ca{M}_{S^\sharp}$.  
\end{definition}
In Case \ref{case-gen-log-sch}, we can also define small and big log diamonds.
\begin{definition}\label{def-small-big-dia-gen-scheme}
Let $(X,\ca{M}_X)$ be an fs log scheme over $\bb{Z}_p$. We define a \textbf{log big diamond} functor $(X,\ca{M}_X)^{\Diamond}$ (or $X^{\log \Diamond}$) as
$$S=\Spa(A,A^+)\in \Perf\mapsto \{(S^\sharp=\Spa(A^\sharp,A^{\sharp,+}),\iota,\ca{M}_{S^\sharp},f)\}/\simeq,$$
where all but the last term in the tuple are as above, and $f: \Spa(A^\sharp,A^{\sharp,+})\to \Spec A^\sharp\to X$ induces a morphism $f^*\ca{M}_X\to \ca{M}_{S^\sharp}$; note that $f$ induces a morphism between {\'e}tale sites $S^\sharp_{\et}\to (\Spec A^\sharp)_{\et}\to X_\et$, and that all perfectoid spaces are {\'e}tale sheafy by \cite[Cor. A.11]{DLLZ23}.\par
We define a \textbf{log small diamond} functor $(X,\ca{M}_X)^{\dia}$ (or $X^{\log \dia}$) in a similarly-written form $S\in\Perf\mapsto \{(S^\sharp=\Spa(A^\sharp,A^{\sharp,+}),\iota,\ca{M}_{S^\sharp},f)\}/\simeq$, but here $f$ is a morphism $f:\Spa(A^\sharp,A^{\sharp,+})\to \Spec A^{\sharp,+}\to X$.
\end{definition}
In Case \ref{case-gen-for-sch}, we can associate a log diamond functor; we denote it by $(X,\ca{M}_X)^\dia$ or $X^{\log \dia}$.
\begin{definition}\label{def-small-dia-gen-for-sch}
Let $(X,\ca{M}_X)$ be an fs log formal scheme over $\Spf \bb{Z}_p$. The functor $(X,\ca{M}_X)^\Dia$ sends $S=\Spa(A,A^+)\in \Perf$ to 
$$\{(S^\sharp=\Spa(A^\sharp,A^{\sharp,+}),\iota,\ca{M}_{S^\sharp},f)\}/\simeq,$$
where $f: S^\sharp\to X$ is a morphism that induces $f^*\ca{M}_{X}\to \ca{M}_{S^\sharp}$ and other components are the same as above.
\end{definition}
This symbol is compatible with the log small diamond for Case \ref{case-gen-log-sch}: Let $\Spf(A,I)$ be an affine formally of finite type formal scheme. Suppose that $A$ is also $p$-adic complete. Let $X=\Spec A$; and we equip $X$ with an fs log structure $\ca{M}$. Then $X^{\log \diamond}$ is by definition the same as $(\Spf (A,(p)), \ca{M})^{\Dia}$; here the latter $\ca{M}$ is the log structure pulled back from $X_\et$.\par
In the case of (bounded) fs log $p$-adic formal schemes, Definition \ref{def-small-dia-gen-for-sch} essentially comes from the category $(X,\ca{M}_X)_{\Prism}^{\perf}$ in \cite[Def. 7.34 and Thm. 7.35]{KY25}.\par
\begin{rk}[{cf. Remark \ref{rk-condition-log-necessary}}]\label{rk-condition-log-notnecessary}
As we can see later from Lemma \ref{lem-can-log-str}, one might remove the $p$-finite generation condition in the definition of log diamonds given above. But it is this more restrictive class of log structures that is useful to our purpose.   
\end{rk}
Assume that $X$ is a $\bb{Z}_p$-scheme satisfying (SF). The diamond functors in Case \ref{case-gen-log-sch} also admit constructions via the diamond functor in the log adic space case.
\begin{lem}\label{lem-reconst-dia}
Under the assumption above, the log big diamond functor can be constructed as a functor
sending any affinoid perfectoid $S=\Spa(A,A^+)\in \ob\Perf$ to the isomorphism classes
$$\{(S^\sharp,\iota,\ca{M}_{S^\sharp},f)\}/\simeq,$$ 
where $f:(S^\sharp,\ca{M}_{S^\sharp})\to(X^\ad,\nu^{\ad,-1}\ca{M}_X)$ a morphism. More precisely, $f$ is a morphism between adic spaces $f:S^{\sharp}\to X^\ad$ that induces a morphism of sites 
$$S^{\sharp}_{\et}\xrightarrow{f} X^\ad_{\et}\xrightarrow{\nu^\ad} X_{\et}$$
together with a morphism $f^{-1}\circ \nu^{\ad,-1}\ca{M}_X\to \ca{M}_{S^{\sharp}}$.\par
Similarly, the log small diamond functor can be constructed as a functor the same as above, but changing $f$ to $f:(S^{\sharp},\ca{M}_{S^{\sharp}})\to (\wat{X}^\ad,\wat{\nu}^{-1}\ca{M}_X)$, which, more precisely, is a morphism of adic spaces $f:S^{\sharp}\to \wat{X}^\ad$ that induces a morphism of sites: 
$$S^{\sharp}_{\et}\xrightarrow{f} \wat{X}^\ad_{\et}\xrightarrow{\wat{\nu}} X_{\et}$$
together with a morphism $f^{-1}\circ \wat{\nu}^{-1}\ca{M}_X\to \ca{M}_{S^{\sharp}}$.
\end{lem}
\begin{proof}The fact that one can construct diamonds via $X^\ad$ and $\wat{X}^\ad$ is in \cite[\S 2.2]{anschutz2022p}. By \cite[Cor. A.11]{DLLZ23}, $X^\ad$ and $\wat{X}^\ad$ are {\'e}tale sheafy.
\end{proof}
The main theorem of \S\ref{subsec-log-diamond} is the following:
\begin{thm}\label{thm-log-diamond-v-sheaves} The functors $X^{\log \Diamond}$ and $X^{\log \dia}$ in Definition \ref{def-log-diamonds-adic-space}, Definition \ref{def-small-big-dia-gen-scheme} and Definition \ref{def-small-dia-gen-for-sch} are $v$-sheaves over $\Perf$.
\end{thm}
We will complete the proof of this theorem in \S\ref{subsubsec-proof-v-sheaf} and \S\ref{subsubsec-end-proof-v-sheaves}.
\subsubsection{}\label{subsubsec-proof-v-sheaf}
Let $Y$ be a perfectoid space. Let $\nu:Y_v\to Y_{\et}$ be the natural projection from the $v$-site to the {\'e}tale site of $Y$. There is a pullback functor $\nu^{-1}: Y_{\et}^\sim\to Y_v^\sim$ and a pushforward functor $\nu_*:Y_v^\sim\to Y_{\et}^\sim$ between topoi.\par 
Let $\ca{LOG}_?$ (resp. $\ca{LOG}_?^\mrm{int}$) be the category (resp. the full subcategory) of log structures (resp. integral log structures) on $Y_?$ for $?=\et$ and $v$. Let $\ca{M}$ be a log structure on $Y_\et$. Define the pullback of $\ca{M}$ to $Y_v$ as $\nu^*\ca{M}:=(\nu^{-1}\ca{M})^a$.\par
The following statement can be viewed as an analogue of \cite[App. A]{Ols03} and \cite[III. Prop. 1.4.1]{Ogu18} for perfectoid spaces.
\begin{lem}\label{lem-v-log-str}
With the conventions above, the pullback $\nu^*:\ca{LOG}_{Y_{\et}}\to \ca{LOG}_{Y_v}$ restricts to a functor from $\ca{LOG}^\mrm{int}_{Y_\et}$ to $\ca{LOG}_{Y_v}^{\mrm{int}}$, and this restriction induces an equivalence from the category of fine perfectoid log structures on $Y_{\et}$ to the category of fine perfectoid log structures on $Y_v$.\par 
Moreover, it induces an equivalence between the category of saturated and fine perfectoid log structures on $Y_{\et}$ and the category of saturated and fine perfectoid log structures on $Y_v$.
\end{lem}
We will show Lemma \ref{lem-v-log-str} in the rest of \S\ref{subsubsec-proof-v-sheaf}. The procedures are similar to the two references above.\par 
We first show the full faithfulness of $\nu^*|_{\ca{LOG}^{\mrm{int}}_{Y_{\et}}}$ and show that the essential image of this functor lies in $\ca{LOG}_{Y_v}^{\mrm{int}}$. 
\begin{lem}\label{lem-v-log-str-fullfaithful}
The functor $\nu^*:\ca{LOG}_{Y_{\et}}^{\mrm{int}}\to \ca{LOG}_{Y_v}$ is fully faithful. That is, we have
$$\Hom_{\ca{LOG}_{Y_{\et}}}(\ca{L}_1,\ca{L}_2)\iso \Hom_{\ca{LOG}_{Y_v}}(\nu^*\ca{L}_1,\nu^*\ca{L}_2)$$
for integral log structures $\ca{L}_1,\ca{L}_2\in\ca{LOG}_{Y_{\et}}^{\mrm{int}}$.
\end{lem}
\begin{proof}
 Denote by $\ca{LOG}^\pre_{Y_?}$ the category of prelog structures on the $?$ topology of $Y$. We have $\Hom_{\ca{LOG}_{Y_v}}(\nu^*\ca{L}_1,\nu^*\ca{L}_2)\iso \Hom_{\ca{LOG}_{Y_v}^\pre}(\nu^{-1}\ca{L}_1,\nu^*\ca{L}_2)\iso\Hom_{\ca{LOG}_{Y_{\et}}^\pre}(\ca{L}_1,\nu_*\nu^*\ca{L}_2)$ by adjunctions. The last one is isomorphic to $\Hom_{\ca{LOG}_{Y_{\et}}}(\ca{L}_1,\nu_*\nu^*\ca{L}_2)$. It suffices to show that the last term is canonically isomorphic to $\ca{L}_2$, which we will do in Lemma \ref{lem-v-log-fullfaithful-adjunction}.
\end{proof}
\begin{lem}\label{lem-v-log-fullfaithful-adjunction}
The canonical morphism between log structures $\ca{L}\to \nu_*\nu^*\ca{L}$ is an isomorphism for any $\ca{L}\in \ca{LOG}^{\mrm{int}}_{Y_{\et}}$.
\end{lem}
\begin{proof}
There is a functor defined on $Y_v$ by $$(f:Y'\to Y)\in Y_v\mapsto f^*\ca{L}:=(f^{-1}\ca{L})^a\in\ca{LOG}_{Y_\et'}.$$ 
We claim that the presheaf $$(f:Y'\to Y)\mapsto f^*\ca{L}(Y')$$
is a sheaf on $Y_v$. 
If we know this, then $\nu^*\ca{L}$ is this (pre)sheaf, and the desired result follows from the definition of $\nu_*$. \par
We now show the claim in the last paragraph. Let $\ca{L}':=f^*\ca{L}$ and $\ca{L}''=p_1^*f^*\ca{L}=p_2^*f^*\ca{L}$.
We show that the composition
\begin{equation*}
\begin{tikzcd} 
\ca{L}(Y)\arrow[r]&\mrm{eq}(\ca{L}'(Y')\arrow[r,shift left=2,"p_1^*" description]\arrow[r,shift right=2,"p_2^*" description]&\ca{L}''(Y'\times_Y Y')).
\end{tikzcd}
\end{equation*}
is an isomorphism. \par 
First of all, $f^*\ca{L}$ is integral: By \cite[Lem. 2.2.4]{DLLZ23} and \cite[I. Prop. 1.3.4]{Ogu18}, we check pointwise that $f^{-1}\ca{L}$ is integral and $f^*\ca{L}$ is integral.\par
Consider the commutative diagram
\begin{equation}\label{diag-log-str-v-sheaf}
\begin{tikzcd} 
\ca{O}^\times(Y)\arrow[r]\arrow[d]&\mrm{eq}(\ca{O}^\times(Y')\arrow[r,shift left=2,"p_1^*" description]\arrow[r,shift right=2,"p_2^*" description]\arrow[d]&\ca{O}^\times(Y'\times_Y Y'))\arrow[d]\\
\ca{L}(Y)\arrow[r]\arrow[d]&\mrm{eq}(\ca{L}'(Y')\arrow[r,shift left=2,"p_1^*" description]\arrow[r,shift right=2,"p_2^*" description]\arrow[d]&\ca{L}''(Y'\times_Y Y'))\arrow[d]\\
\overline{\ca{L}}(Y)\arrow[r]&\mrm{eq}(\overline{\ca{L}}'(Y')\arrow[r,shift left=2,"p_1^*" description]\arrow[r,shift right=2,"p_2^*" description]&\overline{\ca{L}}''(Y'\times_Y Y')).
\end{tikzcd}
\end{equation}
It follows from \cite[Thm. 8.7]{scholze2017etale} that $\ca{O}_Y$ is a $v$-sheaf, so is $\ca{O}^\times_Y$. We have the isomorphism of the first row of (\ref{diag-log-str-v-sheaf}).
The last row is also an isomorphism. This follows from \cite[Prop. 14.7 and 14.8]{scholze2017etale}.\par
We now examine the middle row of (\ref{diag-log-str-v-sheaf}). 
Then it follows from a diagram chasing that $\ca{L}(V)\to \ca{L}'(V\times_Y Y')$ is injective for any $V\in Y_{\et}$: Indeed, suppose that there are sections $s_1,s_2\in\ca{L}(V)$ mapping to the same element $s\in \ca{L}'(V\times_YY')$. We have that $s_1=s_2+u$ for $u\in\ca{O}^\times(V)$. Then $s=s+u$ in $\ca{L}'(V\times_YY')$, and therefore $u=0$ in $\ca{L}'(V\times_YY')$ by integrality. Then $u=0$ in $\ca{O}^\times(V)\to\ca{L}(V)$ by the paragraph above, as desired. \par
On the other hand, for any $s\in \mrm{eq}(\ca{L}'(Y')\rightrightarrows\ca{L}''(Y'\times_Y Y'))$, there is {\'e}tale locally $s'\in\ca{L}(U)$ for some $U\in Y_{\et}$ mapping to $s\in \ca{L}'(Y'\times_Y U)$ by diagram chasing. 
The element $s'$ on the {\'e}tale cover $U$ of $Y$ glues to a section on $Y$ since we can check it over the cover $Y'\times_Y U$ by the injectivity we just proved. We now have the isomorphism of the middle arrow, as desired.
\end{proof}
In the proof, we have shown that
\begin{cor}
    The pullback functor $\nu^*$ restricts to a fully faithful functor 
    $$\nu^*:\ca{LOG}_{Y_{\et}}^{\mrm{int}}\lra \ca{LOG}_{Y_{v}}^{\mrm{int}}.$$
\end{cor}
\begin{cor}\label{cor-fully-faithful-fine-perf}
The functor $\nu^*:\ca{LOG}_{Y_{\et}}^{\mrm{int}}\lra \ca{LOG}_{Y_{v}}^{\mrm{int}}$ sends fine perfectoid objects to fine perfectoid objects.
\end{cor}
\begin{proof}
    All conditions we need follow from the construction in Lemma \ref{lem-v-log-fullfaithful-adjunction}.
\end{proof}
We then show the essential surjectivity. Fix $\ca{M}\in\ca{LOG}_{Y_v}$ that is fine perfectoid. We show that  
\begin{lem}\label{lem-adjuction-iso}
    The adjunction morphism $\nu^*\nu_*\ca{M}\to \ca{M}$ is an isomorphism.
\end{lem}
\begin{proof}
There is an exact sequence $1\to \ca{O}^\times_{Y_{\et}}\to\nu_*\ca{M}\to \overline{\nu_*\ca{M}}\to 0$. Applying $\nu^{-1}$, we get
$1\to \nu^{-1}\ca{O}^\times_{Y_\et}\to \nu^{-1}\nu_*\ca{M}\to\nu^{-1}\overline{\nu_*\ca{M}}\to 0$. Moreover, we have a commutative diagram of exact sequences
\begin{equation}\label{eq-diag-chasing}
    \begin{tikzcd}
1\arrow[r]&\nu^*\ca{O}^\times_{Y_\et}\arrow[r]\arrow[d,"\simeq"]&\nu^*\nu_*\ca{M}\arrow[r]\arrow[d]& \nu^{-1}\overline{\nu_*\ca{M}}\arrow[r]\arrow[d]&0\\
1\arrow[r]&\ca{O}^\times_{Y_v}\arrow[r]&\ca{M}\arrow[r]&\overline{\ca{M}}\arrow[r]&0.
    \end{tikzcd}
\end{equation}
We first show that the composition
\begin{equation}\label{eq-adjuction-iso}
\nu^{-1}\overline{\nu_*\ca{M}}\to\nu^{-1}\nu_*\overline{\ca{M}}\to \overline{\ca{M}}
\end{equation} 
is injective. 
By the proof of Lemma \ref{lem-v-log-fullfaithful-adjunction}, for any $v$-cover $(f: X\to Y)\in Y_v$, $\nu^*\nu_*\ca{M}(X)$ (resp. $\nu^{-1}\overline{\nu_*\ca{M}}(X)$) is given by $f^*\nu_*\ca{M}(X)$ (resp. $f^{-1}\overline{\nu_*\ca{M}}(X)$) with a surjection $f^*\nu_*\ca{M}\to f^{-1}\overline{\nu_*\ca{M}}$ obtained by pulling back $\nu_*\ca{M}\to \overline{\nu_*\ca{M}}$. 
Fix such a morphism $f$. Suppose that there are $\overline{s}_1,\overline{s}_2\in f^{-1}\overline{\nu_*\ca{M}}(X)$ mapping to the same element $\overline{s}\in \overline{\ca{M}}(X)$. 
Up to replacing $X$ with an {\'e}tale cover, there is a cover $U\in X_{\et}$ that is the pullback of an {\'e}tale cover $U_Y$ of $Y$ via $f$, such that $\overline{s}_1, \overline{s}_2\in \overline{\nu_*\ca{M}}(U_Y)=\nu_*\overline{\ca{M}}(U_Y)=\overline{\ca{M}}(U_Y)$. The equation $\overline{\nu_*\ca{M}}=\nu_*\overline{\ca{M}}$ is proved in Lemma \ref{lem-nubarM-iso} below.
Moreover, they lift to $s_1,s_2\in \ca{M}(U_Y')$ for some cover $U_Y'\in U_{Y,\et}$, respectively. 
On the other hand, there is a $v$-cover $V\to U':=U_Y'\times_Y X$ with $u\in \ca{O}^\times(V)$ such that $s_1=s_2+u$ in $\ca{M}(V)$. The descent data of $s_1$ and $s_2$ from $V$ to $U'$, along with the integrality of $\ca{M}$, induces a descent datum of $u$ from $V$ to $U'$. So $u\in \ca{O}^\times(U')$ and $\overline{s}_1=\overline{s}_2$ in $\overline{\ca{M}}(U')$. We have shown the desired claim since $U'\to U_Y$ is a $v$-cover.\par
A similar argument as above implies that $\nu^*\nu_*\ca{M}\to \ca{M}$ is injective. For $Y'\in Y_{\et}$ and $s_1,s_2\in \nu^*\nu_*\ca{M}(Y')$ mapping to $s\in\ca{M}(Y')$. By diagram chasing, $s_1$ and $s_2$ project to the same $\overline{s}\in \nu^{-1}\overline{\nu_*\ca{M}}(Y')$ and $s_1=s_2+u'$ in $\nu^*\nu_*\ca{M}(Y'')$, where $Y''$ is a $v$-cover of $Y'$ and $u\in\ca{O}^\times(Y'')$. By the integrality of $\nu^*\nu_*\ca{M}$, $u$ descends to $\ca{O}^\times(Y')$. Mapping to $\ca{M}(Y')$, by the integrality of $\ca{M}$ and $\nu^*\ca{O}^\times_{Y_\et}=\ca{O}^\times_{Y_v}$, we have that $u=0$, as desired.\par
We next show that (\ref{eq-adjuction-iso}) is an isomorphism. In fact, by assumption, for any geometric point $\overline{x}\in X_{\et}$, there is an {\'e}tale neighborhood $U_{\overline{x}}$ of it such that $\ca{M}$ admits a chart modeled on a weakly $p$-finitely generated $\md{P}$ denoted by $\theta:\md{P}_{U_{\overline{x}}}\to\ca{M}|_{U_{\overline{x}}}$, and such that $\theta$ factors through $\nu^*\nu_*\ca{M}|_{U_{\overline{x}}}$. 
Then there is a morphism $\md{P}\to \nu^{-1}\overline{\nu_*\ca{M}}\to\nu^{-1}\nu_*\overline{\ca{M}}\to \overline{\ca{M}}$ over $U_{\overline{x}}$ whose composition is surjective. Then (\ref{eq-adjuction-iso}) is surjective, as desired.\par
Finally, the remainder of the lemma follows from a diagram chasing for (\ref{eq-diag-chasing}).
\end{proof}
In fact, we have
\begin{lem}\label{lem-nubarM-iso}
The natural morphism $\overline{\nu_*\ca{M}}\to \nu_*\overline{\ca{M}}$ is an isomorphism. 
\end{lem}
\begin{proof}
It is an injection, as there is a commutative diagram of exact sequences
\begin{equation*}\label{eq-diag-chasing-2}
    \begin{tikzcd}
1\arrow[r]&\ca{O}^\times_{Y_\et}\arrow[r]\arrow[d,"\simeq"]&\nu_*\ca{M}\arrow[r]\arrow[d]& \overline{\nu_*\ca{M}}\arrow[r]\arrow[d]&0\\
1\arrow[r]&\nu_*\ca{O}^\times_{Y_v}\arrow[r]&\nu_*\ca{M}\arrow[r]&\nu_*\overline{\ca{M}}&
    \end{tikzcd}
\end{equation*}
and the injectivity is proved by diagram chasing. We show that this is a surjection. It can be checked that, for $U\in Y_{\et}$, $\overline{\nu_*\ca{M}}(U)=$
$$\{s\in \overline{\ca{M}}(U)|s \text{ can be lifted to }s'\in \ca{M}(V)\text{ for some {\'e}tale cover }V\to U\}.$$
Fix $s\in \overline{\ca{M}}(U)$ for $U\in Y_{\et}$. The functor $\overline{\ca{M}}_s$ sending $V\in U_v$ to $\{s'\in\ca{M}(V)|s'\mapsto s|_V\}$ is a $\bb{G}_m$-torsor on $U_v$. By \cite[Lem. 17.1.8]{SW20} and \cite[Thm. 2.5.8]{KL19}, $\overline{\ca{M}}_s$ is the pullback of an {\'e}tale $\bb{G}_m$-torsor $\overline{\ca{N}}_s=\nu_*\overline{\ca{M}}_s$. So $\overline{\ca{N}}_s$ {\'e}tale locally has a section, and $s$ {\'e}tale locally admits a lifting. 
\end{proof}
We now show the second paragraph of Lemma \ref{lem-v-log-str}.
\begin{lem}\label{lem-pullback-saturated}
Let $\ca{L}\in\ca{LOG}^{\mrm{int}}_{Y_{\et}}$. Then $\ca{L}$ is fine perfectoid if and only if $\nu^*\ca{L}$ is. In this case, $\ca{L}$ is saturated if and only if $\nu^*\ca{L}$ is so. 
\end{lem}
\begin{proof}
For fine perfectoidness, the ``only if'' part is Corollary \ref{cor-fully-faithful-fine-perf}. 
For the ``only if'' part of saturatedness, if $\ca{L}$ is saturated, then, by the construction in Lemma \ref{lem-v-log-fullfaithful-adjunction}, we see that $f^{-1}\overline{\ca{L}}$ is saturated by \cite[Lem. 2.2.4]{DLLZ23}, and $f^*\ca{L}$ is saturated by \cite[I. Prop. 1.3.5]{Ogu18} and \cite[Lem. 2.2.4]{DLLZ23}.\par 
The ``if'' parts are immediate from Lemma \ref{lem-v-log-fullfaithful-adjunction} and the definitions.
\end{proof}
This completes the proof of Lemma \ref{lem-v-log-str}.
\subsubsection{End of the proof}\label{subsubsec-end-proof-v-sheaves}
From Lemma \ref{lem-v-log-str}, we obtain the following proposition:
\begin{prop}\label{prop-log-diamond-v-sheaves-bis}
Let $(Y,\ca{M}_Y)$ be a log perfectoid space (cf. \cite[Def. 2.2.2(9)]{DLLZ23}) that is fs. Then the $v$-descent data of the functor sending $(Y',f:Y'\to Y)\in \Perfd_{/Y}$ to the groupoid $$\{\text{saturated and fine perfectoid log structures }\ca{M} \text{ on }Y'_{\et} \text{ with }f^*\ca{M}_Y\to \ca{M}\}$$ are effective.  
\end{prop}
\begin{proof}
By Lemma \ref{lem-v-log-str}, the functor in the proposition is isomorphic to the functor sending $(Y',f:Y'\to Y)\in \Perfd_{/Y}$ to the groupoid $$\{\text{saturated and fine perfectoid log structures }\ca{M} \text{ on }Y_v' \text{ with }\nu^*f^*\ca{M}_Y\to \ca{M}\}.$$
Let $h:Y''\to Y'$ be a $v$-cover. Let $(\ca{M}'',\alpha'')$ be a saturated and fine perfectoid log structure on $Y''_v$ equipped with a descent datum $\sigma:p_1^*\ca{M}''\simeq p_2^*\ca{M}''$ on $Y''\times_{Y'}Y''$ which satisfies the following conditions:
\begin{itemize}
\item $p_2^*\alpha''\circ\sigma=\sigma_0 \circ p_1^*\alpha''$, where $\sigma_0$ is the canonical descent datum associated with pulling back $\ca{O}_{Y'}$ to $Y''$ via $h$. 
\item There is a morphism $h^{\sharp}:\nu^* h^*f^*\ca{M}_Y\to \ca{M}''$ such that $(p_2^*h^\sharp)\circ \sigma_Y=\sigma\circ(p_1^*h^\sharp)$, where $\sigma_Y$ is the canonical descent datum of $\nu^*h^*f^*\ca{M}_Y=h^*\nu^*f^*\ca{M}_Y$ induced by pulling back $\nu^*f^*\ca{M}_Y$ via $h$. 
\end{itemize}
This descent datum determines and is determined by a log structure $(\ca{M}',\alpha')$ on $Y'_v$ together with a morphism $\nu^*f^*\ca{M}_Y\to \ca{M}'$. 
This log structure $\ca{M}'$ is saturated and fine perfectoid. Taking $\nu_*$, we get the desired morphism $f^*\ca{M}_Y\to \nu_* \ca{M}'$ (cf. Lemma \ref{lem-v-log-fullfaithful-adjunction}). Then the Proposition follows from Lemma \ref{lem-v-log-str}.
\end{proof}
We then have the following corollary:
\begin{cor}\label{cor-log-diamond-v-sheaves-bis}
Let $(X,\ca{M}_X)$ be one of the Cases \ref{case-gen-log-adic-space}-\ref{case-gen-for-sch}.\par 
Then the functor ``$\Hom(-,(X,\ca{M}_X))$'' on $\Perfd$ defined as follows is a $v$-sheaf:
$$Y\in \Perfd\longmapsto \{f:(Y,\ca{M})\to (X,\ca{M}_X)|\ca{M}\text{ is saturated and fine perfectoid}\}/\simeq.$$
\end{cor}
\begin{proof}
By \cite[Thm. 17.1.3 and Thm. 18.1.1]{SW20}, $\Hom(-,X)$ is a $v$-sheaf on $\Perfd$. We then get the desired result by applying Proposition \ref{prop-log-diamond-v-sheaves-bis} where $\ca{M}_Y=f^*\ca{M}_X$.    
\end{proof}
\begin{proofof}[Theorem \ref{thm-log-diamond-v-sheaves}]
Now we can show Theorem \ref{thm-log-diamond-v-sheaves} with these facts. 
Indeed, since $\mrm{Untilt}$ is a $v$-sheaf by \cite[Lem. 15.1]{scholze2017etale}, we only have to show that $\Hom(-,(X,\ca{M}_X))$ is a $v$-sheaf on $\Perfd$. But this is Corollary \ref{cor-log-diamond-v-sheaves-bis}. 
\end{proofof}
\subsubsection{}Suppose that $X$ is of finite type over $\bb{Z}_p$ and satisfies (SF). 
It follows from the definition that there is a natural ``structural morphism'' $X^{\log\diamond}\to \Spd \bb{Z}_p$.
Moreover, 
\begin{lem}\label{lem-v-sheaf-map}
There is a natural injective morphism between $v$-sheaves $$X^{\log\diamond}\times_{\Spd\bb{Z}_p}\Spd\bb{Q}_p\lra(X_{\bb{Q}_p})^{\log\Diamond}.$$
\end{lem}
\begin{proof}
For any $S:=\Spa(A,A^+)\in\Perf$, $X^{\log\diamond}\times_{\Spd\bb{Z}_p}\Spd\bb{Q}_p(S)$ parametrizes isomorphism classes of $\{S^\sharp=\Spa(A^\sharp,A^{\sharp,+}),f: (S^\sharp,\ca{M}_{S^\sharp})\to(\wat{X}^\ad,\wat{\nu}^{-1}_X\ca{M}_X)\}$ where $A^\sharp$'s are $\bb{Q}_p$-algebras. To define a class in $(X_{\bb{Q}_p})^{\log\Diamond}(S)$, it remains to define a morphism 
$f_{\bb{Q}_p}:(S^\sharp,\ca{M}_{S^\sharp})\to((X_{\bb{Q}_p})^\ad,\nu^{\ad,-1}_{X_{\bb{Q}_p}}\ca{M}_{X_{\bb{Q}_p}})$. 
For the morphism $f_{\bb{Q}_p}: S^\sharp\to (X_{\bb{Q}_p})^\ad$, it is defined as $S^{\sharp}\xrightarrow{f} \wat{X}^\ad_\eta\hookrightarrow (X_{\bb{Q}_p})^\ad$; for the morphism between sheaves $\nu^{\ad,-1}_{X_{\Qp}}\ca{M}_{X_{\bb{Q}_p}}\to \ca{M}_{S^\sharp}$, note that there is a commutative diagram of sites
\begin{equation*}
    \begin{tikzcd}
(\wat{X}^\ad_\eta)_{\et}\arrow[r,"i_1"]\arrow[d,"j_1"]\arrow[ddr,"\wat{\nu}_{X_{\bb{Q}_p}}"description]&(\wat{X}^\ad)_{\et}\arrow[d,"j_2"]\arrow[ddr,"\wat{\nu}_X"description]&\\
((X_{\bb{Q}_p})^\ad)_{\et}\arrow[r,"i_2"]\arrow[dr,"\nu^\ad_{X_{\bb{Q}_p}}"description]& (X^\ad)_{\et}\arrow[dr,"\nu^\ad_X"description]&\\
{}&X_{\bb{Q}_p,\et}\arrow[r,"i_3"]&X_\et.
    \end{tikzcd}
\end{equation*}
For a fixed $f:S^\sharp\to \wat{X}_\eta^\ad$, a morphism $f^{-1}i_1^{-1}\wat{\nu}_X^{-1}\ca{M}_X\to \ca{M}_{S^{\sharp}}$ corresponds to a morphism $\alpha_0:f^{-1}j_1^{-1}\nu^{\ad,-1}_{X_{\bb{Q}_p}} i_3^{-1}\ca{M}_X\to \ca{M}_{S^\sharp}$.
The latter one induces a morphism $\alpha:f^{-1}j_1^{-1}\nu^{\ad,-1}_{X_{\bb{Q}_p}} i_3^*\ca{M}_X\to \ca{M}_{S^\sharp}$, where $i_3^*\ca{M}_X = \ca{M}_{X_{\Qp}}$. Conversely, $\alpha$ also uniquely determines $\alpha_0$. So the injectivity follows from the injectivity of $X^\diamond\times_{\Spd \bb{Z}_p}\Spd \bb{Q}_p\to X_{\bb{Q}_p}^\Diamond.$
\end{proof}
\begin{definition}\label{def-slashed-diamond} Let $(X,\ca{M}_X)$ be an fs log scheme over $\Spec \bb{Z}_p$.
Define the slashed log diamond $(X,\ca{M}_X)^{\Diamond/}$ (or denoted by $X^{\log\Diamond/}$) as the quotient $v$-sheaf
$$(X,\ca{M}_X)^{\Diamond/}:=X^{\log\diamond}\disju_{X^{\log\diamond}\times_{\Spd\bb{Z}_p}\Spd\bb{Q}_p}(X_{\bb{Q}_p})^{\log\Diamond}.$$
Note that the definition makes sense by Theorem \ref{thm-log-diamond-v-sheaves}.
\end{definition}
\begin{lem}\label{lem-proper-equivalent}
When $X$ is proper over $\bb{Z}_p$, the natural injections $X^{\log\diamond}\hookrightarrow X^{\log\Diamond/}\hookrightarrow X^{\log\Diamond}$ are isomorphisms. 
\end{lem}
\begin{proof}
This follows from the construction and the fact that $X^\ad\iso \wat{X}^\ad$ when $X$ is proper over $\bb{Z}_p$.  
\end{proof}
Let us also remark that there is another way of defining a log diamond, which is a generalization of \cite[Def. 8.3.1]{SW20}.
\begin{rk}\label{rk-log-dia-alt-def}
Let $\mrm{FSPerf}$ be the category of saturated and fine perfectoid log perfectoid spaces of characteristic $p$. We then define a log diamond $X$ as a pro-{\'e}tale sheaf on $\mrm{FSPerf}$ such that: (1) There is a surjective map between sheaves $Y\to X$, where $Y$ is representable by an object in $\mrm{FSPerf}$; (2) the fiber product of sheaves $Y\times_X Y$ is also representable by an object in $\mrm{FSPerf}$, denoted by $Y\times_X^{sat} Y$, such that the two projections $p_1$, $p_2:Y\times_X^{sat} Y\to Y$ are (non-log) pro-{\'e}tale. \par
We expect that Definition \ref{def-log-diamonds-adic-space} should satisfy this definition when $X$ is analytic over $\Spa \bb{Z}_p$ (cf. Lemma \ref{lem-maxproket-perfectoid}). We do not explore this direction in the current paper.\par
One can also change the definition of log diamonds by changing the category $\mrm{FSPerf}$ to other categories with a different class of log structures. This definition may also be viewed as an analogue of the definition of a log algebraic space in the second sense due to Kajiwara-Kato-Nakayama (see \cite[10.1]{KKN15}) in perfectoid geometry. 
The definition of log diamonds in \cite[Def. 7.3]{KY25} may also be viewed as an analogue of log algebraic spaces in the first sense in \cite{KKN15} for some category of log perfectoid spaces.
\end{rk}
\subsection{Log shtukas}\label{subsec-log-shtuka}
We define a notion of $p$-adic shtukas in log geometry.\par 
Let $(X,\ca{M}_X)$ be one of the cases in \S\ref{subsubsec-log-diamonds-setup}. 
\subsubsection{}\label{subsubsec-log-diamond-category}
As \S\ref{subsec-shtuka-nonred}, let $G$ be a connected linear algebraic group over $\bb{Q}$ that has a quasi-parahoric model $\ca{G}$ over $\bb{Z}_p$ in the sense of Definition \ref{def: quasi-parahoric for non-reductive group}. Note that we still use the symbol $G$ and $\G$ here, but $G$ might not be reductive.\par
Let $\mu:\bb{G}_{m,\overline{\bb{Q}}_p}\to G_{\overline{\bb{Q}}_p}$ be a minuscule cocharacter in the sense of Definition \ref{def-minuscule-parabolic} (2), and denote by $\{\mu\}$ the $G(\overline{\bb{Q}}_p)$-conjugacy class of it. Let $E$ be the field of definition of $\{\mu\}$.\par
In what follows, $?$ denotes $\Dia$ or $\dia$. Since $(X,\ca{M}_X)^?$ is a functor from $\Perf$ to isomorphism classes, we can and will regard $(X,\ca{M}_X)^?$ as a functor $(X,\ca{M}_X)^?:\Perf^{\mrm{op}}\to \md{Categories}$ that sends $S\in \Perf$ to the collection of tuples $(S^\sharp,\ca{M}_{S^\sharp},\alpha_{S^\sharp},f)$ viewed as a category, which is still denoted by $(X,\ca{M}_X)^?(S)$. More precisely, the objects of $(X,\ca{M}_X)^?(S)$ are those $(S^\sharp,\ca{M}_{S^\sharp},\alpha_{S^\sharp},f)$, while the morphisms between $(S^\sharp_1,\ca{M}_{S^\sharp_1},\alpha_{S^\sharp_1},f_1)$ and $(S^\sharp_2,\ca{M}_{S^\sharp_2},\alpha_{S^\sharp_2},f_2)$ are the morphisms between log adic spaces $$g:(S^\sharp_1,\ca{M}_{S^\sharp_1},\alpha_{S^\sharp_1})\to (S^\sharp_2,\ca{M}_{S^\sharp_2},\alpha_{S^\sharp_2})$$ such that $f_2\circ g=f_1$. \par 
By Lemma \ref{lem-v-log-str}, Theorem \ref{thm-log-diamond-v-sheaves} and Proposition \ref{prop-log-diamond-v-sheaves-bis}, this functor is a $v$-stack, which we abusively denote by $$p^?:(X,\ca{M}_X)^?\to \Perf.$$ 
\subsubsection{}\label{subsec-log-shtuka-def}
On the other hand, let $$\Sht_{\G,\mu}:\Perf^{\mrm{op}}\to \md{Groupoids}$$
be the $v$-stack over $\Perf$ sending $S\in \Perf$ to the groupoid with objects $(S^\sharp;(\mathscr{P},\phi_\mathscr{P}))$, where $S^\sharp$ is an untilt of $S$ and $(\mathscr{P},\phi_\mathscr{P})$ is a $\G$-shtuka over $S$ with one leg at $S^\sharp$ bounded by $\mu$ (see \cite[Def. 2.4.3]{PR24}). For a fixed $S^\sharp$, let $\Sht_{\G,\mu}(S^\sharp)$ be the groupoid of $\G$-shtukas over $S$ with one leg at $S^\sharp$ bounded by $\mu$.\par
Write the fibered category in groupoids corresponding to $\Sht_{\G,\mu}$ by $p_{\Sht}:\ca{SHT}_{\G,\mu}\to \Perf$.\par 
Note that there is a commutative diagram of fibered categories
\begin{equation}\label{diag-sht-def-setup}
    \begin{tikzcd}
       \ca{SHT}_{\G,\mu/(X,\ca{M}_X)^?}\arrow[rr,"p_1"]\arrow[drr,"{p^\sharp((X,\ca{M}_X)^?)}"description]&& \ca{SHT}_{\G,\mu}\arrow[rr,"p_{\Sht}"]\arrow[drr,"p^\sharp"]&&\Perf\\
       && (X,\ca{M}_X)^?\arrow[rr,"p_2"]&&\mrm{Untilt}.\arrow[u,"\flat"]
    \end{tikzcd}
\end{equation}
In the diagram above, $p^\sharp$ (resp. $p_2$) is a functor sending $(S^\sharp;(\mathscr{P},\phi_\mathscr{P}))$ (resp. $(S^\sharp,\ca{M}_{S^\sharp},f)$) to $S^\sharp$ with morphisms also projected to the morphisms of the first factor.
The fibered category $p^\sharp((X,\ca{M}_X)^?):\ca{SHT}_{\G,\mu/(X,\ca{M}_X)^?}\to (X,\ca{M}_X)^?$ is defined by a functor $\Sht_{\G,\mu/(X,\ca{M}_X)^?}$ sending $(S^\sharp,\ca{M}_{S^\sharp},f)$ to $\Sht_{\G,\mu}(S^\sharp)$. The functor $p_1$ is the natural projection.
In conclusion, this is the situation to which Appendix \ref{subsec-limits-def} applies. When $?=\Diamond/$, by Definition \ref{def-slashed-diamond}, $(X,\ca{M}_X)^{\Diamond/}$ also satisfies the diagram above because it is true for both small and big diamonds.
\begin{definition}\label{def-log-shtuka}
Suppose that $(X,\ca{M}_X)$ is in one of the Cases \ref{case-gen-log-adic-space}-\ref{case-gen-for-sch} in \S\ref{subsubsec-log-diamonds-setup}.  
Define the groupoid of \textbf{log $\G$-shtukas bounded by $\mu$} on $(X,\ca{M}_X)^?$ as (see Definition \ref{def-2-lim})
$$\Sht_{\ca{G},\mu}^?(X,\ca{M}_X):=\twolim\limits_{(S^{\sharp},\ca{M}_{S^\sharp},f_{S^\sharp})\in ((X,\ca{M}_X)^{?})^{\mrm{op}}}\Sht_{\ca{G},\mu}(S^\sharp),$$
where $?=\Diamond$ or $\diamond$.\par
Similarly, we define the groupoid of \emph{$p$-adic $\G$-shtukas} on $(X,\ca{M}_X)^?$ as 
$$\Sht_{\ca{G}}^?(X,\ca{M}_X):=\twolim\limits_{(S^{\sharp},\ca{M}_{S^\sharp},f_{S^\sharp})\in ((X,\ca{M}_X)^{?})^{\mrm{op}}}\Sht_{\ca{G}}(S^\sharp).$$
\end{definition}
The following statement follows immediately from Lemma \ref{lem-mor-lim-correspond}.
\begin{prop}\label{prop-def-shtukas}
    With the definitions above, the following two groupoids are canonically \emph{isomorphic}:
    \begin{enumerate}
    \item The groupoid $\Sht^?_{\ca{G},\mu}(X,\ca{M}_X)$;
    \item The groupoid of $1$-morphisms between $v$-stacks on $\Perf$
        $$\mathscr{P}^{\log?}:(X,\ca{M}_X)^{?}\lra \Sht_{\ca{G},\mu}$$
    that are sections of $p^\sharp((X,\ca{M}_X)^?)$.
    \end{enumerate}
\end{prop}
We call an object in either of the two groupoids a \emph{log shtuka}.
\begin{definition}
In general, we can also define $\Sht_{\G,\mu}(Y)$ for any $v$-stack $Y$ mapping to $\mrm{Untilt}$ satisfying the diagram (\ref{diag-sht-def-setup}) as the groupoid of $1$-morphisms $Y\to \Sht_{\G,\mu}$.\par 
Denote $\Sht^{\Dia/}_{\G,\mu}(X,\ca{M}_X):=\Sht_{\G,\mu}((X,\ca{M}_X)^{\Dia/})$. 
\end{definition}
For example, we can define $\Sht_{\G,\mu}((X,\ca{M}_X)^\dia\times_{\Spd \bb{Z}_p}\Spd \bb{Q}_p)$ as a $2$-limit over $(X,\ca{M}_X)^\dia\times_{\Spd \bb{Z}_p}\Spd \bb{Q}_p$ as in Definition \ref{def-log-shtuka} or the groupoid of $1$-morphisms from $(X,\ca{M}_X)^\dia\times_{\Spd \bb{Z}_p}\Spd \bb{Q}_p$ to $\Sht_{\G,\mu}$.\par
The following lemma is immediate from definition:
\begin{lem}\label{lem-str-etale-descent}
In all cases for $(X,\ca{M}_X)$ in \S\ref{subsubsec-log-diamonds-setup}, (strict) {\'e}tale descent for log shtukas is effective. 
\end{lem}
\begin{proof}
Let $X'\to X$ be an {\'e}tale cover. Denote by $S^{\sharp,\prime}\to S^\sharp$ the pullback of $X'\to X$ via $f:S^\sharp\to X$. Since the descent for the {\'e}tale cover $S^{\sharp,\prime}\to S^\sharp$ is effective for $\Sht_{\G,\mu}$ (see \cite[Prop. 19.5.3]{SW20}), we have the desired assertion by taking limits.
\end{proof}
\subsubsection{}
\begin{lem}\label{lem-can-log-str}
Let $(X,\ca{M}_X)$ be defined as in \S\ref{subsubsec-log-diamonds-setup}. Let $S\in \Perf$. Suppose that $X^{\log \Diamond}(S)$ is nonempty. Let $(S^\sharp,\ca{M}_{S^\sharp},\alpha_{S^\sharp},f)\in X^{\log \Diamond}(S)$. Then there is a saturated and strongly fine perfectoid log structure $(\ca{M}_{S^\sharp}^{\mrm{can}},\alpha^{\mrm{can}}_{S^\sharp})$ that satisfies the following universal property:\par
For any object $(S^\sharp,\ca{M}_{S^\sharp}',\alpha'_{S^\sharp},f)\in X^{\log \Diamond}$ (for the same $S^\sharp$ but varying saturated and fine perfectoid log structures), there is a uniquely determined morphism between log perfectoid spaces $c_{\ca{M}_{S^\sharp}'}:(S^\sharp,\ca{M}_{S^\sharp}',\alpha'_{S^\sharp})\to (S^\sharp,\ca{M}_{S^\sharp}^{\mrm{can}},\alpha^{\mrm{can}})$ extending the identity morphism of the underlying perfectoid spaces. \par
Moreover, for any morphism $g:(S^\sharp_1,\ca{M}_{S^\sharp_1},\alpha_{S^\sharp_1})\to (S^\sharp_2,\ca{M}_{S^\sharp_2},\alpha_{S^\sharp_2})$ in the fibered category $(X,\ca{M}_X)^{\Diamond}$ over $\Perf$, there is a unique morphism $$g^\can:(S^\sharp_1,\ca{M}_{S^\sharp_1}^\can,\alpha_{S^\sharp_1}^\can)\to (S^\sharp_2,\ca{M}_{S^\sharp_2}^\can,\alpha_{S^\sharp_2}^\can)$$
such that $g^\can\circ c_{\ca{M}_{S^\sharp_1}}=c_{\ca{M}_{S^{\sharp}_2}}\circ g$.\par
Similar results hold if we replace ``$\log \Diamond$'' with ``$\log \diamond$''.
\end{lem}
\begin{proof}
We only show the assertions in Case \ref{case-gen-log-sch} for ``$\log \Diamond$''; other cases are proved in the same way. Let us fix $(S^\sharp,\ca{M}_{S^\sharp},\alpha_{S^\sharp},f)$.\par 
The morphism $f^\sharp:(\ca{M}_{S^\sharp},\alpha_{S^\sharp})\to (\ca{M}_X,\alpha_X)$ gives a commutative diagram:
\begin{equation*}
    \begin{tikzcd}
&\ca{M}_{S^\sharp}\arrow[rr,"\alpha_{S^\sharp}"]\arrow[dd]&&\ca{O}_{S^\sharp_\et}\arrow[dd]\\
f^{-1}\ca{M}_X\arrow[rr,"\alpha"description]\arrow[dd]\arrow[ur,"f^\sharp"]&&f^{-1}\ca{O}_{X_\et}\arrow[ur,"f^\sharp"]\arrow[dd]&\\
&\overline{\ca{M}}_{S^\sharp}\arrow[rr,"\overline{\alpha}_{S^\sharp}"]&&\ca{O}_{S^\sharp_\et}/\ca{O}^\times_{S^\sharp_\et}\\
f^{-1}\overline{\ca{M}}_X\arrow[rr,"\overline{\alpha}_X"]\arrow[ur,"\overline{f}^\sharp"]&&f^{-1}(\ca{O}_{X_\et}/\ca{O}_{X_\et}^\times).\arrow[ur,"\overline{f}^\sharp"]&
    \end{tikzcd}
\end{equation*}   
Let $\overline{\ca{M}}^\mrm{canperf}_{S^\sharp}:=(f^{-1}\overline{\ca{M}}_X)^{\mrm{perf}}=(\overline{f^*\ca{M}_X})^{\mrm{perf}}$. 
As $\overline{\ca{M}}_{S^\sharp}$ is uniquely $p$-divisible, there is a unique homomorphism $\overline{\alpha}^{\can}:\overline{\ca{M}}^{\mrm{canperf}}_{S^\sharp}\to \overline{\ca{M}}_{S^\sharp}$.
Denote by $\ca{M}_{S^\sharp}^{\mrm{canperf}}$ the pullback of $\ca{M}_{S^\sharp}$ via $\overline{\alpha}^{\can}$.
That is, $\ca{M}_{S^\sharp}^{\mrm{can}}$ fits into the following commutative diagram as an extension of $\overline{\ca{M}}_{S^\sharp}^{\mrm{canperf}}$ by $\ca{O}^\times_{S^\sharp_\et}$
\begin{equation*}
    \begin{tikzcd}
    1\arrow[r]& \ca{O}^\times_{S^\sharp_\et}\arrow[r]\arrow[d]&\ca{M}^{\mrm{can}}_{S^\sharp}\arrow[r]\arrow[d,"c_{\ca{M}_{S^\sharp}}"]&\overline{\ca{M}}_{S^\sharp}^{\mrm{canperf}}\arrow[r]\arrow[d,"{\overline{\alpha}}^{\mrm{can}}"]&0\\
    1\arrow[r]&\ca{O}^\times_{S^\sharp_\et}\arrow[r]&\ca{M}_{S^\sharp}\arrow[r]&\overline{\ca{M}}_{S^\sharp}\arrow[r]&0.
    \end{tikzcd}
\end{equation*}
Let $\alpha^\can_{S^\sharp}:=\alpha_{S^\sharp}\circ c_{S^\sharp}$. We claim that $(\ca{M}^\can_{S^\sharp},\alpha^\can)$ is the desired pair of log structure.\par
To check the universal property, it suffices to show that the construction is independent of the choice of $f^\sharp:(\ca{M}_{S^\sharp},\alpha_{S^{\sharp}})\to (\ca{M}_X,\alpha_X)$ for a fixed $f:S^\sharp\to X$.
Given two such morphisms $f_1^\sharp:(\ca{M}^1_{S^\sharp},\alpha^1_{S^\sharp})\to (\ca{M}_X,\alpha_X)$ and $f_2^\sharp:(\ca{M}^2_{S^\sharp},\alpha^2_{S^\sharp})\to (\ca{M}_X,\alpha_X)$ that are compatible with a fixed $f:S^\sharp\to X$, we show the definition of $c_{S^\sharp}$ is independent of this choice.
Indeed, there is a commutative diagram
\begin{equation}
\begin{tikzcd}
    && \overline{\ca{M}}^1_{S^\sharp}\arrow[drr,"{\overline{\alpha}^1_{S^\sharp}}"]&&\\
    f^{-1}\overline{\ca{M}}_X\arrow[urr,"{\overline{f}^{1,\sharp}}"]\arrow[drr,"{\overline{f}^{2,\sharp}}"]&&&&\ca{O}_{S^\sharp_\et}/\ca{O}^\times_{S^\sharp_{\et}}\\
    &&\overline{\ca{M}}^2_{S^\sharp}.\arrow[urr,"{\overline{\alpha}^2_{S^\sharp}}"]&&
    \end{tikzcd}
\end{equation}
The map $\overline{f}^{i,\sharp}$ uniquely determines a homomorphism $\overline{\alpha}^{\can,i}:\overline{\ca{M}}_{S^\sharp}^{\mrm{canperf}}\to\overline{\ca{M}}^i_{S^\sharp}$ for $i=1,2$ and 
$\ca{M}_{S^\sharp}^\can$ is constructed by pulling back $\ca{M}^i_{S^\sharp}\to \overline{\ca{M}}^i_{S^\sharp}$ via $\overline{\alpha}^{\can,i}$. 
We define $(\ca{M}^3_{S^\sharp},\alpha^3_{S^\sharp})$ as $((\ca{M}^1_{S^\sharp},\alpha^1_{S^\sharp})\oplus (\ca{M}^2_{S^\sharp},\alpha^2_{S^\sharp}))/\ca{O}^\times_{S^\sharp}$, which is saturated by \cite[I. Prop. 1.3.4]{Ogu18} and \cite[I. Prop. 1.3.5]{Ogu18}. Moreover, $\overline{\ca{M}}^3_{S^\sharp}\iso \overline{\ca{M}}^1_{S^\sharp}\oplus \overline{\ca{M}}^2_{S^\sharp}$ is uniquely $p$-divisible. We have the following commutative diagram 
\begin{equation}
\begin{tikzcd}
    && \ca{M}^1_{S^\sharp}\arrow[drr,"{\alpha^1_{S^\sharp}}"]&&\\
    \ca{M}^\can_{S^\sharp}\arrow[urr,"{c_{\ca{M}_{S^\sharp}^1}}"]\arrow[rr,"{c_{\ca{M}_{S^\sharp}^3}}"]\arrow[drr,"{c_{\ca{M}_{S^\sharp}^2}}"']&&\ca{M}^3_{S^\sharp}\arrow[u,"p_1"]\arrow[d,"p_2"]&&\ca{O}_{S^\sharp_\et}\\
    &&\ca{M}^2_{S^\sharp},\arrow[urr,"{\alpha^2_{S^\sharp}}"]&&
    \end{tikzcd}
\end{equation}
as desired. The functoriality in the third paragraph follows from the construction of $(\ca{M}^\can_{S^\sharp},\alpha^\can_{S^\sharp})$, $c_{\ca{M}_{S^\sharp}}$ and the unique $p$-divisibility of $\overline{\ca{M}}_{S_1^\sharp}$ and $\overline{\ca{M}}_{S_2^\sharp}$.\par
The resulting $\ca{M}_{S^\sharp}^{\mrm{can}}$ is saturated and strongly fine perfectoid. Indeed, since $(X,\ca{M}_X)$ is fs, by \cite[Prop. 2.3.13]{DLLZ23} and \cite[I. 1.3.6]{Ogu18}, $\ca{M}_{S^\sharp}^\can$ admits saturated, uniquely $p$-divisible, and $p$-finitely generated charts. By \cite[I. Prop. 1.3.4, Prop. 1.3.5]{Ogu18}, $\ca{M}_{S^\sharp}^\can$ is saturated; it suffices to prove perfectoidness. This follows again from the construction of $\ca{M}^\can_{S^\sharp}$ as a pullback of $\ca{M}_{S^\sharp}$ via $\overline{\alpha}^\can$ and the perfectoidness of $\ca{M}_{S^\sharp}$.
\end{proof}
\begin{prop}\label{prop-shu-can-obj}
Denote by $\mathfrak{f}: X^{\log\Diamond} \to X^\Diamond$ the natural projection. With the conventions as in Lemma \ref{lem-can-log-str}, the object $(S^\sharp,\ca{M}_{S^\sharp}^{\mrm{can}},\alpha^{\mrm{can}},f)$ is the final object of the fiber of $\mathfrak{f}$ at $(S^\sharp,f)$ (if the fiber is non-empty). Then a log shtuka $$\mathscr{P}^{\log\Diamond}:X^{\log\Diamond}\lra \Sht_{\G,\mu}$$ 
is uniquely determined by its restriction to the objects $(S^\sharp,\ca{M}_{S^\sharp}^\mrm{can},\alpha^\mrm{can})$ up to isomorphisms. The same is true for ``$\diamond$''.
\end{prop}
\begin{proof}
    This follows from Lemma \ref{lem-can-log-str} and Lemma \ref{lem-equiv-subcat-sec}.
\end{proof}
\begin{cor}\label{cor-trivial-log}
If $X$ is equipped with the trivial log structure $\ca{M}_X=\ca{O}^*_{X_\et}$, then a log shtuka
$$\mathscr{P}^{\log\Diamond}: X^{\log\Dia}\lra \Sht_{\G,\mu}$$
determines and is determined by its restriction to the objects $(S^\sharp,\ca{O}^*_{S^\sharp})$, i.e., a nonlog shtuka: The projection $X^{\log\Dia}\to X^\Dia$ admits a canonical section, and therefore determines a nonlog shtuka. The same is true for ``$\dia$'' and ``$\Dia/$''.
\end{cor}
\begin{proof}
    This follows from the same argument as that of Proposition \ref{prop-shu-can-obj}. Note that, in this case, the canonical log structures are all trivial. The case of ``$\Diamond/$'' follows from Definition \ref{def-slashed-diamond} as the corollary is true for both small and big diamonds.
\end{proof}
\begin{lem}\label{lem-compatible-charts}
Continuing with the conventions in Proposition \ref{prop-shu-can-obj} with $X$ satisfying (SF), we consider the small diamond $\mathfrak{f}:X^{\log \diamond}\to X^\diamond$. 
The map $f:(S^\sharp:=\Spa(R^\sharp,R^{\sharp,+}), \ca{M}_{S^\sharp}^\can,\alpha^\can)\to (X,\ca{M}_X)$ by Definition \ref{def-small-big-dia-gen-scheme} factors through a map $f^+:(\Spec R^{\sharp,+},\ca{M}_{R^{\sharp,+}},\alpha_{R^{\sharp,+}})\to (X,\ca{M}_X)$ such that the pullback log structure $\wat{\ca{M}}^{\ad}_{R^{\sharp,+}}$ of $\ca{M}_{R^{\sharp,+}}$ to $\Spa(R^{\sharp,+},R^{\sharp,+})$ via $\Spa(R^{\sharp,+},R^{\sharp,+})\to \Spec R^{\sharp,+}$ induces the log structure $\ca{M}_{S^\sharp}^{\can}$ on $S^\sharp$ via the pullback through $S^\sharp\to \Spa(R^{\sharp,+},R^{\sharp,+})$. 
Moreover, $f$ {\'e}tale locally admits a chart of the form $\md{P}\hookrightarrow \md{P}[\frac{1}{p}]$.   
\end{lem}
\begin{proof}
The proof of the first assertion follows from the same argument as Lemma \ref{lem-can-log-str} after the following changes: Here $X$ is a scheme; we replace $\overline{\ca{M}}_{S^\sharp}^{\mrm{canperf}}$ with the log structure induced by $\overline{\ca{M}}_{R^{\sharp,+}}:=(f^{+,-1}\overline{\ca{M}}_X)^{\perf}$ by pulling back via $S^\sharp\to \Spec R^{\sharp,+}$; the map $\overline{\alpha}^\can$ still exists by Lemma \ref{lem-satfineperf-implies-chart}, since the chart, by definition, factors through $R^{\sharp,+}$. The rest of the proof remains unchanged.\par
For the second assertion, assume that $X$ is affine and admits a global fs sharp chart $\md{P}$. By construction, there is a map $\theta: \md{P}\to \ca{M}_{R^{\sharp,+}}\to R^{\sharp,+}$. Since $\ca{M}_{R^{\sharp,+}}$ is integral and perfectoid, by \cite[Rmk. 2.25 and Rmk. 2.26]{KY25}, the induced log structure $\ca{M}_{R^{\sharp,+}/p}$ on the mod-$p$ log ring is uniquely $p$-divisible, and $\theta/p: \md{P}\to \ca{M}_{R^{\sharp,+}/p}\to R^{\sharp,+}/p$ can be extended to $\wdtd{\theta/p}: \md{P}[\frac{1}{p}]\to \ca{M}_{R^{\sharp,+}/p}\to R^{\sharp,+}/p$. This, in turn, induces a map 
    $$\xi\circ[\wdtd{\theta/p}^\flat]:\md{P}[\frac{1}{p}]\to W(R^+)\to R^{\sharp,+},$$
which induces the log structure $\ca{M}_{R^{\sharp,+}}$. The rest of the assertion follows from the construction in the last paragraph.
\end{proof}
\subsection{Equivalence of categories on generic fiber}\label{subsec-eq-gen}
We discuss a generalization of \cite[\S 2.5]{PR24} in log geometry.
\subsubsection{}\label{subsubsec-locsys}
To formulate the equivalence of categories, we need a notion of \emph{(pro-)$p$-Kummer {\'e}tale local systems} (see \cite{IKY26})\footnote{
We are grateful to Inoue and Koshikawa for suggesting the name ``(pro-)$p$-Kummer''; a similar one ``$p$-primary Kummer'' was also suggested to one of us (PW) during his visit to ZJU. In an update \cite{IKY26} to \cite{KY25}, Inoue, Koshikawa, and Yao have discussed various equivalent definitions for this notion in more detail. Therefore, the definition used here can be viewed as an adaptation of their work to the context of \cite{DLLZ23}.}.\par
Let $(X,\ca{M}_X)$ be a locally Noetherian fs log adic space. Denote by $\Loc(X_{\ket})$ the category of Kummer {\'e}tale local systems on $X$ with values in finite sets. If we fix a finite set or finite abelian group $\Lambda$, we denote by $\Lambda$-$\Loc(X_\ket)$ the category of Kummer {\'e}tale local systems on $X$ with values in $\Lambda$. Let $\zeta$ be a log geometric point on $X$ defined as in \cite[4.4.2-4.4.3]{DLLZ23}. Let $\pi_1^{\ket}(X,\zeta)$ be the Kummer {\'e}tale fundamental group at $\zeta$. The following result is a logarithmic analogue of the classical theory that {\'e}tale local systems correspond to finite {\'e}tale covers.
\begin{prop}[{\cite[Thm. 4.4.15-Cor. 4.4.18]{DLLZ23}}]\label{prop-ket-loc-fket-cover}
Let $(X,\ca{M}_X)$ be a connected locally Noetherian fs log adic space. There is a natural equivalence of categories:
$$X_{\fket}\iso\Loc(X_{\ket})\xrightarrow{\sim} \pi_1^{\ket}(X,\zeta)\text{-}\mrm{Fsets}.$$
Let $\Lambda$ be a finite discrete abelian group. Then the equivalence above restricts to an equivalence
$$\Lambda\text{-}\Loc(X_{\ket})\xrightarrow{\sim}\mrm{Rep}_{\Lambda}(\pi_1^{\ket}(X,\zeta)).$$
The category $\pi_1^{\ket}(X,\zeta)\text{-}\mrm{Fsets}$ (resp. $\mrm{Rep}_{\Lambda}(\pi_1^{\ket}(X,\zeta))$) consists of finite discrete sets with continuous $\pi_1^{\ket}(X,\zeta)$-actions (resp. continuous $\pi_1^{\ket}(X,\zeta)$-representations of $\Lambda$).
\end{prop}
\begin{proof}
Note that the first statement is exactly \emph{loc. cit.}. For any local system $L\in \Lambda$-$\Loc(X_{\ket})$, the corresponding finite Kummer {\'e}tale cover $Y$ is Kummer {\'e}tale locally represented by $X^\Lambda$, and the corresponding $\pi_1^{\ket}(X,\zeta)$-finite set is given by $\Hom_X(\zeta,Y)$. So the fundamental group $\pi_1^{\ket}(X,\zeta)$ acts through $\mrm{Aut}(\Lambda)$. 
\end{proof}
The log geometric point $\zeta$ is constructed from a geometric point $\xi=\Spec l$ of $X$. By \cite[Cor. 4.4.22]{DLLZ23}, we have that $\pi_1^{\ket}(\xi,\zeta)\iso \Hom (\overline{M}^{\gp},\wat{\bb{Z}}'(1)(l))$, where $\wat{\bb{Z}}'(1)(l)=\varprojlim_n\mu_n(l)$, in which $\mu_n(l)$ is the group of $n$-th roots of unity in $l$; the limit runs over all $n$ invertible in $l$.\par 
Now, assume that $X$ is defined over $\Spa(\bb{Q}_p,\bb{Z}_p)$. Then, in this case, we have $\wat{\bb{Z}}'(1)(l)=\wat{\bb{Z}}(1)$.
\begin{definition}\label{def-p-neat-fin-loc-sys}
Let $L$ be a Kummer {\'e}tale local system on $X$ with values in $\Lambda$. For any log geometric point $\zeta$ constructed from a geometric point $\xi$, $L$ corresponds to a representation $\rho_L\in\mrm{Rep}_{\Lambda}(\pi_1^{\ket}(X,\zeta))$, and the action of $\pi_1^{\ket}(\xi,\zeta)$ on $L$ is given by composing this representation with the natural map $p_\xi^X:\pi_1^{\ket}(\xi,\zeta)\to \pi_1^{\ket}(X,\zeta)$. Set $\overline{M}:=\overline{\ca{M}}_{X,\xi}$. 
In this setting, we say that $L$ is \textbf{$p$-Kummer at $\xi$} if $\rho\circ p^X_\xi$ is trivial when restricted to $\Hom (\overline{M}^{\gp},\wat{\bb{Z}}^p(1))$. We say that $L$ is \textbf{$p$-Kummer} if it is $p$-Kummer at every geometric point.
\end{definition}
Next, we consider a pro-Kummer-{\'e}tale $\wat{\bb{Z}}_p$-local system $\wat{L}$ on $X$. By \cite[Lem. 6.3.3]{DLLZ23}, $\wat{L}=\varprojlim_n L_n$, where $L_n$ is a $\bb{Z}/p^n\bb{Z}$-local system on $X_{\ket}$ such that $L_{n+1}/p^nL_{n+1}\iso L_n$. Then, at a log geometric point $\zeta$ over a geometric point $\xi$, the Kummer {\'e}tale fundamental group $\pi_1^{\ket}(X,\zeta)$ has a compatible action on $\{L_{n,\zeta}\}_n$. We obtain a continuous $\bb{Z}_p$-representation $\rho_{\wat{L}}:=\varprojlim_n \rho_{L_n}$.
\begin{definition}\label{def-p-neat-zp-loc-sys} Let $\wat{L}$ be a pro-Kummer {\'e}tale $\wat{\bb{Z}}_p$-local system on $X$.
In the setting of Definition \ref{def-p-neat-fin-loc-sys}, 
we say that $\wat{L}$ is \textbf{pro-$p$-Kummer at $\xi$} if $\rho_{\wat{L}}\circ p^X_\xi$ is trivial when restricted to $\Hom (\overline{M}^{\gp},\wat{\bb{Z}}^p(1))$. We say that $\wat{L}$ is \textbf{pro-$p$-Kummer} if it is pro-$p$-Kummer at every geometric point. Denote by $\wat{\bb{Z}}_p\text{-}\Loc_p(X_{\proket})$ the category of pro-$p$-Kummer {\'e}tale $\wat{\bb{Z}}_p$-local systems on $X$.
\end{definition}
\begin{rk}\label{rk-neat-unip}
Suppose that $\wat{L}$ is torsion free. 
Note that the image of $\Hom(\overline{M}^{\gp},\wat{\bb{Z}}^p(1))$ in $\mrm{Aut}(\wat{L}_{\zeta})$ is always finite because a continuous map from the first term (which is a formed by copies of $\wat{\bb{Z}}^p$) to $\mrm{Aut}(\wat{L}_{\zeta})\iso\mrm{GL}_N(\bb{Z}_p)$ always has a finite image. If $\wat{L}$ has unipotent boundary monodromy (i.e., the image of $\rho_{\wat{L}}\circ p_\xi^X$ is contained in some $U(\bb{Z}_p)\sbst GL_N(\bb{Z}_p)$ for a unipotent subgroup $U\sbst \mrm{GL}_{N,\bb{Z}_p}$), it is always pro-$p$-Kummer.
\end{rk}
Finally, let $\ca{G}$ be the quasi-parahoric $\bb{Z}_p$-model of a connected linear algebraic group $G$ over $\bb{Q}_p$ in the sense of Definition \ref{def: quasi-parahoric for non-reductive group}. 
\begin{definition}\label{def-torsor-p-neat}
    We say a pro-Kummer {\'e}tale $\ul{\G(\bb{Z}_p)}$-torsor $\ca{E}$ is \textbf{pro-$p$-Kummer} if, for any continuous torsion-free $\G(\bb{Z}_p)$-representation $V$, the corresponding pro-Kummer {\'e}tale local system $\ul{V}:=\ca{E}\times^{\ul{\G(\bb{Z}_p)}} V$ is pro-$p$-Kummer. Denote by $\ul{\G(\bb{Z}_p)}$-$\Loc_p(X_{\proket})$ the category of pro-$p$-Kummer {\'e}tale $\ul{\G(\bb{Z}_p)}$-torsors on $X$.
\end{definition}
\subsubsection{}Assume that $(X,\ca{M}_X)$ is either an fs log scheme locally of finite type over $\bb{Q}_p$, or a locally Noetherian fs log adic space over $\Spa(\bb{Q}_p,\bb{Z}_p)$. 
Let $G$, $\mu$, $E$, and $\G$ be as in \S\ref{subsec-log-shtuka}.
\begin{prop}\label{prop-limit-loc-system}
In the case of locally Noetherian fs log adic spaces over $\Spa \bb{Q}_p$, we have a natural bi-exact tensor equivalence of categories 
\begin{equation}\label{eq-lim-locsys-1}
\wat{\bb{Z}}_p\text{-}\Loc_p(X_{\proket})\iso\twolim\limits_{(S^{\sharp},\ca{M}_{S^\sharp}^\can,f_{S^\sharp})\in ((X,\ca{M}_X)^\Diamond)^{\mrm{op}}}\wat{\bb{Z}}_p\text{-}\Loc(S^\sharp_{\proet}).
\end{equation}
By the Tannakian formalism, we have a natural equivalence of categories
\begin{equation}\label{eq-lim-locsys-2}
\ul{\G(\bb{Z}_p)}\text{-}\Loc_p(X_{\proket})\iso\twolim\limits_{(S^{\sharp},\ca{M}_{S^\sharp}^\can,f_{S^\sharp})\in ((X,\ca{M}_X)^{\Diamond})^{\mrm{op}}}\ul{\G(\bb{Z}_p)}\text{-}\Loc(S^\sharp_{\proet}).\end{equation}
In the case of locally of finite type schemes, replacing the categories on the left by $\wat{\bb{Z}}_p\text{-}\Loc_p(X^\ad_{\proket})$ and $\ul{\G(\bb{Z}_p)}\text{-}\Loc_p(X^\ad_{\proket})$, we have the same results.
\end{prop}
\begin{lem}\label{lem-maxproket-perfectoid}
Let $(U,\ca{M}_U)$ be an affinoid Tate fs log adic space over $\Spa(\bb{Q}_p,\bb{Z}_p)$ that admits a global sharp fs chart $\md{P}$.
\begin{enumerate}
\item\label{proket-perfectoid-1} There is an inverse limit $\varprojlim_i(U_i,\ca{M}_i)$ of affinoid finite Kummer {\'e}tale covers $(U_i,\ca{M}_i)$ of $(U,\ca{M}_U)$ that is associated with an affinoid perfectoid log adic space $(\wat{U}_\infty,\wat{\ca{M}}_\infty)$, with $\overline{\wat{\ca{M}}}_\infty$ being uniquely divisible. Moreover, $\wat{\ca{M}}_{\infty}$ is perfectoid.
\item\label{proket-perfectoid-2}
Let $\md{P}[1/p]:=\varinjlim_{n}\frac{1}{p^n}\md{P}$. 
One can define another perfectoid cover $(\wat{U}_{p^{-\infty}},\wat{\ca{M}}_{p^{-\infty}})$ in $U_{\proket}$.
In fact, there is an inverse limit $U_{p^{-\infty}}:=\varprojlim_J \Spa(R_j,R_j^+)$ of affinoid finite {\'e}tale covers of $U\times_{U\langle\md{P}\rangle}U\langle\md{P}[1/p]\rangle$, which is associated with a log perfectoid space $\wat{U}_{p^{-\infty}}$ equipped with the log structure $\wat{\ca{M}}_{p^{-\infty}}$ induced by $\md{P}[1/p]$.
\end{enumerate}
\end{lem} 
\begin{proof}
We prove Part \ref{proket-perfectoid-1}. 
This is essentially \cite[Prop. 5.3.12]{DLLZ23} and \cite[Lem. 15.3]{scholze2017etale}. The argument in \cite[Prop. 5.3.12]{DLLZ23} works for both two parts; let us explain more.\par 
In fact, by \cite[Prop. 4.8]{Sch13} and \cite[Lem. 15.3]{scholze2017etale}, for any $U=\Spa(R,R^+)$ as in the Lemma, there is an affinoid perfectoid object $\varprojlim_{I}U_i=\Spa(R_i,R_i^+)$ in $U_{\proet}$ where $R_i$ are finite {\'e}tale algebras of $R$, and $\varprojlim_I U_i\sim \wat{U}$, that is, is associated with a perfectoid space $\wat{U}$. 
Then, by the argument in \cite[Prop. 5.3.12]{DLLZ23}, the inverse limits 
$\varprojlim_{(i,n)\in I\times \bb{N}}U_i\langle \frac{1}{n}\md{P}\rangle$
and $\varprojlim_{(i,n)\in I\times \bb{N}}U_i\langle \frac{1}{p^n}\md{P}\rangle$
are associated with affinoid log perfectoid spaces $\wdtd{U}_\infty$ and $\wdtd{U}_{p^{-\infty}}$.
Since $U\to U\langle\md{P}\rangle$ is a strict closed immersion, the pullback of $\wdtd{U}_\infty\to U\langle\md{P}\rangle$ (resp. $\wdtd{U}_{p^{-\infty}}\to U\langle\md{P}\rangle$) via $U\to U\langle \md{P}\rangle$ is an affinoid perfectoid log adic space $\wat{U}_\infty$ (resp. $\wat{U}_{p^{-\infty}}$) in $U_{\proket}$. The log structures are as desired by construction.
\end{proof}
\begin{lem}[{cf. \cite[Prop. 5]{IKY26}}]\label{lem-locp-split-cover} Let $(U,\ca{M}_U)$ be a Noetherian affinoid fs log adic space over $\Spa(\bb{Q}_p,\bb{Z}_p)$ that admits a global sharp fs chart $\md{P}$.\par
Let $\wat{L}\in \wat{\bb{Z}}_p\text{-}\Loc(U_{\proket})$. Then $\wat{L}$ is in $\wat{\bb{Z}}_p\text{-}\Loc_p(U_{\proket})$ if and only if $\wat{L}$ splits over a strict pro-finite {\'e}tale cover $\wat{V}_{p^{-\infty}}$ of $\wat{U}_{p^{-\infty}}$ in the notation of Lemma \ref{lem-maxproket-perfectoid}(2), such that $\wat{V}_{p^{-\infty}}\to \wat{U}_{p^{-\infty}}$ is associated with a pro-finite {\'e}tale $V_{p^{-\infty}}\to U_{p^{-\infty}}$ in $U_{\proket}$. In fact, in this case, $\wat{L}$ splits over $\wat{V}_{p^{-\infty}}$.
\end{lem}
\begin{proof}
Without loss of generality, assume that $U$ is connected. 
Since $\md{P}$ is fs and sharp, for any geometric point $\xi\to U$ of $U$, there is a surjective map $\md{P}\to \overline{\ca{M}}_{U,\xi}$ between finitely generated, torsion-free and sharp monoids. Then this induces a surjective map $\pi:\md{P}[1/p]^{\gp}/\md{P}^{\gp}\to \overline{\ca{M}}^{\gp}_{U,\xi,p^{-\infty}}/\overline{\ca{M}}_{U,\xi}^\gp$. Note that $\overline{\ca{M}}^{\gp}_{U,\xi,p^{-\infty}}/\overline{\ca{M}}_{U,\xi}^\gp\iso \overline{\ca{M}}_{U,\xi}^\gp\otimes \bb{Q}_p/\bb{Z}_p\iso \Hom(\overline{\ca{M}}_{U,\xi}^\gp,\bb{Z}_p(1))$ is a quotient of $\pi_1^{\ket}(\xi,\zeta)\iso \Hom(\overline{\ca{M}}^{\gp}_{U,\xi},\wat{\bb{Z}}(1))$.
Furthermore, we have the following commutative diagram
\begin{equation}\label{diag-ket-funda}
    \begin{tikzcd}
        \md{P}^{\gp}_\infty/\md{P}^{\gp}\arrow[r,two heads, "\pi_\infty"]\arrow[d,two heads]&\overline{\ca{M}}^{\gp}_{U,\xi,\infty}/\overline{\ca{M}}^{\gp}_{U,\xi}\arrow[d,two heads]\\
        \md{P}[1/p]^{\gp}/\md{P}^{\gp}\arrow[r,two heads,"\pi"]&\overline{\ca{M}}^{\gp}_{U,\xi,p^{-\infty}}/\overline{\ca{M}}_{U,\xi}^\gp.
    \end{tikzcd}
\end{equation}

Take a log geometric point $\zeta$ over $\xi$. By \cite[Cor. 4.4.18]{DLLZ23}, the fiber $\wat{L}_\zeta$ of $\wat{L}$ at $\zeta$ is a continuous representation of $\pi_1^{\ket}(U,\zeta)$. There is an exact sequence 
$$1\to \pi_1^{\ket}(\xi,\zeta)\to \pi_1^{\ket}(U,\zeta)\to \pi_1^{\et}(U,\xi)\to 1,$$
where $\pi_1^{\et}(U,\xi)$ corresponds to the inverse limit of all finite {\'e}tale covers of $U$ via the equivalence of categories (4.4.19) in \cite[Cor. 4.4.18]{DLLZ23}. Furthermore, the quotient of $\pi_1^{\ket}(U,\zeta)$ by the subgroup $\Hom(\overline{\ca{M}}^{\gp}_{U,\xi},\wat{\bb{Z}}^p(1))\sbst \pi_1^{\ket}(\xi,\zeta)$ corresponds to the cover $\xi_{p^{-\infty}}:=U(\xi)\times_{U(\xi)\langle \overline{\ca{M}}_{U,\xi}\rangle}U(\xi)\langle \overline{\ca{M}}_{U,\xi,p^{-\infty}}\rangle$.\par
Then $\wat{L}$ splits over $\wat{V}_{p^{-\infty}}$ if and only if $\wat{L}|_\xi$ splits over $\xi\times_U V_{p^{-\infty}}$ for any $\xi$. 
Since $\xi_{p^{-\infty}}\to \xi$ factors through $\xi\times_U V_{p^{-\infty}}$, we have the ``$\impliedby$'' direction.\par
For the other direction, 
the action of $\pi_1^{\ket}(\xi,\zeta)$ factors through the image of $\pi$. By (\ref{diag-ket-funda}), the action of $\md{P}^{\gp}_\infty/\md{P}^{\gp}_{p^{-\infty}}$ on $\wat{L}|_{U(\xi)\times_{U(\xi)\langle \md{P}\rangle}U(\xi)\langle \md{P}_{\infty}\rangle}$ is trivial by the pro-$p$-Kummer {\'e}tale condition at $(\xi,\zeta)$. Thus, $\wat{L}|_{U(\xi)}$ becomes trivial over $U(\xi)\times_{U(\xi)\langle \md{P}\rangle}U(\xi)\langle \md{P}[1/p]\rangle $ by \cite[Prop. 4.4.9]{DLLZ23}. By writing $\wat{L}$ as an inverse limit of $\bb{L}_n$, we find a pro-finite {\'e}tale cover $V$ of $U$ such that $\wat{L}$ is trivialized over $V\times_{V\langle\md{P}\rangle}V\langle\md{P}[1/p]\rangle$. Hence, the pullback of $\wat{L}$ to the perfectoid cover $\wat{U}_{p^{-\infty}}$ is trivialized over a pro-finite {\'e}tale cover of $\wat{U}_{p^{-\infty}}$.\par
The last sentence can be seen from the proof.
\end{proof}
\begin{proofof}[{Proposition \ref{prop-limit-loc-system}}]We only need to show the first equivalence (\ref{eq-lim-locsys-1}) as the bi-exactness is clear in the construction. 
Without loss of generality, we assume that we are in the case of log adic spaces since $(X^\ad,\nu^{\ad,*}\ca{M}_X)$ is locally Noetherian and fs by construction in the case of log schemes. We also assume the connectedness of $X$.\par 
Since both sides satisfy {\'e}tale descent, we assume that $X$ is affinoid and that $\ca{M}_X$ admits a global chart $\md{P}$.\par
The RHS of (\ref{eq-lim-locsys-1}) is equivalent to $\twolim\limits_{(S^{\sharp},\ca{M}_{S^\sharp},f_{S^\sharp})\in ((X,\ca{M}_X)^\Diamond)^{\mrm{op}}}\wat{\bb{Z}}_p\text{-}\Loc(S^\sharp_{\proet})$. Let $\{L_{(S^\sharp,\ca{M}_{S^\sharp},f)}\}$ be an object in this limit. This limit gives an object $L_{p^{-\infty}}:=L_{(\wat{X}_{p^{-\infty}},\wat{\ca{M}}_{p^{-\infty}})}$ by evaluating on the log perfectoid space with a perfectoid log structure $(\wat{X}_{p^{-\infty}},\wat{\ca{M}}_{p^{-\infty}})$ formed by Lemma \ref{lem-maxproket-perfectoid}(2). After replacing $\wat{X}_{p^{-\infty}}$ with a pro-finite {\'e}tale cover of it, we can assume that $L_{p^{-\infty}}$ is trivialized. 
The perfectoid space $\wat{X}_{p^{-\infty}}$ is, by construction, associated with a perfectoid object $X_{p^{-\infty}}$ in $X_{\proket}$. 
Note that $X_{p^{-\infty}}\to X$ is a pro-Kummer {\'e}tale cover by construction. 
Then the Kummer {\'e}tale fundamental group of $(X,\ca{M}_X)$ acts on $X_{p^{-\infty}}$, on its associated perfectoid space, and on the trivial local system $L_{p^{-\infty}}$ (which is encoded in the $2$-limit). 
This induces a pro-Kummer {\'e}tale $\wat{\bb{Z}}_p$-local system $\wat{L}$ in $\wat{\bb{Z}}_p\text{-}\Loc(X_{\proket})$. Furthermore, $\wat{L}$ is in $\wat{\bb{Z}}_p\text{-}\Loc_p(X_{\proket})$ by Lemma \ref{lem-locp-split-cover}.\par
Let us show the other direction. By Lemma \ref{lem-locp-split-cover}, an object $\wat{L}\in \wat{\bb{Z}}_p\text{-}\Loc_p(X_{\proket})$ is trivialized over a perfectoid object $X_{p^{-\infty}}$ in $X_{\proket}$. Moreover, $\wat{L}$ is equipped with an action of $\pi_1^{\ket}(X,\zeta)$ that factors through $\mrm{Aut}(X_{p^{-\infty}}/X)$. By applying \cite[Lem. 5.3.8]{DLLZ23} and taking an inverse limit, the saturated product $X_{p^{-\infty}}\times_X^{sat}(X\times_{X\langle\md{P}\rangle}X\langle\md{P}[1/p]\rangle)$ is pro-finite {\'e}tale over $X\times_{X\langle\md{P}\rangle}X\langle\md{P}[1/p]\rangle$. Since, for every object $(S^\sharp,\ca{M}^\can_{S^\sharp},f)\in (X,\ca{M}_X)^\Diamond$, the morphism $f:S^\sharp\to X$ factors through $X\times_{X\langle\md{P}\rangle}X\langle\md{P}[1/p]\rangle$, the saturated product $S^\sharp_{p^{-\infty}}:=X_{p^{-\infty}}\times_X S^\sharp$ is pro-finite {\'e}tale and surjective over $S^\sharp$. The $\pi_1^{\ket}(X,\zeta)$-action on $\wat{L}$ is pulled back to the action of $\mrm{Aut}(S^\sharp_{p^{-\infty}}/S^\sharp)$ on $(f^{-1}\wat{L})(S^\sharp_{p^{-\infty}})$. Hence, $f^{-1}\wat{L}$ is a pro-{\'e}tale local system on $S^\sharp$. For any morphism $g:(S_1^\sharp,\ca{M}_{S^\sharp_1},f_1)\to (S_2^\sharp,\ca{M}_{S^\sharp_2},f_2)$, it is easy to check that the pullbacks satisfy the desired functoriality required by a $2$-limit.
\end{proofof}
\subsubsection{}\label{subsubsec-equi-cat-gen}
In this subsection, we show the following theorem:
\begin{thm}[{cf. \cite[Prop. 2.5.3]{PR24}}]\label{thm-equi-cat-gen}
Let $(X,\ca{M}_X)$ be a locally Noetherian fs log adic space over $\Spa(\bb{Q}_p,\bb{Z}_p)$ (resp. an fs log scheme of finite type over $\Spec\bb{Q}_p$). 
With this assumption, there is a natural equivalence between the following two categories:
  \begin{enumerate}
      \item The category of log shtukas $\Sht_{\ca{G},\mu}^\Diamond(X,\ca{M}_X)$;
      \item The category $HT_{\G,\mu}(X,\ca{M}_X)$ of pairs $(\bb{P},\pi_{HT})$, where $\bb{P}$ is a pro-$p$-Kummer-{\'e}tale $\ul{\ca{G}(\bb{Z}_p)}$-torsor defined over $X$ (resp. $X^\ad$) and a Hodge-Tate period map $\pi_{HT}: \bb{P}\to \ca{F}_{G,\mu^{-1}}^\Diamond$. 
  \end{enumerate}
The meaning of the last morphism is the following: The torsor $\bb{P}$ can be viewed as a limit over $(X,\ca{M}_X)^\Diamond$ of pro-{\'e}tale $\wat{\bb{Z}}_p$-torsors $\ca{E}_{S^\sharp}$ by Proposition \ref{prop-limit-loc-system}. The last morphism is given by a functorial assignment of $\G(\bb{Z}_p)$-equivariant maps between $v$-sheaves $HT(\ca{E}_{S^\sharp}):\ca{E}_{S^\sharp}\to \mrm{Gr}_{G,\Spd E,\mu^{-1}}$ for objects in $(X,\ca{M}_X)^\Diamond$ followed by the Bialynicki-Birula isomorphism (see Lemma \ref{lem: BB map}) $\mrm{Gr}_{G,\Spd E,\mu^{-1}}\to \ca{F}^\Diamond_{G,\mu^{-1}}$.
\end{thm}
\begin{rk}\label{rk-comp-KoshikawaYao}
See \cite[Thm. 7]{IKY26} (and also \cite[Thm. 7.36]{KY25}) for an analogous theorem in the log prismatic theory. In fact, our proof of Theorem \ref{thm-equi-cat-gen} is also similar to the one presented there.
\end{rk}
\begin{rk}\label{rk-adic-generic}
The following statement is included in the case of log adic spaces in Theorem \ref{thm-equi-cat-gen}: 
Let $X$ be a separated $\bb{Z}_p$-scheme that is flat and of finite type. Let $\wat{X}$ be the $p$-adic completion and $Y:=\wat{X}^\ad_\eta$ be the adic generic fiber. If $Y$ is equipped with an fs log structure $\ca{M}$, then applying the theorem above gives an equivalence over $(Y,\ca{M})^\Dia$.
\end{rk}
\begin{proofof}[{Theorem \ref{thm-equi-cat-gen}}]
By definition, we have $\Sht^\Diamond_{\G,\mu}(X,\ca{M}_X)=$ $$\twolim\limits_{(S^{\sharp},\ca{M}_{S^\sharp},f_{S^\sharp})\in ((X,\ca{M}_X)^\Diamond)^{\mrm{op}}}\Sht_{\ca{G},\mu}(S^\sharp).$$
By \cite[Prop. 2.5.2]{PR24}, this expression is equivalent to  $$\twolim\limits_{(S^{\sharp},\ca{M}_{S^\sharp},f_{S^\sharp})\in ((X,\ca{M}_X)^{\Dia})^{\mrm{op}}}HT_{\G,\mu}(S^\sharp).$$
The last $2$-limit is equivalent to $\twolim\limits_{(S^{\sharp},\ca{M}_{S^\sharp}^\can,f_{S^\sharp})\in ((X,\ca{M}_X)^{\Dia})^{\mrm{op}}}HT_{\G,\mu}(S^\sharp)$ by Lemma \ref{lem-can-log-str} and Lemma \ref{lem-equiv-subcat-sec}. Finally, combining Proposition \ref{prop-limit-loc-system} (and the explanation in Theorem \ref{thm-equi-cat-gen}), we have the desired equivalence.
\end{proofof}

\begin{cor}\label{cor-equi-cat-gen}
For an fs log scheme $(X,\ca{M}_X)$ over $\bb{Q}_p$ in \S\ref{subsubsec-log-diamonds-setup} Case \ref{case-gen-log-sch}, we have an equivalence between $\Sht^\Diamond_{\G,\mu}(X,\ca{M}_X)$ and $HT_{\G,\mu}^\Diamond(X,\ca{M}_X)$, the category of pairs $(\bb{P},\pi_{HT})$ where $\bb{P}$ is an object in $$\twolim\limits_{(S^{\sharp},\ca{M}_{S^\sharp},f_{S^\sharp})\in ((X,\ca{M}_X)^{\Diamond})^{\mrm{op}}}\ul{\G(\bb{Z}_p)}\text{-}\Loc(S^\sharp_{\proet}),$$
and $\pi_{HT}$ is as above.
\end{cor}
\begin{proof}
Write as $2$-limits and apply \cite[Prop. 2.5.2]{PR24}. 
\end{proof}
\begin{cor}[Canonical extension of shtukas in characteristic zero]\label{cor-canonical-extensions}
\begin{enumerate}

\item\label{can-ext-cor-2} Let $E$ be a $p$-adic field. Assume that $X$ is a smooth rigid analytic variety over $\Spa E$ with a normal crossings divisor $D$ such that $U:=X\bss D$ is open and dense in $X$. Denote by $\ca{M}_X$ the log structure induced by $D\hookrightarrow X\xhookleftarrow{j} U$. Then the restriction functor
$$\Res^X_U: \Sht^\Diamond_{\G,\mu}(X,\ca{M}_X)\to \Sht^\Diamond_{\G,\mu}(U)$$
is fully faithful.
\item\label{can-ext-cor-3} With the same assumptions as in Part \ref{can-ext-cor-2}, for any object $\mathscr{P}_U\in \Sht_{\G,\mu}^\Diamond(U)$, $\mathscr{P}_U$ corresponds to a $\ul{\G(\bb{Z}_p)}$-torsor $\bb{P}_U$ and a Hodge-Tate map $\pi_{HT}: \bb{P}_U\to \ca{F}_{G,\mu^{-1}}^\Diamond$. 
Then $\mathscr{P}_U$ is in the essential image of $\Res^X_U$ \textbf{only if} $\bb{P}:=j_*\bb{P}_U$ is pro-$p$-Kummer {\'e}tale. 
\item\label{can-ext-cor-1} (Diao-Lan-Liu-Zhu) With the same assumptions as in Parts \ref{can-ext-cor-2} and \ref{can-ext-cor-3}, we further assume that $\bb{P}_U$ is \textbf{de Rham} in the sense that, for any $V\in \mrm{Rep}(\ca{G}(\bb{Z}_p))$, $\ul{V}_{\bb{Q}_p}:=(\bb{P}_U\times^{\ul{\ca{G}(\bb{Z}_p)}}V)\otimes_{\bb{Z}_p}\bb{Q}_p$ is de Rham. 
If $\ul{V}_{\bb{Q}_p}$ has \textbf{unipotent} geometric monodromy along the boundary in the sense of \cite[Def. 6.3.7]{DLLZ23}, then there is a canonical way of associating a Hodge-Tate map $\pi_{HT,X}:\bb{P}\to \ca{F}_{G,\mu^{-1}}^\Dia$ to $\bb{P}$, and therefore defining a log shtuka $\mathscr{P}\in\Sht_{\G,\mu}^\Dia(X,\ca{M}_X)$ via Theorem \ref{thm-equi-cat-gen}. This association is compatible with the one in \cite[Prop. 2.6.3]{PR24} after the restriction back to $U^\Dia$.
\end{enumerate}
\end{cor}
\begin{proof}
   Consider the diagram 
   \begin{equation}
       \begin{tikzcd}
        \Sht^\Diamond_{\G,\mu}(X,\ca{M}_X)\arrow[rr]\arrow[d]&& \Sht^\Diamond_{\G,\mu}(U)\arrow[d]\\
        HT_{\G,\mu}(X,\ca{M}_X)\arrow[rr]&&HT_{\G,\mu}(U).
       \end{tikzcd}
   \end{equation}
This diagram is commutative. The vertical maps are equivalences of categories by Theorem \ref{thm-equi-cat-gen} and \cite[Prop. 2.5.3]{PR24}. Hence, for Part \ref{can-ext-cor-2}, it suffices to show that the bottom arrow is fully faithful.\par
By rigid Abhyankar's Lemma \cite[Prop. 4.2.1]{DLLZ23} and by \cite[Cor. 6.3.4]{DLLZ23}, there is an equivalence 
\begin{equation}\label{eq-equi-loc-sys-res}\wat{\bb{Z}}_p\text{-}\Loc(X_{\proket})\iso \wat{\bb{Z}}_p\text{-}\Loc(U_{\proet})\end{equation}
induced by $\wat{L}\mapsto \wat{L}|_{U_{\proet}}$; the quasi-inverse is defined by the pushforward of $\wat{L}|_{U_{\proet}}$ from $U$ to $X$. By the Tannakian formalism, there is an equivalence
$$\Res^{\G-\Loc}:\ul{\G(\bb{Z}_p)}\text{-}\Loc(X_{\proket})\iso \ul{\G(\bb{Z}_p)}\text{-}\Loc(U_{\proet})$$
induced by restricting the pro-Kummer {\'e}tale $\ul{\G(\bb{Z}_p)}$-torsor to $U$.
Since the functor $\ul{\G(\bb{Z}_p)}$-$\Loc_p(X_{\proket})\to \ul{\G(\bb{Z}_p)}$-$\Loc(X_{\proket})$ is fully faithful, there is a fully faithful functor induced by restriction to $U$:
$$i:\ul{\G(\bb{Z}_p)}\text{-}\Loc_p(X_{\proket})\to\ul{\G(\bb{Z}_p)}\text{-}\Loc(U_{\proet}).$$
The rest of Part \ref{can-ext-cor-2} follows from Lemma \ref{lem-ext-ht} below.\par
The quasi-inverse of (\ref{eq-equi-loc-sys-res}) is given by $\wat{L}\mapsto j_*\wat{L}$. By the construction above, $j_*\bb{P}_U$ in the statement of Part \ref{can-ext-cor-3} is pro-$p$-Kummer if and only if $\bb{P}_U$ is in the essential image of $\ul{\G(\bb{Z}_p)}\text{-}\Loc_p(X_{\proket})$ under the composition $i\circ \Res^{\G-\Loc}$.
So Part \ref{can-ext-cor-3} follows.\par
Part \ref{can-ext-cor-1} follows from the main theorems of \cite{DLLZ} (see also \cite[Thm. 4.2.1]{RCam26}).
In fact, by the proof of \cite[Thm. 3.2.12]{DLLZ} written above \cite[\S 3.5]{DLLZ}, there is a canonical isomorphism
$$\mu^* D_{\dr,\log}(\ul{V}){\otimes}_{\ca{O}_{X_{\proket}}}\ca{O}\bb{B}_{\dr,\log}\xrightarrow{\sim}\ul{V}\otimes_{\wat{\bb{Z}}_p}\ca{O}\bb{B}_{\dr,\log},$$
where $\mu: X_{\proket}\to X_{an}$, $\ca{O}\bb{B}_{\dr,\log}$ is a peroid sheaf with filtration and log connection extending $\ca{O}\bb{B}_{\dr}$, and $D_{\dr,\log}$ is the arithmetic log Riemann-Hilbert functor (see \cite[2.2 and 3.2.6]{DLLZ}). This isomorphism is compatible with filtrations and log connections on both sides.\par
We then define $\bb{M}_{0,\log}:=(D_{\dr,\log}(\ul{V})\otimes_{\ca{O}_X}\ca{O}\bb{B}_{\dr,\log}^+)^{\nabla^{\log}=0}$, and $\bb{M}_{\log}:=(\mrm{Fil}^0(D_{\dr,\log}(\ul{V})\otimes_{\ca{O}_X}\ca{O}\bb{B}_{\dr,\log}))^{\nabla^{\log}=0}$; they are $B_{\dR}^+$-modules due to log Poincar{\'e}'s lemma \cite[Cor. 2.4.2]{DLLZ}. Moreover, the filtrations on them are finite projective modules by \cite[Cor. 3.4.22]{DLLZ}. Then, using the same argument as in \cite[\S 2.6.1]{PR24}, we obtain a Hodge-Tate map $\pi_{HT,X}:\bb{P}\to \mrm{Gr}_{G,\mu^{-1}}$. It is clear that this construction is compatible with and extends the one in \emph{loc. cit.}.
\end{proof}
\begin{lem}\label{lem-ext-ht}
With the conventions and assumptions as in Corollary \ref{cor-canonical-extensions}, every Hodge-Tate map $\pi_{HT}: \bb{P}_U\to \ca{F}^\Diamond_{G,\mu^{-1}}$ on $U^\Diamond$ for a fixed $\ul{\G(\bb{Z}_p)}$-torsor $\bb{P}_U$ admits \textbf{at most one} extension to $\pi_{HT,X}:\bb{P}:=j_*\bb{P}_U\to \ca{F}^\Diamond_{G,\mu^{-1}}$ on $(X,\ca{M}_X)^{\Diamond}$.
\end{lem}
\begin{proof}Denote by $\pi_1$ and $\pi_2$ two extensions of $\pi_{HT}$. As in the proof of \cite[Prop. 4.2.1]{DLLZ23}, let us assume that both $X=\Spa(R,R^+)$ and $U=\Spa(R_0,R_0^+)$ are affinoid.   
By rigid Abhyankar \cite[Prop. 4.2.1]{DLLZ23} and by Lemma \ref{lem-maxproket-perfectoid}(1), we can choose an affinoid perfectoid object $\wdtd{U}:=\Spa(A,A^+)$ such that $\wdtd{U}\to U$ extends to an affinoid perfectoid object $\wdtd{X}:=\Spa(B,B^+)\to X$ in $X_{\proket}$. We then take the associated perfectoid spaces $\wat{X}\sim \wdtd{X}$ and $\wat{U}\sim \wdtd{U}$, and take a strictly totally disconnected cover $\wat{X}^{str}\to \wat{X}\to X$ where $\wat{X}^{str} \to \wat{X}$ is universally open by \cite[Lem. 7.18]{scholze2017etale}. Let $\wat{U}^{str}:=\wat{X}^{str}\times_{\wat{X}}\wat{U}$; this is also a strictly totally disconnected space, since it is an open subspace of $\wat{X}^{str}$. Note that the log structure on $\wat{X}^{str}$ is induced by the pullback of the log structure on $\wat{X}$. Thus, $\wat{X}^{str}$ with its log structure is an object in $(X,\ca{M}_X)^\Diamond$.\par 
By assumption, $|U|$ is open dense in $|X|$. In fact, $|\wat{U}|$ is dense in $|\wat{X}|$. 
Indeed, since $|\wat{X}|\iso\varprojlim_{j\in J} |X_j|$ by \cite[Prop. 6.4]{scholze2017etale}, it suffices to show that each $U\times_X X_j$ is open dense in $X_j$. By \cite[Prop. 2.8]{Han20}, and since the constructions in \cite{Han20} and \cite[Prop. 4.2.1]{DLLZ23} are compatible, the pullback of $D$ through $X_j\to X$ is nowhere-dense. This implies that $U\times_X X_j$ is dense in $X_j$. Furthermore, we deduce from this fact that $\wat{U}^{str}$ is also dense in $\wat{X}^{str}$ since $\wat{X}^{str}\to \wat{X}$ is universally open.\par
If there are two maps $\pi_1(\wat{X}^{str})$ and $\pi_2(\wat{X}^{str}):\bb{P}_{\wat{X}^{str}}\to \ca{F}_{G,\mu^{-1},\wat{X}^{str}}$ extending $\pi_{HT,\wat{U}^{str}}$, since $|\wat{U}^{str}|$ is dense in $|\wat{X}^{str}|$ and $\ca{F}_{G,\mu^{-1}}$ is separated over $E$, we have $\pi_1=\pi_2$. \par
Pick an object $(S^\sharp, \ca{M}_{S^\sharp},f)\in (X,\ca{M}_X)^\Diamond$. 
Let $\pi_1(S^\sharp)$ and $\pi_2(S^\sharp):\bb{P}_{S^\sharp}\to \ca{F}_{G,\mu^{-1},S^\sharp}$ be two Hodge-Tate maps. We now pull them back to the perfectoid space $\wat{S}^\sharp_{p^{-\infty}}$ associated with the saturated product $S^\sharp\times^{sat}_{X}\wdtd{X}$, and to $\wat{S}_1^\sharp:=\wat{S}^\sharp_{p^{-\infty}}\times_{\wat{X}}\wat{X}^{str}$. By functoriality of $2$-limits, the two pullbacks $\pi_1(\wat{S}^\sharp_1)$ and $\pi_2(\wat{S}_1^\sharp)$ are equal. For the same reason, the $v$-descent data of the two maps from $\wat{S}^\sharp_1$ to $S^\sharp$ are canonically identified. Hence, we have $\pi_1(S^\sharp)=\pi_2(S^\sharp)$ by descent.
\end{proof}
\subsection{Extending log shtukas}\label{subsec-ext-logsh}
Let $(X,\ca{M}_X)$ be an fs log scheme, where $X$ is a normal scheme that is separated, flat and of finite type over $\bb{Z}_p$ or $X=\Spec A$ for an excellent Noetherian normal domain $A$ that is flat over $\bb{Z}_p$. 
Denote by $X_{\bb{Q}_p}$ the generic fiber of $X$. 
Let $\ca{M}_{X_{\bb{Q}_p}}=i^*\ca{M}_{X}$ be the log structure induced by the natural inclusion $i:X_{\bb{Q}_p}\hookrightarrow X$.\par
We show that
\begin{thm}[{cf. \cite[Thm. 2.7.7]{PR24}}]\label{thm-ext-shu-gen}
In the situation above, the restriction functor
$$\Res:=\Res^{X}_{X_{\bb{Q}_p}}: \Sht^{\Diamond/}_{\ca{G},\mu}(X,\ca{M}_X)\lra\Sht^\Diamond_{\ca{G},\mu}(X_{\bb{Q}_p},\ca{M}_{X_{\bb{Q}_p}})$$
is fully faithful.\par
Moreover, the restriction 
$$\Res: \Sht^{\diamond}_{\ca{G},\mu}(X,\ca{M}_X)\lra\Sht_{\ca{G},\mu}((X,\ca{M}_{X})^\dia\times_{\Spd\bb{Z}_p}\Spd\bb{Q}_p)$$
is fully faithful.
\end{thm}
\subsubsection{}\label{subsubsec-ext-logsh-proof}
Let us introduce some conventions.\par
Let $\mathscr{P}_1$ and $\mathscr{P}_2$ be two objects in $\Sht^{\Diamond/}_{\G,\mu}(X,\ca{M}_X)$. It follows from the construction in Definition \ref{def-log-shtuka} and Definition \ref{def-2-lim} that $\Hom_{\Sht^{\Diamond/}_{\G,\mu}(X,\ca{M}_X)}(\mathscr{P}_1,\mathscr{P}_2)=$
$$\lim\limits_{(S^{\sharp},\ca{M}_{S^\sharp}^\can,f_{S^\sharp})\in ((X,\ca{M}_X)^{\Diamond/})^{\mrm{op}}}\Hom_{\Sht_{\G,\mu}(S^\sharp)}(\mathscr{P}_{1,(S^\sharp,f)},\mathscr{P}_{2,(S^\sharp,f)}).$$
The morphisms in $\Sht^\Diamond_{\G,\mu}(X_{\bb{Q}_p},\ca{M}_{X_{\bb{Q}_p}})$ are formed in a similar way. \par
Writing $\overline{\mathscr{P}}_i:=\mathscr{P}_{i}|_{X_{\bb{Q}_p}^{\log\Diamond}}$ the restriction of $\mathscr{P}_i$ to $X_{\bb{Q}_p}^{\log\Diamond}$ for $i=1,2$, a morphism $\mbf{H}\in \ob \Hom(\overline{\mathscr{P}}_1,\overline{\mathscr{P}}_2)$ is then represented by $\mbf{H}=\{\mbf{H}_{(S^\sharp,\ca{M}^{\can}_{S^\sharp},f)}\}$, where $\mbf{H}_{(S^\sharp,\ca{M}^{\can}_{S^\sharp},f)}\in \Hom_{\Sht_{\G,\mu}(S^\sharp)}(\mathscr{P}_{1,(S^\sharp,f)},\mathscr{P}_{2,(S^\sharp,f)})$ and the collection runs over the objects $(S^\sharp,\ca{M}_{S^\sharp}^{\can},f)$ in $X_{\bb{Q}_p}^{\log\Diamond}$ that satisfy the condition in Definition \ref{def-2-lim}(2).\par
Similarly, let $\ca{H}\in \Hom_{\Sht^{\Diamond/}_{\G,\mu}(X,\ca{M}_X)}(\mathscr{P}_1,\mathscr{P}_2)$ be a morphism. Denote by $\ca{H}_{(S^\sharp,\ca{M}^{\can}_{S^\sharp},f)}$ the morphism between $\mathscr{P}_{1,(S^\sharp,\ca{M}^{\can}_{S^\sharp},f)}$ and $\mathscr{P}_{2,(S^\sharp,\ca{M}_{S^\sharp}^{\can},f)}$ assigned to $\ca{H}$. \par
In what follows, we will frequently use the abbreviations $\ca{H}_{(S^\sharp,\ca{M}^{\can}_{S^\sharp},f)}=\ca{H}_{S^\sharp}$ and $\mbf{H}_{(S^\sharp,\ca{M}^{\can}_{S^\sharp},f)}=\mbf{H}_{S^\sharp}$ if there is nothing confusing.\par
\begin{proofof}[{Theorem \ref{thm-ext-shu-gen}}] The proof can be adapted from the proofs of \cite[Prop. 2.7.6]{PR24} and \cite[Thm. 2.7.7]{PR24}. First, note that the argument in \emph{loc. cit.} did not use that the scheme is of finite type after restricting the question locally to a normal domain; moreover, the second claim will follow from the same proof as the first. We assume that we are in the first situation for $X$.\par
With the conventions as above, given $\mbf{H}\in \Hom(\overline{\mathscr{P}}_1,\overline{\mathscr{P}}_2)$, we will show that $\mbf{H}$ extends uniquely to a morphism $\ca{H}\in \Hom_{\Sht^{\Diamond/}_{\G,\mu}(X,\ca{M}_X)}(\mathscr{P}_1,\mathscr{P}_2)$.\par
Without loss of generality, we assume that $X=\Spec A^+$ is affine and admits a global fs sharp chart $\md{P}$. (Indeed, we can do this by taking an {\'e}tale neighborhood around each geometric point of $X$ and by taking the log structure on the {\'e}tale neighborhood the pullback of $\ca{M}_X$. After proving this case, one can pass to a Cech nerve of an {\'e}tale cover.) We write $\wat{A}^+$ the ($p$-adic) completion and write $\wat{A}:=\wat{A}^+[\frac{1}{p}]$. Write $\wat{X}:=\mrm{Spf}\wat{A}^+$ and $\wat{X}^\ad:=\Spa(\wat{A}^+,\wat{A}^+)$.\par
Let $f:(Y:=\Spa(B,B^+),\ca{M}_Y)\to (\wat{X}^{\ad},\ca{M}_{\wat{X}^\ad})$ be an object in $(X,\ca{M}_{X})^{\diamond}$. Here, we also assume that $Y=\Spa(B,B^+)$ is affinoid perfectoid and assume that $\ca{M}_Y=\ca{M}_Y^{\can}$.\par
Firstly, we consider the case where $(B,B^+)=(C,C^+)$ is an algebraically closed, (nonarchimedean, complete) perfectoid field. As in the proof of \cite[Prop. 2.7.6]{PR24}, if the morphism $f$ is adic, then $C$ is of characteristic zero, but, in this case, we have that $f$ is already an object in $(X_{\bb{Q}_p},\ca{M}_{X_{\bb{Q}_p}})^{\Diamond}$, so we set $\ca{H}_{\Spa(C,C^+)}=\mbf{H}_{\Spa(C,C^+)}$.\par
If the morphism $f$ is not adic, then $C$ is of characteristic $p$. 
In this case, we claim that $f:\Spa(C,C^+)\to \wat{X}^\ad$ factors through an integral perfectoid space $X_\infty$ over $\wat{X}^\ad$. In fact, let $\wat{X}_0:=\Spa (\wat{A},\wat{A}^+)$. Then $\wat{X}_0$ is an affinoid Tate adic space equipped with a morphism to the adic generic fiber of $\wat{X}^\ad$ and to $\wat{X}^\ad$, $i:\wat{X}_0\to \wat{X}^\ad_\eta\to \wat{X}^\ad$.
Let $R_{p^{-\infty}}$ be the inverse limit of finite {\'e}tale algebras over $\wat{A}_{p^{-\infty}}:=\wat{A}\otimes_{\wat{A}\langle\md{P}\rangle}\wat{A}\langle \md{P}[1/p]\rangle$. 
Let $R^+_{p^{-\infty}}$ be the integral closure in $R_{p^{-\infty}}$ of $\wat{A}^+$. By Lemma \ref{lem-maxproket-perfectoid}(2), the $p$-adic completion $(\wat{R}_{p^{-\infty}},\wat{R}^+_{p^{-\infty}})$ is a complete affinoid perfectoid Huber pair. Now, since $(C,C^+)$ is an algebraically closed perfectoid field, $C^+$ is a valuation ring and $R^+_{p^{-\infty}}$ is formed by an inverse limit of normalizations in finite ring extensions, the map $\wat{A}^+\to C^+$ factors through $\wat{A}^+\to \wat{R}^+_{p^{-\infty}}\to C^+$. Set $X_\infty:=\Spa(\wat{R}^+_{p^{-\infty}})$. The claim is shown.\par
By Lemma \ref{lem-maxproket-perfectoid}(2), $X_\infty$ is equipped with the log structure $\md{P}[1/p]^a$ induced by $\md{P}[1/p]$; this is a saturated and fine perfectoid log structure by construction (see also \cite[Rmk. 2.26]{KY25}). So the associated log perfectoid space $X_{\infty,\eta}:=\Spa(\wat{R}_{p^{-\infty}},\wat{R}^+_{p^{-\infty}})$ with its natural morphism to $\wat{X}^\ad$ is an object in $X^{\log \diamond}$. By \cite[Prop. 2.7.6]{PR24}, $\mbf{H}_{X_{\infty,\eta}}$ extends uniquely to a morphism $\ca{H}_{\Spd(\wat{R}^+_{p^{-\infty}})}:\mathscr{P}_1|_{\Spd(\wat{R}^+_{p^{-\infty}})}\to\mathscr{P}_2|_{\Spd(\wat{R}^+_{p^{-\infty}})}$, and $\ca{H}_{\Spa C}$ is determined by pulling back $\ca{H}_{\Spd(\wat{R}^+_{p^{-\infty}})}$ to $\Spd C$. (To apply the proposition there, one needs that $\wat{R}^+_{p^{-\infty}}=\wat{R}^\circ_{p^{-\infty}}$. This follows from the normalization construction.) \par
Finally, we consider general $Y=\Spa(B,B^+)$. This follows from a $v$-descent argument. 
As in the proof of \emph{loc. cit.}, we choose a collection of maps $\{\pi_i:\Spa(C_i,C_i^+)\to \Spa(B,B^+)\}$ covering $|Y|$, where $\Spa(C_i,C_i^+)$ are algebraically closed perfectoid fields. We then form a product $Z:=\Spa(\prod_iC_i,\prod_iC_i^+)$ with a $v$-cover $\pi:=\prod_i\pi_i: Z\to Y$.
Then $Z$ maps to $\wat{X}^\ad$ via $Z\xrightarrow{\ \pi\ }Y\xrightarrow{\ f\ }\wat{X}^\ad$, and the log structure on $Z$ is defined by pulling back $\ca{M}_Y$. 
We can construct $\ca{H}_Z$: By the proof of \cite[Prop. 2.7.6]{PR24}, the morphisms on $W (C_i^+)$ for the corresponding BKF modules. We then take the product of these morphisms on $\prod_iW(C_i^+)$ and restrict to get a morphism $\ca{H}_Z$ between shtukas.\par
To show that the morphism $\ca{H}_{Z}$ descends to a unique morphism $\ca{H}_Y$, it suffices to check that there is a descent datum $p^*_1\ca{H}_Z\iso p^*_2\ca{H}_Z$ on $\wdtd{Z}:=Z\times_Y Z$. 
As explained in the proof of \emph{loc. cit.}, it suffices to check this pointwisely.
For any point $s:\Spa(C,C^+)\to \wdtd{Z}$, we have $\pi\circ p_1\circ s=\pi\circ p_2\circ s$.
So the two sides are the equal after post-composing $f$. It suffices to deal with the case that $C$ is of characteristic $p$. By the last three paragraphs, we know that $f\circ \pi\circ p_i\circ s$ can also be factored as $x:\Spa C\to \Spd C^+\to \Spd \wat{R}^+_{p^{-\infty}}\to X^{\log \diamond}$, and this factorization is \emph{unique} by the first paragraph of the proof of \emph{loc. cit.} (see \cite[2.7.2]{PR24} and \cite{gleason2025specialization}). Thus the last two paragraphs gives a unique way of defining $\ca{H}_{\Spa C,x}$, as desired. \par
We now check the uniqueness (as there might be a different choice of $v$-cover) and functoriality of this assignment. Given a morphism $g: (Y,\ca{M}_{Y}^\can,f)\to (Y',\ca{M}_{Y'}^\can,f')$ such that $f'\circ g=f$, we check that, taking products of points $Z\to Y$ and $Z'\to Y'$ and constructing $\ca{H}_{Y}$ and $\ca{H}_{Y'}$ using them respectively, we have $g^*\ca{H}_{Y'}=\ca{H}_Y$. 
Upon replacing $Z$ with the disjoint union of it with a $v$-cover of $Y\times_{Y'}Z'$, we assume that there is a morphism $\wdtd{g}:Z\to Z'$ covering $g$. Then, by the same argument as the last paragraph, the $\ca{H}_{\Spa C,x}$ at any point of $Z$ constructed by either pulling back from $Z$ or pulling back from $Z\xrightarrow{\wdtd{g}} Z'$ are the same as the one factoring through $\Spd \wat{R}^+_{p^{-\infty}}$. So $\wdtd{g}^*\ca{H}_{Z'}=\ca{H}_Z$. Next, we consider the fiber product $\wdtd{Z}$ as the last paragraph, and the same argument shows that the morphisms descend and are equal on $Y$.\par
For the case of affine excellent normal flat domains, all arguments above work; note that the normalization argument goes through by \cite[part2, 7.8.3(5)]{EGA4}.
\end{proofof}
\subsubsection{}\label{subsubsec-ext-logsh-cor}
Let us state an immediate corollary of Theorem \ref{thm-ext-shu-gen}.
\begin{cor}[{cf. \cite[Rmk. 2.7.9]{PR24}}]\label{cor-ext-shu-gen}
Let $X$ be a normal scheme that is flat, separated and of finite type over $\bb{Z}_p$. Let $D\sbst X$ be a relative Cartier divisor such that $U:=X\bss D$ is open dense in $X$. 
Let $\ca{M}_X$ be the log structure defined by $D\hookrightarrow X\hookleftarrow U$.
\begin{enumerate}
\item Suppose that the generic fiber $X_{\bb{Q}_p}$ is a smooth variety with the normal crossings divisor $D_{\bb{Q}_p}$. Then the restriction functor 
    $$\Res^X_U:\Sht^{\Diamond/}_{\G,\mu}(X,\ca{M}_X)\to \Sht^{\Diamond/}_{\G,\mu}(U)$$ is fully faithful.
\item Suppose that $X_{\bb{Q}_p}$ is a smooth variety, and that $\wat{U}^\ad_\eta\hookrightarrow \wat{X}^\ad_\eta$ has a normal crossings complement. Then the restriction functor
$$\Res^X_U:\Sht^{\diamond}_{\G,\mu}(X,\ca{M}_X)\to \Sht_{\G,\mu}^\dia(U)$$
is fully faithful.
\end{enumerate}
\end{cor}
\begin{proof}We have the following commutative diagram
\begin{equation*}
    \begin{tikzcd}
    \Sht^{\Diamond/}_{\G,\mu}(X,\ca{M}_X)\arrow[r,"\Res^X_U"]\arrow[d] &\Sht^{\Diamond/}_{\G,\mu}(U)\arrow[d]\\
    \Sht^{\Diamond}_{\G,\mu}(X_{\bb{Q}_p},\ca{M}_{X_{\bb{Q}_p}})\arrow[r,"\Res^{X_{\bb{Q}_p}}_{U_{\bb{Q}_p}}"]& \Sht^{\Diamond}_{\G,\mu}(U_{\bb{Q}_p}).
    \end{tikzcd}
\end{equation*}
Then the assertion follows from combining Theorem \ref{thm-ext-shu-gen}, \cite[Thm. 2.7.7]{PR24} and Corollary \ref{cor-canonical-extensions}. Indeed, we have known that the functors represented by the bottom arrow and the vertical arrows are fully faithful.\par
The second part follows similarly from Remark \ref{rk-adic-generic}, Corollary \ref{cor-canonical-extensions} and Theorem \ref{thm-ext-shu-gen}. 
\end{proof}

\subsection{Restriction to special fiber}\label{subsec-res-special}
Suppose that $(X,\ca{M}_X)$ is in Case \ref{case-gen-log-sch} of log schemes. 
\subsubsection{}Assume that $X$ is defined over $\bb{Z}_p$ or just $\bb{F}_p$.
\begin{prop}\label{prop-perfect-nonlog}
The projection $F^{\log}:(X_{\F_p},\ca{M}_{X_{\F_p}})^\diamond\to X_{\F_p}^\diamond$ admits a canonical section $i^{\log}$. 
This section is given by assigning to any $(S^\sharp,f)\in X_{\F_p}^\diamond(S)$ an object in log diamond $(S^\sharp,(f^*\ca{M}_{X_{\F_p}})^{\perf},f)\in X_{\F_p}^{\log \diamond}.$ 
The same holds after replacing $\diamond$ with $\Diamond$.
\end{prop}
\begin{proof}
    When $S^\sharp$ is an affinoid perfectoid space of characteristic $p$, it is perfect. Then the morphism between sheaves of monoids $f^*\ca{M}_{X_{\F_p}}\to \ca{O}_{S^\sharp_\et}$ induces a morphism $(f^*\ca{M}_{X_{\F_p}})^{\perf}\to\ca{O}_{S^\sharp_\et}$. 
\end{proof}
\begin{cor}\label{cor-special-fiber-nonlog}
Suppose that there is a log $\G$-shtuka bounded by $\mu$, 
$$\mathscr{P}:(X,\ca{M}_X)^{\diamond}\to \Sht_{\G,\mu}.$$
Then the restriction of $\mathscr{P}$ to $(X_{\F_p},\ca{M}_{X_{\F_p}})^\diamond$ gives a morphism 
$$\mathscr{P}|_{X_{\F_p}}:(X_{\F_p},\ca{M}_{X_{\F_p}})^\diamond\to \Sht_{\G,\mu}.$$
This morphism is equivalent to a morphism 
$$\mathscr{P}|_{X_{\F_p}}:X_{\F_p}^\diamond\to \Sht_{\G,\mu}.$$
The same holds if we replace $\diamond$ with $\Diamond$.\par
Moreover, there is a morphism $\mathscr{P}^{\red}:X_{\bb{F}_p}^{\perf}\to \Sht^W_{\G,\mu}$.
\end{cor}
\begin{proof}
The construction of $\ca{M}_{S^{\sharp}}^{\can}$ in Lemma \ref{lem-can-log-str} coincides with the perfection of $f^*\ca{M}_X$ in Proposition \ref{prop-perfect-nonlog}. By Proposition \ref{prop-shu-can-obj}, we have a morphism $X_{\F_p}^\diamond\to \Sht_{\GG, \mu}$. Applying the reduction functor in the sense of \cite{gleason2025specialization}, and noting that $(\Sht_{\GG, \mu})_{\red}$ is represented by $\Sht_{\GG, \mu}^W$, $X_{\F_p}^\diamond = (X_{\F_p}^{\perf})^\diamond$ by \cite[Prop. 18.3.1]{SW20} 
and that $((X_{\F_p}^{\perf})^\diamond)_{\red}$ is represented by $X_{\F_p}^{\perf}$ by \cite[Prop. 3.16]{gleason2025specialization}, we get the desired morphism.
\end{proof}
For simplicity, we will use charts to refer to log structures in the remaining part of \S\ref{subsec-res-special}.
\begin{rk}\label{rk: specialization map of log shtukas}
    Let $K$ be a finite field extension of $\rQ_p$, and let $\varpi \in K$ be a uniformizer. One can consider the absolute log prismatic site $(\OO_K, \varpi^{\N})_{\Prism}$. A log prismatic $F$-crystal $\EE$ on $(\OO_K, \varpi^{\N})_{\Prism}$ is equivalent to a semistable $G_K$-representation $\mathbb{L}$ in $\Rep_{\Z_p}^{\st}(G_K)$ by \cite{yao2023mathbb}. One can read off the monodromy action of the semistable representation $\mathbb{L}[1/p]$ by evaluating the associated log prismatic $F$-crystal $\EE$ on Breuil-Kisin/Hyodo–Kato log prisms. If we evaluate $\EE$ at the cover $(\OO_K, \varpi^{\N}) \to (\OO_C, \varpi^{\N[\frac{1}{p}]})$, one obtains a Breuil-Kisin-Fargues module (see \cite[Ex. 3.1]{yao2023mathbb}). Let $S = \Spa(R, R^+) \in \Perf$. The objects $(S^{\sharp}, \ca{M}_{S^{\sharp}}, f) \in (\OO_K, \varpi^{\N})^{\diamond}(S)$ require $S^{\sharp}$ to be perfectoid and $\ca{M}_{S^{\sharp}}$ to be saturated and fine perfectoid; thus $\overline{\ca{M}}_{S^{\sharp}}$ is uniquely $p$-divisible (see Lemma \ref{lem-satfineperf-satfineperfbar}). In particular, under the equivalence of categories between perfect log prisms and perfectoid log rings in \cite[Prop. 2.39]{KY25}, one only sees the \emph{perfect} objects in $(\OO_K, \varpi^{\N})^{\diamond}$. Therefore, the realization functor from log prismatic $F$-crystals to log shtukas forgets the monodromy action when we consider the specialization map as in \cite[Prop. 2.4.6]{PR24}.
\end{rk}

\subsubsection{Log $p$-divisible groups}

In this subsection, we give an example that might be relevant to Corollary \ref{cor-special-fiber-nonlog}. Proposition \ref{prop: log p-divisible group and p-divisible group} will not be used in the remaining of this paper. 

Let $K$ be a finite field extension of $\rQ_p$, let $\OO_K$ be its ring of integers, and let $k$ be its residue field. Let $\varpi \in K$ be a uniformizer, and set $S = (\Spec \OO_K, \varpi^{\N})$, $S_0 = (\Spec k, 0^{\N})$. By the main theorem of \cite{bertapelle2023log}, taking the generic fiber $(\cdot)_K$ gives an equivalence of categories between the category $\mathbf{BT}^{\log}_{S, d}$ of dual-representable log $p$-divisible groups over $S$ and the category $\mathbf{BT}_K^{\st}$ of $p$-divisible groups over $K$ with semistable reduction. Moreover, taking the Tate module $T_p$ gives an equivalence of categories between $\mathbf{BT}_K^{\st}$ and the category $\Rep_{\Z_p}^{\st, \lrbracket{0, 1}}(G_K)$ of semistable $G_K$-representations with Hodge--Tate weights in $\lrbracket{0, 1}$.

Let $\HH \in \mathbf{BT}^{\log}_{S, d}$. On the one hand, as in Remark \ref{rk: specialization map of log shtukas}, one can attach to it a log prismatic $F$-crystal and then a log shtuka in $(\OO_K, \varpi^{\N})^{\dia}$; by restricting to the special fiber and by Corollary \ref{cor-special-fiber-nonlog}, we obtain a usual shtuka over $\Spec k$. On the other hand, the special fiber $\HH_k$ of $\HH$ is a log $p$-divisible group over $S_0$, which need not be a classical $p$-divisible group. Nevertheless, one can still produce a usual shtuka from $\HH_k$ that is compatible with the above one, by Lemma \ref{lem-can-log-str} and the following result:

\begin{prop}\label{prop: log p-divisible group and p-divisible group}
    Let $\HH$ be a dualizable log $p$-divisible group over $(\Spec k, 0^{\N})$. Then its pullback to $(\Spec k, 0^{\N[\frac{1}{p}]})$ becomes a classical $p$-divisible group.
\end{prop}
\begin{proof}
    This follows from Kato's classification of log $p$-divisible groups \cite[Thm. 3.1]{Kat23}; cf. \cite[Thm. 3.8, Cor. 3.10]{wurthen2024log}. We briefly recall some notation and results from Kato's work on log $p$-divisible groups, following \cite[\S 3]{wurthen2024log}. Let $T = \Spec A$, where $A$ is a Noetherian henselian local ring with the residue characteristic $p$, and suppose that $T$ admits a global chart $P \to \ca{M}_T$ such that the induced map $\md{P} \to \ca{M}_{T,\overline{t}}/\OO_{T, \overline{t}}^{\times}$ is an isomorphism at the geometric point $\overline{t}$ of $T$. Let $(\mathrm{fin}/T)_d$ be the category of finite Kummer flat log group schemes over $T$ defined as in \cite[1.6]{Kat23}; see also \cite[Def. 2.3, 2.6]{wurthen2024log} for details. For $\GG \in (\mathrm{fin}/T)_d$, there is a unique exact sequence $0 \to \GG^{\circ} \to \GG \to \GG^{\et} \to 0$, which restricts to the classical connected-\'etale sequence over any finite Kummer flat cover where $\GG$ becomes a classical $p$-divisible group. Moreover, $\GG^{\circ}$ and $\GG^{\et}$ are classical finite flat group schemes. From Kato's classification, any $\GG \in (\mathrm{fin}/T)_d$ that is $p$-power torsion determines, and is uniquely determined by, a pair $(\GG^{\mathrm{cl}}, \beta)$, where $\GG^{\mathrm{cl}}$ is a classical extension of $\GG^{\circ}$ and $\GG^{\et}$, and $\beta \in \Hom(\GG^{\et}(1), \GG^{\circ}) \otimes_{\Z} \md{P}^{\gp}$.
    
    By definition, $\HH = \varinjlim \HH_n$, where $\HH_n = \Ker(\HH \stackrel{p^n}{\to} \HH)$ is an object in $(\mathrm{fin}/S_0)_d$. Therefore, $\HH_n$ determines and is determined by $(\HH_n^{\mathrm{cl}}, \beta_n)$. Since $\Hom(\HH_n^{\et}(1), \HH_n^{\circ})$ is killed by $p^n$, $\beta \in \Hom(\HH_n^{\et}(1), \HH_n^{\circ}) \otimes_{\Z} \md{P}^{\gp}/p^n\md{P}^{\gp}$. In particular, after passing to $S_{0, \infty} = (\Spec k, 0^{\N[\frac{1}{p}]})$, $\beta$ is trivialized and $\HH_n$ becomes a classical $p$-divisible group. Therefore, $\HH$ is classical over $(\Spec k, 0^{\N[\frac{1}{p}]})$.
\end{proof}

\section{Local systems and shtukas on mixed Shimura varieties}\label{sec-shtukas-mixed-sh}

In this section, we work with the theory in characteristic $0$.

\subsection{Group-theoretic lemmas}\label{subsec: pass from P to P^c}
We collect the results used later regarding various group-theoretic constructions in compactification theory, Bruhat-Tits theory, and their intersection.
\subsubsection{}\label{subsubsec-gen-mix-sh-gp-lem}
Let $P$ be any linear algebraic group. Define $P^c := P/Z(P)_{ac}$, where $Z(P)_{ac}$ is the anti-cuspidal part of the multiplicative connected center $Z(P)^\circ$ of $P$.\par 
Given a mixed Shimura datum $(P, \ca{X})$, let $G$ be the Levi quotient of $P$.
Then the natural projection $P \to G$ induces an embedding $Z(P)\hookrightarrow Z(G)$, and $Z(P)_{ac}$ is sent to $Z(G)_{ac}$ by definition. Note that $P^c \to G^c$ factors through the Levi quotient $G^{\prime, c}$ of $P^c$, and $G^{\prime, c} \to G^c$ is a homomorphism with central kernel.
\begin{lem}\label{lem-levi-c}
For a mixed Shimura datum $(P,\ca{X})$, we have $G^{\prime,c}=G^c$.
\end{lem}
\begin{proof}
Recall that $(P,\ca{X})$ satisfies \cite[2.1 (viii)]{Pin89}. That is, we require that, for any Levi subgroup $\wdtd{G}\sbst P$ lifting $G$, the adjoint action of the center $Z(\wdtd{G})$ of $\wdtd{G}$ on $\lie W$ factors through a cuspidal quotient; $W$ denotes the unipotent radical of $P$. Note that the definition of $Z(\wdtd{G})$ is independent of the choice of lifting $\wdtd{G}$.
From this condition, we see that 
$$Z(G)_{ac}\xrightarrow{\sim} Z(\wdtd{G})_{ac}\sbst \mrm{Cent}_{Z(\wdtd{G})}(\lie W)= Z(P).$$
So the embedding $Z(P)_{ac}\hookrightarrow Z(G)_{ac}$ admits a section by the displayed expression above. We then have $Z(P)_{ac}=Z(G)_{ac}$, and the lemma follows.
\end{proof}
\subsubsection{}\label{subsubsec-bd-mx-sh-gp-lem}
Now let $(G,X)$ be a Shimura datum. 
Let $\Phi=(Q_{\Phi}, X_{\Phi}^+, g_{\Phi})$ be a cusp label representative of $(G, X)$ (see \cite[2.1.7]{Mad19}). Recall that one associates to $\Phi$ a mixed Shimura datum $(P_\Phi,D_\Phi)$. Recall that $Q_\Phi\sbst G$ is an admissible $\bb{Q}$-parabolic that contains a normal subgroup $P_\Phi$ and is equipped with a $Q_\Phi(\R)$-equivariant morphism
\begin{equation*}
   \tau: X \to \pi_0(X) \times \Hom(\DS_{\CC}, P_{\Phi,\CC}),\quad x\mapsto ([x], u_x^{Q_\Phi}\circ h_{\infty}).
\end{equation*}
Moreover, $D_\Phi$ is defined as the $P_\Phi(\R)U_\Phi(\CC)$-orbit of a fixed point $([x], u_x^{Q_\Phi}\circ h_{\infty})$ such that $x \in X^+_\Phi$, where $X_\Phi^+$ is a connected component of $X$; it depends only on $X^+_\Phi$, not on the choice of $x$ in $X^+_\Phi$. Denote by $W_\Phi$ the unipotent radical of $P_\Phi$ and $U_\Phi$ the center of $W_\Phi$. Let $K$ be an open compact subgroup of $G(\A)$. Define $K_{\Phi} = P_\Phi(\A) \cap g_{\Phi}Kg_{\Phi}^{-1}$; define $K_{\Phi,p}$ and $K_\Phi^p$ similarly.
 \begin{lem}\label{lem-ptog-emb}
     The natural embedding $P_\Phi \to G$ induces a finite map $P^c_\Phi \to G^c$. Let $ZP_\Phi:=(Z_G\cdot P_\Phi)^\circ$. Then there is an injective map $ZP_\Phi^c\to G^c$. 
 \end{lem}
 \begin{proof}
Since the center of $ZP_\Phi/ZW_\Phi$ is isogenous to a product of split tori and $\bb{R}$-anisotropic tori (cf. \cite[Cor. 4.10]{Pin89} and \cite[Lem. 1.14]{Wu25}), we know that $Z(P_\Phi)_{ac}\sbst Z(ZP_\Phi)_{ac}=Z(ZW_\Phi)_{ac}=Z(G)_{ac}$. So there is a well-defined quotient map $P_\Phi^c=P_\Phi/Z(P_\Phi)_{ac}\to G^c=G/Z(G)_{ac}$.
For $ZP_\Phi$, by the last paragraph, we see that the embedding $ZP_\Phi\hookrightarrow G$ induces an embedding $ZP^c_\Phi\hookrightarrow G^c$. 
 \end{proof}

 \begin{lem}\label{lem: quasi-parahoric}
   Let $\GG$ be a quasi-parahoric group scheme of $G$, and let $\GG_{\Phi}$ be the $g_{\Phi}$-conjugate of $\GG$ (i.e. $\GG_{\Phi}(\bZ_p) = g_{\Phi}\GG(\bZ_p)g_{\Phi}^{-1}$). Let $\QQ_{\Phi}$ (resp. $\mathcal{W}_{\Phi}$) be the closure of $Q_{\Phi}$ (resp. $W_{\Phi}$) in $\GG_{\Phi}$, and let $\PP_{\Phi}$ be the smoothing of the closure of $P_{\Phi}$ in $\GG_{\Phi}$.
   \begin{enumerate}
       \item When $\GG$ is stabilizer quasi-parahoric (resp. quasi-parahoric, parahoric), then $\QQ_{\Phi}$ is stabilizer quasi-parahoric (resp. quasi-parahoric, parahoric) in the sense of \ref{def: quasi-parahoric for non-reductive group}. Moreover, when $\GG$ is parahoric, then $\QQ_{\Phi} \hookrightarrow \GG_{\Phi}$ is a parabolic embedding in the sense of \ref{def: parabolic embedding}.
       \item When $\GG$ is stabilizer quasi-parahoric (resp. quasi-parahoric, parahoric), then $\PP_{\Phi}$ is stabilizer quasi-parahoric (resp. quasi-parahoric, quasi-parahoric). Moreover, $(\PP_{\Phi}, \mu_{\Phi})$ comes from boundary in the sense of \ref{def: PP, mu comes from boundary}.
   \end{enumerate}

   Let $\GG$ be quasi-parahoric. Let $\UU_{\Phi}$ be the closure of $U_{\Phi}$ in $\mathcal{W}_{\Phi}$, and let $\mathcal{V}_{\Phi}$ be the quotient $\mathcal{W}_{\Phi}/\UU_{\Phi}$. Then $\UU_{\Phi}$ and $\mathcal{V}_{\Phi}$ are affine smooth group schemes with connected fibers, and we have the decomposition $\mathcal{W}_{\Phi} = \UU_{\Phi} \rtimes \mathcal{V}_{\Phi}$.
 \end{lem}
 \begin{proof}
     This is essentially \cite[Prop. 2.67]{Mao25b}. We start from $\GG_{\Phi}$ and omit the index $\Phi$ for simplicity (we won't mention the initial $\GG$). First of all, assume $K_p:=K_{\Phi, p}$ and $Q:=Q_{\Phi}$ are in good position, i.e. the associated (generic) point $x \in B_{\red}(G, \rQ_p)$ of the facet $\FF \subset B_{\red}(G, \rQ_p)$ ($\GG^{\circ} = \GG_{\FF}$) is contained in an apartment $A_{\red}(G, T)$, where $T$ is a maximal $\bQ$-split torus of $G$ defined over $\rQ_p$ such that $T \subset L \subset Q$. By the arguments below Definition \ref{def: quasi-parahoric for non-reductive group}, the arguments above Definition \ref{def: PP, mu comes from boundary}, and by \cite[Prop. 2.67]{Mao25b} with the help of Lemma \ref{lem: quasi-parahoric, intermediate}, we see that both statements $(1)$ and $(2)$ are true. In this case, $\LL$ (resp. $\GG_{h}$) is the smoothening of the closure of $L$ (resp. $G_h$) in $\GG$.

     In general, we first show $(1)$: apply arguments in \cite[Rmk. 2.45]{Mao25b}: we can always find an element $g \in G(\rQ_p)$ such that $K' := g^{-1}Kg$ and $Q$ are in good position, $K = gK'g^{-1}$. We factor $g = qnk$, $q \in Q(\rQ_p)$, $n \in N(\rQ_p)$, $k \in P_x^0$, where $N$ is the normalizer of $T$ over $\rQ_p$, and $P_x^0 = \GG_x^{\circ}(\Z_p) = (K')^{\circ}$. The action $k$ fixes $x$, $n$ moves $x$ to $nx \in A_{\red}(G, T)$, and $q$ shifts the apartment $A_{\red}(G, T)$ in $B_{\red}(G, \rQ_p)$. Let $K_1 = nkK'(nk)^{-1}$, then $K_1$ and $Q$ are in good position, $Q(\rQ_p) \cap K_1 = (W(\rQ_p) \cap K_1) \rtimes (L(\rQ_p) \cap K_1)$ by the first paragraph. Since $W$ is normal in $Q$,
     \[ Q(\rQ_p) \cap K = q(Q(\rQ_p) \cap K_1)q^{-1} = (W(\rQ_p) \cap K) \rtimes (qL(\rQ_p)q^{-1} \cap K),  \]
     we still have $\QQ = \mathcal{W}\rtimes \LL$, and $\mathcal{W}$ is smooth with connected fibers since it is conjugated to the one cut out by $K_1$. But in this case $\LL$ is no longer the closure of $L$ in $\GG$, we need to conjugate the section $L \to Q$ by $q$. When we study the levels of mixed Shimura varieties coming from boundary, this is exactly what we do: we define the level $K_{\Phi, L} \subset L_{\Phi}(\A)$ as the image of $K_{\Phi, Q} \subset Q_{\Phi}(\A)$ under the projection $Q_{\Phi} \to L_{\Phi}$, and we do not prescribe (and do not need) a section $L_{\Phi} \to Q_{\Phi}$. 
     
     For the second statement $(2)$, in the arguments above Definition \ref{def: PP, mu comes from boundary}, we use $K_1$ in place of $K$, and note that the conjugation by $q$ shifts apartments but does not change the correspondence between $x_h$ and $x_L$.

     For the last statement, when $K_p$ and $Q$ are in good position, fix a splitting $W = U \times V$ where both $U$ and $V$ are products of root groups, by Bruhat-Tits theory, we have $\mathcal{W} = \mathcal{U} \times \mathcal{V}$, where both $\mathcal{U}$ and $\mathcal{V}$ are the closures of $U$ and $V$ in $\mathcal{W}$ respectively. $\mathcal{W}$, $\mathcal{U}$ and $\mathcal{V}$ are affine smooth group schemes with connected fibers. In general, since $U$ is the center of $W$, then $Q$ normalizes $U$, as in second paragraph one can conjugate $1 \to \mathcal{U} \to \mathcal{W} \to \mathcal{V} \to 1$ by elements in $Q$ and still get the wanted exact sequence.
 \end{proof}
 \begin{lem}\label{lem: quasi-parahoric, intermediate}
     Let $x \in A_{\red}(G, T)$, $T \subset Q \subset G$, $\KK^{\circ} = \GG_x^{\circ}(\bZ_p)$, $\KK = \GG_x(\bZ_p)$. Let $\KK_1$ be a quasi-parahoric group controlled by $x$, i.e., $\KK^{\circ} \subset \KK_1 \subset \KK$. Then
     \[ Q(\bQ) \cap \KK_1 = (W(\bQ) \cap \KK_1) \rtimes (L(\bQ) \cap \KK_1). \]
 \end{lem}
 \begin{proof}
     We omit the notation $(\bQ)$. Since
     \[ Q \cap \KK^{\circ} = (W \cap \KK^{\circ}) \rtimes (L \cap \KK^{\circ}), \quad Q \cap \KK = (W \cap \KK) \rtimes (L \cap \KK), \]
     and $W \cap \KK^{\circ} = W \cap \KK$ (by the triviality of $\pi_1(W)$), we have $W \cap \KK_1 = W \cap \KK^{\circ}$. On the other hand,
     \[ (Q \cap \KK_1)/(Q \cap \KK^{\circ}) = \pi(\KK_1)/\pi(\KK^{\circ}), \quad (Q \cap \KK)/(Q \cap \KK_1) = \pi(\KK)/\pi(\KK_1). \]
     Here $\pi: Q \to L$ is the projection. Since $\pi(\KK)/\pi(\KK^{\circ}) = (L \cap \KK)/(L \cap \KK^{\circ}) \subset \pi_1(L)_{I, \mathrm{tor}}$ is finite abelian, there is a unique subgroup of $Q \cap \KK$ containing $Q \cap \KK^{\circ}$ with image $\pi(\KK_1)/\pi(\KK^{\circ})$ in $\pi(\KK)/\pi(\KK^{\circ})$, this forces $Q \cap \KK_1 = (W \cap \KK_1) \rtimes (L \cap \KK_1)$ and $L \cap \KK_1 = \pi(\KK_1)$.
 \end{proof}
 This proof also shows the following.
 \begin{lem}\label{lem: bound of quasi-parahoric}
     Let $\PP_x = \UU \rtimes \GG_x$ be a stabilizer quasi-parahoric group scheme, and $\PP_x^{\circ} = \UU \rtimes \GG_x^{\circ}$. Let $\KK = \PP_x(\bZ_p)$, $\KK_G = \GG_x(\bZ_p)$, $\KK^{\circ} = \PP_x^{\circ}(\bZ_p)$, $\KK_G^{\circ} = \GG_x^{\circ}(\bZ_p)$. Given a $\sigma$-invariant closed subgroup $\KK_1 \subset \KK$ that contains $\KK^{\circ}$, let $\KK_{1, G} \subset G(\bZ_p)$ be the image of $\KK_1$, then $\KK_{1, G}$ and $\KK_1$ are quasi-parahoric and there exist unique quasi-parahoric models $\PP_1 = \UU \rtimes \GG_1$ with $\PP_1(\bZ_p) = \KK_1$, $\GG_1(\bZ_p) = \KK_{1, G}$.
 \end{lem}
 
 \begin{lem}\label{lem: induced compatibility on quasi-group scheme}
     Let $P_1 = U_1 \rtimes G_1$, $P_2 = U_2 \rtimes G_2$, $P_1 \to P_2$ be an embedding that is compatible with $U_1 \rtimes G_1 \to U_2 \rtimes G_2$ in the sense of \ref{def: compatible, P and G}. Let $\PP_i = \UU_i \rtimes \GG_i$ be quasi-parahoric group schemes of $P_i$ such that 
     \[  \PP_1(\bZ_p) = \PP_2(\bZ_p) \cap \PP_1(\bQ),\quad \GG_1(\bZ_p) = \GG_2(\bZ_p) \cap \GG_1(\bQ),    \]
     then the induced morphism $\PP_1 \to \PP_2$ is compatible with $\GG_1 \to \GG_2$ in the sense of \ref{def: compatible, PP and GG}.
 \end{lem}
 \begin{proof}
     The decomposition $\PP_2 = \UU_2 \rtimes \GG_2$ fixes a section $G_2 \to P_2$. In characteristic $0$, we can choose a section $G_1 \to P_1$ such that $G_1 \to P_1 \to P_2$ factors through $G_2$. By assumption,
     \[ \GG_1(\bZ_p) = G_1(\bQ) \cap \GG_2(\bZ_p) = G_1(\bQ) \cap \PP_2(\bZ_p) = G_1(\bQ) \cap \PP_1(\bZ_p), \]
     this gives a section $\GG_1 \to \PP_1$ that induces $\PP_1 = \UU_1 \rtimes \GG_1$.
 \end{proof}
 
\subsubsection{Associating quasi-parahoric group schemes}\label{subsubsec-ass-gp-sch}

Recall:
$$G(\bQ)^0 = \ker (G(\bQ) \stackrel{\Tilde{\kappa}_G}{\to} \pi_1(G)_I),\quad G(\bQ)^1 = \Ker(G(\bQ) \stackrel{\Tilde{\nu}_G}{\to} \pi_1(G)_{I} \otimes \rQ),$$ 
Given any quasi-parahoric group scheme $\GG$ over $\Z_p$, one can find $x \in \Bui_{\red}(G, \rQ_p)$ such that
\begin{equation}\label{eq: cut-off quasi-parahoric}
    \Stab_{G(\bQ)}(x) \cap G(\bQ)^0 \subset \GG(\bZ_p) \subset \Stab_{G(\bQ)}(x) \cap G(\bQ)^1,
\end{equation}
moreover, $x$ is in the generic position of the facet $\FF \subset \Bui_{\red}(G, \bQ)$ determined by the parahoric group scheme $\GG^{\circ}$. We have $\GG^{\circ}(\bZ_p) = \Stab_{G(\bQ)}(x) \cap G(\bQ)^0$.

Also, recall that, if $G \to G'$ is a surjection with a central kernel, there is a canonical ($G(\bQ)\to G'(\bQ)$)-equivariant bijection $\Bui_{\red}(G, \bQ) \to \Bui_{\red}(G', \bQ)$. Here we use reduced buildings instead of extended buildings.
\begin{lem}\label{lem: surjection on G^0}
    Let $G$ be a reductive group, $Z \subset G$ be a central torus, and $G' = G/Z$. Then the surjection $G(\bQ) \to G'(\bQ)$ induces a surjection $G(\bQ)^0 \to G'(\bQ)^0$, a surjection $\Stab_{G(\bQ)}(x) \to \Stab_{G'(\bQ)}(x')$, and a morphism $G(\bQ)^1 \to G'(\bQ)^1$ where $x'\in \Bui_{\red}(G', \rQ_p)$ is the image of some $x\in \Bui_{\red}(G, \rQ_p)$. In particular, given a quasi-parahoric (resp. parahoric) group scheme $\GG$ of $G$ associated with a point $x \in \Bui_{\red}(G, \rQ_p)$ in the sense of $(\ref{eq: cut-off quasi-parahoric})$, then the image of $\GG(\bZ_p)$ contains (resp. is) the parahoric group $(\GG^{\prime}_{x'})^{\circ}(\bZ_p)$ and is contained in the stabilizer quasi-parahoric group $\GG^{\prime}_{x'}(\bZ_p)$.
\end{lem}
\begin{proof}
    Since $Z$ is connected, $H^1(\bQ, Z)$ is trivial by Steinberg's theorem, $G(\bQ) \to G'(\bQ)$ is a surjection. Consider the commutative diagram with exact rows:
\[\begin{tikzcd}
	1 & {Z(\bQ)} & {G(\bQ)} & {G'(\bQ)} & 1 \\
	& {\pi_1(Z)_{I}} & {\pi_1(G)_{I}} & {\pi_1(G')_{I}} & {1,}
	\arrow[from=1-1, to=1-2]
	\arrow[from=1-2, to=1-3]
	\arrow["{\tilde{\kappa}_Z}", from=1-2, to=2-2]
	\arrow[from=1-3, to=1-4]
	\arrow["{\tilde{\kappa}_G}", from=1-3, to=2-3]
	\arrow[from=1-4, to=1-5]
	\arrow["{\tilde{\kappa}_{G'}}", from=1-4, to=2-4]
	\arrow[from=2-2, to=2-3]
	\arrow[from=2-3, to=2-4]
	\arrow[from=2-4, to=2-5]
\end{tikzcd}\]
here vertical maps are Kottwitz maps and are surjective (\cite[\S 7.1]{kottwitz1997isocrystals}). Chasing the diagram, we have a surjective $G(\bQ)^0 \to G'(\bQ)^0$.

Let $g \in \Stab_{G'(\bQ)}(x')$ and $g' \in G(\bQ)$ be a point in its preimage. Since $\Bui_{\red}(G, \bQ) \rightiso \Bui_{\red}(G', \bQ)$ is $G(\bQ)$-$G'(\bQ)$-equivariant, $g'$ fixes $x$; thus, $\Stab_{G(\bQ)}(x) \to \Stab_{G'(\bQ)}(x')$ is surjective and moreover, $\Stab_{G(\bQ)}(x)$ is the preimage of $\Stab_{G'(\bQ)}(x')$ under $\pi: G(\bQ) \to G'(\bQ)$. This forces
$$ \pi(\Stab_{G(\bQ)}(x) \cap G(\bQ)^0) = \pi(\Stab_{G(\bQ)}(x)) \cap \pi(G(\bQ)^0) = \Stab_{G'(\bQ)}(x') \cap G'(\bQ)^0$$.
\end{proof}

By Lemmas \ref{lem: surjection on G^0} and \ref{lem: bound of quasi-parahoric}, we make the following definitions:
\begin{Definition and Lemma}\label{def-ass-gp-sch-mix-sh}
Let $\PP = \UU \rtimes \GG$ be a quasi-parahoric group scheme of $P$, and let $\KK_p = \PP(\bZ_p)$. Then $\GG^{\circ}$ is the parahoric group scheme $\GG_x^{\circ} = \GG_{\FF}$ for some $x \in \FF \subset \Bui_{\red}(G, \rQ_p)$.
\begin{enumerate}
    \item Let $\KK^c_p \subset P^c(\bQ)$ be the image of $\KK_p$ under the projection $P \to P^c$. Then $\KK^c_p$ is a quasi-parahoric group, and there is a unique quasi-parahoric group scheme $\PP^c = \UU \rtimes \GG^c$ of $P^c$ with $\PP^c(\bZ_p) = \KK^c_p$, where $(\GG^c)^{\circ}$ is the parahoric group scheme associated with $x \in \Bui_{\red}(G^c, \rQ_p)$. We call $K_p^c := \ca{P}^c(\Z_p)$ and $\ca{P}^c$ the quasi-parahoric group (scheme) \textbf{associated with} $K_p$, and $P \to P^c$ extends to $\PP \to \PP^c$.
    \item In general, let $P \to P'$ be a quotient with central multiplicative kernel. Let $\PP^{\prime, \circ}_x = \UU \rtimes \GG^{\prime, \circ}_x$ (resp. $\PP^{\prime}_x = \UU \rtimes \GG^{\prime}_x$) be the parahoric (resp. stabilizer quasi-parahoric) group scheme associated with $x \in \Bui_{\red}(G', \rQ_p)$. Let $\Breve{K}'_p \subset \PP^{\prime}_x(\bZ_p)$ be a $\sigma$-invariant subgroup that contains $\im(\KK_p)$ and $\PP^{\prime, \circ}_x(\bZ_p)$ and that stabilizes $\UU(\bZ_p)$. Then it is a quasi-parahoric group, and there is a unique quasi-parahoric group scheme $\PP' = \UU \rtimes \GG'$ with $\PP'(\bZ_p) = \KK'_p$ and $(\GG')^{\circ} = \GG^{\prime, \circ}_x$; moreover, $P \to P'$ extends to $\PP \to \PP'$.
\end{enumerate}
\end{Definition and Lemma}
\begin{rk}\label{rk-explain-definition}
       In fact, the subgroup generated by $\im(\KK_p)$ and $\ca{P}^{\prime, \circ}(\bZ_p)$ is just $\im(\KK_p)\cdot \ca{P}^{\prime, \circ}(\bZ_p)$ and is open compact: $\im(\breve{K}_p) \subset \ca{P}_x^{\prime, \circ}(\bZ_p)$, thus its conjugation stabilizes $\ca{P}^{\prime, \circ}(\bZ_p)$. 
\end{rk}


Now let $(G,X)$ be a Shimura datum, and let $\Phi$ be a cusp label representative.
\begin{definition}\label{def-pstar}
Let $Y_\Phi$ be a connected subgroup satisfying $P_\Phi\sbst Y_\Phi\sbst ZP_\Phi$. Denote $Y_\Phi^*:=Z(G)_{ac}\cdot Y_\Phi/Z(G)_{ac}$.\par
By definition, we have a central isogeny $Y_\Phi^c\to Y_\Phi^*\sbst G^c$; in addition, $ZP_\Phi^*=ZP_\Phi^c$.
\end{definition}
Note that $P_\Phi^*$ can also be obtained as follows. Pick $\Phi\in\ca{CLR}(G,X)$. This cusp label representative maps to a cusp label representative $\Phi^*\in\ca{CLR}(G^c,X^c)$. Then $P^*_\Phi=P_{\Phi^*}$.
\begin{definition}\label{def-quasi-parahoric-big}
Let $\KK_{\Phi,p}^*=\KK_p^c\cap P_\Phi^*(\bQ)$ and let $\ca{P}^*_\Phi$ be the quasi-parahoric group scheme whose group of $\bZ_p$-points is $\KK_{\Phi,p}^*$ (see Lemma \ref{lem: quasi-parahoric}). The group $\KK_p^c = \GG^c(\bZ_p)$ is defined in \ref{def-ass-gp-sch-mix-sh}, viewing $G$ here as the $P$ there.
\end{definition}
Then there is, by definition, a map $\KK_{\Phi,p}^c\to \KK_{\Phi,p}^*$; the first group is defined by Definition \ref{def-ass-gp-sch-mix-sh} with input $(P, \KK_p)=(P_\Phi, \KK_{\Phi,p})$. We have $\ca{P}_\Phi^c\to \ca{P}_\Phi^*$.
\begin{definition}\label{def-mixed-sh-tower-*}
The mixed Shimura datum $(P_\Phi,D_\Phi)$ and the map $P_\Phi\to P_\Phi^*$ induce a mixed Shimura datum $(P_\Phi^*,D_\Phi^*)$. We denote the corresponding tower with respect to the weight filtration by
$(P_\Phi^*,D_\Phi^*)\to({\overline{P}_\Phi^*},{\overline{D}_\Phi^*}) \to (G_{\Phi,h}^*,D^*_{\Phi,h})$. Now, the images of $\KK^*_{\Phi,p}$ in $\overline{P}^*_\Phi(\breve{\bb{Q}}_p)$ and $G_{\Phi,h}^*(\breve{\bb{Q}}_p)$ are denoted by $\overline{\KK}^*_{\Phi,p}$ and $\KK^*_{\Phi,h,p}$, respectively.
\end{definition}

By Lemma \ref{lem: quasi-parahoric}, we have quasi-parahoric group schemes $\ovl{\PP}^*_{\Phi}$ and $\GG^*_{\Phi, h}$ such that
$$ \PP^*_{\Phi} = \mathcal{W}_{\Phi} \rtimes \GG^*_{\Phi, h},\quad \ovl{\PP}^*_{\Phi} = \mathcal{V}_{\Phi} \rtimes \GG^*_{\Phi, h},\quad \GG^*_{\Phi, h}(\bZ_p) = \KK^*_{\Phi,h,p},\quad \ovl{\PP}^*_{\Phi}(\bZ_p) = \overline{\KK}^*_{\Phi,p}.  $$

Moreover, the natural morphisms $\ovl{P}_{\Phi}^c \to \ovl{P}_{\Phi}^*$, $G_{\Phi, h}^c \to G_{\Phi, h}^*$ induce morphisms of quasi-parahoric group schemes $\ovl{\PP}_{\Phi}^c \to \ovl{\PP}_{\Phi}^*$, $\GG_{\Phi, h}^c \to \GG_{\Phi, h}^*$. We denote $K^*_{\Phi,h,p} = \GG^*_{\Phi, h}(\Z_p)$, $\ovl{K}^*_{\Phi,p} = \ovl{\PP}^*_{\Phi}(\Z_p)$ as usual.

\subsection{Mixed Shimura varieties}\label{subsec: betti local system}
In this subsection, we quickly recall some facts about mixed Shimura varieties. \par
Let $(P,\ca{X})$ be a mixed Shimura datum in the sense of \cite{Pin89}. According to the weight filtration, there is a tower of mixed Shimura data $(P,\ca{X})\to (\overline{P},\overline{\ca{X}})\to (G_h, X_h)$, where $G_h$ is the Levi quotient of $P$ and $\overline{P}=P/U$, where $U$ is a normal subgroup contained in the unipotent radical of $P$.\par
Let $K$ be a neat open compact subgroup of $P(\A)$. Define $\overline{K}$ (resp. $K_h$) as the image of $K$ in $\overline{P}(\A)$ (resp. $P_h(\A)$). Thus, we have morphisms between mixed Shimura varieties $\sh_K(P,\ca{X})\to \sh_{\overline{K}}(\overline{P},\overline{\ca{X}})\to \sh_{K_h}(G_h,\ca{X}_h)$. The reflex fields of these mixed Shimura varieties are the same.\par

If we assume that $(P, \ca{X})=(P_\Phi,D_\Phi)$ comes from a cusp label representative $\Phi \in \ca{CLR}(G,X)$ of a Shimura datum $(G,X)$, we call such mixed Shimura data/varieties the boundary mixed Shimura data/varieties. See \cite[\S 1.1.3]{Wu25}. 
In this case, $U$ is the center of the unipotent radical of $P$.
Note that there exist mixed Shimura varieties that do not arise in this manner.\par
Write the cusp label representative $\Phi$ as $(Q_\Phi,X_\Phi^+,g_\Phi)$. Consider an open immersion:
\begin{equation}\label{eq-open-emb}
    U({\Phi}) = P_\Phi(\rQ)_+\backslash X^+_\Phi \times P(\A)/K_{\Phi} \subset \shu{K_{\Phi}}(P_\Phi, D_\Phi)(\CC).
\end{equation}
Let 
\begin{equation*}
    \beta_{\Phi}: U({\Phi}) \to \shu{K}(G, X)(\CC),\quad [(x, p)] \mapsto [(x, pg)].
\end{equation*}
By varying $\Phi$, we obtain a covering of $\shu{K}(G, X)(\CC)$ via $\beta_{\Phi}$. The (analytic) construction of $\shuc{K}{\Sigma}(G, X)(\CC)$ uses compactifications of $U({\Phi})$.\par
Note that we can view $\beta_{\Phi}$ as the composition of a natural projection
\begin{equation}\label{eq: beta, sep, proj}
  p:   P_\Phi(\rQ)_+ \backslash X_\Phi^+ \times P_\Phi(\A)/K_{\Phi} \to G(\rQ)_+ \backslash X^+_\Phi \times G(\A)/g_\Phi Kg^{-1}_{\Phi}
\end{equation} 
followed by a right action
\begin{equation}\label{eq: beta, sep, hecke}
    [\cdot g_\Phi]:  G(\rQ)_+ \backslash X^+_\Phi \times G(\A)/g_\Phi K g_\Phi^{-1} \to G(\rQ)_+ \backslash X^+_\Phi \times G(\A)/K.
\end{equation}

Given $x\in X$, denote by $\mu_x$ the Hodge cocharacter at $x$. From now on, fix homomorphisms $h_0$ and $h_\infty$ from $\bb{S}$ to $H_0$ such that their pre-compositions with $\mu:\bb{G}_{m,\bb{C}}\to \bb{S}_\bb{C}$, sending $z$ to $(z,1)$, are identical after conjugation by a matrix $\mbf{c}\in \mrm{GL}_2(\bb{C})$; the choice of $\mbf{c}$ depends only on the choice of $h_0$ and $h_\infty$, but not on the value of $z$. See, for example, \cite[4.3]{Pin89} or the beginning of \cite[3.6]{Pin92}.

Therefore, $\mrm{int}(\mbf{c})\circ \mu_x=u_x^{Q_\Phi}\circ h_\infty\circ\mu$ factors through $P_{\Phi,\bb{C}}$. We denote $\mu_{\Phi,x}:=\mrm{int}(\mbf{c})\circ\mu_x$ and $\mu_{\Phi,h,x}:\bb{G}_{m,\bb{C}}\xrightarrow{\mu_{\Phi,x}} P_{\Phi,\bb{C}}\to G_{\Phi,h}$. The assignment $\mu_x\mapsto \mu_{\Phi,x}$ is well defined.\par
We have the Borel embedding $X \hookrightarrow \check{X}\iso G/P_{\mu_x}(\bb{C})$ sending a point $x$ to the Hodge filtration defined by $\mu_{x}$. Here, $P_{\ast}$ denotes the subgroup determined by the cocharacter $\ast$.\par 
We then have the following commutative diagram:
\begin{equation*}
\begin{tikzcd}
	X & {x\in X^+} & {D_{\Phi}} & {D_{\Phi, h} := W_\Phi(\R)U_\Phi(\CC)\backslash D_{\Phi}} \\
	{\breve{X}\iso G/P_{\mu_x}(\bb{C})} && {\breve{D}_{\Phi}\iso P_{\Phi}/P_{\mu_{\Phi,x}}(\bb{C})} & {\breve{D}_{\Phi, h} \iso G_{\Phi, h}/P_{\mu_{\Phi,h,x}}(\bb{C}).}
	\arrow[hook, from=1-1, to=2-1]
	\arrow[from=1-2, to=1-1]
	\arrow[from=1-2, to=1-3,"\tau"]
	\arrow[from=1-3, to=1-4]
	\arrow[hook, from=1-3, to=2-3]
	\arrow[hook, from=1-4, to=2-4]
	\arrow[from=2-3, to=2-1,"{\mbf{c}^{-1}(-)\mbf{c}}"']
	\arrow[from=2-3, to=2-4]
\end{tikzcd}
\end{equation*}
\subsection{p-adic local systems}\label{subsec: generic fiber, p-adic local system}
Fix a mixed Shimura datum $(P,\ca{X})$ in this subsection. The subgroup $K\sbst P(\A)$ is assumed to be a neat open compact subgroup. If we assume that $K=K_pK^p$, then this always means that $K_p$ is open compact in $P(\bb{Q}_p)$ and $K^p\sbst P(\Ap)$ is neat open compact. Sometimes, we assume that $(P,\ca{X})=(G,X)$ is a Shimura datum and consider boundary mixed Shimura data.
\subsubsection{}\label{subsubsec-proet-gxk}

Let $(P, \ca{X})$ be a mixed Shimura datum. By \cite[Lem. 3.7(b)]{Pin89}, for compact subgroups $K'\sbst K \subset P(\A)$ with $K=K_pK^p\sbst P(\A)$ neat open compact, we have
\begin{equation*}
   \varprojlim_{K'\sbst K} \shu{K'}(P, \ca{X}) = \varprojlim_{K'\sbst K} P(\rQ) \bss \ca{X} \times P(\A) /K' =\varprojlim_{K'\sbst K} P(\rQ) \backslash \ca{X} \times P(\A) / Z(P)(\rQ)^{-}K',
\end{equation*}
where $Z(P)(\rQ)^{-}$ is the closure of $Z(P)(\bb{Q})^\circ:=\lrbracket{z \in Z(P)(\rQ)| z|_{\ca{X}} = \identity}$ in $Z(P)(\A)$. Similarly, we denote by $Z(P)(\rQ)^{-}_K$ the closure of $Z(P)(\rQ)^\circ \cap K$ in $K$.
We let $K'$ run over the neat open compact subgroups in the form $K'=K_p'K^p$. Taking the inverse limit induces a natural pro-\'etale torsor
 \begin{equation}\label{eq-proetale torsor-mixed}
      \varprojlim_{K_{p}' \subset K_{p}}  \shu{K_{p}'K^p}(P, \ca{X}) \to \shu{K}(P, \ca{X})
 \end{equation}
 under the group $\mathscr{G}_p(P):=K/Z(P)(\rQ)^-_KK^p$.\par
\begin{definition}\label{def-proetale-torsor-aut} 
Let $K_p = \PP(\Z_p)$ be a quasi-parahoric subgroup and define
$$\bb{P}_K(P,\ca{X}):= \varprojlim_{K_{p}' \subset K_{p}}  \shu{K_{p}'K^p}(P, \ca{X})\times^{\ul{K/Z(P)(\rQ)^-_KK^p}}\ul{K^c_p},$$ where $K_p^c=\ca{P}^c(\bb{Z}_p)$ is the one defined by Definition \ref{def-ass-gp-sch-mix-sh} for $(P,K_p)$.\par
Note that $Z(P)_K^{\overline{\ }}= Z(P)_{ac}(\bb{Q})^{\overline{\ }}\cap K$ as $K$ is neat, the quotient map $K\to K_p^c$ factors through $K/Z(P)(\bb{Q})_K^{\overline{\ }}K^p$.
\end{definition}
In the situation where $(P,\ca{X})=(G,X)$ is a Shimura datum, fix a cusp label representative $\Phi=(Q_\Phi,X^+_\Phi,g_\Phi)\in \ca{CLR}(G,X)$. Then, in this context, (\ref{eq-proetale torsor-mixed}) defines a pro-{\'e}tale torsor
 \begin{equation*}
      \varprojlim_{K_{\Phi, p}' \subset K_{\Phi,p}}  \shu{K_{\Phi, p}'K^p_{\Phi}}(P_\Phi, D_{\Phi}) \to \shu{K_{\Phi}}(P_\Phi, D_{\Phi})
 \end{equation*}
 under $K_{\Phi}/Z(P_\Phi)(\rQ)^-_{K_\Phi}K_{\Phi}^p$. Let $\PP_\Phi$ be the smoothing of the closure of $P_\Phi$ in $\GG_{\Phi}$, where $\GG_{\Phi}$ is the conjugate of $\GG$ by $g = g_{\Phi,p} \in G(\Qp)$, the $p$-factor of $g_\Phi$. Let $\PP^{c}_\Phi$ be the quasi-parahoric group scheme associated with $\PP_\Phi$.
 
From Definition \ref{def-proetale-torsor-aut}, we can construct a pro-\'etale torsor $\PPp_{K_{\Phi}} \to \shu{K_{\Phi}}(P_\Phi, D_{\Phi})$ under $\ul{\PP^{c}_\Phi(\Z_p)}$. Note that we can further push out the torsor to $\ca{P}^*_\Phi$; we will in fact do this in the final step.
\subsubsection{} Let us go back to the general setup.
 Let $\rho: P(\Qp) \to P^c(\Qp) \to \GL(W_{\Qp})$ be a finite-dimensional $\Qp$-representation. Since $K_{p}$ is compact, there exists a $\Z_p$-lattice $W_{\Z_p} \subset W_{\Qp}$ such that $\rho(\PP(\Z_p)) \subset \GL(W_{\Z_p})$.

 We construct $\ls_{\rho, W_{\Qp}}$ as follows:  Let
 \begin{equation*}
     K_{ p}^{(n)} = K_{ p} \cap \rho^{-1}(\lrbracket{g\in\GL(W_{\Z_p})|\ g\equiv \identity \mod p^n}).
 \end{equation*}
We then have an \'etale $\Z_p/p^n\Z_p$-local system $\ls_{\rho, W_{\Z_p}, n}$ on $\shu{K}$ defined as
 \begin{equation}\label{eq: construction of local system from torsor, mixed}
    \shu{K_{ p}^{(n)}K^p} \times^{K_{ p}/K_{ p}^{(n)}} W_{\Z_p}/p^n,
 \end{equation}
 and we have 
 \begin{equation*}
     \ls_{\rho, W_{\Z_p}} = \varprojlim_n \ls_{\rho, W_{\Z_p}, n},\quad  \ls_{\rho, W_{\Qp}} =  \ls_{\rho, W_{\Z_p}} \otimes \rQ.
 \end{equation*}
 
 We say $\bb{P}_K:=\lp_{K}(P,\ca{X})$ is \emph{de Rham} if, for any such pair $(\rho, W_{\Z_p})$, $\ls_{\rho, W_{\Qp}}$ is de Rham. It suffices to check this on a single faithful representation of $P^c(\Qp)$. This is independent of the choice of $W_{\Z_p} \subset W_{\Qp}$.
 
 \begin{prop}[{Liu-Zhu}]\label{prop: mixed, de Rham}
     The pro-\'etale torsor $\lp_{K}$ under $\ul{\PP^{c}(\Z_p)}$ over $\shu{K}^\ad$ is de Rham.
 \end{prop}
 \begin{proof}
   By the rigidity theorem in \cite[Thm. 1.1]{liu2017rigidity}, it suffices to show that, on every geometrically connected component of $\shu{K}$, there exists a closed point where the torsor is de Rham.
In \cite[Def. 11.5]{Pin89}, for each morphism of mixed Shimura data $\iota: (T, Y) \to (P, \ca{X})$ (which is also called a special point when $(T, Y)$ is a pure Shimura datum such that $T$ is a torus embedded into $P$) and for $K_T \subset T(\A)$, $K \subset P(\A)$ such that $\iota(K_T) \subset K$, the induced morphism $\shu{K_T}(T, Y)_{\CC} \to \shu{K}(P, \ca{X})_{\CC}$ descends to a map between canonical models $\shu{K_T}(T, Y) \to \shu{K}(P, \ca{X}) \times_{E(P, \ca{X})} E(T, Y)$. Then \cite[Lem. 4.4]{liu2017rigidity} implies that $\ls_{\rho, W_{\Qp}}$ is de Rham at special points. Note that, although the notion of special points for mixed Shimura varieties is slightly generalized in the sense that $(T,Y)$ can be a finite cover of a usual Shimura datum, the argument in \emph{loc. cit.} goes through verbatim since it only involves class field theory. Since we have a special point on each geometrically connected component by \cite[Lem. 11.6]{Pin89} and the transitivity of the action of $P(\A)$ on geometrically connected components, the proposition follows. 
 \end{proof}
Replacing $P$ by $\ovl{P}$ (resp. by $G_{h}$), we have a pro-\'etale torsor $\overline{\lp}_{K} \to \shu{\ovl{K}}(\ovl{P}, \overline{\ca{X}})$ (resp. $\lp_{K_{h}} \to \shu{K_{h}}(G_{h}, \ca{X}_{h})$) under $\ul{\ovl{\PP}^{c}(\Z_p)}$ (resp. $\ul{\ca{G}_{h}^c(\bb{Z}_p)}$). Since $P \to \ovl{P}\to G_{h}$ induces $P^c \to \ovl{P}^c \to G_{h}^c$ by Lemma \ref{lem-levi-c}, it follows directly from the construction that we have a commutative diagram
\begin{equation*}
\begin{tikzcd}
	{\lp_{K}} & {\ovl{\lp}_{K}} & {\lp_{K_{h}}} \\
    {\shu{K}(P, \ca{X})} & {\shu{\ovl{K}}(\ovl{P}, \overline{\ca{X}})} & {\shu{K_{h}}(G_{h}, \ca{X}_{h}).}
	\arrow[from=1-1, to=1-2]
	\arrow[from=1-1, to=2-1]
	\arrow[from=1-2, to=1-3]
	\arrow[from=1-2, to=2-2]
	\arrow[from=1-3, to=2-3]
	\arrow[from=2-1, to=2-2]
	\arrow[from=2-2, to=2-3]
\end{tikzcd}
\end{equation*}
Exactly the same argument as Proposition \ref{prop: mixed, de Rham} shows that $\overline{\bb{P}}_K$ and $\bb{P}_{K_h}$ are also de Rham.
\subsubsection{Hodge-Tate period maps}\label{subsec: HT period map, 1}

Recall that the pro-\'etale torsor $\lp_{K}$ over $\shu{K}(P, \ca{X})$ under $\underline{\ca{P}^c(\Z_p)}$ is de Rham. Then there is a $\underline{\ca{P}^c(\Z_p)}$-equivariant Hodge-Tate period map
\begin{equation}\label{eq: HT, G, G_h}
   \HT_{K}: \lp_K \to \Gra{P^c, \mu^{c, -1}}.
\end{equation}
Here the superscript ``$(-)^c$'' means the projection to $P^c$. Similarly, we have
$$ \HT_{K_h}: \lp_{\ovl{K}} \to \Gra{\ovl{P}^c, \ovl{\mu}^{c, -1}}, \quad \HT_{K_h}: \lp_{K_h} \to \Gra{G_{h}^c, \mu_{h}^{c, -1}}.$$

\begin{lem}\label{lem-red-torsors-HT-maps}
    We have the following commutative diagram:
\begin{equation}\label{eq: reduction of torsors and HT maps}
\begin{tikzcd}
	{\lp_{K}} & {\ovl{\lp}_{K}} & {\lp_{K_{h}}} \\
	{ \Gra{P^c, \mu^{c, -1}}} & { \Gra{\ovl{P}^c, \ovl{\mu}^{c, -1}}} & {\Gra{G_{h}^c, \mu_{h}^{c, -1}}.}
	\arrow[from=1-1, to=1-2]
	\arrow["{\HT_{K}}"', from=1-1, to=2-1]
	\arrow[from=1-2, to=1-3]
	\arrow["{\HT_{\ovl{K}}}"', from=1-2, to=2-2]
	\arrow["{\HT_{K_{h}}}"', from=1-3, to=2-3]
	\arrow[from=2-1, to=2-2]
	\arrow[from=2-2, to=2-3]
\end{tikzcd}
\end{equation}
\end{lem}
\begin{proof}
     We prove the commutativity of the left square, and the right one follows from the same reason. 

By construction, the push-out torsor $\lp_{K}':=\lp_{K} \times^{\underline{\PP^{c}_\Phi(\Z_p)}}\underline{\ovl{\PP}^{c}_\Phi(\Z_p)}$ is the pullback of $\ovl{\lp}_{K}$ along $\shu{K}(P, \ca{X}) \to \shu{\ovl{K}}(\ovl{P}, \ovl{\ca{X}})$.
The Hodge-Tate period map of $\bb{P}'_{K}$ is given by the composition $\lp_{K} \to \Gra{P^c, \mu^{c, -1}} \to \Gra{\ovl{P}^c, \ovl{\mu}^{c, -1}}$. Since the Hodge-Tate period maps are defined using the Hodge filtrations associated with the de Rham pro-\'etale torsors, $\lp_{K}' \to \Gra{\ovl{P}^c, \overline{\mu}^{c, -1}}$ factors through $\ovl{\lp}_{K}$.
\end{proof}
By \cite{liu2017rigidity} (cf. Proposition \ref{prop: mixed, de Rham}) and \cite[\S 2.6]{PR24}, the de Rham local system $\bb{P}_K$ (resp. $\overline{\bb{P}}_K$ and resp. $\bb{P}_{K_h}$) canonically induces a shtuka in $\Sht_{\ca{P}^c,\mu^c,\Spd E}$ (resp. $\Sht_{\overline{\ca{P}}^c,\overline{\mu}^c,\Spd E}$ and resp. $\Sht_{\ca{G}_h^c,\mu_h^c,\Spd E}$) on $\sh_K(P,\ca{X})^\Dia$ (resp. $\sh_{\overline{K}}(\overline{P},\overline{\ca{X}})^\Dia$ and resp. $\sh_{K_h}(G_h,\ca{X}_h)^\Dia$). 
Together with Corollary \ref{cor: quasi-parahoric}, we know that:
\begin{cor}\label{cor: reduction, generic fiber, general}
    We have the following commutative diagram:
\begin{equation}\label{eq: reduction, generic fiber, general}
\begin{tikzcd}
	{\shu{K}(P, \ca{X})^{\Dia}} & {\shu{\ovl{K}}(\ovl{P}, \ovl{\ca{X}})^{\Dia}} & {\shu{K_{h}}(G_{ h}, \ca{X}_{h})^{\Dia}} \\
	{\Sht_{\PP^c, \mu^c, \delta = 1, \Spd E}} & {\Sht_{\ovl{\PP}^c, \bar{\mu}^c, \delta = 1, \Spd E}} & {\Sht_{\GG_{h}^c, \mu_{h}^c, \delta = 1, \Spd E}.}
	\arrow[from=1-1, to=1-2]
	\arrow[from=1-1, to=2-1]
	\arrow[from=1-2, to=1-3]
	\arrow[from=1-2, to=2-2]
	\arrow[from=1-3, to=2-3]
	\arrow[from=2-1, to=2-2]
	\arrow[from=2-2, to=2-3]
\end{tikzcd}
\end{equation}
\end{cor}

\subsubsection{Some sections}

Let $(P, \mathcal{X})$ be a mixed Shimura datum, $1 \to W \to P_1 \to P \to 1$ be an extension where $W$ is a unipotent group. Assume there exists a mixed Shimura datum $(P_1, \mathcal{X}_1)$ such that $(P_1, \mathcal{X}_1)/W \cong (P, \mathcal{X})$. Assume $P_1 = W \rtimes P$, we fix a section $(P, \mathcal{X}) \to (P_1, \mathcal{X}_1)$ (\cite[Prop. 2.17(b)]{Pin89}). Let $K_1 \subset P_1(\A)$ be a neat open compact subgroup such that $K_1 = K_W \rtimes K$, where $K_W = K_1 \cap W(\A) \subset W(\A)$, $K \subset P(\A)$ is the image of $K_1$, then the natural projection $p_K: \shu{K_1}(P_1, \mathcal{X}_1) \to \shu{K}(P, \mathcal{X})$ has a section $e_K: \shu{K}(P, \mathcal{X}) \to \shu{K_1}(P_1, \mathcal{X}_1)$. Given $K' \subset K$ such that $K_1' = K_W' \rtimes K'$, then $e_{K'}$ and $p_{K'}$ are compatible with $e_K$ and $p_K$ respectively. Passing to $p^{\infty}$-limit, we then have a commutative diagram:
\begin{equation*}
\begin{tikzcd}
	{\lp_{K}} & {\lp_{K_1}} & {\lp_{K}} \\
	{\shu{K}(P, \mathcal{X})} & {\shu{K_1}(P_1, \mathcal{X}_1)} & {\shu{K}(P, \mathcal{X}).}
	\arrow["e", from=1-1, to=1-2]
	\arrow[from=1-1, to=2-1]
	\arrow["p", from=1-2, to=1-3]
	\arrow[from=1-2, to=2-2]
	\arrow[from=1-3, to=2-3]
	\arrow["{e_K}", from=2-1, to=2-2]
	\arrow["{p_K}", from=2-2, to=2-3]
\end{tikzcd}
\end{equation*}

Let $\PP$ be a quasi-parahoric model of $P$, and let $\mathcal{W}$ be a smooth affine model of $W$ with connected fibers such that $\PP$ normalizes $\mathcal{W}$. Let $\PP_1 = \mathcal{W} \rtimes \PP$, $K_W = \mathcal{W}(\Z_p)$, $K_1 = \PP_1(\Z_p)$, and $K = \PP(\Z_p)$; then $K_1 = K_W \rtimes K$. The fixed section $(P, \mathcal{X}) \to (P_1, \mathcal{X}_1)$ also fixes a section of the conjugacy class of Hodge cocharacters $\lrbracket{\mu} \to \lrbracket{\mu_1}$. Define $\PP^c$, $\PP_1^c$ as in Definition \ref{def-ass-gp-sch-mix-sh}; then, by Lemma \ref{lem-levi-c}, $P_1^c = W \rtimes P^c$, and $\PP_1^c = \mathcal{W} \rtimes \PP^c$.

\begin{lem}\label{lem-sec-torsors-HT-maps}
    We have a commutative diagram:
\[\begin{tikzcd}
	{\lp_{K}} & {\lp_{K_1}} & {\lp_{K}} \\
	{\Gra{P^c, \mu^{c, -1}}} & {\Gra{P_1^c, \mu_1^{c, -1}}} & {\Gra{P^c, \mu^{c, -1}}}
	\arrow[from=1-1, to=1-2]
	\arrow["{\HT_K}"', from=1-1, to=2-1]
	\arrow[from=1-2, to=1-3]
	\arrow["{\HT_{K_1}}"', from=1-2, to=2-2]
	\arrow["{\HT_K}"', from=1-3, to=2-3]
	\arrow[from=2-1, to=2-2]
	\arrow[from=2-2, to=2-3]
\end{tikzcd}\]
\end{lem}
\begin{proof}
    The right diagram commutes due to Lemma \ref{lem-red-torsors-HT-maps}. For the left diagram, fix a representation $\rho: P^c(\rQ_p) \to P^c_1(\rQ_p) \stackrel{\rho_1}{\to} \GL(W_{\rQ_p})$ and a lattice $W_{\Z_p} \subset W_{\rQ_p}$ stabilized by $\PP_1^c(\Z_p)$. The induced local system $\rho_*\lp_K$ on $\shu{K}(P, \mathcal{X})$ is the pullback of the induced local system $\rho_{1*}\lp_{K_1}$ on $\shu{K_1}(P_1, \mathcal{X}_1)$; thus, the Hodge-Tate period maps commute.
\end{proof}

\begin{lem}\label{lem: section of shtukas, generic fiber}
    We have a commutative diagram:
\[\begin{tikzcd}
	{\shu{K}(P, \mathcal{X})} & {\shu{K_1}(P_1, \mathcal{X}_1)} & {\shu{K}(P, \mathcal{X})} \\
	{\Sht_{\PP^c, \mu^c}} & {\Sht_{\PP_1^c, \mu_1^c}} & {\Sht_{\PP^c, \mu^c}.}
	\arrow["{e_K}", from=1-1, to=1-2]
	\arrow[from=1-1, to=2-1]
	\arrow["{p_K}", from=1-2, to=1-3]
	\arrow[from=1-2, to=2-2]
	\arrow[from=1-3, to=2-3]
	\arrow[from=2-1, to=2-2]
	\arrow[from=2-2, to=2-3]
\end{tikzcd}\]
\end{lem}

We can apply this result to $(P_{\Phi}, D_{\Phi}) \to (\ovl{P}_{\Phi}, \ovl{D}_{\Phi})$ and to $(\ovl{P}_{\Phi}, \ovl{D}_{\Phi}) \to (G_{\Phi, h}, D_{\Phi, h})$.

\begin{lem}\label{lem: section of reductions, generic fiber}
    Assume $K_{\Phi} = (K_{\Phi} \cap U_{\Phi}(\A)) \rtimes \ovl{K}_{\Phi}$ and $\ovl{K}_{\Phi} = (\ovl{K}_{\Phi} \cap V_{\Phi}(\A)) \rtimes K_{\Phi, h}$. Then the reduction (\ref{eq: reduction, generic fiber, general}) admits a section:
\[\begin{tikzcd}
	{\shu{K_{\Phi}}(P_{\Phi}, D_{\Phi})^{\Dia}} & {\shu{\ovl{K}_{\Phi}}(\ovl{P}_{\Phi}, \ovl{D}_{\Phi})^{\Dia}} & {\shu{K_{\Phi, h}}(G_{\Phi, h}, D_{\Phi, h})^{\Dia}} \\
	{\Sht_{\PP_{\Phi}^c, \mu_{\Phi}^c, \delta = 1, \Spd E}} & {\Sht_{\ovl{\PP}_{\Phi}^c, \bar{\mu}_{\Phi}^c, \delta = 1, \Spd E}} & {\Sht_{\GG_{\Phi, h}^c, \mu_{\Phi, h}^c, \delta = 1, \Spd E}.}
	\arrow[from=1-1, to=2-1]
	\arrow[from=1-2, to=1-1]
	\arrow[from=1-2, to=2-2]
	\arrow[from=1-3, to=1-2]
	\arrow[from=1-3, to=2-3]
	\arrow[from=2-2, to=2-1]
	\arrow[from=2-3, to=2-2]
\end{tikzcd}\]
\end{lem}

\subsection{Torus action on the tower}\label{subsec-torus-action-tower}
We re-organize the materials in \cite[3.13-3.16]{Pin89} and \cite[2.1.11]{Mad19}, and describe the torus action in an explicit and functorial way. This exposition is helpful in \S\ref{subsubsec-tower-computation}, where we compute the torus action on a tower of mixed Shimura varieties.\par
Throughout the subsection, we pick any mixed Shimura datum $(P,\ca{X})$ and denote its reflex field by $\bb{E}(P,\ca{X})$; the symbol $E$ denotes its completion at a place $v$ over $p$. Sometimes we use the shorthand $\bb{E}:=\bb{E}(P,\ca{X})$. 
\subsubsection{}\label{subsubsec-can-sec-uipotent-part}
Let $U\iso \bb{G}_a^n$ be a connected additive group over $\bb{Q}$. The group $\bb{G}_m$ acts on $U$ by scalar multiplication. 
Denote by $(\bb{G}_m,\ca{H}_0)$ the ``Siegel Shimura datum of dimension zero'', where $\ca{H}_0=\{\pm\mrm{id}\}$ is the two-point set parameterizing the isomorphisms $\bb{Z}\iso \bb{Z}(1)$.
Let $P_0 = U \rtimes \Gm$, where $\Gm$ acts on $U$ by this scalar multiplication. We can write $P_0$ as a blocked mirabolic subgroup 
\begin{equation*}
    \begin{pmatrix}
        \bb{G}_m& U\\
        0& 1
    \end{pmatrix}.
\end{equation*}
In what follows, denote $U:=\begin{pmatrix}
    1& U\\0&1
\end{pmatrix}$ by abuse of notation. \par
Let $X_0$ be $\pi_0(\ca{H}_0)\times$
$$(P_0(\bb{R})U(\bb{C})\text{-orbit of the homomorphism }(z_1,z_2)\in\bb{S}(\bb{C})\mapsto \mrm{diag}\{z_1z_2,1\}\in P_0(\bb{C})).$$ 
Then $(P_0,X_0)$ is a mixed Shimura datum with a natural projection $\lambda_0:(P_0,X_0)\to(\bb{G}_m,\ca{H}_0)$ obtained by the quotient by $U$. See also \cite[2.24]{Pin89}.\par
Let $K_0$ be a neat open compact subgroup of $P_0(\A)$. 
By \cite[3.13]{Pin89}, $\sh_{K_0}(P_0,X_0)(\bb{C})$ is a torus torsor under $\mbf{E}_{K_0}(\bb{C})$, where $\mbf{E}_{K_0}(\bb{C})\iso U(\bb{C})/\Lambda_0$ for $\Lambda_0:=(U(\bb{Q})\cap K_0)(-1)$. More precisely, fix any $p\in P_0(\A)$, the fiber of $\sh_{K_0}(P_0,X_0)(\bb{C})\to \sh_{\lambda_0(K_0)}(\bb{G}_m,\ca{H}_0)(\bb{C}):=\bb{Q}^\times_{>0}\bss \{1\}\times \A^\times/K_0$ at $[(1,\lambda_0(p))]$ determined by $\lambda_0(p)\in \bb{G}_m(\A)$ is $U(\bb{C})/(U(\bb{Q})\cap pK_0p^{-1})(-1)$. Note that $U(\bb{Q})\cap pK_0p^{-1}$ and $U(\bb{Q})\cap K_0$ differ by homothety; let us explain this in the next paragraph.\par
Recall that there is a unique factorization of multiplicative groups $\A^\times=\bb{Q}^\times_{>0}\times \wat{\bb{Z}}^\times$. For $g\in \A^\times$, denote by $r(g)$ the $\bb{Q}^\times_{>0}$-factor. Denote $$r(p):=r\circ \lambda_0(p)=\begin{pmatrix}r(p)&0\\0&1\end{pmatrix}.$$ 
Denote $\Lambda_{0,\A}:=U(\A)\cap K_0$. Writing 
\begin{equation}\label{eq-torus-exp}
\mbf{E}_{K_0}(\bb{C})\iso U(\bb{Q})\bss U(\bb{C})\times U(\A)/\Lambda_{0,\A},\end{equation}
the homothety is given by the left conjugate by $r(p)$:
$$\mbf{t}(p):U(\bb{Q})\bss U(\bb{C})\times U(\A)/\Lambda_{0,\A}\xrightarrow{\sim} U(\bb{Q})\bss U(\bb{C})\times U(\A)/p\Lambda_{0,\A}p^{-1},$$
sending $[(u,u_f)]_{\Lambda_{0,\A}}\mapsto [(r(p)ur(p)^{-1},r(p)u_fr(p)^{-1})]_{r(p)\Lambda_{0,\A}r(p)^{-1}}.$ Note that $r(p)\Lambda_{0,\A}r(p)^{-1}={p\Lambda_{0,\A}p^{-1}}$. This isomorphism between the affine complex tori canonically descends to an isomorphism between the algebraic tori.\par
Define an action
$$\bb{E}_{K_0}(P_0,X_0)(\bb{C}): \mbf{E}_{K_0}(\bb{C})\times \sh_{K_0}(P_0,X_0)(\bb{C})\to \sh_{K_0}(P_0,X_0)(\bb{C})$$ 
as follows:\par
Fix an isomorphism $\mrm{sgn}:\ca{H}_0\iso \{\pm 1\}$. Let $\delta(x)=\mrm{diag}\{\mrm{sgn}(x),1\}$ for $x\in \ca{H}_0$. Picking any $p\in P_0(\A)$, set $[(u,u_f)]_{\Lambda_{0,\A}}\in \mbf{E}_{K_0}(\bb{C})\mapsto$ 
\begin{equation}\label{eq-torus-action-def}
([(x,p)]_{K_0}\mapsto[(\mrm{int}(\delta(x)r(p)ur(p)^{-1}\delta(x)^{-1})x, \delta(x)r(p)u_fr(p)^{-1}\delta(x)^{-1}p)]_{K_0}).
\end{equation}
\begin{lem}\label{lem-ck0-well-defined}
The map $\bb{E}_{K_0}(P_0,X_0)(\bb{C})$ is well defined and algebraic. Moreover, it descends to $\bb{Q}$.
\end{lem}
\begin{proof}
Fix any $u_0\in \Lambda_{0,\A}$. There is $u_0'\in \Lambda_{0,\A}$ such that $pu_0'p^{-1}=r(p)u_0r(p)^{-1}$, because the conjugation of $\wat{\bb{Z}}^\times$ and $U(\A)$ stabilizes $\Lambda_{0,\A}$. Assume that $\mrm{sgn}(x)=1$ without loss of generality. Then the assignment $\bb{E}_{K_0}(P_0,X_0)(\bb{C})(u,u_f\cdot u_0)$ sends $[(x,p)]_{K_0}$ to 
\begin{equation*}
\begin{split}
    &[(\mrm{int}(r(p)ur(p)^{-1})x,r(p)u_fu_0r(p)^{-1}\cdot p)]_{K_0} \\
    &=[(\mrm{int}(r(p)ur(p)^{-1})x,r(p)u_fr(p)^{-1}r(p)u_0r(p)^{-1}\cdot p)]_{K_0}\\
    &=[(\mrm{int}(r(p)ur(p)^{-1})x,r(p)u_fr(p)^{-1}pu_0'p^{-1}\cdot p)]_{K_0}\\
    &=[(\mrm{int}(r(p)ur(p)^{-1})x,r(p)u_fr(p)^{-1}\cdot p)]_{K_0}.
    \end{split}
\end{equation*}
It is easy to check that this map does not change when multiplying $(u,u_f)$ to $(u'u,u'u_f)$ by $u'\in U(\bb{Q})$, and that the compatibility of this map with the multiplication on $\mbf{E}_{K_0}$. So this is a well-defined action. Since $\mbf{t}(p)$ is algebraic and since it suffices to check the algebraicity on each individual connected component, we find that $\bb{E}_{K_0}(P_0,X_0)(\bb{C})$ is algebraic.\par
Now we consider the Galois action. It suffices to check on a dense subset. So we may assume that $u\in U(\bb{Q})$, $u_f=1$ by strong approximation, and $[(x,p)]_K$ is contained in a special point $(T,\ca{Y})\sbst (P_0,X_0)$ such that 
$$\ca{Y}\sbst \pi_0(\ca{H}_0)\times P_0(\bb{Q})\text{-orbit of the homomorphism }h_0:(z_1,z_2)\in \bb{S}(\bb{C})\mapsto \mrm{diag}\{z_1z_2,1\}\in P_0(\bb{C}).$$ 

Then, by taking a conjugation in $P_0(\bb{Q})$, we may consider the case where $T=\mrm{diag}\{\bb{G}_m,1\}$ and $\ca{Y}=\{\pm h_0\}$. 
Again, we assume that $\mrm{sgn}(x)=1$ and the case of $x=-1$ is similar.
Now, for $\tau\in \gal(\bb{Q}^{ab}/\bb{Q})$, we denote by $[(x,p)]\mapsto [(x,r_{(T,\ca{Y})}(\tau)p)]$ the Galois action on the points of $\sh_{K_T}(T,\ca{Y})(\overline{\bb{Q}})$ for a suitable $K_T$. Write $d:=r_{(T,\ca{Y})}(\tau)\in \pi_0(\bb{G}_m(\bb{A})/\bb{G}_m(\bb{Q}))$; we choose a representative of it in $\wat{\bb{Z}}^\times$. Writing $u=\begin{pmatrix}
    1&\ul{u}\\
    0&1
\end{pmatrix}$, we check that $\tau\cdot [(u,1)]\cdot [(x,p)]$
\begin{equation}\label{computation-gal-1}
    \begin{split}
    &=\tau\cdot [(\mrm{int}(r(p)ur(p)^{-1})x,p)]\\
    &=\tau\cdot[(x,r(p)u^{-1}r(p)^{-1}\cdot p)]\\
    &=\tau\cdot [(x,\begin{pmatrix}
        p&-r(p)\ul{u}\\
        0&1
    \end{pmatrix})]\\
    &=[(x, \begin{pmatrix}
        p&-r(p)\ul{u}\\
        0&1
    \end{pmatrix}\begin{pmatrix}
        d&0\\
        0&1
    \end{pmatrix})]\\
    &=[(x,\begin{pmatrix}
      dp& -r(p)\ul{u}\\
      0& 1
    \end{pmatrix})],
    \end{split}
\end{equation} 
and that $[(u,1)]\cdot \tau\cdot [(x,p)]$
\begin{equation}\label{computation-gal-2}
    \begin{split}
    &=[(u,1)]\cdot[(x,d \cdot p)]\\
    &=[(\mrm{int}(r(d)r(p)ur(p)^{-1}r(d)^{-1})x,d\cdot p)]\\
    &=[(x,r(d)r(p)u^{-1}r(p)^{-1}r(d)^{-1}dp)]\\
    &=[(x,r(p)u^{-1}r(p)^{-1}dp)]\\
    &=[(x,\begin{pmatrix}
        dp&-r(p)\ul{u}\\
        0&1
    \end{pmatrix})].
    \end{split}
\end{equation}
Hence, we see that (\ref{computation-gal-1})$=$(\ref{computation-gal-2}). Note that the multiplication of $d$ induced by the Galois action in (\ref{computation-gal-1}) is on the right because one can check that this is the only way such that the Galois action determining the canonical model is compatible with right Hecke actions. 
By density, the Galois action commutes with the action of $\mbf{E}_{K_0}(\bb{C})$ and $\bb{E}_{K_0}(P_0,X_0)(\bb{C})$ descends to $\bb{Q}$.
\end{proof}
Denote the algebraic action over $\bb{Q}$ by $\bb{E}_{K_0}(P_0,X_0)$.
\begin{convention}
Note that an isomorphism $\mrm{sgn}$ is fixed. If we change $\mrm{sgn}$ to $-\mrm{sgn}$, the $\mbf{E}_{K_0}$-action in (\ref{eq-torus-action-def}) will be reversed. Moreover, the isomorphism (\ref{eq-torus-exp}) is induced by an exponential map, which also depends on a choice of $\sqrt{-1}$, and here we fixed $\bb{C}$ and chose $\sqrt{-1}$ to be $i\in\bb{C}$. Nevertheless, the torsor structure given by (\ref{eq-torus-action-def}) is canonical up to unique isomorphism.
\end{convention}
From the explicit construction above, we see that
\begin{lem}\label{lem-functorial-ck0}
The morphism $\bb{E}_{K_0}(P_0,X_0)$ is functorial in $(U,K_0)$. That is, for any homomorphism $f:U_1\to U_2$ between additive groups as in the beginning of \S\ref{subsubsec-can-sec-uipotent-part}, we produce $(P_0^1,X_0^1)$ (resp. $(P_0^2,X_0^2)$) for $U_1$ (resp. $U_2$) as in the first paragraph. Then there is a morphism between mixed Shimura data $\wdtd{f}:(P_0^1,X_0^1)\to (P_0^2,X_0^2)$ induced by $f$. Picking neat open compact subgroups $K_0^1$ and $K_0^2$ such that $\wdtd{f}(K^1_0)\sbst K^2_0$, we obtain a morphism between mixed Shimura varieties $\wdtd{f}:\sh_{K^1_0}(P_0^1,X_0^1)\to\sh_{K^2_0}(P_0^2,X_0^2)$ over $\bb{Q}$ that fits into the commutative diagram
\begin{equation*}
    \begin{tikzcd}
        \mbf{E}_{K_0^1}\times\sh_{K_0^1}(P_0^1,X_0^1)\arrow[rr,"{\bb{E}_{K_0^1}(P_0^1,X_0^1)}"]\arrow[d,"{[f]\times\wdtd{f}}"]&&\sh_{K_0^1}(P_0^1,X_0^1)\arrow[d,"\wdtd{f}"]\\
        \mbf{E}_{K_0^2}\times \sh_{K_0^2}(P_0^2,X_0^2)\arrow[rr,"{\bb{E}_{K_0^2}(P_0^2,X_0^2)}"]&&\sh_{K_0^2}(P_0^2,X_0^2).
    \end{tikzcd}
\end{equation*}
The map $[f]$ of the homomorphism between tori is induced by $f$ and exponential maps.
\end{lem}
\begin{proof}
This follows from the explicit formula (\ref{eq-torus-action-def}) above.
\end{proof}
\subsubsection{}
Now, we consider a mixed Shimura datum $(P,\ca{X})$. There is a weight filtration on $\lie P$ given by the definition of $(P,\ca{X})$ being a mixed Shimura datum. 
Let $U$ be the unipotent group whose Lie algebra corresponds to the weight-($-2$)-part of $\lie P$ and $W$ be the unipotent radical of $P$. Recall that the adjoint representation of $P$ on $\lie U\iso U$ factors through $G_h:=P/W$.\par 
For our purpose, we assume that this action factors as $$\mrm{Ad}: P\xrightarrow{\lambda} \bb{G}_m\to \mrm{GL}(U),$$ 
where $\bb{G}_m$ acts as scalar multiplication.\par
As in \S\ref{subsubsec-can-sec-uipotent-part}, we define $r(p)=r\circ\lambda(p)$ for $p\in P(\A)$.\par
Moreover, we can extend $\lambda$ to a morphism between mixed Shimura data $\lambda: (P, \ca{X}) \to (\Gm, \HH_0)$; picking any connected component $\ca{X}^+\sbst \ca{X}$, this extension depends on the choice of $\ca{X}^+\to \HH_0$, which has two options.\par
By \cite[3.13]{Pin89}, for any neat open compact subgroup $K\sbst P(\A)$, the morphism $\sh_{K}(P,\ca{X})(\bb{C})\to \sh_{\overline{K}}(\overline{P},\overline{\ca{X}})$ is a torus torsor under $\mbf{E}_K(\bb{C}):=U(\bb{C})/\Lambda$. Here $\Lambda=[p_2(Z(P)(\bb{Q})^\circ\times U(\bb{Q})\cap K)](-1)$, where $Z(P)(\bb{Q})^\circ$ is the centralizer in $Z(P)(\bb{Q})$ of $\ca{X}$. \par
As (\ref{eq-torus-action-def}), we define the action $\bb{E}_K(P,\ca{X})(\bb{C})$ of $[(u,u_f)]_{\Lambda_{\A}}$: 
\begin{equation}\label{eq-torus-action-gen}
  ([(x,p)]_{K}\mapsto[(\mrm{int}(\delta(x)r(p)ur(p)^{-1}\delta(x)^{-1})x, \delta(x)r(p)u_fr(p)^{-1}\delta(x)^{-1}p)]_{K}).  
\end{equation}
\begin{lem}\label{lem-torus-action-gen}
The action (\ref{eq-torus-action-gen}) is well defined.\par 
We can write the action as a $4$-step composition:
\begin{equation}\label{eq-action-4-step}
\begin{split}
&\bb{E}_K(P,\ca{X})(\bb{C}):\mbf{E}_{K}\times \sh_K(P,\ca{X})(\bb{C})=\\
&\mbf{E}_{K}\times \sh_{\lambda(K)}(\bb{G}_m,\ca{H}_0)(\bb{C})\times_{\sh_{\lambda(K)}(\bb{G}_m,\ca{H}_0)(\bb{C})}\sh_K(P,\ca{X})(\bb{C})\to \\
&\mbf{E}_K\times \sh_{K_0}(P_0,X_0)(\bb{C})\times_{\sh_{\lambda(K)}(\bb{G}_m,\ca{H}_0)(\bb{C})}\sh_K(P,\ca{X})(\bb{C})\xrightarrow{\ \bb{E}_{K_0}(P_0,X_0)(\bb{C})\ }\\
&\sh_{K_0}(P_0,X_0)(\bb{C})\times_{\sh_{\lambda(K)}(\bb{G}_m,\ca{H}_0)(\bb{C})}\sh_K(P,\ca{X})(\bb{C})\xrightarrow{\ m_K(P,\ca{X})(\bb{C})\ } \\
&\sh_K(P,\ca{X})(\bb{C}).
\end{split}
\end{equation}
Here, $(P_0,X_0)$ is the one defined by $U$ in \S\ref{subsubsec-can-sec-uipotent-part} and $K_0=[p_2(Z(P)(\bb{Q})^\circ \times U(\A)\cap K)]\rtimes \lambda(K)$. The first equality is from the definition of fiber products, the second arrow is induced by the diagonal section from $(\bb{G}_m,\ca{H}_0)$ to $(P_0,X_0)$ by construction, and the last arrow is given by the multiplication operation in \cite[2.22]{Pin89} (in which it was denoted by ``$\mu$'').\par
Moreover, the composition $\bb{E}_K(P,X)(\bb{C})$ in (\ref{eq-action-4-step}) descends to an algebraic morphism over $\bb{E}(P,\ca{X})$. 
\end{lem}
\begin{proof}
The first paragraph is verified exactly the same way as Lemma \ref{lem-ck0-well-defined}. The second paragraph is proved in \cite[2.22 and Cor. 3.12]{Pin89}, combined with Lemma \ref{lem-ck0-well-defined}. The third paragraph then follows from the second sentence of Lemma \ref{lem-ck0-well-defined}, the functoriality of canonical models of mixed Shimura varieties induced by morphisms between mixed Shimura data \cite[Ch. 11]{Pin89} and the identification  $(P_0,X_0)\times_{(\bb{G}_m,\ca{H}_0)}(P,\ca{X})=(U\rtimes P,\ca{X}')$ for some $\ca{X}'$.
\end{proof}
We denote the algebraic morphism above by $\bb{E}_K(P,\ca{X})$.\par
\begin{lem}\label{lem-functorial-c-action}
The formalism of the action $\bb{E}_K(P,\ca{X})$ is functorial in $(P,\ca{X},K)$.
\end{lem}
\begin{proof}
This follows from the fact that all steps in the composition (\ref{eq-action-4-step}) in the second paragraph of Lemma \ref{lem-torus-action-gen} are functorial. See Lemma \ref{lem-functorial-ck0} and \cite[Ch. 11]{Pin89}.
\end{proof}
\begin{rk}
As mentioned above, changing the choice of $\ca{X}^+\to \ca{H}_0$ will change $\bb{E}_K(P,\ca{X})$ to $\bb{E}_K(P,\ca{X})\circ ((-\mrm{id})\times \mrm{id})$, but the $\mbf{E}_K$-torsor structure of $\sh_K(P,\ca{X})$ is canonical up to this isomorphism.
\end{rk}
\subsubsection{}\label{subsubsec-tower-computation}
Fix a choice of $\ca{X}^+\to \ca{H}_0$. 
For any pair of neat open compact subgroups $K'\sbst K$, we have a commutative diagram over the reflex field $\bb{E}$ by Lemma \ref{lem-functorial-c-action}:
  \begin{equation}\label{diag-change-K}
\begin{tikzcd}
	{\mbf{E}_{K'}\times \shu{K'}(P, \ca{X})}\arrow[rr,"{\bb{E}_{K'}(P,\ca{X})}"]\arrow[d] && { \shu{K'}(P, \ca{X})} \arrow[d]\\
	{\mbf{E}_{K}\times \shu{K}(P, \ca{X})} \arrow[rr, "{\bb{E}_{K}(P,\ca{X})}"]&& { \shu{K}(P, \ca{X}).}
\end{tikzcd}
  \end{equation}
Assume that $K=K_pK^p$. Taking the limit over all $K'=K_p'K^p\sbst K$, we denote $\mbf{E}_{\infty,K}:=\varprojlim_{K_p'\sbst K_p}\mbf{E}_{K'}$. Note that the Galois group of $\mbf{E}_{\infty,K,E}$ over $\mbf{E}_{K,E}$ is 
$$\mathscr{G}_p(\mbf{E}_K):= \varprojlim_{K'_p\sbst K_p} \frac{p_2(Z(P)(\bb{Q})^\circ\times U(\bb{Q})\cap K_pK^p)}{p_2(Z(P)(\bb{Q})^\circ\times U(\bb{Q})\cap K_p'K^p)}.$$

When $P=P^c$, from the neatness of $K$, the group $\mathscr{G}_p(\mbf{E}_K)=$
$$\varprojlim_{K'_p\sbst K_p} \frac{U(\bb{Q})\cap K_pK^p}{U(\bb{Q})\cap K_p'K^p}=\varprojlim_{K_p'\sbst K_p} \frac{(U(\A)\cap K_pK^p)\cap U(\bb{Q})}{(U(\A)\cap K_p'K^p)\cap U(\bb{Q})}.$$
The group $U$ is unipotent, and $U(\A)\cap K_p K^p=(U(\bb{Q}_p)\cap K_p)(U(\Ap)\cap K^p)$. Therefore, the group $\mathscr{G}_p(\mbf{E}_K)$ is pro-$p$, and $\mbf{E}_{\infty,K}$ is the inverse limit of $\cdots \xrightarrow{p}\mbf{E}_K\xrightarrow{p}\mbf{E}_K\xrightarrow{p}\mbf{E}_K$. Taking the inverse limit of (\ref{diag-change-K}) over $K'=K_p'K^p$, we have
\begin{equation}\label{diag-lim-K_p-cuspidal}
\begin{tikzcd}
    \mbf{E}_{\infty,K}\times \bb{P}_K\arrow[d]\arrow[rr,"{\bb{E}_{K^p}(P,\ca{X})}"]&&\bb{P}_K\arrow[d]\\
    \mbf{E}_{K}\times \sh_K(P,\ca{X})\arrow[rr,"{\bb{E}_K(P,\ca{X})}"]&&\sh_K(P,\ca{X}),
\end{tikzcd}
\end{equation}
where $\bb{E}_{K^p}(P,\ca{X})=\varprojlim_{K_p'\sbst K_p}\bb{E}_{K'}(P,\ca{X})$.
\par
Now, we consider a general $P$. Let $K_p^c$ be defined as in Definition \ref{def-ass-gp-sch-mix-sh}, and let $K^{c,p}$ be the image of $K^p$ in $P^c(\Ap)$. 
Note that $K_p^c$ is open compact and $K^{c,p}$ is neat open compact.\par
By Lemma \ref{lem-functorial-c-action}, there is a commutative diagram induced by the morphism $(P,\ca{X},K)\to (P^c,\ca{X}^c,K^c)$:
\begin{equation*}
    \begin{tikzcd}
    {\mbf{E}_{K}\times \shu{K}(P, \ca{X})}\arrow[rr,"{\bb{E}_{K}(P,\ca{X})}"]\arrow[d] && { \shu{K}(P, \ca{X})} \arrow[d]\\
	{\mbf{E}_{K^c}\times \shu{K^c}(P^c, \ca{X}^c)} \arrow[rr, "{\bb{E}_{K^c}(P^c,\ca{X}^c)}"]&& { \shu{K^c}(P^c, \ca{X}^c).}
    \end{tikzcd}
\end{equation*}
\begin{lem}\label{lem-proet-torsor-commutative-diag}
Taking the inverse limit with respect to all open compact subgroups $K_p'\sbst K_p$, we have a commutative diagram
\begin{equation*}\label{diag-functoriality-kc-kinf}
    \begin{tikzcd}
    {\mbf{E}_{\infty,K}\times \varprojlim_{K_p'\sbst K_p}\sh_{K'}(P,\ca{X})}\arrow[r,"{\bb{E}_{K^p}(P,\ca{X})}"]\arrow[dr]\arrow[dd]&{\varprojlim_{K_p'\sbst K_p}\sh_{K'}(P,\ca{X})}\arrow[dr]\arrow[dd]&\\
    &{\mbf{E}_K\times\sh_K(P,\ca{X})}\arrow[r,"{\bb{E}_K(P,\ca{X})}"]\arrow[dd,shift right=1ex]&\sh_K(P,\ca{X})\arrow[dd]\\
    {\mbf{E}_{\infty,K^c}\times \varprojlim_{K_p^{\prime,c}\sbst K_p^c}\sh_{K^{\prime,c}}(P^c,\ca{X}^c)}\arrow[r,"{\bb{E}_{K^{c,p}}(P^c,\ca{X}^c)}"]\arrow[dr]&{\varprojlim_{K_p^{\prime,c}\sbst K_p^c} \sh_{K^{\prime,c}}(P^c,\ca{X}^c)}\arrow[dr]\\
    &{\mbf{E}_{K^c}\times \sh_{K^c}(P^c,\ca{X}^c)}\arrow[r,"{\bb{E}_{K^c}(P^c,\ca{X}^c)}"]&\sh_{K^c}(P^c,\ca{X}^c).
    \end{tikzcd}
\end{equation*}
Moreover, the commutative diagram 
\begin{equation}\label{diag-equi-torus-action}
    \begin{tikzcd}
    {\mbf{E}_{\infty,K}\times \varprojlim_{K_p'\sbst K_p}\sh_{K'}(P,\ca{X})}\arrow[rr,"{\bb{E}_{K^p}(P,\ca{X})}"]\arrow[d]&&{\varprojlim_{K_p'\sbst K_p}\sh_{K'}(P,\ca{X})}\arrow[d]\\
    {\mbf{E}_{\infty,K^c}\times \varprojlim_{K_p^{\prime,c}\sbst K_p^c}\sh_{K^{\prime,c}}(P^c,\ca{X}^c)}\arrow[rr,"{\bb{E}_{K^{c,p}}(P^c,\ca{X}^c)}"]&&{\varprojlim_{K_p^{\prime,c}\sbst K_p^c} \sh_{K^{\prime,c}}(P^c,\ca{X}^c)}
    \end{tikzcd}
\end{equation}
is equivariant under the commutative diagram
\begin{equation*}\label{eq-group-comm-diag}
    \begin{tikzcd}
    {\{1\}\times K/Z(P)(\bb{Q})^{\overline{\ }}_KK^p}\arrow[rr,"\mrm{id}"]\arrow[d]&&{ K/Z(P)(\bb{Q})^{\overline{\ }}_KK^p}\arrow[d]\\
    {\{1\}\times K^c_p}\arrow[rr,"\mrm{id}"]&& K^c_p
    \end{tikzcd}
\end{equation*}
over
\begin{equation*}
    \begin{tikzcd}
    {\mbf{E}_K\times\sh_K(P,\ca{X})}\arrow[rr,"{\bb{E}_K(P,\ca{X})}"]\arrow[d]&&\sh_K(P,\ca{X})\arrow[d]\\
    {\mbf{E}_{K^c}\times \sh_{K^c}(P^c,\ca{X}^c)}\arrow[rr,"{\bb{E}_{K^c}(P^c,\ca{X}^c)}"]&&\sh_{K^c}(P^c,\ca{X}^c).
    \end{tikzcd}
\end{equation*}
\end{lem}
\begin{proof}
The first diagram follows from Lemma \ref{lem-functorial-c-action}. 
Let us show the equivariance of (\ref{diag-equi-torus-action}). From the functoriality between mixed Shimura varieties, the two vertical arrows are $(K/Z(P)(\bb{Q})^{\overline{\ }}_KK^p\to K_p^c)$-equivariant.\par
We check the $( K/Z(P)(\bb{Q})_K^{\overline{\ }}K^p\to K/Z(P)(\bb{Q})_K^{\overline{\ }}K^p)$-equivariance of $\bb{E}_{K^p}(P,\ca{X})$. Fix $K'_p\sbst K_p$ and assume that it is a normal subgroup. It suffices to check this over $\bb{C}$. 
Pick $[(u,1)]\in \mbf{E}_{K'}$ and $k\in K_p$. On the one hand, we compute that $$ [(u,1)]\cdot k\cdot [(x,p)]_{K'}=[(u,1)]\cdot[(x,pk)]_{K'}=[(x,\mrm{int}(\delta(x)r(pk))(u^{-1})pk)]_{K'};$$
and, on the other hand, $k\cdot[(u,1)]\cdot[(x,p)]_{K'}=[(x,\mrm{int}(\delta(x)r(p))(u^{-1})pk)]_{K'}$. The two expressions are the same since $r(k)=1$.
\end{proof}
\begin{lem}\label{lem-torus-action-tower-pro-p-gen}
With the conventions in Lemma \ref{lem-proet-torsor-commutative-diag}, the pushout $$\bb{P}_K:=\varprojlim_{K_p'\sbst K_p}\sh_{K'}(P,\ca{X})\times^{\ul{K/Z(P)(\bb{Q})^{\overline{\ }}_KK^p}}\ul{K_p^c}$$ is the pullback of $\bb{P}_{K^c}$ on $\sh_{K^c}(P^c,\ca{X}^c)$. 
The action $\bb{E}_{K^p}(P,\ca{X})$ induces a $(\mbf{E}_{\infty,K}'\to \mbf{E}_K)$-equivariant action $\bb{E}'_{K^p}(P,\ca{X})$ on $\bb{P}_K\to \sh_K(P,\ca{X})$ such that $\mbf{E}'_{\infty,K}$ is an inverse limit of $p$-power covers of $\mbf{E}_K$, i.e., $\mbf{E}'_{\infty,K}=\varprojlim_{i\geq 0}\mbf{E}_K^{(i)}$, where $\mbf{E}_K^{(0)}=\mbf{E}_K$ and $\mbf{E}_K^{(i)}\to \mbf{E}_K$ are $p$-power isogenies between tori for positive integers $i$. \par
In other words, there is a commutative diagram
\begin{equation*}
\begin{tikzcd}
    {\mbf{E}_{\infty,K}'\times \bb{P}_K}\arrow[rr,"{\bb{E}_{K^p}'(P,\ca{X})}"]\arrow[d]&&{\bb{P}_K}\arrow[d]\\
    {\mbf{E}_{K}\times \sh_K(P,\ca{X})}\arrow[rr,"{\bb{E}_{K}(P,\ca{X})}"]&&{\sh_K(P,\ca{X}).}
\end{tikzcd} 
\end{equation*}
\end{lem}
\begin{proof}
The first sentence follows from the equivariance statement in Lemma \ref{lem-proet-torsor-commutative-diag}. 
Fix an open compact normal subgroup $K_p'\sbst K_p$. Consider the pullback of 
\begin{equation*}
    \begin{tikzcd}
    {\mbf{E}_{K^{\prime,c}}\times\sh_{K^{\prime,c}}(P^c,\ca{X}^c)}\arrow[rr,"{\bb{E}_{K^{\prime,c}}(P^c,\ca{X}^c)}"]\arrow[d]&&{\sh_{K^{\prime,c}}(P^c,\ca{X}^c)}\arrow[d]\\
    {\mbf{E}_{K^c}\times\sh_{K^c}(P^c,\ca{X}^c)}\arrow[rr,"{\bb{E}_{K^c}(P^c,\ca{X}^c)}"]&&{\sh_{K^c}(P^c,\ca{X}^c)}
    \end{tikzcd}
\end{equation*}
along $\mbf{E}_K\times \sh_K(P,\ca{X})\to \mbf{E}_{K^c}\times\sh_{K^c}(P^c,\ca{X}^c)$. It suffices to show that the connected component of the fiber product $\mbf{E}_{K,\bb{E}}\times_{\mbf{E}_{K^c,\bb{E}}}\mbf{E}_{K^{\prime,c},\bb{E}}$ is a torus $\mbf{E}$ such that the morphism $\mbf{E}\to \mbf{E}_{K,\bb{E}}$ induced by the projection to the first factor of this fiber product is a $p$-power isogeny. This is Lemma \ref{lem-torus-isog-p-ele}.
\end{proof}
\begin{lem}\label{lem-torus-isog-p-ele}
Let $\mbf{E}_0$, $\mbf{E}_1$ and $\mbf{E}_2$ be split tori of the same rank over a field $k$ of characteristic zero. Let $f:\mbf{E}_1\to\mbf{E}_0$ and $g:\mbf{E}_2\to \mbf{E}_0$ be isogenies between tori. Assume that $g$ is a $p$-power isogeny. Then the identity component $\mbf{E}_3$ of $\mbf{E}_1\times_{f,\mbf{E}_0,g}\mbf{E}_2$ is a split torus equipped with a $p$-power isogeny $\mbf{E}_3\to\mbf{E}_1$.
\end{lem}
\begin{proof}
This fiber product is smooth, so irreducible components are exactly the connected components. If we know that $\mbf{E}_3$ is a split torus, then $\mbf{E}_3\to \mbf{E}_1$ is a $p$-power isogeny by construction. 
We show that $\mbf{E}_3$ is a split torus. 
Write the character groups of $\mbf{E}_i$ by $X_i$ for $i=0,1,2$. Then there are inclusions $X_0\sbst X_1\sbst X_0\otimes\bb{Q}$ and $X_0\sbst X_2\sbst X_0\otimes\bb{Q}$. Let $X'=X_1\cap X_2\supset X_0$ be the intersection of the lattices $X_1$ and $X_2$ in $X_0\otimes \bb{Q}$. Let $X$ be the abelian group generated by $X_1$ and $X_2$ in $X_0\otimes \bb{Q}$.\par
We have a map $X_1\oplus_{X'} X_2\to X$ given by addition in $X_0\otimes \bb{Q}$. This is surjective by the definition of $X$ and is injective by the definition of $X'$. Hence, we see that $\mbf{E}:=\spec k[X]\to \spec k[X_1]\times_{\spec k[X']}\spec k[X_2]$ is an isomorphism. On the other hand, the map $\mbf{E}\to \spec k[X_1]\times_{\spec k[X_0]}\spec k[X_2]$ is a closed embedding by the separatedness of $\spec k[X']\to \spec k[X_0]$, and is a morphism between varieties of the same dimension. We then conclude that $\mbf{E}\iso \mbf{E}_3:=(\mbf{E}_1\times_{\mbf{E}_0}\mbf{E}_2)^\circ$.
\end{proof}

\subsubsection{Action on shtukas}\label{subsubsec-action-shtukas}
Now, we work over $E$, the completion of $\bb{E}$ as the beginning of this subsection.
 Let $\PP$ be a quasi-parahoric model of $P$ and $K_p = \PP(\Z_p)$. Fix a neat level $K^p \subset P(\A^p)$, and write $K = K_p K^p$. Let $E_K = \mathbf{E}_K \otimes \rQ_p$, and $E_{\infty,K}' = \varprojlim(\cdots \stackrel{p}{\to} E_K \stackrel{p}{\to} E_K)$. 

 Note that $(E'_{\infty,K})^{\Dia}$ is represented by a perfectoid space in the sense that, when base-changed to $\Spa(C, \OO_C)$, $(E_{\infty, K}')^{\Dia}_C$ is the $\Dia$ of $E_{\infty, K, C}' \sim \varprojlim(\cdots \stackrel{p}{\to} E_{K, C}^{\ad} \stackrel{p}{\to} E_{K, C}^{\ad})$; see \cite[Lem. 2.14]{blakestad2018perfectoid}.
 
 Let $(\PPs, \phi_{\PPs})$ be the $\PP^c$-shtuka on $\shu{K}(P, \mathcal{X})^{\Dia}$ associated with $(\lp_K, \HT_K)$. Let $F$ be a finite field extension of $E$ or $\breve{E}$, and $C$ be a non-archimedean perfectoid field containing $F$. Given $\gamma \in E_K(F)$, by taking $p^{\infty}$-roots in $C$, we can lift $\gamma$ to some $\Tilde{\gamma} \in E_{\infty, K}'(C)$.
\begin{lem}\label{lem: torus action on shtuka, generic fiber}
   $\Tilde{\gamma}$ induces an isomorphism of $\PP^c$-shtukas $(\PPs, \phi_{\PPs})$ and $\gamma^*(\PPs, \phi_{\PPs})$ over $\shu{K}(P, \mathcal{X})_F$.
\end{lem}
\begin{proof}
   In Lemma \ref{lem-torus-action-tower-pro-p-gen}, we see that the $\Tilde{\gamma}$-action on $\lp_{K, C}$ is equivariant with the $\gamma$-action on $\shu{K}(P, \mathcal{X})_C$. This gives an isomorphism of pro-\'etale torsors $\Tilde{\gamma}_*: \lp_{K, C} \to \gamma^*\lp_{K, C}$. We need to show that $\Tilde{\gamma}_*$ is compatible with the Hodge--Tate period maps (i.e. $\HT_{K}(\Tilde{\gamma}(x)) = \HT_K(x)$ for any $x \in \lp_{K, C}$). The action $E_{\infty, K, C}' \times \lp_{K, C} \to \lp_{K, C}$ is the base change of $E_{\infty, K}' \times \lp_{K} \to \lp_{K}$. Since the action of elements in $E_{\infty, K}'$ is compatible with the Hodge--Tate period maps (see the proof of Proposition \ref{prop: Delta acts trivially on generic fiber}), $\Tilde{\gamma}$ is compatible with the Hodge--Tate period maps.
\end{proof}
\begin{rk} 
We do not have a canonical lifting of $\gamma \in E_K$ to $E_{\infty,K}'$, and there is no equivariant action of $E_K$ on $\shu{K}(P, \mathcal{X}) \to \Sht_{\PP^c, \mu^c}$. Thus, the $\ca{P}^c$-shtuka on $\sh_K(P,\ca{X})$ does not descend to $\sh_{\overline{K}}(\overline{P},\overline{\ca{X}})$. The reader can compare this situation with Proposition \ref{prop: Delta acts trivially on generic fiber} later, where the action induces descent data.
\end{rk}
\subsection{Abelian scheme action on the tower}

We focus on the abelian-scheme torsor $\shu{\ovl{K}}:=\shu{\ovl{K}}(\ovl{P}, \ovl{\mathcal{X}}) \to \shu{K_h}:=\shu{K_h}(G_h, \mathcal{X}_h)$. Here $G_h:= P/W = \ovl{P}/V$, $K_h = \pi(\ovl{K})$, and $\pi: \ovl{P} \to G_h$ is the projection.

Recall that the fiber of $\shu{\ovl{K}}(\CC) \to \shu{K_h}(\CC)$ at $(x, q) \in \shu{K_h}(\CC)$ can be computed explicitly (see \cite[\S 3.13]{Pin89}). Let $e: (G_h, X_h) \to (\ovl{P}, \ovl{\mathcal{X}})$ be a splitting. Then the map $[w] \mapsto [w\cdot e(x), q]$ gives an isomorphism from the fiber over $(x, q)$ to $\Gamma_{\ovl{K}}\backslash V(\R)$, where $\Gamma_{\ovl{K}}$ is the image of $(\lrbracket{z \in Z(\ovl{P})(\rQ)|\ z|_{\ovl{\mathcal{X}}} = \identity} \rtimes V(\rQ)) \cap q\ovl{K}q^{-1}$ under the projection $Z(\ovl{P}) \times V \to V$. When we fix the away-from-$p$ part, the embedding $\Gamma_{\ovl{K}(p^n)} \to \Gamma_{\ovl{K}(p^m)}$ ($m \leq n$) has $p$-torsion cokernel and gives the projection
\begin{equation}\label{eq: morphism between abelian variety}
    p_{n, m}:\Gamma_{\ovl{K}(p^n)}\backslash V(\R) \to \Gamma_{\ovl{K}(p^m)}\backslash V(\R).
\end{equation}
By \cite[\S 3.22]{Pin89}, $\Gamma_{\ovl{K}}\backslash V(\R) = \Gamma_{\ovl{K}}\backslash V(\CC)/\exp(\Fil^0(\Lie(V))_{\CC})$ is an abelian variety, and (\ref{eq: morphism between abelian variety}) is a $p$-isogeny of abelian varieties. This picture descends over the reflex field.

Let $K_V$ be the image of $(\lrbracket{z \in Z(\ovl{P})(\rQ)|\ z|_{\ovl{\mathcal{X}}} = \identity} \rtimes V(\A)) \cap \ovl{K}$ under the projection $Z(\ovl{P}) \times V \to V$, and let $K_h^* \subset G_h(\A)$ be an open compact subgroup that contains $K_h$ and normalizes $K_V$; let $\wdt{K}^* = K_V \rtimes K_h^*$. By \cite[Cor. 3.12]{Pin89} and \cite[\S 10]{Pin89}, $\shu{\ovl{K}} \to \shu{K_h}$ is canonically a torsor under the family of abelian varieties $\shu{\wdt{K}^*} \to \shu{K_h^*}$.

To simplify the calculation, we assume that $K_h$ normalizes $K_V$ (and we can, and will, always do this in later sections), and take $K_h^* = K_h$. In particular, for any pair of neat subgroups $\ovl{K}' \subset \ovl{K}$ such that $K_h'$ normalizes $K_V'$, we have a commutative diagram over the reflex field:
\begin{equation}\label{eq: abelian scheme diagram}
\begin{tikzcd}
	{\shu{\wdt{K}^{\prime}} \times_{\shu{K_h^{\prime}}} \shu{\ovl{K}'}} & {\shu{\ovl{K}'}} \\
	{\shu{\wdt{K}} \times_{\shu{K_h}} \shu{\ovl{K}}} & {\shu{\ovl{K}}.}
	\arrow["{\mu_{\ovl{K}'}}", from=1-1, to=1-2]
	\arrow[from=1-1, to=2-1]
	\arrow[from=1-2, to=2-2]
	\arrow["{\mu_{\ovl{K}}}", from=2-1, to=2-2]
\end{tikzcd}
\end{equation}

For our purpose, we directly work over $E$ rather than $\mathbb{E}$. We fix a system of normal subgroups $\ovl{K}' \subset \ovl{K}$ such that $\ovl{K}^{\prime p} = \ovl{K}^p$; then $\wdt{K}^{\prime} \subset \wdt{K}$ are normal and $\wdt{K}^{\prime p} = \wdt{K}^{p}$. Fix $x:=x_{\ovl{K}} \in \shu{K_h}(F)$, where $F$ is a finite field extension of $E$ or $\breve{E}$. Lift $x$ to $\tilde{x} \in \shu{K_h^p}(\ovl{\rQ}_p) := \varprojlim_{K_{h, p}' \subset K_{h, p}} \shu{K_{h, p}'K_h^p}(\ovl{\rQ}_p)$, and let $x_{\ovl{K}'} \in \shu{K_h^{\prime}}(\ovl{\rQ}_p)$ be its image. Denote the abelian scheme $\shu{\wdt{K}'} \to \shu{K_h'}$ by $A_{\ovl{K}'} \to \shu{K_h'}$, the fiber at $x_{\ovl{K}'}$ by $A_{x_{\ovl{K}'}}$, and $A_{x_{\ovl{K}}}$ by $A_x$. Also denote the fiber of $\shu{\ovl{K}'} \to \shu{K_h'}$ at $x_{\ovl{K}'}$ by $\shu{\ovl{K}', x_{\ovl{K}'}}$, and $\shu{\ovl{K}, x_{\ovl{K}}}$ by $\shu{\ovl{K}, x}$. The projection $A_{\ovl{K}'} \to A_{\ovl{K}}$ induces a $p$-isogeny $A_{x_{\ovl{K}'}} \to A_{x}$, by (\ref{eq: morphism between abelian variety}).

Passing to the $p^{\infty}$-level, denote
\[ \mathcal{S}_{\wdt{K}^p} = \varprojlim_{\wdt{K}_p' \subset \wdt{K}_p} \shu{\wdt{K}_p'\wdt{K}^p},\quad \ab_{\ovl{K}^p} = \varprojlim_{\ovl{K}_p' \subset \ovl{K}_p} A_{\ovl{K}_p'\ovl{K}^p},\quad \mathcal{S}_{\wdt{K}^p} = \ab_{\ovl{K}^p}. \]
\[ \mathcal{S}_{\ovl{K}^p} = \varprojlim_{\ovl{K}_p' \subset \ovl{K}_p} \shu{\ovl{K}_p'\ovl{K}^p},\quad \mathcal{S}_{K_h^p} = \varprojlim_{K_{h, p}' \subset K_{h, p}} \shu{K_{h, p}'K_h^p},\quad \ab_{\infty, \tilde{x}} =  \varprojlim_{\ovl{K}_p' \subset \ovl{K}_p} A_{x_{\ovl{K}'_p\ovl{K}^{p}}}. \]

\begin{rk}
    In abelian-type case, $\mathcal{S}_{K_h^p}^{\Dia}$ is represented by a perfectoid space (see \cite{shen2017perfectoid}), one can show that $\mathcal{S}_{\wdt{K}^p}^{\Dia}$ and $\mathcal{S}_{\ovl{K}^p}^{\Dia}$ are represented by perfectoid spaces using the fact that the fibers $\ab_{\infty, \tilde{x}}^{\Dia}$ are represented by perfectoid spaces (see \cite[Thm. 1]{blakestad2018perfectoid}). We do not need this result, so we do not spell out the proof.
\end{rk}

Let $\tilde{x}^c$ be the image of $\tilde{x}$ in $\shu{K_h^{c, p}}(\ovl{\rQ}_p) := \varprojlim_{K_{h, p}^{c, \prime} \subset K_{h, p}^c} \shu{K_{h, p}^{c, \prime}K_h^{c, p}}(\ovl{\rQ}_p)$. We can similarly define these objects using $\shu{\ovl{K}^c}(\ovl{P}^c, \ovl{\mathcal{X}}^c) \to \shu{K_h^c}(G_h^c, X_h^c)$. Then diagram (\ref{eq: abelian scheme diagram}) induces a commutative diagram:
\begin{equation}\label{eq: abelian scheme diagram, 2}
\begin{tikzcd}[sep=tiny]
	{A_{x_{\ovl{K}'}} \times_{x_{\ovl{K}'}} \shu{\ovl{K}', x_{\ovl{K}'}}} && {\shu{\ovl{K}', x_{\ovl{K}'}}} & \\
	& {A_{x_{\ovl{K}^{c, \prime}}} \times_{x_{\ovl{K}^{c, \prime}}} \shu{\ovl{K}^{c, \prime}, x_{\ovl{K}^{c, \prime}}}} && {\shu{\ovl{K}^{c, \prime}, x_{\ovl{K}^{c, \prime}}}} \\
	{A_{x} \times_{x} \shu{\ovl{K}, x}} && {\shu{\ovl{K}, x}} \\
	& {A_{x^c} \times_{x^c} \shu{\ovl{K}^{c}, x^c}} && { \shu{\ovl{K}^{c}, x^{c}}.}
	\arrow["{\mu_{x_{\ovl{K}'}}}", from=1-1, to=1-3]
	\arrow[from=1-1, to=2-2]
	\arrow[from=1-1, to=3-1]
	\arrow[from=1-3, to=2-4]
	\arrow[from=1-3, to=3-3]
	\arrow["{\mu_{x_{\ovl{K}^{c, \prime}}}}", from=2-2, to=2-4]
	\arrow[from=2-2, to=4-2]
	\arrow[from=2-4, to=4-4]
	\arrow["{\mu_{x}}", from=3-1, to=3-3]
	\arrow[from=3-1, to=4-2]
	\arrow[from=3-3, to=4-4]
	\arrow["{\mu_{x^c}}", from=4-2, to=4-4]
\end{tikzcd}
\end{equation}
By passing to $p^{\infty}$-level, we have a commutative diagram:
\begin{equation}\label{eq: abelian scheme diagram, 3}
\begin{tikzcd}[sep=tiny]
	{\mathcal{A}_{\infty, \tilde{x}} \times_{\tilde{x}} \mathcal{S}_{\ovl{K}^p, \tilde{x}}} && { \mathcal{S}_{\ovl{K}^p, \tilde{x}}} & \\
	& {\mathcal{A}_{\infty, \tilde{x}^c} \times_{\tilde{x}^c} \mathcal{S}_{\ovl{K}^{c, p}, \tilde{x}^c}} && {\mathcal{S}_{\ovl{K}^{c, p}, \tilde{x}^c}} \\
	{A_{x} \times_{x} \shu{\ovl{K}, x}} && {\shu{\ovl{K}, x}} \\
	& {A_{x^c} \times_{x^c} \shu{\ovl{K}^{c}, x^c}} && { \shu{\ovl{K}^{c}, x^{c}},}
	\arrow["{\mu_{\tilde{x}}}", from=1-1, to=1-3]
	\arrow[from=1-1, to=2-2]
	\arrow[from=1-1, to=3-1]
	\arrow[from=1-3, to=2-4]
	\arrow[from=1-3, to=3-3]
	\arrow["{\mu_{\tilde{x}^c}}", from=2-2, to=2-4]
	\arrow[from=2-2, to=4-2]
	\arrow[from=2-4, to=4-4]
	\arrow["{\mu_{x}}", from=3-1, to=3-3]
	\arrow[from=3-1, to=4-2]
	\arrow[from=3-3, to=4-4]
	\arrow["{\mu_{x^c}}", from=4-2, to=4-4]
\end{tikzcd}
\end{equation}
where $\mathcal{S}_{\ovl{K}^p, \tilde{x}}$ is the fiber of $\mathcal{S}_{\ovl{K}^p} \to \mathcal{S}_{K_h^p}$ at $\tilde{x}$. Similarly for $\mathcal{S}_{\ovl{K}^{c, p}, \tilde{x}^c}$.

We introduce lemmas analogous to Lemma \ref{lem-torus-action-tower-pro-p-gen} and Lemma \ref{lem-torus-isog-p-ele}.
\begin{lem}\label{lem-ab-action-tower-pro-p-gen}
    With the above conventions, the pushout
    \[ \ovl{\PPp}_K := \varprojlim_{\ovl{K}'_p \subset \ovl{K}_p} \shu{\ovl{K}}(\ovl{P}, \ovl{\mathcal{X}}) \times^{\underline{\ovl{K}/Z(\ovl{P})(\rQ)^{-}_{\ovl{K}}\ovl{K}^p}} \underline{\ovl{K}^c_p} \]
    is the pullback of $\ovl{\PPp}_{K^c}$ on $\shu{\ovl{K}^c}(\ovl{P}^c, \ovl{\mathcal{X}}^c)$. The action $\mu_{\tilde{x}}$ induces a $(\ab_{\infty, \tilde{x}}' \to A_x)$-equivariant action $\mu_{\tilde{x}}'$ on $\ovl{\PPp}_{K, \tilde{x}} \to \shu{\ovl{K}}(\ovl{P}, \ovl{\mathcal{X}})_x$, where $\ab_{\infty, \tilde{x}}'$ is an inverse limit of $p$-power covers of $\ab_x$, and $\ovl{\PPp}_{K, \tilde{x}}$ is the fiber of $\ovl{\PPp}_{K} \to \PPp_{K_h}$ at $\tilde{x}$.
\end{lem}
\begin{proof}
    Same as the proof of Lemma \ref{lem-torus-action-tower-pro-p-gen} (using Lemma \ref{lem-ab-isog-p-ele}).
\end{proof}
\begin{lem}\label{lem-ab-isog-p-ele}
    Let $A_0, A_1, A_2$ be abelian varieties over a field of characteristic $0$, and let $f: A_1 \to A_0$ and $g: A_2 \to A_0$ be isogenies. Assume that $g$ is a $p$-power isogeny. Then the identity component $A_3$ of $A_1 \times_{f, A_0, g} A_2$ is an abelian variety equipped with a $p$-power isogeny $A_3 \to A_1$.
\end{lem}
\begin{proof}
    $A_1 \times_{f, A_0, g} A_2$ is a closed subgroup of $A_1 \times A_2$, hence projective; its identity component $A_3$ is therefore an abelian variety by Cartier's theorem since all varieties are defined over a field of characteristic zero. By construction, $A_3 \to A_1$ is a $p$-power isogeny.
\end{proof}

Let $(\ovl{\PPs}, \phi_{\ovl{\PPs}})$ be the $\ovl{\PP}^c$-shtuka on $\shu{\ovl{K}}^{\Dia}$ associated with $(\ovl{\lp}_{K}, \ovl{\HT}_{K})$. 
\begin{lem}\label{lem: abelian scheme action on shtuka, generic fiber}
   Let $x \in \shu{K_h}(F)$. Given $\gamma \in A_{x}(F)$, there is an isomorphism of $\ovl{\PP}^c$-shtukas $(\ovl{\PPs}, \phi_{\ovl{\PPs}})$ and $\gamma^*(\ovl{\PPs}, \phi_{\ovl{\PPs}})$ over $\shu{\ovl{K}, x}^{\Dia}$.
\end{lem}
\begin{proof}
  This follows similarly from the proof of Lemma \ref{lem: torus action on shtuka, generic fiber}. Let $\shu{\ovl{K}', x}$ (resp. $A_{\ovl{K}', x}'$) be the fiber of $x$ along $\shu{\ovl{K}'} \to \shu{K_h}$ (resp. $A_{\ovl{K}'}' \to \shu{K_h'}$; here we use the identity component of $A_{\ovl{K}} \times_{A_{\ovl{K}^c}} A_{\ovl{K}^{c, \prime}}$, where we pull back these abelian schemes over $\shu{K_h'}$). Let $S:=\mathcal{S}_{K_h^p, x}$. Let $\ab_{\infty, x}' = \varprojlim_{\ovl{K}_p' \subset \ovl{K}_p} A'_{\ovl{K}'_p\ovl{K}^{p}, x}$; it is a group object over $S$. Given any $x' \in S$, the fiber of $\ab_{\infty, x}'$ at $x'$ is $\ab_{\infty, \tilde{x}'}'$ for some $\tilde{x}' \in \mathcal{S}_{\ovl{K}^p}$ lifting $x'$. Consider the projection $\ab_{\infty, x}' \to A_x \times_x S$; this is a morphism of group objects over $S$. Let $T_p(A_{S})$ be its kernel; it is a $v$-sheaf over $S$. Recall that for any $x' \in S$, the fiber $T_p(A_S)_{x'}$ is the kernel of $\ab'_{\infty, \tilde{x}'} \to A_x$, which is isomorphic to $T_p(A_x)$. In other words, $T_p(A_{S})$ is a torsor under $T_p(A_x)_S$.

  Given $\gamma \in A_x(F)$, it has image $x \in \shu{K_h}(F)$. We lift $\gamma$ to $\gamma_S \in A_x(S)$. Take a $v$-cover $T \to S$ that trivializes $T_p(A_S)$; we can lift $\gamma_S$ to $\Tilde{\gamma}_T \in (\ab'_{\infty, x})_T$ and obtain an action $(\ab'_{\infty, x})_T \times_T (\lp_{\ovl{K}, x})_T \to (\lp_{\ovl{K}, x})_T$ lifting the group action $A_x \times_{x} \shu{\ovl{K}, x} \to \shu{\ovl{K}, x}$. In particular, the $\Tilde{\gamma}_T$-action on $(\lp_{\ovl{K}, x})_T$ lifts the $\gamma$-action on $\shu{\ovl{K}, x}$. Moreover, the action of $(\ab'_{\infty, x})_T$ on $(\lp_{\ovl{K}, x})_T$ essentially comes from the action of $\ab'_{\infty, x}$ on $\lp_{\ovl{K}, x}$ over $S$; the $\Tilde{\gamma}_T$-action is a global lifting of $\gamma$ and is compatible with the Hodge--Tate period map, thus it induces an isomorphism of $\ovl{\PP}^c$-shtukas $(\ovl{\PPs}, \phi_{\ovl{\PPs}})$ and $\gamma^*(\ovl{\PPs}, \phi_{\ovl{\PPs}})$ over $\shu{\ovl{K}, x}^{\Dia}$.

\end{proof}

\section{Canonicity and functoriality}\label{sec-canonical}
\subsection{Canonical integral models}\label{subsec: canonical integral models}

Mimicking the axioms in \cite{PR24}, we define canonical integral models for mixed Shimura varieties coming from the boundary.

Let $(G, X)$ be a pure Shimura datum, let $\GG$ be a quasi-parahoric group scheme, and let $K_p = \GG(\Z_p)$. We fix a cusp label $[\Phi] = [(Q_{\Phi}, X_{\Phi}^+, g_{\Phi})] \in \Cusp_K(G, X)$. Here we can vary $K^p$, and do not distinguish $\Phi$ for different levels $K^p$ once we have prescribed $g = g_{\Phi} \in G(\A)$. Recall that, given $\GG$, we can attach quasi-parahoric group schemes $\PP_{\Phi}$ (resp. $\PP_{\Phi}^c$, $\PP_{\Phi}^*$) to $P_{\Phi}$ (resp. $P_{\Phi}^c$, $P_{\Phi}^*$) as in Section \ref{subsec: pass from P to P^c}.

Recall the Pappas--Rapoport integral models:
\begin{axiom}[{\cite[Conj. 4.2.2]{PR24}, \cite[Def. 4.1.2]{daniels2024conjecture}, \cite[Def. 4.3]{DY25}}]\label{def: canonical model for pure Shimura data}
     Consider a system $\lrbracket{\Shum{K}(G, X)}_{K^p}$ of normal flat schemes over $\OO_E$ with generic fiber $\lrbracket{\shu{K}(G, X)}_{K^p}$, with $K^p$ varying over all sufficiently small compact open subgroups of $G(\A^p)$. We say $\lrbracket{\Shum{K}(G, X)}_{K^p}$ is a canonical integral model of $\lrbracket{\shu{K}(G, X)}_{K^p}$ if the following properties are satisfied.
    \begin{enumerate}
        \item For every discrete valuation ring $R$ of mixed characteristic over $\OO_E$, we have
        \begin{equation*}
            \shu{K}(G, X)(R[1/p]) = (\varprojlim_{K^p} \Shum{K}(G, X))(R).
        \end{equation*}
        \item For every $K^p \subset G(\A^p)$ and $K^{\prime p} \subset G(\A^p)$ and an $h \in G(\A^p)$ with $h^{-1}K^{\prime p}h \subset K^p$, there are finite \'etale morphisms 
        \[ t_{K', K}(h): \Shum{K'}(G, X) \to \Shum{K}(G, X) \]
        extending the generic fiber.
        \item The $\GG^c$-shtuka $\PPs_{K, E}$ extends to a $\GG^c$-shtuka $\PPs_{K}$ on $\Shum{K}(G, X)^{\Dia/}$ for every sufficiently small $K^p$.
        \item Consider $x \in \Shum{K}(G, X)(k)$, with corresponding $b_x \in G^c(\bQ)$. Let $x_0$ be the base point in $\MM^{\intg}_{\GG^c, b_x, \mu^c}(k)$. Then there exists an isomorphism of $v$-sheaves:
        \begin{equation*}
            \Theta_x: (\MM^{\intg}_{\GG^c, b_x, \mu^c})^{\wedge}_{/x_0} \rightiso (\Shum{K}(G, X)^{\wedge}_{/x})^{\Dia},
        \end{equation*}
        such that $\Theta_x^*(\PPs_{K})$ is the universal $\GG^c$-shtuka on $\MM^{\intg}_{\GG^c, b_x, \mu^c}$.
    \end{enumerate}
\end{axiom}

We mimic the above definition:
\begin{axiom}\label{def: canonical model for mixed Shimura data}
    Consider a system $\lrbracket{\Shum{K_{\Phi}}(P_{\Phi}, D_{\Phi})}_{K_{\Phi}^p}$ of normal flat schemes over $\OO_E$ with generic fiber $\lrbracket{\shu{K_{\Phi}}(P_{\Phi}, D_{\Phi})}_{K_{\Phi}^p}$, with $K_{\Phi}^p$ varying over all sufficiently small compact open subgroups of $P_{\Phi}(\A^p)$. We say $\lrbracket{\Shum{K_{\Phi}}(P_{\Phi}, D_{\Phi})}_{K_{\Phi}^p}$ is a canonical integral model of $\lrbracket{\shu{K_{\Phi}}(P_{\Phi}, D_{\Phi})}_{K_{\Phi}^p}$ if the following properties are satisfied.
    \begin{enumerate}
        \item For every discrete valuation ring $R$ of mixed characteristic over $\OO_E$, we have
        \begin{equation*}
            \shu{K_{\Phi, p}}(P_{\Phi}, D_{\Phi})(R[1/p]) = (\varprojlim_{K_{\Phi}^p} \Shum{K_{\Phi, p}K_{\Phi}^p}(P_{\Phi}, D_{\Phi}))(R).
        \end{equation*}
        \item For every $K_{\Phi}^p \subset P_{\Phi}(\A^p)$ and $K_{\Phi}^{\prime p} \subset P_{\Phi}(\A^p)$ and an $h \in P_{\Phi}(\A^p)$ with $h^{-1}K_{\Phi}^{\prime p}h \subset K_{\Phi}^p$, there are finite \'etale morphisms 
        \[ t_{K'_{\Phi}, K_{\Phi}}(h): \Shum{K_{\Phi}'}(P_{\Phi}, D_{\Phi}) \to \Shum{K_{\Phi}}(P_{\Phi}, D_{\Phi}) \]
        extending the generic fiber.
        \item The $\PP_{\Phi}^*$-shtuka $\PPs_{K_{\Phi}, E}$ extends to a $\PP_{\Phi}^*$-shtuka $\PPs_{K_{\Phi}}$ on $\Shum{K_{\Phi}}(P_{\Phi}, D_{\Phi})^{\Dia/}$ for every sufficiently small $K_{\Phi}^p$.
        \item Consider $x \in \Shum{K_{\Phi}}(P_{\Phi}, D_{\Phi})(k)$, with corresponding $b_{\Phi, x} \in P_{\Phi}^*(\bQ) \hookrightarrow G^c(\bQ)$ (see the first paragraph of subsection \ref{subsec: Newton strata and central leaves}). Let $x_0$ be the base point in $\MM^{\intg}_{\GG^c_{\Phi}, b_{\Phi, x}, \mu^c_{\Phi}}(k)$. Then there exists an isomorphism of $v$-sheaves:
        \begin{equation*}
            \Theta_x: (\MM^{\intg}_{\GG^c_{\Phi}, b_{\Phi, x}, \mu^c_{\Phi}})^{\wedge}_{/x_0} \rightiso (\Shum{K_{\Phi}}(P_{\Phi}, D_{\Phi})^{\wedge}_{/x})^{\Dia},
        \end{equation*}
        such that $\GG_{\Phi}^c \times^{\PP_{\Phi}^*} \Theta_x^*(\PPs_{K_{\Phi}})$ is the universal $\GG^c_{\Phi}$-shtuka on $\MM^{\intg}_{\GG^c_{\Phi}, b_{\Phi, x}, \mu^c_{\Phi}}$.
    \end{enumerate}

    In general, let $Y_\Phi$ be a connected subgroup of $G$ satisfying $P_\Phi\sbst Y_\Phi\sbst ZP_\Phi$, define $Y_{\Phi}^*$ as in Definition \ref{def-pstar}, we set up the axioms for the canonical integral model $\lrbracket{\Shum{K_{\Phi}^Y}(Y_{\Phi}, D_{\Phi}^Y)}_{K_{\Phi}^{Y, p}}$ of $\lrbracket{\shu{K_{\Phi}^Y}(Y_{\Phi}, D_{\Phi}^Y)}_{K_{\Phi}^{Y, p}}$ in exactly same way.
\end{axiom}
\begin{rk}
     We use the convention that when $(P_{\Phi}, D_{\Phi}) = (G, X)$, $P_{\Phi}^c = P_{\Phi}^* = G^c$ and $\PP_{\Phi}^c = \PP_{\Phi}^* = \GG^c$, so Axiom \ref{def: canonical model for mixed Shimura data} generalizes Axiom \ref{def: canonical model for pure Shimura data}.
\end{rk}
\begin{rk}\label{rk: why use star not c}
    We can consider the following condition $(3)'$ (resp. $(4)'$) compared with $(3)$ (resp. $(4)$): in $(3)$ (resp. $(4)$), replace $\PP_{\Phi}^*$ with $\PP_{\Phi}^c$ in the statement.
    
    By devissage \ref{lem: Devissage}, since $\PP^c_{\Phi} \to \GG^c_{\Phi}$ factors through $\PP^*_{\Phi} \hookrightarrow \GG^c_{\Phi}$, $(4)$ and $(4)'$ are equivalent. However, $(3)'$ is stronger than $(3)$. In the rest of the article, we use $(3)$ instead of $(3)'$, mainly due to two reasons:
    \begin{enumerate}
        \item The $\Delta_{\Phi, K}^{\circ}$-action on $\shu{K_{\Phi}}(P_{\Phi}, D_{\Phi})$ fixes the $\PP^*_{\Phi}$-shtuka rather than the $\PP^c_{\Phi}$-shtuka,
        \item $P_{\Phi}^c \to G^c$ might not be an embedding, in practice, we need an embedding $P_{\Phi}^* \to G^c$ to apply the funtoriality result Proposition \ref{prop: functoriality of canonical integral models} (which is needed in proving the canonicity of integral models of abelian-type, see Theorem \ref{thm-ext-cim-ab}).
    \end{enumerate}
\end{rk}

In order to prove functoriality for such canonical integral models, we need further assumptions when $G^c \neq G$.

\begin{definition}\label{def: well-adapted}
    Let $(P, \mathcal{X})$ be a Shimura datum, let $Z'$ be a central multiplicative group in $P$ that contains $Z(P)_{ac}$, and let $(P', \mathcal{X}'):=(P, \mathcal{X})/Z'$. Let $\PP$ be a quasi-parahoric group scheme of $P$, and let $\PP'$ be a quasi-parahoric group scheme of $P'$ defined in Definition \ref{def-ass-gp-sch-mix-sh}. Let $K_p = \PP(\Z_p)$ and $K_p' = \PP'(\Z_p)$. Consider a system $\lrbracket{\Shum{K}(P, \mathcal{X})}_{K^p}$ (resp. $\lrbracket{\Shum{K'}(P', \mathcal{X}')}_{K^{\prime, p}}$) of normal flat schemes over $\OO_E$ with generic fiber $\lrbracket{\shu{K}(P, \mathcal{X})}_{K^p}$ (resp. $\lrbracket{\shu{K'}(P', \mathcal{X}')}_{K^{\prime, p}}$), with $K^p$ (resp. $K^{\prime, p}$) varying over all sufficiently small compact open subgroups of $P(\A^p)$ (resp. $P'(\A^p)$).

    Let $\pi: (P, \mathcal{X}) \to (P', \mathcal{X}')$ be the projection. We say $\lrbracket{\Shum{K}(P, \mathcal{X})}_{K^p}$ and $\lrbracket{\Shum{K'}(P', \mathcal{X}')}_{K^{\prime, p}}$ are \textbf{adapted} if, for any pair $(K^p, K^{\prime, p})$ such that $\pi(K^p) \subset K^{\prime, p}$, there is a finite morphism
    \begin{equation}\label{eq: proj K to K'}
        \Shum{K}(P, \mathcal{X}) \to  \Shum{K'}(P', \mathcal{X}')
    \end{equation}
    extending the generic fiber $\shu{K}(P, \mathcal{X}) \to \shu{K'}(P', \mathcal{X}')$ induced by $\pi$.

    We say $\lrbracket{\Shum{K}(P, \mathcal{X})}_{K^p}$ is \textbf{adapted with respect to $P \to P'$} if there exists $\lrbracket{\Shum{K'}(P', \mathcal{X}')}_{K^{\prime, p}}$ such that $\lrbracket{\Shum{K}(P, \mathcal{X})}_{K^p}$ and $\lrbracket{\Shum{K'}(P', \mathcal{X}')}_{K^{\prime, p}}$ are adapted.
\end{definition}
\begin{definition}\label{def: well-adapted, 2}
    Keep notation in Definition \ref{def: well-adapted}. Let $(G, X):=(P, \mathcal{X})/W$. We say a system of normal flat models $\lrbracket{\Shum{K}(P, \mathcal{X})}_{K^p}$ of $\lrbracket{\shu{K}(P, \mathcal{X})}_{K^p}$ adapted with respect to $P \to P'$ is moreover \textbf{adapted with respect to $P \to G$} if there exists a system of normal flat models $\lrbracket{\Shum{K_G}(G, X)}_{K^{p}_G}$ of $\lrbracket{\shu{K_G}(G, X)}_{K^{p}_G}$ that is adapted with $G \to G'$, where $G' = G/Z'$ (by Lemma \ref{lem-levi-c}, $Z(G)_{ac} = Z(P)_{ac} \subset Z'$) and $\PP' = \UU \rtimes \GG'$, such that when $K^p_G$ and $K^{\prime, p}_G$ contain the image of $K^p$ and $K^{\prime, p}$ respectively, the induced morphisms on generic fibers
    \[ \shu{K}(P, \mathcal{X}) \to \shu{K_G}(G, X),\quad \shu{K'}(P', \mathcal{X}') \to \shu{K_G'}(G',  X') \]
    extend to morphisms on integral models
    \[ \Shum{K}(P, \mathcal{X}) \to \Shum{K_G}(G, X),\quad \Shum{K'}(P', \mathcal{X}') \to \Shum{K_G'}(G', X'), \]
    respectively.
\end{definition}

 Recall the following fact:
\begin{lem}
    Given a morphism between Shimura data $\iota: (P_1, \mathcal{X}_1) \to (P_2, \mathcal{X}_2)$, then $\iota(W_1) \subset W_2$ and $\iota(U_1) \subset U_2$, and $\iota$ induces $G_1 \to G_2$ and $V_1 \to V_2$. In particular, the morphism $P_1 \to P_2$ is compatible with $U_1 \rtimes G_1 \to U_2 \rtimes G_2$ in the sense of \ref{def: compatible, P and G}.
\end{lem}
\begin{proof}
    Let $h_1 \in \mathcal{X}_1$, and set $h_2 = h_1\circ\iota \in \mathcal{X}_2$. By definition, the weight filtration on $\Lie P_i$ induced by $\Ad_{P_i}\circ h_i$ is $W_{-1}(\Lie P_i) = \Lie W_i$, $W_{-2}(\Lie P_i) = \Lie U_i$. By construction, $\Lie P_1 \to \Lie P_2$ is a morphism between mixed Hodge structures, and it is strict; thus it maps $\Lie W_1$ to $\Lie W_2$ and $\Lie U_1$ to $\Lie U_2$.
\end{proof}

\begin{lem}\label{lem: two cases}
   Consider one of the following two cases:
   \begin{enumerate}
       \item Let $\iota: (P_1, \mathcal{X}_1) \to (P_2, \mathcal{X}_2)$ be an embedding of mixed Shimura data, let $(G_i, X_i) = (P_i, \mathcal{X}_i)/W$ be the induced pure Shimura datum. Let $E_1$, $E_2$ be the completion of reflex fields. Let $\PP_1, \PP_2$ be quasi-parahoric group schemes of $P_1$ and $P_2$ respectively, such that
       \begin{equation}\label{eq: cpt groups P induces G}
           \PP_2(\bZ_p) \cap P_1(\bQ) = \PP_1(\bZ_p),\quad \GG_2(\bZ_p) \cap G_1(\bQ) = \GG_1(\bZ_p).
       \end{equation}
       Let $K_2^p \subset P_2(\A^p)$ be a sufficiently small subgroup. One can always find $K_1^p \subset P_1(\A^p)$ such that the induced morphisms $\shu{K_1}(P_1, \mathcal{X}_1) \to \shu{K_2}(P_2, \mathcal{X}_2) \otimes_{E_2} E_1$ and $\shu{K_{G_1}}(G_1, X_1) \to \shu{K_{G_2}}(G_2, X_2) \otimes_{E_2} E_1$ are closed embeddings (such levels exist; see \cite[Lem. 1.22]{Wu25}), where $K_i = \PP_i(\Z_p)K^p_i$, and $K_{G_i} = \GG_i(\Z_p)K_{G_i}^p$, where $K_{G_i}^p \subset G_i(\A^p)$ are the images of $K_i^p$.

       Let $\PP_i'$ be a quasi-parahoric group scheme of $P_i'$ as in Definition \ref{def-ass-gp-sch-mix-sh}, such that $P_1' \to P_2'$ induces $\PP_1' \to \PP_2'$. We \emph{assume} $P_1' \to P_2'$ is a closed embedding and moreover the smoothing $\wdt{\PP}_1'$ of the closure of $P_1'$ in $\PP_2'$ is a quasi-parahoric group scheme with $(\wdt{\PP}_1')^{\circ} = \PP_1^{\prime\circ}$.
       \item Let $(P, \mathcal{X}):=(P_1, \mathcal{X}_1) = (P_2, \mathcal{X}_2)$ be a mixed Shimura datum and $(G, X) = (G_i, X_i) = (P_i, \mathcal{X}_i)/W_i$ be the induced pure Shimura datum, let $E_1 = E_2$ be the completion of the reflex field. Let $\PP_1 \to \PP_2$ be quasi-parahoric group schemes of $P$ such that $\UU_1 = \UU_2$ and $\PP_1^{\circ} = \PP_2^{\circ}$ (i.e. $\GG_1^{\circ} = \GG_2^{\circ}$). Define $\PP_1'$, $\PP_2'$ as in Definition \ref{def-ass-gp-sch-mix-sh} such that we have a morphism $\PP_1' \to \PP_2'$ (which automatically induces $(\PP_1')^{\circ} = (\PP_2')^{\circ}$ by construction). Let $K_i = \PP_i(\Z_p)K^p$ for some neat $K^p \subset P(\A^p)$, and $K_{G_i} = \GG_i(\Z_p)K_G^p$, where $K_G^p \subset G(\A^p)$ are the images of $K^p$. Define $K_i'$, $K_{G_i}'$ similarly, for some $K^{\prime, p} \subset P'(\A^p)$.
   \end{enumerate}
    
    Assume
    \[ \lrbracket{\shu{K_2}(P_2, \mathcal{X}_2)}_{K_2^p},\quad \lrbracket{\shu{K_{G_2}}(G_2, X_2)}_{K_{G_2}^p},\quad \lrbracket{\shu{K_2'}(P'_2, \mathcal{X}_2')}_{K_2^{\prime, p}},\quad \lrbracket{\shu{K_{G_2}'}(G'_2, X_2')}_{K_{G_2}^{\prime, p}} \]
    have integral models
    \begin{equation}\label{eq: 123}
        \lrbracket{\Shum{K_2}(P_2, \mathcal{X}_2)}_{K_2^p},\quad \lrbracket{\Shum{K_{G_2}}(G_2, X_2)}_{K_{G_2}^p},\quad \lrbracket{\Shum{K_2'}(P'_2, \mathcal{X}_2')}_{K_2^{\prime, p}},\quad \lrbracket{\Shum{K_{G_2}'}(G'_2, X_2')}_{K_{G_2}^{\prime, p}}
    \end{equation}
    respectively that are normal and flat over $\OO_{E_2}$ and are adapted with $P_2 \to P_2'$ and $P_2 \to G_2$ in the sense of Definitions \ref{def: well-adapted} and \ref{def: well-adapted, 2}, such that the shtukas on generic fibers extend over these integral models, i.e., for $(P_2, \mathcal{X}_2)$ (and similarly for $(P_2', \mathcal{X}_2')$), given the commutative diagram on the right, we have horizontal morphisms on the left diagram that make the diagram commute:
    \[
\begin{tikzcd}
	{\Shum{K_2}(P_2, \mathcal{X}_2)^{\Dia/}} & {\Sht_{\PP_2', \mu', \delta = 1}} & {\shu{K_2}(P_2, \mathcal{X}_2)^{\Dia}} & {\Sht_{\PP_2', \mu', \delta = 1, \Spd E}} \\
	{\Shum{K_{G_2}}(G_2, X_2)^{\Dia/}} & {\Sht_{\GG_2', \mu', \delta = 1},} & {\shu{K_{G_2}}(G_2, X_2)^{\Dia}} & {\Sht_{\GG_2', \mu', \delta = 1, \Spd E}.}
	\arrow[from=1-1, to=1-2]
	\arrow[from=1-1, to=2-1]
	\arrow[from=1-2, to=2-2]
	\arrow[from=1-3, to=1-4]
	\arrow[from=1-3, to=2-3]
	\arrow[from=1-4, to=2-4]
	\arrow[from=2-1, to=2-2]
	\arrow[from=2-3, to=2-4]
\end{tikzcd}
    \]
    Let
    \[ \lrbracket{\Shum{K_1}(P_1, \mathcal{X}_1)}_{K_1^p},\quad \lrbracket{\Shum{K_{G_1}}(G_1, X_1)}_{K_{G_1}^p},\quad \lrbracket{\Shum{K_1'}(P'_1, \mathcal{X}_1')}_{K_1^{\prime, p}},\quad \lrbracket{\Shum{K_{G_1}'}(G'_1, X_1')}_{K_{G_1}^{\prime, p}} \] 
    be the relative normalizations of ((\ref{eq: 123})$\otimes \OO_{E_1}$) in
    \[ \lrbracket{\shu{K_1}(P_1, \mathcal{X}_1)}_{K_1^p},\quad \lrbracket{\shu{K_{G_1}}(G_1, X_1)}_{K_{G_1}^p},\quad \lrbracket{\shu{K_1'}(P'_1, \mathcal{X}_1')}_{K_1^{\prime, p}},\quad \lrbracket{\shu{K_{G_1}'}(G'_1, X_1')}_{K_{G_1}^{\prime, p}}, \]
    then these integral models are normal and flat over $\OO_{E_1}$ and are adapted with $P_1 \to P_1'$ and $P_1 \to G_1$, such that the shtukas on generic fibers extend over these integral models, i.e., for $(P_1, \mathcal{X}_1)$ (and similarly for $(P_1', \mathcal{X}_1')$), given the commutative diagram on the right, we have horizontal morphisms on the left diagram that make the diagram commute:
\[
\begin{tikzcd}
	{\Shum{K_1}(P_1, \mathcal{X}_1)^{\Dia/}} & {\Sht_{\PP_1', \mu', \delta = 1}} & {\shu{K_1}(P_1, \mathcal{X}_1)^{\Dia}} & {\Sht_{\PP_1', \mu', \delta = 1, \Spd E}} \\
	{\Shum{K_{G_1}}(G_1, X_1)^{\Dia/}} & {\Sht_{\GG_1', \mu', \delta = 1},} & {\shu{K_{G_1}}(G_1, X_1)^{\Dia}} & {\Sht_{\GG_1', \mu', \delta = 1, \Spd E}.}
	\arrow[from=1-1, to=1-2]
	\arrow[from=1-1, to=2-1]
	\arrow[from=1-2, to=2-2]
	\arrow[from=1-3, to=1-4]
	\arrow[from=1-3, to=2-3]
	\arrow[from=1-4, to=2-4]
	\arrow[from=2-1, to=2-2]
	\arrow[from=2-3, to=2-4]
\end{tikzcd}
\]
\end{lem}
\begin{proof}
    We prove these two cases separately.
    \begin{enumerate}
    \item Since $P_1' \to P_2'$ is a closed embedding, $G_1' \to G_2'$ is also a closed embedding. Let $\wdt{\GG}_1'$ be the smoothing of the closure of $G_1'$ in $\GG_2'$. We claim that $\wdt{\GG}_1'$ is a quasi-parahoric group scheme and $(\wdt{\GG}_1')^{\circ} = (\GG_1')^{\circ}$: condition (\ref{eq: cpt groups P induces G}) shows that $\PP_1 \to \PP_2$ is compatible with $\GG_1 \to \GG_2$ by Lemma \ref{lem: induced compatibility on quasi-group scheme}; in particular, we have compatible sections $G_1 \to P_1$ and $G_2 \to P_2$ by taking generic fibers. Consider the induced compatible sections $G_1' \to P_1'$ and $G_2' \to P_2'$, $\PP_1' \to \PP_2'$ is compatible with $\GG_1' \to \GG_2'$ under such sections. Then the closure of $G_1'$ in $\GG_2'$ is the same as the closure of $G_1'$ in $\PP_2'$, hence also the closure of $G_1'$ in $\wdt{\PP}_1'$. Since $(\wdt{\PP}_1')^{\circ} = \PP_1^{\prime\circ}$, we have $\wdt{\PP}_1' = \UU_1 \rtimes \wdt{\GG}_1'$, and $(\wdt{\GG}_1')^{\circ} = (\GG_1')^{\circ}$.
    
    Let $\wdt{K}_{G_1, p}' = \wdt{\GG}_1'(\Z_p)$ and $\wdt{K}_{G_1}' = \wdt{K}_{G_1, p}'K_{G_1}^p$. Then $\shu{\wdt{K}_{G_1}'}(G_1', X_1')^{\Dia} \to \Sht_{\wdt{\GG}_1', \mu', \delta = 1, \Spd E}$ extends to $\Shum{\wdt{K}_{G_1}'}(G'_1, X_1')^{\Dia/} \to \Sht_{\wdt{\GG}_1', \mu', \delta = 1}$, by \cite[Thm. 4.3.1, 4.5.2]{PR24} (with modifications in the quasi-parahoric case using \cite[Thm. 4.1.8]{daniels2024conjecture}). Since $(\wdt{\GG}_1')^{\circ} = (\GG_1')^{\circ}$, we pass from $\wdt{\GG}_1'$ to $\GG_1'$ using \cite[Thm. 4.1.15]{daniels2024conjecture}; $\shu{K_{G_1}'}(G_1', X_1')^{\Dia} \to \Sht_{\GG_1', \mu', \delta = 1, \Spd E}$ extends to $\Shum{K_{G_1}'}(G'_1, X_1')^{\Dia/} \to \Sht_{\GG_1', \mu', \delta = 1}$.
     
     Let $P_3 = U_2 \rtimes G_1$, where $G_1$ acts on $U_2$ via the embedding $G_1 \to G_2$, and let $\PP_3 = \UU_2 \rtimes \GG_1$. Moreover, let $P_3' = U_2 \rtimes G_1'$ and $\PP_3' = \UU_2 \rtimes \GG_1'$. We can always do this since $\PP_1 \to \PP_2$ (resp. $\PP_1' \to \PP_2'$) is compatible with $\GG_1 \to \GG_2$ (resp. $\GG_1' \to \GG_2'$). Consider $\PP_3' \to \PP_2'$ that is compatible with $\GG_3' = \GG_1' \to \GG_2'$, by Corollary \ref{cor: quasi-parahoric, fiber product}, we have a dashed morphism that makes the diagrams commute:
\[
\begin{tikzcd}
	{\Shum{K_1'}(P_1', \mathcal{X}_1')^{\Dia/}} && {\Shum{K_2'}(P_2', \mathcal{X}_2')^{\Dia/}} \\
	& {\Sht_{\PP_3', \mu', \delta = 1}} & {\Sht_{\PP_2', \mu', \delta = 1}} \\
	{\Shum{K_{G_1}'}(G_1', X_1')^{\Dia/}} & {\Sht_{\GG_1', \mu', \delta = 1}} & {\Sht_{\GG_2', \mu', \delta = 1}.}
	\arrow[from=1-1, to=1-3]
	\arrow[dashed, from=1-1, to=2-2]
	\arrow[from=1-1, to=3-1]
	\arrow[from=1-3, to=2-3]
	\arrow[from=2-2, to=2-3]
	\arrow[from=2-2, to=3-2]
	\arrow["\lrcorner"{anchor=center, pos=0.125}, draw=none, from=2-2, to=3-3]
	\arrow[from=2-3, to=3-3]
	\arrow[from=3-1, to=3-2]
	\arrow[from=3-2, to=3-3]
\end{tikzcd}
\]
    Apply Proposition \ref{prop: generic fiber reduce, then integral reduce} to the exact sequence $1 \to \PP_1' \to \PP_3' \to \UU_2/\UU_1 \to 1$ (the Part $(2)$ in the proof of Proposition \ref{prop: generic fiber reduce, then integral reduce} does not need $U_1$ being normal in $U_3$). Then the dashed morphism $\Shum{K_1'}(P_1', \mathcal{X}_1')^{\Dia/} \to \Sht_{\PP_3', \mu', \delta = 1}$ factors through $\Sht_{\PP_1', \mu'}$, hence through $\Sht_{\PP_1', \mu', \delta = 1}$; see Corollary \ref{cor: quasi-parahoric}. Compose with
    $$\Shum{K_1}(P_1, \mathcal{X}_1) \to \Shum{K_1'}(P_1', \mathcal{X}_1')\quad (\textit{resp.}\ \Shum{K_{G_1}}(G_1, X_1) \to \Shum{K_{G_1}'}(G_1', X_1')),$$ 
    we have
    $$\Shum{K_1}(P_1, \mathcal{X}_1)^{\Dia/} \to \Sht_{\PP_1', \mu', \delta = 1}\quad (\textit{resp.}\ \Shum{K_{G_1}}(G_1, X_1)^{\Dia/} \to \Sht_{\GG_1', \mu', \delta = 1}) $$
    extending the one on the generic fiber. The commutativity of the diagram follows from \cite[Cor. 2.7.10]{PR24}.

 \item First of all, since $(\GG_1')^{\circ} = (\GG_2')^{\circ}$ by construction, $\shu{K_{G_1}'}(G', X')^{\Dia} \to \Sht_{\GG_1', \mu', \delta = 1, \Spd E}$ extends to $\Shum{K_{G_1}'}(G', X')^{\Dia/} \to \Sht_{\GG_1', \mu', \delta = 1}$, by \cite[Thm. 4.1.15]{daniels2024conjecture}. Composing with $\Shum{K_{G_1}}(G, X) \to \Shum{K_{G_1}'}(G', X')$, we obtain $\Shum{K_{G_1}}(G, X)^{\Dia/} \to \Sht_{\GG_1', \mu', \delta = 1}$.

   By Lemma \ref{lem-levi-c}, $P' = U \rtimes G'$. By Corollary \ref{cor: quasi-parahoric, fiber product}, we obtain the desired dashed morphism that extends the one on the generic fiber and makes the diagrams commute:
\[
\begin{tikzcd}
	{\Shum{K_1'}(P',\mathcal{X}')^{\Dia/}} && {\Shum{K'_2}(P', \mathcal{X}')^{\Dia/}} \\
	& {\Sht_{\PP_1', \mu', \delta = 1}} & {\Sht_{\PP_2', \mu', \delta = 1}} \\
	{\Shum{K'_{G_1}}(G', X')^{\Dia/}} & {\Sht_{\GG_1', \mu', \delta = 1}} & {\Sht_{\GG_2', \mu', \delta = 1}.}
	\arrow[from=1-1, to=1-3]
	\arrow[dashed, from=1-1, to=2-2]
	\arrow[from=1-1, to=3-1]
	\arrow[from=1-3, to=2-3]
	\arrow[from=2-2, to=2-3]
	\arrow[from=2-2, to=3-2]
	\arrow["\lrcorner"{anchor=center, pos=0.125}, draw=none, from=2-2, to=3-3]
	\arrow[from=2-3, to=3-3]
	\arrow[from=3-1, to=3-2]
	\arrow[from=3-2, to=3-3]
\end{tikzcd}
\]
    Similarly for $\Shum{K_1}(P, \mathcal{X})^{\Dia/} \to \Sht_{\PP_1', \mu', \delta = 1}$. The commutativity of the diagram follows from \cite[Cor. 2.7.10]{PR24}.
    \end{enumerate}
\end{proof}

\begin{prop}[{\cite[Thm. 4.3.1, 4.5.2]{PR24}, \cite[Thm. 4.1.15]{daniels2024conjecture}}]\label{prop: functoriality of canonical integral models}
    Consider one of the following two cases:
    \begin{enumerate}
        \item Let $\iota: (G_1, X_1) \to (G_2, X_2)$ be an embedding of pure Shimura data, and let $E_1$, $E_2$ be completions of reflex fields. Let $\GG_1$, $\GG_2$ be quasi-parahoric group schemes of $G_1$ and $G_2$ respectively, such that $\GG_2(\bZ_p) \cap G_1(\bQ) = \GG_1(\bZ_p)$. Let $K_2^p \subset G_2(\A^p)$ be any sufficiently small subgroup; one can always find $K_1^p \subset G_1(\A^p)$ such that the induced morphism $\shu{K_1}(G_1, X_1) \to \shu{K_2}(G_2, X_2)$ is a closed embedding, where $K_i = \GG_i(\Z_p)K_i^p$.
        Let $[\Phi_1] = [(Q_{\Phi_1}, X_{\Phi_1}^+, g_{\Phi_1})] \in \Cusp_{K_1}(G_1, X_1)$, and $[\Phi_2] = [\iota_*\Phi_1] = [(Q_{\Phi_2}, X_{\Phi_2}^+, g_{\Phi_2})] \in \Cusp_{K_2}(G_2, X_2)$ (with $g_{\Phi_1} = g_{\Phi_2}$). Here we can vary $K_1^p$ (resp. $K_2^p$), and do not distinguish $\Phi_1$ (resp. $\Phi_2$) for different levels $K_1^p$ (resp. $K_2^p$) once we have prescribed $g = g_{\Phi_1} = g_{\Phi_2} \in G(\A)$. Assume
        \begin{enumerate}
            \item $G_1^c \to G_2^c$ is an embedding, and the smoothing $\wdt{\GG}_1^c$ of the closure of $G_1^c$ in $\GG_2^c$ is a quasi-parahoric group scheme with $(\wdt{\GG}_1^c)^{\circ} = (\GG_1^c)^{\circ}$.
            \item $\PP_{\Phi_2}(\bZ_p) \cap P_{\Phi_1}(\bQ) = \PP_{\Phi_1}(\bZ_p)$ induces $\GG_{\Phi_2, h}(\bZ_p) \cap G_{\Phi_1, h}(\bQ) = \GG_{\Phi_1, h}(\bZ_p)$.
        \end{enumerate}

        \item Let $(G, X) := (G_1, X_1) = (G_2, X_2)$ be a pure Shimura datum, and let $E_1 = E_2$ be the completion of the reflex field. Let $\GG_1 \to \GG_2$ be quasi-parahoric group schemes of $G$ such that $\GGc_1 = \GGc_2$. Fix $(Q_{\Phi}, X_{\Phi}^+, g_{\Phi}) := (Q_{\Phi_1}, X_{\Phi_1}^+, g_{\Phi_1}) = (Q_{\Phi_2}, X_{\Phi_2}^+, g_{\Phi_2})$ as above. Let $K_i = \GG_i(\Z_p)K^p$ for some neat $K^p \subset G(\A^p)$.
    \end{enumerate}
    
    Assume $\lrbracket{\shu{K_2}}_{K^p_2}$ and $\lrbracket{\shu{K_2^c}}_{K^{c, p}_2}$ have canonical integral models $\lrbracket{\Shum{K_2}}_{K^p_2}$ and $\lrbracket{\Shum{K_2^c}}_{K^{c, p}_2}$ respectively that are adapted with $G_2 \to G_2^c$. Let $\lrbracket{\Shum{K_1}}_{K^p_1}$ and $\lrbracket{\Shum{K_1^c}}_{K^{c, p}_1}$ be the relative normalizations of $\lrbracket{\Shum{K_2}}_{K^p_2}$ and $\lrbracket{\Shum{K_2^c}}_{K^{c, p}_2}$ in $\lrbracket{\shu{K_1}}_{K^p_1}$ and $\lrbracket{\shu{K_1^c}}_{K^{c, p}_1}$ respectively. In case $(2)$, we \textbf{further assume $\lrbracket{\Shum{K_1} \to \Shum{K_1^c}}_{K^p_1}$ are \'etale}. Then $\lrbracket{\Shum{K_1}}_{K^p_1}$ and $\lrbracket{\Shum{K_1^c}}_{K^{c, p}_1}$ are canonical integral models of $\lrbracket{\shu{K_1}}_{K^p_1}$ and $\lrbracket{\shu{K_1^c}}_{K^{c, p}_1}$ respectively, adapted with $G_1 \to G_1^c$.

    Moreover, assume 
    \begin{equation}\label{eq: ca-1}
        \lrbracket{\shu{K_{\Phi_2}}}_{K_{\Phi_2}^p},\quad \lrbracket{\shu{K_{\Phi_2}^*}}_{K_{\Phi_2}^{*, p}},\quad \lrbracket{\shu{K_{\Phi_2, h}}}_{K_{\Phi_2, h}^p},\quad \lrbracket{\shu{K_{\Phi_2, h}^*}}_{K_{\Phi_2, h}^{*, p}}
    \end{equation}
    have canonical integral models
    \begin{equation}\label{eq: ca-2}
        \lrbracket{\Shum{K_{\Phi_2}}}_{K_{\Phi_2}^p},\quad \lrbracket{\Shum{K_{\Phi_2}^*}}_{K_{\Phi_2}^{*, p}},\quad \lrbracket{\Shum{K_{\Phi_2, h}}}_{K_{\Phi_2, h}^p},\quad \lrbracket{\Shum{K_{\Phi_2, h}^*}}_{K_{\Phi_2, h}^{*, p}}
    \end{equation}
    respectively, adapted with $P_{\Phi_2} \to P_{\Phi_2}^*$ and $P_{\Phi_2} \to G_{\Phi_2, h}$. Let
    \begin{equation}\label{eq: ca-4}
        \lrbracket{\Shum{K_{\Phi_1}}}_{K_{\Phi_1}^p},\quad \lrbracket{\Shum{K_{\Phi_1}^*}}_{K_{\Phi_1}^{*, p}},\quad \lrbracket{\Shum{K_{\Phi_1, h}}}_{K_{\Phi_1, h}^p},\quad \lrbracket{\Shum{K_{\Phi_1, h}^*}}_{K_{\Phi_1, h}^{*, p}} 
    \end{equation}
    be the relative normalizations of (\ref{eq: ca-2}) in
    \begin{equation}\label{eq: ca-3}
        \lrbracket{\shu{K_{\Phi_1}}}_{K_{\Phi_1}^p},\quad \lrbracket{\shu{K_{\Phi_1}^*}}_{K_{\Phi_1}^{*, p}},\quad \lrbracket{\shu{K_{\Phi_1, h}}}_{K_{\Phi_1, h}^p},\quad \lrbracket{\shu{K_{\Phi_1, h}^*}}_{K_{\Phi_1, h}^{*, p}}
    \end{equation}
    respectively. In case $(2)$, we \textbf{further assume $\lrbracket{\Shum{K_{\Phi_1}} \to \Shum{K_{\Phi_1}^*}}_{K_{\Phi_1}^p}$ and $\lrbracket{\Shum{K_{\Phi_1, h}} \to \Shum{K_{\Phi_1, h}^*}}_{K_{\Phi_1, h}^p}$ are \'etale}. Then (\ref{eq: ca-4}) are canonical integral models of (\ref{eq: ca-3}) respectively, adapted with $P_{\Phi_1} \to P_{\Phi_1}^*$ and $P_{\Phi_1} \to G_{\Phi_1, h}$.
\end{prop}
\begin{proof}
    Axiom $(1)$ is automatic; it comes from the construction of the relative normalization, see the first paragraph in the proof of \cite[Thm. 4.1.15]{daniels2024conjecture}. Axiom $(2)$ can be verified using axiom $(4)$, since $\Dia$ is full-faithful on the category of flat and normal formal schemes locally formally of finite type over $\OO_E$; see \cite[Prop. 18.4.1]{SW20}. We only need to consider axioms $(3)$ and $(4)$.
    
    For $(G_1, X_1)$:
    \begin{enumerate}
        \item We first verify that $\lrbracket{\Shum{K_1^c}(G_1^c, X_1^c)}_{K^{c, p}_1}$ is canonical in the sense of the Pappas--Rapoport conjectural framework; this is \cite[Thm. 4.3.1, 4.5.2]{PR24} (with quasi-parahoric modifications using \cite[Thm. 4.1.8]{daniels2024conjecture}, and with modifications in the $G \neq G^c$ case using the second paragraph of the proof of Lemma \ref{lem: two cases}, case $(1)$). Next, we show $\lrbracket{\Shum{K_1}(G_1, X_1)}_{K^{ p}_1}$ is canonical. The extension of shtukas is given by composing $\Shum{K_1} \to \Shum{K_1^c}$ with $\Shum{K_1^c} \to \Sht_{\GG_1^c, \mu_1^c, \delta = 1}$. Let $x \in \Shum{K_1}(k)$ and $\Bar{x} \in \Shum{K_1^c}(k)$ be its image. The existence of an isomorphism $\Theta_x: ((\Shum{K_1})^{\wedge}_{/x})^{\Dia} \cong (\mathcal{M}^{\intg}_{\GG_1^c, b_x, \mu_1^c})^{\wedge}_{/x_0}$ follows from the same arguments as in \cite[\S 4.7, 4.8]{PR24}, and $\Theta_x$ is compatible with $\Theta_{\Bar{x}}$ by construction. This in turn shows that the finite morphism $\Shum{K_1} \to \Shum{K_1^c}$ is \'etale, and $\lrbracket{\Shum{K_1}}_{K^{ p}_1}$ is adapted with $G_1 \to G_1^c$.
        \item The canonicity of $\lrbracket{\Shum{K_1^c}(G_1^c, X_1^c)}_{K^{c, p}_1}$ comes from \cite[Thm. 4.1.15]{daniels2024conjecture}. Let us show the canonicity of $\lrbracket{\Shum{K_1}(G_1, X_1)}_{K^{p}_1}$. Note that $\Shum{K_1} \to \Shum{K_2}$ is finite \'etale since $\Shum{K_1} \to \Shum{K_1^c}$ and $\Shum{K_2} \to \Shum{K_2^c}$ are finite \'etale by assumption, and $\Shum{K_1^c}\to \Shum{K_2^c}$ is finite {\'e}tale by axiom $(4)$; this implies axiom $(4)$. The extension of shtukas comes from Lemma \ref{lem: two cases}, case $(2)$.
    \end{enumerate}

    For $(P_{\Phi_1}, D_{\Phi_1})$ and $(G_{\Phi_1, h}, D_{\Phi_1, h})$, apply Lemma \ref{lem: two cases}, and let $P_i' = P_{\Phi_i}^*$, $\PP_i' = \PP_{\Phi_i}^*$. Verify the conditions in Lemma \ref{lem: two cases}:
    \begin{enumerate}
        \item We need to show that $P_1' \to P_2'$ is a closed embedding and that the smoothing $\wdt{\PP}_1'$ of the closure of $P_1'$ in $\PP_2'$ is a quasi-parahoric group scheme with $(\wdt{\PP}_1')^{\circ} = (\PP_1')^{\circ}$. Since $G_1^c \hookrightarrow G_2^c$ and $P_i' \hookrightarrow G_i^c$, we have $P_1' \hookrightarrow P_2'$. We have
        $$ \wdt{\PP}_1'(\bZ_p) = P_1'(\bQ) \cap \PP_2'(\bZ_p) = P_1'(\bQ) \cap \GG_2^c(\bZ_p) = P_1'(\bQ) \cap \GG_1^c(\bZ_p) = \PP_1'(\bZ_p). $$
        Since a smooth affine group scheme of a given group is uniquely determined by its set of $\bZ_p$-points, we have $\wdt{\PP}_1' = \PP_1'$.
        \item We need to show that $\GG_1 \to \GG_2$ induces $\PP_{\Phi_1} \to \PP_{\Phi_2}$ such that $\PP_{\Phi_1}^{\circ} = \PP_{\Phi_2}^{\circ}$. This is standard, since $(\PP_{\Phi_i}, \mu_{\Phi_i})$ comes from the boundary in the sense of Definition \ref{def: PP, mu comes from boundary}; see Lemma \ref{lem: quasi-parahoric}. On the other hand, $\GG_1^c \to \GG_2^c$ induces $\PP_1' \to \PP_2'$ such that $(\PP_1')^{\circ} = (\PP_2')^{\circ}$.
    \end{enumerate}

    Now let us show canonicity of integral models. For $(G_{\Phi_1, h}, D_{\Phi_1, h})$:
    \begin{enumerate}
        \item This is similar to $(G_1, X_1)$, but we need to handle the case where $G_{\Phi_1, h} \to G_{\Phi_1, h}^*$ has non-connected kernel. By assumption, the smoothing $\wdt{\GG}_{\Phi_1, h}^*$ of the closure of $G_{\Phi_1, h}^*$ in $\GG_{\Phi_2, h}^*$ is a quasi-parahoric group scheme and $(\wdt{\GG}_{\Phi_1, h}^*)^{\circ} = (\GG_{\Phi_1, h}^*)^{\circ}$; see the first paragraph of the proof of Lemma \ref{lem: two cases}. For the extension of shtukas, we apply Lemma \ref{lem: two cases}, case $(1)$. For axiom $(4)$, we apply results $(1)$ and $(2)$ successively from the first two paragraphs of this proof.
        \item This is similar to $(G_1, X_1)$. 
    \end{enumerate}
    For $(P_{\Phi_1}, D_{\Phi_1})$:
    \begin{enumerate}
        \item For extension of shtukas, we apply Lemma \ref{lem: two cases} Case $(1)$. Note that we can adjust level away from $p$ such that $\shu{K_1} \to \shu{K_2}$, $\shu{K_{\Phi_1}} \to \shu{K_{\Phi_2}}$ and $\shu{K_{\Phi_1, h}} \to \shu{K_{\Phi_2, h}}$ are all closed embeddings by \cite[Lem. 1.22]{Wu25}.
        Axiom $(4)$ follows from the same arguments as in \cite[\S 4.7, 4.8]{PR24}; we only explain where the arguments need to be adjusted or verified in our case.
        We use notation from \cite[\S 4.7, 4.8]{PR24}. The base point $x_0$ in $\MM^{\intg}_{\GG^{c}_{\Phi_2}, b_{\Phi_2, x}, \mu_{\Phi_2}^c}$ is the image of the base point in $\MM^{\intg}_{\GG^{c}_{\Phi_1}, b_{\Phi_1, x}, \mu_{\Phi_1}^c}$ (with $b_{\Phi_2, x} = b_{\Phi_1, x}$). The arguments in \cite[\S 4.7.2]{PR24} (with modifications in \cite[\S 4.8]{PR24}) imply that there is a $\Spd \hat{R}_x$-point (where $\Spf \hat{R}_x = \Shum{K_{\Phi_1}}(P_{\Phi_1}, D_{\Phi_1})^{\wedge}_{/x}$) of $\MM^{\intg}_{\GG^{c}_{\Phi_1}, b_{\Phi_1, x}, \mu_{\Phi_1}^c}$ lifting the base point $x_0$, such that the corresponding $\GG^{c}_{\Phi_1}$-shtuka is equal to the pushout $\GG^{c}_{\Phi_1} \times^{\PP^*_{\Phi_1}}$ of the $\PP^*_{\Phi_1}$-shtuka coming from $\PPs_{K_{\Phi_1}}$. In the arguments below \cite[Prop. 4.7.1]{PR24}, note that $\Spd (\hat{R}_x)_{\eta}$ and $(\MM^{\intg}_{\GG^{c}_{\Phi_1}, b_{\Phi_1, x}, \mu^c_{\Phi_1}})^{\wedge}_{/x_0}$ have the same dimension (see Remark \ref{rk: same dimension for mixed Shimura varieties}). Similarly for $\Shum{K_{\Phi_1}^*}(P_{\Phi_1}^*, D_{\Phi_1}^*)$, and $\Theta_{\Bar{x}}$ is compatible with $\Theta_{x}$, where $\Bar{x} \in \Shum{K_{\Phi_1}^*}(k)$ is the projection of $x$. This in particular shows that the finite morphism $\Shum{K_{\Phi_1}} \to \Shum{K_{\Phi_1}^*}$ is \'etale, and $\lrbracket{\Shum{K_{\Phi_1}}}_{K_{\Phi_1}^p}$ is adapted with $P_{\Phi_1} \to P_{\Phi_1}^*$.
        \item For the extension of shtukas, we apply Lemma \ref{lem: two cases}, case $(2)$. Axiom $(4)$ for $\Shum{K_{\Phi_1}^*}$: we need to show $\Shum{K_{\Phi_1}^*} \to \Shum{K_{\Phi_2}^*}$ is \'etale, and this follows from the same process as the proof of \cite[Thm. 4.1.15]{daniels2024conjecture}: Note that $\shu{K_{\Phi_1}^*} \to \shu{K_{\Phi_2}^*}$ is a finite \'etale torsor under the finite abelian group $K_{\Phi_2, p}^*/K_{\Phi_1, p}^* = \pi_0(\PP_{\Phi_2}^*)^{\phi}/\pi_0(\PP_{\Phi_1}^*)^{\phi}$ since $(P_i^*)^c = P_i^*$, and by \cite[Prop. 2.3.1]{daniels2024conjecture} and Corollary \ref{cor: quasi-parahoric, fiber product}, we have a cartesian diagram:
        \[
\begin{tikzcd}
	{\Shum{K_{\Phi_1}^*}(P_{\Phi_1}^*, D_{\Phi_1}^*)^{\Dia/}} & {\Sht_{\PP_{\Phi_1}^*, \mu_1^*, \delta=1}} \\
	{\Shum{K_{\Phi_2}^*}(P_{\Phi_2}^*, D_{\Phi_2}^*)^{\Dia/}} & {\Sht_{\PP_{\Phi_2}^*, \mu_2^*, \delta=1},}
	\arrow[from=1-1, to=1-2]
	\arrow[from=1-1, to=2-1]
	\arrow[from=1-2, to=2-2]
	\arrow[from=2-1, to=2-2]
\end{tikzcd}
        \]
        where $\Sht_{\PP_{\Phi_1}^*, \mu_1^*, \delta=1, E} \to \Sht_{\PP_{\Phi_2}^*, \mu_2^*, \delta=1, E}$ is a finite \'etale torsor under $\pi_0(\PP_{\Phi_2}^*)^{\phi}/\pi_0(\PP_{\Phi_1}^*)^{\phi}$, and $\Shum{K_{\Phi_1}^*} \to \Shum{K_{\Phi_2}^*}$ defined via relative normalization is also a a finite \'etale torsor under $\pi_0(\PP_{\Phi_2}^*)^{\phi}/\pi_0(\PP_{\Phi_1}^*)^{\phi}$. Axiom $(4)$ for $\Shum{K_{\Phi_1}}$ then follows from the \'etaleness of $\Shum{K_{\Phi_1}} \to \Shum{K_{\Phi_1}^*}$ (which also implies $\Shum{K_{\Phi_1}} \to \Shum{K_{\Phi_2}}$ is \'etale since $\Shum{K_{\Phi_2}} \to \Shum{K_{\Phi_2}^*}$ and $\Shum{K_{\Phi_1}^*} \to \Shum{K_{\Phi_2}^*}$ are \'etale).
    \end{enumerate}
\end{proof}

\begin{prop}[{\cite[Thm. 4.2.4]{PR24} and \cite[Prop. 4.1.10]{daniels2024conjecture}}]\label{prop: morphisms extend to integral models}
    Let $\iota: (G_1, X_1) \to (G_2, X_2)$ be a morphism of Shimura data (not necessarily an embedding), and let $E_1$ and $E_2$ be completions of reflex fields. Let $\GG_1$, $\GG_2$ be quasi-parahoric group schemes of $G_1$ and $G_2$ respectively, such that $G_1 \to G_2$ induces $\GG_1 \to \GG_2$. Let $[\Phi_1] = [(Q_{\Phi_1}, X_{\Phi_1}^+, g_{\Phi_1})] \in \Cusp_{K_1}(G_1, X_1)$, and $[\Phi_2] = [\iota_*\Phi_1] = [(Q_{\Phi_2}, X_{\Phi_2}^+, g_{\Phi_2})] \in \Cusp_{K_2}(G_2, X_2)$. Assume $\lrbracket{\shu{K_{\Phi_i}}(P_{\Phi_i}, D_{\Phi_i})}_{K_i^p}$ have canonical models $\lrbracket{\Shum{K_{\Phi_i}}(P_{\Phi_i}, D_{\Phi_i})}_{K_{\Phi_i}^p}$ that are adapted with $P_{\Phi_i} \to P_{\Phi_i}^*$ and $P_{\Phi_i} \to G_{\Phi_i, h}$ for $i = 1, 2$. Then the induced morphism $\shu{K_{\Phi_1}}(P_{\Phi_1}, D_{\Phi_1}) \to \shu{K_{\Phi_2}}(P_{\Phi_2}, D_{\Phi_2}) \otimes_{E_1} E_2$ extends canonically to
    \[ \Shum{K_{\Phi_1}}(P_{\Phi_1}, D_{\Phi_1}) \to \Shum{K_{\Phi_2}}(P_{\Phi_2}, D_{\Phi_2}) \otimes_{\OO_{E_1}} \OO_{E_2}. \]
\end{prop}
\begin{proof}
    Consider the diagonal embedding of Shimura data $(G_1, X_1) \hookrightarrow (G_1 \times G_2, X_1 \times X_2)$. Then $[\Phi_1] = [(Q_{\Phi_1}, X_{\Phi_1}^+, g_{\Phi_1})] \in \Cusp_{K_1}(G_1, X_1)$ induces $[\Phi_1 \times \Phi_2] = [(Q_{\Phi_1} \times Q_{\Phi_2}, X_{\Phi_1}^+ \times X_{\Phi_2}^+, g_{\Phi_1} \times g_{\Phi_2})] \in \Cusp_{K_1 \times K_2}(G_1 \times G_2, X_1 \times X_2)$. Note that $\GG_1(\bZ_p) = G_1(\bQ) \cap (\GG_1(\bZ_p) \times \GG_2(\bZ_p))$, $(G_1 \times G_2)^c = G_1^c \times G_2^c$, the closure of $G_1^c$ in $\GG_1^c \times \GG_2^c$ is exactly $\GG_1^c$. Moreover, $\GG_{\Phi_1, h}(\bZ_p) = G_{\Phi_1, h}(\bQ) \cap (\GG_{\Phi_1, h}(\bZ_p) \times \GG_{\Phi_2, h}(\bZ_p))$. Then we apply Proposition \ref{prop: functoriality of canonical integral models} case $(1)$, and get a canonical model $\lrbracket{\Shum{K_{\Phi_1}}(P_{\Phi_1}, D_{\Phi_1})'}_{K_{\Phi_1}^p}$ of $\lrbracket{\shu{K_{\Phi_1}}(P_{\Phi_1}, D_{\Phi_1})}_{K_{\Phi_1}^p}$. We need to show that the projection to the first factor $\Shum{K_{\Phi_1}}(P_{\Phi_1}, D_{\Phi_1})' \to \Shum{K_{\Phi_1}}(P_{\Phi_1}, D_{\Phi_1})$ is an isomorphism. We apply arguments from the proof of \cite[Prop. 4.1.10]{daniels2024conjecture}, and use the pushout $\GG_{\Phi_i}^{c}$-shtukas instead of $\PP_{\Phi_i}^*$-shtukas.
\end{proof}

\begin{cor}[{\cite[Thm. 4.2.4]{PR24} and \cite[Cor. 4.1.13]{daniels2024conjecture}}]\label{cor: uniqueness of canonical mixed}
    A canonical integral model $\lrbracket{\Shum{K_{\Phi}}(P_{\Phi}, D_{\Phi})}_{K_{\Phi}^p}$ of $\lrbracket{\shu{K_{\Phi}}(P_{\Phi}, D_{\Phi})}_{K_{\Phi}^p}$ that is adapted with $P_{\Phi} \to P_{\Phi}^*$ and $P_{\Phi} \to G_{\Phi, h}$ is unique up to a unique isomorphism, if it exists.
\end{cor}
\begin{proof}
    This follows from Proposition \ref{prop: morphisms extend to integral models} (see the proof of \cite[Cor. 4.1.13]{daniels2024conjecture}). 
\end{proof}

\section{Canonical extensions on toroidal compactifications}\label{sec-can-ext}

\subsection{Axiomatic setup of good compactifications}\label{subsec-gluing-setup}
Let $(G,X)$ be a Shimura datum. 
Fix an open compact subgroup $K\sbst G(\A)$, and assume that $K=K_pK^p$, where $K^p\sbst G(\Ap)$ is neat open compact and $K_p$ is open compact. 
Choose an admissible (rational polyhedral) smooth projective cone decomposition $\Sigma$ (without self-intersections).
Denote the toroidal compactification by  $\sh_K^{\Sigma}:=\sh_K^{\Sigma}(G,X)$ and the minimal compactification by $\sh_K^\mmin:=\sh_K^\mmin(G,X)$; they are defined over the reflex field $\bb{E}:=\bb{E}(G,X)$. 
There is a proper morphism from $\sh_K^\Sigma$ to $\sh_K^\mmin$ that is compatible with the stratifications on the source and the target.\par 
For integral models, a similar story is expected. 
The properties of these integral models are summarized as \cite[Prop. 2.1.2]{lan2018compactifications} and \cite[Thm. 4.1.5]{Mad19}:
\begin{axiom}[{Qualitative descriptions of good compactifications; \cite[Prop. 2.2]{LS18i} and \cite[Thm. 4.1.5 and Thm. 5.2.11]{Mad19}}]\label{axiom-good-compactification}Fix a prime number $p$ and a place $v$ of $\bb{E}$ over $p$. Set $E:=\bb{E}_v$.
For the cone decomposition $\Sigma$ above\footnote{In practice, it is usually harmless to refine the cone decompositions, as long as the refinements are still smooth and projective or are refinements of some smooth projective cone decompositions.}, there is a normal, proper, flat model $\mathscr{S}_K^\Sigma$ for $\sh^\Sigma_{K,E}$ over $\ca{O}_{E}$, and also a normal, projective, flat model $\mathscr{S}_K^\mmin$ for $\sh_{K,E}^\mmin$ over $\ca{O}_{E}$. The following properties hold for $\mathscr{S}_{K}^{\Sigma}$ and $\mathscr{S}_K^\mmin$:
\begin{enumerate}[label=(\textrm{\ref{axiom-good-compactification}}.\arabic*)]
\item\label{axiom-1} There is a proper surjective morphism $\oint_K^\Sigma: \mathscr{S}_K^{\Sigma}\to \mathscr{S}_K^\mmin$ with geometrically connected fibers. The map $\oint_K^\Sigma$ extends the one constructed in the characteristic zero theory. There are open dense embeddings $J^\Sigma:\mathscr{S}_K\to \mathscr{S}_K^\Sigma$ and $J^\mmin:\mathscr{S}_K\to \mathscr{S}_K^\mmin$ such that $\oint^\Sigma_K\circ J^{\Sigma}=J^\mmin$.
\item\label{axiom-2} There is a stratification of locally closed subschemes for $\mathscr{S}_K^\Sigma$, 
$$\mathscr{S}_K^\Sigma:=\disju_{\Upsilon:=[(\Phi,\sigma)]\in\cusp_K(G,X,\Sigma)}\ca{Z}_{[(\Phi,\sigma)],K},$$
and a stratification of locally closed subschemes for $\mathscr{S}_K^\mmin$,
$$\mathscr{S}_K^\mmin:=\disju_{[\Phi]\in\cusp_K(G,X)}\ca{Z}_{[\Phi],K}$$
extending the stratifications for $\sh_K^\Sigma$ and $\sh_K^\mmin$, respectively. Each stratum is normal and is flat over $\ca{O}_{E}$. The same partial orders among strata and the same expressions of closures of strata as the characteristic zero theory hold for $\{\ca{Z}_{[(\Phi,\sigma)],K}\}$ and $\{\ca{Z}_{[\Phi],K}\}$. The preimage $\oint^\Sigma_K(\ca{Z}_{[\Phi],K})$ is the union of the strata labeled by preimage of $[\Phi]$ in $\cusp_K(G,X,\Sigma)$ under the projection $\cusp_K(G,X,\Sigma)\to \cusp_K(G,X)$. For the definitions of $\cusp_K(G,X)$ and $\cusp_K(G,X,\Sigma)$, see, e.g., \cite[\S 1.1 and \S 1.3]{Wu25}.\par 
For any cusp label representative $\Phi$, there is a cone decomposition $\Sigma^+(\Phi)$ associated, which decomposes an open self-adjoint nondegenerated cone $\mbf{P}^+_\Phi$; the $\sigma$'s in $[(\Phi,\sigma)]$ are taken from $\Sigma^+(\Phi)$.
\item\label{axiom-3} There is a partial-ordered set $\ca{CLR}(G,X)$ of cusp label representatives from characteristic zero theory; the quotient of $\ca{CLR}(G,X)$ by the equivalence relations induces $\cusp_K(G,X)$. For each $\Phi\in\ca{CLR}(G,X)$, one can associate a mixed Shimura datum $(P_\Phi, D_\Phi)$ in the sense of \cite[Def. 2.1]{Pin89}. The mixed Shimura variety $\sh_{K_\Phi}:=\sh_{K_\Phi}(P_\Phi,D_\Phi)$ admits a tower $$\sh_{K_\Phi}(P_\Phi,D_\Phi)\to {\sh}_{\overline{K}_\Phi}(\overline{P}_\Phi,\overline{D}_\Phi)\to \sh_{K_{\Phi,h}}(G_{\Phi,h},D_{\Phi,h}),$$ 
where the first map is a torsor under a split torus $\mbf{E}_{K_\Phi}$ and the second map is a torsor under an abelian scheme over $\sh_{K_{\Phi,h}}$. The variety $\sh_{K_{\Phi,h}}$ is a pure Shimura variety over $E$.\par 
There is a tower 
\begin{equation}\label{eq-tower-int-mixed-sh}
\mathscr{S}_{K_\Phi}\to \mathscr{S}_{\overline{K}_\Phi}\to \mathscr{S}_{K_{\Phi,h}}
\end{equation}
of normal flat schemes of finite type over $\ca{O}_{E}$, such that the tower extends the one displayed above. The first map is also a torsor under the same split torus. The second map is proper and surjective.\par
For any equivalence $\Phi\xrightarrow{\sim}\Phi'$ in $\ca{CLR}(G,X)$, the corresponding towers of integral models as in (\ref{eq-tower-int-mixed-sh}) are canonically isomorphic.
\item\label{axiom-4} There is a commutative diagram
\begin{equation}
    \begin{tikzcd}
    {\mathscr{S}^\Sigma_K\supset \ca{Z}_{[(\Phi,\sigma)],K}}\arrow[dd]&\mathscr{S}_{K_\Phi,\sigma}\arrow[l,"{/\Delta^\circ_{\Phi,K}}"]\arrow[r,hook]&\mathscr{S}_{K_\Phi}(\sigma)\arrow[d]\\
    &&\overline{\mathscr{S}}_{K_\Phi}\arrow[d]\\
    {\mathscr{S}^\mmin_K\supset \ca{Z}_{[\Phi],K}}&&\mathscr{S}_{K_{\Phi,h}}.\arrow[ll,"{/\Delta_{\Phi,K}}"]
    \end{tikzcd}
\end{equation}
In the diagram above, there are groups $\Delta^\circ_{\Phi,K}\triangleleft \Delta_{\Phi,K}$ acting on $\mathscr{S}_{K_\Phi}(\sigma)$ and $\mathscr{S}_{K_{\Phi,h}}$, respectively. Their actions factor through finite quotients, and the quotient schemes have canonical isomorphisms $\ca{Z}_{[(\Phi,\sigma)],K}\iso \Delta^\circ_{\Phi,K}\bss \mathscr{S}_{K_\Phi,\sigma}$ and $\ca{Z}_{[\Phi],K}\iso \Delta_{\Phi,K}\bss \mathscr{S}_{K_{\Phi,h}}$. The scheme $\mathscr{S}_{K_\Phi,\sigma}$ is the $\sigma$-stratum of the (relative) affine toric embedding $\mathscr{S}_{K_\Phi}\hookrightarrow\mathscr{S}_{K_\Phi}(\sigma)$ with respect to the cone $\sigma$.\par 
For the definitions of $\Delta^\circ_{\Phi,K}$ and $\Delta_{\Phi,K}$, see \cite[\S 1.3]{Wu25}. In fact, there is an equivariant action of $\Delta_{\Phi,K}$ on (\ref{eq-tower-int-mixed-sh}), and $\Delta^\circ_{\Phi,K}$ is the normal subgroup fixing the induced action on the cone decomposition (which is independent of the choice of cones).\par
Moreover, the quotient $\mathscr{S}^*_{K_\Phi}:=\Delta_{\Phi,K}^\circ\bss \mathscr{S}_{K_\Phi}\to \overline{\mathscr{S}}^*_{K_\Phi}:=\Delta_{\Phi,K}^\circ\bss \overline{\mathscr{S}}_{K_\Phi}$ is also a torsor under a split torus ${\mbf{E}}_{\wdtd{K}_\Phi}$.
\item\label{axiom-5} There is a strata-preserving isomorphism
$$\mathfrak{X}_{\Upsilon,K}:=\mathfrak{X}_{[(\Phi,\sigma)],K}:=\cpl{\mathscr{S}_K^\Sigma}{\ca{Z}_{[(\Phi,\sigma)],K}}\iso \Delta_{\Phi,K}^\circ\bss \cpl{\mathscr{S}_{K_\Phi}(\sigma)}{\mathscr{S}_{K_\Phi,\sigma}}.$$
\item\label{axiom-6} There is a stronger strata-preserving isomorphism
$$\mathfrak{X}_{\Upsilon,K}^\circ\iso \Delta^\circ_{\Phi,K}\bss \cpl{\mathscr{S}_{K_\Phi}(\sigma)}{\mathscr{S}_{K_\Phi,\sigma}^+}.$$
The scheme $\mathscr{S}^+_{K_\Phi,\sigma}$ denotes $\cup_{\tau\subset \sigma,\tau\in \Sigma^+(\Phi)}\mathscr{S}_{K_\Phi,\tau}$ and $\mathfrak{X}^\circ_{\Upsilon,K}:=\cpl{\mathscr{S}_K^\Sigma}{\ca{Z}^+_{\Upsilon,K}}$, where $\ca{Z}^+_{\Upsilon,K}:=\cup_{\tau\subset \sigma,\tau\in\Sigma^+(\Phi)}\ca{Z}_{[(\Phi,\tau)],K}$.
\end{enumerate}
\end{axiom}
\begin{convention}
Denote $\mathscr{S}^{+,*}_{K_\Phi,\sigma}:=\Delta^\circ_{\Phi,K}\bss \mathscr{S}^{+}_{K_\Phi,\sigma}$. When $\Phi$ and $K$ are clear, we denote $\mathfrak{X}_{\Upsilon,K}$, $\mathfrak{X}^\circ_{\Upsilon,K}$, $\mathscr{S}^{+,(*)}_{K_\Phi,\sigma}$ and $\ca{Z}^+_{\Upsilon,K}$ by $\mathfrak{X}_{\sigma}$, $\mathfrak{X}_{\sigma}^\circ$, $\mathscr{S}^{+,(*)}_\sigma$ and $\ca{Z}^+_{\sigma}$, respectively.\par
Over the generic fiber, good compactifications satisfying the above axioms was proved by Pink \cite{Pin89} based on \cite{AMRT10}. The notation system for the generic fiber can be obtained by replacing $\mathscr{S}$ with $\sh$ and replacing $\ca{Z}$ with $\mrm{Z}$. 
\end{convention}
\begin{rk}
The main difference between Axiom \ref{axiom-good-compactification} and \cite[Prop. 2.2]{LS18i} is that we need to treat the case where $\Delta_{\Phi,K}^\circ$ is nontrivial. The item (9) in \emph{loc. cit.} is in fact a consequence of having a strata-preserving isomorphism between completions as in \ref{axiom-5} (see, e.g., \cite[Prop. 4.53]{Wu25}).
\end{rk}
Axiom \ref{axiom-6} was proved in the Hodge-type case in \cite[Prop. 2.1.3]{lan2018compactifications}.
In fact, Axiom \ref{axiom-6} is also true in the abelian-type case. The proof is the same as \cite[Prop. 4.32 and Lem. 4.49]{Wu25} with only some symbols changed; let us record it below.
\begin{prop}\label{prop-ab-axiom-6}
Axiom \ref{axiom-6} is true for integral models $\mathscr{S}^{\Sigma_2'}_{K_2}$ stated in \cite[Thm. 4.39]{Wu25}.
\end{prop}
\begin{proof}In this proof, we will use some conventions in \cite{Wu25} without explanation; the readers are recommended to consult the list of symbols there.\par
By replacing ``$\ca{S}_{\K_\Phi,\sigma}$'' with $\mathscr{S}_{\K_\Phi,\sigma}^+:=\cup_{\tau\subset \sigma,\tau\in \Sigma^+(\Phi)}\mathscr{S}_{\K_\Phi,\tau}$ and replacing ``$\ca{Z}_{[ZP^b(\Phi,\sigma)],K}$'' with $\ca{Z}^+_{[ZP^b(\Phi,\sigma)],K}:=\cup_{\tau\subset \sigma,\tau\in \Sigma^+(\Phi)}\ca{Z}_{[ZP^b(\Phi,\tau)],K}$, the argument in \cite[Prop. 4.32]{Wu25} goes through verbatim with the following exception: We need to see that $[\sigma]_{ZP}$ and $[\tau]_{ZP}$ are of the same cardinality. But $[\tau]_{ZP}$ is defined to be the $\Delta$-orbit in $\Delta_{ZP,K}$-orbit of $\tau$, where $\Delta\sbst \Delta^{ZP}_{\Phi,K}$ is some subgroup independent of the choice of $\tau$. Also, the stabilizer $\Delta^{ZP,\circ}_{\Phi,K}$ of the $\Delta^{ZP}_{\Phi,K}$-action is independent of the choice of $\tau$. We now have 
$$\cpl{\mathscr{S}_{K,\ca{O}_{K_Z}}^{\Sigma}}{\ca{Z}^+_{[ZP^b(\Phi,\sigma)],K,\ca{O}_{K_Z}}}\iso\cpl{\mathscr{S}_{\K_\Phi}(\sigma)_{\ca{O}_{K_Z}}}{\mathscr{S}^+_{\K_\Phi,\sigma,\ca{O}_{K_Z}}}.$$
Now, with the same replacements, the proof of \cite[Lem. 4.49]{Wu25} goes through verbatim. 
\end{proof}
In summary,
\begin{thm}[{\cite{FC90}, \cite{Lan16b}, \cite{Mad19}, and \cite{Wu25}}]\label{thm-abelian-type-axiom}
For abelian-type Shimura data, there are integral models satisfying Axiom \ref{axiom-good-compactification} for all $K$ at the beginning of \S\ref{subsec-gluing-setup}.
\end{thm}
\begin{rk}\label{rk-context}
In \S\ref{subsec-can-ext-generic-fiber}, we work over the generic fiber, so all arguments work without assumptions. In \S\ref{subsec-can-ext-int}, we assume Axiom \ref{axiom-good-compactification} holds for certain compactifications; however, the only condition we will need is a stratification with formal completions described by \textbf{toric schemes}. In \S\ref{subsec-can-mod-II}, we will further impose Axiom \ref{def: canonical model for mixed Shimura data} to define canonical integral models of compactifications; as we just mentioned, one does not have to require the stronger conditions in \S\ref{subsec-can-mod-II} to show the canonical extension of shtukas to integral models of toroidal compactifications.
\end{rk}
\subsection{Canonical extensions on generic fiber}\label{subsec-can-ext-generic-fiber}
Let $(G,X)$ be a Shimura datum and $\Sigma$ be smooth and projective with respect to $K$. On $\sh^{\Sigma}_K(G,X)$, we will explain that the pro-Kummer-\'etale $\ul{\ca{G}^c(\bb{Z}_p)}$-torsor induces a log shtuka with a good description at the boundary. 
\subsubsection{}
Let $\bb{P}_K:=\bb{P}_K(G,X)$ and $\bb{P}_{K_\Phi}:=\bb{P}_{K_\Phi}(P_\Phi,D_\Phi)$ be the pro-{\'e}tale torsor defined as in Definition \ref{def-proetale-torsor-aut}.
\begin{lem}\label{lem: delta-action}
    Let $\Phi_1 \stackrel{(\gamma, q_2)_K}{\longrightarrow} \Phi_2$ be an equivalence between cusp label representatives in $\ca{CLR}(G,X)$ (see \cite[\S 2.1.14]{Mad19}). Then we have a canonical isomorphism $\lp_{K_{\Phi_1}} \to \lp_{K_{\Phi_2}}$ that makes the diagram commute:
\[\begin{tikzcd}
	{\lp_{K_{\Phi_1}}} & {\lp_{K_{\Phi_2}}} \\
	{\shu{K_{\Phi_1}}(P_{\Phi_1}, D_{\Phi_1})} & {\shu{K_{\Phi_2}}(P_{\Phi_2}, D_{\Phi_2}).}
	\arrow["\cong", from=1-1, to=1-2]
	\arrow[from=1-1, to=2-1]
	\arrow[from=1-2, to=2-2]
	\arrow["\cong", from=2-1, to=2-2]
\end{tikzcd}\]
Similar statements hold for $\lp_{\ovl{K}_{\Phi_i}} \to \shu{\ovl{K}_{\Phi_i}}(\ovl{P}_{\Phi_i}, \ovl{D}_{\Phi_i})$ and $\lp_{K_{\Phi_i, h}} \to \shu{K_{\Phi_i, h}}(G_{\Phi_i, h}, D_{\Phi_i, h})$.
\end{lem}
\begin{proof}
    It suffices to show that, at each normal subgroup $K_p' \subset K_p$, we have a canonical dashed morphism fitting into the commutative diagram
\[\begin{tikzcd}
	{\shu{K'_{\Phi_1}}(P_{\Phi_1}, D_{\Phi_1})} & {\shu{K'_{\Phi_2}}(P_{\Phi_2}, D_{\Phi_2})} \\
	{\shu{K_{\Phi_1}}(P_{\Phi_1}, D_{\Phi_1})} & {\shu{K_{\Phi_2}}(P_{\Phi_2}, D_{\Phi_2}),}
	\arrow["\cong", dashed, from=1-1, to=1-2]
	\arrow[from=1-1, to=2-1]
	\arrow[from=1-2, to=2-2]
	\arrow["\cong", from=2-1, to=2-2]
\end{tikzcd}\]
where $K' = K_p'K^p$ and $K = K_pK^p$. Recall that $\gamma g_1 = q_2 g_2 k$ for some $k \in K$, we decompose $(\gamma, q_2)_K$ as
\begin{equation}\label{eq: decomposed of gamma, q_2}
    \Phi_1 \stackrel{(\gamma, 1)_K}{\longrightarrow} \Phi_1' \stackrel{(1, q_2)_K}{\longrightarrow} \Phi_2'  \stackrel{(1, k)_K}{\longrightarrow} \Phi_2,
\end{equation}
where $\Phi_1' = (Q_{\Phi_2}, X_{\Phi_2}^+, \gamma g_1)$, $\Phi_2' = (Q_{\Phi_2}, X_{\Phi_2}^+, g_2k)$.

$(\gamma, 1)_K$: Conjugation by $\gamma$ gives an isomorphism $[\Int(\gamma)]: \shu{K_{\Phi_1}}(P_{\Phi_1}, D_{\Phi_1}) \to \shu{K_{\Phi_1'}}(P_{\Phi_1'}, D_{\Phi_1'})$, similarly at $K'$-level.

$(1, q_2)_K$: Since $(\gamma g_1) = (q_2)(g_2k)$, $K_{\Phi_2'} = q_2^{-1}K_{\Phi_1'}q_2$, right multiplication by $q_2$ gives an isomorphism $[\cdot q_2]: \shu{K_{\Phi_1'}}(P_{\Phi_1'}, D_{\Phi_1'}) \to \shu{K_{\Phi_2'}}(P_{\Phi_2'}, D_{\Phi_2'})$, similarly at $K'$-level since $K_{\Phi_2'}' = q_2^{-1}K_{\Phi_1'}'q_2$.

$(1, k)_K$: Note that $K_{\Phi_2'} = K_{\Phi_2}$. The identity morphism induces $\shu{K_{\Phi_2'}}(P_{\Phi_2'}, D_{\Phi_2'}) \cong \shu{K_{\Phi_2}}(P_{\Phi_2}, D_{\Phi_2})$, similarly at $K'$-level since $K'$ is normal in $K$ and $K'_{\Phi_2'} = K'_{\Phi_2}$.

To summarize, we have a commutative diagram induced by (\ref{eq: decomposed of gamma, q_2}):
\begin{equation}\label{eq: decomposition of gamma, q_2 on Shimura}
\begin{tikzcd}
	{\shu{K'_{\Phi_1}}(P_{\Phi_1}, D_{\Phi_1})} & {\shu{K'_{\Phi_1'}}(P_{\Phi_1'}, D_{\Phi_1'})} & {\shu{K'_{\Phi_2'}}(P_{\Phi_2'}, D_{\Phi_2'})} & {\shu{K'_{\Phi_2}}(P_{\Phi_2}, D_{\Phi_2})} \\
	{\shu{K_{\Phi_1}}(P_{\Phi_1}, D_{\Phi_1})} & {\shu{K_{\Phi_1'}}(P_{\Phi_1'}, D_{\Phi_1'})} & {\shu{K_{\Phi_2'}}(P_{\Phi_2'}, D_{\Phi_2'})} & {\shu{K_{\Phi_2}}(P_{\Phi_2}, D_{\Phi_2}).}
	\arrow["{[\Int(\gamma)]}", from=1-1, to=1-2]
	\arrow[from=1-1, to=2-1]
	\arrow["{[\cdot q_2]}", from=1-2, to=1-3]
	\arrow[from=1-2, to=2-2]
	\arrow["\identity", from=1-3, to=1-4]
	\arrow[from=1-3, to=2-3]
	\arrow[from=1-4, to=2-4]
	\arrow["{[\Int(\gamma)]}", from=2-1, to=2-2]
	\arrow["{[\cdot q_2]}", from=2-2, to=2-3]
	\arrow["\identity", from=2-3, to=2-4]
\end{tikzcd}
\end{equation}
The statements for $\lp_{\ovl{K}_{\Phi_i}} \to \shu{\ovl{K}_{\Phi_i}}(\ovl{P}_{\Phi_i}, \ovl{D}_{\Phi_i})$ and $\lp_{K_{\Phi_i, h}} \to \shu{K_{\Phi_i, h}}(G_{\Phi_i, h}, D_{\Phi_i, h})$ are proved in the same way.
\end{proof}

Take $\Phi_1 = \Phi_2 = \Phi$. Then $\gamma$ belongs to
\[ \Delta_{\Phi,K} := (\stb_{Q(\rQ)}(D_\Phi) \cap P(\A)g_\Phi Kg_\Phi^{-1})/P(\rQ). \]
For each $\gamma \in \Delta_{\Phi,K}$, we find a $q_2 \in P(\A)$ such that $\Phi \stackrel{(\gamma, q_2)_K}{\longrightarrow} \Phi$. In particular, $\Delta_{\Phi,K}$ naturally acts on the tower
\begin{equation}\label{eq: Delta-action on tower}
    \shu{K_{\Phi}}(P_{\Phi}, D_{\Phi}) \to \shu{\ovl{K}_{\Phi}}(\ovl{P}_{\Phi}, \ovl{D}_{\Phi}) \to \shu{K_{\Phi, h}}(G_{\Phi, h}, D_{\Phi, h}),
\end{equation}
and by Lemma \ref{lem: delta-action}, $\Delta_{\Phi,K}$-action naturally lifts to the tower
\begin{equation}\label{eq: Delta-action on torsors}
    \lp_{K_{\Phi}} \to \lp_{\ovl{K}_{\Phi}} \to \lp_{K_{\Phi, h}}.
\end{equation}

\begin{prop}\label{prop: Delta acts trivially on generic fiber}
    The commutative diagrams (\ref{eq: reduction of torsors and HT maps}) and (\ref{eq: reduction, generic fiber, general}) in Lemma \ref{lem-red-torsors-HT-maps} and Corollary \ref{cor: reduction, generic fiber, general} are equivariant under the $\Delta_{\Phi,K}$-action.
\end{prop}
\begin{proof}
   We work with $\shu{K_{\Phi}}(P_{\Phi}, D_{\Phi})$; the other two follow in the same way. This essentially follows from Lemma~\ref{lem: delta-action}, but we spell out the details for completeness. We abbreviate $\shu{K_{\Phi}}(P_{\Phi}, D_{\Phi})$ to $\shu{K_{\Phi}}$. Let $x \in \shu{K_{\Phi}}$ be a closed point. The supported shtuka $\shu{K_{\Phi}}^{\Dia} \to \Sht_{\PP^c_{\Phi}, \mu_{\Phi}^c, \delta = 1, \Spd E}$ at $x$ is determined by the Hodge–Tate period map $\HT_{\Phi, x}: \lp_{K_{\Phi}, x} \to \Gra{P^c_{\Phi}, \mu_{\Phi}^c}$. This map is determined by the de Rham $\underline{\PP_{\Phi}^c(\Z_p)}$-torsor $\lp_{K_{\Phi}, x}$ itself, which in turn is determined by the Galois representation $\rho_x: \Gal(\ovl{k(x)}|k(x)) \to \PP_{\Phi}^c(\Z_p)$. We now recall the construction of $\rho_x$, following \cite[\S 3.1]{klevdal2023compatibility}.

    Fix a lifting $x \in \shu{\emptyset_{\Phi}}(\ovl{\rQ}_p):= \mathrm{Sh}(P_{\Phi}, D_{\Phi})(\ovl{\rQ}_p)$, with image $xK_{\Phi}' \in \shu{K_{\Phi}'}(\ovl{\rQ}_p)$ at each level $K_{\Phi}' \subset K_{\Phi}$. Let $S_{K_{\Phi}', x}$ be the geometrically connected component of $\shu{K_{\Phi'}}$ containing $xK_{\Phi}'$; it is defined over some abelian extension $E_{K_{\Phi}', x}$ of $E$. For an open normal subgroup $K_{\Phi}' \subset K_{\Phi}$, we have Cartesian squares of Galois coverings under the right action of $K_{\Phi}/Z(P_{\Phi})(\rQ)^-_{K_{\Phi}}K_{\Phi}'$:
    \begin{equation}\label{eq: Galois covering}
\begin{tikzcd}
	{\shu{K_{\Phi}'}} & {\shu{K_{\Phi}'} \otimes_E E_{K_{\Phi}, x}} & {\shu{K_{\Phi}'}\otimes_{\shu{K_{\Phi}}} S_{K_{\Phi}, x}} & {\wdt{S}_{K_{\Phi}', x}} \\
	{\shu{K_{\Phi}}} & {\shu{K_{\Phi}} \otimes_E E_{K_{\Phi}, x}} & {S_{K_{\Phi}, x},}
	\arrow[from=1-1, to=2-1]
	\arrow[from=1-2, to=1-1]
	\arrow[from=1-2, to=2-2]
	\arrow[hook', from=1-3, to=1-2]
	\arrow[from=1-3, to=2-3]
	\arrow[hook', from=1-4, to=1-3]
	\arrow[from=1-4, to=2-3]
	\arrow[from=2-2, to=2-1]
	\arrow[hook', from=2-3, to=2-2]
\end{tikzcd}
    \end{equation}
      where $\wdt{S}_{K_{\Phi}', x} \subset \shu{K_{\Phi}'} \otimes_E E_{K_{\Phi}, x}$ is the connected component that contains $xK_{\Phi}'$.  Let $\ovl{x}$ be a geometric point over $x$, then the above diagram induces a morphism
      \[ \pi_1(S_{K_{\Phi}, x}, \bar{x}) \to \Aut(\wdt{S}_{K_{\Phi}', x}/S_{K_{\Phi}, x})^{\mathrm{op}} \subset K_{\Phi}/Z(P_{\Phi})(\rQ)^-_{K_{\Phi}}K_{\Phi}'.  \]
      Taking the inverse limit, and specializing it to $(\Spec E_{K_{\Phi}, x}, \bar{s})$, we have
      \[ \wdt{\rho}_x: \Gal(\ovl{k(x)}|k(x)) \to \pi_1(S_{K_{\Phi}, x}, \bar{x}) \to K_{\Phi}/Z(P_{\Phi})(\rQ)^-_{K_{\Phi}}.  \]
      The projection of $\wdt{\rho}_x$ to $p$-factor $\PP^c_{\Phi}(\Z_p)$ is $\rho_x$.

      Let $\gamma \in \Delta_{\Phi,K}$. Let $y = \gamma x \in \shu{\emptyset_{\Phi}}(\ovl{\rQ}_p)$, as in Lemma \ref{lem: delta-action}, $[\gamma] $ acts compatibly on the diagram (\ref{eq: Galois covering}), thus $\wdt{\rho}_x$ and $\wdt{\rho}_y$ are conjugated by $[\gamma]$:
      \[
\begin{tikzcd}
	{\tilde{\rho}_x: \Gal(\ovl{k(x)}|k(x))} & {\pi_1(S_{K_{\Phi}, x}, \bar{x})} & {\Aut(\wdt{S}_{K_{\Phi}', x}/S_{K_{\Phi}, x})^{\mathrm{op}}} & {K_{\Phi}/Z(P_{\Phi})(\rQ)^-_{K_{\Phi}}K_{\Phi}'} \\
	{\tilde{\rho}_y: \Gal(\ovl{k(y)}|k(y))} & {\pi_1(S_{K_{\Phi}, y}, \bar{y})} & {\Aut(\wdt{S}_{K_{\Phi}', y}/S_{K_{\Phi}, y})^{\mathrm{op}}} & {K_{\Phi}/Z(P_{\Phi})(\rQ)^-_{K_{\Phi}}K_{\Phi}'.}
	\arrow[from=1-1, to=1-2]
	\arrow[from=1-2, to=1-3]
	\arrow["{[\gamma]}"', from=1-2, to=2-2]
	\arrow[from=1-3, to=1-4]
	\arrow["{[\gamma]}"', from=1-3, to=2-3]
	\arrow["{\Int(q_2^{-1}\gamma)}"', from=1-4, to=2-4]
	\arrow[from=2-1, to=2-2]
	\arrow[from=2-2, to=2-3]
	\arrow[from=2-3, to=2-4]
\end{tikzcd}
      \]
      We fix a representation $\rho_{\Phi}: P_{\Phi}(\rQ_p) \to P_{\Phi}^c(\rQ_p) \to \GL(W_{\rQ_p})$, and a lattice $W_{\Z_p} \subset W_{\rQ_p}$ such that $\rho_{\Phi}(\PP_{\Phi}^c(\Z_p)) \subset \GL(W_{\Z_p})$, this produces a de Rham local system $\ls_{\rho_{\Phi}, W_{\Z_p}}$ as in (\ref{eq: construction of local system from torsor, mixed}). For simplicity, we denote $\ls:=\ls_{\rho_{\Phi}, W_{\Z_p}}$. The canonical isomorphism $\gamma: \lp_{K_{\Phi}, x} \to \lp_{K_{\Phi}, y}$ induces an isomorphism $\gamma: \ls_x \to \ls_y$. We identify the underlying $\Z_p$-local systems $\ls_x \to \ls_y$, then the Galois action on $\ls_y$ is twisted by $[\gamma]$-conjugation from the one on $\ls_x$ since $\wdt{\rho}_x$ and $\wdt{\rho}_y$ are conjugated by $[\gamma]$. By the Cartan-Leray spectral sequence, the identification $\ls_x \cong \ls_y$ induces a canonical identification $(D_{\dR}(\ls_x), \nabla_{\ls_x}, \Fil_{\ls_x}) \cong (D_{\dR}(\ls_y), \nabla_{\ls_y}, \Fil_{\ls_y})$. In particular, in the construction \cite[Prop. 2.6.3]{PR24}, $[\gamma]$-action identifies $\mathrm{DRT}(\ls_x)$ and $\mathrm{DRT}(\ls_y)$. Therefore, $[\gamma]$ induces an equivariant action on $\shu{K_{\Phi}}(P_{\Phi}, D_{\Phi})^{\Dia} \to [\Gra{P^c_{\Phi}, \mu_{\Phi}^{c, -1}}/\underline{\PP^c_{\Phi}(\Z_p)}]$ where it acts trivially on the target. Finally, since $\Delta_{\Phi,K}$-action on $\shu{K_{\Phi}}$ lifts canonically to a $\Delta_{\Phi,K}$-action on $\lp_{K_{\Phi}}$ as in Lemma \ref{lem: delta-action}, we have a cocycle condition for the identifications among $\lrbracket{\mathrm{DRT}(\ls_{\gamma(x)})}_{\gamma \in \Delta_{\Phi,K}}$, $\Delta_{\Phi,K}$ acts on the diagram $(\ref{eq: reduction, generic fiber, general})$.
\end{proof}
\subsubsection{}\label{*-local-system}
Recall that we have
 \begin{equation*}
     \mathfrak{X}_{\sigma,\eta}^{\circ}:= \Delta_{\Phi,K}^\circ\backslash (\sh_{K_\Phi}(\sigma))_{\underset{\tau\in\Sigma^+(\Phi),\ \tau\subset\sigma}{\cup}\sh_{K_\Phi,\tau}}^{\wedge} \cong(\shuc{K}{\Sigma})_{\underset{\tau\in \Sigma^+(\Phi),\ \tau\subset\sigma}{\cup} \mrm{Z}_{[(\Phi, \tau)]}}^{\wedge}.
 \end{equation*}
 Here $\Delta_{\Phi,K}^\circ \subset \Delta_{\Phi,K}$ is the subgroup that stabilizes $\sigma$; when $K$ is neat, it is independent of the choice of $\sigma$ (see \cite[Thm. 6.19]{Pin89} and also\cite[\S 2.1.19]{Mad19}).\par
 We first consider the torsor defined by the tower $\varprojlim_{K_p'}\Delta^\circ_{\Phi,K'}\bss\sh_{K'_\Phi}$.
Denote by $ZP_\Phi$ the connected component of $Z_G\cdot P_\Phi$, and by $ZP_\Phi(\bb{Q})_1:=\stb_{ZP_\Phi(\bb{Q})}(D_\Phi)$. When $K$ is neat, by \cite[Lem. 1.11]{Wu25}, we have that $$\Delta^\circ_{\Phi,K}\bss \sh_{K_\Phi}(\bb{C})=ZP_\Phi(\bb{Q})_1\bss D_\Phi\times ZP_\Phi(\bb{Q})_1P_\Phi(\A)/ZP_\Phi(\bb{Q})_1P_\Phi(\A)\cap g_\Phi K g_\Phi^{-1}.$$
Then the Galois group is $\mathscr{G}_p^\circ(\Phi,K):=\Gal(\varprojlim_{K_p'\sbst K_p}\Delta^\circ_{\Phi,K'}\bss \sh_{K_\Phi'}/\Delta^\circ_{\Phi,K}\bss \sh_{K_\Phi})$
\begin{equation}\label{eq-galois-zp}
    \begin{split}
    &=\Gal(\frac{\varprojlim_{K_p'\sbst K_p}ZP_\Phi(\bb{Q})_1\bss D_\Phi\times ZP_\Phi(\bb{Q})_1P_\Phi(\A)/ZP_\Phi(\bb{Q})_1P_\Phi(\A)\cap g_\Phi K' g_\Phi^{-1}}{ZP_\Phi(\bb{Q})_1\bss D_\Phi\times ZP_\Phi(\bb{Q})_1P_\Phi(\A)/ZP_\Phi(\bb{Q})_1P_\Phi(\A)\cap g_\Phi K g_\Phi^{-1}})\\
    &=\frac{ZP_\Phi(\bb{Q})_1P_\Phi(\A)\cap g_\Phi K_pg_\Phi^{-1}}{(Z(\bb{Q})^{\overline{\ }}\cap g_\Phi Kg_\Phi^{-1})(ZP_\Phi(\bb{Q})_1P_\Phi(\A)\cap g_\Phi K^pg_\Phi^{-1})}\\
    &=\frac{ZP_\Phi(\bb{Q})_1P_\Phi(\A)g_\Phi K^p g_\Phi^{-1}\cap g_\Phi K_p g_\Phi^{-1}}{(Z(\bb{Q})^{\overline{ }}\cap g_\Phi Kg_\Phi^{-1})g_\Phi K^pg_\Phi^{-1}\cap g_\Phi K_pg_\Phi^{-1}}\\
    &=g_\Phi\frac{(g_\Phi^{-1}ZP_\Phi(\bb{Q})_1P_\Phi(\A)g_\Phi)K^p\cap K_p}{(Z(\bb{Q})^{\overline{ }}\cap K)K^p\cap K_p}g_\Phi^{-1}.
    \end{split}
\end{equation}
The image of the numerator of $\mathscr{G}_p^\circ(\Phi,K)$ in $ZP^c_\Phi(\bb{Q}_p)=(ZP_\Phi)^c(\bb{Q}_p)$ lies in $P^*_\Phi(\bb{Q}_p)$; this is also due to neatness and the fact that $ZP^c_\Phi/P_\Phi^*$ is a cuspidal torus.
\begin{definition}\label{def-*-local-system}
Denote $\bb{P}_{K_\Phi}^*:=\varprojlim_{K_{\Phi,p}'} \Delta^\circ_{\Phi,K'}\bss \sh_{K'_\Phi}\times^{\mathscr{G}_p^\circ(\Phi,K)}\ul{\ca{P}^*_\Phi(\bb{Z}_p)}$.
\end{definition}
Then, by Proposition \ref{prop: Delta acts trivially on generic fiber}, the quotient induced by the action of $\Delta_{\Phi,K}^\circ$ descends $\bb{P}_{K_\Phi}\times^{\ul{\ca{P}^c_\Phi(\bb{Z}_p)}} \ul{\ca{P}^*_\Phi(\bb{Z}_p)}$ to $\bb{P}^*_{K_\Phi}$; the Hodge-Tate map $\mrm{HT}_{K_\Phi}$ descends to 
$$\mrm{HT}_{K_\Phi}^*:\bb{P}^*_{K_\Phi}\to \Gra{P^*_{\Phi},\mu_\Phi^{*,-1}},$$
where $\mu_\Phi^*$ is the projection of $\mu_\Phi$ to $P^*_\Phi$. 
More precisely,
\begin{lem}\label{lem: reduction, generic fiber, general, star}
   The commutative diagram (\ref{eq: reduction, generic fiber, general}) descends to the commutative diagram:
    \begin{equation}\label{eq: reduction, generic fiber, general, star}
\begin{tikzcd}
	{\Delta_{\Phi, K}^{\circ}\backslash\shu{K_{\Phi}}(P_{\Phi}, D_{\Phi})^{\Dia}} & {\Delta_{\Phi, K}^{\circ}\backslash\shu{\ovl{K}_{\Phi}}(\ovl{P}_{\Phi}, \ovl{D}_{\Phi})^{\Dia}} & {\Delta_{\Phi, K}^{\circ}\backslash\shu{K_{\Phi, h}}(G_{\Phi, h}, D_{\Phi, h})^{\Dia}} \\
	{\Sht_{\PP_{\Phi}^*, \mu_{\Phi}^*, \delta = 1, \Spd E}} & {\Sht_{\ovl{\PP}_{\Phi}^*, \bar{\mu}_{\Phi}^*, \delta = 1, \Spd E}} & {\Sht_{\GG_{\Phi, h}^*, \mu_{\Phi, h}^*, \delta = 1, \Spd E}.}
	\arrow[from=1-1, to=1-2]
	\arrow[from=1-1, to=2-1]
	\arrow[from=1-2, to=1-3]
	\arrow[from=1-2, to=2-2]
	\arrow[from=1-3, to=2-3]
	\arrow[from=2-1, to=2-2]
	\arrow[from=2-2, to=2-3]
\end{tikzcd}
    \end{equation}
\end{lem}
\begin{proof}
This follows from Proposition \ref{prop: Delta acts trivially on generic fiber}, Lemma \ref{lem: fiber product of shtukas} and Proposition \ref{prop: shtuka and loc system, general}.
\end{proof}

 Let $\mathfrak{W} = \Spf( R,I) \subset \mathfrak{X}_{\sigma,\eta}^{\circ}$ be an affine open formal subscheme, we can consider the flat morphisms 
 $W = \Spec R \to \shuc{K}{\Sigma}$ and $W \to \Delta_{\Phi,K}^\circ\backslash\sh_{K_\Phi}(\sigma)$.
Let $W^0 \subset W$ be the common open subscheme associated with $\shu{K} \subset \shuc{K}{\Sigma}$ and $\Delta_{\Phi,K}^\circ\bss\sh_{K_\Phi} \subset \Delta_{\Phi,K}^\circ\bss\sh_{K_\Phi}(\sigma)$. We have flat morphisms $W^0 \to \shu{K}(G, X)$ and $W^0 \to \Delta_{\Phi,K}^\circ\bss\shu{K_{\Phi}}(P_\Phi, D_{\Phi})$.

 We pull $\PPp_K$ and $\PPp_{K_{\Phi}}^*$ back to $W^0$, and denote them by $\PPp_{K, W^0}$ and $\PPp_{K_{\Phi}, W^0}$, respectively. Now we compare these two pro-\'etale torsors.

 \begin{lem}\label{lem: torsors reduction over W^0, generic}
          $\PPp_{K, W^0} \cong \ul{\GG^c(\Z_p)} \times^{\ul{\PP^{*}_{\Phi}(\Z_p)}} \PPp^*_{K_{\Phi}, W^0}$, where $\PP^*_{\Phi}(\bb{Z}_p) \to \GG^c(\bb{Z}_p)$ is the composition of $\ca{P}_\Phi^*(\bb{Z}_p)\xrightarrow{\Int(g^{-1}_\Phi)} g_\Phi^{-1}\ca{P}^*_\Phi(\bb{Z}_p) g_\Phi\hookrightarrow \ca{G}^c(\bb{Z}_p)$.
 \end{lem}
\begin{proof} Note that $g_\Phi^{-1}\ca{P}^*_\Phi(\bb{Z}_p)g_\Phi\to \ca{G}^c(\bb{Z}_p)$ is injective by Definition \ref{def-quasi-parahoric-big}. 
Let $K' \subset K$. The cone decomposition of the toroidal compactification at level $K'$ is the cone decomposition induced by $\Sigma$, which is denoted by the same symbol.
There is a transition map $\pi_{K',K}^\Sigma:\sh_{K'}^\Sigma\to \sh_K^\Sigma$. Pick a stratum $\mrm{Z}_{[(\Phi,\sigma)],K}$. The preimage of $\mrm{Z}_{[(\Phi,\sigma)],K}$ under $\pi^\Sigma_{K',K}$ is the disjoint union
$$\disju_{[(\Phi',\sigma')]\to[(\Phi,\sigma)]}\mrm{Z}_{[(\Phi',\sigma')],K'},$$
where $[(\Phi',\sigma')]$ are cusp labels with cones in $\cusp_{K'}(G,X,\Sigma)$ that map to $[(\Phi,\sigma)]$ in $\cusp_K(G,X,\Sigma)$.\par
Let $\mathfrak{W} = \Spf R \subset \mathfrak{X}^\circ_{\sigma,\eta}$ be an open formal subscheme. Let $\mathfrak{W}' =\Spf R'$ be the pullback via $\pi^\Sigma_{K',K}$. Denote $W:=\spec R$ and $W':=\spec R'$. We then have a commutative diagram
\begin{equation*}
    \begin{tikzcd}
    \sh_{K'}^\Sigma\arrow[d]& W'\arrow[l]\arrow[r]\arrow[d]&{ \disju\limits_{[(\Phi',\sigma')]\to[(\Phi,\sigma)]}\Delta^\circ_{\Phi',K'}\bss\sh_{K_{\Phi'}}(\sigma')}\arrow[d]\\
    \sh_K^\Sigma&W\arrow[l]\arrow[r]&\Delta^\circ_{\Phi,K}\bss\sh_{K_\Phi}(\sigma),
    \end{tikzcd}
\end{equation*}
where both squares are Cartesian.\par
Recall that, fixing $Q_\Phi$, the set $I_K(Q_\Phi,\Sigma)$ (see \cite[Def. 1.32]{Wu25}) consists of the cusp labels with cones $[(\Phi_1,\sigma_1)]$ equivalent to $[(Q_\Phi,X_\Phi^{+,\prime},g_\Phi',\sigma')]$ for some $(X_\Phi^{+,\prime},g_\Phi',\sigma')$ in $\cusp_K(G,X,\Sigma)$; we add subscript $K$ here to emphasize the role of $K$ in the equivalence relation. By \cite[Prop. 1.4]{Wu25} and the paragraph above \cite[Def. 1.16]{Wu25}, 
$$I_K(Q_\Phi,\Sigma)\iso \stb_{Q_\Phi(\bb{Q})}(D_\Phi)\bss(\Sigma^+(\Phi)\times P_\Phi(\A)\bss G(\A)/K).$$

From this, we see that $K$ permutes all $[(\Phi',\sigma')]$ mapping to $[(\Phi,\sigma)]$.
Restricting to $W^0$, the diagram above gives
\begin{equation*}
    \begin{tikzcd}
    \sh_{K'}\arrow[d]& W^{\prime,0}\arrow[l]\arrow[r]\arrow[d]&{ \disju\limits_{[(\Phi',\sigma')]\to[(\Phi,\sigma)]}\Delta^\circ_{\Phi',K'}\bss\sh_{K_{\Phi'}}}\arrow[d]\\
    \sh_K&W^0\arrow[l]\arrow[r]&\Delta^\circ_{\Phi,K}\bss\sh_{K_\Phi},
    \end{tikzcd}
\end{equation*}
where both squares are Cartesian.
Now, we shrink $K_p$ and take the inverse limit. The pushout of the inverse limit of $[(\Phi,\sigma)]\to[(\Phi,\sigma)]$ on the right vertical arrow gives rise to $\bb{P}^*_{K_\Phi}$.\par 
We need to see that on $W^0$, $\bb{P}_{K,W^0}\iso \ul{\ca{G}^c(\bb{Z}_p)}\times^{\ul{\ca{P}_\Phi^*(\bb{Z}_p)}}\bb{P}^*_{K_\Phi,W^0}$. 
For this, we compute that $K_p$ acts transitively on $\varprojlim_{K_p'} I_{K'}(Q_\Phi,\Sigma)$ with stabilizer $(g_\Phi^{-1}\stb_{Q_\Phi(\bb{Q})}(D_\Phi,\sigma)P_\Phi(\A)g_\Phi K^p)\cap K_p$ at $[(\Phi,\sigma)]$. When $K$ is neat, this intersection is independent of the choice of $\sigma\in \Sigma^+(\Phi)$ and is equal to $(g_\Phi^{-1} ZP_\Phi(\bb{Q})_1P_\Phi(\A)g_\Phi K^p)\cap K_p$. \par
We finally compute that on $W^0$: $[\varprojlim\Delta^\circ_{\Phi,K}\bss \sh_{K_\Phi}]\times^{\mathscr{G}^\circ_{p}(\Phi,K)}\ul{\ca{G}^c(\bb{Z}_p)}$
\begin{equation*}
    \begin{split}
        &=[\varprojlim\Delta^\circ_{\Phi,K}\bss \sh_{K_\Phi}]\times^{\mathscr{G}^\circ_{p}(\Phi,K)}\ul{K_p/(K_p\cap K^pZ(\bb{Q})_K^{\overline{\ }})}\times^{ \ul{ K_p/(K_p\cap K^pZ(\bb{Q})_K^{\overline{\ }})}}\ul{\ca{G}^c(\bb{Z}_p)}\\
        &=[\varprojlim_{K_p'\sbst K_p}\disju_{[(\Phi',\sigma')]\to[(\Phi,\sigma)]}\Delta^\circ_{\Phi',K'}\bss \sh_{K_{\Phi'}}]\times^{ \ul{ K_p/(K_p\cap K^pZ(\bb{Q})_K^{\overline{\ }})}}\ul{\ca{G}^c(\bb{Z}_p)}\\
        &=\varprojlim_{K_p'\sbst K_p}\sh_{K'}\times^{ \ul{ K_p/(K_p\cap K^pZ(\bb{Q})_K^{\overline{\ }})}}\ul{\ca{G}^c(\bb{Z}_p)}.
    \end{split}
\end{equation*}
The only nontrivial part is the computation from the second line to the third line. This is done by comparing (\ref{eq-galois-zp}) with the stabilizer in the last paragraph; note that there is a $g_\Phi$-conjugation difference in the construction.
\end{proof}
 \begin{lem}\label{lem: HT match}
     We have the following commutative diagram
\[\begin{tikzcd}
	{\PPp_{K_{\Phi}, W^0}^*} & {\PPp_{K, W^0}} \\
	{\Gra{P_{\Phi}^*, \mu_{\Phi}^{*, -1}}} & {\Gra{G^c, \mu^{c, -1}}.}
	\arrow[from=1-1, to=1-2]
	\arrow["{\HT_{K_{\Phi}}^*}"', from=1-1, to=2-1]
	\arrow["{\HT_K}"', from=1-2, to=2-2]
	\arrow["{\mrm{Int}(g_{\Phi}^{-1})}", from=2-1, to=2-2]
\end{tikzcd}\]
 \end{lem}
 \begin{proof}       
Recall that, in (\ref{eq: construction of local system from torsor, mixed}), for any representation $\rho: G^c(\Qp) \to \GL(V_{\Qp})$ and a lattice $V_{\Z_p} \subset V_{\Qp}$ such that $\rho(\GG^c(\Z_p)) \subset \GL(V_{\Z_p})$, $\rho$ induces a morphism $\rho_{\Phi}: P_{\Phi}^*(\Qp) \to \GL(V_{\Qp})$ such that $\rho_\Phi(\PP_{\Phi}^{*}(\Z_p))=\rho(g_\Phi^{-1}(\ca{P}^*_\Phi(\bb{Z}_p))g_\Phi) \subset \GL(V_{\Z_p})$. Since $\ls_{{\rho_\Phi},V_{\bb{Z}_p},W^0}=\ls_{\rho(g_{\Phi}^{-1}(-)g_\Phi),V_{\bb{Z}_p},W^0}$ is de Rham by \cite{liu2017rigidity}, the association of Hodge-Tate period map is intrinsic (see \cite[Prop. 2.6.3]{PR24}). This gives the desired commutative diagram.  
 \end{proof}
 \begin{cor}\label{cor: shtuka comparison, over W, generic fiber}
     We have the following commutative diagram
     \begin{equation}\label{eq: shtuka comparison, over W}
\begin{tikzcd}
	{\shu{K}(G, X)^{\Dia}} & {W^{0, \Diamond}} & {\Delta_{\Phi,K}^\circ\bss\shu{K_{\Phi}}(P_{\Phi}, D_{\Phi})^{\Dia}} \\
	{\Sht_{\GG^{c}, \mu^c}} & {} & {\Sht_{\PP_{\Phi}^{*}, \mu^*_{\Phi}}.}
	\arrow[from=1-1, to=2-1]
	\arrow[from=1-2, to=1-1]
	\arrow[from=1-2, to=1-3]
	\arrow[from=1-3, to=2-3]
	\arrow["{\mrm{Int}(g_{\Phi}^{-1})}"', from=2-3, to=2-1]
\end{tikzcd}
     \end{equation}
 \end{cor} 
  \begin{cor}\label{cor-nilp-monodromy}
For any representation $(\rho,V_{\bb{Z}_p})$ of $\ca{G}^c(\bb{Z}_p)$ where $V_{\bb{Z}_p}$ is a finite free module over $\bb{Z}_p$, the pro-Kummer {\'e}tale local system $J^\Sigma_*\ls_{\rho,V_{\bb{Z}_p},W^0}$ has unipotent geometric monodromy in the sense of Remark \ref{rk-neat-unip} (cf. \cite[Def. 6.3.7]{DLLZ23}). Hence,  by \cite{liu2017rigidity} and Corollary \ref{cor-canonical-extensions} (2) and (3), there is a log shtuka $\mathscr{P}^{\can}_\eta$ on $(\sh^{\Sigma}_K)^{\log \Diamond}$ uniquely extending the shtuka $\mathscr{P}_\eta$ on $(\sh_K(G,X))^\Diamond$ associated with $\bb{P}_K$ under \cite[Prop. 2.5.3]{PR24}.
 \end{cor}
 \begin{proof}
Note that it suffices to consider the mixed Shimura varieties in the form of 
$$\varprojlim_{K_p'\sbst K_p}\sh_{\K_{\Phi}'}(ZP_\Phi,ZP_\Phi(\bb{Q})D_\Phi),$$
where $\K_\Phi':=ZP_\Phi(\A)\cap g_\Phi K' g_\Phi^{-1}$. Indeed, this follows from the description of the boundary using $ZP$-cusps (see \cite[(1.32)]{Wu25}), and from the open and closed immersion $\Delta^\circ_{\Phi,K'}\bss\sh_{K_\Phi'}\hookrightarrow\sh_{\K_\Phi'}(ZP_\Phi,ZP_\Phi(\bb{Q})D_\Phi)$ for all $K'_p$ by \cite[Lem. 1.42 (1)]{Wu25} and \cite[Lem. 1.11]{Wu25}.
For unipotence, combine Lemma \ref{lem: torsors reduction over W^0, generic} and Lemma \ref{lem-torus-action-tower-pro-p-gen} (the mixed Shimura variety in the lemmas is $\sh_{\K_\Phi}(ZP_\Phi,ZP_\Phi(\bb{Q})D_\Phi)$).
Then the action of the geometric Kummer {\'e}tale fundamental group factors through $\mbf{E}_\infty'$ on the pro-Kummer {\'e}tale torsor $J(\sigma)_*\bb{P}^*_{K_\Phi}$ on $\Delta_{\Phi,K}^\circ\bss\sh_{K_\Phi}(P_\Phi,D_\Phi)(\sigma)\sbst \sh_{\K_\Phi}(ZP_\Phi,ZP_\Phi(\bb{Q})D_\Phi)(\sigma)$. Here $J(\sigma): \Delta^\circ_{\Phi,K}\bss\sh_{K_\Phi}(P_\Phi,D_\Phi)\hookrightarrow \Delta^\circ_{\Phi,K}\bss \sh_{K_\Phi}(P_\Phi,D_\Phi)(\sigma)$ is the twisted toric embedding. 
By the proof of Lemma \ref{lem-torus-action-tower-pro-p-gen} and since we are considering $P^*_\Phi$, which is a quotient from $P^c_\Phi$, we can further reduce the problem to the case where $P=P^c$. The torus action is defined by (\ref{eq-torus-action-gen}) and (\ref{eq-action-4-step}). When $K_p$ varies, the action of the torus part is given by the action of a conjugate of a unipotent element $u_f\in U(\A)$ as explained at the beginning of \S\ref{subsubsec-tower-computation}.  
The other statement is self-explanatory.
 \end{proof}
\subsection{Canonical extensions on integral models}\label{subsec-can-ext-int}
We formulate and prove two general extension results for log shtukas. The techniques we use are consequences of toric charts but not the geometry of integral models of (open) Shimura varieties; this, in particular, enables us to show the main theorems in satisfactory generality.
\subsubsection{Gluing lemma}\label{subsubsec-glue-thm}
The rough idea of the following lemma is to glue the log shtukas using toric charts. In particular, we do not use the geometry of $\mathscr{S}_K$. Let $E$ be the completion of the reflex field of the Shimura datum $(G,X)$ that we consider.
\begin{lem}\label{lem-gluing-shtukas}We still assume that $\Sigma$ is smooth projective. 
Suppose that there is an integral model $\mathscr{S}^\Sigma_K$ satisfying Axiom \ref{axiom-good-compactification} for $(G,X,K)$. (In fact, we only need the assumptions for toroidal compactifications.) \par
Choose any affine formal open subscheme $\mathfrak{W}=\Spf(A,I)\sbst \mathfrak{X}_{\Upsilon}$. Let ${W}:=\Spec A$. The open embedding $j:\mathscr{S}_K\hookrightarrow \mathscr{S}_K^\Sigma$ induces an open embedding $j_{W}: {W}^0:= \mathscr{S}_K\times_{j,\mathscr{S}_K^\Sigma}W \hookrightarrow W$, which fits into the following commutative diagram
\begin{equation*}
    \begin{tikzcd}
    {W}^0\arrow[rr,hook,"j_{W}"]\arrow[d,"i^0"]&& {W}\arrow[d,"i"]\\
    \mathscr{S}_K\arrow[rr,hook,"j"]&& \mathscr{S}_K^\Sigma.
    \end{tikzcd}
\end{equation*}
Denote $\mathscr{S}:=\mathscr{S}_K$ and $\ca{X}:=\mathscr{S}_K^\Sigma$ for simplicity. We further assume:
\begin{itemize}
\item There is a $\G^c$-shtuka $\mathscr{P}_{\mathscr{S}}$ on $\mathscr{S}^{\Diamond/}$ extending the one $\mathscr{P}_\eta:=\mathscr{P}_{\mathscr{S}}|_{\mathscr{S}_{E}}$ on the generic fiber. Assume that $\mathscr{P}_\eta$ extends (uniquely) to a log shtuka $\mathscr{P}_\eta^{\can}$ on $\ca{X}_{E}^{\log \Diamond}$. 
Then there is a canonical isomorphism $\theta_\eta: i^{0,*}_{E}\mathscr{P}_\eta\iso i^{0,*}_{E}j^*_{E}\mathscr{P}_\eta^{\can}\xrightarrow{\sim} j_{W,E}^*i^*_{E}\mathscr{P}_\eta^{\can}$ on $W^{0,\Dia}_E$.
\item For any $\mathfrak{W}$, there is a log $\G^c$-shtuka $\mathscr{P}_{W}$ on ${W}^{\log \diamond}$ extending $i^*_{E}\mathscr{P}_\eta^\can$, such that there is a unique isomorphism 
$$\theta: i^{0,*}\mathscr{P}_{\mathscr{S}}\xrightarrow{\sim} j^*_{W}\mathscr{P}_{W}$$
on $W^{0,\dia}$ extending $\theta_\eta$.
\end{itemize}
Then there is a unique log $\G^c$-shtuka $\mathscr{P}^{\can}$ on $\ca{X}^{\log \diamond}=\ca{X}^{\log \Diamond}$ extending both $\mathscr{P}_\eta^{\can}$ and $\mathscr{P}_\mathscr{S}$.  
\end{lem}
\begin{proof}
Induction on the dimension of $\sigma$ in the pair $\Upsilon=[(\Phi,\sigma)]$. Note that $\dim \sigma=\mrm{codim}\ \ca{Z}_{\Upsilon,K}$. When $\sigma$ is trivial, by assumption, there is a shtuka $\mathscr{P}_{\mathscr{S}}$.\par
Assume that there is an extension of log shtuka $\mathscr{P}^k$ on $\ca{X}^{k,\log\Diamond/}$, whose underlying scheme $\ca{X}^k$ is set-theoretically the union of all strata $\ca{Z}_{\Upsilon}$ of codimension $\leq k$. The strata of codimension exactly $k+1$ are disjoint in $\ca{X}$. So we add one more stratum $\Upsilon':=[(\Phi',\sigma')]$ to $\ca{X}^k$ and show that there is an extension of the log shtuka on $(\ca{X}^k\cup \ca{Z}_{\Upsilon'})^{\log \Diamond/}$ or $(\ca{X}^k\cup \ca{Z}_{\Upsilon'})^{\log \diamond}$. (Here $\ca{X}^k\cup \ca{Z}_{\Upsilon'}$ denotes the normal subscheme in $\ca{X}$.)\par
Pick any affine open subscheme $V$ of $\ca{X}$ that is contained in $\ca{X}^k\cup \ca{Z}_{\Upsilon'}$ and intersects with $\ca{Z}_{\Upsilon'}$. 
After taking a suitable {\'e}tale cover and taking an affine open subscheme again, there is an affine scheme $U=\spec B$ and {\'e}tale morphisms $\mathsf{e}_1: U\to V=\spec A$ and $\mathsf{e}_2:U\to T= \spec R\times \mbf{E}_{K_{\Phi'}}(\sigma')$ for some normal flat $R$, such that the pullbacks of the stratifications on $V$ and $T$ via $\mathsf{e}_1$ and $\mathsf{e}_2$ coincide on $U$. Indeed, this follows from a refinement of Artin's approximation theorem (see \cite[Prop. 4.53]{Wu25}). We can choose finitely many such $U$ covering $\ca{Z}_{\Upsilon'}$.\par
Denote by $\tau_1$, ..., $\tau_{k+1}$ the $1$-dimensional faces of $\sigma'$. Denote by $D_i$ the divisor on $U$ defined by $\mathsf{e}_2^{-1}(T_{\tau_i})$ for $1\leq i\leq (k+1)$; denote by $U_i$ the complement $U\bss D_i$. Denote $U':=\bigcup_{i=1}^{k+1}U_i$.\par 
By construction, $D_i$ is defined by a principal ideal $I_i=(s_i)\sbst B$. Denote by $I$ the ideal generated by $I_i$ for all $i$. Denote by $\mathbf{V}$ the completion of $V$ along the closed stratum $Z\sbst V$ corresponding to $T_{\sigma'}$; the stratum $Z$ corresponds to an ideal $J$. Then $\spec \wat{A}_{J}\times_V U=\spec \wat{B}_{I}$ and $\mathsf{e}_1^{-1}V(J)=V(I)=\mathsf{e}_2^{-1}T_{\sigma'}$. Denote $\mbf{U}:=\spec \wat{B}_I$ and $\mbf{U}':=\spec \wat{B}_I\bss V(I)$.\par
We now claim that there is an isomorphism between $v$-sheaves
\begin{equation}\label{eq-iso-v-sheaf-gluing}(U')^{\log \diamond}\disju_{(\mbf{U}')^{\log\diamond}} \mbf{U}^{\log \diamond}\iso U^{\log \diamond}.\end{equation}
Here, the LHS is defined as the $v$-sheaf gluing $(U')^{\log\diamond}$ and $\mbf{U}^{\log \diamond}$. The map from left to right is injective. Since both side are $v$-sheaves by Theorem \ref{thm-log-diamond-v-sheaves}, for any affinoid $S^\sharp\in \Perfd$, it suffices to check surjectivity over a $v$-cover. 
For this, let $Z=\Spa(\prod_{t\in\ca{T}} C_t^+[1/\varpi_t],\prod_{t\in \ca{T}}C_t^+)\to S^\sharp$ be a $v$-cover constructed by a product of points, where $C_t^+[1/\varpi_t]$ are complete algebraically closed non-archimedean fields and $C_t^+$ are open and bounded valuation rings (see \cite[Ex. 1.1, Def. 1.2]{gleason2025specialization}). This cover is a strictly totally disconnected perfectoid space (see \cite[Prop. 1.6]{gleason2025specialization}), so all points are analytic.\par 
For each $C_t^+$, if there is a map $f:\spec C_t^+\to  U$, then the map factors through either $\mbf{U}$ or $U'$ according to the valuations of $\{s_i\}$ are all in the maximal ideal of $C^+_t$ or not: Indeed, for any $(C_t^+[1/\varpi_t],C^+_t)$, the valuation $v_t:C_t^+[1/\varpi_t]\to \bb{R}_{\geq 0}$ is given by the $\varpi_t$-topology by \cite[Lem. 4.2.2, Def. 4.2.3]{SW20}. Therefore, for any $s_i\in C^+_t$, $s_i$ lies in the maximal ideal if and only if $v_t(s_i)<1$ if and only if the whole $\wat{B}_{(s_i)}$ maps into $C_t^+$. If $v_t(s_i)=1$ then $s_i^{-1}\in C^+_t$. \par
For any $(S^\sharp,\ca{M}_{S^\sharp})\in U^{\log\dia}$, we pullback the log structure to $Z$ and denote it by $\ca{M}$. Then the log structure will also factor through $U'$ or $\mbf{U}$, since the log structures on $U'$ and $\mbf{U}$ are pulled back from the one on $U$. The claim is proved.\par
Note that there is a log shtuka $\mathscr{P}_{\mbf{U}}$ on $(\mbf{U})^{\log \diamond}$ given by pulling back the log shtuka $\mathscr{P}_{\mbf{V}}$ on $(\mbf{V})^{\log \diamond}$ and a log shtuka $\mathscr{P}^k_{U'}$ given by pulling back $\mathscr{P}^k$. The restrictions of both $\mathscr{P}^k_{U'}$ and $\mathscr{P}_{\mbf{U}}$ to $(\mbf{U}')^{\log \diamond}$ both extend the pullback of $\mathscr{P}^\can_\eta$ to $(\mbf{U}_{\bb{Q}_p}')^{\log \Diamond}$. Now it follows from Theorem \ref{thm-ext-shu-gen} that there is a canonical isomorphism between $\mathscr{P}^k_{U'}$ and $\mathscr{P}_{\mbf{U}}$ on $(\mbf{U}')^{\log \diamond}$. (Here we can use Theorem \ref{thm-ext-shu-gen} because $\mbf{U}'$ is a finite union of spectra of excellent Noetherian normal flat domains; for the fact that the completion is still excellent in this case, see \cite[Appendix A]{KSG21}, especially the table in the end.)
Combining (\ref{eq-iso-v-sheaf-gluing}) in the last paragraph with Proposition \ref{prop-def-shtukas}, we extend the log shtuka to $U$. Taking an {\'e}tale cover of $U$'s constructed as the third paragraph and by Theorem \ref{thm-ext-shu-gen} and the {\'e}tale descent of shtukas, we obtain an extension $\mathscr{P}^{k,\Upsilon'}$ on $(\ca{X}^k\cup \ca{Z}_{\Upsilon'})^{\log \diamond}$.\par 
Now, we choose another stratum $\Upsilon''=[(\Phi'',\sigma'')]$ with codimension $(k+1)$ and repeat the construction above. Since there are only finitely many such strata and finitely many choice of $k$, we finally construct an extension on $(\ca{X}^{k+1})^{\log \diamond}$, and we are done by induction.
\end{proof}
 \subsubsection{Extension along an affine toric embedding}\label{subsubsec-ext-affine-emb}
  Fix $\mbf{E}$ a split torus over $\bb{Z}$. Let $N=\mbf{X}_*(\mbf{E})$ and $M = \mbf{X}^*(\mbf{E})$. Let $\sigma \subset N \otimes \R$ be a smooth convex polyhedral cone, and $\sigma^{\vee} \subset M \otimes \R$ be its dual. We consider the toric embedding $\mbf{E} \hookrightarrow \mbf{E}(\sigma) = \Spec \Z[\sigma^{\vee} \cap M]$. 
We can naturally endow $\mbf{E}(\sigma)$ with a log structure associated with $D_{\sigma} = \mbf{E}(\sigma) \backslash \mbf{E} \hookrightarrow \mbf{E}(\sigma)$, and we denote it by $(\mbf{E}(\sigma), \ca{M}_{\sigma})$, with chart $\md{P}_{\sigma} \to \ca{M}_{\sigma}$ associated with the monoid $\sigma^{\vee} \cap M/(
\sigma^\vee\cap M)^\times$.\par
Let $Y$ be a normal scheme that is separated, flat, and of finite type over $\bb{Z}_p$; we also assume that the generic fiber of $Y$ is smooth. 

Let $X$ be an $\mbf{E}$-torsor over $Y$. The twisted toric embedding $X\hookrightarrow X(\sigma):=X\times^{\mbf{E}}\mbf{E}(\sigma)$ is equipped with a log structure associated with the chart $\md{P}_\sigma$.\par
Let $\mbf{E}_\infty':=\varprojlim_{x\mapsto x^p}\mbf{E}=\varprojlim \mbf{E}_n$, where $\mbf{E}_n$ denotes $\Spec \bb{Z}[\frac{1}{p^n}M]$. Write $\mbf{E}_n(\sigma)=\Spec \bb{Z}[\sigma^\vee\cap \frac{1}{p^n}M]$ and $\mbf{E}'_\infty(\sigma)=\Spec \bb{Z}[\sigma^\vee\cap M[\frac{1}{p}]]$.\par
Note that $X\hookrightarrow X(\sigma)$ is equipped with an equivariant $\mbf{E}$-action over $Y$. We write this action as $\bb{E}:\mbf{E}\times X\to X$ and $\bb{E}(\sigma):\mbf{E}\times X(\sigma)\to X(\sigma)$. Pre-composing with $p_\infty\times \mrm{id}: \mbf{E}_\infty'\times X\to \mbf{E}\times X$ (resp. $p_\infty\times \mrm{id}: \mbf{E}_\infty'\times X(\sigma)\to \mbf{E}\times X(\sigma)$), where $p_\infty:\mbf{E}'_\infty\to \mbf{E}$ is the canonical projection, we obtain an $\mbf{E}'_\infty$-action $\bb{E}_\infty: \mbf{E}'_\infty\times X\to X$ (resp. $\bb{E}_\infty(\sigma): \mbf{E}'_\infty\times X(\sigma)\to X(\sigma)$); the actions $\bb{E}_\infty$ and $\bb{E}_\infty(\sigma)$ are equivariant with respect to $X\hookrightarrow X(\sigma)$.\par
\begin{construction}\label{const-einf-action}\upshape
The cover $X\to Y$ is faithfully flat. There is a canonical isomorphism between schemes $i:\mbf{E}\times X\xrightarrow{\sim}X\times_Y X$; for any $\bb{Z}_p$-algebra $R$, $e\in \mbf{E}(R)$, $x\in X(R)$, and $y:=e\cdot x$, $i$ sends $(e,x)$ to $(x,y)$.
The descent datum $\varphi: X\times_Y X\xrightarrow{\sim} X\times_Y X$ for the cover $X\to Y$ is an isomorphism $\varphi=((-\mrm{id})\circ p_1\times \bb{E}): \mbf{E}\times X\xrightarrow{\sim}\mbf{E}\times X$ after pulled back via $i$.\par 
Let $i_\infty:=i\circ (p_\infty\times \mrm{id})$. We also pull back $\varphi$ to an isomorphism of $\mbf{E}_\infty'\times X$; i.e., we define $\varphi_\infty:=((-\mrm{id})\circ p_1\times \bb{E}_\infty): \mbf{E}_\infty'\times X\xrightarrow{\sim}\mbf{E}_\infty'\times X$.\par
For toric embeddings, we can also define 
$$\varphi_\infty(\sigma):=((-\mrm{id})\circ p_1\times \bb{E}_\infty(\sigma)):\mbf{E}'_\infty\times X(\sigma)\xrightarrow{\sim} \mbf{E}'_\infty\times X(\sigma).$$
Similarly, we define $\varphi(\sigma)$ by projection to $\mbf{E}$. 
The action of $\mbf{E}$ on $X(\sigma)$ is not free if $\sigma^\vee\cap M/(\sigma^\vee\cap M)^\times$ is nontrivial; in this case $\varphi(\sigma)$ is not an fpqc descent datum.\hfill$\square$
\end{construction}
\begin{definition}\label{def-einf-action-shtuka}
Let $(\ca{G},\mu)$ be as in \S\ref{subsec-shtuka-nonred}, where $\ca{G}$ is a quasi-parahoric model of a possibly non-reductive linear algebraic group (it was denoted as $\ca{P}$ in \S\ref{subsec-shtuka-nonred}).
Let $(\mathscr{P},\phi_{\mathscr{P}})$ be a $\ca{G}$-shtuka (bounded by $\mu$) on $X^\diamond$. We say that there is an (equivariant) $\mbf{E}_\infty'$-action on $\mathscr{P}$ if there is an isomorphism $$\varphi(\mathscr{P}):p_2^*(\mathscr{P},\phi_{\mathscr{P}})\xrightarrow{\sim}\bb{E}_\infty^*(\mathscr{P},\phi_{\mathscr{P}})$$ 
on $(\mbf{E}_\infty'\times X)^\diamond$ that satisfies the following conditions:
\begin{enumerate}
    \item Let $s: X \to \mbf{E}'_\infty\times X$ be the identity section of $\mbf{E}'_\infty$. Then $s^*\varphi(\mathscr{P})=\mrm{id}_{(\mathscr{P},\phi_\mathscr{P})}$;
    \item There is a commutative diagram (we omit $\phi_\mathscr{P}$ below):
    \begin{equation}\label{eq-associativity-shtuka}
        \begin{tikzcd}
        p_{23}^*p_2^*\mathscr{P}\arrow[r,"p_{23}^*\varphi(\mathscr{P})"]\arrow[d,equal]&p_{23}^*\bb{E}_\infty^*\mathscr{P}=(\mrm{id}\times \bb{E}_\infty)^*p_2^*\mathscr{P}\arrow[rr,"(\mrm{id}\times\bb{E}_\infty)^*\varphi(\mathscr{P})"]&&(\mrm{id}\times\bb{E}_\infty)^*\bb{E}^*_\infty\mathscr{P}\arrow[d,equal]\\
        (m\times\mrm{id})^*p_2^*\mathscr{P}\arrow[rrr,"(m\times\mrm{id})^*\varphi(\mathscr{P})"]&&&(m\times\mrm{id})^*\bb{E}_\infty^*\mathscr{P}
        \end{tikzcd}
    \end{equation}
over the commutative diagram
    \begin{equation*}
        \begin{tikzcd}
            \mbf{E}'_\infty\times \mbf{E}'_\infty\times X\arrow[rr,"(\mrm{id}\times \bb{E}_\infty)"]\arrow[d,"m\times\mrm{id}"]&&\mbf{E}'_\infty\times X\arrow[d,"\bb{E}_\infty"]\\
            \mbf{E}'_\infty\times X\arrow[rr,"\bb{E}_\infty"]&&X.
        \end{tikzcd}
    \end{equation*}
    
\end{enumerate}
\end{definition}
\begin{definition}\label{def-einf-action-shtuka-log}
Let $(\mathscr{P}(\sigma),\phi_{\mathscr{P}(\sigma)})$ be a $\ca{G}$-shtuka (bounded by $\mu$) on $X(\sigma)^{\log\diamond}$. We say that there is an $\mbf{E}_\infty'$-action on $\mathscr{P}(\sigma)$ if there is an isomorphism $$\varphi(\mathscr{P}(\sigma)):p_2^*(\mathscr{P}(\sigma),\phi_{\mathscr{P}(\sigma)})\xrightarrow{\sim}\bb{E}_\infty(\sigma)^*(\mathscr{P}(\sigma),\phi_{\mathscr{P}(\sigma)})$$ 
on $(\mbf{E}_\infty'\times X(\sigma))^{\log\diamond}$ that satisfies similar conditions above. 
\end{definition}
The following proposition was inspired by \cite[4.1.1]{Har89} (cf. \cite[Prop. 1.3.5]{HZ01}).
\begin{thm}\label{thm-ext-toric-equiv}
Let $(\ca{G},\mu)$ be as in \S\ref{subsec-shtuka-nonred}. With the conventions above, there is a natural equivalence of categories between
\begin{itemize}
    \item The category $\Sht^\diamond_{\ca{G},\mu,\mbf{E}_\infty'}(X)$ of $\ca{G}$-shtukas $(\mathscr{P},\phi_{\mathscr{P}})$ with one leg bounded by $\mu$ on $X^\diamond$ that are equipped with equivariant $\mbf{E}_\infty'$-actions;
    \item The category $\Sht^\diamond_{\ca{G},\mu,\mbf{E}_\infty'}(X(\sigma),\ca{M}_\sigma)$ of log $\ca{G}$-shtukas $(\mathscr{P}(\sigma),\phi_{\mathscr{P}(\sigma)})$ with one leg bounded by $\mu$ on $X^{\log \diamond}$ that are equipped with equivariant $\mbf{E}_\infty'$-actions.
\end{itemize}
The morphisms in the two categories are those preserving the $\mbf{E}_\infty'$-actions. In fact, we show that every object in the first category can be uniquely extended to the second one.
\end{thm}
Note that the last statement in the theorem is stronger than saying the restriction functor is essentially surjective.\par 
We prove this theorem in \S\ref{subsubsec-ext-affine-toric-continued}. An immediate corollary is:
\begin{cor}\label{cor-ext-toric-emb}
Let $(G,X)$ be any Shimura datum. If the $\ca{P}_\Phi^*$-shtuka $\mathscr{P}_E$ determined by $\bb{P}^*_{K_\Phi}$ extends to a $\ca{P}_\Phi^*$-shtuka $\mathscr{P}$ on $(\Delta^\circ_{\Phi,K}\bss \mathscr{S}_{K_\Phi})^\diamond$, it uniquely extends to a log shtuka $\mathscr{P}^\can$ on $(\Delta^\circ_{\Phi,K}\bss\mathscr{S}_{K_\Phi}(\sigma))^{\log \diamond}$.
\end{cor}
\begin{proof}
There is a map between $\ul{\ca{P}^*_\Phi(\bb{Z}_p)}$-torsors $$\bb{P}_{K_\Phi}^*\to \bb{P}_{K_{\Phi^*}}(P_\Phi^*,D_{\Phi^*})$$
over a map
$$\Delta^\circ_{\Phi,K}\bss \sh_{K_\Phi}(P_\Phi,D_\Phi)\to \sh_{K_{\Phi^*}}(P_\Phi^*,D_{\Phi^*}),$$
where $\Phi^*$ is the one defined below Definition \ref{def-pstar}. \par
The first scheme is a torus torsor under $\mbf{E}_{\K_\Phi}$, which corresponds to a lattice (see \cite[1.3.3]{Wu25})
$$\bm{\Lambda}_{K_\Phi}:=p_2(Z(\bb{Q})\times U_\Phi(\A)\cap g_\Phi Kg_\Phi^{-1});$$
since $K$ is neat, this lattice maps to 
$$\bm{\Lambda}_{K_\Phi}^*:=p_2((1\times U_\Phi(\A))\cap g_\Phi K^c g_\Phi^{-1}).$$
Thus, the map $$\Delta^\circ_{\Phi,K}\bss \sh_{K_\Phi}(P_\Phi,D_\Phi)\to \sh_{K_{\Phi^*}}(P_\Phi^*,D_{\Phi^*})$$ is equivariant under the isogeny between tori $\mbf{E}_{\K_\Phi}\to \mbf{E}_{K_{\Phi^*}}$.\par
Now we use the proof of Lemma \ref{lem-torus-action-tower-pro-p-gen} {verbatim}, replacing $\sh_K(P,\ca{X})\to \sh_{K^c}(P^c,\ca{X}^c)$ with $\Delta^\circ_{\Phi,K}\bss \sh_{K_\Phi}(P_\Phi,D_\Phi)\to \sh_{K_{\Phi^*}}(P_\Phi^*,D_{\Phi^*})$. We then see that there is an $\mbf{E}'_{\K_\Phi,\infty}:=\varprojlim_{x\mapsto x^p}\mbf{E}_{\K_\Phi}$-action. That is, there is a diagram
\begin{equation*}
    \begin{tikzcd}
        \mbf{E}'_{\K_\Phi,\infty}\times\bb{P}^*_{K_\Phi}\arrow[r]\arrow[d]&\bb{P}^*_{K_\Phi}\arrow[d]\\
        \mbf{E}_{\K_\Phi}\times \Delta^\circ_{\Phi,K}\bss\sh_{K_\Phi}\arrow[r]&\Delta^\circ_{\Phi,K}\bss\sh_{K_\Phi}.
    \end{tikzcd}
\end{equation*}

Thus, by \cite[Thm. 2.7.7]{PR24} and the Tannakian formalism, $\mathscr{P}$ is equipped with an $\mbf{E}_{\K_\Phi,\infty}'$-action. Note that although $\mbf{E}'_{\K_\Phi,\infty}\times \mathscr{S}$ is not of finite type, the argument in \emph{loc. cit.} still works. Indeed, denote $\mathscr{S}:=\Delta^\circ_{\Phi,K}\bss\mathscr{S}_{K_\Phi}$. It suffices to work with an affine open $\Spec A\sbst \mathscr{S}$. By the proof of \emph{loc. cit.}, we find a perfectoid cover $\Spa(\wdtd{R},\wdtd{R}^+)\to \Spa(\wat{A}[1/p],\wat{A})$ by taking the completion of finite {\'e}tale extensions of $\wat{A}[1/p]$ in its fraction field. Then we can replace $\Spa(\wdtd{R},\wdtd{R}^+)$ with the perfectoid space $\Spa(\wdtd{R}\langle\md{R}\rangle,\wdtd{R}^+\langle\md{R}\rangle)$ (see \cite[Lem. 2.2.15]{DLLZ23}), where $\md{R}$ is the character group of $\mbf{E}'_{\K_\Phi,\infty}$; and the rest of the proof in \emph{loc. cit.} works for $\mbf{E}'_{\K_\Phi,\infty}\times \Spec A$.\par 
We then obtain $\mathscr{P}(\sigma)$ extending $\mathscr{P}$ by applying Theorem \ref{thm-ext-toric-equiv}. The uniqueness follows from Corollary \ref{cor-ext-shu-gen}.
\end{proof}
\subsubsection{Proof of Theorem \ref{thm-ext-toric-equiv}}\label{subsubsec-ext-affine-toric-continued}
There is a natural restriction functor 
$$\mrm{Res}: \Sht^\diamond_{\ca{G},\mu,\mbf{E}_\infty'}(X(\sigma),\ca{M}_\sigma)\to \Sht^\diamond_{\ca{G},\mu,\mbf{E}_\infty'}(X).$$
It is fully faithful by Corollary \ref{cor-ext-shu-gen}. We show that the last statement in the theorem is true; in particular, the restriction functor is essentially surjective. \par
Let $\md{P}:=\md{P}_\sigma$ and $\md{Q}:=\sigma^\vee\cap M$. Let $\md{Q}_\infty:=\md{Q}[\frac{1}{p}]$. 
There is a natural injective homomorphism between monoids $j:\md{Q}\hookrightarrow \md{Q}_\infty$. Also, there is an injective homomorphism
$\psi: \md{Q}\to \md{Q}^{\gp}\oplus \md{Q}_\infty$ given by the embeddings $\md{Q}\hookrightarrow\md{Q}^{\gp}$ and $j$.\par
\begin{construction}\label{const-coaction}\upshape
There is a natural coaction of $\md{Q}^\gp_\infty$ on $\md{Q}^{\gp}\oplus \md{Q}_\infty$.\par 
In fact, the action of $\mbf{E}_{\infty}'$ on $\mbf{E}_\infty'(\sigma)$ induces a coaction of $\md{Q}^{\gp}_\infty$ on $\md{Q}_\infty$. Explicitly, this is given by 
$$\md{Q}_\infty\to \md{Q}_\infty\oplus \md{Q}_\infty^{\gp};\ x\mapsto (x,x).$$
Similarly, the action of $\mbf{E}_\infty'$ on $\mbf{E}$ given by $\mbf{E}_\infty'\times \mbf{E}\xrightarrow{p_\infty\times \mrm{id}}\mbf{E}\times\mbf{E}\xrightarrow{\text{multiplication}}\mbf{E}$ also induces a comultiplication $\md{Q}^{\gp}\to \md{Q}^{\gp}\oplus \md{Q}_\infty^{\gp}$.\par
Combining the two (co)actions, we obtain a map between monoids 
$$d_\infty(\sigma):\md{Q}^{\gp}\oplus \md{Q}_\infty\to \md{Q}^{\gp}\oplus \md{Q}_\infty\oplus \md{Q}_\infty^{\gp};\ (x,y)\mapsto (x,y,-x+y).$$
This corresponds to an action $d_\infty(\sigma):\mbf{E}_\infty'\times\mbf{E}\times\mbf{E}_\infty'(\sigma)\to \mbf{E}\times \mbf{E}'_\infty(\sigma)$ that pointwisely sends $(g,a,b)$ to $(g^{-1}a,gb)$. Note that the quotient of $\mbf{E}\times \mbf{E}'_\infty(\sigma)$ by this action is $\mbf{E}\times^{\mbf{E}'_\infty}\mbf{E}'_\infty(\sigma)\iso \mbf{E}(\sigma)$.\par
Denote by $d_\infty$ the $\mbf{E}'_\infty$-action $d_\infty(\sigma)$ restricted to $\mbf{E}\times \mbf{E}_\infty'(\sigma)$. \hfill$\square$
\end{construction}
We need a lemma on saturated fiber products.
\begin{lem}\label{lem-sat-product-abs}
We have that $(\md{Q}_\infty\oplus_{j,\md{Q},\psi}(\md{Q}^{\gp}\oplus \md{Q}_\infty))^{\sat}\iso \md{Q}_\infty\oplus \md{Q}_\infty^{\gp}$. Moreover, the $\md{Q}_\infty^{\gp}$-coaction on the second factor $\md{Q}^{\gp}\oplus \md{Q}_\infty$ induces the comultiplication on the right-hand side on the second factor of $\md{Q}_\infty\oplus \md{Q}_\infty^\gp$.
\end{lem}
\begin{proof}
We first compute $(\md{Q}_\infty\oplus_{j,\md{Q},\psi}(\md{Q}^\gp\oplus \md{Q}_\infty))^{\mrm{int}}$. 
Since $\md{Q}_\infty\oplus \md{Q}^{\gp}\oplus \md{Q}_\infty$ is integral and $(\md{Q}_\infty\oplus \md{Q}^{\gp}\oplus \md{Q}_\infty)^{\gp}=\md{Q}_\infty^{\gp}\oplus \md{Q}^{\gp}\oplus \md{Q}_\infty^{\gp}$, by \cite[I. Prop. 1.3.4]{Ogu18}, $(\md{Q}_\infty\oplus_{j,\md{Q},\psi}(\md{Q}^\gp\oplus \md{Q}_\infty))^{\mrm{int}}$ is the image of $\md{Q}_\infty\oplus \md{Q}^{\gp}\oplus \md{Q}_\infty$ in $\md{Q}_\infty^{\gp}\oplus_{j^\gp,\md{Q}^{\gp},\psi^\gp}(\md{Q}^\gp\oplus \md{Q}_\infty^\gp)$. 
Therefore, 
\begin{equation}\label{eq-int-mono}(\md{Q}_\infty\oplus_{j,\md{Q},\psi}(\md{Q}^\gp\oplus \md{Q}_\infty))^{\mrm{int}}\sbst(\md{Q}_\infty+\md{Q}^\gp)\oplus_{j^{\gp},\md{Q}^\gp,\psi^\gp}(\md{Q}^\gp\oplus (\md{Q}_\infty+\md{Q}^\gp)),\end{equation}
where $(\md{Q}_\infty+\md{Q}^{\gp})$ denotes the monoid in $\md{Q}_\infty^{\gp}$ generated by $\md{Q}_\infty$ and $\md{Q}^{\gp}$.\par
Denote by $(a,b,c)$ an element in $\md{Q}_\infty\oplus \md{Q}^{\gp}\oplus \md{Q}_\infty$. Then $(a,b,c)\sim (a+b,0,c-b)$ in the RHS of (\ref{eq-int-mono}) since $(0,b,b)\sim (b,0,0)$ by amalgamated product.
Hence, $(\md{Q}_\infty\oplus_{j,\md{Q},\psi}(\md{Q}^\gp\oplus \md{Q}_\infty))^{\mrm{int}}$ is generated by $(a,0,0)$, $(0,0,c)$, where $a,c\in \md{Q}_\infty$, and $(b,0,-b)$, where $b\in \md{Q}^{\gp}$. \par
Thus, $(\md{Q}_\infty\oplus_{j,\md{Q},\psi}(\md{Q}^\gp\oplus \md{Q}_\infty))^{\mrm{sat}}=\md{Q}_{\infty}\oplus \md{Q}_\infty^{\gp}$. Indeed, for any $(x,0,y)\in (\md{Q}_\infty^\gp,0,\md{Q}_\infty^\gp)$ such that $(nx,0,ny)=(a+b,0,c-b)$, by enlarging $n$, we assume that $nx$ and $ny$ are in $\md{Q}^{\gp}$; and therefore, $a$ and $c$ are in $\md{Q}^\gp\cap \md{Q}_\infty=\md{Q}$. 
We then write $(a+b,0,c-b)=(a+b+c-c,0,c-b)=(a+c,0,0)+(b-c,0,c-b)$. Hence, $y$ can be any element in $\md{Q}_\infty^\gp$ and $x+y$ can be any element in $\md{Q}_\infty$.
Note that the generators of $(\md{Q}_\infty\oplus_{j,\md{Q},\psi}(\md{Q}^\gp\oplus \md{Q}_\infty))^{\mrm{sat}}$ are $(u,0,0)\in (\md{Q}_\infty,0,0)$ and $v\in \md{Q}^{\gp}_\infty\xrightarrow{(-v,v)}(\md{Q}_\infty^\gp,0,\md{Q}^\gp_\infty)$. \par
Now we check the last sentence. The induced $\md{Q}^\gp_\infty$-coaction on the RHS of (\ref{eq-int-mono}) sends $(u,v,w)\in (\md{Q}_\infty+\md{Q}^\gp)\oplus_{j^{\gp},\md{Q}^\gp,\psi^\gp}(\md{Q}^\gp\oplus (\md{Q}_\infty+\md{Q}^\gp))$ to $(u,v,w,-v+w)$; note that this is compatible with the amalgamated product as the coactions on the last two factors are canceled and therefore coincide with the trivial action on the first factor. Hence, the induced coaction on $(-v,0,v)$ is the comultiplication on the third factor, while the coaction on $(u,0,0)$ is trivial. 
\end{proof}
We now go back to the proof. 
Note that 
\begin{lem}\label{lem-ess-surj-etale-cover}
It suffices to show the last statement of Theorem \ref{thm-ext-toric-equiv} over an {\'e}tale cover $\wdtd{Y}$ of $Y$.
\end{lem}
\begin{proof}
Let $\wdtd{X}:=\wdtd{Y}\times_Y X$ and $\wdtd{X}(\sigma):=\wdtd{Y}\times_Y X(\sigma)$. 
Suppose that there is a log shtuka $\wdtd{\mathscr{P}}(\sigma)$ in $\Sht^\dia_{\G,\mu,\mbf{E}'_\infty}(\wdtd{X}(\sigma),\ca{M}_{\wdtd{X}(\sigma)})$ that extends some $\wdtd{\mathscr{P}}$ in $\Sht_{\G,\mu,\mbf{E}'_\infty}^\diamond(X)$ to $\wdtd{X}(\sigma)^{\log\dia}$.
Since $\wdtd{\mathscr{P}}$ is the pullback of some $\mathscr{P}$ in $\Sht^\dia_{\G,\mu,\mbf{E}'_\infty}(X)$. 
By Corollary \ref{cor-ext-shu-gen}, the descent datum of $\wdtd{\mathscr{P}}$ from $\wdtd{X}$ to $X$ extends uniquely to a strict {\'e}tale descent datum of $\wdtd{\mathscr{P}}(\sigma)$ from $\wdtd{X}(\sigma)$ to $X(\sigma)$, as $\sigma$ is a smooth cone decomposition and the generic fiber of $X$ is smooth.
By strict {\'e}tale descent of log shtukas (see Lemma \ref{lem-str-etale-descent}), there is a log shtuka $\mathscr{P}(\sigma)$ over $(X(\sigma))^{\log\dia}$ descending $\wdtd{\mathscr{P}}(\sigma)$.
\end{proof}
By Lemma \ref{lem-ess-surj-etale-cover}, we can assume that $X$ is a trivial torsor over $Y$, and that $X\times \mbf{E}_\infty'(\sigma)\to X(\sigma)$ admits a chart $\md{Q}\to \md{Q}^{\gp}\oplus \md{Q}_\infty$.\par
Let $(S^{\sharp}=\Spa(R,R^+), \ca{M}_{S^{\sharp}}, f) \in (X(\sigma), \ca{M}_{\sigma})^{\log\diamond}(S)$, where $$f: (S^\sharp,\ca{M}_{S^\sharp})\to(S^{\sharp,+}:=\Spec R^{+}, \ca{M}_{S^{\sharp,+}}) \xrightarrow{f^+} (X(\sigma), \ca{M}_{\sigma}).$$ 
We assume that $\ca{M}_{S^\sharp}=\ca{M}_{S^\sharp}^\can$ by Proposition \ref{prop-shu-can-obj}. And $f$ corresponds to an $R^+$-point of $X(\sigma)$ denoted by $f^+:(\Spec R^+,\ca{M}_{S^{\sharp,+}})\to (X(\sigma),\ca{M}_\sigma)$ that admits a chart $\md{Q}\to \md{Q}_\infty$ by Lemma \ref{lem-compatible-charts}.

  Consider the fiber product in the category of saturated log schemes
  \begin{equation}
\begin{tikzcd}
	{\tilde{T}^{\sharp,+}} & {X \times \mbf{E}_{\infty}'(\sigma)} \\
	{\Spec R^{+}} & {X(\sigma).}
	\arrow[from=1-1, to=1-2]
	\arrow[from=1-1, to=2-1]
	\arrow[from=1-2, to=2-2,"\alpha_\infty"]
	\arrow[from=2-1, to=2-2]
\end{tikzcd}
  \end{equation}
  \begin{lem}\label{lem-sat-fiber-prod-perfd}
     The $p$-adic completion of $\wdt{T}^{\sharp,+}$ is the spectrum of the integral perfectoid ring $R^{+}\langle \md{Q}_\infty^\gp\rangle$ equipped with an $\mbf{E}_{\infty}'$-action.
  \end{lem}
  \begin{proof}
Recall that we can work {\'e}tale locally over $Y$ by Lemma \ref{lem-ess-surj-etale-cover}, and assume that $X=Y\times \mbf{E}$ and that $Y$ is affinoid. Set $X_\infty(\sigma):=X(\sigma)\times_{\Spec \bb{Z}_p[\md{Q}]}\Spec \bb{Z}_p[\md{Q}_\infty].$
Then 
\begin{equation*}
    \begin{split}
&(S^{\sharp,+},\ca{M}_{S^{\sharp,+}})\times^{\sat}_{X(\sigma)}(X\times \mbf{E}_\infty'(\sigma))\\
       & =(\Spec R^{+},\ca{M}_{S^{\sharp,+}})\times_{X_\infty(\sigma)} X_\infty(\sigma)\times^{\sat}_{X(\sigma)}(X\times \mbf{E}_\infty'(\sigma))\\
      &  =(\Spec R^{+},\ca{M}_{S^{\sharp,+}})\times_{X_\infty(\sigma)} X_\infty(\sigma)[ \md{Q}_\infty^\gp].
    \end{split}
\end{equation*}
The second line follows from Lemma \ref{lem-compatible-charts} and the third line follows from Lemma \ref{lem-sat-product-abs}.
  \end{proof}
We now finish the proof. We omit $\phi_\mathscr{P}$ in the proof. Suppose that we are given $\mathscr{P}$ in $\Sht^\diamond_{\ca{G},\mu,\mbf{E}'_\infty}(X)$. 
Set $\wdtd{T}^\sharp:=\Spa(R\langle\md{Q}_\infty^\gp\rangle,R^+\langle\md{Q}_\infty^\gp\rangle)$, which is affinoid perfectoid by \cite[Lem. 2.2.5]{DLLZ23}. \par
Note that there is a commutative diagram
\begin{equation}\label{eq-commu-Einf-X}
\begin{tikzcd}
	{\mathbf{E}_\infty'\times\mathbf{E}_\infty'\times X} && {\mathbf{E}_{\infty}' \times X} && X & \\
	{\mathbf{E}_\infty'\times\mathbf{E}_\infty'\times X} && {\mathbf{E}_{\infty}' \times X} && {X.} & {}
	\arrow["{d_{\infty}':=d_{\infty}\circ(i\times \identity \times \identity)}"{pos=0.6}, shift left=2, from=1-1, to=1-3]
	\arrow["{p_{23}}"'{pos=0.6}, shift right=2, from=1-1, to=1-3]
	\arrow["{\identity \times \varphi_{\infty}}", from=1-1, to=2-1]
	\arrow["{\mathbb{E}_{\infty}}", from=1-3, to=1-5]
	\arrow["{\varphi_{\infty}}", from=1-3, to=2-3]
	\arrow[equals, from=1-5, to=2-5]
	\arrow["{m \times \identity}", shift left=2, from=2-1, to=2-3]
	\arrow["{p_{23}}"', shift right=2, from=2-1, to=2-3]
	\arrow["{p_2}", from=2-3, to=2-5]
\end{tikzcd}
\end{equation}
Here, $d_\infty$ also denotes the action of the first factor $\mbf{E}_\infty'$ on $\mbf{E}_\infty'\times X$ that is a (positive) multiplication on the first factor $\mbf{E}_\infty'$ and $-\bb{E}_\infty$ on the second factor $X$; $m$ (resp. $i$) denotes the multiplication (resp. inverse) on $\mbf{E}_\infty'$; $d_{\infty}'$ maps $(a, b, x)$ to $(a^{-1}b, ax)$.\par

Let $\mathscr{P}_2$ be the pullback of $\mathscr{P}$ to $(\mbf{E}'_\infty\times X)^\diamond$ along the projection $p_2$ to the second factor. Note that $p_2\circ d_{\infty}' = p_2\circ\varphi_{\infty}\circ p_{13}=\bb{E}_\infty\circ p_{13}$. From this, we obtain a descent datum $\varphi(\mathscr{P}_2):p_{23}^*\mathscr{P}_2\xrightarrow{\sim} d_{\infty}^{\prime*}\mathscr{P}_2$ by pulling back $\varphi(\mathscr{P}): p_2^*\mathscr{P} \rightiso \varphi_{\infty}^*p_2^*\mathscr{P}=\bb{E}_\infty^*\mathscr{P}$ along $p_{13}$. The cocycle condition is given by the group action $\varphi(\mathscr{P})$. 
Therefore, $\mathscr{P}_2$ descends to $\mathscr{P}'$ on $X$ with this descent datum. \par
On the other hand, denote by $\mathscr{P}_3 := \mathbb{E}_{\infty}^*\mathscr{P} = \varphi_{\infty}^*p_2^*\mathscr{P}$, we have a standard descent datum $\varphi(\mathscr{P}_3): p_{23}^*\mathscr{P}_3\xrightarrow{\sim}d_{\infty}^{\prime*}\mathscr{P}_3$ given by the pullback of $\mathscr{P}$ along $\bb{E}_\infty$. 
Note that $\varphi(\mathscr{P})$ gives an isomorphism between descent data $(\mathscr{P}_2, \varphi(\mathscr{P}_2)) \to (\mathscr{P}_3, \varphi(\mathscr{P}_3))$, and it descends to an isomorphism $\mathscr{P}' \rightiso \mathscr{P}$. Indeed, define $\bm{\alpha}: \mbf{E}'_\infty\times\mbf{E}'_\infty\times X\xrightarrow{\sim}\mbf{E}'_\infty\times \mbf{E}'_\infty\times X$ by the assignment $(a,b,x)\mapsto (a^{-1}b,a,x)$. Then $p_{23}\circ\bm{\alpha}=p_{13}$, $(\mrm{id}\times\bb{E}_\infty)\circ\bm{\alpha}=d'_\infty$ and $(m\times \mrm{id})\circ\bm{\alpha}=p_{23}$. Pulling back (\ref{eq-associativity-shtuka}) along $\bm{\alpha}$, we obtain a commutative diagram
\begin{equation}\label{eq-associativity-pullback}
    \begin{tikzcd}
       p_{13}^*p_2^*\mathscr{P}\arrow[r,"p_{13}^*\varphi(\mathscr{P})"]\arrow[d,equal]&p_{13}^*\bb{E}_\infty^*\mathscr{P}=d^{\prime,*}_\infty p_2^*\mathscr{P}\arrow[rr,"d^{\prime,*}_\infty\varphi(\mathscr{P})"]&&d^{\prime,*}_\infty\bb{E}^*_\infty\mathscr{P}\arrow[d,equal]\\
p_{23}^*p_2^*\mathscr{P}\arrow[rrr,"p_{23}^*\varphi(\mathscr{P})"]&&&p_{23}^*\bb{E}_\infty^*\mathscr{P}.  
    \end{tikzcd}
\end{equation}
We write (\ref{eq-associativity-pullback}) above in the following form:
\begin{equation*}
    \begin{tikzcd}
        p^*_3\mathscr{P}\arrow[rr,"p_{13}^*\varphi(\mathscr{P})"]\arrow[d,"p_{23}^*\varphi(\mathscr{P})"]&&p_{13}^*\bb{E}_{\infty}^*\mathscr{P}=d^{\prime,*}_\infty p_2^*\mathscr{P}\arrow[d,"d^{\prime,*}_\infty\varphi(\mathscr{P})"]\\
        p_{23}^*\bb{E}_{\infty}^*\mathscr{P}=(\mrm{id}\times\varphi_\infty)^*p_3^*\mathscr{P}\arrow[rr,equal]&&(\mrm{id}\times \varphi_\infty)^*(m\times \mrm{id})^*p_2^*\mathscr{P}=d^{\prime,*}_\infty\bb{E}_{\infty}^*\mathscr{P};
    \end{tikzcd}
\end{equation*}
this implies the claim that $\varphi(\mathscr{P})$ is an isomorphism between the descent data of $\mathscr{P}_2$ and $\mathscr{P}_3$.\par
Moreover, we show that $\mathscr{P}'\xrightarrow{\sim} \mathscr{P}$ is an automorphism of $\mathscr{P}$. Consider the canonical section $s: X \to \mbf{E}_{\infty}' \times X $ sending $x$ to $(e, x)$, it is a section of both $\mathbb{E}_{\infty}$ and $p_2$. Then, 
$$\mathscr{P}' = s^*\mathbb{E}_{\infty}^*\mathscr{P}' = s^*\mathscr{P}_2 = s^*p_2^*\mathscr{P} = \mathscr{P}.$$
Also, the descended automorphism $\mathscr{P} \rightiso \mathscr{P}$ is an identity, as it can be realized as the pullback $s^*\varphi(\mathscr{P}): \mathscr{P} \rightiso \mathscr{P}$ of the morphism $\varphi(\mathscr{P})$ between the descent data; here we use the fact that the following composition is the identity:
$$ X \stackrel{s}{\to} \mbf{E}_{\infty}' \times X \stackrel{\varphi_{\infty}}{\to} \mbf{E}_{\infty}' \times X \stackrel{p_2}{\to} X. $$
This automorphism is an identity by the first condition in Definition \ref{def-einf-action-shtuka}.

Next, we extend the first row of (\ref{eq-commu-Einf-X}) to 
\begin{equation}\label{eq-esigma-x}
\begin{tikzcd}
	{\mathbf{E}_\infty'\times\mathbf{E}_\infty'(\sigma)\times X} && {\mathbf{E}_{\infty}'(\sigma) \times X} && {X(\sigma).} & \\
	&&&&& {}
	\arrow["{d_{\infty}'(\sigma):=d_{\infty}(\sigma)\circ(i\times \identity \times \identity)}"{pos=0.6}, shift left=2, from=1-1, to=1-3]
	\arrow["{p_{23}}"'{pos=0.6}, shift right=2, from=1-1, to=1-3]
	\arrow["{\mathbb{E}_{\infty}(\sigma)'}", from=1-3, to=1-5]
\end{tikzcd}
\end{equation}
Here, $\bb{E}_\infty(\sigma)'$ is defined by $\mbf{E}'_\infty(\sigma)\times X\to \mbf{E}(\sigma)\times X\to \mbf{E}(\sigma)\times^{\mbf{E}}X=X(\sigma)$ and it extends the map $\bb{E}_\infty$.
Again, we equip $\mbf{E}_\infty'(\sigma)$ with the positive multiplication of $\mbf{E}_\infty'$ and $X$ with the negative action. Thus, $d_\infty(\sigma)$ extends $d_\infty$ (resp. $d_\infty'(\sigma)$ extends $d_\infty'$).\par
Denote by $\mathscr{Q}$ the pullback of $\mathscr{P}$ to $(\mbf{E}'_\infty(\sigma)\times X)^{\log\dia}$ via the projection to the second factor $X$.\par
For any $f^+:(\Spec R^{+},\ca{M}_{R^{+}})\to X(\sigma)$ as above, the evaluation $\mathscr{Q}_{\wdtd{T}^{\sharp}}$ of $\mathscr{Q}$ at $\wdtd{T}^\sharp\to\wdtd{T}^{\sharp,+}\to X\times\mbf{E}'_\infty(\sigma)$ is equipped with an $\mbf{E}_\infty'$-action given by the pullback of $d_\infty'(\sigma)$ by construction and by Lemma \ref{lem-sat-fiber-prod-perfd}.
By $v$-descent of shtukas (see \cite[Prop. 19.5.3]{SW20}), the shtuka $\mathscr{Q}_{\wdtd{T}^\sharp}$ descends to a shtuka $\mathscr{Q}_{S^\sharp}$.\par 
This assignment can be easily checked to be functorial, and we denote the obtained shtuka on $X(\sigma)^{\log \dia}$ by $\mathscr{P}(\sigma)$; the $\mbf{E}_\infty'$-action on $\mathscr{P}(\sigma)$ is descended from that of $\mathscr{Q}$ by the multiplication action of $\mbf{E}_\infty'$ on $\mbf{E}'_\infty(\sigma)$. Moreover, when we restrict this construction to $(\mbf{E}'_\infty\times X)^{\dia}$, we have seen that we obtain a shtuka $\mathscr{P}' = \mathscr{P}$. 
Now we have completed the proof of last statement in Theorem \ref{thm-ext-toric-equiv}.       \hfill$\square$
\subsection{Canonical integral models of compactifications}\label{subsec-can-mod-II}
\subsubsection{Assumptions}
As mentioned in Remark \ref{rk-context}, the extension of shtukas to toroidal compactifications $\mathscr{S}_K^\Sigma$ can be achieved without understanding the geometry of the interior $\mathscr{S}_K$.
To highlight this point, we present a sequence of assumptions and definitions at slightly different levels of generality.\par
Let $(G, X)$ be a Shimura datum, let $\GG$ be a quasi-parahoric model of $G = G_{\rQ_p}$, and let $K_p = \GG(\Z_p)$. Let $\lrbracket{\Shum{K_pK^p}(G, X)}_{K^p \subset G(\A^p)}$ be a family of integral models of $\lrbracket{\shu{K_pK^p}(G, X)}_{K^p \subset G(\A^p)}$.

To study the boundary stratification of the integral model of the toroidal compactification and to construct log shtukas on it, we assume:
\begin{assumption}\label{ass-well-position_mixed}
  For each neat $K^p \subset G(\A^p)$,
   \begin{enumerate}
       \item $\Shum{K}(G, X)$ has a good compactification theory as in Axiom \ref{axiom-good-compactification}. (In fact, we only need the assumptions for toroidal compactifications.)
       \item For any $\Phi\in\ca{CLR}(G,X)$, the $\PP^*_{\Phi}$-shtuka $(\PPs_{\Phi, E}, \phi_{\PPs_{\Phi, E}})$ over $\shu{K_{\Phi}}(P_{\Phi}, D_{\Phi})^{\Dia}$ with one leg bounded by $\mu^*_{\Phi}$ associated with the push-out of the de Rham pro-\'etale $\underline{\PP^{c}_{\Phi}(\Z_p)}$-torsor $\PPp_{K_{\Phi}} \to \shu{K_{\Phi}}(P_{\Phi}, D_{\Phi})$ via $\PP^c_{\Phi} \to \PP^*_{\Phi}$ extends to a $\PP^*_{\Phi}$-shtuka $(\PPs_{\Phi}, \phi_{\PPs_{\Phi}})$ over $\Shum{K_{\Phi}}(P_{\Phi}, D_{\Phi})^{\Dia/}$.
   \end{enumerate}
   \end{assumption}

   To consider the boundary stratification on the integral model of the minimal compactification and prove the well-positionedness of various strata on the special fiber, we further assume:
   \begin{assumption}\label{ass-well-position}
    For each neat $K^p \subset G(\A^p)$, assumptions \ref{ass-well-position_mixed} hold. Moreover, for each $[\Phi] \in \ca{CLR}_K(G, X)$, $\Delta_{\Phi, K}^{\circ}\backslash \Shum{\ovl{K}_{\Phi}}(\ovl{P}_{\Phi}, \ovl{D}_{\Phi}) \to \Delta_{\Phi, K}^{\circ}\backslash \Shum{K_{\Phi, h}}(G_{\Phi, h}, D_{\Phi, h})$ is an abelian scheme torsor.
   \end{assumption}

   To consider functoriality and uniqueness of integral models of toroidal (resp. minimal compactifications), we make the following definition:
   \begin{definition}\label{def-PR-int-mod}
      Fix $K'=K_p K^{\prime,p}$ for $K_p$ a quasi-parahoric subgroup and $K^{\prime,p}$ neat open compact. Assume $\Shum{K'}:=\Shum{K'}(G, X)$ has a good compactification theory as in Axiom \ref{axiom-good-compactification} and let $\Sigma$ be a projective smooth cone decomposition for $(G,X,K')$.\par 
      We say that $\{\mathscr{S}_{K}^{\Sigma}\}_{K^p}$ (resp. $\{\mathscr{S}_{K}^{\min}\}_{K^p}$) is a system of \textbf{canonical integral models} (in the sense of Pappas-Rapoport) of toroidal compactifications $\{\shu{K}^{\Sigma}\}_{K^p}$ (resp. minimal compactifications $\{\shu{K}^{\min}\}_{K^p}$) if, for any $\Phi\in\ca{CLR}(G,X)$, the integral model of mixed Shimura variety $\{\mathscr{S}_{K_\Phi}\}_{K^p}$ (resp. $\{\mathscr{S}_{K_{\Phi, h}}\}_{K^p}$) is a system of canonical integral models in the sense of Axiom \ref{def: canonical model for mixed Shimura data} (resp. Axiom \ref{def: canonical model for pure Shimura data}).
      The inverse system runs over a cofinal collection of neat open compact subgroups $K^{p}\sbst K^{\prime,p}$ with the cone decomposition the induced one (see \cite[Def. 1.18(2)]{Wu25}).
   \end{definition}
\begin{thm}\label{thm-can-int-mod-cpt-summary}
Let $(G, X)$ be any abelian-type Shimura datum, let $p > 0$ be any prime, and let $\GG$ be any quasi-parahoric model with $K_p=\ca{G}(\bb{Z}_p)$. Then $\lrbracket{\shu{K_pK^p}(G, X)}_{K^p \subset G(\A^p)}$ has a \emph{canonical integral model} $\lrbracket{\Shum{K_pK^p}}_{K^p\sbst G(\Ap)}$ in the sense of Axiom \ref{def: canonical model for pure Shimura data} that satisfies Assumption \ref{ass-well-position}. Moreover, $\lrbracket{\Shum{K_pK^p}}_{K^p \subset G(\A^p)}$ has canonical integral models of toroidal compactifications and minimal compactifications for any given neat open compact $K^{\prime,p}$.
\end{thm}
\begin{proof}
In Theorem \ref{thm-ext-cim-ab}, we will check Assumption \ref{ass-well-position_mixed}(2) and the second assertion using the construction in \cite{Wu25} (specializing to Case ($\mrm{STB}_1$) and its normalization models with corresponding quasi-parahoric levels; see \emph{loc. cit.} Sec. 4.2). Assumption \ref{ass-well-position_mixed}(1) (i.e., Axiom \ref{axiom-good-compactification}) is Theorem \ref{thm-abelian-type-axiom}. The second statement of Assumption \ref{ass-well-position} is \cite[Prop. 4.59]{Wu25}; as one will see later, this, in fact, requires Proposition \ref{prop-always-comp-in-derived-part}. 
\end{proof}
\subsubsection{}
We mention some important consequences. The following proposition means that the shtukas on the strata \emph{automatically} glue together.
\begin{prop}\label{prop-degeneration-int}
Under Assumption \ref{ass-well-position_mixed}, there is a unique morphism
$$(\mathscr{S}_K^{\Sigma})^{\log \Diamond}\to \Sht_{\ca{G}^c,\mu^c,\delta=1}$$
extending $\mathscr{S}_K^{\Diamond/}\to \Sht_{\ca{G}^c,\mu^c}$. Moreover, there is a commutative diagram
\begin{equation}\label{eq-diag-degeneration-int}
    \begin{tikzcd}
    {(\mathscr{S}_K^{\Sigma})^{\log \Diamond}}\arrow[d]&W^{\log \Diamond/}\arrow[l]\arrow[r]&(\Delta^\circ_{\Phi,K}\bss\mathscr{S}_{K_\Phi}(\sigma))^{\log\Diamond/}\arrow[d]\\
    \Sht_{\ca{G}^c,\mu^c,\delta=1}&&\Sht_{\ca{P}_\Phi^*,\mu_\Phi^*,\delta=1}.\arrow[ll,"\Int(g_\Phi^{-1})"]
    \end{tikzcd}
\end{equation}
\end{prop}
By Lemma \ref{lem-proper-equivalent}, there is no difference of using $\log \Dia$ and $\log \Diamond/$ on the top-left corner of (\ref{eq-diag-degeneration-int}) since $\mathscr{S}_K^\Sigma$ is proper.
\begin{proof}
Let $\mathscr{S}_K$ be an integral model satisfying the Assumption \ref{ass-well-position_mixed}. By Corollary \ref{cor: shtuka comparison, over W, generic fiber} and Theorem \ref{thm-ext-shu-gen}, there is a unique $\ca{P}^*_\Phi$-shtuka $\mathscr{P}_\Phi^*$ on $\Delta^\circ_{\Phi,K}\bss \mathscr{S}_{K_\Phi}^{\Diamond/}$ such that the commutative diagram (\ref{eq: shtuka comparison, over W}) extends integrally. Then the first sentence follows from applying Corollary \ref{cor-ext-toric-emb} and Lemma \ref{lem-gluing-shtukas} in order. The diagram (\ref{eq-diag-degeneration-int}) follows from the proof of Lemma \ref{lem-gluing-shtukas}.
\end{proof}
Combing Proposition \ref{prop-degeneration-int} and Corollary \ref{cor-special-fiber-nonlog}, we see that
\begin{cor}\label{cor-degeneration-sp}
The log shtuka in Proposition \ref{prop-degeneration-int} induces a non-log shtuka on the special fiber 
$$\mathscr{S}_{K,\ca{O}_{E}/\mathfrak{p}_v}^{\Sigma,\Diamond}=\mathscr{S}_{K,\ca{O}_E/\mathfrak{p}_v}^{\Sigma,\dia}\to \Sht_{\ca{G}^c,\mu^c,\delta=1}$$
with a commutative diagram
\begin{equation}\label{eq-diag-degeneration-sp}
    \begin{tikzcd}
    {(\mathscr{S}_{K,\ca{O}_{E}/\mathfrak{p}_v}^{\Sigma})^{\perf}}\arrow[d]&W^{\perf}_{\ca{O}_{E}/\mathfrak{p}_v}\arrow[l]\arrow[r]&(\Delta^\circ_{\Phi,K}\bss\mathscr{S}_{K_\Phi}(\sigma))^{\perf}_{\ca{O}_{E}/\mathfrak{p}_v}\arrow[d]\\
    \Sht_{\ca{G}^c,\mu^c,\delta=1}^W&&\Sht_{\ca{P}_\Phi^*,\mu_\Phi^*,\delta=1}^W.\arrow[ll,"\Int(g_\Phi^{-1})"]
    \end{tikzcd}
\end{equation}
\end{cor}

\begin{cor}\label{cor-big-diamond}
Under Assumption \ref{ass-well-position_mixed}, for any $\Upsilon=[(\Phi,\sigma)]\in\cusp_K(G,X,\Sigma)$, there is a map 
$$(\ca{Z}_{\Upsilon,K})^{\log ?}\to \Sht_{\ca{G}^c,\mu^c,\delta=1},$$
where $?=\Diamond, \diamond$ or $\Diamond/$ and the log structure on $\ca{Z}_{\Upsilon,K}$ is defined by pulling back the log structure on $\mathscr{S}_K^\Sigma$.
In particular, there is a morphism
$$\mathscr{S}_K^\Diamond\to \Sht_{\ca{G}^c,\mu^c,\delta=1}$$
extending $\mathscr{S}_K^{\Diamond/}\to \Sht_{\ca{G}^c,\mu^c,\delta=1}.$
\end{cor}
\begin{proof}
Since $\mathscr{S}^\Sigma_K$ is proper, $(\mathscr{S}_K^{\Sigma})^{\log \Diamond}=(\mathscr{S}_K^\Sigma)^{\log \Diamond/}=(\mathscr{S}_K^\Sigma)^{\log\diamond}$ (see Lemma \ref{lem-proper-equivalent}). For a map between fs log schemes $(X,\ca{M}_X)\to(Y,\ca{M}_Y)$ that are separated and of finite type, by definition, there is a map between log diamonds $(X,\ca{M}_X)^{\log ?}\to(Y,\ca{M}_Y)^{\log ?}$. For the last sentence, note that the pullback log structure to the interior is the trivial log structure. Then it follows from Corollary \ref{cor-trivial-log}.
\end{proof}
\subsubsection{}

    \begin{lem}\label{lem: left square}
         Under Assumption \ref{ass-well-position}, for any $\Phi\in\ca{CLR}(G,X)$, the $\ovl{\PP}^*_{\Phi}$-shtuka $(\ovl{\PPs}_{\Phi, E}, \phi_{\ovl{\PPs}_{\Phi, E}})$ over $\shu{\ovl{K}_{\Phi}}^{\Dia}$ with one leg bounded by $\ovl{\mu}^*_{\Phi}$ associated with the push-out of the de Rham pro-\'etale $\underline{\ovl{\PP}^{c}_{\Phi}(\Z_p)}$-torsor $\ovl{\PPp}_{K_{\Phi}} \to \shu{\ovl{K}_{\Phi}}$ via $\ovl{\PP}^c_{\Phi} \to \ovl{\PP}^*_{\Phi}$ extends to a $\ovl{\PP}^*_{\Phi}$-shtuka $(\ovl{\PPs}_{\Phi}, \phi_{\ovl{\PPs}_{\Phi}})$ over $\Shum{\ovl{K}_{\Phi}}^{\Dia/}$. In other words, we have a morphism $(\Delta_{\Phi, K}^{\circ}\backslash\Shum{\ovl{K}_{\Phi}})^{\Dia/} \to \Sht_{\ovl{\PP}_{\Phi}^*, \bar{\mu}_{\Phi}^*}$ extending the one on generic fiber. Similarly, we have a morphism $(\Delta_{\Phi, K}^{\circ}\backslash\Shum{K_{\Phi, h}})^{\Dia/} \to \Sht_{\GG_{\Phi, h}^*, \mu_{\Phi, h}^*}$ extending the one on generic fiber.
    \end{lem}
    \begin{proof}
        Zariski locally over $\Shum{\ovl{K}_{\Phi}} =: C(\Phi)$, $\Shum{K_{\Phi}} \to \Shum{\ovl{K}_{\Phi}}$ is the projection $\mathbf{E}(\Phi) \times_{\Spec \Z} C(\Phi) \to C(\Phi)$. It admits the zero section $C(\Phi) \to \mathbf{E}(\Phi) \times C(\Phi)$. We have the following commutative diagrams:
\[
\begin{tikzcd}
	{(\mathbf{E}(\Phi)_{\eta} \times C(\Phi)_{\eta})^{\Dia}} & {\Sht_{\PP_{\Phi}^*, \mu_{\Phi}^*}} & {(\mathbf{E}(\Phi) \times C(\Phi))^{\Dia/}} & {\Sht_{\PP_{\Phi}^*, \mu_{\Phi}^*}} \\
	{C(\Phi)_{\eta}^{\Dia}} & {\Sht_{\ovl{\PP}_{\Phi}^*, \bar{\mu}_{\Phi}^*},} & {C(\Phi)^{\Dia/}} & {\Sht_{\ovl{\PP}_{\Phi}^*, \bar{\mu}_{\Phi}^*}.}
	\arrow[from=1-1, to=1-2]
	\arrow[from=1-1, to=2-1]
	\arrow[from=1-2, to=2-2]
	\arrow[from=1-3, to=1-4]
	\arrow[from=1-3, to=2-3]
	\arrow[from=1-4, to=2-4]
	\arrow[shift left=3, from=2-1, to=1-1]
	\arrow[from=2-1, to=2-2]
	\arrow[shift left=3, from=2-3, to=1-3]
	\arrow[dashed, from=2-3, to=2-4]
\end{tikzcd}
\]
       The composition  $C(\Phi)^{\Dia/} \to \Sht_{\ovl{\PP}_{\Phi}^*, \bar{\mu}_{\Phi}^*}$ exists Zariski locally on $C(\Phi)$. Since, over generic fiber, we have a globally defined morphism $C(\Phi)_{\eta}^{\Dia} \to \Sht_{\ovl{\PP}_{\Phi}^*, \bar{\mu}_{\Phi}^*}$, and the commutativity of the left square implies that the morphism is independent of the choice of (Zariski local) section $C(\Phi) \to \mathbf{E}(\Phi) \times C(\Phi)$, we obtain a globally defined morphism $C(\Phi)^{\Dia/} \to \Sht_{\ovl{\PP}_{\Phi}^*, \bar{\mu}_{\Phi}^*}$ using \cite[Cor. 2.7.10]{PR24}. This further descends to $(\Delta_{\Phi, K}^{\circ}\backslash\Shum{\ovl{K}_{\Phi}})^{\Dia/}$ using the descent data coming from the generic fiber. We can apply this argument to $\Shum{K_{\Phi, h}}$ since $\Shum{\ovl{K}_{\Phi}} \to \Shum{K_{\Phi, h}}$ is an abelian scheme torsor by assumption.
    \end{proof}

     In particular, we have
 \begin{prop}\label{prop: boundary shtukas extend over integral base}
     Under Assumption \ref{ass-well-position}, the commutative diagram (\ref{eq: reduction, generic fiber, general, star}) in Lemma \ref{lem: reduction, generic fiber, general, star} extends over integral base. That is to say, for any $[\Phi]$, we have the following commutative diagram
\[
\begin{tikzcd}
	{\Delta_{\Phi, K}^{\circ}\backslash\Shum{K_{\Phi}}(P_{\Phi}, D_{\Phi})^{\Dia/}} & {\Delta_{\Phi, K}^{\circ}\backslash\Shum{\ovl{K}_{\Phi}}(\ovl{P}_{\Phi}, \ovl{D}_{\Phi})^{\Dia/}} & {\Delta_{\Phi, K}^{\circ}\backslash\Shum{K_{\Phi, h}}(G_{\Phi, h}, D_{\Phi, h})^{\Dia/}} \\
	{\Sht_{\PP_{\Phi}^*, \mu_{\Phi}^*, \delta = 1}} & {\Sht_{\ovl{\PP}_{\Phi}^*, \bar{\mu}_{\Phi}^*, \delta = 1}} & {\Sht_{\GG_{\Phi, h}^*, \mu_{\Phi, h}^*, \delta = 1}.}
	\arrow[from=1-1, to=1-2]
	\arrow[from=1-1, to=2-1]
	\arrow[from=1-2, to=1-3]
	\arrow[from=1-2, to=2-2]
	\arrow[from=1-3, to=2-3]
	\arrow[from=2-1, to=2-2]
	\arrow[from=2-2, to=2-3]
\end{tikzcd}
\]
 \end{prop}

    \begin{prop}\label{prop: shtuka comparison, over W, integral base}
Under Assumption \ref{ass-well-position}, the commutative diagram (\ref{eq: shtuka comparison, over W}) in Corollary \ref{cor: shtuka comparison, over W, generic fiber} extends over the integral base. That is to say, for all such coverings $W^0$ (where $\mathfrak{W}$ are open affine coverings of $\mathfrak{X}^{\circ}_{[\Phi, \sigma], K}$, here for each $\Phi$, we let $\sigma$ run over all $\Sigma(\Phi)^+$), we have following commutative diagram
        \begin{equation}\label{eq: main diagram}
\begin{tikzcd}
	{\Shum{K}(G, X)^{\dia}} & {W^{0, \dia}} & {\Delta_{\Phi, K}^{\circ}\backslash\Shum{K_{\Phi}}(P_{\Phi}, D_{\Phi})^{\dia}} & {\Delta_{\Phi, K}^{\circ}\backslash\Shum{K_{\Phi, h}}(G_{\Phi, h}, D_{\Phi, h})^{\dia}} \\
	{\Sht_{\GG^c, \mu^c, \delta = 1}} & {} & {\Sht_{\PP_{\Phi}^*, \mu_{\Phi}^*, \delta = 1}} & {\Sht_{\GG_{\Phi, h}^*, \mu_{\Phi, h}^*, \delta = 1}.}
	\arrow[from=1-1, to=2-1]
	\arrow[from=1-2, to=1-1]
	\arrow[from=1-2, to=1-3]
	\arrow[from=1-3, to=1-4]
	\arrow[from=1-3, to=2-3]
	\arrow[from=1-4, to=2-4]
	\arrow["{\Int(g_{\Phi}^{-1})}"', from=2-3, to=2-1]
	\arrow[from=2-3, to=2-4]
\end{tikzcd}
        \end{equation}

    Taking the special fibers (apply the reduction functor), we have
    \begin{equation}\label{eq: special fiber of main diagram}
\begin{tikzcd}
	{\Shum{K}(G, X)^{\perf}_{\bar{s}}} & {W^{0, \perf}_{\bar{s}}} & {\Delta_{\Phi, K}^{\circ}\backslash\Shum{K_{\Phi}}(P_{\Phi}, D_{\Phi})^{\perf}_{\bar{s}}} & {\Delta_{\Phi, K}^{\circ}\backslash\Shum{K_{\Phi, h}}(G_{\Phi, h}, D_{\Phi, h})^{\perf}_{\bar{s}}} \\
	{\Sht_{\GG^c, \mu^c, \delta = 1}^W} & {} & {\Sht_{\PP_{\Phi}^*, \mu_{\Phi}^*, \delta = 1}^W} & {\Sht_{\GG_{\Phi, h}^*, \mu_{\Phi, h}^*, \delta = 1}^W.}
	\arrow[from=1-1, to=2-1]
	\arrow[from=1-2, to=1-1]
	\arrow[from=1-2, to=1-3]
	\arrow[from=1-3, to=1-4]
	\arrow[from=1-3, to=2-3]
	\arrow[from=1-4, to=2-4]
	\arrow["{\Int(g_{\Phi}^{-1})}"', from=2-3, to=2-1]
	\arrow[from=2-3, to=2-4]
\end{tikzcd}
    \end{equation}
    \end{prop}
    Note that $W^0$ is not of finite type over $\Z_p$; nevertheless, since $W^0$ is affine, noetherian and flat over $\Z_p$, one can still apply \cite[Cor. 2.7.10]{PR24} (see the first paragraph of the proof of Theorem \ref{thm-ext-shu-gen}).

\subsubsection{}
The compactifications satisfying Definition \ref{def-PR-int-mod} are fully functorial.
\begin{prop}\label{prop-full-functoriality}
Let $f:(G_1,X_1,\ca{G}_1, K_1^p,\Sigma_1)\to (G_2,X_2,\ca{G}_2,K_2^p,\Sigma_2)$ be a morphism. Here, the meaning of a morphism is: (1) $f:(G_1,X_1)\to (G_2,X_2)$ is a morphism between Shimura data; (2) $\ca{G}_1\to\ca{G}_2$ is a morphism between quasi-parahoric group schemes whose generic fiber is $f:G_1\to G_2$; (3) $K_1^p\to K_2^p$ is a morphism between neat open compact subgroups compatible with $f: G_1(\Ap)\to G_2(\Ap)$; (4) $\Sigma_1$ and $\Sigma_2$ are smooth projective and are compatible in the sense of \cite[Def. 1.18(3)]{Wu25}.\par
Let $\mathscr{S}_{K_i}^{\Sigma_i}(G_i,X_i)$ be a canonical integral model of $\sh_{K_i}^{\Sigma_i}(G_i,X_i)$ for $i=1,2$. Then there is a unique morphism $f_{\Sigma_1,\Sigma_2}:\mathscr{S}_{K_1}^{\Sigma_1}\to \mathscr{S}_{K_2}^{\Sigma_2}$ extending the one $f:\sh_{K_1}\to \sh_{K_2}$ over the generic fiber by functoriality between Shimura data, and makes the following diagram commutes:
\[
\begin{tikzcd}
	{(\mathscr{S}_{K_1}^{\Sigma_1})^{\log\Dia}} & {\Sht_{\GG^c_1, \mu^c_1, \delta = 1}} \\
	{(\mathscr{S}_{K_2}^{\Sigma_2})^{\log\Dia}} & {\Sht_{\GG^c_2, \mu^c_2, \delta = 1}.}
	\arrow[from=1-1, to=1-2]
	\arrow[from=1-1, to=2-1]
	\arrow[from=1-2, to=2-2]
	\arrow[from=2-1, to=2-2]
\end{tikzcd}
\]
\end{prop}
\begin{proof}
By Proposition \ref{prop: morphisms extend to integral models}, for any cusp label representative $\Phi_1$ mapping to $\Phi_2$, there is a morphism $f_{\Phi_1}:\mathscr{S}_{K_{\Phi_1}}\to \mathscr{S}_{K_{\Phi_2}}$ extending the one on the generic fiber. Note that this morphism is automatically a morphism between torus torsors equivariant under $(\mbf{E}_{K_{\Phi_1}}\to \mbf{E}_{K_{\Phi_2}})$ by separatedness and normality. Hence, we obtain a morphism of integral models $f_{\Phi_1}(\sigma_1):\mathscr{S}_{K_{\Phi_1}}(\sigma_1)\to \mathscr{S}_{K_{\Phi_2}}(\sigma_2)$ for $\sigma_1$ mapping to $\sigma_2$; we also have a morphism between closed strata $f_{\Phi_1,\sigma_1}$. Axiom \ref{axiom-4} requires that $\Delta_{\Phi,K}^\circ$ acts on $\mathscr{S}_{K_\Phi}$. By separatedness and normality again, there are morphisms $f^\circ_{\Phi_1}:\Delta^\circ_{\Phi_1,K_1}\bss\mathscr{S}_{K_{\Phi_1}}\to \Delta^\circ_{\Phi_2,K_2}\bss \mathscr{S}_{K_{\Phi_2}}$, $f^\circ_{\Phi_1}(\sigma_1):\Delta^\circ_{\Phi_1,K_1}\bss\mathscr{S}_{K_{\Phi_1}}(\sigma_1)\to \Delta^\circ_{\Phi_2,K_2}\bss \mathscr{S}_{K_{\Phi_2}}(\sigma_2)$ and also a morphism between their corresponding closed strata. We have the desired morphism between compactifications by \cite[Lem. A.3.4]{Mad19}. The commutativity of the diagram is guaranteed by Theorem \ref{thm-ext-shu-gen}, and by the commutativity of the diagram on the generic fiber.
\end{proof}
In particular,
\begin{cor}\label{cor-fully-functorial-uniqueness}
Fix $(G,X,\Sigma)$. The inverse system of canonical integral models $\{\mathscr{S}_K^\Sigma\}_{K^p}$ is unique up to a unique isomorphism.
\end{cor}
Upon restricting the maps to the open dense strata, we obtain that
\begin{cor}
    Fix $(G, X)$. Let $\GG_1 \to \GG_2$ be a morphism of quasi-parahoric group schemes, then there exists a \emph{proper} morphism $\Shum{K_1} \to \Shum{K_2}$ extending $\shu{K_1} \to \shu{K_2}$.
\end{cor}

Similarly, one has canonicity and functoriality of canonical integral models of minimal compactifications.
\begin{prop}\label{prop-full-functoriality-min}
Let $f:(G_1,X_1,\ca{G}_1, K_1^p)\to (G_2,X_2,\ca{G}_2,K_2^p)$ be a morphism as in Proposition \ref{prop-full-functoriality}. 
Let $\mathscr{S}_{K_i}^{\mmin}(G_i,X_i)$ be a canonical integral model of $\sh_{K_i}^\mmin(G_i,X_i)$ for $i=1,2$. Then there is a unique morphism $f^\mmin:\mathscr{S}_{K_1}^{\mmin}\to \mathscr{S}_{K_2}^{\mmin}$ extending $f:\sh_{K_1}\to \sh_{K_2}$ over the generic fiber between Shimura varieties induced by functoriality between Shimura data.  \par
In particular, the inverse system of canonical integral models of minimal compactifications $\{\mathscr{S}_K^\mmin\}_{K^p}$ is unique up to a unique isomorphism.
\end{prop}
\begin{proof}
The argument as in \cite[Lem. A.3.4]{Mad19} still works. See \cite[Lem. 5.5]{Wu25} and \cite[4.1.8-4.1.9]{Mao25b}.
\end{proof}

\section{Shtukas on integral models of boundary mixed Shimura varieties of abelian type}\label{sec-can-ext-ab}

In this section, we construct shtukas on integral models of boundary mixed Shimura varieties of abelian type, and verify the axioms \ref{def: canonical model for mixed Shimura data}.

\subsection{Hodge-type case}

\subsubsection{Set up}

Let $(G, X)$ be a Hodge-type Shimura datum, and let $\GG$ be a stabilizer quasi-parahoric model of $G$, with $K_p = \GG(\Z_p)$. 
We fix an adjusted Hodge embedding $(G, X, K_p) \hookrightarrow (G^{\dd}, X^{\dd}, K_p^{\dd})$ in the sense of \cite[Def. 3.16, Rmk. 3.17]{Mao25b}, where $G^{\dd} = \GSp(V, \psi)$, $\GSP = \GSp(V_{\Z_p})$ is a hyperspecial model of $G^{\dd}$, $K_p^{\dd} = \GSP(\Z_p)$ is the stabilizer group of the self-dual lattice $V_{\Z_p}$, and $K_p = G(\rQ_p) \cap K_p^{\dd}$. Explicitly, we work with the construction in \cite{KPZ24}, combine the chain of lattices into a single lattice and then apply the Zarhin's trick to get a self-dual lattice in the Siegel side. We fix such an adjusted Hodge embedding instead of a random one as in \cite{PR24} in order to have a good theory of compactifications as in \cite{Mad19} and to ensure good compatibility of levels at the boundary; see \cite[Prop. 1.2]{Mao25b}. 
The induced morphism $\GG \to \GSP$ factors through a smoothening map $\GG \to \ovl{\GG}$, and $\ovl{\GG}\hookrightarrow \GSP$ is a closed embedding. Let $(s_{\alpha}) \in V_{\Z_p}^{\otimes}$ be the Hodge tensors that cut out the closed subscheme $\ovl{\GG}$. Note that $s_{\alpha}$ can be defined in $V_{\Z_{(p)}}^{\otimes}$.

\begin{rk}\label{rk: Hodge-Type case, no c}
In the Hodge-type case, the center $Z(G)^{\circ}$ is an almost product of a split $\rQ$-torus and a compact-type $\rQ$-torus; thus $Z(G)_{ac}$ is trivial, $G = G^c$, $P_{\Phi} = P^c_{\Phi} = P^*_{\Phi}$, and $\Delta_{\Phi, K}^{\circ}$ is trivial. See \cite[Lem. 2.1.20]{Mad19}.
\end{rk}

On the boundary, since the level group $K_{\Phi, G} = gKg^{-1}$ is twisted by $g = g_{\Phi} \in G(\A)$, we need to twist the Hodge tensors. Let $(s_{\alpha, \Phi}) \in V_{\Z_{(p)}}^{g, \otimes}$ (where $V^g_{\Z_{(p)}} = gV_{\Z_p} \cap V_{\rQ}$) be the collection of Hodge tensors that cut out the closed subgroup scheme $\ovl{\GG}_{\Phi} \to \GSP_{\Phi^{\dd}}$; here $\ovl{\GG}_{\Phi}$ and $\GSP_{\Phi^{\dd}}$ are the left $g$-conjugates of $\ovl{\GG}$ and $\GSP$, respectively. Indeed, we take $g_p$-conjugation, where $g_p$ is the $p$-component of $g$; we do not distinguish $g$ and $g_p$ when the difference is clear from the context. Let $\GG_{\Phi} \to \ovl{\GG}_{\Phi}$ be the dilated morphism; $\GG_{\Phi}$ is a smooth affine model of $G$. 

 Over $\Z_p$, $(s_{\alpha})$ and $(s_{\alpha, \Phi})$ are related as follows:
 \begin{equation*}
\begin{tikzcd}
	{V^{\otimes n}_{\Z_p}} & {V^{g, \otimes n}_{\Z_p}} \\
	{V^{\otimes m}_{\Z_p}} & {V^{g, \otimes m}_{\Z_p}.}
	\arrow["g", from=1-1, to=1-2]
	\arrow["{s_{\alpha}}", from=1-1, to=2-1]
	\arrow["{s_{\alpha, \Phi}}", from=1-2, to=2-2]
	\arrow["g", from=2-1, to=2-2]
\end{tikzcd}
 \end{equation*}

\subsubsection{Reduction}
 Recall that over $\shu{K}(G, X)$, the family of Hodge tensors $(s_{\alpha}) \in V^{\otimes}_{\Z_p}$ defines a family of \'etale tensors $(t_{\alpha, \et})$ on the Tate module $T_p(A)$, where $A$ is the universal abelian scheme over $\shu{K}(G, X)$, pulled back from $\shu{K^{\dd}}(G^{\dd}, X^{\dd})$. We denote by $\ls_{\rho, V_{\Z_p}}$ the local system given by $T_p(A)$, where $\rho: G \to \GSp(V) \to \GL(V)$. Let $\lp_K$ be the pro-\'etale torsor under $\underline{K_p}$ defined by
 \begin{equation}\label{eq: def of lp_K}
     \Isom_{(t_{\alpha, \et}), (s_{\alpha})}(\ls_{\rho, V_{\Z_p}}, V_{\Z_p, \shu{K}(G, X)}),
 \end{equation}
 this coincides with the pro-\'etale torsor (\ref{eq-proetale torsor-mixed}).

 One can repeat the above process with $A$ replaced by the universal $1$-motive $\QQ$ over $\shu{K_{\Phi}}:=\shu{K_{\Phi}}(P_{\Phi}, D_{\Phi})$, pulled back from $\shu{K_{\Phi^{\dd}}}:= \shu{K_{\Phi^{\dd}}}(P_{\Phi^{\dd}}, D_{\Phi^{\dd}})$. Here $\Phi^{\dd} \in \ca{CLR}(G^{\dd}, X^{\dd})$ is the cusp-label representative induced by $\Phi$.

 By the work of Brylinski, all Hodge classes on $1$-motives over $\CC$ are absolute Hodge cycles; thus $(s_{\alpha, \Phi})$ induces $(t_{\alpha, \Phi, \et}) \in T_p(\QQ)^{\otimes}$ over $\shu{K_{\Phi}}$. We denote by $\ls_{\rho_{\Phi}, K_{\Phi}}$ the local system given by the Tate module $T_p(\QQ)$. Consider the following torsor under $\underline{\GG_{\Phi}(\Z_p)}$ over $\shu{K_{\Phi}}$:
 \begin{equation}\label{eq: def of lp_KPhi}
    \lp_{K_{\Phi, G}}:= \Isom_{(t_{\alpha, \Phi, \et}), (s_{\alpha, \Phi})}(\ls_{\rho_{\Phi}, K_{\Phi}}, V^g_{\Z_p, \shu{K_{\Phi}}}).
 \end{equation}
Consider the following torsor under $\underline{K_{\Phi, Q, p}} = \underline{\QQ_{\Phi}(\Z_p)}$ over $\shu{K_{\Phi}}$:
 \begin{equation*}
    \lp_{K_{\Phi, W^{\bullet}}}:= \Isom_{(t_{\alpha, \Phi, \et}), (s_{\alpha, \Phi})}(W_{\bullet}\ls_{\rho_{\Phi}, K_{\Phi}}, W_{\bullet}V^g_{\Z_p, \shu{K_{\Phi}}}),
 \end{equation*}
 where the filtration $W_{\bullet}\ls_{\rho_{\Phi}, K_{\Phi}}$ comes from the weight filtration $W_{\bullet}T_p(\mQ)$ on $T_p(\mQ)$. Since the weight filtration is constant, both the Hodge tensors and the weight filtration are trivialized when $T_p(\mQ)$ is trivialized. Therefore, we have a canonical embedding $\lp_{K_{\Phi, W^{\bullet}}} \to \lp_{K_{\Phi, G}}$, which induces
 \begin{equation}\label{eq: G reduce to Q, generic}
     \underline{\GG_{\Phi}(\Z_p)} \times^{\underline{\QQ_{\Phi}(\Z_p)}} \lp_{K_{\Phi, W_{\bullet}}} = \lp_{K_{\Phi, G}}.
 \end{equation}

  On the other hand, there is a natural pro-\'etale torsor $\lp_{K_{\Phi}}$ under $K_{\Phi}$ over $\shu{K_{\Phi}}$, as defined in Definition \ref{def-proetale-torsor-aut}.

  \begin{lem}\label{lem: generic fiber, G reduces to P}
     $T_p(\QQ)$ is trivialized over $\varprojlim_{K_{\Phi}' \subset K_{\Phi}} \shu{K_{\Phi}'K_{\Phi}^p}$.
 \end{lem}
 \begin{proof}
     It suffices to show that for each $n$, $T_p(\QQ)/p^n$ is trivialized over $\shu{K_{\Phi, p}^{(n)}K_{\Phi}^p}$. It also suffices to work with the Siegel case, where the result follows from the moduli interpretation; see \cite[\S 2.2.14]{Mad19}.
 \end{proof}
 \begin{cor}\label{cor: generic fiber, G reduces to P}
     $\lp_{K_{\Phi, G}}$ reduces to $\lp_{K_{\Phi}}$.
 \end{cor}
 \begin{proof}
     Along the standard representation $\PP_{\Phi} \to \GL(V^g_{\Z_p})$, the push-out torsor of $\lp_{K_{\Phi}}$ is the frame bundle associated with $T_p(\QQ)$. By Lemma \ref{lem: generic fiber, G reduces to P}, we have a morphism
     \begin{equation*}
         \lp_{K_{\Phi}} \to \Isom(\ls_{\rho_{\Phi}, K_{\Phi}}, V^g_{\Z_p, \shu{K_{\Phi}}}),
     \end{equation*}
     that factors through $\lp_{K_{\Phi, G}}$ and thus factors through $\lp_{K_{\Phi, W_{\bullet}}}$.
 \end{proof}

 \subsubsection{Hodge-Tate period map}

 The de Rham pro-\'etale torsor $\lp_K$ under $\underline{K_p}$ over $\shu{K}$ produces a Hodge-Tate period map (\ref{eq: HT, G, G_h}), which is the sheaf-theoretic version of the Hodge-Tate period map
 \begin{equation*}
    \shu{K^p} \sim \varprojlim_{K_p}\shu{K_pK^p} \to \fl_{G, \mu^{-1}}
 \end{equation*}
 from the perfectoid space $\shu{K^p}$ to the flag variety $\fl_{G, \mu^{-1}}$. Here $(\shu{K^p})^{\Dia} = \varprojlim_{K_p}\shu{K_pK^p}^{\Dia}$. Recall that given a geometric point $x \in \shu{K^p}(C, C^+)$, we have an abelian variety $A_x$ over $C$ together with a trivialization $(T_p(A_x), (t_{\alpha, \et, x})) \cong (V_{\Z_p}, (s_{\alpha}))$, then $\HT_{K}(x) \in \fl_{G, \mu^{-1}}(C, C^+)$ corresponds to the Hodge-Tate filtration:
 \begin{equation*}
     0 \to (\Lie(A_x))(1) \to T_p(A_x) \otimes_{\Z_p} C  \to (\Lie A_x^{\vee})^{\vee} \to 0, 
 \end{equation*}
 see \cite[Thm. 2.1.3]{caraiani2017generic}.

 Similarly, let $G_{\Phi}$ be the $\rQ_p$-group $gG_{\rQ_p}g^{-1}$. One has a perfectoid space together with a Hodge-Tate period map
  \begin{equation*}
      \HT_{K_{\Phi, G}}: \shu{K_{\Phi}^p} \to \fl_{G_{\Phi}, \mu_{\Phi}^{-1}},\quad \shu{K_{\Phi}^p} \sim \varprojlim_{K_{\Phi, p}} \shu{K_{\Phi, p}K_{\Phi}^p},\quad  \shu{K_{\Phi}^p}^{\Dia} = \varprojlim_{K_{\Phi, p}} \shu{K_{\Phi, p}K_{\Phi}^p}^{\Dia},
  \end{equation*}
 defined as follows: given $x \in \shu{K_{\Phi}}(C, C^+)$, by Corollary \ref{cor: generic fiber, G reduces to P}, we have a $1$-motive $\QQ_x$ over $C$ together with a trivialization $(T_p(\QQ_x), (t_{\alpha, \Phi, \et, x})) \cong (V_{\Z_p}^g, (s_{\alpha, \Phi}))$. Then $\HT_{K_{\Phi, G}}(x) \in \fl_{G_{\Phi}, \mu_{\Phi}^{-1}}(C, C^+)$ corresponds to the Hodge-Tate filtration
 \begin{equation*}
     0 \to (\Lie(\QQ_x))(1) \to T_p(\QQ_x) \otimes_{\Z_p} C  \to (\Lie \QQ_x^{\vee})^{\vee} \to 0. 
 \end{equation*}

 By Corollary \ref{cor: generic fiber, G reduces to P}, $\HT_{K_{\Phi, G}}$ factors through $\HT_{K_{\Phi}}: \shu{K_{\Phi}^p} \to \fl_{P_{\Phi}, \mu_{\Phi}^{-1}}$ and is the sheaf-theoretic version of the Hodge-Tate period map (\ref{eq: HT, G, G_h}).

 \begin{cor}\label{cor: generic fiber, shtuka, G reduces to P}
     Over $\shu{K_{\Phi}}^{\Dia}$, the $\GG_{\Phi}$-shtuka associated with $(\lp_{K_{\Phi, G}}, \HT_{K_{\Phi, G}})$ with a leg bounded by $\mu_{\Phi}$ reduces to the $\PP_{\Phi}$-shtuka associated with $(\lp_{K_{\Phi}}, \HT_{K_{\Phi}})$ with a leg bounded by $\mu_{\Phi}$.
 \end{cor}

\subsubsection{Comparison}

 Let $\shu{K}(G, X)^{\Dia} \to \Sht_{\GG, \mu}$ (resp. $\shu{K_{\Phi}}(P_{\Phi}, D_{\Phi})^{\Dia} \to \Sht_{\PP_{\Phi}, \mu_{\Phi}}$) be the shtuka associated with $\lp_K$ (resp. $\lp_{K_{\Phi}}$). We use the notation from subsection \ref{subsec-can-ext-generic-fiber}.

 \begin{lem}[{See Corollary \ref{cor: shtuka comparison, over W, generic fiber}}]
     We have a commutative diagram:
     \[ 
\begin{tikzcd}
	{\shu{K}(G, X)^{\Dia}} & {W^{0, \Dia}} & {\shu{K_{\Phi}}(P_{\Phi}, D_{\Phi})^{\Dia}} \\
	{\Sht_{\GG, \mu}} & {} & {\Sht_{\PP_{\Phi}, \mu_{\Phi}}.}
	\arrow[from=1-1, to=2-1]
	\arrow[from=1-2, to=1-1]
	\arrow[from=1-2, to=1-3]
	\arrow[from=1-3, to=2-3]
	\arrow["{\Int(g_{\Phi}^{-1})}"', from=2-3, to=2-1]
\end{tikzcd}
     \]
 \end{lem}
 \begin{rk}
     This is a special case of Corollary \ref{cor: shtuka comparison, over W, generic fiber}; see Remark \ref{rk: Hodge-Type case, no c}. We give a more direct and intuitive proof in the Hodge-type case.
 \end{rk}
 \begin{proof}
    By Corollary \ref{cor: generic fiber, shtuka, G reduces to P}, it suffices to work with $\shu{K_{\Phi}}(P_{\Phi}, D_{\Phi})^{\Dia} \to \Sht_{\GG_{\Phi}, \mu_{\Phi}}$. In \cite[Prop. 3.1.6]{Mad19} (using fpqc descent from $\bigsqcup V \to W^0$; see \cite[\S 3.1.5]{Mad19} for the notation $V$), we have comparisons over $W^0$:
      \begin{equation*}
    \mQ[p^{\infty}] \cong A[p^{\infty}],\quad T_p(\mQ) \cong T_p(A),\quad (t_{\alpha, \et}) \cong (t_{\alpha, \Phi, \et}),
 \end{equation*}
    the tensor comparisons are compatible with $(s_{\alpha}) \cong (s_{\alpha, \Phi})$, in the sense that the diagram commutes:
\[ 
\begin{tikzcd}
	{V_{\Z_p} ^{\otimes}} & {V^{g, \otimes}_{\Z_p}} & {T_p(\mQ_{W^0})^{\otimes}} & {T_p(A_{W^0})^{\otimes}} \\
	{V_{\Z_p} ^{\otimes}} & {V^{g, \otimes}_{\Z_p}} & {T_p(\mQ_{W^0})^{\otimes}} & {T_p(A_{W^0})^{\otimes}.}
	\arrow["{\stackrel{g^{\otimes}}{\rightiso}}", from=1-1, to=1-2]
	\arrow["{(s_{\alpha})}"', from=1-1, to=2-1]
	\arrow["\cong", from=1-2, to=1-3]
	\arrow["{(s_{\alpha, \Phi})}"', from=1-2, to=2-2]
	\arrow["\cong", from=1-3, to=1-4]
	\arrow["{(t_{\alpha, \Phi, \et})}"', from=1-3, to=2-3]
	\arrow["{(t_{\alpha, \et})}"', from=1-4, to=2-4]
	\arrow["{\stackrel{g^{\otimes}}{\rightiso}}", from=2-1, to=2-2]
	\arrow["\cong", from=2-2, to=2-3]
	\arrow["\cong", from=2-3, to=2-4]
\end{tikzcd}
\]
This gives the isomorphism $\Int(g_{\Phi}^{-1}): (\lp_K, \HT_K) \cong (\lp_{K_{\Phi}}, \HT_{K_{\Phi}})$.
 \end{proof}

 \subsubsection{$\GG_{\Phi}$-shtukas}

  We move on to integral models.

  \begin{prop}\label{prop: ext of G shtukas, Phi}
      The $\GG_{\Phi}$-shtuka $(\PPs_{\Phi, G, E}, \phi_{\PPs_{\Phi, G, E}})$ over $\shu{K_{\Phi}}^{\Dia}\to\Spd(E)$ with one leg bounded by $\mu_{\Phi}$ extends (uniquely) to a $\GG_{\Phi}$-shtuka $(\PPs_{\Phi, G}, \phi_{\PPs_{\Phi, G}})$ over $\Shum{K_{\Phi}}^{\Dia/} \to \Spd(\OO_E)$ with one leg bounded by $\mu_{\Phi}$.
  \end{prop}
   This proposition follows almost verbatim from \cite[\S 4.6.3]{PR24}. We only need to replace the BKF module of the universal $p$-divisible group $A[p^{\infty}]$ over $\Shum{K}$ by the universal $p$-divisible group $\mQ[p^{\infty}]$ over $\Shum{K_{\Phi}}$. The uniqueness of the extension follows from \cite[Cor. 2.7.10]{PR24}. Let us write down some key steps in the construction in order to fix notation.

  \begin{proof}
      By Tannakian formalism, under 
   \begin{equation*}
       \rho_{\Phi}: \PP_{\Phi} \to \GSP(V_{\Z_p}^g) \to \GL(V_{\Z_p}^g), 
   \end{equation*}
   the $\GG_{\Phi}$-shtuka $(\PPs_{\Phi, G, E}, \phi_{\PPs_{\Phi, G, E}})$ over $\shu{K_{\Phi}}^{\Dia}\to\Spd(E)$ gives a vector shtuka $(\VVs_{\Phi, E}, \phi_{\VVs_{\Phi, E}})$ over $\shu{K_{\Phi}}^{\Dia}$, associated with the de Rham local system $\ls_{\rho_{\Phi}, K_{\Phi}}$ as in (\ref{eq: def of lp_KPhi}). We have a family of tensors $(t_{\alpha, \Phi, E}) \in(\VVs_{\Phi, E}, \phi_{\VVs_{\Phi, E}})^{\otimes}$, which can be viewed as shtukas homomorphisms over $\shu{K_{\Phi}}^{\Dia}$:
   \begin{equation}\label{eq: tensors, Phi, generic}
       t_{\alpha, \Phi, E}: (\oplus_i\VVs_{\Phi, E}^{\otimes m_i}, \phi_{\oplus_i\VVs_{\Phi, E}^{\otimes m_i}}) \to (\oplus_i\VVs_{\Phi, E}^{\otimes n_i}, \phi_{\oplus_i\VVs_{\Phi, E}^{\otimes n_i}}).
   \end{equation}

   Given $S = \Spa(R, R^+) \in \Perf$, a map $S \to \Shum{K_{\Phi}}^{\dia}$ is given by an untilt $S^{\sharp} = \Spa(R^{\sharp}, R^{\sharp+}) \to (\wdh{\Shum{K_{\Phi}}})^{\ad}$. We pull back the universal $p$-divisible group $\QQ[p^{\infty}]$ to $\Spec R^{\sharp+}$, and let $M_{\inf, \Phi}(R^{\sharp+})$ be the BKF-module of $\QQ[p^{\infty}]$ with one leg at $\phi(\xi) = 0$. Denote $(\VVs_{\Phi, S}, \phi_{\VVs_{\Phi, S}})$ by the corresponding minuscule vector shtuka with height $2n$ and dimension $n$ over $S$ with one leg at $S^{\sharp}$, given by the restriction of $(\phi^{-1})^*M_{\inf, \Phi}(R^{\sharp+})$ to $\Spa(W(R^+))\backslash \lrbracket{[\varpi] = 0}$. By \cite[Thm. 2.7.7]{PR24}, the tensors $(t_{\alpha, \Phi, E})$ in (\ref{eq: tensors, Phi, generic}) extend uniquely to tensors $(t_{\alpha, \Phi}) \in (\VVs_{\Phi, S}, \phi_{\VVs_{\Phi, S}})^{\otimes}$. 
   
   By exactly the same arguments as in Steps $1, 2, 3$ in \cite[\S 4.6.3]{PR24}, the $v$-sheaf over $S$ defined by
   \begin{equation}\label{eq: ovlP, G torsor}
       \ovl{\PPs}_{\Phi, G, S}:= \Isom_{(t_{\alpha, \Phi}), (s_{\alpha, \Phi}) \otimes \identity}(\VVs_{\Phi, S}, V^g_{\Z_p} \otimes_{\Z_p} \OO_{S\dottimes\Z_p})
   \end{equation}
    is induced by a $\GG_{\Phi}$-torsor $\PPs_{\Phi, G, S}$ over $S\dottimes\Z_p$ via $\ovl{\PPs}_{\Phi, G, S} = \ovl{\GG}_{\Phi} \times^{\GG_{\Phi}} \PPs_{\Phi, G, S}$. The Frobenius structure $\phi_{\PPs_{\Phi, G, S}}$ comes from the Frobenius structure $\phi_{\VVs_{\Phi, S}}$ of $\VVs_{\Phi, S}$. Then $(\PPs_{\Phi, G, S}, \phi_{\Phi, G, S})$ gives a $\GG_{\Phi}$-shtuka over $S$. It has a leg bounded by $\mu_{\Phi}$ by its nature over the generic fiber: one can reduce it to $(C, \OO_C)$-points and then apply arguments in \cite[\S 3.3.7]{PR24}. By varying $S$, and using descent data from the generic fiber, we finish the proof.
  \end{proof}

  \subsubsection{$\QQ_{\Phi}$-shtukas}

We follow the notation from the proof of Proposition \ref{prop: ext of G shtukas, Phi}. The BKF module $M_{\inf, \Phi}(R^{\sharp+})$ associated with the universal $p$-divisible group $\mQ[p^{\infty}]$ has a canonical filtration coming from the weight filtration on $\mQ[p^{\infty}]$; therefore, the associated vector shtuka $(\VVs_{\Phi, S}, \phi_{\VVs_{\Phi, S}})$ is equipped with a canonical weight filtration $W_{\bullet}(\VVs_{\Phi, S}, \phi_{\VVs_{\Phi, S}})$. Consider the $v$-sheaf over $S$ defined by
   \begin{equation*}
       \ovl{\PPs}_{\Phi, W_{\bullet}, S}:= \Isom_{(t_{\alpha, \Phi}), (s_{\alpha, \Phi}) \otimes \identity}(W_{\bullet}\VVs_{\Phi, S}, W_{\bullet}V^g_{\Z_p} \otimes_{\Z_p} \OO_{S\dottimes\Z_p}).
   \end{equation*}

   Note that $\ovl{\PPs}_{\Phi, W_{\bullet}, S}$ has a natural $\QQ_{\Phi}$-action factoring through $\ovl{\QQ}_{\Phi}$.

   \begin{lem}\label{lem: shtuka, reduces to Q torsor}
       The $\GG_{\Phi}$-torsor $\PPs_{\Phi, G, S}$ reduces to a $\QQ_{\Phi}$-torsor $\PPs_{\Phi, W_{\bullet}, S}$ (i.e. $\PPs_{\Phi, G, S} = \GG_{\Phi} \times^{\QQ_{\Phi}} \PPs_{\Phi, W_{\bullet}, S}$), and $\ovl{\PPs}_{\Phi, W_{\bullet}, S} = \ovl{\QQ}_{\Phi}\times^{\QQ_{\Phi}}\PPs_{\Phi, W_{\bullet}, S}$. Moreover, the $\GG_{\Phi}$-torsor isomorphism $\phi_{\PPs_{\Phi, G, S}}$ reduces to a $\QQ_{\Phi}$-torsor isomorphism
       \begin{equation*}
           \phi_{\PPs_{\Phi, W_{\bullet}, S}}: \Frob_{S}^*(\PPs_{\Phi, W_{\bullet}, S})|_{S\dottimes \Z_p \setminus S^{\sharp}} \rightiso \PPs_{\Phi, W_{\bullet}, S}|_{S\dottimes \Z_p \setminus S^{\sharp}},
       \end{equation*}
   \end{lem}
   \begin{proof}
       This follows from the fact that the weight filtration on $\mQ$ is constant. To be more precise, in the construction \cite[\S 4.6.3]{PR24}, one first deals with geometric points, and then deals with products of geometric points. Let $T \to S$ be a $v$-cover such that $T$ is a product of geometric points. By \cite[Prop. 19.5.3]{SW20}, the category of $\GG_{\Phi}$-torsors (resp. $\QQ_{\Phi}$-torsors) on $S \dottimes \Z_p$ is equivalent to the category of $\GG_{\Phi}$-torsors (resp. $\QQ_{\Phi}$-torsors) on $T\dottimes\Z_p$ with suitable descent data (resp. descent data for $\GG_{\Phi}$-torsors plus descent data for the weight filtration).
       
       Denote $T = \Spa(B, B^+)$, where $B^+ = \prod_j C_j^+$. The $\GG_{\Phi}$-torsor on $T\dottimes\Z_p$ comes from the restriction of the (trivial) $\GG_{\Phi}$-torsor on $A_{\inf}(B^+) = W(B^+)$ given by the (trivial) $\GG_{\Phi}$-BKF-module $\mathcal{T}_{\Phi, G}$ constructed in \cite[Lem. 4.6.6]{PR24}. Recall that the weight filtration on $M_{\inf, \Phi}$ is the filtration of BKF modules of $0 \subset \mT[p^{\infty}] \subset \mG[p^{\infty}] \subset \mQ[p^{\infty}]$. Over $W(B^+) = \prod_j W(C_j^+)$, this filtration becomes constant, thus the $\GG_{\Phi}$-torsor $\mathcal{T}_{\Phi, G}$ reduces to a $\QQ_{\Phi}$-torsor $\mathcal{T}_{\Phi, W_{\bullet}}$. We have descent data of the weight filtration since the weight filtration over $T$ is the one pulled back from $S$, the $\GG_{\Phi}$-torsor $\PPs_{\Phi, G, S}$ reduces to a $\PP_{\Phi}$-torsor $\PPs_{\Phi, W_{\bullet}, S}$, and $\ovl{\PPs}_{\Phi, W_{\bullet}, S} = \ovl{\QQ}_{\Phi}\times^{\QQ_{\Phi}}\PPs_{\Phi, W_{\bullet}, S}$.

       The Frobenius structure $\phi_{\PPs_{\Phi, G, S}}$ comes from the Frobenius structure $\phi_{\VVs_{\Phi, S}}$ on $\VVs_{\Phi, S}$, which comes from the restriction of the Frobenius structure $\phi_{M_{\inf, \Phi}}$ on the BKF module $M_{\inf, \Phi}(R^{\sharp+})$ of $\mQ[p^{\infty}]$. Since the Frobenius structure on $M_{\inf, \Phi}(R^{\sharp+})$ preserves the weight filtration, $\phi_{\PPs_{\Phi, G, S}}$ reduces to the $\QQ_{\Phi}$-torsor isomorphism $\phi_{\PPs_{\Phi, W_{\bullet}, S}}$.
        \end{proof}

   \begin{lem}\label{lem: shtuka, reduces to Q shtukas}
      The $\GG_{\Phi}$-shtuka $(\PPs_{\Phi, G}, \phi_{\PPs_{\Phi, G}})$ over $\Shum{K_{\Phi}}^{\Dia/} \to \Spd(\OO_E)$ with one leg bounded by $\mu_{\Phi}$ reduces to a $\QQ_{\Phi}$-shtuka $(\PPs_{\Phi, W^{\bullet}}, \phi_{\PPs_{\Phi, W^{\bullet}}})$ with one leg bounded by $\mu_{\Phi}$.
   \end{lem}
   \begin{proof}
       By Lemma \ref{lem: shtuka, reduces to Q torsor}, it suffices to show that the leg of the $\QQ_{\Phi}$-shtuka $(\PPs_{\Phi, W^{\bullet}}, \phi_{\PPs_{\Phi, W^{\bullet}}})$ is bounded by $\mu_{\Phi}$. It suffices to check the geometric points. Over generic fiber, this follows from Corollary \ref{cor: generic fiber, shtuka, G reduces to P}. Over special fiber, let $D$ be an algebraically closed field, and $x: \Spec D \to \Shum{K_{\Phi}}$ be a geometric point on the special fiber, we need to show the pullback $\QQ_{\Phi}$-shtuka $x^*(\PPs_{\Phi, W_{\bullet}}, \phi_{\PPs_{\Phi, W_{\bullet}}})$ is bounded by $\mu_{\Phi}$. As in \cite[\S 4.6.3, Step(3)]{PR24}, we lift $\Spa D$ to a point $\Tilde{x}: \Spa(C, \OO_C) \to (\wdh{\Shum{K_{\Phi}}})^{\ad}$, where $C$ is a complete non-archimedean algebraically closed field of characteristic $0$, such that $\OO_C/\mathfrak{m}_C = D$. By \cite[\S 4.6.3, Step(2)]{PR24}, the $\GG_{\Phi}$-shtuka $\Tilde{x}^*(\PPs_{\Phi, G}, \phi_{\PPs_{\Phi, G}})$ on $\Spa(C, \OO_C)$ extends to a $\GG_{\Phi}$-BKF module $(\mathcal{T}_{\Phi, G}, \phi_{\mathcal{T}_{\Phi, G}})$ on $W(\OO_{C^{\flat}})$ and then to a $\GG_{\Phi}$-BKF module $(\ovl{\mathcal{T}}_{\Phi, G}, \phi_{\ovl{\mathcal{T}}_{\Phi, G}})$ on $W(D)$. By the correspondence of $\GG_{\Phi}$-BKF modules on $W(D)$ and $\GG_{\Phi}$-shtukas on $\Spec D$ (\cite[Thm. 2.3.8]{PR24}), and by \cite[Prop. 2.4.6]{PR24}, the $\GG_{\Phi}$-shtuka $x^*(\PPs_{\Phi, G}, \phi_{\PPs_{\Phi, G}})$ on $\Spec D$ is the specialization of the $\GG_{\Phi}$-shtuka $\Tilde{x}^*(\PPs_{\Phi, G}, \phi_{\PPs_{\Phi, G}})$ on $\Spa(C, \OO_C)$. In the proof of Lemma \ref{lem: shtuka, reduces to Q torsor}, we see that the $\GG_{\Phi}$-torsor $\mathcal{T}_{\Phi, G}$ on $W(\OO_{C^{\flat}})$ reduces to a $\QQ_{\Phi}$-torsor $\mathcal{T}_{\Phi, W_{\bullet}}$, thus the $\GG_{\Phi}$-torsor $\ovl{\mathcal{T}}_{\Phi, G}$ on $W(D)$ reduces to a $\QQ_{\Phi}$-torsor $\ovl{\mathcal{T}}_{\Phi, W_{\bullet}}$. In particular, the $\QQ_{\Phi}$-shtuka $x^*(\PPs_{\Phi, W_{\bullet}}, \phi_{\PPs_{\Phi, W_{\bullet}}})$ on $\Spec D$ is the specialization of the $\QQ_{\Phi}$-shtuka $\Tilde{x}^*(\PPs_{\Phi, W_{\bullet}}, \phi_{\PPs_{\Phi, W_{\bullet}}})$ on $\Spa(C, \OO_C)$.
       
       In particular, the relative position of $\phi_{\PPs_{\Phi, W^{\bullet}}}(\Frob^*(\PPs_{\Phi, W^{\bullet}}))$ and $\PPs_{\Phi, W^{\bullet}}$ at $\Tilde{x}$ gives a $\Spa(C, \OO_C)$-point in $\Gra{Q_{\Phi}, \mu_{\Phi}^{-1}}$, which lifts to a $\Spa(\OO_C)$-point in $\Gra{\QQ_{\Phi}}$ whose specialization gives the relative position of $\phi_{\PPs_{\Phi, W^{\bullet}}}(\Frob^*(\PPs_{\Phi, W^{\bullet}}))$ and $\PPs_{\Phi, W^{\bullet}}$ at $x$.
       
       By construction, $\vM_{\QQ_{\Phi}, \mu_{\Phi}} = \underline{|\Gra{Q_{\Phi}, \mu_{\Phi}^{-1}}|^{wgc}} \times_{\underline{|\Gra{\QQ_{\Phi}, \OO_E}|}} \Gra{\QQ_{\Phi}, \OO_E}$ (where $|\Gra{Q_{\Phi}, \mu_{\Phi}^{-1}}|^{wgc}$ is the weakly generalizing closure of $|\Gra{Q_{\Phi}, \mu_{\Phi}^{-1}}|$ in $|\Gra{\QQ_{\Phi}, \OO_E}|$). Since $|\Spa(C, \OO_C)|$ is open and dense in $|\Spa(\OO_C)|$ and the image of $|\Spa(\OO_C)|$ in $|\Gra{\QQ_{\Phi}, \OO_E}|$ is weakly generalizing, then $\Spa(\OO_C) \to \Gra{\QQ_{\Phi}}$ factors through $\vM_{\QQ_{\Phi}, \mu_{\Phi}}$, the $\QQ_{\Phi}$-shtuka $x^*(\PPs_{\Phi, W_{\bullet}}, \phi_{\PPs_{\Phi, W_{\bullet}}})$ on $\Spec D$ is bounded by $\mu_{\Phi}$.
   \end{proof}
   
\subsubsection{$\PPs_{\Phi}$-shtukas}

    \begin{prop}\label{prop: shtuka, reduces to P shtukas, integral}
        The $\GG_{\Phi}$-shtuka $(\PPs_{\Phi, G}, \phi_{\PPs_{\Phi, G}})$ over $\Shum{K_{\Phi}}^{\Dia/} \to \Spd(\OO_E)$ with one leg bounded by $\mu_{\Phi}$ reduces to a $\PP_{\Phi}$-shtuka $(\PPs_{\Phi}, \phi_{\PPs_{\Phi}})$ with one leg bounded by $\mu_{\Phi}$. This extends Corollary \ref{cor: generic fiber, shtuka, G reduces to P} over the integral base.
    \end{prop}
    \begin{proof}
        By Lemma \ref{lem: shtuka, reduces to Q torsor} and \ref{lem: shtuka, reduces to Q shtukas}, the $\GG_{\Phi}$-shtuka $(\PPs_{\Phi, G}, \phi_{\PPs_{\Phi, G}})$ with one leg bounded by $\mu_{\Phi}$ reduces to the $\QQ_{\Phi}$-shtuka $(\PPs_{\Phi, W^{\bullet}}, \phi_{\PPs_{\Phi, W^{\bullet}}})$ with one leg bounded by $\mu_{\Phi}$. Corollary \ref{cor: generic fiber, shtuka, G reduces to P} says that the $\QQ_{\Phi}$-shtuka $(\PPs_{\Phi, W^{\bullet}, E}, \phi_{\PPs_{\Phi, W^{\bullet}, E}})$ with one leg bounded by $\mu_{\Phi}$ reduces to a $\PP_{\Phi}$-shtuka $(\PPs_{\Phi, E}, \phi_{\PPs_{\Phi, E}})$ with one leg bounded by $\mu_{\Phi}$. Then the result follows from Proposition \ref{prop: generic fiber reduce, then integral reduce}.
    \end{proof}

\subsubsection{Rest of axioms}

\begin{prop}\label{prop: hodge type satisfies axioms}
   Let $(G, X)$ be a Hodge-type Shimura datum, and let $\GG$ be any quasi-parahoric model of $G_{\bb{Q}_p}$. Then there exists a family of integral models $\lrbracket{\Shum{K_{\Phi}}(P_{\Phi}, D_{\Phi})}_{K_{\Phi}^p}$ that is canonical in the sense of Definition \ref{def: canonical model for mixed Shimura data}, and is adapted with $P_{\Phi} \to G_{\Phi, h}$ in the sense of Definition \ref{def: well-adapted, 2}.
\end{prop}

As in Remark \ref{rk: Hodge-Type case, no c}, $P_{\Phi} = P_{\Phi}^c$, so it is trivially true that the integral model is adapted with $P_{\Phi} \to P_{\Phi}^c$ in the sense of \ref{def: well-adapted}.
\begin{proof}
   It suffices to deal with stabilizer quasi-parahoric $\GG$ by the functoriality result (Proposition \ref{prop: functoriality of canonical integral models}, case $(2)$). Let $\GG$ be a stabilizer quasi-parahoric. We show that the family of integral models $\lrbracket{\Shum{K_{\Phi}}(P_{\Phi}, D_{\Phi})}_{K_{\Phi}^p}$ defined via relative normalizations is canonical. This follows from the Siegel case (Lemma \ref{lem: Siegel type satisfies axioms}) and the functoriality result (Proposition \ref{prop: functoriality of canonical integral models}, case $(1)$). Note that the condition $(b)$ is not automatic, we verified it in \cite[Prop. 1.2(2)]{Mao25b} under a fixed adjusted Hodge embedding. Also note that, for the integral canonical model $\Shum{K}(G, X)$ of a Siegel-type Shimura variety $\shu{K}(G, X)$ with $K_p$ being hyperspecial, $\Shum{K_{\Phi, h}}(G_{\Phi, h}, D_{\Phi, h})$ is again the integral canonical model of the Siegel-type Shimura variety $\shu{K_{\Phi, h}}(G_{\Phi, h}, D_{\Phi, h})$ with $K_{\Phi, h}$ being hyperspecial, for each $\Phi \in \ca{CLR}(G, X)$, see \cite[Thm. 3.58(2)]{Mao25b}.
\end{proof}

\begin{lem}\label{lem: Siegel type satisfies axioms}
    Proposition \ref{prop: hodge type satisfies axioms} holds when $(G, X) = (G^{\dd}, X^{\dd})$ is a Siegel Shimura datum and when $\GG$ is hyperspecial.
\end{lem}
\begin{proof}
    In this case, we have moduli interpretations for the smooth integral models $\Shum{K_{\Phi}}(P_{\Phi}, D_{\Phi})$ and $\Shum{K_{\Phi, h}}(G_{\Phi, h}, D_{\Phi, h})$ (\cite[\S 2.2.14, 2.2.15]{Mad19}), and a projection $\Shum{K_{\Phi}}(P_{\Phi}, D_{\Phi}) \to \Shum{K_{\Phi, h}}(G_{\Phi, h}, D_{\Phi, h})$ defined via the forgetful functor (\cite[\S 2.2.16]{Mad19}); then $\Shum{K_{\Phi}}(P_{\Phi}, D_{\Phi})$ is adapted with $P_{\Phi} \to G_{\Phi, h}$. Now we verify those axioms in Definition \ref{def: canonical model for mixed Shimura data}. 

    Axiom $(1)$ follows from the N\'eron-Ogg-Shafarevich criterion for $1$-motives; see \cite[Appendix, Lem. A.3.5]{Mad19}. Indeed, recall that $\Shum{K_{\Phi}}(P_{\Phi}, D_{\Phi})$ represents the following moduli functor: given a $\Z_p$-scheme $S$, $\Shum{K_{\Phi}}(P_{\Phi}, D_{\Phi})(S)$ is the set of isomorphism classes of tuples $(\mQ, \lambda, \alpha, \alpha^{\vee}, \epsilon)$ over $S$, where $(\mQ, \lambda)$ is a polarized $1$-motive over $S$, equipped with isomorphisms of sheaves of $\Z$-modules over $S$:
    \[ \alpha: \gr_0^W\underline{V}^g(\Z) \cong \mQ^{\et},\quad \alpha^{\vee}: \gr_0^W\underline{V}^{\vee, g}(\Z)(\nu) \cong \mQ^{\vee, \et} = \mQ^{\mult}, \]
    where $V^{\vee}(\nu)$ is the dual representation $V^{\vee}$ twisted by the similitude character $\nu: \GSp(V, \psi) \to \Gm$, and $\epsilon \in H^0(S, \Isom(V^g(\hat{\Z}^p), \hat{T}^p(\mQ))/K_{\Phi, G}^p)$ is a $K_{\Phi, G}^p$-level structure, where $K_{\Phi, G}^p = g^pK^pg^{p, -1}$. See \cite[\S 2.2.14]{Mad19} for conventions. 

    Let $x \in \Shum{K_{\Phi, p}}(P_{\Phi}, D_{\Phi})(R[1/p])$, it gives a tuple $(\mQ_{\eta}, \lambda_{\eta}, \alpha_{\eta}, \alpha^{-1}_{\eta}, \epsilon_{\eta})$ over $\Spec R[1/p]$. By \cite[Appendix, Lem. A.3.5]{Mad19}, $\mQ_{\eta}$ (resp. $\mQ_{\eta}^{\vee}$) extends to a $1$-motif $\mQ$ (resp. $\mQ^{\vee}$) over $\Spec R$. Note that $\alpha_{\eta}$, $\alpha^{-1}_{\eta}$, $\epsilon_{\eta}$ extend to $\alpha$, $\alpha^{-1}$, $\epsilon$ respectively by the \'etaleness, and $\lambda_{\eta}$ extends to $\lambda$ using the rigidity theorem for homomorphisms of semi-abelian schemes (thus of 1-motives), see \cite[I, Prop. 2.7]{FC90}. The tuple $(\mQ, \lambda, \alpha, \alpha^{-1}, \epsilon)$ over $\Spec R$ gives the wanted lifting.
    
    Axiom $(2)$: This follows from moduli interpretations. Also see \cite[Lem. 4.25]{Wu25}.
    
    Axiom $(3)$: This is a special case of Proposition \ref{prop: shtuka, reduces to P shtukas, integral}.

    Axiom $(4)$: We replace the Serre-Tate theorem for abelian schemes with the Serre-Tate theorem for $1$-motives. Fix a $k$-point $x \in \Shum{K_{\Phi}}(P_{\Phi}, D_{\Phi})(k)$. We have a tuple $(\mQ_0, \lambda_0, \alpha_0, \alpha^{\vee}_0, \epsilon_0)$ over $\Spec k$, the formal completion $\hat{S}_x:=(\Shum{K_{\Phi}}(P_{\Phi}, D_{\Phi}))^{\wedge}_{/x}$ represents the functor $\mathrm{Def}: \mathrm{Art}_{W(k)} \to \mathrm{Set}$ sending $S$ to isomorphism classes of deformations over $S$ of the tuple $(\mQ_0, \lambda_0, \alpha_0, \alpha^{\vee}_0, \epsilon_0)$.

    On the other hand, $\mathcal{M}^{\intg}:=\mathcal{M}^{\intg}_{\GG_{\Phi, b_{\Phi, x}, \mu_{\Phi}}}$ is a Rapoport-Zink space, the deformation space of $p$-divisible groups, for example, see \cite[Thm. 25.1.2]{SW20}. Let $(X_1, \lambda_1)$ be a polarized $p$-divisible group in the isogeny class determined by $[b_{\Phi, x}] \in B(G, \lrbracket{\mu^{-1}}) = B(G_{\Phi}, \lrbracket{\mu_{\Phi}^{-1}})$, given a $\Z_p$-scheme $S$, $\mathcal{M}^{\intg}(S)$ is the set of isomorphism classes of $(X, \lambda, \rho)$ over $S$, where $(X, \lambda)$ is a polarized $p$-divisible group over $S$ and $\rho: (X, \lambda) \times_S S/p \to (X_1, \lambda_1) \times_k S/p$ is a quasi-isogeny. Let $x_0 \in \mathcal{M}^{\intg}(k)$ be the base-point associated with $(X_1, \lambda_1, \identity)$. Since the quasi-isogeny $\rho$ over $S/p$ is rigid in the sense that there is a unique $\rho$ lifting the isomorphism $\rho \otimes S/\mathfrak{m}_S$, for example, see \cite[\S 2.3]{kim2018rapoport}. Therefore, $(\mathcal{M}^{\intg})^{\wedge}_{/x_0}$ represents the functor $\mathrm{Def}_{(X_1, \lambda_1)}: \mathrm{Art}_{W(k)} \to \mathrm{Set}$ sending $S$ to isomorphism classes of deformations over $S$ of the polarized $p$-divisible group $(X_1, \lambda_1)$.

    Let $(X_0, \lambda_0)$ be the $p$-divisible group $(\mQ_0, \lambda_0)[p^{\infty}]$, then we have an isogeny $\rho_{01}: (X_0, \lambda_0) \to (X_1, \lambda_1)$. $(X_0, \lambda_0, \rho)$ gives a $k$-point in $\mathcal{M}^{\intg}$, which we denote by $x_1$. 
    
    By the main theorem of \cite{bertapelle2019deformations}, the Serre-Tate theorem holds for $1$-motives: let $R \in \mathrm{Art}_{W(k)}$, $\mathcal{M}_1(R)$ be the category of $1$-motives over $R$. Let $M$ (resp. $B$) be a $1$-motive (resp. a $p$-divisible group) over $R$, $M_0$ (resp. $B_0$) be its reduction over $\Spec k$. Let $\mathrm{Def}(R, k)$ be the category of triples $(M_0, B, e)$, where $M_0$ is a $1$-motive over $\Spec k$, $B$ is a $p$-divisible group over $R$, and $e: B_0 \to M_0[p^{\infty}]$ is an isomorphism of $p$-divisible groups. There is a natural equivalence of categories: 
    \[ \Delta_R: \mathcal{M}_1(R) \to \mathrm{Def}(R, k),\quad M \to (M_0, M[p^{\infty}], \textrm{natural}\ e). \]

    Consider the morphism $\hat{S}_x \to (\mathcal{M}^{\intg})^{\wedge}_{/x_1}$ defined by sending $(\mQ, \lambda, \alpha, \alpha^{\vee}, \epsilon)$ to $((\mQ, \lambda)[p^{\infty}], \identity)$. We claim this is an equivalence of categories. By Serre-Tate for $1$-motives, and by the fact that $(\alpha, \alpha^{\vee})$ and $\epsilon$ are rigid by the \'etaleness (thus determined by their reductions $(\alpha_0, \alpha_0^{\vee})$ and $\epsilon_0$ over $\Spec k$ respectively), it suffices to show that, given a polarization of $p$-divisible groups $\lambda[p^{\infty}]: \mQ[p^{\infty}] \to \mQ^{\vee}[p^{\infty}]$, we can lift it canonically to a polarization of $1$-motives $\lambda: \mQ \to \mQ^{\vee}$. Recall that a polarization $\lambda: \mQ \to \mQ^{\vee}$ is a morphism between $1$-motives such that $\lambda^{\mathrm{ab}}: \mQ^{\mathrm{ab}} \to \mQ^{\vee, \mathrm{ab}}$ is a polarization of abelian schemes, $\lambda^{\et}: \mQ^{\et} \to \mQ^{\mult}$ is injective, and $(\lambda^{\et})^{\vee} = \lambda^{\mult}$. The full-faithfulness of $\Delta_R$ gives a unique morphism $\lambda: \mQ \to \mQ^{\vee}$ and its proof also shows that $\lambda^{\mathrm{ab}}$ is a polarization, $\lambda^{\et}$ and $\lambda^{\mult}$ satisfy the properties since they are \'etale in nature and $\lambda_0^{\et}$ and $\lambda_0^{\mult}$ do satisfy the properties.

    Therefore, we have an isomorphism $\hat{S}_x \cong (\mathcal{M}^{\intg})^{\wedge}_{/x_1}$. We show that the isomorphism $\hat{S}_x \cong (\mathcal{M}^{\intg})^{\wedge}_{/x_1}$ matches the supported shtukas. By \cite[Thm. 25.1.2]{SW20}, the identification of the Rapoport-Zink space with $\mathcal{M}^{\intg}$ is to evaluate the universal $p$-divisible group at perfectoid space $S = \Spa(R, R^+) \in \Perf$ to get a BKF module over $A_{\inf}(R, R^+)$, and then restrict the BKF module to $S\dot{\times}\Z_p$ to get the desired shtuka; also, the quasi-isogeny $\rho$ provides a framing. Under this identification, since $\hat{S}_x \cong (\mathcal{M}^{\intg})^{\wedge}_{/x_1}$ matches the universal $p$-divisible groups, it matches the supported shtukas. Finally, after changing the base point, $(\mathcal{M}^{\intg})^{\wedge}_{x_1} \cong (\mathcal{M}^{\intg})^{\wedge}_{x_0}$, we are done.
\end{proof}

Finally, to prove the axiom $(1)$ in abelian-type case, we need following enhanced result:
\begin{lem}\label{lem: enhanced (1)}
    Keep the notation from Proposition \ref{prop: hodge type satisfies axioms}. Let $l$ be a prime $\neq p$, $K_{\Phi}^{l,p} =K_\Phi^p\cap P_\Phi(\A^{p,l})$, $K_{\Phi}^{l} =K_{\Phi,p}K_\Phi^{p,l}$. For every discrete valuation ring $R$ of mixed characteristic over $\OO_E$, we have
        \begin{equation*}
            \shu{K_{\Phi}^l}(P_{\Phi}, D_{\Phi})(R[1/p]) =\Shum{K_{\Phi}^l}(P_{\Phi}, D_{\Phi})(R).
        \end{equation*}
\end{lem}
\begin{proof}
    This follows from the proof of the Siegel case (Lemma \ref{lem: Siegel type satisfies axioms}) and the proof of the functoriality result (Proposition \ref{prop: functoriality of canonical integral models}, case $(1)$). Note that in \cite[Appendix, Lem. A.3.5]{Mad19}, one only needs to trivialize $T^l(\mQ_{\eta})$, rather than the whole $\hat{T}^p(\mQ_{\eta})$.
\end{proof}

\subsection{Abelian-type case}\label{subsec-ab}
\subsubsection{Construction}\label{subsubsection}
Let us recall the construction in \cite{Wu25}. The key steps are listed in Construction \ref{const-step-1-gp}, \ref{const-gp-fromab-to-hodge}, and \ref{const-scheme-from-hodge-to-ab} for convenience.
\begin{construction}\label{const-step-1-gp}\upshape
Let $(G_2,X_2)$ be an abelian-type Shimura datum with an associated Hodge-type Shimura datum $(G,X)$; note that, in \cite{Wu25}, it was denoted by $(G_0,X_0)$. There is a central isogeny $\pi^\der:G^\der\to G^{\der}_2$ with kernel $C^\der$. \par
Let $G':=G_2\times^{G_2^\der} (G/C^\der)$. 
There is an embedding $\pi^a: G_2\to G'$ and a map $\pi^b: G\to G'$ such that they induce morphisms between Shimura data $\pi^a:(G_2,X_2)\to (G',X_a)$ and $\pi^b:(G,X)\to (G',X_b)$. 
By replacing $(G,X)$ with a conjugate by $g\in G^\ad(\bb{Q})$, we assume that the images of $X_2$ and $X$ in $X^\ad$ have nontrivial intersection. We then have that the images of $X_a$ and $X_b$ in $X^\ad$ coincide.
Then $\pi^a$ and $\pi^b$ satisfy the conditions in \cite[\S 1.4.3]{Wu25}. See \cite[\S 1.4.3 and Const. 1.46]{Wu25}.
\end{construction}
We recall more group-theoretic preparations.
\begin{construction}\label{const-gp-fromab-to-hodge}\upshape
We start with $K_2=K_{2,p}K^p_2\sbst G_2(\A)$ and an admissible cone decomposition $\Sigma_2$ for $(G_2,X_2,K_2)$. 
\begin{itemize}
\item By \cite[Prop. 1.35, Prop. 1.37 and Prop. 1.47]{Wu25}, we can find $K':=K_p'K^{\prime,p}\sbst G'(\A)$ containing $K_2$, and an admissible smooth projective $ZP$-invariant cone decomposition $\Sigma$ for $(G',X_a,K')$ and $(G',X_b,K')$ inducing a cone decomposition $\Sigma_2'$ refining $\Sigma_2$ for $(G_2,X_2,K_2)$, such that the strata between toroidal (with cones chosen as above) and minimal compactifications of $(G_2,X_2)$ and $(G',X_a)$ are all open and closed embeddings, and such that the map between mixed Shimura varieties for any cusp label representative $\Phi_2\in \ca{CLR}(G_2,X_2)$ mapping to a cusp label representative $\Phi\in \ca{CLR}(G',X_a)$ is an isomorphism. (Note that $K'_p$ is chosen such that $K'_p\cap G_2(\bb{Q}_p)=K_{2,p}$.)  There is an identification $\ca{CLR}(G',X_a)=\ca{CLR}(G',X_b)$ by definition and construction.
\item Denote $I_{G'/G,K'}:=\stb_{G'(\bb{Q})}(X)\pi^b(G(\A))\bss G'(\A)/K'$.
Choose a complete system $\{g_\alpha\}_{\alpha\in I_{G'/G,K'}}$ of representatives of $I_{G'/G,K'}$ in $G'(\A)$.
Let $K^\alpha_p:=\pi^{b,-1}(g_\alpha K'_p g_\alpha^{-1})\sbst G(\bb{Q}_p)$. Choose suitable neat open compact $K^{\alpha,p}\sbst G(\Ap)$, and denote $K^\alpha:=K^\alpha_pK^{\alpha,p}$.
We obtain the induced cone decompositions $\Sigma^\alpha:=\pi^{b,*}(\Sigma)$. 
\item In our case, we assume that $K_p'$ is a quasi-parahoric subgroup, which is the $\bb{Z}_p$-points of a quasi-parahoric group scheme $\ca{G}'$ corresponding to ${x}\in \ca{B}(G',\bb{Q}_p)$. \par
Following the construction in Case (STB) of \cite[\S 4.2]{Wu25}, there is a Siegel-type Shimura datum $(G^\ddag,X^\ddag)=(\mrm{GSp}(V,\psi),X^\ddag)$ with an embedding $\rho:(G,X)\hookrightarrow (G^\ddag,X^\ddag)$ obtained as follows:\par 
For any $\alpha$, we choose $K_p^{\ddag,\alpha}\sbst G^\ddag(\bb{Q}_p)$ and neat open compact $K^{\ddag,\alpha,p}\sbst G^\ddag(\Ap)$ such that
    $$\rho^\alpha:(G,X,K_p^\alpha)\hookrightarrow (G^{\ddag,\alpha},X^{\ddag,\alpha},K_p^{\ddag,\alpha})$$
is an adjusted Hodge embedding satisfying the setup in both \cite[\S 3.1]{Mad19} and \cite{daniels2024conjecture}.\par 
In fact, $G^{\ddag,\alpha}:=\mrm{GSp}(W^\alpha,\psi^\alpha)$,  $K_p^{\ddag,\alpha,\sharp}=\stb_{G^{\ddag,\alpha}(\bb{Q}_p)}(W_{\bb{Z}_p}^\alpha)$ for some self-dual $\bb{Z}_p$-lattice $W_{\bb{Z}_p}^\alpha$ in $W_{\bb{Q}_p}^\alpha$. Moreover, $W_{\bb{Z}_p}^\alpha$ comes from a self-dual $\bb{Z}$-lattice $W_{\bb{Z}}^\alpha$.\par
Furthermore, there are positive integers $n_\alpha$ such that $(W^\alpha,\psi^\alpha)^{\perp n_\alpha}=(V,\psi)$ for a symplectic space $(V,\psi)$ for all $\alpha$. By defining $V^{\alpha}_?:=W_?^{\alpha,\perp n_\alpha}$ for $?=\bb{Q},\bb{Z},\bb{Z}_p$, $(G^\ddag,X^\ddag)$ as the Siegel Shimura datum defined by $(V,\psi)$, and $K_p^{\ddag,\alpha}:=\stb_{G^\ddag(\bb{Q}_p)}(V_{\bb{Z}_p}^\alpha)$, we obtain from $\rho^\alpha$ the Hodge embeddings 
$$\iota^\alpha:(G,X,K_p^\alpha)\hookrightarrow (G^\ddag,X^\ddag,K_p^{\ddag,\alpha}).$$
After this slight modification, \cite[Thm. 3.58]{Mao25b} still holds.
\item By \cite[Prop. 1.47]{Wu25} again, cone decompositions $\Sigma^{\ddag,\alpha}$ can be chosen for $(G^\ddag,X^\ddag,K^{\ddag,\alpha})$ such that, for any $\alpha$, $\Sigma^\alpha$ is induced by both $\Sigma$ via $\pi^b(g_\alpha)$ and $\Sigma^{\ddag,\alpha}$ via $\iota^\alpha$, and such that $\mathscr{S}_{K^\alpha}^{\Sigma^\alpha}$ can be constructed satisfying Axiom \ref{axiom-good-compactification} according to \cite{Mad19}.
\end{itemize}
\end{construction}
We now recall some schematic constructions in \cite{Wu25}.
\begin{construction}\label{const-scheme-from-hodge-to-ab}\upshape
We recall the construction of integral models of boundary mixed Shimura varieties $\mathscr{S}_{K_{\Phi_2}}$ for $\Phi_2\in\ca{CLR}(G_2,X_2)$ following \cite[Sec. 4.2]{Wu25} in Case (STB) and its normalization (DL) case to parahoric levels. 
\begin{enumerate}
\item We introduce more group-theoretic setup related to Bruhat-Tits theory. In general, there are maps between groups $\pi^a:K_{2,p}\hookrightarrow K_p'$ and $\pi^b(g_\alpha): K_p^\alpha\to g_\alpha K_p' g_\alpha^{-1}$. Starting with $\Phi_2\in \ca{CLR}(G_2,X_2)$, let $\Phi$ be the image of $\Phi_2$ in $\ca{CLR}(G',X_a)\iso \ca{CLR}(G',X_b)$. \par
Fix a point in the reduced building $x\in \ca{B}(G_2,\bb{Q}_p)$. We denote by $K_{2,p}^\circ$ (resp. $K_{2,p}^{\mrm{stb}}$ and $K_{2,p}$) the parahoric (resp. the stabilizer quasi-parahoric and resp. a general quasi-parahoric) subgroup of $G_2(\bb{Q}_p)$ corresponding to $x$. For $\Phi_2$, define $K_{\Phi_2,p}^\mmin:=P_{\Phi_2}(\bb{Q}_p)\cap g_{\Phi_2}K_{2,p}^\circ g_{\Phi_2}^{-1}$ (resp. $K_{\Phi_2,p}^{\mrm{stb}}:=P_{\Phi_2}(\bb{Q}_p)\cap g_{\Phi_2} K_{2,p}^{\mrm{stb}}g_{\Phi_2}^{-1}$ and resp. $K_{\Phi_2,p}:=P_{\Phi_2}(\bb{Q}_p)\cap g_{\Phi_2} K_{2,p}g_{\Phi_2}^{-1} $); note that $K_{\Phi_2,p}^{\mrm{stb}}$ is still a stabilizer quasi-parahoric subgroup. Let $K_p^{\prime,\mrm{stb}}$ be the stabilizer quasi-parahoric subgroup in $G'(\bb{Q}_p)$ corresponding to $x$. Similarly, define $K_{\Phi,p}^{\mrm{stb}}:=P_\Phi(\bb{Q}_p)\cap g_{\Phi_2} K_p^{\prime,\mrm{stb}}g_{\Phi_2}^{-1}$ and $\K_{\Phi,p}^{\mrm{stb}}:=ZP_\Phi(\bb{Q}_p)\cap g_{\Phi_2} K_p^{\prime,\mrm{stb}}g_{\Phi_2}^{-1}$. Note that we can still use Definition \ref{def: P comes from boundary} and Lemma \ref{lem: quasi-parahoric} to get that $\K_{\Phi,p}^{\mrm{stb}}$ is a stabilizer quasi-parahoric subgroup. The prime-to-$p$ neat open compact subgroups are chosen and written in an obvious pattern.
\item Denote by $(P_\Phi,D_\Phi^b)$ (resp. $(P_\Phi,D^a_\Phi)$) the mixed Shimura datum associated with $\Phi\in\ca{CLR}(G,X_b)$. We construct $\{\mathscr{S}_{\K^{\mrm{stb}}_{\Phi}}(ZP_\Phi,ZP_\Phi(\bb{Q})D^b_\Phi)\}_{K^{\prime,p}}$ and $\{\mathscr{S}_{\K^{\mrm{stb}}_{\Phi}}(ZP_\Phi,ZP_\Phi(\bb{Q})D^a_\Phi)\}_{K^{\prime,p}}$. 
Note that the Shimura datum $(G,X)$ is chosen such that any place $v'$ of $\bb{E}':=\bb{E}(G,X)\cdot \bb{E}(G_2,X_2)$ over a place $v_2$ of $\bb{E}_2:=\bb{E}(G_2,X_2)$ over $p$ splits completely over $v_2$.\par
We define $\mathscr{S}_{\K_\Phi^{\mrm{stb}}}(ZP_\Phi,ZP_\Phi(\bb{Q})D^b_\Phi)$ over $\ca{O}':=\ca{O}_{\bb{E}'}\otimes_{\ca{O}_{\bb{E}}}\ca{O}_{\bb{E},(v)}$ for any place $v$ of $\bb{E}:=\bb{E}(G,X)$. 
For any cusp label $\Phi\in\ca{CLR}(G',X_b)$ extending to a $ZP$-cusp label $[ZP^b(\Phi)]$, by \cite[(4.22)]{Wu25}, we construct 
\begin{equation}\label{eq-def-zp-cusp-1}\mathscr{S}_{\K_\Phi^{\mrm{stb}}}(ZP_\Phi,ZP_\Phi(\bb{Q})D^b_\Phi):=\disju_{\alpha\in I_{G'/G,K^{\prime,\mrm{stb}}}}\disju_{\pi^b(g_0^\alpha)\alpha\sim g^b}\disju_{[\sigma^\alpha]\in [\sigma]_{ZP}}\mathscr{S}_{K_{\Phi_0^\alpha}}/\Delta_{\Phi_0^\alpha,K^{\prime,\mrm{stb}}}(ZP_\Phi).\end{equation}
The group $\Delta_{\Phi_0^\alpha,K^{\prime,\mrm{stb}}}(ZP_\Phi)$ is a group that acts on $\mathscr{S}_{K_{\Phi_0^\alpha}}$ through a finite group; this group is defined and functorial in \emph{all} neat $K^{\prime,\mrm{stb}}$; its quotient is finite {\'e}tale by \cite[Lem. 4.25]{Wu25}. Writing in this form a priori depends on the choice of $\mathscr{S}_{K^\alpha}^{\Sigma^\alpha}$ and, in particular, on the choice of $\Sigma^\alpha$. But we can also re-label the disjoint union in (\ref{eq-def-zp-cusp-1}) as 
\begin{equation}\label{eq-def-zp-cusp-2}\mathscr{S}_{\K_\Phi^{\mrm{stb}}}(ZP_\Phi,ZP_\Phi(\bb{Q})D^b_\Phi):=\disju_{\beta\in I_{\Phi,K^{\prime,\mrm{stb}}}} \mathscr{S}_{K_{\Phi_\beta}}/\Delta_{\Phi_\beta,K^{\prime,\mrm{stb}}}(ZP_\Phi),\end{equation}
where $\Phi_\beta$ are cusp label representatives in $ \ca{CLR}(G,X)$ that occurred in (\ref{eq-def-zp-cusp-1}) and the index set $I_{\Phi,K^{\prime,\mrm{stb}}}$ is defined in \cite[4.2.6]{Wu25}; it depends only on $\Phi$, $K^{\prime,\mrm{stb}}$, the Shimura data $(G,X)$ and the group $G'$, satisfying obvious functoriality for varying $K'\sbst G'(\A)$. The level group $K_{\Phi_{\beta}}$ is defined as $P_{\Phi_\beta}(\A)\cap g_{\Phi_\beta}K^\alpha g^{-1}_{\Phi_\beta}$, where $K_p^\alpha=\pi^{b,-1}(g_\alpha K_p^{\prime,\mrm{stb}}g_\alpha^{-1})$.\par 
Replacing $\mathscr{S}_{\K_\Phi^{\mrm{stb}}}(ZP_\Phi,ZP_\Phi(\bb{Q})D^b_\Phi)$ with $\mathscr{S}_{\overline{\K}_\Phi^{\mrm{stb}}}(ZP_\Phi,ZP_\Phi(\bb{Q})D^b_\Phi)$ (resp. $\mathscr{S}_{\K_{\Phi,h}^{\mrm{stb}}}(ZP_\Phi,ZP_\Phi(\bb{Q})D^b_\Phi)$), the integral model is constructed by replacing the Hodge-type mixed Shimura varieties $\mathscr{S}_{K_{\Phi_\beta}}$ in (\ref{eq-def-zp-cusp-2}) with $\mathscr{S}_{\overline{K}_{\Phi_\beta}}$ (resp. $\mathscr{S}_{K_{\Phi_\beta,h}}$). Again, the group action of $\Delta_{\Phi_\beta}$ factors through a finite group that acts freely on $\mathscr{S}_{\overline{K}_{\Phi_\beta}}$ (resp. $\mathscr{S}_{K_{\Phi_\beta,h}}$). Moreover, there is a sequence of morphisms
\begin{equation}\label{eq-stb-b-mix-sequence}\mathscr{S}_{\K_\Phi^{\mrm{stb}}}(ZP_\Phi,ZP_\Phi(\bb{Q})D^b_\Phi)\to \mathscr{S}_{\overline{\K}_\Phi^{\mrm{stb}}}(ZP_\Phi,ZP_\Phi(\bb{Q})D^b_\Phi)\to \mathscr{S}_{\K_{\Phi,h}^{\mrm{stb}}}(ZP_\Phi,ZP_\Phi(\bb{Q})D^b_\Phi).\end{equation}
The first morphism is a torsor under a split torus $\mbf{E}_{\K_\Phi}$, and the second morphism is an abelian scheme torsor by \cite[Prop. 4.30]{Wu25}.\par
To construct $\mathscr{S}_{\K_\Phi^{\mrm{stb}}}(ZP_\Phi,ZP_\Phi(\bb{Q})D^a_\Phi)$, we first base change $\mathscr{S}_{\K_\Phi^{\mrm{stb}}}(ZP_\Phi,ZP_\Phi(\bb{Q})D^b_\Phi)$ to $\ca{O}^{ur}$ such that $\ca{O}^{ur}=\ca{O}_{\bb{E}''}\otimes_{\ca{O}_{\bb{E}}}\ca{O}_{\bb{E},(v)}$ where $\bb{E}''$ is finite extension of $\bb{E}'$ that is unramified over $p$, and
that the generic fiber of the base change $\mathscr{S}_{\K_\Phi^{\mrm{stb}}}(ZP_\Phi,ZP_\Phi(\bb{Q})D^b_\Phi)_{\ca{O}^{ur}}$ is $\sh_{\K_\Phi^{\mrm{stb}}}(ZP_\Phi,ZP_\Phi(\bb{Q})D^b_\Phi)_{\bb{E}''}$. Thus, the integral model $\mathscr{S}_{\K_\Phi^{\mrm{stb}}}(ZP_\Phi,ZP_\Phi(\bb{Q})D^a_\Phi)$ is obtained from desending the model $\mathscr{S}_{\K_\Phi^{\mrm{stb}}}(ZP_\Phi,ZP_\Phi(\bb{Q})D^b_\Phi)_{\ca{O}^{ur}}$ with the descend datum given by the Shimura datum $(G,X_a)$. There is also a sequence of morphisms as (\ref{eq-stb-b-mix-sequence}) replacing ``$b$'' with ``$a$'' having the same property.
\item Let $K_p'$ be a quasi-parahoric subgroup associated with $K_{p}^{\prime,\mrm{stb}}$. Denote $\K_{\Phi,p}'=ZP_\Phi(\A)\cap g_\Phi K'_p g_\Phi^{-1}$. We construct $\{\mathscr{S}_{\K_\Phi}(ZP_\Phi,ZP_\Phi(\bb{Q})D^b_\Phi)\}_{K^{\prime,p}}$ as the normalization from the tower $\{\mathscr{S}_{\K_\Phi^{\mrm{stb}}}(ZP_\Phi,ZP_\Phi(\bb{Q})D^b_\Phi)\}_{K^{\prime,p}}$. 
More explicitly, we also construct $\mathscr{S}_{\K_\Phi}(ZP_\Phi,ZP_\Phi(\bb{Q})D^b_\Phi)$ as $\disju_{\beta'\in I_{\Phi,K'}}\mathscr{S}_{K_{\Phi_{\beta'}}}/\Delta_{\Phi_{\beta',K'}}(ZP_\Phi)$. Each $\mathscr{S}_{K_{\Phi_{\beta'}}}$ is constructed as follows: For Any $\Phi_{\beta'}\in \ca{CLR}(G,X)$, the integral model $\mathscr{S}_{K^{\mrm{stb}}_{\Phi_{\beta'}}}$ at stabilizer quasi-parahoric level is constructed as before. Then $\mathscr{S}_{K_{\Phi_{\beta'}}}$, the integral model assciated with the same cusp label representative but at the quasi-parahoric level, is defined by taking relative normalization from $\mathscr{S}_{\K^{\mrm{stb}}_{\Phi_{\beta'}}}(ZP_\Phi,ZP_\Phi(\bb{Q})D^b_\Phi)$. There is a finite map $$\mathscr{S}_{K_{\Phi_{\beta'}}}/\Delta_{\Phi_{\beta'},K'}(ZP_\Phi)\to \mathscr{S}_{K_{\Phi_{\beta'}}^{\mrm{stb}}}/\Delta_{\Phi_{\beta'},K^{\prime,\mrm{stb}}}(ZP_\Phi).$$ 
There is also a sequence as (\ref{eq-stb-b-mix-sequence}) for $\mathscr{S}_{\K_\Phi}(ZP_\Phi,ZP_\Phi(\bb{Q})D^b_\Phi)$.\par 
We can also construct $\mathscr{S}_{\K_\Phi}(ZP_\Phi,ZP_\Phi(\bb{Q})D^a_\Phi)$ with the same method as above. See \cite[Prop. 4.30 and Const. 4.33]{Wu25}.
\item We construct $\{\mathscr{S}_{K_{\Phi_2}^{\mmin}}\}_{K_2^p}$, $\{\mathscr{S}_{K_{\Phi_2}}\}_{K_2^p}$ and $\{\mathscr{S}_{K_{\Phi_2}^{\mrm{stb}}}\}_{K_2^p}$. In fact, they are constructed by taking the relative normalization of $\{\mathscr{S}_{\K^{\mrm{stb}}_{\Phi}}(ZP_\Phi,ZP_\Phi(\bb{Q})D^a_\Phi)\}_{K^{\prime,p}}$ in $\{\sh_{K_{\Phi_2}^{\mmin}}\}_{K_2^p}$, $\{\sh_{K_{\Phi_2}}\}_{K_2^p}$ and $\{\sh_{K_{\Phi_2}^{\mrm{stb}}}\}_{K_2^p}$, respectively.
\end{enumerate}
\hfill $\square$
\end{construction}
\subsubsection{Choosing accessible Hodge-type liftings}\label{subsubsec-accessible-hodge}
\begin{definition}\label{def-compatible-in-derived-part}
Fix a prime number $p$.
Let $(G,X)$ and $(G_1,X_1)$ be Shimura data with a central isogeny $\pi: G^\der\to G^{\der}_1$ such that $\pi$ induces an isomorphism between the two adjoint Shimura data. Let $\ca{G}^{\circ}_{1,x}$ be a parahoric group scheme corresponding to $x$ in the building $\ca{B}(G_1,\bb{Q}_p)$, and $\ca{G}^\circ_{x}$ be the corresponding parahoric group scheme of $G_{\bb{Q}_p}$.\par
We say that $\breve{K}_p^\circ:=\ca{G}_{x}^\circ(\breve{\bb{Z}}_p)$ \textbf{is accessible to} $\breve{K}^{\circ}_{1,p}:=\ca{G}_{1,x}^{\circ}(\breve{\bb{Z}}_p)$ if the map $\pi: G^\der(\breve{\bb{Q}}_p)\to G_1^{\der}(\breve{\bb{Q}}_p)$ restricts to a map 
\begin{equation}\label{eq-finer}\breve{K}_p^\circ\cap G^\der(\breve{\bb{Q}}_p)\to\breve{K}_{1,p}^{\circ}\cap G^{\der}_1(\breve{\bb{Q}}_p).\end{equation}
We say that $(G,X)$ \textbf{is accessible to} $(G_1,X_1)$ (at $p$ with respect to $\pi$) if, for any point $x$ in the reduced building of $G_{1,\bb{Q}_p}$, $\breve{K}_p^\circ$ is always accessible to $\breve{K}^{\circ}_{1,p}$.
\end{definition}
\begin{convention}\label{conv-kernel-der-pi1}
Let $H$ and $G$ be two connected reductive groups over $\bb{Q}$ or $\bQ$ with a map $f:H\to G$. Denote $\ker(f:\pi_1(H_{\bQ})_I\to \pi_1(G_{\bQ})_I)$ by $\pi(H,G)$.
\end{convention}
Given $f: H\to G$ and $g: G\to G'$, we have, by definition, an inclusion 
$\pi(H,G)\sbst \pi(H,G')$.\par
It follows from the definition of (quasi-)parahoric subgroups that
\begin{lem}\label{lem-trivial-reduction}
To show that $(G,X)$ is accessible to $(G_1,X_1)$, it suffices to show that 
$$\pi(G^\der, G)\sbst \pi(G^\der,G_1).$$
\end{lem}
\begin{prop}\label{prop-comp-der-consequence}
Fix a prime $p$. Let $(G_2,X_2)$ be an abelian-type Shimura datum. We assume that there is always a Hodge-type lifting $(G,X)$ of $(G_2,X_2)$ such that
\begin{enumerate}
    \item For any place $v'$ of $\bb{E}'$ over a place $v_2$ of $\bb{E}_2$ over $p$, $v'$ splits completely over $v_2$;
    \item $(G,X)$ is accessible to $(G_2,X_2)$ at $p$.
\end{enumerate}
Then, for any cusp label representative $\Phi_2\in\ca{CLR}(G_2,X_2)$ and any quasi-parahoric subgroup $K_{\Phi_2,p}$ as in Construction \ref{const-scheme-from-hodge-to-ab}(1), the maps
$\mathscr{S}_{K_{\Phi_2}}\to \mathscr{S}_{K_{\Phi_2}^{\mrm{stb}}}$, 
$\mathscr{S}_{\overline{K}_{\Phi_2}}\to \mathscr{S}_{\overline{K}^{\mrm{stb}}_{\Phi_2}}$
and 
$\mathscr{S}_{K_{\Phi_2,h}}\to \mathscr{S}_{K_{\Phi_2,h}^{\mrm{stb}}}$
are finite {\'e}tale. 
\end{prop}
\begin{proof}
With these two assumptions, the equality $\pi^b(K_{p}^\circ)\pi^a(K_{2,p})\cap G_2(\bb{Q}_p)=K_{2,p}$ holds. 
Indeed, it suffices to show that  $\pi^b(K_p^\circ)\cap G_2(\bb{Q}_p)\sbst K_{2,p}$. Since $\pi^b(G(\bb{Q}_p))\cap G_2(\bb{Q}_p)\sbst G_2^\der(\bb{Q}_p)$, we show that $\pi^b(K_p^\circ)\cap G^\der_2(\bb{Q}_p)\sbst K_{2,p}\cap G^\der_2(\bb{Q}_p)$. Since the kernel of $\pi^b$ is in the center of $G^\der$, we know that if $x\in G(\bb{Q}_p)$ maps to $G^{\prime,\der}(\bb{Q}_p)=G_2^\der(\bb{Q}_p)$, then $x\in G^\der(\bb{Q}_p)$. Thus, $\pi^b(K_p^\circ\cap G^\der(\bb{Q}_p))=\pi^b(K_p^\circ)\cap G^{\prime,\der}(\bb{Q}_p)$. Finally, we know from the assumption that $\pi^b(K_p^\circ\cap G^\der(\bb{Q}_p))=\pi^b(K_p^\circ)\cap G^{\prime,\der}(\bb{Q}_p)\sbst K_{2,p}\cap G_2^\der(\bb{Q}_p)$.\par
When this equality holds, for any quasi-parahoric $K_{2,p}\sbst G_2(\bb{Q}_p)$, we can choose $K_p'\sbst G'(\bb{Q}_p)$ such that $K_p'\cap G_2(\bb{Q}_p)=K_{2,p}$. Moreover, $\pi^{b,-1}(K_p')$ are quasi-parahoric subgroups. See \cite[Lem. 4.57]{Wu25}.\par
Now, $\mathscr{S}_{K_{\Phi_2}^{\mrm{stb}}}$ can be constructed by choosing the corresponding stabilizer quasi-parahoric subgroup $\K_{\Phi,p}^{\mrm{stb}}$ of $ZP_{\Phi}(\bb{Q}_p)$, and $\mathscr{S}_{K_{\Phi_2}^\mmin}$ can be constructed from $\mathscr{S}_{\K_\Phi}$ such that $\K_{\Phi,p}$ is an open compact subgroup whose preimage in $G(\bb{Q}_p)$ contains the corresponding parahoric subgroups of $P_{\Phi_0^\alpha}(\bb{Q}_p)$. Then, by the statement that the quotient of each $\Delta_{\Phi_{\beta'},K'}(ZP_\Phi)$ is finite {\'e}tale proved in \cite[Lem. 4.25]{Wu25}, we are reduced to showing the corresponding statement for Hodge-type Shimura data. This follows from Proposition \ref{prop: functoriality of canonical integral models}(2), noting that the assumption there is trivial by \cite[Cor. 1.9]{Wu25}.
\end{proof}
The crucial part is to show that
\begin{prop}\label{prop-always-comp-in-derived-part}
There always exists a Hodge-type lifting $(G,X)$ such that the assumptions in Proposition \ref{prop-comp-der-consequence} are true. 
\end{prop}
\begin{construction}\label{const-min-del-const-ref}\upshape
For an abelian-type $(G_2,X_2)$, we associate a Hodge-type $(G,X)$ such that it satisfies (1) of Proposition \ref{prop-comp-der-consequence}; this can be done by the construction as in \cite[Prop. 2.3.10]{Del79} and \cite{KPZ24} for any prime $p$. We now modify the Shimura datum to make it to satisfy (2) of Proposition \ref{prop-comp-der-consequence}.\par 
Let $(T,h)\sbst (G,X)$ be a special point with $T$ a maximal torus. We first consider $G^{db}:=G\times_{G^{ab}} G$ where the map $G\to G^{ab}=G/G^\der$ is the natural one.  Denote $T^{ab}:=T/(T\cap G^\der)$. Define a refined construction $G^{rf}:=G\times_{G^{ab}}T$.\par 
Since $(G,X)$ is of Hodge-type, we can decompose $G$ to an almost-direct product $G=G^\der\cdot Z^c\cdot\bb{G}_m$, where $Z^c$ is an $\bb{R}$-anisotropic torus. Then $G^{db}\iso (G^\der\times G^\der)\cdot Z^c\cdot \bb{G}_m\sbst G\times G$.\par
In the second expression of the above line, $Z^c\cdot \bb{G}_m$ maps to $G\times G$ diagonally into $(Z^c\cdot \bb{G}_m)\times (Z^c\cdot \bb{G}_m)$. Alternatively, $G^{db}\iso ((G^\der\times G^\der)\times Z^c\cdot \bb{G}_m)/Z^c\cdot \bb{G}_m\cap G^\der$, where the last intersection maps to the first factor diagonally and to the second factor naturally.\par
Denote the corresponding Shimura data by $(G^{db},X^{db})$ and $(G^{rf},X^{rf})$, respectively. There is a natural map $(G^{rf},X^{rf})\to (G^{db},X^{db})$. Then both of them are of Hodge type because they are both contained in the Hodge-type Shimura datum defined by $G\times_{\eta,\bb{G}_m,\eta}G$.
\end{construction}
\begin{prop}\label{prop-acc-lifting-is-true}
Proposition \ref{prop-always-comp-in-derived-part} is true for any prime $p$. Given $(G_2,X_2)$ of abelian type, the accessible Hodge-type lifting is $(G^{rf},X^{rf})$, as given by Construction \ref{const-min-del-const-ref}.\par
More precisely, $\pi(G^{rf,\der},G^{rf})$ is trivial. In particular, the intersection of any parahoric subgroup $\breve{K}_p$ of $G^{rf}(\bQ)$ with $G^{rf,\der}(\bQ)$ is parahoric.
\end{prop}
\begin{proof}
It suffices to show that $\pi(G^{rf,\der},G^{rf})$ is trivial.
We first consider $\pi(G^{db,\der},G^{db})$. By Lemma \ref{lem: db-gen} below, we have that 
$$ \pi(G^{db,\der},G^{db})=\diag_{\{1,2\}} \{\pi(G^\der,G)\},$$
where the index $\{1,2\}$ labels the two factors of $G^{db,\der}=G^\der\times G^\der$.\par
Note that $G^{rf,\der}=G^\der\times\{1\}\sbst G^{db,\der}$. Then $\pi(G^{\der}\times\{1\},G^{db})\sbst \pi(G^{db,\der},G^{db})$ is trivial. So $\pi(G^{rf,\der},G^{rf})\sbst \pi(G^{rf,\der},G^{db})$ is trivial.
\end{proof}
\begin{lem}\label{lem: db-gen}
    Let $G$ be any reductive group, $\pi: G \to G^{ab}$, $G^{db} = G \times_{G^{ab}} G$. Then $$\pi(G^{db, \der}, G^{db}) = \diag_{\{1, 2\}}\pi(G^{\der}, G) \subset \pi(G^{\der} \times G^{\der}, G \times G).$$
\end{lem}
\begin{proof}
    Consider the commutative diagram of exact sequences:
    \[
\begin{tikzcd}
	& {G^{\der} \times G^{\der}} & {G^{\der} \times G^{\der}} && \\
	0 & {G^{db}} & {G \times G} & {G^{ab}} & 0 \\
	0 & {G^{ab}} & {G^{ab} \times G^{ab}} & {G^{ab}} & {0.}
	\arrow[equals, from=1-2, to=1-3]
	\arrow[hook, from=1-2, to=2-2]
	\arrow[hook, from=1-3, to=2-3]
	\arrow[from=2-1, to=2-2]
	\arrow["i", from=2-2, to=2-3]
	\arrow["\pi", two heads, from=2-2, to=3-2]
	\arrow["{(\pi, \pi^{-1})}", from=2-3, to=2-4]
	\arrow[two heads, from=2-3, to=3-3]
	\arrow[from=2-4, to=2-5]
	\arrow[equals, from=2-4, to=3-4]
	\arrow[from=3-1, to=3-2]
	\arrow["\Delta", from=3-2, to=3-3]
	\arrow["{(\identity, \identity^{-1})}", from=3-3, to=3-4]
	\arrow[from=3-4, to=3-5]
\end{tikzcd}
    \]
    Take the long exact sequences of it,
    \[
\begin{tikzcd}
	{H_1(I, \pi_1(G^{ab})_I)} & {H_1(I, \pi_1(G^{ab})_I) \times H_1(I, \pi_1(G^{ab})_I)} \\
	{\pi_1(G^{\der})_I \times \pi_1(G^{\der})_I} & {\pi_1(G^{\der})_I \times \pi_1(G^{\der})_I} \\
	{\pi_1(G^{db})_I} & {\pi_1(G)_I \times \pi_1(G)_I}
	\arrow["\Delta", hook, from=1-1, to=1-2]
	\arrow["{\Delta(\delta)}", from=1-1, to=2-1]
	\arrow["{(\delta, \delta)}", from=1-2, to=2-2]
	\arrow["{=}", no head, from=2-1, to=2-2]
	\arrow[from=2-1, to=3-1]
	\arrow[from=2-2, to=3-2]
	\arrow["i", from=3-1, to=3-2]
\end{tikzcd}
    \]
    From the diagram, we see that $\pi(G^{db, \der}, G^{db}) = \Image \Delta(\delta)$,  and $\pi(G^{\der} \times G^{\der}, G \times G) = \Image \delta \times \delta$. It follows from definition that $\Image \Delta(\delta) \subset \Image \delta \times \delta$ is exactly $\diag_{\{1, 2\}}\Image \delta = \diag_{\{1, 2\}}\pi(G^{\der}, G)$.
\end{proof}
\begin{proofof}[Proposition \ref{prop-always-comp-in-derived-part}]
Now let $(G, X)$ be any Hodge lifting of $(G_2, X_2)$. Since $G^{db} \subset G \times_{\eta, \Gm, \eta} G$, where $\eta: G \to \Gm$ is the similitude character (here $\eta$ factors through $G^{ab}$), $(G^{db}, X^{db})$ is a Hodge-type Shimura datum. Therefore, $(G^{rf}, X^{rf})$ is a Hodge-lifting of $(G_2, X_2)$ and satisfies the second condition of Proposition \ref{prop-comp-der-consequence}.
\end{proofof}
\begin{rk}\label{rk-mistake}
We remark that, in our construction, $G=G^{rf}$ might not be $R$-smooth when $p=2$, and that $Z_G$ might not be connected. Some supplementary examples are provided in a note \cite{MW26note}.
\end{rk}
\subsubsection{Main theorem for boundary mixed Shimura varieties of abelian type}\label{subsubsec-ab-main-result}
\begin{thm}\label{thm-ext-cim-ab}
Let $(G_2,X_2)$ be an abelian-type Shimura datum, and let $\GG_2$ be any quasi-parahoric model of $G_{2,\bb{Q}_p}$. Then the construction above gives a family of canonical integral models $$\lrbracket{\Shum{K_{\Phi_2}}(P_{\Phi_2}, D_{\Phi_2})}_{K_{\Phi_2}^p}$$ that is adapted with $P_{\Phi_2} \to G_{\Phi_2, h}$ and $G_2\to G_2^c$. Here we can take $\{K_\Phi^p\}$ to be the collection of all neat open compact subgroups in $G_2(\Ap)$. 
\end{thm}
\begin{proof}
As the first step, we show that $\{\mathscr{S}_{\K_\Phi^{\mrm{stb}}}:=\mathscr{S}_{\K_\Phi^{\mrm{stb}}}(ZP_\Phi,ZP_\Phi(\bb{Q})D^b_\Phi)\}_{K^{\prime,p}}$ is a canonical integral model adapted with $ZP_\Phi\to ZG_{\Phi, h}:=ZP_\Phi/W_\Phi$ and $ZP_\Phi\to ZP_\Phi^c$ in the sense of Definition \ref{def: well-adapted, 2} and Definition \ref{def: well-adapted}, respectively.\par 
Note that there is a morphism $\pi^c_{\K_\Phi^{\mrm{stb}}}:\mathscr{S}_{\K_\Phi^{\mrm{stb}}}\to \mathscr{S}_{\K_\Phi^{\mrm{stb},c}}(ZP_\Phi^c,(ZP_\Phi(\bb{Q})D^b_\Phi)^c)$ by construction. Indeed, in Construction \ref{const-scheme-from-hodge-to-ab}, we replace $G'$ in $(G',X_b)$ with $G^{\prime,c}$. Since $G'=G_2\times^{G^\der_2}(G/C^\der)$ and $G=G^c$, we see that $\pi^{b,c}:(G,X)\to (G^{\prime,c},X_b^c)$ also satisfies the setup in Construction \ref{const-step-1-gp}. 
The morphism $\pi^c_{\K_\Phi^{\mrm{stb}}}$ is finite {\'e}tale since we can write the map as
$$\disju_{\beta\in I_{\Phi,K^{\prime,\mrm{stb}}}} \mathscr{S}_{K_{\Phi_\beta}}/\Delta_{\Phi_{\beta},K^{\prime,\mrm{stb}}}(ZP_\Phi)\to \disju_{\beta^c\in I_{\Phi^c,K^{\prime,\mrm{stb},c}}} \mathscr{S}_{K_{\Phi_{\beta^c}}}/\Delta_{\Phi_{\beta^c},K^{\prime,\mrm{stb},c}}(ZP_\Phi^c),$$
and for each $\beta$ mapping to $\beta'$, there is a commutative diagram 
\begin{equation*}
    \begin{tikzcd}
\mathscr{S}_{K_{\Phi_\beta}}\arrow[r]\arrow[d]&\mathscr{S}_{K_{\Phi_{\beta^c}}}\arrow[d]\\
\mathscr{S}_{K_{\Phi_\beta}}/\Delta_{\Phi_{\beta},K^{\prime,\mrm{stb}}}(ZP_\Phi)\arrow[r]&\mathscr{S}_{K_{\Phi_{\beta^c}}}/\Delta_{\Phi_{\beta^c},K^{\prime,\mrm{stb},c}}(ZP_\Phi^c),
    \end{tikzcd}
\end{equation*}
such that all but the bottom arrow are already known to be finite {\'e}tale.

We check that the models $\{\mathscr{S}_{\K_\Phi^{\mrm{stb}}}\}_{K^{\prime,p}}$ are Pappas-Rapoport canonical integral models in the sense of Axiom \ref{def: canonical model for mixed Shimura data}. For (2), it follows from the construction and \cite[Lem. 4.25]{Wu25}. For (1), we choose any discrete valuation ring $R$ with fraction field $F$. 
Choose any prime $l\neq p$. Choose a $F$-point $x$ of $\mathscr{S}_{\K_\Phi^{\mrm{stb},l}}$.
Then it follows from Construction \ref{const-scheme-from-hodge-to-ab} (2) (and (3)) and the proof of \cite[Prop. 4.30]{Wu25} that the action of $\varprojlim_{K_l'}\Delta_{\Phi_\beta,K^{\prime,\mrm{stb}}}(ZP_\Phi)$ on $\varprojlim_{K_l}\mathscr{S}_{K_{\Phi_\beta}}$ factors through a \emph{finite} group. Thus, by Lemma \ref{lem: enhanced (1)}, there is a finite field extension $F'$ of $F$, such that $x$ lifts to an $F'$-point $\wdtd{x}$ of some $\mathscr{S}_{K_{\Phi_\beta}^l}$.\par 
Denote by $R'$ the normalization of $R$ in $F'$; this is still a \emph{discrete} valuation ring. Then we can apply Lemma \ref{lem: enhanced (1)} in the Hodge-type case to get a unique extension of $\wdtd{x}$ to an $R'$-point $\wdtd{y}$ of $\mathscr{S}_{K_{\Phi_\beta}^l}$. This point maps to an $R'$-point $y$ of $\mathscr{S}_{\K_\Phi^{\mrm{stb},l}}$. The uniqueness of this lifting follows from separatedness. Now (a stronger version of) (1) follows because $R'\cap F=R$.\par
We now check (3) and (4). Denote by $(\mathscr{P}_\beta,\phi_\beta)$ the $\ca{P}_{\Phi_\beta}$-shtuka on $\mathscr{S}_{K_{\Phi_\beta}}^{\diamond}$ given by Proposition \ref{prop: hodge type satisfies axioms}. 
We push out $\mathscr{P}_\beta$ via $\ca{P}_{\Phi_\beta}\to g_{\Phi_\beta}\ca{ZP}_\Phi^cg_{\Phi_\beta}^{-1}\xrightarrow{\sim}\ca{ZP}_\Phi^c$; we denote by $(\mathscr{Q}_\beta,\psi_\beta)$ the pushout. 
By checking over the generic fiber, the $\Delta_{\Phi_\beta}$-action lifts to an action on $(\mathscr{Q}_\beta,\psi_\beta)$. By {\'e}tale descent, there is a $\ca{ZP}_\Phi^c$-shtuka $(\wdtd{\mathscr{P}},\wdtd{\phi})$ on $\mathscr{S}_{\K_\Phi^{\mrm{stb}}}^\diamond$.
This checks (3). For (4), it follows again from corresponding (4) in the Hodge-type case and \cite[Lem. 4.25]{Wu25} and from the d{\'e}vissage of integral local Shimura varieties by \cite[Prop. 5.3.1]{pappas2022integral}. \par
We can also verify Axiom \ref{def: canonical model for mixed Shimura data} for $\{\mathscr{S}_{\K^{\mrm{stb}}_{\Phi,h}}(ZP_\Phi,ZP_\Phi(\bb{Q})D^b_\Phi)\}$, $\{\mathscr{S}_{\K_\Phi^{\mrm{stb},c}}(ZP_\Phi^c,(ZP_\Phi(\bb{Q})D^b_\Phi)^c)\}$ and $\{\mathscr{S}_{\K_{\Phi,h}^{\mrm{stb},c}}(ZP_\Phi^c,(ZP_\Phi(\bb{Q})D^b_\Phi)^c)\}$ in exactly the same way. We now have finished the first step.\par
For the second step, we show that above properties for integral models associated with $(G', X_b, \K_\Phi^{\mrm{stb}})$ can be transferred to the integral models associated with $(G', X_a, \K_\Phi^{\mrm{stb}})$: All geometric properties required here can be checked over an {\'e}tale cover. We base change to $\ca{O}^{ur}$ and descend the schemes to $\mathscr{S}_{\K_\Phi^{\mrm{stb}}}(ZP_\Phi,ZP_\Phi(\bb{Q})D^a_\Phi)\to \mathscr{S}_{\K_{\Phi,h}^{\mrm{stb}}}(ZP_\Phi,ZP_\Phi(\bb{Q})D^a_\Phi)$ and $\mathscr{S}_{\K_\Phi^{\mrm{stb}}}(ZP_\Phi,ZP_\Phi(\bb{Q})D^a_\Phi)\to \mathscr{S}_{\K_\Phi^{\mrm{stb}}}(ZP_\Phi^c,(ZP_\Phi(\bb{Q})D^a_\Phi)^c)$. The latter map is again finite {\'e}tale by {\'e}tale descent. Note that the Weil descent datum of shtukas extends from the generic fiber by \cite[Cor. 2.7.10]{PR24}.\par 
For the third step, we show that $\{\mathscr{S}_{K_{\Phi_2}^{\mrm{stb}}}\}$, $\{\mathscr{S}_{K_{\Phi_2, h}^{\mrm{stb}}}\}$, $\{\mathscr{S}_{K_{\Phi_2}^{\mrm{stb}, c}}\}$, $\{\mathscr{S}_{K_{\Phi_2, h}^{\mrm{stb}, c}}\}$ satisfy Axiom \ref{def: canonical model for mixed Shimura data} and has the desired properties. We apply Lemma \ref{lem: two cases}(1) to the embedding of mixed Shimura data $(P_{\Phi_2},D_{\Phi_2}) = (P_{\Phi}, D_{\Phi}) \to (ZP_\Phi,ZP_\Phi(\bb{Q})D_{\Phi}^a)$, to verify axiom $(3)$ (and also axiom $(1)$, $(2)$, as in the proof of Proposition \ref{prop: functoriality of canonical integral models} (1)). Note that, since the whole tower $\{\mathscr{S}_{\K_\Phi^{\mrm{stb}}}(ZP_\Phi,ZP_\Phi(\bb{Q})D^a_\Phi)\}_{K^{\prime,p}}$ is constructed, we can adjust prime-to-$p$ levels so that the cited lemma can be applied. In fact, we adjust $K^{\prime,p}$ so that both $\mathscr{S}_{K_{\Phi_2}^{\mrm{stb}}}\to\mathscr{S}_{\K_\Phi^{\mrm{stb}}} (ZP_\Phi,ZP_\Phi(\bb{Q})D^a_\Phi)$ and $\mathscr{S}_{K_{\Phi_2}^{\mrm{stb},c}}\to\mathscr{S}_{\K_\Phi^{\mrm{stb},c}} (ZP_\Phi^c,(ZP_\Phi(\bb{Q})D^a_\Phi)^c)$ are open and closed embeddings. We see from this that $\mathscr{S}_{K_{\Phi_2}^{\mrm{stb}}}\to\mathscr{S}_{K_{\Phi_2}^{\mrm{stb},c}}$ is finite {\'e}tale for any $K_{\Phi_2}^{\mrm{stb},p}$. In particular, axiom $(4)$ follows easily, using the \emph{d{\'e}vissage} of integral local Shimura varieties (\cite[Prop. 5.3.1]{pappas2022integral}).\par
Now, we consider a general quasi-parahoric level $K_{2,p}$ and finish by applying Proposition \ref{prop: functoriality of canonical integral models}(2) to $K_{2,p}\to K_{2,p}^{\mrm{stb}}$. We just need to verify the condition highlighted there. We adopt the conventions in Construction \ref{const-scheme-from-hodge-to-ab}. In fact, $\mathscr{S}_{K_{\Phi_2}^{c}}\to\mathscr{S}_{K_{\Phi_2}^{\mrm{stb},c}}$ is finite {\'e}tale by applying Proposition \ref{prop: functoriality of canonical integral models}(2), whose condition is trivial for this case. Now, it suffices to show that $\mathscr{S}_{K_{\Phi_2}}\to \mathscr{S}_{K_{\Phi_2}^{\mrm{stb}}}$, 
$\mathscr{S}_{\overline{K}_{\Phi_2}}\to \mathscr{S}_{\overline{K}^{\mrm{stb}}_{\Phi_2}}$
and 
$\mathscr{S}_{K_{\Phi_2,h}}\to \mathscr{S}_{K_{\Phi_2,h}^{\mrm{stb}}}$
are finite {\'e}tale. This is done by Proposition \ref{prop-always-comp-in-derived-part} and Proposition \ref{prop-comp-der-consequence}.
\end{proof}

\section{Well-positioned subschemes}\label{sec-well-position}

   Fix a Shimura datum $(G, X)$ and a quasi-parahoric model $\GG$ of $G_{\rQ_p}$. For each $\Phi \in \ca{CLR}(G, X)$, we denote $G_{\Phi} = gG_{\rQ_p}g^{-1}$, where $g = g_{\Phi, p} \in G(\rQ_p)$. In this subsection, we generalize the results in \cite{boxer2015torsion}, \cite{lan2018compactifications} and \cite{Mao25}. In order to present the nature of those well-positioned subschemes, \textbf{we work under Assumption \ref{ass-well-position}}.

   Let us recall the definition of a well-positioned subscheme. 
    \begin{Definition and Proposition}[{\cite[Def. 2.2.1 and Lem. 2.2.2]{lan2018compactifications}}]\label{def: well-positioned}
    	Let $T$ be a locally noetherian scheme over $\OO_{E}$. A locally closed subset (resp. subscheme) $Y$ of $\Shum{K, T}:=(\Shum{K})_T$ is called \emph{well positioned}, if, for each $\Phi \in \ca{CLR}(G, X)$, there exists a locally closed subset (resp. subscheme) $Y^{\natural}(\Phi) \subset \Zb(\Phi)_T \rightiso \mathcal{Z}([\Phi])_T$ such that for some (thus for all) cone decompositions $\Sigma$, and for each $\sigma\in\Sigma(\Phi)^+$, there are some (thus for all) open affine coverings $\mathfrak{W}$ of $\mathfrak{X}^{\circ}_{\sigma}$ satisfying the following property: the pullback of $Y^{\natural}(\Phi)$ along 
        \[W^0_T \to \Delta_{\Phi, K}^{\circ}\backslash\Xi(\Phi)_T \to \Delta_{\Phi, K}^{\circ}\backslash C(\Phi)_T \to \Delta_{\Phi, K}^{\circ}\backslash\Zb^{\bigsur}(\Phi)_T \to \Zb(\Phi)_T\] coincides with the pullback of $Y$ along $W^0_T \to \Shum{K, T}$ as a subset (resp. subscheme). If this is the case, we say $Y$ is well positioned with respect to $Y^{\natural}:=\{Y^{\natural}(\Phi)\}_{\Phi}$, and $Y^{\natural}$ is associated with $Y$.
    \end{Definition and Proposition}
     In practice, we usually take $T = s = \Spec k_E$ or $T = \bar{s} = \Spec \Bar{k}_E$.
     \begin{rk}
         Compared with \cite[Def. 2.2.1 and Lem. 2.2.2]{lan2018compactifications}, we need to consider the $\Delta_{\Phi, K}^{\circ}$-action. Note that $\Delta_{\Phi, K}^{\circ}$ is trivial in the Hodge-type case, and thus was not considered in \cite{lan2018compactifications}. With this modification, if we replace everything
         \begin{equation}\label{eq: boundary schemes}
             \Zb^{\bigsur}(\Phi),\ C(\Phi),\ \Xi(\Phi),\ \Xi(\Phi)(\sigma),\ \Xi(\Phi)_{\sigma},\ \ovl{\Xi}_{\Sigma},\ \mathfrak{X}_{[\Phi, \sigma]},\ \mathfrak{X}_{[\Phi, \sigma]}^{\circ},\ \mathfrak{X}_{\Sigma}
         \end{equation}
         by their quotients under the $\Delta_{\Phi, K}^{\circ}$-action, then the statements and proofs in \cite{lan2018compactifications} and \cite{Mao25} work without change.
     \end{rk}
     \begin{rk}
         The quotient map $\Zb^{\bigsur}(\Phi) \to \Zb(\Phi)$ is finite \'etale by \cite[Cor. 4.27]{Wu25}. It follows that $\Delta_{\Phi, K}^{\circ}\backslash\Zb^{\bigsur}(\Phi) \to \Zb(\Phi)$ is also finite \'etale. Under Assumption \ref{ass-well-position}, $\Delta_{\Phi, K}^{\circ}\backslash C(\Phi) \to \Zb(\Phi)$ is flat with geometrically reduced fibers, which verifies \cite[Assumption 2.7]{Mao25}. In this situation, well-positioned subschemes enjoy more satisfactory properties.
     \end{rk}
     \begin{definition}[{\cite[Def. 2.3.1]{lan2018compactifications}}]
    	Let $Y$ be a locally closed subscheme of $\Shum{K, T}$. Let $\overline{Y}$ be the closure of $Y$ in $\Shum{K, T}$, $Y_0=\overline{Y}\backslash Y$ the complement of $Y$ in $\overline{Y}$. Let $\overline{Y}^{\min}$ and $Y_0^{\min}$ be the closure of $\overline{Y}$ and $Y_0$ in $\Shumm{K, T}$ respectively, $Y^{\min}$ be the complement of $Y_0^{\min}$ in $\overline{Y}^{\min}$. Similarly we define $\overline{Y}^{\Sigma}$, $Y_{0}^{\Sigma}$, $Y^{\Sigma}$. We call $Y^{\min}$ (resp. $Y^{\Sigma}$) the \emph{partial minimal compactification} (resp. the \emph{partial toroidal compactification} with respect to $\Sigma$) of $Y$.
    \end{definition}

    If $Y \subset \Shum{K, T}$ is well positioned with respect to $Y^{\natural}:=\lrbracket{Y^{\natural}(\Phi)}_{\Phi}$, then $\oint_{K, T}^{\Sigma}: \Shum{K, T}^{\Sigma} \to \Shum{K, T}^{\min}$ induces a morphism $\oint_{Y}^{\Sigma}: Y^{\Sigma} \to Y^{\min}$. Moreover, $(\oint_{K, T}^{\Sigma})^{-1}(Y^{\min}) = Y^{\Sigma}$, and the compactifications $Y^{\Sigma} \to Y^{\min}$ of $Y$ satisfy Axiom \ref{axiom-good-compactification} except for the flatness and normality of $Y^{\Sigma}$, of $Y^{\min}$, and of their boundary strata. See \cite[Thm. 2.3.2]{lan2018compactifications}. Note that, in Axiom \ref{axiom-good-compactification}, for any scheme (or formal scheme) in (\ref{eq: boundary schemes}), we replace $(?)$ by $Y_{(?)}^{\natural}$, which is defined as the pullback of $Y^{\natural}(\Phi)$ along $(?) \to \Zb(\Phi)$. The identification $\Zb(\Phi) \rightiso \mathcal{Z}([\Phi])$ induces a canonical morphism $Y^{\natural}(\Phi) \to Y_{\mathcal{Z}([\Phi])}:=(Y^{\min} \cap \mathcal{Z}([\Phi]))$, which induces a bijection on the underlying sets.

    \subsection{Newton strata and central leaves}\label{subsec: Newton strata and central leaves}

   We define Newton strata and central leaves following \cite{PR24}. Let $x \in \Shum{K}(G, X)(k)$. Pulling back the universal $\GG^{c}$-shtuka $(\PPs, \phi_{\PPs})$ over $\Shum{K}(G, X)^{\Dia/}$, we get a $\GG^{c}$-shtuka $(\PPs_x, \phi_{\PPs_x}):=x^*(\PPs, \phi_{\PPs})$ over $\Spec (k)$, which is associated with a $\GG^{c}$-torsor on $\Spec W(k)$ together with an isomorphism $\phi_{\PPs_x}: \phi^*(\PPs_x)[1/p] \rightiso \PPs_x[1/p]$. A choice of trivialization of $\PPs_x$ defines an element $b_x \in G^c(\bQ)$; a different choice of trivialization gives $\GG^{c}(\bZ_p)$-$\sigma$-conjugation of $b_x$. Let $C(\GG^c)=G^c(\bQ)/\GG^{c}(\bZ_p)_{\sigma}$ and $B(G^c)=G^c(\bQ)/G^{c}(\bQ)_{\sigma}$; we have
   \begin{equation}\label{eq: Newton, k-points}
       \Upsilon_K(k): \Shum{K}(G, X)(k) \to C(\GG^c),\quad \delta_K(k): \Shum{K}(G, X)(k) \to B(G^c).
   \end{equation}

   Since $(\PPs, \phi_{\PPs})$ is bounded by $\mu^c$, the $\sigma$-conjugation class of $b_x$ sits inside $B(G^c, \lrbracket{\mu^{c, -1}}) \subset B(G^c)$.

   One can define $\Upsilon_K$, $\delta_K$ globally, without referring to $k$-points. In fact, in subsection \ref{subsec: isocrystals}, we recall a natural projection $\Sht_{\GG^c, \mu^c, \delta = 1}^W \to G^c\textit{-}\Isoc_{\mu^{c, -1}}$. By construction, $\Upsilon_K(k)$ and $\delta_K(k)$ in (\ref{eq: Newton, k-points}) are the $k$-points of $\Shum{K}(G, X)_{\bar{s}}^{\perf} \to \Sht^W_{\GG^c, \mu^c, \delta = 1}$ and $\Shum{K}(G, X)_{\bar{s}}^{\perf} \to G^c\textit{-}\Isoc_{\mu^{c, -1}}$, respectively. Note that perfection does not change the underlying topological space, so we have a morphism
   \begin{equation*}
       \delta_K: \Shum{K}(G, X)_{\bar{s}} \to B(G^c, \lrbracket{\mu^{c, -1}})
   \end{equation*}
   whose $k$-points are $\delta_K(k)$. Given $[b] \in B(G^c, \lrbracket{\mu^{c, -1}})$ with $b \in G^c(\bQ)$, let $\NE^{[b]} \subset \Shum{K}(G, X)_{\bar{s}}$ be the preimage $\delta_K^{-1}([b])$; we call it the \emph{Newton stratum} associated with $[b]$. \cite[Thm. 3.6]{rapoport1996classification} says that $\Shum{K}(G, X)_{\bar{s}} \to B(G^c, \lrbracket{\mu^{c, -1}})$ is semi-lower continuous; therefore, Newton strata are locally closed. We endow them with the unique reduced subscheme structures.

   Now we define $\Upsilon_K$ globally. Given $[[b]] \in C(\GG^c)$ with $b \in G^c(\bQ)$, let $\CE^{[[b]]}(k) \subset \Shum{K}(G, X)(k)$ (resp. $\NE^{[b]}(k) \subset \Shum{K}(G, X)(k)$) be the preimage $\Upsilon_K(k)^{-1}([[b]])$ (resp. $\delta_K(k)^{-1}([b])$).

   \begin{lem}\label{lem: central leaves are closed in Newton strata}
      $\CE^{[[b]]}(k) \subset \NE^{[b]}(k)$ is closed.
   \end{lem}
   \begin{proof}
      Consider the diagram
       \begin{equation}\label{eq: auxilary shtuka fiber product}
\begin{tikzcd}
	{\Shum{K^{\circ}}(G, X)^{\Dia/}} & {\Shum{K}(G, X)^{\prime\Dia/}} & {\Shum{K}(G, X)^{\Dia/}} \\
	& {\Sht_{\GG^{\circ, c}, \mu^c, \delta = 1}} & {\Sht_{\GG^{c}, \mu^c, \delta = 1}.}
	\arrow["g", from=1-1, to=1-2]
	\arrow[from=1-1, to=2-2]
	\arrow["f", from=1-2, to=1-3]
	\arrow[from=1-2, to=2-2]
	\arrow["\square"{description}, draw=none, from=1-2, to=2-3]
	\arrow[from=1-3, to=2-3]
	\arrow[from=2-2, to=2-3]
\end{tikzcd}
  \end{equation}
       where $\Shum{K}(G, X)^{\prime\Dia/}$ is the fiber product, and $\Shum{K^{\circ}}(G, X)$ is the relative normalization of $\Shum{K}(G, X)$ in $\shu{K^{\circ}}(G, X)$. Since $\Sht_{\GG^{c, \circ}, \mu^c, \delta = 1} \to \Sht_{\GG^{c}, \mu^c, \delta = 1}$ is a finite \'etale torsor under the abelian group $\pi_0(\GG^c)^{\phi}$, $\Shum{K}(G, X)^{\prime\Dia/}$ is represented by a normal flat scheme $\Shum{K}(G, X)'$, and $\Shum{K}(G, X)' \to \Shum{K}(G, X)$ is a finite \'etale torsor under $\pi_0(\GG^c)^{\phi}$, thanks to \cite[Prop. 2.3.1]{daniels2024conjecture}. We take the special fiber of the diagram (\ref{eq: auxilary shtuka fiber product}), and denote the Newton strata and central leaves in $\Shum{K}(G, X)'(k)$ by $\CE^{[[b]]}(k)$ and $\NE^{[b]}(k)$, respectively, using $\Sht_{\GG^{c, \circ}, \mu^c, \delta = 1}^W$. Since $\NE^{[b]}(k)$ is locally closed in $\Shum{K}(G, X)'(k)$ by \cite{rapoport1996classification}, we claim that $\CE^{[[b]]}(k)$ is closed in $\NE^{[b]}(k)$. If this is true, by using torsors, a central leaf in $\Shum{K}(G, X)(k)$ is the image of some $\CE^{[[b]]}(k)$ under $f$, and its preimage under $f$ is a topologically disjoint union of some $\CE^{[[b_i]]}(k)$; then the lemma follows from the fact that $f$ is finite \'etale.

       Let $x \in \CE^{[[b]]}(k)$. In \cite[\S 2.14]{hamacher2025point}, given a perfect scheme $S$ and a $\GG^c$-shtuka $(\PPs, \phi_{\PPs})$ on $S$, the authors define
       \[ \CE^{[[b]]}_{(\PPs, \phi_{\PPs})}:= \left\{ s \in S |\ \bar{s}^*(\PPs, \phi_{\PPs}) \cong (\PPs_x, \phi_{\PPs_x})|_{\ovl{\kappa(s)}} \right\}, \]
       and similarly define $\NE^{[b]}_{(\PPs, \phi_{\PPs})}$ using isocrystals. \cite[Prop. 2.15(3)]{hamacher2025point} says that $\CE^{[[b]]}_{(\PPs, \phi_{\PPs})} \subset \NE^{[b]}_{(\PPs, \phi_{\PPs})}$ is closed. Now let $S$ be the perfection of $\NE^{[b]}$; then $\NE^{[b]}_{(\PPs, \phi_{\PPs})}$ is the whole space $\NE^{[b], \perf}$. Thus, $\CE^{[[b]]}(k) = \CE^{[[b]]}_{(\PPs, \phi_{\PPs})}(k) \subset \NE^{[b]}(k)$ is closed.
   \end{proof}

  By some standard arguments on Jacobson schemes, (e.g., \cite[\S 3.2]{Mao25b}), there exists a unique closed subscheme $\CE^{[[b]]} \subset \NE^{[b]}$ with induced reduced subscheme structure such that the set of $k$-points of $\CE^{[[b]]}$ is $\CE^{[[b]]}(k)$. Finally, let
  \begin{equation*}
      C(\GG^c, \lrbracket{\mu^{c, -1}}):= \Sht_{\GG^c, \mu^c, \delta = 1}^W(k),
  \end{equation*}
  then we obtain a globally defined $\Upsilon_K$, which is determined by its values on $k$-points $\Upsilon_K(k)$:
\begin{equation*}
      \Upsilon_K: \Shum{K}(G, X)_{\bar{s}} \to C(\GG^c, \lrbracket{\mu^{c, -1}}),
  \end{equation*}
  and we call the fibers \emph{central leaves}. In the proof of Lemma \ref{lem: central leaves are closed in Newton strata}, fix a $k$-point $s \in \CE^{[[b]]}$, we see that
  $$ \CE^{[[b]]} = \left\{ x \in \Shum{K}(G, X)_{\bar{s}} |\ (\PPs_s, \phi_{\PPs_s})|_{\ovl{\kappa(x)}} \cong (\PPs_x, \phi_{\PPs_x})|_{\ovl{\kappa(x)}} \right\}. $$

  \begin{rk}
      The above arguments and definitions also apply to $(G^*_{\Phi, h}, \GG^*_{\Phi, h})$. The central leaves on $\Shum{K_{\Phi, h}}(G_{\Phi, h}, D_{\Phi, h})_{\bar{s}}$, defined using $\Sht_{\GG^c_{\Phi, h}, \mu^c_{\Phi, h}, \delta = 1}^W$, are finer than those defined using $\Sht_{\GG^*_{\Phi, h}, \mu^*_{\Phi, h}, \delta = 1}^W$. Nevertheless, the index set of Newton strata (and of KR strata, of EKOR strata, recalled in the next subsection) depends only on the adjoint group, so there is no difference in whether one uses $\GG^c_{\Phi, h}$ or $\GG^*_{\Phi, h}$.
  \end{rk}
   
   \begin{prop}\label{prop: Newton strata are well-positioned}
       Newton strata are well positioned. Moreover, let $\NE^{[b]}$ be a Newton stratum on $\Shum{K}(G, X)_{\bar{s}}$ with some $[b] \in B(G^c, \lrbracket{\mu^{c, -1}})$. Then, for each $\Phi \in \ca{CLR}(G, X)$, $(\NE^{[b]})_{\Zb^{\bigsur}(\Phi)}^{\natural}$ is either empty or a Newton stratum $\NE^{[b_{\Phi, h}]}$ on $\Zb^{\bigsur}(\Phi)_{\bar{s}} \cong \Shum{K_{\Phi, h}}(G_{\Phi, h}, D_{\Phi, h})_{\bar{s}}$, for some $[b_{\Phi, h}] \in B(G_{\Phi, h}^*, \lrbracket{\mu_{\Phi, h}^{*, -1}})$. The relation between $[b]$ and $[b_{\Phi, h}]$ is given in the proof.
   \end{prop}
   \begin{proof}
       We pass the commutative diagram (\ref{eq: special fiber of main diagram}) to isocrystals:
\[
\begin{tikzcd}
	{\Shum{K}(G, X)^{\perf}_{\bar{s}}} & {W^{0, \perf}_{\bar{s}}} & {\Delta_{\Phi,K}^{\circ}\backslash\Shum{K_{\Phi}}(P_{\Phi}, D_{\Phi})^{\perf}_{\bar{s}}} & {\Delta_{\Phi,K}^{\circ}\backslash\Shum{K_{\Phi, h}}(G_{\Phi, h}, D_{\Phi, h})^{\perf}_{\bar{s}}} \\
	{G^{c}\textit{-}\Isoc_{\mu^{c, -1}}} && {P^*_{\Phi}\textit{-}\Isoc_{\mu_{\Phi}^{*, -1}}} & {G_{\Phi, h}^*\textit{-}\Isoc_{\mu_{\Phi,h}^{*, -1}}.}
	\arrow[from=1-1, to=2-1]
	\arrow[from=1-2, to=1-1]
	\arrow[from=1-2, to=1-3]
	\arrow[from=1-3, to=1-4]
	\arrow[from=1-3, to=2-3]
	\arrow[from=1-4, to=2-4]
	\arrow["{\Int(g_{\Phi}^{-1})}"', from=2-3, to=2-1]
	\arrow[from=2-3, to=2-4]
\end{tikzcd}
\]
       We further take the underlying topological spaces and get the commutative diagram
\[
\begin{tikzcd}
	{\Shum{K}(G, X)_{\bar{s}}} & {W^{0}_{\bar{s}}} & {\Delta_{\Phi,K}^{\circ}\backslash\Shum{K_{\Phi}}(P_{\Phi}, D_{\Phi})_{\bar{s}}} & {\Delta_{\Phi,K}^{\circ}\backslash\Shum{K_{\Phi, h}}(G_{\Phi, h}, D_{\Phi, h})_{\bar{s}}} \\
	{B(G^c, \lrbracket{\mu^{c, -1}})} && {B(P_{\Phi}^*, \lrbracket{\mu_{\Phi}^{*, -1}})} & {B(G_{\Phi, h}^*, \lrbracket{\mu_{\Phi, h}^{*, -1}}).}
	\arrow[from=1-1, to=2-1]
	\arrow[from=1-2, to=1-1]
	\arrow[from=1-2, to=1-3]
	\arrow[from=1-3, to=1-4]
	\arrow[from=1-3, to=2-3]
	\arrow[from=1-4, to=2-4]
	\arrow["{\Int(g_{\Phi}^{-1})}"', from=2-3, to=2-1]
	\arrow[from=2-3, to=2-4]
\end{tikzcd}
\]
      Let $[b] \in B(G^c, \lrbracket{\mu^{c, -1}})$. If $\NE^{[b]} \subset \Shum{K}(G, X)_{\bar{s}}$ is non-empty over $W_{\bar{s}}^0$, then $[b]$ is in the image of $\Int(g_{\Phi}^{-1}): B(P_{\Phi}^*, \lrbracket{\mu_{\Phi}^{*, -1}}) \to B(G^c, \lrbracket{\mu^{c, -1}})$ for some $[b_{\Phi}] \in B(P_{\Phi}^*, \lrbracket{\mu_{\Phi}^{*, -1}})$.
      
      On the other hand, let $[b_{\Phi, h}]\in B(G_{\Phi, h}^*, \lrbracket{\mu_{\Phi, h}^{*, -1}})$ be the image of $[b_{\Phi}]$. Since $B(P_{\Phi}^*, \lrbracket{\mu_{\Phi}^{*, -1}}) \rightiso B(G_{\Phi, h}^*, \lrbracket{\mu_{\Phi, h}^{*, -1}})$ is a bijection, the preimage of $\NE^{[b]}$ over $W_{\bar{s}}^0$ is the finite union of preimages of those Newton strata $\NE^{[b_{\Phi, h}]}$ on $\Shum{K_{\Phi, h}}(G_{\Phi, h}, D_{\Phi, h})_{\bar{s}}$, where $[b_{\Phi, h}]$ are the preimages of $[b] \in B(G^c, \lrbracket{\mu^{c, -1}})$ under $B(G_{\Phi, h}^*, \lrbracket{\mu_{\Phi, h}^{*, -1}}) \rightiso B(P_{\Phi}^*, \lrbracket{\mu_{\Phi}^{*, -1}}) \to B(G^c, \lrbracket{\mu^{c, -1}})$.

      We further claim that $B(P_{\Phi}^*, \lrbracket{\mu_{\Phi}^{*, -1}}) \to B(G^c, \lrbracket{\mu^{c, -1}})$ is injective; then there is a unique $[b_{\Phi, h}]$ associated with $[b]$. By \cite[Lem. 2.1]{he2024affine}, we have injectivity:
      \[ B(L_{\Phi}^*, \lrbracket{\mu_{\Phi}^{*, -1}}) = B(Q_{\Phi}^*, \lrbracket{\mu_{\Phi}^{*, -1}}) \to B(G_{\Phi}^c, \lrbracket{\mu_{\Phi}^{c, -1}}) = B(G^c, \lrbracket{\mu^{c, -1}}). \]
      Also, note that $B(P_{\Phi}^*, \lrbracket{\mu_{\Phi}^{*, -1}}) \to B(Q_{\Phi}^*, \lrbracket{\mu_{\Phi}^{*, -1}})$ is bijective:
      \[ B(P_{\Phi}^*, \lrbracket{\mu_{\Phi}^{*, -1}}) = B(G_{\Phi, h}^*, \lrbracket{\mu_{\Phi, h}^{*, -1}}) = B(G_{\Phi, h}^{\ad}, \lrbracket{\mu_{\Phi, h}^{\ad, -1}}), \]
      \[ B(Q_{\Phi}^*, \lrbracket{\mu_{\Phi}^{*, -1}}) = B(L_{\Phi}^*, \lrbracket{\mu_{\Phi, h}^{*, -1}}) = B(L_{\Phi}^{\ad}, \lrbracket{\mu_{\Phi, h}^{\ad, -1}}), \]
      and 
      \[ B(L_{\Phi}^{\ad}, \lrbracket{\mu_{\Phi, h}^{\ad, -1}}) = B(G_{\Phi, h}^{\ad}, \lrbracket{\mu_{\Phi, h}^{\ad, -1}}) \times B(G_{\Phi, l}^{\ad}, \lrbracket{\identity}) = B(G_{\Phi, h}^{\ad}, \lrbracket{\mu_{\Phi, h}^{\ad, -1}}). \]
      
      Finally, to show that Newton strata are well positioned, we pass to $\Zb(\Phi) = \Delta_{\Phi,K}\backslash \Zb^{\bigsur}(\Phi)$. It suffices to show that $\Delta_{\Phi,K}$ stabilizes each Newton stratum; then there is a unique subscheme $(\NE^{[b]})^{\natural} \subset \Zb(\Phi)_{\bar{s}}$ whose preimage in $\Zb^{\bigsur}(\Phi)_{\bar{s}}$ is $(\NE^{[b]})^{\natural}_{\Zb^{\bigsur}(\Phi)} = \NE^{[b_{\Phi, h}]}$. Such $(\NE^{[b]})^{\natural} \subset \Zb(\Phi)_{\bar{s}}$ is automatically locally closed if it exists; see \cite[Lem. 2.3.10]{lan2018compactifications}. Consider the $\Delta_{\Phi,K}$-action on Newton strata; now we apply Proposition \ref{prop: boundary shtukas extend over integral base}.
   \end{proof}
   \begin{rk}
       As explained in \cite[Lem. 2.16, 2.20, Prop. 2.30]{Mao25}, the connected components and closures of Newton strata, as well as $\NE^{\leq [b]}$ for any $[b] \in B(G^c, \lrbracket{\mu^{c, -1}})$, are all well-positioned. The same applies to central leaves, KR strata, and EKOR strata proved in later sections, and we do not repeat the arguments there.
   \end{rk}

   The following proof was explained to us by Sian Nie.
   \begin{lem}\label{lem: map of basic elements}
       Let $G$ be any reductive group over $\rQ_p$, $M$ be any proper Levi subgroup of $G$, $\mu$ be a cocharacter of $M$ that is non-central in $G$. The image of $B(M, \lrbracket{\mu}) \to B(G, \lrbracket{\mu})$ does not contain the basic element of $B(G, \lrbracket{\mu})$.
   \end{lem}
   \begin{proof}
       Given $[b_1] \leq [b_2]$ in $B(M, \lrbracket{\mu})$, then $[b_1] \leq [b_2]$ in $B(G, \lrbracket{\mu})$. Since basic elements are minimal, it suffices to show that given the basic element $[b] \in B(M, \lrbracket{\mu})$, its image in $B(G, \lrbracket{\mu})$ can not be basic.
       
       We fix a pinning and assume $\mu$ is dominant in $G$. Let $\mu^{\dia}$ be the $\sigma$-average of $\mu$, recall that 
       \[ B(M, \lrbracket{\mu}) = \left\{ [b] \in B(M) |\ \kappa([b]) \leq \kappa(\mu), \nu_{b} \leq \mu^{\dia}  \right\}  \]
       Let $[b] \in B(M, \lrbracket{\mu})$ be basic, then $\nu_b = \mu^{\dia} - h$ is central in $M$, where $h$ is a non-negative linear combination of positive coroots of $M$. Recall that $[b]$ is basic in $B(G, \lrbracket{\mu})$ if and only if $\nu_b$ is central in $G$ if and only if $\lrangle{\alpha, \nu_b} = 0$ for all simple roots in $G$. For those simple roots $\alpha$ not in $M$, $\lrangle{\alpha, \mu^{\dia}} \geq 0$ since $\mu^{\dia}$ is dominant in $G$, and $\lrangle{\alpha, h} \leq 0$ by the Cartan matrix, then 
       \[ \lrangle{\alpha, \nu_b} = \lrangle{\alpha, \mu^{\dia}} - \lrangle{\alpha, h} \geq 0, \]
       it takes equality if and only if $\lrangle{\alpha, \mu^{\dia}} = \lrangle{\alpha, h} = 0$. This implies that $h = 0$. Since $\nu_b = \mu^{\dia}$ is central in $M$, then $\mu^{\dia}$ is central in $G$, we get a contradiction.
   \end{proof}
   
   \begin{cor}\label{cor: basic newton stratum has no boundary}
       Basic Newton stratum has no boundary.
   \end{cor}
   \begin{proof}
       As in the proof of Proposition \ref{prop: Newton strata are well-positioned}, it suffices to show that the basic element $[b_0] \in B(G^c, \lrbracket{\mu^{c, -1}})$ does not come from any Levi subgroup; this is proved in Lemma \ref{lem: map of basic elements}. Note that the Hodge cocharacter $\mu$ is central in $G$ if and only if $G = T$ is a torus, and the Shimura variety associated with $(T, \lrbracket{h})$ has no boundary.
   \end{proof}

   \begin{definition}
       Let $Y$ be a scheme and $\lrbracket{Y_i}_{i \in I}$ be subschemes of $Y$. We say that $Y$ is a \textbf{topologically disjoint union} of $\lrbracket{Y_i}_{i \in I}$ if $Y = \bigsqcup_{i \in I} Y_i$ and all $Y_i \subset Y$ are open and closed.
   \end{definition}

   \begin{prop}\label{prop: central leaves are well-positioned}
       Central leaves are well positioned. Moreover, let $\CE^{[[b]]}$ be a central leaf on $\Shum{K}(G, X)_{\bar{s}}$ with some $[[b]] \in C(\GG^c, \lrbracket{\mu^{c, -1}})$. Then, for each $\Phi \in \ca{CLR}(G, X)$, $(\CE^{[[b]]})_{\Zb^{\bigsur}(\Phi)}^{\natural}$ is either empty or a topologically disjoint union of central leaves $\CE^{[[b_{\Phi, h}]]}$ on $\Zb^{\bigsur}(\Phi)_{\bar{s}} \cong \Shum{K_{\Phi, h}}(G_{\Phi, h}, D_{\Phi, h})_{\bar{s}}$, for some collection of $[[b_{\Phi, h}]] \in C(\GG_{\Phi, h}^*, \lrbracket{\mu_{\Phi, h}^{*, -1}})$. The relation between $[[b]]$ and the collection of $[[b_{\Phi, h}]]$ is given in the proof.
   \end{prop}
   \begin{proof}
       Let $[[b]] \in C(\GG^c, \lrbracket{\mu^{c, -1}})$. If $\CE^{[[b]]} \subset \Shum{K}(G, X)_{\bar{s}}$ is non-empty over $W_{\bar{s}}^0$, then $[[b]]$ is in the image of some $[[b_{\Phi}]]$ under $\Int(g_{\Phi}^{-1}): C(\PP_{\Phi}^*, \lrbracket{\mu_{\Phi}^{*, -1}}) \to C(\GG^c, \lrbracket{\mu^{c, -1}})$; here $C(\PP_{\Phi}^*, \lrbracket{\mu_{\Phi}^{*, -1}}) := \Sht_{\PP_{\Phi}^*, \mu_{\Phi}^*}^W(k)$. We let $[[b_{\Phi, h}]] \in C(\GG_{\Phi, h}^*, \lrbracket{\mu_{\Phi, h}^{*, -1}})$ be the image of $[[b_{\Phi}]]$. 
       
       We claim that there is at most one $[[b_{\Phi}]] \in C(\PP_{\Phi}^*, \lrbracket{\mu_{\Phi}^{*, -1}})$ in the image of $\Shum{K_{\Phi}}(k) \to C(\PP_{\Phi}^*, \lrbracket{\mu_{\Phi}^{*, -1}})$ that has fixed image $[[b_{\Phi, h}]] \in C(\GG_{\Phi, h}^*, \lrbracket{\mu_{\Phi, h}^{*, -1}})$. If this is the case, the preimage of $\CE^{[[b]]}$ over $W_{\bar{s}}^0$ is the finite union of preimages of those central leaves $\CE^{[[b_{\Phi, h}]]}$ on $\Shum{K_{\Phi, h}, \bar{s}}$, where $[[b_{\Phi, h}]]$ are the preimages of $[[b]] \in C(\GG^c, \lrbracket{\mu^{c, -1}})$ under 
       \[ C(\GG_{\Phi, h}^*, \lrbracket{\mu_{\Phi, h}^{*, -1}}) \leftarrow C(\PP_{\Phi}^*, \lrbracket{\mu_{\Phi}^{*, -1}}) \stackrel{\Int(g^{-1}_{\Phi})}{\to}  C(\GG^c, \lrbracket{\mu^{c, -1}}). \]

       Let us show the claim. We first work with abelian scheme torsor part. Let $x \in \CE^{[[b_{\Phi, h}]]}(k) \subset \Shum{K_{\Phi, h}}(k)$, $y_1, y_2 \in \Shum{\ovl{K}_{\Phi}}(k)$ be in the preimage of $x$. Since $\Shum{\ovl{K}_{\Phi}} \to \Shum{K_{\Phi, h}}$ is a torsor under the abelian scheme $\ab_{K}(\Phi) \to \Shum{K_{\Phi, h}}$, there exists $\gamma \in \ab_{K}(\Phi)(k)$ with image $x$ such that $\gamma y_1 = y_2$. We could lift $y_1$, $y_2$, $x$, $\gamma$ to some $\OO_F$-points $\Tilde{y}_1$, $\Tilde{y}_2$, $\Tilde{x}$, $\Tilde{\gamma}$ respectively, where $F$ is a finite extension of $\bQ$, such that $\Tilde{\gamma}\Tilde{y}_1 = \Tilde{y}_2$ with image $\tilde{x}$. By Lemma \ref{lem: abelian scheme action on shtuka, generic fiber} and \cite[Corollary 2.7.10]{PR24}, there exists an isomorphism of $\ovl{\PP}_{\Phi}^*$-shtukas $(\ovl{\PPs}_{\Phi}, \phi_{\ovl{\PPs}_{\Phi}})$ and $\tilde{\gamma}^*(\ovl{\PPs}_{\Phi}, \phi_{\ovl{\PPs}_{\Phi}})$ over $\Shum{\ovl{K}_{\Phi}, \tilde{x}}$, thus there exists an isomorphism of $\ovl{\PP}_{\Phi}^*$-shtukas $y_1^*(\ovl{\PPs}_{\Phi}, \phi_{\ovl{\PPs}_{\Phi}})$ and $y_2^*(\ovl{\PPs}_{\Phi}, \phi_{\ovl{\PPs}_{\Phi}})$. Next, we show that given different $x, x' \in \CE^{[[b_{\Phi, h}]]}(k)$, one can take some points in the fiber $y, y' \in \Shum{\ovl{K}_{\Phi}}(k)$ respectively such that $y^*(\ovl{\PPs}_{\Phi}, \phi_{\ovl{\PPs}_{\Phi}})$ is isomorphic to $y^{\prime *}(\ovl{\PPs}_{\Phi}, \phi_{\ovl{\PPs}_{\Phi}})$. \'Etale locally $T \to \Shum{K_{\Phi, h}}$, there is a section $T \to \Shum{\ovl{K}_{\Phi}, T} \cong \ab_K(\Phi)_T$, we can take it to be the zero section of the abelian scheme, such that its generic fiber $T_{\eta} \to \ab_K(\Phi)_{T_{\eta}}$ coincides with the base change of the zero section $\shu{K_{\Phi, h}} \to \shu{K_{\Phi, V} \rtimes K_{\Phi, h}}$. This can be done using the group structure of $\ab_K(\Phi)_T$. By Lemma \ref{lem: section of reductions, generic fiber} and \cite[Cor. 2.7.10]{PR24}, we have a unique extension of morphism of shtukas:
\[
\begin{tikzcd}
	{T^{\Dia/}} & {\Shum{\ovl{K}_{\Phi}, T}^{\Dia/}} \\
	{\Sht_{\GG_{\Phi, h}^*, \mu_{\Phi, h}^*}} & {\Sht_{\ovl{\PP}_{\Phi}^*, \bar{\mu}_{\Phi}^*}.}
	\arrow[from=1-1, to=1-2]
	\arrow[from=1-1, to=2-1]
	\arrow[from=1-2, to=2-2]
	\arrow[from=2-1, to=2-2]
\end{tikzcd}
\]
       In particular, for any $x \in \CE^{[[b_{\Phi, h}]]}(k)$, we can find a preimage $y \in \Shum{\ovl{K}_{\Phi}}(k)$ such that the $\ovl{\PP}^*_{\Phi}$-shtuka $y^*(\ovl{\PPs}_{\Phi}, \phi_{\ovl{\PPs}_{\Phi}})$ is isomorphic to the pushforward of the $\GG_{\Phi, h}^*$-shtuka $x^*(\ovl{\PPs}_{\Phi, h}, \phi_{\ovl{\PPs}_{\Phi, h}})$. This finishes the claim for the abelian scheme torsor part. Similar arguments can be made for the torus torsor part $\Shum{K_{\Phi}} \to \Shum{\ovl{K}_{\Phi}}$, with Lemma \ref{lem: abelian scheme action on shtuka, generic fiber} replaced by Lemma \ref{lem: torus action on shtuka, generic fiber}.
       
       To show central leaves are well positioned, we need to pass to $\Zb(\Phi) = \Delta_{\Phi,K}\backslash\Zb^{\bigsur}(\Phi)$ and show that $\Delta_{\Phi,K}$ stabilizes each central leaf. See the last paragraph of the proof of Proposition \ref{prop: Newton strata are well-positioned}.

       Finally, to show \emph{topological disjointness}, note that central leaves are closed in Newton strata; this follows from Proposition \ref{prop: Newton strata are well-positioned} and the fact that if $Y_1 \subset Y_2$ is a closed embedding of well-positioned subschemes, then $Y^{\natural}_{1, \Zb^{\bigsur}(\Phi)} \subset Y^{\natural}_{2, \Zb^{\bigsur}(\Phi)}$ are closed embeddings for all $\Phi \in \ca{CLR}(G, X)$. See \cite[Lem. 2.17, 2.18, 2.19]{Mao25}.
   \end{proof}

\subsubsection{Toroidal compactifications}

Recall that we have globally defined shtukas on the special fiber of the integral model of the toroidal compactification by Proposition \ref{prop-degeneration-int} and Corollary \ref{cor-degeneration-sp}:
\begin{equation}\label{eq: special fiber of toroidal compactification to shtuka}
   \Shumc{K}{\Sigma}(G, X)_{\bar{s}}^{\perf} \to \Sht_{\GG^c, \mu^c, \delta = 1}^W \to G^c\textit{-}\Isoc_{\mu^{c, -1}},
\end{equation}
which induces
\begin{equation*}
    \Upsilon_K^{\Sigma}(k): \Shumc{K}{\Sigma}(G, X)(k) \to C(\GG^c, \lrbracket{\mu^{c, -1}}),\quad \delta_K^{\Sigma}: \Shumc{K}{\Sigma}(G, X) \to B(G^c, \lrbracket{\mu^{c, -1}}).
\end{equation*}

\begin{prop}\label{prop: toroidal of Newton strata}
    Let $[b] \in B(G^c, \lrbracket{\mu^{c, -1}})$, and let $\NE^{[b]} \subset \Shum{K}(G, X)_{\bar{s}}$ be the Newton stratum. Then its partial toroidal compactification $(\NE^{[b]})^{\Sigma}$ (see \cite[Definition 2.3.1]{lan2018compactifications}) is the fiber $\delta_K^{\Sigma, -1}([b])$.
\end{prop}
\begin{proof}
    By restricting the diagram (\ref{eq-diag-degeneration-int}) to the special fiber, we have a commutative diagram:
    \begin{equation}\label{eq: graph over W}
\begin{tikzcd}
	{\Shumc{K}{\Sigma}(G, X)^{\perf}_{\bar{s}}} & {W^{\perf}_{\bar{s}}} & {\Delta_{\Phi, K}^{\circ}\backslash\Shum{K_{\Phi}}(P_{\Phi}, D_{\Phi})(\sigma)^{\perf}_{\bar{s}}} & {\Delta_{\Phi, K}^{\circ}\backslash\Shum{K_{\Phi, h}}(G_{\Phi, h}, D_{\Phi, h})^{\perf}_{\bar{s}}} \\
	{\Sht_{\GG^c, \mu^c, \delta = 1}^W} && {\Sht_{\PP_{\Phi}^*, \mu_{\Phi}^*, \delta = 1}^W} & {\Sht_{\GG_{\Phi, h}^*, \mu_{\Phi, h}^*, \delta = 1}^W} \\
	{G^{c}\textit{-}\Isoc_{\mu^{c, -1}}} && {P^*_{\Phi}\textit{-}\Isoc_{\mu_{\Phi}^{*, -1}}} & {G_{\Phi, h}^*\textit{-}\Isoc_{\mu_{\Phi, h}^{*, -1}}.}
	\arrow[from=1-1, to=2-1]
	\arrow[from=1-2, to=1-1]
	\arrow[from=1-2, to=1-3]
	\arrow[from=1-3, to=1-4]
	\arrow[from=1-3, to=2-3]
	\arrow[from=1-4, to=2-4]
	\arrow[from=2-1, to=3-1]
	\arrow["{\Int(g_{\Phi}^{-1})}"', from=2-3, to=2-1]
	\arrow[from=2-3, to=2-4]
	\arrow[from=2-3, to=3-3]
	\arrow[from=2-4, to=3-4]
	\arrow["{\Int(g_{\Phi}^{-1})}"', from=3-3, to=3-1]
	\arrow[from=3-3, to=3-4]
\end{tikzcd}
    \end{equation}
    By Proposition \ref{prop: Newton strata are well-positioned}, $\NE^{[b]}$ is well positioned, and its preimage in $W^0_{\bar{s}}$ is the preimage of some $(\NE^{[b]})^{\natural}_{\Zb^{\bigsur}(\Phi)} \subset \Zb^{\bigsur}(\Phi)_{\bar{s}} = \Shum{K_{\Phi, h}, \bar{s}}$. The preimage of its partial toroidal compactification $(\NE^{[b]})^{\Sigma}$ in $W_{\bar{s}}$ should also be the preimage of the same $(\NE^{[b]})^{\natural}_{\Zb^{\bigsur}(\Phi)}$; thus $(\NE^{[b]})^{\Sigma}$ is contained in the fiber $(\delta^{\Sigma}_K)^{-1}([b])$. On the other hand, Newton strata form a (weak) stratification $\Shum{K,\bar{s}} = \bigsqcup_{[b]} \NE^{[b]}$; thus $\Shumc{K}{\Sigma}(G, X)_{\bar{s}} = \bigsqcup_{[b]} (\NE^{[b]})^{\Sigma}$ (see \cite[Lem. 2.20]{Mao25}). This forces $(\NE^{[b]})^{\Sigma} = (\delta^{\Sigma}_K)^{-1}([b])$.
\end{proof}

\begin{prop}\label{prop: toroidal of central leaves}
    Let $[[b]] \in C(\GG^c, \lrbracket{\mu^{c, -1}})$, and let $\CE^{[[b]]} \subset \Shum{K}(G, X)_{\bar{s}}$ be the central leaf. Then its partial toroidal compactification $(\CE^{[[b]]})^{\Sigma}$ has the set of $k$-points $\Upsilon_K^{\Sigma}(k)^{-1}([[b]])$. In particular, we can upgrade $\Upsilon_K^{\Sigma}(k)$ to $\Upsilon_K^{\Sigma}: \Shumc{K}{\Sigma}(G, X)_{\bar{s}} \to C(\GG^c, \lrbracket{\mu^{c, -1}})$, and $(\CE^{[[b]]})^{\Sigma} = \Upsilon_K^{\Sigma, -1}([[b]])$.
\end{prop}
\begin{proof}
    We cannot directly apply the proof of Proposition \ref{prop: toroidal of Newton strata}, since
    \begin{equation}\label{eq: non-inj}
        C(\PP^*_{\Phi}, \lrbracket{\mu_{\Phi}^{*, -1}}) := \Sht_{\PP_{\Phi}^*, \mu_{\Phi}^*, \delta = 1}^W(k) \to C(\GG^*_{\Phi, h}, \lrbracket{\mu_{\Phi, h}^{*, -1}}) :=  \Sht_{\GG_{\Phi, h}^*, \mu_{\Phi, h}^*, \delta = 1}^W(k)
    \end{equation}
    is not injective.

    Taking the set of $k$-points of the diagram (\ref{eq: graph over W}), we have
\[
\begin{tikzcd}
	{\Shumc{K}{\Sigma}(G, X)(k)} & {W(k)} & {\Delta_{\Phi, K}^{\circ}\backslash\Shum{K_{\Phi}}(P_{\Phi}, D_{\Phi})(\sigma)(k)} & {\Delta_{\Phi, K}^{\circ}\backslash\Shum{K_{\Phi, h}}(G_{\Phi, h}, D_{\Phi, h})(k)} \\
	{C(\GG^c, \lrbracket{\mu^{c, -1}})} && {C(\PP^*_{\Phi}, \lrbracket{\mu_{\Phi}^{*, -1}})} & {C(\GG^*_{\Phi, h}, \lrbracket{\mu_{\Phi, h}^{*, -1}}).}
	\arrow[from=1-1, to=2-1]
	\arrow[from=1-2, to=1-1]
	\arrow[from=1-2, to=1-3]
	\arrow[from=1-3, to=1-4]
	\arrow[from=1-3, to=2-3]
	\arrow[from=1-4, to=2-4]
	\arrow["{\Int(g_{\Phi}^{-1})}"', from=2-3, to=2-1]
	\arrow[from=2-3, to=2-4]
\end{tikzcd}
\]

By Proposition \ref{prop: central leaves are well-positioned}, $\CE^{[[b]]}(k)$ is well positioned, and its preimage in $W(k)$ is the preimage of a union of central leaves $(\CE^{[[b]]})^{\natural}_{\Zb^{\bigsur}(\Phi)}(k) \subset \Zb^{\bigsur}(\Phi)(k) = \Shum{K_{\Phi, h}}(k)$. The preimage of $(\CE^{[[b]]})^{\Sigma}(k)$ in $W(k)$ should also be the preimage of $(\CE^{[[b]]})^{\natural}_{\Zb^{\bigsur}(\Phi)}(k)$. Let $\CE^{[[b_{\Phi, h}]]} \subset (\CE^{[[b]]})^{\natural}_{\Zb^{\bigsur}(\Phi)}$ be a central leaf. The arguments in the proof of Proposition \ref{prop: central leaves are well-positioned} show that the preimage of $\CE^{[[b_{\Phi, h}]]}(k)$ in $\Shum{K_{\Phi}}(k)$ is a central leaf $\CE^{[[b_{\Phi}]]}(k)$ with a unique $[[b_{\Phi}]] \in C(\PP^c_{\Phi}, \lrbracket{\mu_{\Phi}^{c, -1}})$.

Let us denote by $\CE^{[[b_{\Phi}]]}(\sigma)(k)$ (resp. $\NE^{[b_{\Phi}]}(\sigma)(k)$) the preimage of $\CE^{[[b_{\Phi, h}]]}(k)$ (resp. $\NE^{[b_{\Phi, h}]}(\sigma)(k)$) in $\Shum{K_{\Phi}}(\sigma)(k)$. By Lemma \ref{lem: closedness of central leaves, general group}, $\CE^{[[b_{\Phi}]]}(k) \subset \NE^{[b_{\Phi}]}(k)$ and $\CE^{[[b_{\Phi, h}]]}(k) \subset \NE^{[b_{\Phi, h}]}(k)$ are closed. Then $\CE^{[[b_{\Phi}]]}(\sigma)(k)$ is closed in $\NE^{[b_{\Phi}]}(\sigma)(k)$ and has to be the closure of $\CE^{[[b_{\Phi}]]}(k)$ in $\NE^{[b_{\Phi}]}(\sigma)(k)$. Applying Lemma \ref{lem: closedness of central leaves, general group} again, the central leaf in $\Shum{K_{\Phi}}(\sigma)(k)$ associated with $[[b_{\Phi}]]$ is closed in the Newton stratum and contains $\CE^{[[b_{\Phi}]]}(k)$; therefore, it must contain $\CE^{[[b_{\Phi}]]}(\sigma)(k)$. In particular, the preimages of $\CE^{[[b_{\Phi, h}]]}(k)$ in $W(k)$ and in $W^0(k)$ have the same image $[[b_{\Phi}]] \in C(\PP^c_{\Phi}, \lrbracket{\mu_{\Phi}^{c, -1}})$. Therefore, $(\CE^{[[b]]})^{\Sigma}(k)$ is contained in $(\Upsilon^{\Sigma}_K(k))^{-1}([[b]])$. Then we apply the last paragraph of the proof of Proposition \ref{prop: toroidal of Newton strata}.
\end{proof}

\begin{lem}\label{lem: closedness of central leaves, general group}
    Let $S$ be a perfect scheme and let $\PP$ be a quasi-parahoric group scheme. Given $p: S \to \Sht_{\PP}^W \to \Isoc_{P}$, we define (the set of $k$-points of) central leaves to be the fibers of $S(k) \to \Sht_{\PP}^W(k)$, and define Newton strata to be the fibers of $S \to B(P)$. Then Newton strata are locally closed in $S$, and central leaves are closed in the (set of $k$-points of) Newton strata.
\end{lem}
\begin{proof}
     When $\PP$ is parahoric, this follows from \cite[Prop. 2.15]{hamacher2025point}. Note that in \cite[\S 2]{hamacher2025point}, it is not necessary to assume that the generic fiber of $\PP$ is reductive; the setup applies to any flat affine group scheme of finite type over the base with connected fibers. The references \cite{rad2019local}, \cite{arasteh2021uniformizing}, and \cite[Thm. A.14]{imai2024tannakian} also work in this general setting. 
     
     When $\PP$ is quasi-parahoric, we instead apply the arguments in the proof of Lemma~\ref{lem: central leaves are closed in Newton strata}, using \cite[Prop.~2.3.1]{daniels2024conjecture} with $\pi_0(\PP^c)^{\phi}$ in place of $\pi_0(\GG^c)^{\phi}$.
\end{proof}

\subsection{KR strata and EKOR strata}

\subsubsection{Algebraicity}
Let us recall the setting in \cite[\S 3]{gleason2025specialization}. Given a morphism $f \colon S \to T$ of affine perfect schemes, $f$ is a universally subtrusive cover if and only if the induced morphism $f^{\dia} \colon S^{\dia} \to T^{\dia}$ is a $v$-cover.
\begin{lem}\label{lem: subtrusive torsor}
    Let $S$ be a perfect scheme, and let $H$ be the perfection of an affine group scheme of finite type. Then $H^1_{uv}(S, H) = H^1_v(S^{\dia}, H^{\dia})$, where $H^1_{uv}$ is defined under the universally subtrusive topology.
\end{lem}
\begin{proof}
    Since $\dia$ is fully faithful on the category of perfect schemes and preserves surjections, and by definition preserves fiber products, we have a morphism $H^1_{uv}(S, H) \to H^1_v(S^{\dia}, H^{\dia})$. On the other hand, given an affinoid perfectoid space $U = \Spa(A, A^+) \to S^{\dia}$ that trivializes the given $H^{\dia}$-torsor $\FF$ on $S^{\dia}$, $U_{\red}$ is represented by $\Spec A_{\red}$ ($A_{\red} = (A/A\cdot A^{\circ\circ})^{\perf}$, see \cite[Prop. 3.18]{gleason2025specialization}), and $\Spec A_{\red} \to S = \Spec R$ is surjective (indeed, we have the surjective specialization map $\spe_U: \Spa(A, A^+) \to \Spec A^+_{\red}$). We claim that $\FF_{\red}$ is represented by a perfect scheme. Following this, since the reduction functor preserves finite limits and $\dia$ is fully faithful, taking reduction gives a section $H^1_v(S^{\dia}, H^{\dia}) \to H^1_{uv}(S, H)$. Moreover, by adjointness, $(\FF_{\red})^{\dia} \to \FF$ is a morphism of $H^{\dia}$-torsors, which is automatically an isomorphism; then $H^1_v(S^{\dia}, H^{\dia}) \cong H^1_{uv}(S, H)$.

    Let $H_0$ be a linear algebraic group such that $H = H_0^{\perf}$. We claim that a $H_0$-torsor on $S$ is the same as a $H$-torsor on $S$: recall that a $H$-torsor over $S$ can be viewed as a trivialization $U \to S$ with a section $H(U\times_S U)$ that satisfies cocycle conditions. Since for any perfect scheme $T$, $H(T) = H_0(T)$, we have $H^1_{uv}(S, H_0) = H^1_{uv}(S, H)$.
    
    Now we work in the topos $\wdt{\SchPerf}$. Assume $H_0 = \GL_n$. A $\GL_n$-torsor $\PP$ on $S$ can be viewed as a vector bundle on $S$ via $\PP \times^{\GL_n} \OO_S^n$. By \cite[Thm. 4.1]{bhatt2017projectivity}, vector bundles over $S$ form a $v$-stack, thus $\PP \times^{\GL_n} \OO_S^n$ is representable by the full-faithfulness of $\dia$. By taking the framing, $\PP = \Isom_{S}(\OO_S, \PP \times^{\GL_n} \OO_S^n)$ is representable. In general, take a closed embedding $H_0 \to \GL_n$ such that the quotient is quasi-affine (see \cite[Prop. 1.3]{pappas2008twisted}), then $\GL_n^{\perf}/H_0^{\perf}$ is a perfect quasi-affine scheme. Let $G = \GL_n^{\perf}$, $H = H_0^{\perf}$. Consider the push-out torsor $\QQ = \PP \times^{H} G$, $\QQ$ is representable, in particular, $\QQ \to S$ is fpqc. Consider the quotient sheaf (under fpqc topology) $\QQ/H = \QQ \times^{G} G/H =\PP \times^{H} G/H$, it is trivialized over the fpqc cover $\QQ \to S$, thus $\QQ/H$ itself is a perfect quasi-affine scheme by the fpqc descent. Since the $v$-topology on the category of perfect qcqs schemes is subcanonical (see \cite[Remark 4.2]{bhatt2017projectivity}), then $\QQ/H = \PP \times^{H} G/H$ (viewed as a representable $v$-sheaf) has a section over $S$, thus $\PP$ is representable.
\end{proof}

\subsubsection{Kottwitz-Rapoport strata}\label{subsubsec: KR}

Recall that in \cite[\S 4.9]{PR24}, given a quasi-parahoric group scheme $\GG$, the authors constructed a $v$-sheaf theoretical local model diagram:
\begin{equation}\label{eq: shtukas to local model}
    \Sht_{\GG, \mu} \to [\GG^{\Dia}\backslash \vM_{\GG, \mu}^v],
\end{equation}
where $\GG^{\Dia}(S)$ ($S = \Spa(R, R^+) \in \Perf$) consists of pairs $(S^{\sharp}, g)$, $S^{\sharp} = \Spa(R^{\sharp}, R^{\sharp+})$ is an untilt of $S$ and $g \in \GG(R^{\sharp})$. Such constructions are functorial.

On the other hand, we have a local model diagram
\begin{equation}\label{eq: witt shtukas to local model}
    \Sht_{\GG, \mu}^W \to [\GG_0\backslash \lcM_{\GG, \mu}]
\end{equation}
constructed as follows: let $\Spec R \in \PCAlg$, $\Sht_{\GG, \mu}^{W, \square}(R)$ classifies tuples $((\PPs, \phi_{\PPs}), \alpha)$, where $(\PPs, \phi_{\PPs}) \in \Sht_{\GG, \mu}^W(R)$ and $\alpha: \PPs_0 \rightiso \phi^*(\PPs)$ is a trivialization of the $\GG$-torsor $\phi^*(\PPs)$ over $W(R)$. Then $\Sht_{\GG, \mu}^{W, \square}(R) \to \lcM_{\GG, \mu}(R)$ that maps $((\PPs, \phi_{\PPs}), \alpha)$ to $(\PPs, \phi_{\PPs}\circ\alpha)$ is $L^+\GG$-equivariant. Note that the $L^+\GG$-action on $\lcM_{\GG, \mu}$ factors through $\GG_0$, and the trivialization $\alpha$ is uniquely determined by its reduction over $R$ by the smoothness of $\GG$, then we get (\ref{eq: witt shtukas to local model}).

\begin{rk}
    Let us also compare (\ref{eq: witt shtukas to local model}) with the one constructed in \cite[\S 7.2.3]{xiao2017cycles} (cf. \cite[\S 4.2.2]{shen2021ekor}). The local Hecke stack $\mathrm{Hk}_{\GG, \mu}^W(R)$ (resp. $\lcM_{\GG, \mu}(R) \subset \Gra{\GG}^W(R)$) classifies the modifications $\gamma: \overleftarrow{\PPs} \dashrightarrow \overrightarrow{\PPs}$ (resp. $\PPs_0 \dashrightarrow \PPs$) of $\GG$-torsors over $W(R)$ of type $\mu$. By choosing a trivialization $\alpha: \PP_0 \rightiso \overleftarrow{\PPs}$, we have $\mathrm{Hk}_{\GG, \mu}^W = [L^+\GG\backslash \lcM_{\GG, \mu}]$. Also recall the truncated $\mathrm{Hk}_{\GG, \mu}^{W, (1)} = [\GG_0 \backslash \lcM_{\GG, \mu}]$. We have a morphism $\Sht_{\GG, \mu}^W \to \mathrm{Hk}_{\GG, \mu}^W$ that maps $(\PPs, \phi_{\PPs})$ to $(\overleftarrow{\PPs}=\phi^*\PPs, \overrightarrow{\PPs} = \PPs, \gamma = \phi_{\PPs})$. Compose with $\mathrm{Hk}_{\GG, \mu}^W \to \mathrm{Hk}_{\GG, \mu}^{W, (1)}$, we get (\ref{eq: witt shtukas to local model}).
\end{rk}

\begin{lem}
   Let $\GG_0$ be the special fiber of $\GG$. We have a projection $[\GG^{\Dia}\backslash \vM_{\GG, \mu}^v]_{\red} \to [\GG_0\backslash \lcM_{\GG, \mu}]$. In particular, from (\ref{eq: shtukas to local model}), we further have $\Sht_{\GG, \mu}^W \to [\GG_0\backslash \lcM_{\GG, \mu}]$ which coincides with (\ref{eq: witt shtukas to local model}).
\end{lem}
\begin{proof}
    Take the reduction of (\ref{eq: shtukas to local model}). Let $S \in \PCAlg$ be an affine perfect scheme, 
    \[ [\GG^{\Dia}\backslash \vM_{\GG, \mu}^v]_{\red}(S) = \Hom(S^{\dia}, [\GG^{\Dia}\backslash \vM_{\GG, \mu}^v]) \]
    gives $S^{\dia} \stackrel{p}{\leftarrow} \wdt{S} \stackrel{q}{\rightarrow} \vM_{\GG, \mu}^v$, where $p$ is a $\GG^{\Dia}$-torsor, and $q$ is $\GG^{\Dia}$-equivariant. Take reduction again, we have $(S^{\dia})_{\red} \stackrel{\bar{p}}{\leftarrow} \wdt{S}_{\red} \stackrel{q}{\rightarrow} \lcM_{\GG, \mu}$, where $\bar{p}$ is a $(\GG_0^{\dia})_{\red}$-torsor, and $q$ is $(\GG_0^{\dia})_{\red}$-equivariant: given a perfect algebra $R$, by \cite[Theorem 2]{gleason2025specialization},
    \[ (\GG^{\Dia})_{\red}(\Spec R) = \GG^{\Dia}(\Spd R) = \GG(R) = \GG_0(R) = (\GG_0^{\dia})_{\red}(\Spec R). \]
    Since $S \cong (S^{\dia})_{\red}$, $(\GG_0^{\dia})_{\red}$ is represented by $\GG_0^{\perf}$, Lemma \ref{lem: subtrusive torsor} shows that $\wdt{S}_{\red}$ is represented by a perfect scheme $\wdt{S}_0$ and $\wdt{S}_0 \to S$ is a $\GG_0^{\perf}$-torsor. On the other hand, since $\dia$ is fully faithful and $\lcM_{\GG, \mu}$ is represented by a perfect scheme (which we still denote by $\lcM_{\GG, \mu}$), $\wdt{S}_{\red} \to \lcM_{\GG, \mu}$ is represented by a $\GG_0^{\perf}$-equivariant morphism $\wdt{S}_0 \to \lcM_{\GG, \mu}$ between perfect schemes, thus we have a point in $[\GG_0^{\perf} \backslash \lcM_{\GG, \mu}](S)$. This gives a morphism $[\GG^{\Dia}\backslash \vM_{\GG, \mu}^v]_{\red} \to [\GG_0^{\perf} \backslash \lcM_{\GG, \mu}] = [\GG_0 \backslash \lcM_{\GG, \mu}]$. By constructions, these two $\Sht_{\GG, \mu}^W \to [\GG_0\backslash \lcM_{\GG, \mu}]$ coincide.
\end{proof}

We apply this to $\Shum{K}(G, X)^{\Dia/} \to \Sht_{\GG^c, \mu^c}$, and use Lemma \ref{lem: subtrusive torsor} again. We then have a morphism:
\begin{equation}\label{eq: morphism between quotient stacks}
    \Shum{K}(G, X)_{\bar{s}}^{\perf} \to [\GG_0^c \backslash \lcM_{\GG^c, \mu^{c}}].
\end{equation}

When $\GG^c$ is parahoric, it is homeomorphic to 
\begin{align*}
    \Adm_{G^c}(\lrbracket{\mu^{c, -1}})_{\KK^c} &:= \KK^c \backslash \KK^c \Adm_{G^c}(\lrbracket{\mu^{c, -1}}) \KK^c/\KK^c \\ &= W_{K^c}\backslash W_{K^c}\Adm_{G^c}(\lrbracket{\mu^{c, -1}})W_{K^c}/ W_{K^c} \subset W_{K^c} \backslash \wdt{W}^c/W_{K^c},
\end{align*}
where we use the induced Bruhat order (from the affine-Weyl group $\wdt{W}^c$ of $G^c$) on $W_{K^c} \backslash \wdt{W}^c/W_{K^c}$, thus the morphism (\ref{eq: morphism between quotient stacks}) induces a continuous morphism of underlying topological spaces
\[ l_K: \Shum{K}(G, X)_{\bar{s}} \to \Adm_{G^c}(\lrbracket{\mu^{c, -1}})_{\KK^c}.  \]

When $\GG^c$ is quasi-parahoric, since $\vM_{\GG^{c, \circ}, \mu^{c}} \cong \vM_{\GG^c, \mu^{c}}$, we identify 
\[\lcM_{\GG^c, \mu^{c}}(k) = \KK^{c, \circ} \Adm_{G^c}(\lrbracket{\mu^{c, -1}}) \KK^c/\KK^c.\] 
It is equal to $\KK^c \Adm_{G^c}(\lrbracket{\mu^{c, -1}}) \KK^c/\KK^c$ since the conjugation of $\KK^c/\KK^{c, \circ} \hookrightarrow \pi_1(G^c)_I$ acts trivially on $\Adm_{G^c}(\lrbracket{\mu^{c, -1}})$.. Thus 
\[ |[\GG_0^{c, \circ, \perf} \backslash \lcM_{\GG^{c, \circ}, \mu^{c}}]| = \Adm_{G^c}(\lrbracket{\mu^{c, -1}})_{\KK^{c, \circ}} = \Adm_{G^c}(\lrbracket{\mu^{c, -1}})_{\KK^c} = |[\GG_0^c \backslash \lcM_{\GG^c, \mu^{c}}]|. \]
We define the fibers of $l_K$ as \emph{Kottwitz-Rapoport strata}. In particular, KR strata are locally closed. 

By applying the exact same arguments, given a quasi-parahoric group scheme $\PP$ such that $(\PP, \mu)$ comes from boundary (see Definition \ref{def: PP, mu comes from boundary}), and by Corollary \ref{cor: comes from boundary, action factors through G_0}, the $L^+\PP$-action on $\vM_{\PP, \mu}^v$ factors through $\PP^{\Dia}$, we have morphisms
\[  \Sht_{\PP, \mu} \to  [\PP^{\Dia}\backslash \vM_{\PP, \mu}^v],\quad  \Sht_{\PP, \mu}^W \to [\PP_0 \backslash \lcM_{\PP, \mu}].  \]


\begin{prop}\label{prop: KR strata are well-positioned}
       KR strata are well positioned. Moreover, let $\KR_w$ be a KR stratum on $\Shum{K, \bar{s}}$ with some $w \in \Adm_{G^c}(\lrbracket{\mu^{c, -1}})_{\KK^c}$. Then, for each $\Phi \in \ca{CLR}(G, X)$, $(\KR_w)_{\Zb^{\bigsur}(\Phi)}^{\natural}$ is either empty or a finite union of KR strata $\KR_{w_{\Phi, h}}$ on $\Zb^{\bigsur}(\Phi)_{\bar{s}} \cong \Shum{K_{\Phi, h}, \bar{s}}$, for some collection of $w_{\Phi, h} \in \Adm_{G_{\Phi, h}^*}(\lrbracket{\mu_{\Phi, h}^{*, -1}})_{\KK_{\Phi, h}^*}$. The relation between $w$ and the collection of $w_{\Phi, h}$ is given in the proof.
\end{prop}
\begin{proof}
    KR strata are locally closed and are unions of central leaves. Since central leaves are well positioned by Proposition \ref{prop: central leaves are well-positioned}, KR strata are well positioned; see \cite[Lem. 2.20]{Mao25}. We describe their boundaries. Using the main diagram (\ref{eq: main diagram}), we have
   \[
\begin{tikzcd}
	{\Shum{K}(G, X)^{\dia}} & {W^{0, \dia}} & {\Delta_{\Phi, K}^{\circ}\backslash\Shum{K_{\Phi}}(P_{\Phi}, D_{\Phi})^{\dia}} & {\Delta_{\Phi, K}^{\circ}\backslash\Shum{K_{\Phi, h}}(G_{\Phi, h}, D_{\Phi, h})^{\dia}} \\
	{[\GG^{c, \Dia}\backslash \vM_{\GG^c, \mu^{c, -1}}^v]} & {} & {[\PP_{\Phi}^{*, \Dia}\backslash \vM_{\PP_{\Phi}^*, \mu_{\Phi}^{*, -1}}^v]} & {[\GG_{\Phi, h}^{*, \Dia}\backslash \vM_{\GG_{\Phi, h}^*, \mu_{\Phi, h}^{*, -1}}^v].}
	\arrow[from=1-1, to=2-1]
	\arrow[from=1-2, to=1-1]
	\arrow[from=1-2, to=1-3]
	\arrow[from=1-3, to=1-4]
	\arrow[from=1-3, to=2-3]
	\arrow[from=1-4, to=2-4]
	\arrow["{\Int(g_{\Phi}^{-1})}"', from=2-3, to=2-1]
	\arrow[from=2-3, to=2-4]
\end{tikzcd}
   \]
Applying the above arguments, we have
\[
\begin{tikzcd}
	{\Shum{K}(G, X)_{\bar{s}}} & {W^{0}_{\bar{s}}} & {\Delta_{\Phi, K}^{\circ}\backslash\Shum{K_{\Phi}}(P_{\Phi}, D_{\Phi})_{\bar{s}}} & {\Delta_{\Phi, K}^{\circ}\backslash\Shum{K_{\Phi, h}}(G_{\Phi, h}, D_{\Phi, h})_{\bar{s}}} \\
	{\Adm_{G^c}(\lrbracket{\mu^{c, -1}})_{\KK^c}} & {} & {\Adm_{P_{\Phi}^*}(\lrbracket{\mu_{\Phi}^{*, -1}})_{\KK_{\Phi}^*}} & {\Adm_{G_{\Phi, h}^*}(\lrbracket{\mu_{\Phi, h}^{*, -1}})_{\KK^*_{\Phi, h}}.}
	\arrow[from=1-1, to=2-1]
	\arrow[from=1-2, to=1-1]
	\arrow[from=1-2, to=1-3]
	\arrow[from=1-3, to=1-4]
	\arrow[from=1-3, to=2-3]
	\arrow[from=1-4, to=2-4]
	\arrow["{\Int(g_{\Phi}^{-1})}"', from=2-3, to=2-1]
	\arrow["{=(\ref{eq: P-grassmanian factorization})}", from=2-3, to=2-4]
\end{tikzcd}
\]
    The arguments in the proof of Proposition \ref{prop: Newton strata are well-positioned} can be adapted here without much change.
\end{proof}

   \begin{rk}
       We expect that the boundary of a KR stratum should be a single KR stratum, as proved in \cite{lan2018compactifications} for PEL types and in \cite{Mao25} for Hodge types. In other words, we expect that
       \[  \Adm_{G_{\Phi, h}^*}(\lrbracket{\mu_{\Phi, h}^{*, -1}})_{\KK^*_{\Phi, h}} \cong  \Adm_{L_{\Phi}^*}(\lrbracket{\mu_{\Phi, L}^{*, -1}})_{\KK_{\Phi, L}^*} \hookrightarrow \Adm_{G^c}(\lrbracket{\mu^{c, -1}})_{\KK^c} \]
       is an injection.
   \end{rk}

Next, we work with partial toroidal compactifications of KR strata. Equations (\ref{eq: special fiber of toroidal compactification to shtuka}) and (\ref{eq: morphism between quotient stacks}) give
\begin{equation}\label{eq: extension of local model diagram, special fiber}
     \Shumc{K}{\Sigma}(G, X)_{\bar{s}}^{\perf} \to [\GG_0^c \backslash \lcM_{\GG^c, \mu^{c}}],\quad l_K^{\Sigma}: \Shumc{K}{\Sigma}(G, X)_{\bar{s}} \to \Adm_{G^c}(\lrbracket{\mu^{c, -1}})_{\KK^c}. 
\end{equation}
\begin{prop}\label{prop: toroidal of KR strata}
    Let $w \in \Adm_{G^c}(\lrbracket{\mu^{c, -1}})_{\KK^c}$, and let $\KR_{w} \subset \Shum{K}(G, X)_{\bar{s}}$ be the KR stratum. Then its partial toroidal compactification $\KR_{w}^{\Sigma}$ is the fiber $l_K^{\Sigma, -1}(w)$.
\end{prop}
\begin{proof}
    Similar to the proof of Proposition \ref{prop: toroidal of Newton strata}.
\end{proof}

\subsubsection{Relation with schematic local model diagram}
In \cite{anschutz2022p} and \cite{gleason2024tubular}, it was proved that $\vM_{\GG, \mu}^v$ is represented by a normal scheme $\mathcal{M}_{\GG, \mu}$ that is flat and of finite type over $\OO_E$. In \cite[\S 4.9.1]{PR24}, it is conjectured that there exists a schematic local model diagram; that is to say, there exists a smooth morphism
\[ \pi_{\dR, \GG^c}: \Shum{K}(G, X) \to [\GG^c\backslash \mathcal{M}_{\GG^c, \mu^{c}}], \]
such that we have a ($2$-)commutative diagram
\begin{equation}\label{eq: schematic local model diagram}
\begin{tikzcd}
	{\Shum{K}(G, X)^{\Dia/}} & {\Sht_{\GG^c, \mu^c}} \\
	{[\GG^{c, \Dia/}\backslash \vM_{\GG^c, \mu^{c}}^v]} & {[\GG^{c, \Dia}\backslash \vM_{\GG^c, \mu^{c}}^v].}
	\arrow[from=1-1, to=1-2]
	\arrow["{\pi_{\dR, \GG^c}^{\Dia/}}", from=1-1, to=2-1]
	\arrow[from=1-2, to=2-2]
	\arrow[from=2-1, to=2-2]
\end{tikzcd}
\end{equation}
Here $\mathcal{M}_{\GG^c, \mu^{c}}^{\Dia/} = \mathcal{M}_{\GG^c, \mu^{c}}^{\Dia} = \vM_{\GG^c, \mu^{c}}^v$ by properness. In most cases of Hodge type, the Kisin-Pappas-Zhou integral models have schematic local model diagrams; see \cite[Thm. 7.1.3]{KPZ24} and \cite[Appendix A]{daniels2024conjecture}. In these cases, by the (perfect) smoothness of the morphism, we have a closure relation
\begin{equation}\label{eq: closure relation of KR}
    \ovl{\KR_w} = \bigsqcup_{w' \leq w, w' \in \Adm_{G^c}(\lrbracket{\mu^{c, -1}})_{\KK}} \KR_{w'},
\end{equation}
Here the KR strata are defined using the fibers of $\pi_{\dR, \GG^c}$, and this coincides with the fibers of $l_K$ (\ref{eq: morphism between quotient stacks}), by the commutativity of the above diagram (\ref{eq: schematic local model diagram}).
 
We focus on the abelian-type case. Assume $\GG$ is parahoric. In \cite[Thm. 7.2.20 and Rmk. 7.2.22]{KPZ24} (with supplements in \cite{daniels2024conjecture}), the authors showed that, given an abelian-type Shimura datum $(G, X)$, when $p>2$ (we use $(G, X)$ instead of $(G_2, X_2)$ to keep consistency of the notation in this subsection), there exists an integral model $\Shum{K}(G, X)$ of $\shu{K}(G, X)$ that has a list of good properties (see \cite[Thm. 7.2.20]{KPZ24}). Among these properties, we have a smooth morphism
\begin{equation}\label{eq: weak local model diagram}
    \pi_{\dR, \GG^{\ad}}: \Shum{K}(G, X) \to [\GG^{\ad, \circ}\backslash \mathcal{M}_{\GG^{\ad, \circ}, \mu^{\ad}}].
\end{equation}
Since $\Adm_{G^c}(\lrbracket{\mu^{c, -1}})_{\KK^c} \cong \Adm_{G^{\ad}}(\lrbracket{\mu^{\ad, -1}})_{\KK^{\ad, \circ}}$ (see \cite[Lem. 5.1.4]{shen2021ekor}), we have a stratification (\ref{eq: closure relation of KR}), here the KR strata are defined using the fibers of $\pi_{\dR, \GG^{\ad}}$. 

\begin{prop}\label{prop: schematic local model diagram, 2}
    When $\GG$ is parahoric and $p > 2$, we have a ($2$-)commutative diagram
\[
\begin{tikzcd}
	{\Shum{K}(G, X)^{\Dia/}} & {\Sht_{\GG^c, \mu^c}} \\
	{[\GG^{\ad, \circ, \Dia/}\backslash \vM_{\GG^{\ad, \circ}, \mu^{\ad}}^v]} & {[\GG^{\ad, \circ, \Dia}\backslash \vM_{\GG^{\ad, \circ}, \mu^{\ad}}^v].}
	\arrow[from=1-1, to=1-2]
	\arrow["{\pi_{\dR, \GG^{\ad, \circ}}^{\Dia/}}", from=1-1, to=2-1]
	\arrow[from=1-2, to=2-2]
	\arrow[from=2-1, to=2-2]
\end{tikzcd}
\]
\end{prop}
\begin{proof}
    When $(G, X)$ is Hodge-type, this is essentially \cite[Thm. A.3.3, Prop. 4.3.3]{daniels2024conjecture}. We generalize this result to the abelian-type case. We argue as follows: we first show the proposition when $G=G^\ad$, and then show it for a general abelian-type Shimura datum. \par
When $G=G^\ad$, by the construction in \cite{KPZ24}, there is a Hodge-type Shimura datum $(G',X')$ and a parahoric group scheme $\G'$ lifting $(G^\ad,X^\ad,\ca{G}^{\ad,\circ})$ that satisfy Conditions (A)–(E) in \S 7.1 of \emph{loc. cit.}. Consider the diagram:
\begin{equation*}
\begin{tikzcd}
	{\Shum{K'}(G', X')^{\Dia/}}\arrow[rr]\arrow[dr]\arrow[dd,"{\pi^{\Diamond/}_{\dr,\G^{\ad,\circ}}(G',X')}"] && {\Sht_{\GG^{\prime}, \mu^{\prime}}}\arrow[dr]\arrow[dd]& \\
    &{\Shum{K^\ad}{(G^\ad,X^\ad)}^{\Diamond/}}\arrow[rr]\arrow[dl,"{\pi^{\Diamond/}_{\dr,\G^{\ad,\circ}}(G^\ad,X^\ad)}"]&&{\Sht_{\GG^{\ad,\circ},\mu^\ad}}\arrow[dl]\\
	{[\GG^{\ad, \circ, \Dia/}\backslash \vM_{\GG^{\ad, \circ}, \mu^{\ad}}^v]}\arrow[rr] && {[\GG^{\ad,\circ, \Dia}\backslash \vM_{\GG^{\ad,\circ}, \mu^{\ad}}^v].}&
    \end{tikzcd}
    \end{equation*}
    The left triangle is commutative by the construction of the $\G^{\ad,\circ}$ local model diagram from the Hodge-type case in \cite{KPZ24}. The right triangle is canonically commutative. Since we can cover $\mathscr{S}_{K^\ad}(G^\ad,X^\ad)$ by a disjoint union of $\mathscr{S}_{K^{\prime,\alpha}}(G',X')$ satisfying (A)-(E) (where $K^{\prime,\alpha}$ is $K_p^{\prime,\alpha} K^{\prime,\alpha,p}$ for $K^{\prime,\alpha}_p$ the conjugation of $K'_p=\G'(\bb{Z}_p)$ by an element $\alpha\in G'(\bb{Q}_p)$ and $K^{\prime,\alpha,p}$ neat open compact), it suffices to show the commutativity of the bottom square by composition with $\mathscr{S}_{K'}^{\Diamond/}\to \mathscr{S}_{K^\ad}^{\Diamond/}$. Then the desired commutativity for adjoint Shimura data follows from the Hodge-type case and diagram-chasing.\par
    For a general $(G,X)$, the assertion follows from the last paragraph and the commutativity of
    \begin{equation*}
\begin{tikzcd}
	{\Shum{K}(G, X)^{\Dia/}}\arrow[d,"\pi^{\ad,\Dia/}"]\arrow[rr] && {\Sht_{\GG^{c}, \mu^{c}}}\arrow[d] \\
    {\Shum{K^\ad}{(G^\ad,X^\ad)}^{\Diamond/}}\arrow[rr]&&{\Sht_{\G^{\ad,\circ},\mu^\ad}}
    \end{tikzcd}
    \end{equation*}
    by functoriality of canonical integral models. Note that, since the local model diagram of an abelian-type integral model $\mathscr{S}_K(G,X)$ in \cite{KPZ24} is constructed by first pushing out the $\G^{\prime,\circ}$-local model diagram of $\mathscr{S}_{K'}(G',X')$ to a $\G^{\ad,\circ}$-local model diagram and then passing to $\mathscr{S}_K(G,X)$ as in \cite[Cor. 4.6.18]{KP15}, we still have $\pi^{\Dia/}_{\dr,\G^{\ad,\circ}}(G^\ad,X^\ad)\circ\pi^{\ad,\Dia/}=\pi^{\Dia/}_{\dr,\G^{\ad,\circ}}(G,X)$. The proposition is now proved.
\end{proof}

\begin{cor}\label{cor: compatibility of KR strata}
    Such constructed KR strata (using $\pi_{\dR, \GG^{\ad}}$) coincide with the ones defined in the last subsection (using $l_K$).
\end{cor}
\begin{cor}
    Under the setting of Proposition \ref{prop: schematic local model diagram, 2}, we have a closure relation:
\[\ovl{\KR_w^{\Sigma}} = \bigsqcup_{w' \leq w, w' \in \Adm_{G^c}(\lrbracket{\mu^{c, -1}})_{\KK}} \KR_{w'}^{\Sigma}. \]
\end{cor}
\begin{proof}
    This follows from (\ref{eq: closure relation of KR}) and the compatibility in Corollary \ref{cor: compatibility of KR strata} (see \cite[Prop. 7.1]{Mao25}).
\end{proof}
\begin{rk}\label{rk-local-model-diagram}
    In \cite[Thm. A]{DY25}, the authors constructed a shtuka map on Kisin-Pappas-Zhou integral models when $p > 2$ with parahoric level structures; the construction of the integral model and the shtuka map $\Shum{K}(G, X) \to \Sht_{\GG^c, \mu^{c}}$ in this paper coincides with \emph{loc. cit.} under the condition therein by the uniqueness of canonical integral models (see Theorem \ref{thm-ext-cim-ab} and \cite[Thm. B]{DY25}) and by the unique extension of shtukas from generic fibers to integral models (\cite[Cor. 2.7.10]{PR24}). We also remark that the existence of integral models at parahoric levels with \textbf{schematic local model diagrams} is still not known for all abelian-type Shimura data when $p=2$; cf. \cite{Yan25} for recent progress. 
\end{rk}

\subsubsection{Ekedahl-Kottwitz-Oort-Rapoport strata}

It is difficult to study EKOR strata only using shtukas; nevertheless, let us mention some easily deduced properties. By the constructions in \cite[\S 4.1]{shen2021ekor}, for parahoric $\GG$, from $\Shum{K}(G, X)_{\bar{s}}^{\perf} \to \Sht_{\GG^c, \mu^c}^W$, one can define
\[ \Shum{K}(G, X)_{\bar{s}}^{\perf} \to (\Sht_{\GG^c, \mu^c}^W)^{\loc(m, 1)},\quad \nu_K: \Shum{K}(G, X)_{\bar{s}} \to \prescript{\KK^c}{}{\Adm_{G^c}(\lrbracket{\mu^{c, -1}})}, \]
where $(\Sht_{\GG^c, \mu^c}^W)^{\loc(m, 1)}$ is the stack of $(m, 1)$-truncated shtukas. It is an algebraic stack, with underlying topological space homeomorphic to $\prescript{\KK^c}{}{\Adm_{G^c}(\lrbracket{\mu^{c, -1}})}$. We define the preimages of $\nu_K$ as \emph{EKOR strata}. EKOR strata are locally closed.

Similarly, we can consider the partial toroidal compactifications of EKOR strata. (\ref{eq: special fiber of toroidal compactification to shtuka}) gives
\[ \Shumc{K}{\Sigma}(G, X)_{\bar{s}}^{\perf} \to (\Sht_{\GG^c, \mu^c}^W)^{\loc(m, 1)},\quad \nu_K^{\Sigma}: \Shumc{K}{\Sigma}(G, X)_{\bar{s}} \to \prescript{\KK^c}{}{\Adm_{G^c}(\lrbracket{\mu^{c, -1}})}. \]

In the abelian-type case, when $\GG$ is quasi-parahoric, we further define the EKOR strata on $\Shum{K, \bar{s}}$ (resp. $\Shumc{K, \bar{s}}{\Sigma}$) as the images of EKOR strata on $\Shum{K^{\circ}, \bar{s}}$ (resp. $\Shumc{K^{\circ}, \bar{s}}{\Sigma}$) under the finite \'etale surjection $\Shum{K^{\circ}, \bar{s}} \to \Shum{K, \bar{s}}$ (resp. $\Shumc{K^\circ, \bar{s}}{\Sigma} \to \Shumc{K, \bar{s}}{\Sigma}$; the \'etaleness can be easily seen using the \'etaleness of $\Shum{K^{\circ}_{\Phi}} \to \Shum{K_{\Phi}}$ for any $\Phi \in \ca{CLR}(G, X)$).
\begin{prop}\label{prop: toroidal of EKOR strata}
    EKOR strata are well positioned. Let $w \in \prescript{\KK^c}{}{\Adm_{G^c}(\lrbracket{\mu^{c, -1}})}$, $\EKOR_{w} \subset \Shum{K}(G, X)_{\bar{s}}$ be the EKOR strata; then its partial toroidal compactification $\EKOR_{w}^{\Sigma}$ is the fiber $\nu_K^{\Sigma, -1}(w)$.
\end{prop}
\begin{proof}
    EKOR strata are locally closed and are unions of central leaves. Since central leaves are well positioned by Proposition \ref{prop: central leaves are well-positioned}, EKOR strata are well positioned (see \cite[Lem. 2.20]{Mao25}). Write $\EKOR_w = \bigsqcup_{i \in I} \CE^{[[b_i]]}$ as a union of central leaves; then $\EKOR_w^{\Sigma} = \bigsqcup_{i \in I} (\CE^{[[b_i]]})^{\Sigma}$ by \cite[Lem. 2.20]{Mao25}. By Proposition \ref{prop: toroidal of central leaves}, $(\CE^{[[b_i]]})^{\Sigma} \subset \nu_K^{\Sigma, -1}(w)$; then $\EKOR_w^{\Sigma} \subset \nu_K^{\Sigma, -1}(w)$. Since $\Shum{K}(G, X)_{\bar{s}} = \bigsqcup_{w} \EKOR_w$, $\Shumc{K}{\Sigma}(G, X)_{\bar{s}} = \bigsqcup_{w} \EKOR_w^{\Sigma}$ by \cite[Lem. 2.20]{Mao25}. This forces $\EKOR_w^{\Sigma} = \nu_K^{\Sigma, -1}(w)$.
\end{proof}

\begin{rk}
    Finally, let us briefly discuss the Ekedahl–Oort strata. Assume that $\GG$ is hyperspecial. Then, for any $[\Phi] \in \Cusp_K(G, X)$, $\GG_{\Phi, h}$ is hyperspecial. In \cite[\S 5.3]{xiao2017cycles}, there is a natural perfectly smooth morphism
$$ (\Sht_{\GG, \mu}^W)^{\loc(m, 1)} \to G\textit{-}\ZIP_{\mu^{-1}}^{\perf}, $$
where $G\textit{-}\ZIP^{\perf}$ is the perfection of $G\textit{-}\ZIP$. By composing with $(\Shum{K, \bar{s}})^{\perf} \to \Sht_{\GG^c, \mu^c}^W$ and $(\Shum{K, \bar{s}}^{\Sigma})^{\perf} \to \Sht_{\GG^c, \mu^c}^W$, we obtain
$$ \epsilon_K: (\Shum{K, \bar{s}})^{\perf} \to G^c\textit{-}\ZIP_{\mu^{c, -1}}^{\perf},\quad \epsilon_K^{\Sigma}: (\Shum{K, \bar{s}}^{\Sigma})^{\perf} \to G^c\textit{-}\ZIP_{\mu^{c, -1}}^{\perf}.$$
    We can define the EO strata as the fibers of $\epsilon_K$, endowed with induced reduced subscheme structures, and similarly show that the EO strata are well positioned.\par
Thus, we obtain a similar diagram for $\epsilon_K$ and $\epsilon_K^\Sigma$ as those in Proposition \ref{prop: KR strata are well-positioned}, so that we deduce that, when $\Sigma$ is smooth, $\epsilon_K$ is perfectly smooth if and only if $\epsilon_K^\Sigma$ is so (cf. \cite[Cor. 4.51]{Wu25}).\par 
Moreover, the partial toroidal compactification of the EO strata coincides with the EO strata defined using $\epsilon_K^{\Sigma}$. Also, when $\ca{G}$ is hyperspecial, EKOR strata are the same as EO strata; thus, we omit the proof for EO strata here.
    \end{rk}
\begin{rk}
For the KR strata, we expect the existence of \textbf{schematic} local model diagram. For EO strata, one expects that these morphisms arise from a prismatic or syntomic approach, where the coherent and infinitesimal structures can be directly seen; that is, the morphisms $\epsilon_K$ and $\epsilon_K^{\Sigma}$ should respectively be the perfections of smooth morphisms
   $$ \Shum{K, s} \to G^c\textit{-}\ZIP_{\mu^{c, -1}},\quad \Shum{K, s}^{\Sigma} \to G^c\textit{-}\ZIP_{\mu^{c, -1}}. $$
   For the interior morphism $\epsilon_K$, the most general result currently available is given by \cite{MY2026}; let us also mention the previous works \cite{oort2001stratification}, \cite{viehmann2013ekedahl}, \cite{pink2015f}, \cite{zhang2018ekedahl}, \cite{imai2023prismatic}, etc. Indeed, when $(G, X)$ is of abelian type or when $p$ is large enough, Madapusi and Youcis constructed a formally \'etale map of $p$-adic formal stacks over $\Spf \OO_E$: $\wdh{\Shum{K}} \to \mathrm{BT}_{\infty}^{\GG^c, -\mu^c}$, where $\mathrm{BT}_{\infty}^{\GG^c, -\mu^c}$ serves as the stack of '$p$-divisible groups with $\GG^c$-structures' without actually working with $p$-divisible groups. The compositions $\wdh{\Shum{K}} \to \mathrm{BT}_{n}^{\GG^c, -\mu^c}$ are smooth and surjective, for all $n > 0$. There is a smooth surjection $\mathrm{BT}_1^{\GG^c, -\mu^c} \otimes k_E \to G^c\textit{-}\ZIP_{\mu^{c, -1}}$. The composition gives the deperfection of $\epsilon_K$.
\end{rk}

\appendix
\section{Some category theory}\label{app-cat}

This appendix serves as a complement to the categorical language employed throughout the paper, offering a direct reference applicable to the context considered herein.
\subsection{Morphisms between presheaves in categories}\label{subsec-mor-fibered-cat} 
Let $\ca{C}$ be a ($1$-)category. Let $\md{Categories}$ (resp. $\md{Groupoids}$) be the $2$-category (resp. $(2,1)$-category) of categories (resp. groupoids). Let $F_1$ and $F_2$ be two ($2$-)functors (or presheaves) $F_1,F_2:\ca{C}^\mrm{op} \to \md{Categories}$. \par 
As in \cite[\href{https://stacks.math.columbia.edu/tag/02XV}{Ex. 02XV}]{stacks-project}, there are fibered categories $p_1:\ca{X}_1\to \ca{C}$ and $p_2:\ca{X}_2\to \ca{C}$ corresponding to $F_1$ and $F_2$, respectively. More precisely, a category $\ca{X}$ with a functor $p$ to $\ca{C}$ is called a fibered category if, for any $x\in \ca{X}$ with $f:V\to p(x)\in \mrm{Mor}_{\ca{C}}(V,p(x))$, there exists a strongly Cartesian lift $y\to x$ over $f$ (cf. \cite[\href{https://stacks.math.columbia.edu/tag/02XM}{Def. 02XM}]{stacks-project}).\par
A \emph{($1$-)morphism} $\mathscr{F}: (\ca{X}_1,p_1)\to (\ca{X}_2,p_2)$ is a $1$-morphism in the $2$-category of fibered categories over $\ca{C}$. That is, $\mathscr{F}$ is a functor from $\ca{X}_1$ to $\ca{X}_2$ such that $p_2\circ\mathscr{F}=p_1$, and $\mathscr{F}$ sends a strongly Cartesian morphism to a strongly Cartesian morphism (see \cite[\href{https://stacks.math.columbia.edu/tag/02XP}{Def. 02XP}]{stacks-project}).\par 
On the other hand, we define $1$-morphisms between presheaves in categories.
\begin{definition}\label{def-1-mor-presh-cat}
Define a ($1$-)morphism (or a natural transformation) $\mbf{F}:F_1\to F_2$ to be the collection $\{\mbf{F}(U);\phi(U_1,U_2,f)\}_{U,U_1,U_2\in\ob \ca{C},f\in \mrm{Mor}_{\ca{C}}(U_1,U_2)}$ as follows:\par
\begin{enumerate}
\item For any $U\in\ob \ca{C}$, $\mbf{F}(U):F_1(U)\to F_2(U)$ is a $1$-morphism (i.e., a functor) between categories.
\item For $f\in \mrm{Mor}_{\ca{C}}(U_1,U_2)$ with $U_1,U_2\in\ob \ca{C}$, $$\phi(U_1,U_2,f):\mbf{F}(U_1)F_1(f)(F_1(U_2))\to F_2(f)\mbf{F}(U_2)(F_1(U_2))$$ is a natural transformation between functors $\mbf{F}(U_1)F_1(f)$ and $F_2(f)\mbf{F}(U_2)$ mapping from $F_1(U_2)$ to $F_2(U_1)$ such that, 
\begin{enumerate}
\item $\phi(U,U,\mrm{id}_U)(x)=\mrm{id}_{\mbf{F}(U)x}$ for any $U\in \ob \ca{C}$ and $x\in F_1(U)$.
\item For any morphisms in $\ca{C}$, $f:U_1\to U_2$ and $g:U_2\to U_3$, and any $x_3\in\ob F_1(U_3)$, we have $F_2(f)(\phi(U_2,U_3,g))\circ \phi(U_1,U_2,f)(F_1(g)x_3)=\phi(U_1,U_3,g\circ f)(x_3)$.
\item For any pair $(U_1,U_2)$ and $f:U_1\to U_2$, the diagram 
\begin{equation*}
    \begin{tikzcd}
     F_1(U_2)\arrow[rr,"{F_1(f)}"]\arrow[d,"\mbf{F}(U_2)"]&& F_1(U_1)\arrow[d,"\mbf{F}(U_1)"]\\ 
     F_2(U_2)\arrow[rr,"{F_2(f)}"]&&F_2(U_1)
    \end{tikzcd}
\end{equation*}
is commutative up to composing $\phi(U_1,U_2,f)$.
\end{enumerate}
\end{enumerate}
\end{definition}
In fact, the definitions of $1$-morphisms in the two contexts correspond to each other. 
\begin{lem}\label{lem-1-morphism}$\mathscr{F}$ determines and is determined by a morphism $\mbf{F}$ from $F_1$ to $F_2$.
\end{lem}
\begin{proof}
Recall that the fibered category $p_1:\ca{X}_1\to\ca{C}$ is defined as follows: the objects are $(U,x)$ such that $U\in \ob \ca{C}$ and $x\in \ob F_1(U)$; the morphisms are $\mrm{Mor}_{\ca{X}_1}((U_1,x_1),(U_2,x_2))=\{(f,\phi)| f\in \mrm{Mor}_{\ca{C}}(U_1,U_2),\phi\in \mrm{Mor}_{F_1(U_1)}(x_1,F_1(f)x_2)\}$. The composition is defined as $(g,\psi)\circ (f,\phi)=(g\circ f,F_1(f)(\psi)\circ \phi)$.\par 
The $1$-morphism $\mathscr{F}$ maps $(U,x)$ to $(V,y)\in \ob\ca{X}_2$. Since $p_2\circ\mathscr{F}=p_1$, we have that $V=U$ and $y\in \ob F_2(U)$. This determines an assignment $\mathscr{F}(U)$ from $\ob F_1(U)$ to $\ob F_2(U)$.\par 
Hence, the functor $\mathscr{F}$ maps $\mrm{Mor}_{\ca{X}_1}((U_1,x_1),(U_2,x_2))$ to $\mrm{Mor}_{\ca{X}_2}((U_1,\mathscr{F}(U_1)x_1),(U_2,\mathscr{F}(U_2)x_2))$. More precisely, a pair $(f,\phi)$ as above is sent to $(\mathscr{F}(f),\mathscr{F}(\phi))$. But $\mathscr{F}(f)=f$ since $p_2\circ\mathscr{F}=p_1$. So $\mathscr{F}(\phi): \mathscr{F}(U_1)x_1\to F_2(f)(\mathscr{F}(U_2)x_2)$ is a uniquely determined morphism in $F_2(U_1)$ induced by $\mathscr{F}$.\par 
Now let $(U_1,x_1)=(U_1,F_1(f)x_2)$ and let $\phi=\mrm{id}$. Then $\mathscr{F}(\phi)=: \phi_{\mathscr{F}}(U_1,U_2,f)(x_2)$ is a functor $\mathscr{F}(\phi):\mathscr{F}(U_1)(F_1(f)x_2)\to F_2(f)(\mathscr{F}(U_2)x_2)$. 
Moreover, fix any $(t: x_2\to x_2')\in \mrm{Mor}_{F_1(U_2)}(x_2,x_2')$ and consider the commutative diagram
\begin{equation*}
    \begin{tikzcd}
    (U_1,F_1(f)x_2)\arrow[d,"F_1(f)(t)"]\arrow[rr,"{(f,\mrm{id})}"]&&(U_2,x_2)\arrow[d,"t"]\\
    (U_1,F_1(f)x_2')\arrow[rr,"{(f,\mrm{id})}"]&& (U_2,x_2').
    \end{tikzcd}
\end{equation*}
The diagram above is commutative after applying $\mathscr{F}$ as $\mathscr{F}$ is a functor, which implies the commutativity of 
\begin{equation*}
    \begin{tikzcd}
    \mathscr{F}(U_1)F_1(f)x_2\arrow[rr,"{\phi_{\mathscr{F}}(U_1,U_2,f)(x_2)}"]\arrow[d,"\mathscr{F}(U_1)F_1(f)t"]&&F_2(f)\mathscr{F}(U_2)x_2\arrow[d,"F_2(f)\mathscr{F}(t)"]\\
    \mathscr{F}(U_1)F_1(f)x_2'\arrow[rr,"{\phi_{\mathscr{F}}(U_1,U_2,f)(x_2')}"]&& F_2(f)\mathscr{F}(U_2)x_2'.
    \end{tikzcd}
\end{equation*}
This implies that $\phi_{\mathscr{F}}(U_1,U_2,f)$ is a natural transformation.\par
Hence, the diagram of functors
\begin{equation*}
    \begin{tikzcd}
     F_1(U_2)\arrow[rr,"{F_1(f)}"]\arrow[d,"\mathscr{F}(U_2)"]&& F_1(U_1)\arrow[d,"\mathscr{F}(U_1)"]\\ 
     F_2(U_2)\arrow[rr,"{F_2(f)}"]&&F_2(U_1)
    \end{tikzcd}
\end{equation*}
is commutative up to uniquely determined $2$-morphisms. Then it can be checked that $\{\mathscr{F}(U);\phi_{\mathscr{F}}(U_1,U_2,f)\}$ determines a $1$-morphism $\mbf{F}:F_1\to F_2$. The second condition of $\mbf{F}$ listed above follows from the fact that $\mathscr{F}$ is a functor that preserves compositions.\par
Conversely, suppose that there is a $1$-morphism $\mbf{F}:F_1\to F_2$ given by the collection $\{\mbf{F}(U);\phi(U_1,U_2,f)\}$. Then set $\mathscr{F}(U,x_1)=(U,\mbf{F}(U)x_1)$. Let $(f:U_1\to U_2,\phi:x_1\to F_1(f)x_2)\in \mrm{Mor}_{\ca{X}_1}((U_1,x_1),(U_2,x_2))$. Set $\mathscr{F}(f,\phi)=(f,\phi(U_1,U_2,f)(x_2)\circ\mbf{F}(U_1)(\phi))$.\par 
Then it can be checked that $\mathscr{F}$ is a functor. We only check that it preserves compositions because the other conditions are easier. Suppose that we are given $U_1\xrightarrow{f}U_2\xrightarrow{g} U_3$, $x_i\in F_1(U_i)$ for $i=1,2,3$, and $\phi: x_1\to F_1(f)x_2$ and $\psi:x_2\to F_1(g)x_3$. We have $\mathscr{F}(g\circ f,F_1(f)(\psi)\circ\phi)=$
$$(g\circ f, \phi(U_1,U_3,g\circ f)(x_3)\circ \mbf{F}(U_1)(F_1(f)(\psi)\circ\phi)).$$
The last expression is computed as
\begin{equation*}
\begin{split}
&(g\circ f, F_2(f)(\phi(U_2,U_3,g))\circ \phi(U_1,U_2,f)(F_1(g)x_3)\circ\mbf{F}(U_1)(F_1(f)(\psi)\circ \phi))=\\
&(g\circ f, F_2(f)(\phi(U_2,U_3,g))\circ \phi(U_1,U_2,f)(F_1(g)x_3)\circ(\mbf{F}(U_1)(F_1(f)(\psi))\circ \mbf{F}(U_1)(\phi)))=\\
&(g\circ f,F_2(f)(\phi(U_2,U_3,g))\circ F_2(f)\mbf{F}(U_2)(\psi) \circ \phi(U_1,U_2,f)(x_2)\circ \mbf{F}(U_1)(\phi) )=\\
&(g\circ f, F_2(f)(\phi(U_2,U_3,g)\circ\mbf{F}(U_2)(\psi))\circ (\phi(U_1,U_2,f)(x_2)\circ \mbf{F}(U_1)(\phi)))=\\
&(g,\phi(U_2,U_3,g)\circ\mbf{F}(U_2)(\psi))\circ(f,\phi(U_1,U_2,f)(x_2)\circ\mbf{F}(U_1)(\phi))=
\mathscr{F}(g,\psi)\circ\mathscr{F}(f,\phi).
\end{split}
\end{equation*}
The second line to the third line follows from the fact that $\phi(U_1,U_2,f)$ is a natural transformation, from which the diagram
\begin{equation*}
    \begin{tikzcd}
    \mbf{F}(U_1)F_1(f)x_2\arrow[rrr,"{\phi(U_1,U_2,f)(x_2)}"]\arrow[d,"\mbf{F}(U_1)F_1(f)(\psi)"]&&&F_2(f)\mbf{F}(U_2)x_2\arrow[d,"F_2(f)\mbf{F}(U_2)(\psi)"]\\
    \mbf{F}(U_1)F_1(f)(F_1(g)x_3)\arrow[rrr,"{\phi(U_1,U_2,f)(F_1(g)x_3)}"]&&& F_2(f)\mbf{F}(U_2)(F_1(g)x_3)
    \end{tikzcd}
\end{equation*}
commutes. 
\end{proof}
\subsection{2-limits}\label{subsec-limits-def}
\subsubsection{}We explain what a $2$-limit is, specialized to our situation.
\begin{definition}\label{def-2-lim}
Let $p:\ca{X}\to\ca{C}$ be a fibered category in groupoids constructed from a presheaf $F:\ca{C}^{\mrm{op}}\to \md{Groupoids}$ (see \cite[\href{https://stacks.math.columbia.edu/tag/003T}{Def. 003T}]{stacks-project} and \cite[\href{https://stacks.math.columbia.edu/tag/0049}{Ex. 0049}]{stacks-project}). Define 
$$2\text{-}\varprojlim_{\ca{C}^{\mrm{op}}}\ca{X}$$
as a groupoid. More precisely, 
\begin{enumerate}
\item The objects of it are of the form 
$$L:=\{(U,x_U);T_{f_{U,V}}\}_{U\in\ob \ca{C}^{\mrm{op}},f_{U,V}\in\mrm{Mor}_{\ca{C}^{\mrm{op}}}(U,V)}$$ where $x_U\in \ob F(U)$, $T_{f_{U,V}}\in \mrm{Mor}_{\ca{X}^{\mrm{op}}}((U,x_U),(V,x_V))$ such that $p^{\mrm{op}}(T_{f_{U,V}})=f_{U,V}$, $T_{f_{U_2,U_1}}\circ T_{f_{U_3,U_2}}=T_{f_{U_3,U_1}}$ for $f_{U_2,U_1}\circ f_{U_3,U_2}=f_{U_3,U_1}$ in $\ca{C}^{\mrm{op}}$. 
\item The morphisms from $L$ to $L'=\{(U,x'_U);T'_{f_{U,V}}\}$ are of the form $\{C_U\}_{U\in \ob\ca{C}^{\mrm{op}}}$, where $C_U\in\mrm{Mor}_{\ca{X}^{\mrm{op}}}((U,x_U),(U,x_U'))$ such that $T'_{f_{U,V}}\circ C_U=C_V\circ T_{f_{U,V}}$.
\end{enumerate}
\end{definition}
In Definition \ref{def-2-lim}, an object $L$ determines and is determined by a functor $\mathscr{L}:\ca{C}\to \ca{X}$ sending $U$ to $(U,x_U)$ and $f_{U,V}^{\mrm{op}}$ to $T_{f_{U,V}}^{\mrm{op}}$ such that $p\circ\mathscr{L}=\mrm{id}_{\ca{C}}$. \par
Let $F_1$, $F_2$, $(\ca{X}_1,p_1)$, $(\ca{X}_2,p_2)$, $\mbf{F}$, and $\mathscr{F}$ be as in \S\ref{subsec-mor-fibered-cat}. 
Suppose that there is a $1$-morphism $\pi:\ca{X}_2\to\ca{X}_1$ (such that $p_1\circ\pi=p_2$, and) such that $(\ca{X}_2,p)$ is isomorphic to the fibered category in groupoids over $\ca{X}_1$ associated with a functor $\pi^{presh}:\ca{X}_1^{\mrm{op}}\to \md{Groupoids}$. Denote by $\Phi: F_2\to F_1$ the $1$-morphism between presheaves corresponding to $\pi$ according to Lemma \ref{lem-1-morphism}. 
\begin{lem}\label{lem-mor-lim-correspond}
In this setup, an object in $2$-$\varprojlim_{\ca{X}_1^{\mrm{op}}}\ca{X}_2$ determines and is determined by a section $\mathscr{S}:\ca{X}_1\to \ca{X}_2$ of $\pi:\ca{X}_2\to \ca{X}_1$ in $1$-$\mrm{Mor}_{\ca{C}}(\ca{X}_1,\ca{X}_2)$, which, from Lemma \ref{lem-1-morphism}, corresponds to a $1$-morphism $\Psi: F_1\to F_2$ such that $\Phi\circ\Psi=\mrm{id}_{F_1}$. Moreover, the groupoid of such sections $\mathscr{S}$ (with morphisms among them defined by the natural transformations between functors from $\ca{X}_1$ to $\ca{X}_2$) is \emph{isomorphic} to the groupoid $2$-$\varprojlim_{\ca{X}_1^{\mrm{op}}}\ca{X}_2$.
\end{lem}
\begin{proof}
In the context of this lemma, we have $p_2\circ\mathscr{L}=p_1\circ\pi\circ\mathscr{L}=p_1$, which means that $\mathscr{L}$ is automatically a morphism in $1$-$\mrm{Mor}_{\ca{C}}(\ca{X}_1,\ca{X}_2)$; note that every morphism is strongly Cartesian by \cite[\href{https://stacks.math.columbia.edu/tag/003V}{Lem. 003V}]{stacks-project}. The first part of the lemma now follows from Lemma \ref{lem-1-morphism}. The isomorphism of groupoids follows from the construction in Definition \ref{def-2-lim}(2).
\end{proof}
\begin{rk}\label{rk-use-three-notions}
By Lemma \ref{lem-mor-lim-correspond}, the three notions ``an object in $2$-$\varprojlim_{\ca{X}_1^{\mrm{op}}}\ca{X}_2$'', ``a $1$-morphism section of $\pi: \ca{X}_2\to \ca{X}_1$'', and ``a $1$-morphism section of $\Phi:F_2 \to F_1$'' can and should be used interchangeably in this paper.   
\end{rk}
Assume further that $\ca{X}_1^{\lim}:=$ 
$$\twolim\limits_{U\in\ob \ca{C}^{\mrm{op}}}\varinjlim_{F_1(U)}F_1(U)$$
exists. If it exists, this is a groupoid since the inner colimits $\varinjlim_{F_1(U)}F_1(U)$ (that is, the colimit in the category $F_1(U)$) are groupoids, as they are unique up to unique isomorphisms by definition.\par
Denote by $\ca{X}_1'$ the full subcategory of $\ca{X}_1$ consisting of objects of the form $(U,\varinjlim_{F_1(U)} F_1(U))$. Note that there is also a projection $\ell:\ca{X}_1\to \ca{X}_1'$ according to the definition of colimits, which makes natural inclusion $i:\ca{X}_1'\to \ca{X}_1$ a section of $\ell$, i.e., $\ell\circ i=\mrm{id}_{\ca{X}_1}$.
\begin{lem}\label{lem-equiv-subcat-sec}
With the assumptions above, there is a natural equivalence of groupoids
$$ 2\text{-}\varprojlim_{\ca{X}_1^{\mrm{op}}}\ca{X}_2\iso 2\text{-}\varprojlim_{\ca{X}_1^{\prime,\mrm{op}}}\ca{X}_2.$$
The second limit is formed for the projection $\ell\circ\pi:\ca{X}_2\to\ca{X}_1'$.
\end{lem}
\begin{proof}
The functor $E_1$ from the LHS to the RHS is the natural projection; the functor $E_2$ in the other direction is defined by pulling back objects of $\ca{X}_2$ over $\ca{X}_1'\sbst \ca{X}_1$ via the morphisms defined by colimits (since $\pi:\ca{X}_2\to \ca{X}_1$ is a fibered category by assumption). The fact that $E_2 \circ E_1\simeq \mrm{id}$(of LHS) is due to the universal property of strongly Cartesian morphisms and the assumption that $\ca{X}_2$ is fibered in groupoids: Write an object in LHS by $L=\{(x_1,x_2);T_{f_{x_1,x_1'}}\}$. Suppose that $x_1\in \ob F_1(U)$ maps to $x_1^*\in \ob \varinjlim_{F_1(U)} F_1(U)$. Suppose that $p:x_2^*\to x_1^*$ is the object in the data of $L$. Then there is a (strongly Cartesian) lift $g:x_2'\to x_2^*$ of $f:x_1\to x_1^*$. As $\ca{X}_2\to \ca{X}_1$ is fibered in groupoids and there is a isomorphism $T_{f_{x_1,x_1^*}}:x_2\to x_2^*$, there is a unique isomorphism $C_{x_1}:x_2\to x_2'$.\par The commutativity in Definition \ref{def-2-lim}(2) is obtained in a similar way. In fact, let $y_1$, $y_1^*$, $y_2^*$, $y_2'$, and $y_2$ be the objects constructed in the same way for another $V\to U$ in $\ca{C}$. We have a diagram
\begin{equation}
    \begin{tikzcd}
    &x_2\arrow[dr]\arrow[drr]&&\\
    y_2\arrow[ur]\arrow[ddrr,bend left]\arrow[dr]\arrow[ddr]&&x_2'\arrow[r]\arrow[d]&x_2^*\arrow[d]\\
    &y_2'\arrow[ur]\arrow[dr]\arrow[d]&x_1\arrow[r]&x_1^*\\
    &y_1\arrow[dr]\arrow[ur]&y_2^*\arrow[uur]\arrow[d]&\\
    &&y_1^*.\arrow[uur]&
    \end{tikzcd}
\end{equation}
To show the commutativity of the diagram formed by $y_2$, $x_2$, $x_2'$ and $y_2'$, it suffices to show the commutativity of $y_2$, $x_2$, $x_2^*$ and $y_2^*$. This follows from the definition of $2$-limits.
\end{proof}

\bibliographystyle{amsalpha}
\bibliography{references}

@incollection {Del79,
    AUTHOR = {Deligne, P.},
     TITLE = {Vari\'et\'es de {S}himura: interpr\'etation modulaire, et
              techniques de construction de mod\`eles canoniques},
 BOOKTITLE = {Automorphic forms, representations and {$L$}-functions
              ({P}roc. {S}ympos. {P}ure {M}ath., {O}regon {S}tate {U}niv.,
              {C}orvallis, {O}re., 1977), {P}art 2},
    SERIES = {Proc. Sympos. Pure Math.},
    VOLUME = {XXXIII},
     PAGES = {247--289},
 PUBLISHER = {Amer. Math. Soc., Providence, RI},
      YEAR = {1979},
      ISBN = {0-8218-1437-0},
   MRCLASS = {10D20 (14D20 14G25 14K15)},
  MRNUMBER = {546620},
MRREVIEWER = {James\ Milne},
}

@article {KP15,
    AUTHOR = {Kisin, M. and Pappas, G.},
     TITLE = {Integral models of {S}himura varieties with parahoric level
              structure},
   JOURNAL = {Publ. Math. Inst. Hautes \'{E}tudes Sci.},
  FJOURNAL = {Publications Math\'{e}matiques. Institut de Hautes \'{E}tudes
              Scientifiques},
    VOLUME = {128},
      YEAR = {2018},
     PAGES = {121--218},
      ISSN = {0073-8301},
   MRCLASS = {11G18},
  MRNUMBER = {3905466},
MRREVIEWER = {Martin Orr},
       DOI = {10.1007/s10240-018-0100-0},
       URL = {https://doi.org/10.1007/s10240-018-0100-0},
}

@book {Lan13,
    AUTHOR = {Lan, K.-W.},
     TITLE = {Arithmetic compactifications of {PEL}-type {S}himura
              varieties},
    SERIES = {London Mathematical Society Monographs Series},
    VOLUME = {36},
 PUBLISHER = {Princeton University Press, Princeton, NJ},
      YEAR = {2013},
     PAGES = {xxvi+561},
      ISBN = {978-0-691-15654-5},
   MRCLASS = {14G35 (11G18 14D23 14M27)},
  MRNUMBER = {3186092},
MRREVIEWER = {Rolf\ Berndt},
       DOI = {10.1515/9781400846016},
       URL = {https://doi.org/10.1515/9781400846016},
NOTE = {\href{https://www.kwlan.org/articles/cpt-PEL-type-thesis-revision.pdf}{Link}}
}

@article {Mad19,
    AUTHOR = {Madapusi Pera, K.},
     TITLE = {Toroidal compactifications of integral models of {S}himura
              varieties of {H}odge type},
   JOURNAL = {Ann. Sci. \'{E}c. Norm. Sup\'{e}r. (4)},
  FJOURNAL = {Annales Scientifiques de l'\'{E}cole Normale Sup\'{e}rieure.
              Quatri\`{e}me S\'{e}rie},
    VOLUME = {52},
      YEAR = {2019},
    NUMBER = {2},
     PAGES = {393--514},
      ISSN = {0012-9593,1873-2151},
   MRCLASS = {11G18 (14G35)},
  MRNUMBER = {3948111},
MRREVIEWER = {Brandon\ Levin},
       DOI = {10.24033/asens.2391},
       URL = {https://doi.org/10.24033/asens.2391},
}

@book {Ogu18,
    AUTHOR = {Ogus, A.},
     TITLE = {Lectures on logarithmic algebraic geometry},
    SERIES = {Cambridge Studies in Advanced Mathematics},
    VOLUME = {178},
 PUBLISHER = {Cambridge University Press, Cambridge},
      YEAR = {2018},
     PAGES = {xviii+539},
      ISBN = {978-1-107-18773-3},
   MRCLASS = {14D06 (14A20 14M25)},
  MRNUMBER = {3838359},
MRREVIEWER = {Howard\ M.\ Thompson},
       DOI = {10.1017/9781316941614},
       URL = {https://doi.org/10.1017/9781316941614},
}

@book {Pin89,
    AUTHOR = {Pink, R.},
     TITLE = {Arithmetical compactification of mixed {S}himura varieties},
    SERIES = {Bonner Mathematische Schriften [Bonn Mathematical
              Publications]},
    VOLUME = {209},
      NOTE = {Dissertation, Rheinische Friedrich-Wilhelms-Universit\"at
              Bonn, Bonn, 1989},
 PUBLISHER = {Universit\"at Bonn, Mathematisches Institut, Bonn},
      YEAR = {1990},
     PAGES = {xviii+340},
   MRCLASS = {11G18 (14G35)},
  MRNUMBER = {1128753},
MRREVIEWER = {Min\ Ho\ Lee},
}

@misc{Mao25,
      title={Compactifications of subschemes of integral models of {H}odge-type {S}himura Varieties with Parahoric level structures}, 
      author={Mao, S.},
      year={2025},
      eprint={2504.08574},
      archivePrefix={arXiv},
      primaryClass={math.NT},
      note={arXiv:\href{https://arxiv.org/abs/2504.08574}{2504.08574}}
}

@misc{Mao25b,
      title={Boundary structures of integral models of {H}odge-type {S}himura Varieties}, 
      author={Mao, S.},
      year={2025},
      eprint={2504.13911},
      archivePrefix={arXiv},
      primaryClass={math.NT},
      note={arXiv:\href{https://arxiv.org/abs/2504.13911}{2504.13911}}, 
}

@article {DLLZ,
    AUTHOR = {Diao, H. and Lan, K.-W. and Liu, R. and Zhu,
              X.},
     TITLE = {Logarithmic {R}iemann-{H}ilbert correspondences for rigid
              varieties},
   JOURNAL = {J. Amer. Math. Soc.},
  FJOURNAL = {Journal of the American Mathematical Society},
    VOLUME = {36},
      YEAR = {2023},
    NUMBER = {2},
     PAGES = {483--562},
      ISSN = {0894-0347,1088-6834},
   MRCLASS = {14F40 (14D07 14F30 14G22 14G35)},
  MRNUMBER = {4536903},
MRREVIEWER = {Fumio\ Hazama},
       DOI = {10.1090/jams/1002},
       URL = {https://doi.org/10.1090/jams/1002},
}

@incollection {DLLZ23,
    AUTHOR = {Diao, H. and Lan, K.-W. and Liu, R. and Zhu,
              X.},
     TITLE = {Logarithmic adic spaces: some foundational results},
 BOOKTITLE = {{$p$}-adic {H}odge theory, singular varieties, and non-abelian
              aspects},
    SERIES = {Simons Symp.},
     PAGES = {65--182},
 PUBLISHER = {Springer, Cham},
      YEAR = {2023},
      ISBN = {978-3-031-21549-0; 978-3-031-21550-6},
   MRCLASS = {14A21},
  MRNUMBER = {4592580},
MRREVIEWER = {Huy\ Dang},
       DOI = {10.1007/978-3-031-21550-6\_3},
       URL = {https://doi.org/10.1007/978-3-031-21550-6_3},
}

@article{Sch13, title={$p$-ADIC {H}ODGE THEORY FOR RIGID-ANALYTIC VARIETIES}, volume={1}, DOI={10.1017/fmp.2013.1}, journal={Forum of Mathematics, Pi}, publisher={Cambridge University Press}, author={Scholze, P.}, year={2013}, pages={e1}}

@book {SW20,
    AUTHOR = {Scholze, P. and Weinstein, J.},
     TITLE = {Berkeley lectures on {$p$}-adic geometry},
    SERIES = {Annals of Mathematics Studies},
    VOLUME = {207},
 PUBLISHER = {Princeton University Press, Princeton, NJ},
      YEAR = {2020},
     PAGES = {x+250},
      ISBN = {978-0-691-20209-9; 978-0-691-20208-2; 978-0-691-20215-0},
   MRCLASS = {14G45 (14A15 14F30 14G22 14G35 14M15)},
  MRNUMBER = {4446467},
MRREVIEWER = {Lance\ Edward\ Miller},
}

@book {FC90,
    AUTHOR = {Faltings, G. and Chai, C.-L.},
     TITLE = {Degeneration of abelian varieties},
    SERIES = {Ergebnisse der Mathematik und ihrer Grenzgebiete (3) [Results
              in Mathematics and Related Areas (3)]},
    VOLUME = {22},
      NOTE = {With an appendix by David Mumford},
 PUBLISHER = {Springer-Verlag, Berlin},
      YEAR = {1990},
     PAGES = {xii+316},
      ISBN = {3-540-52015-5},
   MRCLASS = {14K10 (11G10 14D20 14K25)},
  MRNUMBER = {1083353},
MRREVIEWER = {Min\ Ho\ Lee},
       DOI = {10.1007/978-3-662-02632-8},
       URL = {https://doi.org/10.1007/978-3-662-02632-8},
}

@article {BT84,
    AUTHOR = {Bruhat, F. and Tits, J.},
     TITLE = {Groupes r\'eductifs sur un corps local. {II}. {S}ch\'emas en
              groupes. {E}xistence d'une donn\'ee radicielle valu\'ee},
   JOURNAL = {Inst. Hautes \'Etudes Sci. Publ. Math.},
  FJOURNAL = {Institut des Hautes \'Etudes Scientifiques. Publications
              Math\'ematiques},
    NUMBER = {60},
      YEAR = {1984},
     PAGES = {197--376},
      ISSN = {0073-8301,1618-1913},
   MRCLASS = {20G25 (14L15)},
  MRNUMBER = {756316},
MRREVIEWER = {James\ E.\ Humphreys},
       URL = {http://www.numdam.org/item?id=PMIHES_1984__60__5_0},
}

@book {AMRT10,
    AUTHOR = {Ash, A. and Mumford, D. and Rapoport, M. and Tai,
              Y.-S.},
     TITLE = {Smooth compactifications of locally symmetric varieties},
    SERIES = {Cambridge Mathematical Library},
   EDITION = {Second},
      NOTE = {With the collaboration of Peter Scholze},
 PUBLISHER = {Cambridge University Press, Cambridge},
      YEAR = {2010},
     PAGES = {x+230},
      ISBN = {978-0-521-73955-9},
   MRCLASS = {14M27 (32J05 32M15)},
  MRNUMBER = {2590897},
       DOI = {10.1017/CBO9780511674693},
       URL = {https://doi.org/10.1017/CBO9780511674693},
}

@misc{stacks-project,
  author       = {The {Stacks project authors}},
  title        = {The Stacks project},
  note = {\href{link}{https://stacks.math.columbia.edu}},
  year         = {2025},
}

@BOOK{EGA4,
    AUTHOR = "Grothendieck, A. and Dieudonn{\'e}, J.",
    TITLE = "{\'E}l{\'e}ments de g{\'e}om{\'e}trie alg{\'e}brique {IV}",
    PUBLISHER = "Institute des {H}autes {\'E}tudes {S}cientifiques.",
    YEAR = "1964-1967",
    SERIES = "Publications {M}ath{\'e}matiques",
    VOLUME = "20, 24, 28, 32"
}

@article {Lan16b,
    AUTHOR = {Lan, K.-W.},
     TITLE = {Compactifications of {PEL}-type {S}himura varieties in
              ramified characteristics},
   JOURNAL = {Forum Math. Sigma},
  FJOURNAL = {Forum of Mathematics. Sigma},
    VOLUME = {4},
      YEAR = {2016},
     PAGES = {Paper No. e1, 98},
      ISSN = {2050-5094},
   MRCLASS = {11G18 (11G15 14D06)},
  MRNUMBER = {3482277},
MRREVIEWER = {Matteo\ Longo},
       DOI = {10.1017/fms.2015.31},
       URL = {https://doi.org/10.1017/fms.2015.31},
}

@article {PR24,
    AUTHOR = {Pappas, G. and Rapoport, M.},
     TITLE = {{$p$}-adic shtukas and the theory of global and local
              {S}himura varieties},
   JOURNAL = {Camb. J. Math.},
  FJOURNAL = {Cambridge Journal of Mathematics},
    VOLUME = {12},
      YEAR = {2024},
    NUMBER = {1},
     PAGES = {1--164},
      ISSN = {2168-0930,2168-0949},
   MRCLASS = {11G18 (14G35 14G45)},
  MRNUMBER = {4701491},
}

@article {Har89,
    AUTHOR = {Harris, M.},
     TITLE = {Functorial properties of toroidal compactifications of locally
              symmetric varieties},
   JOURNAL = {Proc. London Math. Soc. (3)},
  FJOURNAL = {Proceedings of the London Mathematical Society. Third Series},
    VOLUME = {59},
      YEAR = {1989},
    NUMBER = {1},
     PAGES = {1--22},
      ISSN = {0024-6115,1460-244X},
   MRCLASS = {11G18 (11F55 14L32 32M15)},
  MRNUMBER = {997249},
MRREVIEWER = {Gerd\ Faltings},
       DOI = {10.1112/plms/s3-59.1.1},
       URL = {https://doi.org/10.1112/plms/s3-59.1.1},
}

@article {HZ01,
    AUTHOR = {Harris, M. and Zucker, S.},
     TITLE = {Boundary cohomology of {S}himura varieties. {III}. {C}oherent
              cohomology on higher-rank boundary strata and applications to
              {H}odge theory},
   JOURNAL = {M\'{e}m. Soc. Math. Fr. (N.S.)},
  FJOURNAL = {M\'{e}moires de la Soci\'{e}t\'{e} Math\'{e}matique de France.
              Nouvelle S\'{e}rie},
    NUMBER = {85},
      YEAR = {2001},
     PAGES = {vi+116},
      ISSN = {0249-633X,2275-3230},
   MRCLASS = {14G35 (11G18 14D07)},
  MRNUMBER = {1850830},
MRREVIEWER = {James\ Milne},
       DOI = {10.24033/msmf.398},
       URL = {https://doi.org/10.24033/msmf.398},
}

@article {DY25,
    AUTHOR = {Daniels, P. and Youcis, A.},
     TITLE = {Canonical integral models for {S}himura varieties of abelian
              type},
   JOURNAL = {Forum Math. Sigma},
  FJOURNAL = {Forum of Mathematics. Sigma},
    VOLUME = {13},
      YEAR = {2025},
     PAGES = {Paper No. e69},
      ISSN = {2050-5094},
   MRCLASS = {11G18 (11R39 14G35)},
  MRNUMBER = {4888034},
       DOI = {10.1017/fms.2025.27},
       URL = {https://doi.org/10.1017/fms.2025.27},
}

@misc{DvHKZ24ig,
      title={{I}gusa Stacks and the Cohomology of {S}himura Varieties}, 
      author={Daniels, P. and van Hoften, P. and Kim, D. and Zhang, M.},
      year={2024},
      eprint={2408.01348},
      archivePrefix={arXiv},
      primaryClass={math.NT},
      note={arXiv:\href{https://arxiv.org/abs/2408.01348}{2408.01348}}, 
}

@article{daniels2024conjecture,
    AUTHOR = {Daniels, P. and van Hoften, P. and Kim, D. and
              Zhang, M.},
     TITLE = {On a conjecture of {P}appas and {R}apoport},
   JOURNAL = {Math. Ann.},
  FJOURNAL = {Mathematische Annalen},
    VOLUME = {395},
      YEAR = {2026},
    NUMBER = {2},
     PAGES = {31},
      ISSN = {0025-5831,1432-1807},
   MRCLASS = {14},
  MRNUMBER = {5061066},
       DOI = {10.1007/s00208-026-03450-4},
       URL = {https://doi.org/10.1007/s00208-026-03450-4},
}

@misc{KPZ24,
      title={Integral models of {S}himura varieties with parahoric level structure, {II}}, 
      author={Kisin, M. and Pappas, G. and Zhou, R.},
      year={2024},
      eprint={2409.03689},
      archivePrefix={arXiv},
      primaryClass={math.NT},
      note={arXiv:\href{https://arxiv.org/abs/2409.03689}{2409.03689}, to appear in Forum of Math. Pi},
}

@article {LS18i,
    AUTHOR = {Lan, K.-W. and Stroh, B.},
     TITLE = {Nearby cycles of automorphic \'etale sheaves},
   JOURNAL = {Compos. Math.},
  FJOURNAL = {Compositio Mathematica},
    VOLUME = {154},
      YEAR = {2018},
    NUMBER = {1},
     PAGES = {80--119},
      ISSN = {0010-437X,1570-5846},
   MRCLASS = {11G18 (11F75 11G15 14F20 14G35)},
  MRNUMBER = {3719245},
MRREVIEWER = {Giovanni\ Rosso},
       DOI = {10.1112/S0010437X1700745X},
       URL = {https://doi.org/10.1112/S0010437X1700745X},
}

@article {Pin92,
    AUTHOR = {Pink, R.},
     TITLE = {On {$l$}-adic sheaves on {S}himura varieties and their higher
              direct images in the {B}aily-{B}orel compactification},
   JOURNAL = {Math. Ann.},
  FJOURNAL = {Mathematische Annalen},
    VOLUME = {292},
      YEAR = {1992},
    NUMBER = {2},
     PAGES = {197--240},
      ISSN = {0025-5831,1432-1807},
   MRCLASS = {11G18 (11F75 14F20)},
  MRNUMBER = {1149032},
MRREVIEWER = {Min\ Ho\ Lee},
       DOI = {10.1007/BF01444618},
       URL = {https://doi.org/10.1007/BF01444618},
}

@article {liu2017rigidity,
    AUTHOR = {Liu, R. and Zhu, X.},
     TITLE = {Rigidity and a {R}iemann-{H}ilbert correspondence for
              {$p$}-adic local systems},
   JOURNAL = {Invent. Math.},
  FJOURNAL = {Inventiones Mathematicae},
    VOLUME = {207},
      YEAR = {2017},
    NUMBER = {1},
     PAGES = {291--343},
      ISSN = {0020-9910,1432-1297},
   MRCLASS = {14G22 (14G35 14J60)},
  MRNUMBER = {3592758},
MRREVIEWER = {Marco\ A.\ Garuti},
       DOI = {10.1007/s00222-016-0671-7},
       URL = {https://doi.org/10.1007/s00222-016-0671-7},
}

@article {caraiani2017generic,
    AUTHOR = {Caraiani, A. and Scholze, P.},
     TITLE = {On the generic part of the cohomology of compact unitary
              {S}himura varieties},
   JOURNAL = {Ann. of Math. (2)},
  FJOURNAL = {Annals of Mathematics. Second Series},
    VOLUME = {186},
      YEAR = {2017},
    NUMBER = {3},
     PAGES = {649--766},
      ISSN = {0003-486X,1939-8980},
   MRCLASS = {11F75 (11G18 11R23 14G35)},
  MRNUMBER = {3702677},
MRREVIEWER = {Nguy\cftil en Qu\^oc Th\'ang},
       DOI = {10.4007/annals.2017.186.3.1},
       URL = {https://doi.org/10.4007/annals.2017.186.3.1},
}

@article{richarz2020basics,
  title={Basics on affine {G}rassmannians},
  author={Richarz, T.},
  journal={Note},
  year={2019},
  note={Available from \href{https://www.mathematik.tu-darmstadt.de/media/algebra/homepages/richarz/Notes_on_affine_Grassmannians.pdf}{Link}},
}

@misc{scholze2017etale,
      title={{\'E}tale cohomology of diamonds}, 
      author={Scholze, P.},
      year={2026},
      eprint={1709.07343},
      archivePrefix={arXiv},
      primaryClass={math.AG},
      note={arXiv:\href{https://arxiv.org/abs/1709.07343}{1709.07343}, to appear in Ast.}, 
}

@article {KY25,
    AUTHOR = {Koshikawa, T. and Yao, Z.},
     TITLE = {Logarithmic prismatic cohomology {II}},
   JOURNAL = {Adv. Math.},
  FJOURNAL = {Advances in Mathematics},
    VOLUME = {479},
      YEAR = {2025},
     PAGES = {Paper No. 110446},
      ISSN = {0001-8708,1090-2082},
   MRCLASS = {14A21 (14F20 14F30 14F40 14G20)},
  MRNUMBER = {4937354},
       DOI = {10.1016/j.aim.2025.110446},
       URL = {https://doi.org/10.1016/j.aim.2025.110446},
}

@misc{anschutz2022p,
      title={On the $p$-adic theory of local models}, 
      author={Anschütz, J. and Gleason, I. and Lourenço, J. and Richarz, T.},
      year={2022},
      eprint={2201.01234},
      archivePrefix={arXiv},
      primaryClass={math.AG},
      note={arXiv:\href{https://arxiv.org/abs/2201.01234}{2201.01234}, to appear in Ann. of Math.}, 
}

@misc{Wu25,
      title={Arithmetic compactifications of integral models of {S}himura varieties of abelian type}, 
      author={Wu, P.},
      year={2025},
      eprint={2505.09135},
      archivePrefix={arXiv},
      primaryClass={math.NT},
      url={https://arxiv.org/abs/2505.09135}, 
      note={arXiv:\href{https://arxiv.org/abs/2505.09135}{2505.09135}},
}

@incollection {Kat89,
    AUTHOR = {Kato, K.},
     TITLE = {Logarithmic structures of {F}ontaine-{I}llusie},
 BOOKTITLE = {Algebraic analysis, geometry, and number theory ({B}altimore,
              {MD}, 1988)},
     PAGES = {191--224},
 PUBLISHER = {Johns Hopkins Univ. Press, Baltimore, MD},
      YEAR = {1989},
      ISBN = {0-8018-3841-X},
   MRCLASS = {14F30 (14G20)},
  MRNUMBER = {1463703},
MRREVIEWER = {Adolfo\ Quir\'os},
}

@inproceedings{lan2018compactifications,
    AUTHOR = {Lan, K.-W. and Stroh, B.},
     TITLE = {Compactifications of subschemes of integral models of
              {S}himura varieties},
   JOURNAL = {Forum Math. Sigma},
  FJOURNAL = {Forum of Mathematics. Sigma},
    VOLUME = {6},
      YEAR = {2018},
     PAGES = {Paper No. e18, 105},
      ISSN = {2050-5094},
   MRCLASS = {11G18 (11F75 11G15 14G35)},
  MRNUMBER = {3859178},
MRREVIEWER = {Pietro\ Mercuri},
       DOI = {10.1017/fms.2018.20},
       URL = {https://doi.org/10.1017/fms.2018.20},
}

@article{klevdal2023compatibility,
  title={Compatibility of canonical $\ell$-adic local systems on {S}himura varieties},
  author={Klevdal, C. and Patrikis, S.},
  journal={arXiv preprint},
  year={2023},
  note={arXiv:\href{https://arxiv.org/abs/2303.03863}{2303.03863}},
}

@article {Ols03,
    AUTHOR = {Olsson, M.},
     TITLE = {Logarithmic geometry and algebraic stacks},
   JOURNAL = {Ann. Sci. \'Ecole Norm. Sup. (4)},
  FJOURNAL = {Annales Scientifiques de l'\'Ecole Normale Sup\'erieure.
              Quatri\`eme S\'erie},
    VOLUME = {36},
      YEAR = {2003},
    NUMBER = {5},
     PAGES = {747--791},
      ISSN = {0012-9593},
   MRCLASS = {14D20 (14A20)},
  MRNUMBER = {2032986},
MRREVIEWER = {Ivan\ S.\ Kausz},
       DOI = {10.1016/j.ansens.2002.11.001},
       URL = {https://doi.org/10.1016/j.ansens.2002.11.001},
}

@book {Hub96,
    AUTHOR = {Huber, R.},
     TITLE = {\'Etale cohomology of rigid analytic varieties and adic
              spaces},
    SERIES = {Aspects of Mathematics},
    VOLUME = {E30},
 PUBLISHER = {Friedr. Vieweg \& Sohn, Braunschweig},
      YEAR = {1996},
     PAGES = {x+450},
      ISBN = {3-528-06794-2},
   MRCLASS = {14G22 (14F20)},
  MRNUMBER = {1734903},
MRREVIEWER = {Lorenzo\ Ramero},
       DOI = {10.1007/978-3-663-09991-8},
       URL = {https://doi.org/10.1007/978-3-663-09991-8},
}

@article {Hub94,
    AUTHOR = {Huber, R.},
     TITLE = {A generalization of formal schemes and rigid analytic varieties},
   JOURNAL = {Math. Z.},
  FJOURNAL = {Mathematische Zeitschrift},
    VOLUME = {217},
      YEAR = {1994},
    NUMBER = {4},
     PAGES = {513--551},
      ISSN = {0025-5874,1432-1823},
   MRCLASS = {14A20 (13A18 14L15 32P05)},
  MRNUMBER = {1306024},
MRREVIEWER = {W.\ Bartenwerfer},
       DOI = {10.1007/BF02571959},
       URL = {https://doi.org/10.1007/BF02571959},
}

@article{heuer2022g,
    AUTHOR = {Heuer, B.},
     TITLE = {{$G$}-torsors on perfectoid spaces},
   JOURNAL = {\'Epijournal G\'eom. Alg\'ebrique},
  FJOURNAL = {\'Epijournal de G\'eom\'etrie Alg\'ebrique. EPIGA},
    VOLUME = {10},
      YEAR = {2026},
     PAGES = {Art. 2, 33},
      ISSN = {2491-6765},
   MRCLASS = {14G45 (14F20 14G22)},
  MRNUMBER = {5046468},
       DOI = {10.46298/epiga.2026.13796},
       URL = {https://doi.org/10.46298/epiga.2026.13796},
}

@article{zhu2017affine,
    AUTHOR = {Zhu, X.},
     TITLE = {Affine {G}rassmannians and the geometric {S}atake in mixed
              characteristic},
   JOURNAL = {Ann. of Math. (2)},
  FJOURNAL = {Annals of Mathematics. Second Series},
    VOLUME = {185},
      YEAR = {2017},
    NUMBER = {2},
     PAGES = {403--492},
      ISSN = {0003-486X,1939-8980},
   MRCLASS = {14D24 (14L35 14M15 20G25)},
  MRNUMBER = {3612002},
MRREVIEWER = {Rolf\ Berndt},
       DOI = {10.4007/annals.2017.185.2.2},
       URL = {https://doi.org/10.4007/annals.2017.185.2.2},
}

@article {bhatt2017projectivity,
    AUTHOR = {Bhatt, B. and Scholze, P.},
     TITLE = {Projectivity of the {W}itt vector affine {G}rassmannian},
   JOURNAL = {Invent. Math.},
  FJOURNAL = {Inventiones Mathematicae},
    VOLUME = {209},
      YEAR = {2017},
    NUMBER = {2},
     PAGES = {329--423},
      ISSN = {0020-9910,1432-1297},
   MRCLASS = {14F05 (14M15 19G12)},
  MRNUMBER = {3674218},
MRREVIEWER = {Marc-Hubert\ Nicole},
       DOI = {10.1007/s00222-016-0710-4},
       URL = {https://doi.org/10.1007/s00222-016-0710-4},
}

@article{blakestad2018perfectoid,
    AUTHOR = {Blakestad, C. and Gvirtz, D. and Heuer, B. and
              Shchedrina, D. and Shimizu, K. and Wear, P. and Yao,
              Z.},
     TITLE = {Perfectoid covers of abelian varieties},
   JOURNAL = {Math. Res. Lett.},
  FJOURNAL = {Mathematical Research Letters},
    VOLUME = {29},
      YEAR = {2022},
    NUMBER = {3},
     PAGES = {631--662},
      ISSN = {1073-2780,1945-001X},
   MRCLASS = {14G45},
  MRNUMBER = {4516034},
MRREVIEWER = {Lance\ Edward\ Miller},
       DOI = {10.4310/mrl.2022.v29.n3.a2},
       URL = {https://doi.org/10.4310/mrl.2022.v29.n3.a2},
}

@misc{KL19,
      title={Relative p-adic {H}odge theory, {II}: {I}mperfect period rings}, 
      author={Kedlaya, K. S. and Liu, R.},
      year={2019},
      eprint={1602.06899},
      archivePrefix={arXiv},
      primaryClass={math.NT},
      url={https://arxiv.org/abs/1602.06899}, 
 note={arXiv:\href{https://arxiv.org/abs/1602.06899}{1602.06899}}, 
}

@article{shen2021ekor,
    AUTHOR = {Shen, X. and Yu, C.-F. and Zhang, C.},
     TITLE = {E{KOR} strata for {S}himura varieties with parahoric level
              structure},
   JOURNAL = {Duke Math. J.},
  FJOURNAL = {Duke Mathematical Journal},
    VOLUME = {170},
      YEAR = {2021},
    NUMBER = {14},
     PAGES = {3111--3236},
      ISSN = {0012-7094,1547-7398},
   MRCLASS = {14G35 (11G18)},
  MRNUMBER = {4319228},
MRREVIEWER = {Su-ion\ Ih},
       DOI = {10.1215/00127094-2021-0047},
       URL = {https://doi.org/10.1215/00127094-2021-0047},
}

@article{anschutz2022extending,
    AUTHOR = {Ansch\"utz, J.},
     TITLE = {Extending torsors on the punctured {${\rm Spec}(A_{\inf})$}},
   JOURNAL = {J. Reine Angew. Math.},
  FJOURNAL = {Journal f\"ur die Reine und Angewandte Mathematik. [Crelle's
              Journal]},
    VOLUME = {783},
      YEAR = {2022},
     PAGES = {227--268},
      ISSN = {0075-4102,1435-5345},
   MRCLASS = {20G05 (14M15)},
  MRNUMBER = {4373246},
MRREVIEWER = {Nguy\cftil en Qu\^oc Th\'ang},
       DOI = {10.1515/crelle-2021-0077},
       URL = {https://doi.org/10.1515/crelle-2021-0077},
}

@article{he2024affine,
    AUTHOR = {He, X. and Nie, S. and Yu, Q.},
     TITLE = {Affine {D}eligne-{L}usztig varieties with finite {C}oxeter
              parts},
   JOURNAL = {Algebra Number Theory},
  FJOURNAL = {Algebra \& Number Theory},
    VOLUME = {18},
      YEAR = {2024},
    NUMBER = {9},
     PAGES = {1681--1714},
      ISSN = {1937-0652,1944-7833},
   MRCLASS = {11G25 (20F55 20G25)},
  MRNUMBER = {4856606},
MRREVIEWER = {Lei\ Yang},
       DOI = {10.2140/ant.2024.18.1681},
       URL = {https://doi.org/10.2140/ant.2024.18.1681},
}

@article {gleason2025specialization,
    AUTHOR = {Gleason, I.},
     TITLE = {Specialization maps for {S}cholze's category of diamonds},
   JOURNAL = {Math. Ann.},
  FJOURNAL = {Mathematische Annalen},
    VOLUME = {391},
      YEAR = {2025},
    NUMBER = {2},
     PAGES = {1611--1679},
      ISSN = {0025-5831,1432-1807},
   MRCLASS = {14G45},
  MRNUMBER = {4853001},
       DOI = {10.1007/s00208-024-02952-3},
       URL = {https://doi.org/10.1007/s00208-024-02952-3},
}

@article{hamacher2025point,
    AUTHOR = {Hamacher, P. and Kim, W.},
     TITLE = {Point counting on {I}gusa varieties for function fields},
   JOURNAL = {Adv. Math.},
  FJOURNAL = {Advances in Mathematics},
    VOLUME = {480},
      YEAR = {2025},
     PAGES = {Paper No. 110517, 88},
      ISSN = {0001-8708,1090-2082},
   MRCLASS = {14G35 (14H60)},
  MRNUMBER = {4957546},
       DOI = {10.1016/j.aim.2025.110517},
}

@article{rapoport1996classification,
    AUTHOR = {Rapoport, M. and Richartz, M.},
     TITLE = {On the classification and specialization of {$F$}-isocrystals with additional structure},
   JOURNAL = {Compositio Math.},
  FJOURNAL = {Compositio Mathematica},
    VOLUME = {103},
      YEAR = {1996},
    NUMBER = {2},
     PAGES = {153--181},
      ISSN = {0010-437X,1570-5846},
   MRCLASS = {14F30 (22E50)},
  MRNUMBER = {1411570},
MRREVIEWER = {Abdellah\ Mokrane},
       URL = {http://www.numdam.org/item?id=CM_1996__103_2_153_0},
}

@misc{yao2023mathbb,
  title={{$\bb{Z}_p$}-lattices in semistable {G}alois representations},
  author={Yao, Z.},
  eprint={2308.15468},
  year={2023},
  note={arXiv:\href{https://arxiv.org/abs/2308.15468}{2308.15468}}
}

@article{fargues2021geometrization,
  title={Geometrization of the local {L}anglands correspondence},
  author={Fargues, L. and Scholze, P.},
  journal={arXiv preprint},
  year={2021},
  note={\href{https://arxiv.org/2102.13459}{2102.13459}, to appear in Ast.}
}

@article{gleason2024tubular,
    AUTHOR = {Gleason, I. and Louren\c co, J.},
     TITLE = {Tubular neighborhoods of local models},
   JOURNAL = {Duke Math. J.},
  FJOURNAL = {Duke Mathematical Journal},
    VOLUME = {173},
      YEAR = {2024},
    NUMBER = {4},
     PAGES = {723--743},
      ISSN = {0012-7094,1547-7398},
   MRCLASS = {14G45 (11G18 14G35)},
  MRNUMBER = {4734553},
       DOI = {10.1215/00127094-2023-0028},
       URL = {https://doi.org/10.1215/00127094-2023-0028},
}

@article{pappas2022integral,
    AUTHOR = {Pappas, G. and Rapoport, M.},
     TITLE = {On integral local {S}himura varieties},
   JOURNAL = {J. Inst. Math. Jussieu},
  FJOURNAL = {Journal of the Institute of Mathematics of Jussieu. JIMJ.
              Journal de l'Institut de Math\'ematiques de Jussieu},
    VOLUME = {25},
      YEAR = {2026},
    NUMBER = {1},
     PAGES = {375--443},
      ISSN = {1474-7480,1475-3030},
   MRCLASS = {11G18 (14G35)},
  MRNUMBER = {5018888},
       DOI = {10.1017/S1474748025101345},
       URL = {https://doi.org/10.1017/S1474748025101345},
}

@article{xiao2017cycles,
  title={Cycles on {S}himura varieties via geometric {S}atake},
  author={Xiao, L. and Zhu, X.},
  journal={arXiv preprint},
  year={2017},
  note={arXiv:\href{https://arxiv.org/abs/1707.05700}{1707.05700}}
}

@article{kottwitz1997isocrystals,
    AUTHOR = {Kottwitz, R. E.},
     TITLE = {Isocrystals with additional structure. {II}},
   JOURNAL = {Compositio Math.},
  FJOURNAL = {Compositio Mathematica},
    VOLUME = {109},
      YEAR = {1997},
    NUMBER = {3},
     PAGES = {255--339},
      ISSN = {0010-437X,1570-5846},
   MRCLASS = {20G25 (11S25 14F30 14L05)},
  MRNUMBER = {1485921},
MRREVIEWER = {Guy\ Rousseau},
       DOI = {10.1023/A:1000102604688},
       URL = {https://doi.org/10.1023/A:1000102604688},
}

@article{anschutz2019reductive,
    AUTHOR = {Ansch\"utz, J.},
     TITLE = {Reductive group schemes over the {F}argues-{F}ontaine curve},
   JOURNAL = {Math. Ann.},
  FJOURNAL = {Mathematische Annalen},
    VOLUME = {374},
      YEAR = {2019},
    NUMBER = {3-4},
     PAGES = {1219--1260},
      ISSN = {0025-5831,1432-1807},
   MRCLASS = {11S31 (11S37 14H60)},
  MRNUMBER = {3985110},
MRREVIEWER = {Giovanni\ Rosso},
       DOI = {10.1007/s00208-018-1785-6},
       URL = {https://doi.org/10.1007/s00208-018-1785-6},
}

@article{gleason2023meromorphic,
  title={Meromorphic vector bundles on the {F}argues--{F}ontaine curve},
  author={Gleason, I. and Ivanov, A. B.},
  journal={arXiv preprint},
  year={2023},
  note={arXiv:\href{https://arxiv.org/abs/2307.00887}{2307.00887}, to appear in JEMS},
}

@article{zhang2023pel,
    AUTHOR = {Zhang, M.},
     TITLE = {A {PEL}-type {I}gusa stack and the {$p$}-adic geometry of
              {S}himura varieties},
   JOURNAL = {Camb. J. Math.},
  FJOURNAL = {Cambridge Journal of Mathematics},
    VOLUME = {14},
      YEAR = {2026},
    NUMBER = {2},
     PAGES = {377--486},
      ISSN = {2168-0930,2168-0949},
   MRCLASS = {14 (11)},
  MRNUMBER = {5073031},
       DOI = {10.4310/cjm.260514231850},
       URL = {https://doi.org/10.4310/cjm.260514231850},
}

@article{pappas2008twisted,
    AUTHOR = {Pappas, G. and Rapoport, M.},
     TITLE = {Twisted loop groups and their affine flag varieties},
      NOTE = {With an appendix by T. Haines and Rapoport},
   JOURNAL = {Adv. Math.},
  FJOURNAL = {Advances in Mathematics},
    VOLUME = {219},
      YEAR = {2008},
    NUMBER = {1},
     PAGES = {118--198},
      ISSN = {0001-8708,1090-2082},
   MRCLASS = {22E67 (14M15 17B67 20G25)},
  MRNUMBER = {2435422},
MRREVIEWER = {Dmitry\ A.\ Timash\"ev},
       DOI = {10.1016/j.aim.2008.04.006},
       URL = {https://doi.org/10.1016/j.aim.2008.04.006},
}

@article{fargues2020g,
    AUTHOR = {Fargues, L.},
     TITLE = {{$G$}-torseurs en th\'eorie de {H}odge {$p$}-adique},
   JOURNAL = {Compos. Math.},
  FJOURNAL = {Compositio Mathematica},
    VOLUME = {156},
      YEAR = {2020},
    NUMBER = {10},
     PAGES = {2076--2110},
      ISSN = {0010-437X,1570-5846},
   MRCLASS = {11S31 (14F30 14L05 14L24)},
  MRNUMBER = {4179595},
       DOI = {10.1112/s0010437x20007423},
       URL = {https://doi.org/10.1112/s0010437x20007423},
}

@article{viehmann2024newton,
    AUTHOR = {Viehmann, E.},
     TITLE = {On {N}ewton strata in the {$B_{\rm dR}^+$}-{G}rassmannian},
   JOURNAL = {Duke Math. J.},
  FJOURNAL = {Duke Mathematical Journal},
    VOLUME = {173},
      YEAR = {2024},
    NUMBER = {1},
     PAGES = {177--225},
      ISSN = {0012-7094,1547-7398},
   MRCLASS = {11G18 (14G20 14M15)},
  MRNUMBER = {4728690},
MRREVIEWER = {Nguy\cftil en Qu\^oc Th\'ang},
       DOI = {10.1215/00127094-2024-0005},
       URL = {https://doi.org/10.1215/00127094-2024-0005},
}

@article{haines2010satake,
    AUTHOR = {Haines, T. J. and Rostami, S.},
     TITLE = {The {S}atake isomorphism for special maximal parahoric {H}ecke algebras},
   JOURNAL = {Represent. Theory},
  FJOURNAL = {Representation Theory. An Electronic Journal of the American Mathematical Society},
    VOLUME = {14},
      YEAR = {2010},
     PAGES = {264--284},
      ISSN = {1088-4165},
   MRCLASS = {20G25 (11E95 11G18 14G35 22E50)},
  MRNUMBER = {2602034},
MRREVIEWER = {Rainer\ Schulze-Pillot},
       DOI = {10.1090/S1088-4165-10-00370-5},
       URL = {https://doi.org/10.1090/S1088-4165-10-00370-5},
}

@article{bertapelle2019deformations,
    AUTHOR = {Bertapelle, A. and Mazzari, N.},
     TITLE = {On deformations of $1$-motives},
   JOURNAL = {Canad. Math. Bull.},
  FJOURNAL = {Canadian Mathematical Bulletin. Bulletin Canadien de
              Math\'ematiques},
    VOLUME = {62},
      YEAR = {2019},
    NUMBER = {1},
     PAGES = {11--22},
      ISSN = {0008-4395,1496-4287},
   MRCLASS = {14C15 (14L05 14L15)},
  MRNUMBER = {3943763},
MRREVIEWER = {Federico\ Binda},
       DOI = {10.4153/cmb-2017-076-2},
       URL = {https://doi.org/10.4153/cmb-2017-076-2},
}

@article{shen2017perfectoid,
    AUTHOR = {Shen, X.},
     TITLE = {Perfectoid {S}himura varieties of abelian type},
   JOURNAL = {Int. Math. Res. Not. IMRN},
  FJOURNAL = {International Mathematics Research Notices. IMRN},
      YEAR = {2017},
    NUMBER = {21},
     PAGES = {6599--6653},
      ISSN = {1073-7928,1687-0247},
   MRCLASS = {14G35 (11G18 14J28)},
  MRNUMBER = {3719474},
MRREVIEWER = {Shuichiro\ Takeda},
       DOI = {10.1093/imrn/rnw202},
       URL = {https://doi.org/10.1093/imrn/rnw202},
}

@article{kim2018rapoport,
    AUTHOR = {Kim, W.},
     TITLE = {Rapoport-{Z}ink spaces of {H}odge type},
   JOURNAL = {Forum Math. Sigma},
  FJOURNAL = {Forum of Mathematics. Sigma},
    VOLUME = {6},
      YEAR = {2018},
     PAGES = {Paper No. e8, 110},
      ISSN = {2050-5094},
   MRCLASS = {14F30 (11G05 14G35 14L05)},
  MRNUMBER = {3812116},
MRREVIEWER = {Robin\ de Jong},
       DOI = {10.1017/fms.2018.6},
       URL = {https://doi.org/10.1017/fms.2018.6},
}

@article{arasteh2021uniformizing,
    AUTHOR = {Arasteh Rad, E. and Hartl, U.},
     TITLE = {Uniformizing the moduli stacks of global {$G$}-shtukas},
   JOURNAL = {Int. Math. Res. Not. IMRN},
  FJOURNAL = {International Mathematics Research Notices. IMRN},
      YEAR = {2021},
    NUMBER = {21},
     PAGES = {16121--16192},
      ISSN = {1073-7928,1687-0247},
   MRCLASS = {14G35 (14D23)},
  MRNUMBER = {4338216},
MRREVIEWER = {Ariel\ Shnidman},
       DOI = {10.1093/imrn/rnz223},
       URL = {https://doi.org/10.1093/imrn/rnz223},
}

@article{rad2019local,
    AUTHOR = {Arasteh Rad, E. and Habibi, S.},
     TITLE = {Local models for the moduli stacks of global {$\germ
              G$}-shtukas},
   JOURNAL = {Math. Res. Lett.},
  FJOURNAL = {Mathematical Research Letters},
    VOLUME = {26},
      YEAR = {2019},
    NUMBER = {2},
     PAGES = {323--364},
      ISSN = {1073-2780,1945-001X},
   MRCLASS = {14G35 (14D23)},
  MRNUMBER = {3999548},
MRREVIEWER = {David\ Tweedle},
       DOI = {10.4310/MRL.2019.v26.n2.a1},
       URL = {https://doi.org/10.4310/MRL.2019.v26.n2.a1},
}

@article{imai2024tannakian,
  title={A {T}annakian framework for prismatic {$F$}-crystals},
  author={Imai, N. and Kato, H. and Youcis, A.},
  journal={arXiv preprint},
  year={2024},
  note={arXiv:\href{https://arxiv.org/abs/2406.08259}{2406.08259}, to appear in Forum of Math. Sigma}
}

@misc{MY2026,
  title={On canonicity for integral models for Shimura varieties with hyperspecial level},
  author={Madapusi, K. and Youcis, A.},
  journal={preprint},
  year={2026},
  note={Available \href{https://www.keerthimadapusi.com/assets/preprints/syntomic_with_construction.pdf}{here}}
}

@misc{Yan25,
      title={$2$-adic integral models of some {S}himura varieties with parahoric level structure}, 
      author={Yang, J.},
      year={2025},
      eprint={2301.12981},
      archivePrefix={arXiv},
      primaryClass={math.NT},
      url={https://arxiv.org/abs/2301.12981}, 
      note={arXiv:\href{https://arxiv.org/abs/2301.12981}{2301.12981}}
}

@article{Richarz2013,
    AUTHOR = {Richarz, T.},
     TITLE = {Schubert varieties in twisted affine flag varieties and local models},
   JOURNAL = {J. Algebra},
  FJOURNAL = {Journal of Algebra},
    VOLUME = {375},
      YEAR = {2013},
     PAGES = {121--147},
      ISSN = {0021-8693,1090-266X},
   MRCLASS = {14M15},
  MRNUMBER = {2998951},
MRREVIEWER = {Leonardo\ Constantin\ Mihalcea},
       DOI = {10.1016/j.jalgebra.2012.11.013},
       URL = {https://doi.org/10.1016/j.jalgebra.2012.11.013},
}

@article{HainesRicharz2021Test,
    AUTHOR = {Haines, T. J. and Richarz, T.},
     TITLE = {The test function conjecture for parahoric local models},
   JOURNAL = {J. Amer. Math. Soc.},
  FJOURNAL = {Journal of the American Mathematical Society},
    VOLUME = {34},
      YEAR = {2021},
    NUMBER = {1},
     PAGES = {135--218},
      ISSN = {0894-0347,1088-6834},
   MRCLASS = {14G35 (14M15 20G05)},
  MRNUMBER = {4188816},
MRREVIEWER = {Alan\ Koch},
       DOI = {10.1090/jams/955},
       URL = {https://doi.org/10.1090/jams/955},
}

@incollection{oort2001stratification,
    AUTHOR = {Oort, F.},
     TITLE = {A stratification of a moduli space of abelian varieties},
 BOOKTITLE = {Moduli of abelian varieties ({T}exel {I}sland, 1999)},
    SERIES = {Progr. Math.},
    VOLUME = {195},
     PAGES = {345--416},
 PUBLISHER = {Birkh\"auser, Basel},
      YEAR = {2001},
      ISBN = {3-7643-6517-X},
   MRCLASS = {14K10 (14L15)},
  MRNUMBER = {1827027},
MRREVIEWER = {Takashi\ Ichikawa},
       DOI = {10.1007/978-3-0348-8303-0\_13},
       URL = {https://doi.org/10.1007/978-3-0348-8303-0_13},
}

@article{pink2015f,
    AUTHOR = {Pink, R. and Wedhorn, T. and Ziegler, P.},
     TITLE = {{$F$}-zips with additional structure},
   JOURNAL = {Pacific J. Math.},
  FJOURNAL = {Pacific Journal of Mathematics},
    VOLUME = {274},
      YEAR = {2015},
    NUMBER = {1},
     PAGES = {183--236},
      ISSN = {0030-8730,1945-5844},
   MRCLASS = {14F05 (14F40 14G15 18D10)},
  MRNUMBER = {3347958},
MRREVIEWER = {Eva\ Viehmann},
       DOI = {10.2140/pjm.2015.274.183},
       URL = {https://doi.org/10.2140/pjm.2015.274.183},
}

@article{viehmann2013ekedahl,
    AUTHOR = {Viehmann, E. and Wedhorn, T.},
     TITLE = {Ekedahl-{O}ort and {N}ewton strata for {S}himura varieties of
              {PEL} type},
   JOURNAL = {Math. Ann.},
  FJOURNAL = {Mathematische Annalen},
    VOLUME = {356},
      YEAR = {2013},
    NUMBER = {4},
     PAGES = {1493--1550},
      ISSN = {0025-5831,1432-1807},
   MRCLASS = {14G35 (11G18)},
  MRNUMBER = {3072810},
MRREVIEWER = {Marc-Hubert\ Nicole},
       DOI = {10.1007/s00208-012-0892-z},
       URL = {https://doi.org/10.1007/s00208-012-0892-z},
}

@article{zhang2018ekedahl,
    AUTHOR = {Zhang, C.},
     TITLE = {Ekedahl-{O}ort strata for good reductions of {S}himura
              varieties of {H}odge type},
   JOURNAL = {Canad. J. Math.},
  FJOURNAL = {Canadian Journal of Mathematics. Journal Canadien de
              Math\'ematiques},
    VOLUME = {70},
      YEAR = {2018},
    NUMBER = {2},
     PAGES = {451--480},
      ISSN = {0008-414X,1496-4279},
   MRCLASS = {14G35 (11G18)},
  MRNUMBER = {3759007},
MRREVIEWER = {Brandon\ Levin},
       DOI = {10.4153/CJM-2017-020-5},
       URL = {https://doi.org/10.4153/CJM-2017-020-5},
}

@misc{imai2023prismatic,
  title={The prismatic realization functor for {S}himura varieties of abelian type},
  author={Imai, N. and Kato, H. and Youcis, A.},
  journal={arXiv},
  year={2023},
  note={arXiv:\href{https://arxiv.org/abs/2310.08472}{2310.08472}},
}

@misc{IKY26,
author = {Inoue, K. and Koshikawa, T. and Yao, Z.},
title = {CORRIGENDUM TO ``LOGARITHMIC PRISMATIC COHOMOLOGY {II}''},
year = {2026},
note = {preprint, available \href{https://researchmap.jp/7000017226?lang=en}{here}},
}

@misc{MW26note,
author = {Mao, S. and Wu, P.},
title = {Some examples and counterexamples in complement to this paper},
year = {2026},
note = {Not for publishing, available \href{https://peihang-wu.github.io/files/pi1-hodge.pdf}{here}},
}

@inproceedings{Inoue2025LogPD,
  title={Log Prismatic {D}ieudonn{\'e} theory and its application to {S}himura varieties},
  author={Inoue, K.},
  year={2025},
  url={https://api.semanticscholar.org/CorpusID:276929375}
}

@article {KSG21,
    AUTHOR = {Kurano, K. and Shimomoto, K.},
     TITLE = {Ideal-adic completion of quasi-excellent rings (after
              {G}abber)},
   JOURNAL = {Kyoto J. Math.},
  FJOURNAL = {Kyoto Journal of Mathematics},
    VOLUME = {61},
      YEAR = {2021},
    NUMBER = {3},
     PAGES = {707--722},
      ISSN = {2156-2261,2154-3321},
   MRCLASS = {13B35 (13F25 13F40)},
  MRNUMBER = {4301055},
MRREVIEWER = {Cristodor\ Ionescu},
       DOI = {10.1215/21562261-2021-0011},
       URL = {https://doi.org/10.1215/21562261-2021-0011},
}

@misc{RCam26,
      title={Locally analytic completed cohomology}, 
      author={Rodríguez Camargo, J. E.},
      year={2026},
      eprint={2209.01057},
      archivePrefix={arXiv},
      primaryClass={math.NT},
      note={arXiv:\href{https://arxiv.org/abs/2209.01057}{2209.01057}, to appear in JAMS}, 
}

@article {Han20,
    AUTHOR = {Hansen, D.},
     TITLE = {Vanishing and comparison theorems in rigid analytic geometry},
   JOURNAL = {Compos. Math.},
  FJOURNAL = {Compositio Mathematica},
    VOLUME = {156},
      YEAR = {2020},
    NUMBER = {2},
     PAGES = {299--324},
      ISSN = {0010-437X,1570-5846},
   MRCLASS = {14F20 (14F08 14G05 14G22)},
  MRNUMBER = {4045974},
MRREVIEWER = {Alessandra\ Bertapelle},
       DOI = {10.1112/s0010437x19007371},
       URL = {https://doi.org/10.1112/s0010437x19007371},
}

@article {KKN15,
    AUTHOR = {Kajiwara, T. and Kato, K. and Nakayama, C.},
     TITLE = {Logarithmic abelian varieties, {P}art {IV}: {P}roper models},
   JOURNAL = {Nagoya Math. J.},
  FJOURNAL = {Nagoya Mathematical Journal},
    VOLUME = {219},
      YEAR = {2015},
     PAGES = {9--63},
      ISSN = {0027-7630,2152-6842},
   MRCLASS = {14K10 (14D06 14J10)},
  MRNUMBER = {3413572},
MRREVIEWER = {Arvid\ Siqveland},
       DOI = {10.1215/00277630-3140577},
       URL = {https://doi.org/10.1215/00277630-3140577},
}

@misc{boxer2015torsion,
  title={Torsion in the coherent cohomology of Shimura varieties and Galois representations},
  author={Boxer, G.},
  eprint={1507.05922},
  note={arXiv:\href{https://arxiv.org/abs/1507.05922}{1507.05922}},
  year={2015}
}

@misc{bertapelle2023log,
  title={Log $ p $-divisible groups and semi-stable representations},
  author={Bertapelle, A. and Wang, S. and Zhao, H.},
  eprint={2302.11030},
  note={arXiv:\href{https://arxiv.org/abs/2302.11030}{2302.11030}},
  year={2023}
}

@article{wurthen2024log,
  title={Log p-divisible groups associated with log 1-motives},
  author={W{\"u}rthen, M. and Zhao, H.},
  journal={Canadian Journal of Mathematics},
  volume={76},
  number={3},
  pages={946--983},
  year={2024},
  publisher={Canadian Mathematical Society}
}

@misc{Kat23,
      title={Logarithmic {D}ieudonn{\'e} theory}, 
      author={Kato, K.},
      year={2023},
      eprint={2306.13943},
      archivePrefix={arXiv},
      primaryClass={math.AG},
      note={arXiv:\href{https://arxiv.org/abs/2306.13943}{2306.13943}}, 
}
\end{document}